\documentclass[numbers=enddot,12pt,final,onecolumn,notitlepage]{scrartcl}%
\usepackage[headsepline,footsepline,manualmark]{scrlayer-scrpage}
\usepackage[all,cmtip]{xy}
\usepackage{amssymb}
\usepackage{amsmath}
\usepackage{amsthm}
\usepackage{framed}
\usepackage{comment}
\usepackage{color}
\usepackage[breaklinks=True]{hyperref}
\usepackage[sc]{mathpazo}
\usepackage[T1]{fontenc}
\usepackage{needspace}
\usepackage{tabls}
\usepackage{tikz}
\usepackage{bbding}
\usepackage{hyperxmp}
\usepackage[type={CC}, modifier={zero}, version={1.0},]{doclicense}
\providecommand{\U}[1]{\protect\rule{.1in}{.1in}}
\usetikzlibrary{arrows.meta}
\usetikzlibrary{chains}
\usetikzlibrary{shapes}
\usetikzlibrary{decorations.pathmorphing}
\newcounter{exer}
\numberwithin{exer}{section}
\theoremstyle{definition}
\newtheorem{theo}{Theorem}[subsection]
\newenvironment{theorem}[1][]
{\begin{theo}[#1]\begin{leftbar}}
{\end{leftbar}\end{theo}}
\newtheorem{lem}[theo]{Lemma}
\newenvironment{lemma}[1][]
{\begin{lem}[#1]\begin{leftbar}}
{\end{leftbar}\end{lem}}
\newtheorem{prop}[theo]{Proposition}
\newenvironment{proposition}[1][]
{\begin{prop}[#1]\begin{leftbar}}
{\end{leftbar}\end{prop}}
\newtheorem{defi}[theo]{Definition}
\newenvironment{definition}[1][]
{\begin{defi}[#1]\begin{leftbar}}
{\end{leftbar}\end{defi}}
\newtheorem{remk}[theo]{Remark}
\newenvironment{remark}[1][]
{\begin{remk}[#1]\begin{leftbar}}
{\end{leftbar}\end{remk}}
\newtheorem{coro}[theo]{Corollary}
\newenvironment{corollary}[1][]
{\begin{coro}[#1]\begin{leftbar}}
{\end{leftbar}\end{coro}}
\newtheorem{conv}[theo]{Convention}
\newenvironment{convention}[1][]
{\begin{conv}[#1]\begin{leftbar}}
{\end{leftbar}\end{conv}}
\newtheorem{quest}[theo]{Question}
\newenvironment{question}[1][]
{\begin{quest}[#1]\begin{leftbar}}
{\end{leftbar}\end{quest}}
\newtheorem{warn}[theo]{Warning}
\newenvironment{warning}[1][]
{\begin{warn}[#1]\begin{leftbar}}
{\end{leftbar}\end{warn}}
\newtheorem{conj}[theo]{Conjecture}

\newtheorem{exam}[theo]{Example}
\newenvironment{example}[1][]
{\begin{exam}[#1]\begin{leftbar}}
{\end{leftbar}\end{exam}}
\newtheorem{exmp}[exer]{Exercise}
\newenvironment{exercise}[1][]
{\begin{exmp}[#1]\begin{leftbar}}
{\end{leftbar}\end{exmp}}
\newenvironment{statement}{\begin{quote}}{\end{quote}}
\newenvironment{fineprint}{\begin{small}}{\end{small}}
\let\sumnonlimits\sum
\let\prodnonlimits\prod
\let\cupnonlimits\bigcup
\let\capnonlimits\bigcap
\renewcommand{\sum}{\sumnonlimits\limits}
\renewcommand{\prod}{\prodnonlimits\limits}
\renewcommand{\bigcup}{\cupnonlimits\limits}
\renewcommand{\bigcap}{\capnonlimits\limits}
\newenvironment{noncompile}{}{}
\excludecomment{verlong}
\includecomment{vershort}
\excludecomment{noncompile}
\newcommand{\set}[1]{\left\{ #1 \right\}}

\newcommand{\tup}[1]{\left( #1 \right)}
\newcommand{\of}{\text{ of }}
\newtheoremstyle{plainsl}
{8pt plus 2pt minus 4pt}
{8pt plus 2pt minus 4pt}
{\slshape}
{0pt}
{\bfseries}
{.}
{5pt plus 1pt minus 1pt}
{}
\theoremstyle{plainsl}
\ihead{An introduction to graph theory, version June 7, 2025}
\ohead{page \thepage}
\begin{document}

\title{An introduction to graph theory\\{\normalsize (Text for Math 530 in Spring 2022 at Drexel University)}}
\author{Darij Grinberg\thanks{Drexel University, Korman Center, 15 S 33rd Street,
Office \#263, Philadelphia, PA 19104 (USA). //
\texttt{darijgrinberg@gmail.com} //
\texttt{http://www.cip.ifi.lmu.de/\symbol{126}grinberg/}}}
\date{Spring 2025 edition,
June 7, 2025
}
\maketitle

\begin{abstract}
\textbf{Abstract.} This is a graduate-level introduction to graph theory,
corresponding to a quarter-long course. It covers simple graphs, multigraphs
as well as their directed analogues, and more restrictive classes such as
tournaments, trees and arborescences. Among the features discussed are
Eulerian circuits, Hamiltonian cycles, spanning trees, the matrix-tree and
BEST theorems, proper colorings, Turan's theorem, bipartite matching and the
Menger and Gallai--Milgram theorems. The basics of network flows are
introduced in order to prove Hall's marriage theorem.

Around a hundred exercises are included (without solutions).

\end{abstract}
\tableofcontents

\doclicenseThis

\section{\label{chp.pref}Preface}

\subsection{What is this?}

This is a course on \textbf{graphs} -- a rather elementary concept (actually a
cluster of closely related concepts) that can be seen all over mathematics. We
will discuss several kinds of graphs (simple graphs, multigraphs, directed
graphs, etc.) and study their features and properties. In particular, we will
encounter walks on graphs, matchings of graphs, flows on networks (networks
are graphs with extra data), and take a closer look at certain types of graphs
such as trees and tournaments.

The theory of graphs goes back at least to Leonhard Euler, who in
\href{https://scholarlycommons.pacific.edu/cgi/viewcontent.cgi?article=1052&context=euler-works}{a
1736 paper} \cite{Euler36} (see \cite{EulerEn} for an
\href{https://www.imsc.res.in/~sitabhra/teaching/sb15b/ScientificAmerican_1953_Leonhard_Euler_and_the_Koenigsberg_Bridges.pdf}{English
translation}) solved a puzzle about an optimal tour of the town of
K\"{o}nigsberg. It saw some more developments in the 19th century and
straight-up exploded in the 20th; now it is one of the most active fields of
mathematics. There are now dozens (if not hundreds) textbooks available on the
subject, such as

\begin{itemize}
\item the comprehensive works \cite{BonMur08}, \cite{Berge91}, \cite{Ore74},
\cite{Bollob98},\newline\cite{Dieste17}, \cite{ChLeZh16}, \cite{Jungni13}

\item or the more introductory \cite{Ore96}, \cite[Chapters 5--6]{BenWil06},
\cite{Bollob71}, \cite{Griffi21}, \cite{Galvin}, \cite[Chapter 5]{Guicha16},
\cite{Harary69}, \cite{Harju14},\newline\cite[Chapter 1]{HaHiMo08},
\cite{Wilson10}, \cite{Tait21}, \cite[Chapters 10--13]{LeLeMe},
\cite{Merris00}, \cite[Part III]{Bona23}, \cite{Ruohon13}, \cite{KelTro17},
\cite{LoPeVe03}, \cite{West01}, \cite{Verstr21}, \cite{HarRin03},
\cite{Bickle24}.
\end{itemize}

\noindent These texts are written at different levels of sophistication, rigor
and detail, are tailored to different audiences, and (beyond the absolute
basics) often cover different ground (for instance, \cite{Dieste17}
distinguishes itself by treating infinite and random graphs, whereas
\cite{Even12} focusses on algorithms, and \cite{Griffi21} is strong on applications).

The present notes are self-contained and do not follow any existing book.
Nevertheless, I recommend skimming the texts cited above to gain a wider
perspective on graph theory (far beyond what we can cover in an introductory
course), and perhaps marking some for later reading. Our focus in these notes
is on the more discrete and algebraic sides of graph theory (finite graphs of
various kinds, existential results, counting formulas), and they are limited
both by the time constraints (being written for a quarter-long course) and the
limits of my own knowledge.

\subsubsection{Remarks}

\textbf{Prerequisites.} These notes target a graduate-level (or advanced
undergraduate) reader. A certain mathematical sophistication and willingness
to think along (as well as invent one's own examples) is expected. Beyond
that, the main prerequisites are the basic properties of determinants,
polynomials and finite sums. Rings and fields are occasionally mentioned, but
the reader can make do with just the most basic examples thereof ($\mathbb{Q}%
$, $\mathbb{R}$, polynomial rings and matrix rings; also the finite field
$\mathbb{F}_{2}$ in a few places). No analysis (or even calculus) is required
anywhere in this text.\medskip

\textbf{Course websites.} These notes were written for my Math 530 course at
Drexel University in Spring 2022. The website of this course can be found at%
\[
\text{\url{https://www.cip.ifi.lmu.de/~grinberg/t/22s} .}%
\]
An older, but similarly structured course is my Spring 2017 course at the
University of Minnesota. Its website is available at%
\[
\text{\url{https://www.cip.ifi.lmu.de/~grinberg/t/17s} ,}%
\]
and contains some additional materials (such as solutions to some selected
exercises, a few more detailed topics, and a
\href{https://www.cip.ifi.lmu.de/~grinberg/t/17s/nogra.pdf}{stub of a text}
\cite{17s} that covers parts of our Chapter \ref{chp.sg} in more depth). If
you are reading the present notes on the arXiv, then said additional materials
can also be found as ancillary files to this arXiv submission.\medskip

\textbf{Exercises.} These notes include exercises of varying difficulty and
significance. Almost all of the exercises are optional (i.e., they are not
used anywhere in the text, except perhaps in other exercises), but they often
provide practice, context and additional inspiration. Naturally, one person's
inspiration is another's distraction, so I do not recommend assigning too much
importance to any specific exercise; it is usually better to read on than to
dwell for hours. However, a dozen minutes of thought per exercise will likely
not be a waste of time. \medskip

\begin{noncompile}
Exercises from 2017 to be added to future versions:

\begin{itemize}
\item MT3 exercise 4 (Pfaffian).

\item MT3 exercises 5, 6 (acyclic orientations).

\item HW4 exercise 7 (requires solution).

\item HW5 exercises 1--6 (chipfiring, acyclic orientations).
\end{itemize}
\end{noncompile}

\textbf{Acknowledgments.} I have learned a lot from conversations with Joel
Brewster Lewis, Lukas Katth\"{a}n and Victor Reiner. Chiara Libera Carnevale,
Amanda Johnson, Victor Michel and Eric M\"{u}ller corrected errors in previous
versions of these notes. I am indebted to all of the above, and would
appreciate any further input -- please contact
\texttt{darijgrinberg@gmail.com} about any corrections (however small) and suggestions.

\subsection{\label{sec.pref.nots}Notations}

The following notations will be used throughout these notes:

\begin{itemize}
\item We let $\mathbb{N}=\left\{  0,1,2,\ldots\right\}  $. Thus,
$0\in\mathbb{N}$.

\item The size (i.e., cardinality) of a finite set $S$ is denoted by
$\left\vert S\right\vert $.

\item If $S$ is a set, then the \textbf{powerset} of $S$ means the set of all
subsets of $S$. This powerset will be denoted by $\mathcal{P}\left(  S\right)
$.

Moreover, if $S$ is a set, and $k$ is an integer, then $\mathcal{P}_{k}\left(
S\right)  $ will mean the set of all $k$-element subsets of $S$. For instance,%
\[
\mathcal{P}_{2}\left(  \left\{  1,2,3\right\}  \right)  =\left\{  \left\{
1,2\right\}  ,\ \left\{  1,3\right\}  ,\ \left\{  2,3\right\}  \right\}  .
\]

\item For any number $n$ and any $k\in\mathbb{N}$, we define the
\textbf{binomial coefficient} $\dbinom{n}{k}$ to be the number%
\[
\dfrac{n\left(  n-1\right)  \left(  n-2\right)  \cdots\left(  n-k+1\right)
}{k!}=\dfrac{\prod_{i=0}^{k-1}\left(  n-i\right)  }{k!}.
\]

These binomial coefficients have many interesting properties, which can often
be found in textbooks on enumerative combinatorics (e.g., \cite[Chapter
2]{19fco}). Some of the most important ones are the following:

\begin{itemize}
\item The factorial formula: If $n,k\in\mathbb{N}$ and $n\geq k$, then
$\dbinom{n}{k}=\dfrac{n!}{k!\cdot\left(  n-k\right)  !}$.

\item The combinatorial interpretation: If $n,k\in\mathbb{N}$, and if $S$ is
an $n$-element set, then $\dbinom{n}{k}$ is the number of all $k$-element
subsets of $S$ (in other words, $\left\vert \mathcal{P}_{k}\left(  S\right)
\right\vert =\dbinom{n}{k}$).

\item Pascal's recursion: For any number $n$ and any positive integer $k$, we
have%
\[
\dbinom{n}{k}=\dbinom{n-1}{k-1}+\dbinom{n-1}{k}.
\]

\end{itemize}
\end{itemize}

\section{\label{chp.sg}Simple graphs}

\subsection{\label{sec.sg.defs}Definitions}

The first type of graphs that we will consider are the \textquotedblleft
simple graphs\textquotedblright, named so because of their very simple definition:

\begin{definition}
\label{def.sg.sg} A \textbf{simple graph} is a pair $\left(  V,E\right)  $,
where $V$ is a finite set, and where $E$ is a subset of $\mathcal{P}%
_{2}\left(  V\right)  $.
\end{definition}

To remind, $\mathcal{P}_{2}\left(  V\right)  $ is the set of all $2$-element
subsets of $V$. Thus, a simple graph is a pair $\left(  V,E\right)  $, where
$V$ is a finite set, and $E$ is a set consisting of $2$-element subsets of
$V$. We will abbreviate the word \textquotedblleft simple
graph\textquotedblright\ as \textquotedblleft graph\textquotedblright\ in this
chapter, but later (in Chapter \ref{chp.mg}) we will learn some more advanced
and general notions of \textquotedblleft graphs\textquotedblright.

\begin{example}
\label{exa.sg.sg1}Here is a simple graph:%
\[
\left(  \left\{  1,2,3,4\right\}  ,\ \ \left\{  \left\{  1,3\right\}
,\ \ \left\{  1,4\right\}  ,\ \ \left\{  3,4\right\}  \right\}  \right)  .
\]

\end{example}

\begin{example}
\label{exa.sg.sg2}For any $n\in\mathbb{N}$, we can define a simple graph
$\operatorname*{Cop}\nolimits_{n}$ to be the pair $\left(  V,E\right)  $,
where $V=\left\{  1,2,\ldots,n\right\}  $ and
\[
E=\left\{  \left\{  u,v\right\}  \in\mathcal{P}_{2}\left(  V\right)
\ \mid\ \gcd\left(  u,v\right)  =1\right\}  .
\]
We call this the $n$\textbf{-th coprimality graph}.
\end{example}

(Some authors do not require $V$ to be finite in Definition \ref{def.sg.sg};
this leads to \textbf{infinite graphs}. But I shall leave this can of worms
closed for this quarter.) \medskip

The purpose of simple graphs is to encode relations on a finite set --
specifically the kind of relations that are binary (i.e., relate pairs of
elements), symmetric (i.e., mutual) and irreflexive (i.e., an element cannot
be related to itself). For example, the graph $\operatorname*{Cop}%
\nolimits_{n}$ in Example \ref{exa.sg.sg2} encodes the coprimality (aka
coprimeness) relation on the set $\left\{  1,2,\ldots,n\right\}  $, except
that the latter relation is not irreflexive ($1$ is coprime to $1$, but
$\left\{  1,1\right\}  $ is not in $E$; thus, the graph $\operatorname*{Cop}%
\nolimits_{n}$ \textquotedblleft forgets\textquotedblright\ that $1$ is
coprime to $1$). For another example, if $V$ is a set of people, and $E$ is
the set of $\left\{  u,v\right\}  \in\mathcal{P}_{2}\left(  V\right)  $ such
that $u$ has been married to $v$ at some point, then $\left(  V,E\right)  $ is
a simple graph. Even in 2022, marriage to oneself is not a thing, so all
marriages can be encoded as $2$-element subsets.\footnote{The more standard
example for a social graph would be a \textquotedblleft friendship
graph\textquotedblright; here, $V$ is again a set of people, but $E$ is now
the set of $\left\{  u,v\right\}  \in\mathcal{P}_{2}\left(  V\right)  $ such
that $u$ and $v$ are friends. Of course, this only works if you think of
friendship as being automatically mutual (true for facebook friendship,
questionable for the actual thing).} \medskip

The following notations provide a quick way to reference the elements of $V$
and $E$ when given a graph $\left(  V, E \right)  $:

\begin{definition}
\label{def.sg.VE} Let $G=\left(  V,E\right)  $ be a simple graph.

\begin{enumerate}
\item[\textbf{(a)}] The set $V$ is called the \textbf{vertex set} of $G$; it
is denoted by $\operatorname{V}\left(  G\right)  $. (Notice that the letter
\textquotedblleft$\operatorname{V}$\textquotedblright\ in \textquotedblleft%
$\operatorname{V}\left(  G\right)  $\textquotedblright\ is upright, as opposed
to the letter \textquotedblleft$V$\textquotedblright\ in \textquotedblleft%
$\left(  V,E\right)  $\textquotedblright, which is italic. These are two
different symbols, and have different meanings: The letter $V$ stands for the
specific set $V$ which is the first component of the pair $G$, whereas the
letter $\operatorname{V}$ is part of the notation $\operatorname{V}\left(
G\right)  $ for the vertex set of any graph. Thus, if $H=\left(  W,F\right)  $
is another graph, then $\operatorname{V}\left(  H\right)  $ is $W$, not $V$.)

The elements of $V$ are called the \textbf{vertices} (or the \textbf{nodes})
of $G$.

\item[\textbf{(b)}] The set $E$ is called the \textbf{edge set} of $G$; it is
denoted by $\operatorname{E}\left(  G\right)  $. (Again, the letter
\textquotedblleft$\operatorname{E}$\textquotedblright\ in \textquotedblleft%
$\operatorname{E}\left(  G\right)  $\textquotedblright\ is upright, and stands
for a different thing than the \textquotedblleft$E$\textquotedblright.)

The elements of $E$ are called the \textbf{edges} of $G$. When $u$ and $v$ are
two elements of $V$, we shall often use the notation $uv$ for $\left\{
u,v\right\}  $; thus, each edge of $G$ has the form $uv$ for two distinct
elements $u$ and $v$ of $V$. Of course, we always have $uv=vu$.

Notice that each simple graph $G$ satisfies $G=\left(  \operatorname{V}\left(
G\right)  ,\operatorname{E}\left(  G\right)  \right)  $.

\item[\textbf{(c)}] Two vertices $u$ and $v$ of $G$ are said to be
\textbf{adjacent} (to each other) if $uv\in E$ (that is, if $uv$ is an edge of
$G$). In this case, the edge $uv$ is said to \textbf{join} $u$ with $v$ (or
\textbf{connect} $u$ and $v$); the vertices $u$ and $v$ are called the
\textbf{endpoints} of this edge. When the graph $G$ is not obvious from the
context, we shall often say \textquotedblleft adjacent in $G$%
\textquotedblright\ instead of just \textquotedblleft
adjacent\textquotedblright.

Two vertices $u$ and $v$ of $G$ are said to be \textbf{non-adjacent} (to each
other) if they are not adjacent (i.e., if $uv\notin E$).

\item[\textbf{(d)}] Let $v$ be a vertex of $G$ (that is, $v\in V$). Then, the
\textbf{neighbors} of $v$ (in $G$) are the vertices $u$ of $G$ that satisfy
$vu\in E$. In other words, the \textbf{neighbors} of $v$ are the vertices of
$G$ that are adjacent to $v$.
\end{enumerate}
\end{definition}

\begin{example}
Let $G$ be the simple graph
\[
\left(  \left\{  1,2,3,4\right\}  ,\ \ \left\{  \left\{  1,3\right\}
,\ \ \left\{  1,4\right\}  ,\ \ \left\{  3,4\right\}  \right\}  \right)
\]
from Example \ref{exa.sg.sg1}. Then, its vertex set and its edge set are%
\[
\operatorname*{V}\left(  G\right)  =\left\{  1,2,3,4\right\}
\ \ \ \ \ \text{and}\ \ \ \ \ \operatorname*{E}\left(  G\right)  =\left\{
\left\{  1,3\right\}  ,\ \ \left\{  1,4\right\}  ,\ \ \left\{  3,4\right\}
\right\}  =\left\{  13,\ 14,\ 34\right\}
\]
(using our notation $uv$ for $\left\{  u,v\right\}  $). The vertices $1$ and
$3$ are adjacent (since $13\in\operatorname*{E}\left(  G\right)  $), but the
vertices $1$ and $2$ are not (since $12\notin\operatorname*{E}\left(
G\right)  $). The neighbors of $1$ are $3$ and $4$. The endpoints of the edge
$34$ are $3$ and $4$.
\end{example}

\subsection{\label{sect.intro.draw}Drawing graphs}

There is a common method to represent graphs visually: Namely, a graph can be
drawn as a set of points in the plane and a set of curves connecting some of
these points with each other.

More precisely:

\begin{definition}
\label{def.intro.draw} A simple graph $G$ can be visually represented by
\textbf{drawing} it on the plane. To do so, we represent each vertex of $G$ by
a point (at which we put the name of the vertex), and then, for each edge $uv$
of $G$, we draw a curve that connects the point representing $u$ with the
point representing $v$. The positions of the points and the shapes of the
curves can be chosen freely, as long as they allow the reader to unambiguously
reconstruct the graph $G$ from the picture. (Thus, for example, the curves
should not pass through any points other than the ones they mean to connect.)
\end{definition}

\begin{example}
\label{exa.intro.draw} Let us draw some simple graphs. \medskip

\textbf{(a)} The simple graph $\left(  \left\{  1,2,3\right\}  ,\ \left\{
12,23\right\}  \right)  $ (where we are again using the shorthand notation
$uv$ for $\left\{  u,v\right\}  $) can be drawn as follows:
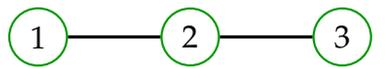
\[%
%
\ \ .
\]
This is (in a sense) the simplest way to draw this graph: The edges are
represented by straight lines. But we can draw it in several other ways as
well -- e.g., as follows:%
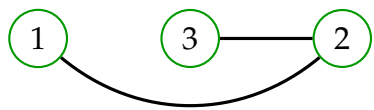
\[%
%
\ \ .
\]
Here, we have placed the points representing the vertices $1,2,3$ differently.
As a consequence, we were not able to draw the edge $12$ as a straight line,
because it would then have overlapped with the vertex $3$, which would make
the graph ambiguous (the edge $12$ could be mistaken for two edges $13$ and
$32$).

Here are three further drawings of the same graph $\left(  \left\{
1,2,3\right\}  ,\ \left\{  12,23\right\}  \right)  $:
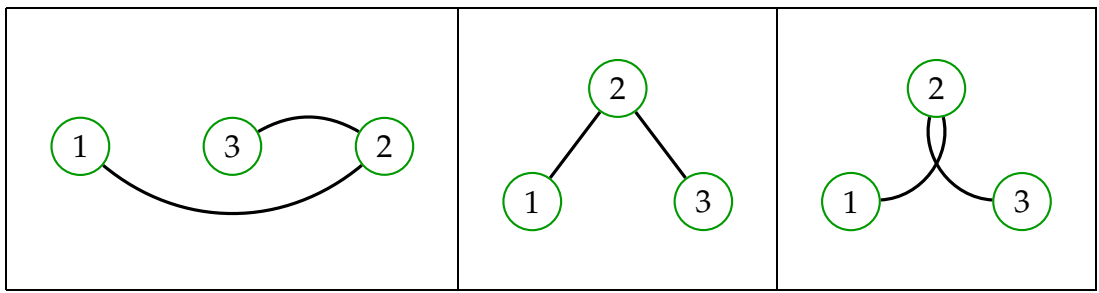
\[%
\begin{tabular}
[c]{|ccc|}\hline
\multicolumn{1}{|c|}{} & \multicolumn{1}{c|}{} & \\
\hline
\end{tabular}
\ \ \ \ .
\]

\textbf{(b)} Consider the $5$-th coprimality graph $\operatorname*{Cop}%
\nolimits_{5}$ defined in Example \ref{exa.sg.sg2}. Here is one way to draw
it:%
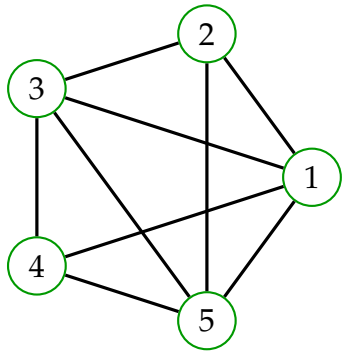
\[%
%
\ \ .
\]
Here is another way to draw the same graph $\operatorname*{Cop}\nolimits_{5}$,
with fewer intersections between edges:
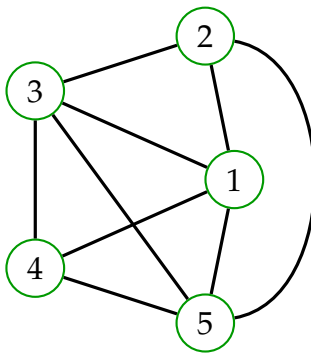
\[%
%
\ \ .
\]
By appropriately repositioning the points corresponding to the five vertices
of $\operatorname*{Cop}\nolimits_{5}$, we can actually get rid of all
intersections and make all the edges straight (as opposed to curved). Can you
find out how?\medskip

\textbf{(c)} Let us draw one further graph: the simple graph $\left(  \left\{
1,2,3,4,5\right\}  ,\ \mathcal{P}_{2}\left(  {\left\{  1,2,3,4,5\right\}
}\right)  \right)  $. This is the simple graph whose vertices are $1,2,3,4,5$,
and whose edges are all possible two-element sets consisting of its vertices
(i.e., each pair of two distinct vertices is adjacent). We shall later call
this graph the \textquotedblleft complete graph $K_{5}$\textquotedblright.
Here is a simple way to draw this graph:
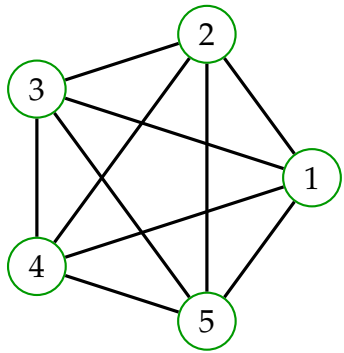
\[%
%
\ \ .
\]
This drawing is useful for many purposes; for example, it makes the abstract
symmetry of this graph (i.e., the fact that, roughly speaking, its vertices
$1,2,3,4,5$ are \textquotedblleft equal in rights\textquotedblright) obvious.
But sometimes, you might want to draw it differently, to minimize the number
of intersecting curves. Here is a drawing with fewer intersections:
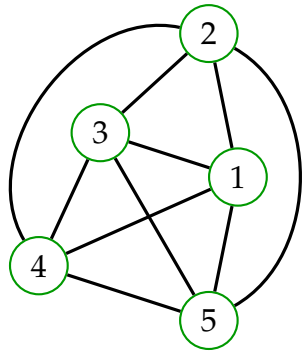
\[%
%
\ \ .
\]
In this drawing, we have only one intersection between two curves left. Can we
get rid of all intersections?

This is a question of topology, not of combinatorics, since it really is about
curves in the plane rather than about finite sets and graphs. The answer is
\textquotedblleft no\textquotedblright. (That is, no matter how you draw this
graph in the plane, you will always have at least one pair of curves
intersect.) This is a classical result (one of the first theorems in the
theory of
\textbf{\href{https://en.wikipedia.org/wiki/Planar_graph}{\textbf{planar
graphs}}}), and proofs of it can be found in various textbooks (e.g.,
\cite[Theorem 4.1.2]{FriFri98}, which is generally a good introduction to
planar graph theory even if it uses terminology somewhat different from ours).
Note that any proof must use some analysis or topology, since the result
relies on the notion of a (continuous) curve in the plane (if curves were
allowed to be non-continuous, then they could \textquotedblleft jump
over\textquotedblright\ one another, so they could easily avoid intersecting!).
\end{example}

\begin{noncompile}
What follows is an OLD version of Example \ref{exa.intro.draw}, with all
graphs drawn in xymatrix:

\begin{example}
Let us draw some simple graphs. \medskip

\textbf{(a)} The simple graph $\left(  \left\{  1,2,3\right\}  ,\ \left\{
12,23\right\}  \right)  $ (where we are again using the shorthand notation
$uv$ for $\left\{  u,v\right\}  $) can be drawn as follows:
\[
\xymatrix{ 1 \are[r] & 2 \are[r] & 3 }.
\]
This is (in a sense) the simplest way to draw this graph: The edges are
represented by straight lines. But we can draw it in several other ways as
well -- e.g., as follows:%
\begin{align*}
&  \xymatrix{ 1 \are@/_2pc/[rr] & 3 \are[r] & 2 }.\\
&  \phantom{a}
\end{align*}
Here, we have placed the points representing the vertices $1,2,3$ differently.
As a consequence, we were not able to draw the edge $12$ as a straight line,
because it would then have overlapped with the vertex $3$, which would make
the graph ambiguous (the edge $12$ could be mistaken for two edges $13$ and
$32$).

Here are three further drawings of the same graph $\left(  \left\{
1,2,3\right\}  ,\ \left\{  12,23\right\}  \right)  $:
\[%
\begin{tabular}
[c]{|c|c|c|}\hline
\ \ \ $\xymatrix{ 1 \are@/_2pc/[rr] & 3 \are@/^1pc/[r] & 2 }\ \ \ $ &
\ \ \ $\xymatrix{ & 2 \are[dl] \are[dr] \\ 1 & & 3 }\ \ \ $ &
\ \ \ $\xymatrix{ & 2 \are@/^2pc/[dl] \are@/_2pc/[dr] \\ 1 & & 3 }\ \ \ $%
\\\hline
\end{tabular}
\ \ \ .
\]

\textbf{(b)} Consider the $5$-th coprimality graph $\operatorname*{Cop}%
\nolimits_{5}$ defined in Example \ref{exa.sg.sg2}. Here is one way to draw
it:%
\[
\xymatrix{ & 2 \are[rr] \are[dl] \are[rdd] & & 3 \are[llld] \are[rd] \are[ldd] \\ 1 \are[rrrr] \are[rrd] & & & & 4 \are[lld] \\ & & 5 }
\]
Here is another way to draw the same graph $G$, with fewer intersections
between edges:
\[
\xymatrix{ & 2 \are[rr] \are[d] \are@/_6pc/[rdd] & & 3 \are[lld] \are[rd] \are[ldd] \\ & 1 \are[rrr] \are[rd] & & & 4 \are[lld] \\ & & 5 }.
\]
By appropriately repositioning the points corresponding to the five vertices
of $G$, we can actually get rid of all intersections and make all the edges
straight (as opposed to curved). Can you find out how?\medskip

\textbf{(c)} Let us draw one further graph: the simple graph $\left(  \left\{
1,2,3,4,5\right\}  ,\mathcal{P}_{2}\left(  {\left\{  1,2,3,4,5\right\}
}\right)  \right)  $. This is the simple graph whose vertices are $1,2,3,4,5$,
and whose edges are all possible two-element sets consisting of its vertices
(i.e., each pair of two distinct vertices is adjacent). We shall later call
this graph the \textquotedblleft complete graph $K_{5}$\textquotedblright.
Here is a simple way to draw this graph:
\[
\xymatrix{ & 2 \are[ld] \are[rr] \are[rrrd] \are[rdd] & & 3 \are[rd] \are[ldd] \are[llld] \\ 1 \are[rrrr] \are[rrd] & & & & 4 \are[lld] \\ & & 5 }
\]
This drawing is useful for many purposes; for example, it makes the abstract
symmetry of this graph (i.e., the fact that, roughly speaking, its vertices
$1,2,3,4,5$ are \textquotedblleft equal in rights\textquotedblright) obvious.
But sometimes, you might want to draw it differently, to minimize the number
of intersecting curves. Here is a drawing with fewer intersections:
\begin{align*}
&
\xymatrix{ & 2 \are[rr] \are@/^4pc/[rrrd] \are@/_6pc/[rrdd] \are[dd] & & 3 \are[rd] \are[lldd] \are[dd] \\ & & & & 4 \are[llld] \\ & 5 \are[rr] & & 1 \are[ru] }.\\
&  \phantom{a}
\end{align*}
In this drawing, we have only one intersection between two curves left. Can we
get rid of all intersections?

This is a question of topology, not of combinatorics, since it really is about
curves in the plane rather than about finite sets and graphs. The answer is
\textquotedblleft no\textquotedblright. (That is, no matter how you draw this
graph in the plane, you will always have at least one pair of curves
intersect.) This is a classical result (one of the first theorems in the
theory of \textbf{planar graphs}), and proofs of it can be found in various
textbooks (e.g., \cite[Theorem 4.1.2]{FriFri98}, which is generally a good
introduction to planar graph theory even if it uses terminology somewhat
different from ours). Note that any proof must use some analysis or topology,
since the result relies on the notion of a (continuous) curve in the plane (if
curves were allowed to be non-continuous, then they could \textquotedblleft
jump over\textquotedblright\ one another, so they could easily avoid intersecting!).
\end{example}
\end{noncompile}

\subsection{\label{sect.intro.R33}A first fact: The Ramsey number $R\left(
3,3 \right)  = 6$}

Enough definitions; let's state a first result:

\begin{proposition}
\label{prop.simple.R33} Let $G$ be a simple graph with $\left|
\operatorname{V}\left(  G \right)  \right|  \geq6$ (that is, $G$ has at least
$6$ vertices). Then, at least one of the following two statements holds:

\begin{itemize}
\item \textit{Statement 1:} There exist three distinct vertices $a$, $b$ and
$c$ of $G$ such that $ab$, $bc$ and $ca$ are edges of $G$.

\item \textit{Statement 2:} There exist three distinct vertices $a$, $b$ and
$c$ of $G$ such that none of $ab$, $bc$ and $ca$ is an edge of $G$.
\end{itemize}
\end{proposition}

In other words, Proposition~\ref{prop.simple.R33} says that if a graph $G$ has
at least $6$ vertices, then we can either find three distinct vertices that
are mutually adjacent\footnote{by which we mean (of course) that any two
\textbf{distinct} ones among these three vertices are adjacent} or find three
distinct vertices that are mutually non-adjacent (i.e., no two of them are
adjacent), or both. Often, this is restated as follows: \textquotedblleft In
any group of at least six people, you can always find three that are
(pairwise) friends to each other, or three no two of whom are
friends\textquotedblright\ (provided that friendship is a symmetric relation).

We will give some examples in a moment, but first let us introduce some
convenient terminology:

\begin{definition}
\label{def.sg.triangle} Let $G$ be a simple graph.

\begin{enumerate}
\item[\textbf{(a)}] A set $\left\{  a,b,c\right\}  $ of three distinct
vertices of $G$ is said to be a \textbf{triangle} (of $G$) if every two
distinct vertices in this set are adjacent (i.e., if $ab$, $bc$ and $ca$ are
edges of $G$).

\item[\textbf{(b)}] A set $\left\{  a,b,c\right\}  $ of three distinct
vertices of $G$ is said to be an \textbf{anti-triangle} (of $G$) if no two
distinct vertices in this set are adjacent (i.e., if none of $ab$, $bc$ and
$ca$ is an edge of $G$).
\end{enumerate}
\end{definition}

Thus, Proposition~\ref{prop.simple.R33} says that every simple graph with at
least $6$ vertices contains a triangle or an anti-triangle (or both).

\begin{example}
\label{exa.simple.R33} Let us show two examples of graphs $G$ to which
Proposition~\ref{prop.simple.R33} applies, as well as an example to which it
does not:

\begin{enumerate}
\item[\textbf{(a)}] Let $G$ be the graph $\left(  V,E\right)  $, where
\begin{align*}
V  &  =\left\{  1,2,3,4,5,6\right\}  \qquad\text{and}\\
E  &  =\left\{  \left\{  1,2\right\}  ,\left\{  2,3\right\}  ,\left\{
3,4\right\}  ,\left\{  4,5\right\}  ,\left\{  5,6\right\}  ,\left\{
6,1\right\}  \right\}  .
\end{align*}
(This graph can be drawn in such a way as to look like a hexagon:%
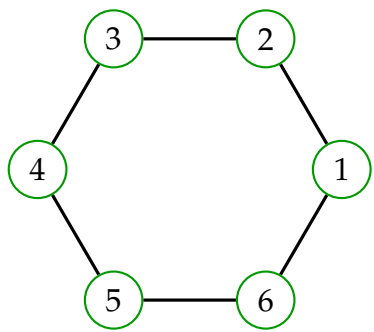
\[%
\begin{tikzpicture}
\begin{scope}[every node/.style={circle,thick,draw=green!60!black}]
\node(A) at (0:2) {$1$};
\node(B) at (60:2) {$2$};
\node(C) at (120:2) {$3$};
\node(D) at (180:2) {$4$};
\node(E) at (240:2) {$5$};
\node(F) at (300:2) {$6$};
\end{scope}
\begin{scope}[every edge/.style={draw=black,very thick}]
\path
[-] (A) edge (B) (B) edge (C) (C) edge (D) (D) edge (E) (E) edge (F) (F) edge (A);
\end{scope}
\end{tikzpicture}%
\ \ .
\]
) This graph satisfies Proposition~\ref{prop.simple.R33}, since $\left\{
1,3,5\right\}  $ is an anti-triangle (or since $\left\{  2,4,6\right\}  $ is
an anti-triangle).

\item[\textbf{(b)}] Let $G$ be the graph $\left(  V,E\right)  $, where
\begin{align*}
V  &  =\left\{  1,2,3,4,5,6\right\}  \qquad\text{and}\\
E  &  =\left\{  \left\{  1,2\right\}  ,\left\{  2,3\right\}  ,\left\{
3,4\right\}  ,\left\{  4,5\right\}  ,\left\{  5,6\right\}  ,\left\{
6,1\right\}  ,\left\{  1,3\right\}  ,\left\{  4,6\right\}  \right\}  .
\end{align*}
(This graph can be drawn in such a way as to look like a hexagon with two
extra diagonals:
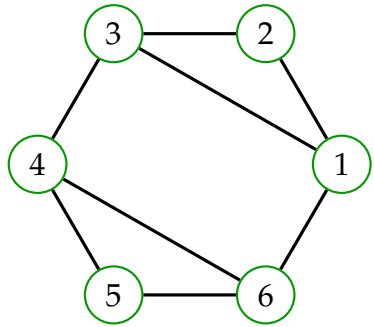
\[%
%
\ \ .
\]
) This graph satisfies Proposition~\ref{prop.simple.R33}, since $\left\{
1,2,3\right\}  $ is a triangle.

\item[\textbf{(c)}] Let $G$ be the graph $\left(  V,E\right)  $, where
\begin{align*}
V  &  =\left\{  1,2,3,4,5\right\}  \qquad\text{and}\\
E  &  =\left\{  \left\{  1,2\right\}  ,\left\{  2,3\right\}  ,\left\{
3,4\right\}  ,\left\{  4,5\right\}  ,\left\{  5,1\right\}  \right\}  .
\end{align*}
(This graph can be drawn to look like a pentagon:
\[%
%
\ \ .
\]
) Proposition~\ref{prop.simple.R33} says nothing about this graph, since this
graph does not satisfy the assumption of Proposition~\ref{prop.simple.R33} (in
fact, its number of vertices $\left\vert \operatorname{V}\left(  G\right)
\right\vert $ fails to be $\geq6$). By itself, this does not yield that the
claim of Proposition~\ref{prop.simple.R33} is false for this graph. However,
it is easy to check that the claim actually \textbf{is} false for this graph:
It has neither a triangle nor an anti-triangle.
\end{enumerate}
\end{example}

\begin{noncompile}
Here is an OLD version of Example \ref{exa.simple.R33}, with graphs drawn
using xymatrix:

\begin{example}
Let us show two examples of graphs $G$ to which
Proposition~\ref{prop.simple.R33} applies, as well as an example to which it
does not:

\begin{enumerate}
\item[\textbf{(a)}] Let $G$ be the graph $\left(  V,E\right)  $, where
\begin{align*}
V  &  =\left\{  1,2,3,4,5,6\right\}  \qquad\text{and}\\
E  &  =\left\{  \left\{  1,2\right\}  ,\left\{  2,3\right\}  ,\left\{
3,4\right\}  ,\left\{  4,5\right\}  ,\left\{  5,6\right\}  ,\left\{
6,1\right\}  \right\}  .
\end{align*}
(This graph can be drawn in such a way as to look like a hexagon:
\[
\xymatrix{ & 1 \are[r] & 2 \are[rd] \\ 6 \are[ru] & & & 3 \are[ld] \\ & 5 \are[lu] & 4 \are[l] }.
\]
) This graph satisfies Proposition~\ref{prop.simple.R33}, since $\left\{
1,3,5\right\}  $ is an anti-triangle (or since $\left\{  2,4,6\right\}  $ is
an anti-triangle).

\item[\textbf{(b)}] Let $G$ be the graph $\left(  V,E\right)  $, where
\begin{align*}
V  &  =\left\{  1,2,3,4,5,6\right\}  \qquad\text{and}\\
E  &  =\left\{  \left\{  1,2\right\}  ,\left\{  2,3\right\}  ,\left\{
3,4\right\}  ,\left\{  4,5\right\}  ,\left\{  5,6\right\}  ,\left\{
6,1\right\}  ,\left\{  1,3\right\}  ,\left\{  4,6\right\}  \right\}  .
\end{align*}
(This graph can be drawn in such a way as to look like a hexagon with two
extra diagonals:
\[
\xymatrix{ & 1 \are[r] \are[rrd] & 2 \are[rd] \\ 6 \are[ru] & & & 3 \are[ld] \\ & 5 \are[lu] & 4 \are[llu] \are[l] }.
\]
) This graph satisfies Proposition~\ref{prop.simple.R33}, since $\left\{
1,2,3\right\}  $ is a triangle.

\item[\textbf{(c)}] Let $G$ be the graph $\left(  V,E\right)  $, where
\begin{align*}
V  &  =\left\{  1,2,3,4,5\right\}  \qquad\text{and}\\
E  &  =\left\{  \left\{  1,2\right\}  ,\left\{  2,3\right\}  ,\left\{
3,4\right\}  ,\left\{  4,5\right\}  ,\left\{  5,1\right\}  \right\}  .
\end{align*}
(This graph can be drawn to look like a pentagon:
\[
\xymatrix{ & 2 \are[rr] & & 3 \are[rd] \\ 1 \are[ru] \are[rrd] & & & & 4 \are[lld] \\ & & 5 }
\]
) Proposition~\ref{prop.simple.R33} says nothing about this graph, since this
graph does not satisfy the assumption of Proposition~\ref{prop.simple.R33} (in
fact, its number of vertices $\left\vert \operatorname{V}\left(  G\right)
\right\vert $ fails to be $\geq6$). By itself, this does not yield that the
claim of Proposition~\ref{prop.simple.R33} is false for this graph. However,
it is easy to check that the claim actually \textbf{is} false for this graph:
It has neither a triangle nor an anti-triangle.
\end{enumerate}
\end{example}
\end{noncompile}

\begin{proof}
[Proof of Proposition~\ref{prop.simple.R33}.]We need to prove that $G$ has a
triangle or an anti-triangle (or both).

Choose any vertex $u\in\operatorname{V}\left(  G\right)  $. (This is clearly
possible, since $\left\vert \operatorname{V}\left(  G\right)  \right\vert
\geq6\geq1$.) Then, there are at least $5$ vertices distinct from $u$ (since
$G$ has at least $6$ vertices). We are in one of the following two cases:

\textit{Case 1:} The vertex $u$ has at least $3$ neighbors.

\textit{Case 2:} The vertex $u$ has at most $2$ neighbors.

Let us consider Case 1 first. In this case, the vertex $u$ has at least $3$
neighbors. Hence, we can find three distinct neighbors $p$, $q$ and $r$ of
$u$. Consider these $p$, $q$ and $r$. If one (or more) of $pq$, $qr$ and $rp$
is an edge of $G$, then $G$ has a triangle (for example, if $pq$ is an edge of
$G$, then $\left\{  u,p,q\right\}  $ is a triangle). If not, then $G$ has an
anti-triangle (namely, $\left\{  p,q,r\right\}  $). Thus, in either case, our
proof is complete in Case 1.

Let us now consider Case 2. In this case, the vertex $u$ has at most $2$
neighbors. Hence, the vertex $u$ has at least $3$ non-neighbors\footnote{The
word \textquotedblleft non-neighbor\textquotedblright\ shall here mean a
vertex that is not adjacent to $u$ and \textbf{distinct from }$u$. Thus, $u$
does not count as a non-neighbor of $u$.} (since there are at least $5$
vertices distinct from $u$ in total). Thus, we can find three distinct
non-neighbors $p$, $q$ and $r$ of $u$. Consider these $p$, $q$ and $r$. If all
of $pq$, $qr$ and $rp$ are edges of $G$, then $G$ has a triangle (namely,
$\left\{  p,q,r\right\}  $). If not, then $G$ has an anti-triangle (for
example, if $pq$ is not an edge of $G$, then $\left\{  u,p,q\right\}  $ is an
anti-triangle). In either case, we are thus done with the proof in Case 2.
Thus, both cases are resolved, and the proof is complete.
\end{proof}

Notice the symmetry between Case 1 and Case 2 in our above proof: the
arguments used were almost the same, except that neighbors and non-neighbors
swapped roles.

\begin{remark}
Proposition~\ref{prop.simple.R33} could also be proved by brute force (using a
computer). Indeed, it clearly suffices to prove it for all simple graphs with
$6$ vertices (as opposed to $\geq6$ vertices), because if a graph has more
than $6$ vertices, then we can just throw away some of them until we have only
$6$ left. However, there are only finitely many simple graphs with $6$
vertices (up to relabeling of their vertices), and the validity of Proposition
\ref{prop.simple.R33} can be checked for each of them. This is, of course,
cumbersome (even a computer would take a moment checking all the $2^{15}$
possible graphs for triangles and anti-triangles) and unenlightening.
\end{remark}

Proposition~\ref{prop.simple.R33} is the first result in a field of graph
theory known as \textbf{Ramsey theory}. I shall not dwell on this field in
this course, but let me make a few more remarks. The first step beyond
Proposition~\ref{prop.simple.R33} is the following generalization (see, e.g.,
\cite[problem 14.11]{Tomesc85} for a proof):

\begin{proposition}
\label{prop.simple.Rrs} Let $r$ and $s$ be two positive integers. Let $G$ be a
simple graph with $\left|  \operatorname{V}\left(  G \right)  \right|
\geq\dbinom{r+s-2}{r-1}$. Then, at least one of the following two statements holds:

\begin{itemize}
\item \textit{Statement 1:} There exist $r$ distinct vertices of $G$ that are
mutually adjacent (i.e., each two distinct ones among these $r$ vertices are adjacent).

\item \textit{Statement 2:} There exist $s$ distinct vertices of $G$ that are
mutually non-adjacent (i.e., no two distinct ones among these $s$ vertices are adjacent).
\end{itemize}
\end{proposition}

Applying Proposition~\ref{prop.simple.Rrs} to $r=3$ and $s=3$, we can recover
Proposition~\ref{prop.simple.R33}.

One might wonder whether the number $\dbinom{r+s-2}{r-1}$ in
Proposition~\ref{prop.simple.Rrs} can be improved -- i.e., whether we can
replace it by a smaller number without making
Proposition~\ref{prop.simple.Rrs} false. In the case of $r=3$ and $s=3$, this
is impossible, because the number $6$ in Proposition~\ref{prop.simple.R33}
cannot be made smaller\footnote{Indeed, we saw in Example~\ref{exa.simple.R33}
\textbf{(c)} that $5$ vertices would not suffice.}. However, for some other
values of $r$ and $s$, the value $\dbinom{r+s-2}{r-1}$ can be improved. (For
example, for $r=4$ and $s=4$, the best possible value is $18$ rather than
$\dbinom{4+4-2}{4-1}=20$.) The smallest possible value that could stand in
place of $\dbinom{r+s-2}{r-1}$ in Proposition~\ref{prop.simple.Rrs} is called
the \textbf{Ramsey number} $R\left(  r,s\right)  $; thus, we have just showed
that $R\left(  3,3\right)  =6$. Finding $R\left(  r,s\right)  $ for higher
values of $r$ and $s$ is a hard computational challenge; here are some values
that have been found with the help of computers:
\begin{align*}
R\left(  3,4\right)   &  =9;\qquad R\left(  3,5\right)  =14;\qquad R\left(
3,6\right)  =18;\qquad R\left(  3,7\right)  =23;\\
R\left(  3,8\right)   &  =28;\qquad R\left(  3,9\right)  =36;\qquad R\left(
4,4\right)  =18;\qquad R\left(  4,5\right)  =25.
\end{align*}
(We are only considering the cases $r\leq s$, since it is easy to see that
$R\left(  r,s\right)  =R\left(  s,r\right)  $ for all $r$ and $s$. Also, the
trivial values $R\left(  1,s\right)  =1$ and $R\left(  2,s\right)  =s+1$ for
$s\geq2$ are omitted.) The Ramsey number $R\left(  5,5\right)  $ is still
unknown (although it is known that $43\leq R\left(  5,5\right)  \leq46$; see
\cite{AngMcK24} for details).

Proposition~\ref{prop.simple.Rrs} can be further generalized to a result
called \textit{Ramsey's theorem}. The idea behind the generalization is to
slightly change the point of view, and replace the simple graph $G$ by a
complete graph (i.e., a simple graph in which every two distinct vertices are
adjacent) whose edges are colored in two colors (say, blue and red). This is a
completely equivalent concept, because the concepts of \textquotedblleft
adjacent\textquotedblright\ and \textquotedblleft
non-adjacent\textquotedblright\ in $G$ can be identified with the concepts of
\textquotedblleft adjacent through a blue edge\textquotedblright\ (i.e., the
edge connecting them is colored blue) and \textquotedblleft adjacent through a
red edge\textquotedblright, respectively. Statements 1 and 2 then turn into
\textquotedblleft there exist $r$ distinct vertices that are mutually adjacent
through blue edges\textquotedblright\ and \textquotedblleft there exist $s$
distinct vertices that are mutually adjacent through red
edges\textquotedblright, respectively. From this point of view, it is only
logical to generalize Proposition~\ref{prop.simple.Rrs} further to the case
when the edges of a complete graph are colored in $k$ (rather than two)
colors. The corresponding generalization is known as Ramsey's theorem. We
refer to the well-written Wikipedia page
\url{https://en.wikipedia.org/wiki/Ramsey's_theorem} for a treatment of this
generalization with proof, as well as a table of known Ramsey numbers
$R\left(  r,s\right)  $ and a self-contained (if somewhat terse) proof of
Proposition~\ref{prop.simple.Rrs}. Ramsey's theorem can be generalized and
varied further; this usually goes under the name \textquotedblleft Ramsey
theory\textquotedblright. For elementary introductions, see the Cut-the-knot
page
\url{http://www.cut-the-knot.org/Curriculum/Combinatorics/ThreeOrThree.shtml}
, the above-mentioned Wikipedia article, as well as the texts by Harju
\cite{Harju14}, Bollobas \cite{Bollob98} and West \cite{West01}. Several other
texts on combinatorics have a chapter on Ramsey theory, e.g., \cite[Chapter
14]{Tomesc85}. \medskip

There is one more direction in which Proposition~\ref{prop.simple.R33} can be
improved a bit: A graph $G$ with at least $6$ vertices has not only one
triangle or anti-triangle, but at least two of them (this can include having
one triangle and one anti-triangle). Proving this makes for a nice exercise:

\begin{exercise}
\label{exa.simple.R33.two} Let $G$ be a simple graph. A
\textbf{triangle-or-anti-triangle} in $G$ means a set that is either a
triangle or an anti-triangle.

\begin{enumerate}
\item[\textbf{(a)}] Assume that $\left\vert \operatorname{V}\left(  G\right)
\right\vert \geq6$. Prove that $G$ has at least two
triangle-or-anti-triangles. (For comparison: Proposition~\ref{prop.simple.R33}
shows that $G$ has at least one triangle-or-anti-triangle.)

\item[\textbf{(b)}] Assume that $\left\vert \operatorname{V}\left(  G\right)
\right\vert =m+6$ for some $m\in\mathbb{N}$. Prove that $G$ has at least $m+1$
triangle-or-anti-triangles. \medskip
\end{enumerate}

[\textbf{Solution:} This is Exercise 1 on homework set \#1 from my Spring 2017
course; see \href{https://www.cip.ifi.lmu.de/~grinberg/t/17s/}{the course
page} for solutions.]
\end{exercise}

\subsection{\label{sec.sg.deg}Degrees}

\subsubsection{Definition and basic properties}

The \textbf{degree} of a vertex in a simple graph just counts how many edges
contain this vertex:

\begin{definition}
\label{def.sg.deg}Let $G=\left(  V,E\right)  $ be a simple graph. Let $v\in V$
be a vertex. Then, the \textbf{degree} of $v$ (with respect to $G$) is defined
to be%
\begin{align*}
\deg v:=  &  \ \left(  \text{the number of edges }e\in E\text{ that contain
}v\right) \\
=  &  \ \left(  \text{the number of neighbors of }v\right) \\
=  &  \ \left\vert \left\{  u\in V\ \mid\ uv\in E\right\}  \right\vert \\
=  &  \ \left\vert \left\{  e\in E\ \mid\ v\in e\right\}  \right\vert .
\end{align*}

\end{definition}

(These equalities are pretty easy to check: Each edge $e\in E$ that contains
$v$ contains exactly one neighbor of $v$, and conversely, each neighbor of $v$
belongs to exactly one edge that contains $v$. However, these equalities are
specific to simple graphs, and won't hold any more once we move on to multigraphs.)

For example, in the graph%
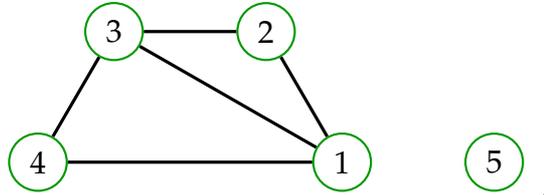
\[%
\begin{tikzpicture}
\begin{scope}[every node/.style={circle,thick,draw=green!60!black}]
\node(A) at (0:2) {$1$};
\node(B) at (60:2) {$2$};
\node(C) at (120:2) {$3$};
\node(D) at (180:2) {$4$};
\node(E) at (4,0) {$5$};
\end{scope}
\begin{scope}[every edge/.style={draw=black,very thick}]
\path[-] (A) edge (B) (B) edge (C) (C) edge (A) (C) edge (D) (D) edge (A);
\end{scope}
\end{tikzpicture}%
\ \ ,
\]
the vertices have degrees%
\[
\deg1=3,\ \ \ \ \ \ \ \ \ \ \deg2=2,\ \ \ \ \ \ \ \ \ \ \deg
3=3,\ \ \ \ \ \ \ \ \ \ \deg4=2,\ \ \ \ \ \ \ \ \ \ \deg5=0.
\]

Here are some basic properties of degrees in simple graphs:

\begin{proposition}
\label{prop.degv.lessn}Let $G$ be a simple graph with $n$ vertices. Let $v$ be
a vertex of $G$. Then,%
\[
\deg v\in\left\{  0,1,\ldots,n-1\right\}  .
\]

\end{proposition}

\begin{proof}
All neighbors of $v$ belong to the $\left(  n-1\right)  $-element set
$\operatorname*{V}\left(  G\right)  \setminus\left\{  v\right\}  $. Thus,
their number is $\leq n-1$.
\end{proof}

\begin{proposition}
[Euler 1736]\label{prop.deg.euler}Let $G$ be a simple graph. Then, the sum of
the degrees of all vertices of $G$ equals twice the number of edges of $G$. In
other words,%
\[
\sum_{v\in\operatorname*{V}\left(  G\right)  }\deg v=2\cdot\left\vert
\operatorname*{E}\left(  G\right)  \right\vert .
\]

\end{proposition}

\begin{proof}
Write the simple graph $G$ as $G=\left(  V,E\right)  $; thus,
$\operatorname*{V}\left(  G\right)  =V$ and $\operatorname*{E}\left(
G\right)  =E$.

Now, let $N$ be the number of all pairs $\left(  v,e\right)  \in V\times E$
such that $v\in e$. We compute $N$ in two different ways (this is called
\textquotedblleft double-counting\textquotedblright):

\begin{enumerate}
\item We can obtain $N$ by computing, for each $v\in V$, the number of all
$e\in E$ that satisfy $v\in e$, and then summing these numbers over all $v$.
Since these numbers are just the degrees $\deg v$, the result will be
$\sum\limits_{v\in V}\deg v$.

\item On the other hand, we can obtain $N$ by computing, for each $e\in E$,
the number of all $v\in V$ that satisfy $v\in e$, and summing these numbers
over all $e$. Since each $e\in E$ contains exactly $2$ vertices $v\in V$, this
result will be $\sum\limits_{e\in E}2=\left\vert E\right\vert \cdot
2=2\cdot\left\vert E\right\vert $.
\end{enumerate}

Since these two results must be equal (because they both equal $N$), we thus
see that $\sum\limits_{v\in V}\deg v=2\cdot\left\vert E\right\vert $. But this
is the claim of Proposition \ref{prop.deg.euler}.
\end{proof}

\begin{corollary}
[handshake lemma]\label{cor.deg.odd-deg-even}Let $G$ be a simple graph. Then,
the number of vertices $v$ of $G$ whose degree $\deg v$ is odd is even.
\end{corollary}

\begin{proof}
Proposition \ref{prop.deg.euler} yields that $\sum\limits_{v\in
\operatorname*{V}\left(  G\right)  }\deg v=2\cdot\left\vert \operatorname*{E}%
\left(  G\right)  \right\vert $. Hence, $\sum\limits_{v\in\operatorname*{V}%
\left(  G\right)  }\deg v$ is even. However, if a sum of integers is even,
then it must have an even number of odd addends. Thus, the sum $\sum
\limits_{v\in\operatorname*{V}\left(  G\right)  }\deg v$ must have an even
number of odd addends. In other words, the number of vertices $v$ of $G$ whose
degree $\deg v$ is odd is even.
\end{proof}

Corollary \ref{cor.deg.odd-deg-even} is often stated as follows: In a group of
people, the number of persons with an odd number of friends (in the group) is
even. It is also known as the \textbf{handshake lemma}.

Here is another property of degrees in a simple graph:

\begin{proposition}
\label{prop.sg.deg.equal-degs}Let $G$ be a simple graph with at least two
vertices. Then, there exist two distinct vertices $v$ and $w$ of $G$ that have
the same degree.
\end{proposition}

\begin{proof}
Assume the contrary. So the degrees of all $n$ vertices of $G$ are distinct,
where $n=\left\vert \operatorname*{V}\left(  G\right)  \right\vert $.

In other words, the map%
\begin{align*}
\deg:\operatorname*{V}\left(  G\right)   &  \rightarrow\left\{  0,1,\ldots
,n-1\right\}  ,\\
v  &  \mapsto\deg v
\end{align*}
is injective. But this is a map between two finite sets of the same size
($n$). When such a map is injective, it has to be bijective (by the pigeonhole
principle). Therefore, in particular, it takes both $0$ and $n-1$ as values.

In other words, there are a vertex $u$ with degree $0$ and a vertex $v$ with
degree $n-1$. Are these two vertices adjacent or not? Yes because of $\deg
v=n-1$; no because of $\deg u=0$. Contradiction!

(Fine print: The two vertices $u$ and $v$ must be distinct, since $0\neq n-1$.
It is here that we are using the \textquotedblleft at least two
vertices\textquotedblright\ assumption!)
\end{proof}

\subsubsection{\label{subsec.sg.deg.mantel}Mantel's theorem}

Here is an application of counting neighbors to proving a fact about graphs.
This is known as \textbf{Mantel's theorem}:

\begin{theorem}
[Mantel's theorem]\label{thm.sg.mantel}Let $G$ be a simple graph with $n$
vertices and $e$ edges. Assume that $e>n^{2}/4$. Then, $G$ has a triangle
(i.e., three distinct vertices that are pairwise adjacent).
\end{theorem}

\begin{example}
\label{exa.sg.mantel.ex1}Let $G$ be the graph $\left(  V,E\right)  $, where%
\begin{align*}
V  &  =\left\{  1,2,3,4,5,6\right\}  ;\\
E  &  =\left\{  12,\ 23,\ 34,\ 45,\ 56,\ 61,\ 14,\ 25,\ 36\right\}  .
\end{align*}
Here is a drawing:
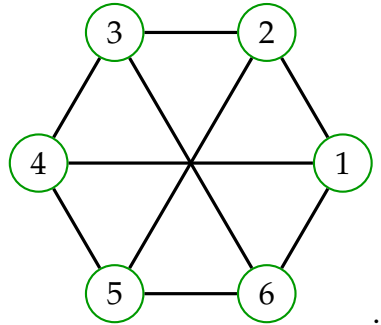
\[%
\begin{tikzpicture}
\begin{scope}[every node/.style={circle,thick,draw=green!60!black}]
\node(1) at (0:2) {$1$};
\node(2) at (60:2) {$2$};
\node(3) at (120:2) {$3$};
\node(4) at (180:2) {$4$};
\node(5) at (240:2) {$5$};
\node(6) at (300:2) {$6$};
\end{scope}
\begin{scope}[every edge/.style={draw=black,very thick}]
\path[-] (1) edge (2) (2) edge (3) (3) edge (4) (4) edge (5);
\path[-] (5) edge (6) (6) edge (1);
\path[-] (1) edge (4) (2) edge (5) (3) edge (6);
\end{scope}
\end{tikzpicture}%
.
\]
This graph has no triangle (which, by the way, is easy to verify without
checking all possibilities: just observe that every edge of $G$ joins two
vertices of different parity, but a triangle would necessarily have two
vertices of equal parity). Thus, by the contrapositive of Mantel's theorem, it
satisfies $e\leq n^{2}/4$ with $n=6$ and $e=9$. This is indeed true because
$9=6^{2}/4$. But this also entails that if we add any further edge to $G$,
then we obtain a triangle.
\end{example}

\begin{noncompile}
Old version of the picture:
$\xymatrix{ & 1 \are[r] \are[rdd] & 2 \are[rd] \are[ldd] \\ 6 \are[rrr] \are[ru] & & & 3 \are[ld] \\ & 5 \are[lu] & 4 \are[l] }$%

\end{noncompile}

\begin{proof}
[Proof of Mantel's theorem.]We will prove the theorem by strong induction on
$n$. Thus, we assume (as the induction hypothesis) that the theorem holds for
all graphs with fewer than $n$ vertices. We must now prove it for our graph
$G$ with $n$ vertices. Let $V=\operatorname*{V}\left(  G\right)  $ and
$E=\operatorname*{E}\left(  G\right)  $, so that $G=\left(  V,E\right)  $.

We must prove that $G$ has a triangle. Assume the contrary. Thus, $G$ has no triangle.

From $e>n^{2}/4\geq0$, we see that $G$ has an edge. Pick any such edge, and
call it $vw$. Thus, $v\neq w$.

Let us now color each edge of $G$ with one of three colors, as follows:

\begin{itemize}
\item The edge $vw$ is colored black.

\item Each edge that contains exactly one of $v$ and $w$ is colored red.

\item All other edges are colored blue.
\end{itemize}

The following picture shows an example of this coloring:%
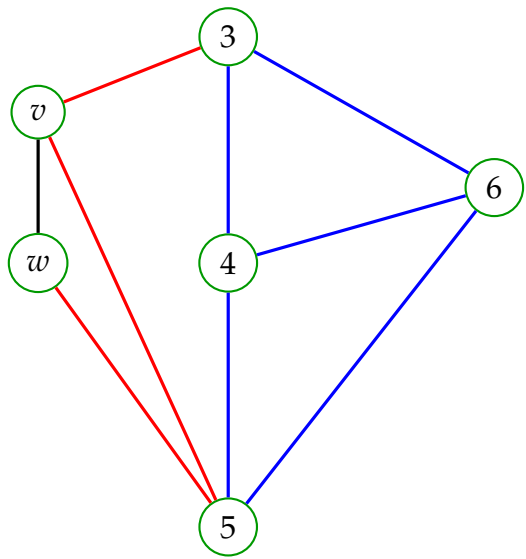
\[%
\begin{tikzpicture}
\begin{scope}[every node/.style={circle,thick,draw=green!60!black}]
\node(A) at (0,1) {$w$};
\node(B) at (0,3) {$v$};
\node(C) at (2.5,4) {$3$};
\node(D) at (2.5,1) {$4$};
\node(E) at (2.5,-2.5) {$5$};
\node(F) at (6,2) {$6$} ;
\end{scope}
\begin{scope}[every edge/.style={draw=black,very thick}]
\path[-] (A) edge[draw=black] (B);
\path[-] (B) edge[draw=red] (C);
\path[-] (D) edge[draw=blue] (C);
\path[-] (A) edge[draw=red] (E);
\path[-] (D) edge[draw=blue] (E);
\path[-] (C) edge[draw=blue] (F);
\path[-] (E) edge[draw=blue] (F);
\path[-] (B) edge[draw=red] (E);
\path[-] (D) edge[draw=blue] (F);
\end{scope}
\end{tikzpicture}%
\ \ .
\]

\begin{noncompile}
(If you are wondering why this doesn't look like our previous graphs: that's
because I'm using tikz rather than xymatrix to draw it. In general, tikz is
probably the better tool, but I'm more used to xymatrix.)
\end{noncompile}

We now count the edges of each color:

\begin{itemize}
\item There is exactly $1$ black edge -- namely, $vw$.

\item How many red edges can there be? I claim that there are at most $n-2$.
Indeed, each vertex other than $v$ and $w$ is connected to at most one of $v$
and $w$ by a red edge, since otherwise it would form a triangle with $v$ and
$w$.

\item How many blue edges can there be? The vertices other than $v$ and $w$,
along with the blue edges that join them, form a graph with $n-2$ vertices;
this graph has no triangles (since $G$ has no triangles). By the induction
hypothesis, however, if this graph had more than $\left(  n-2\right)  ^{2}/4$
edges, then it would have a triangle. Thus, it has $\leq\left(  n-2\right)
^{2}/4$ edges. In other words, there are $\leq\left(  n-2\right)  ^{2}/4$ blue edges.
\end{itemize}

In total, the number of edges is therefore%
\[
\leq1+\left(  n-2\right)  +\left(  n-2\right)  ^{2}/4=n^{2}/4.
\]
In other words, $e\leq n^{2}/4$. This contradicts $e>n^{2}/4$. This is the
contradiction we were looking for, so the induction is complete.
\end{proof}

Quick question: What about equality? Can a graph with $n$ vertices and exactly
$n^{2}/4$ edges have no triangles? Yes (for even $n$). Indeed, for any even
$n$, we can take the graph%
\[
\left(  \left\{  1,2,\ldots,n\right\}  ,\ \left\{  ij\ \mid\ i\not \equiv
j\operatorname{mod}2\right\}  \right)
\]
(keep in mind that $ij$ means the $2$-element set $\left\{  i,j\right\}  $
here, not the product $i\cdot j$). We can also do this for odd $n$, and obtain
a graph with $\left(  n^{2}-1\right)  /4$ edges (which is as close to
$n^{2}/4$ as we can get when $n$ is odd -- after all, the number of edges has
to be an integer). So the bound in Mantel's theorem is optimal (as far as
integers are concerned). \medskip

The following exercise can be regarded as a \textquotedblleft mirror
version\textquotedblright\ of Mantel's theorem:

\begin{exercise}
\label{exe.intro.mantel-co}Let $G$ be a simple graph with $n$ vertices and $e$
edges. Assume that $e<n\left(  n-2\right)  /4$. Prove that $G$ has an
anti-triangle (i.e., three distinct vertices that are pairwise
non-adjacent).\medskip

[\textbf{Solution:} This is Exercise 2 on homework set \#1 from my Spring 2017
course; see \href{https://www.cip.ifi.lmu.de/~grinberg/t/17s/}{the course
page} for solutions.]
\end{exercise}

\bigskip

Mantel's theorem can be generalized:

\begin{theorem}
[Turan's theorem]\label{thm.sg.turan}Let $r$ be a positive integer. Let $G$ be
a simple graph with $n$ vertices and $e$ edges. Assume that%
\[
e>\dfrac{r-1}{r}\cdot\dfrac{n^{2}}{2}.
\]
Then, there exist $r+1$ distinct vertices of $G$ that are mutually adjacent.
\end{theorem}

Mantel's theorem is the particular case for $r=2$. We will see a proof of
Turan's theorem later (Theorem \ref{thm.indep.turan}). Mantel's and Turan's
theorems are two of the simplest results of \textbf{extremal graph theory} --
the study of how inequalities between some graph parameters (in our case: the
numbers of vertices and edges) imply the existence of certain substructures
(in our case: of a triangle or of $r+1$ mutually adjacent vertices). Deeper
introductions to this subject can be found in \cite[Chapters 1 and 5]{Zhao23}
and \cite{Jukna11}.

\begin{exercise}
\label{exe.1.1}Let $G=\left(  V,E\right)  $ be a simple graph. Set
$n=\left\vert V\right\vert $. Prove that we can find some edges $e_{1}%
,e_{2},\ldots,e_{k}$ of $G$ and some triangles $t_{1},t_{2},\ldots,t_{\ell}$
of $G$ such that $k+\ell\leq n^{2}/4$ and such that each edge $e\in
E\setminus\left\{  e_{1},e_{2},\ldots,e_{k}\right\}  $ is a subset of (at
least) one of the triangles $t_{1},t_{2},\ldots,t_{\ell}$. \medskip

[\textbf{Remark:} In other words, this exercise is claiming that all edges of
$G$ can be covered by at most $n^{2}/4$ edge-or-triangles. Here, an
\textbf{edge-or-triangle} means either an edge or a triangle of $G$, and the
word \textquotedblleft covers\textquotedblright\ means that each edge of $G$
is a subset of the chosen edge-or-triangles.] \medskip

[\textbf{Hint:} Imitate the above proof of Mantel's theorem.]
\end{exercise}

\begin{remark}
Exercise \ref{exe.1.1} is a generalization of Mantel's theorem. Indeed, if the
simple graph $G=\left(  V,E\right)  $ has no triangles, then the number $\ell$
in Exercise \ref{exe.1.1} must be $0$, and thus the edges $e_{1},e_{2}%
,\ldots,e_{k}$ must be all edges of $G$, so that we conclude that $\left\vert
E\right\vert =k\leq k+\ell\leq n^{2}/4$.
\end{remark}

\begin{exercise}
\label{exe.1.2}Let $G$ be a simple graph with $n$ vertices and $k$ edges,
where $n>0$. Prove that $G$ has at least $\dfrac{k}{3n}\left(  4k-n^{2}%
\right)  $ triangles.\medskip

[\textbf{Hint:} First argue that for any edge $vw$ of $G$, the total number of
triangles that contain $v$ and $w$ is at least $\deg v + \deg w - n$. Then,
use the inequality $n \left(  a_{1}^{2} + a_{2}^{2} + \cdots+ a_{n}^{2}
\right)  \geq\left(  a_{1} + a_{2} + \cdots+ a_{n} \right)  ^{2}$, which holds
for any $n$ real numbers $a_{1}, a_{2}, \ldots, a_{n}$. (This is a particular
case of the Cauchy--Schwarz inequality or the Chebyshev inequality or the
Jensen inequality -- pick your favorite!)]
\end{exercise}

\begin{remark}
Exercise \ref{exe.1.2} is known as the \textbf{Moon--Moser inequality for
triangles}. It, too, generalizes Mantel's theorem: If $k > n^{2} / 4$, then
$\dfrac{k}{3n} \left(  4k-n^{2} \right)  > 0$, and therefore Exercise
\ref{exe.1.2} entails that $G$ has at least one triangle.
\end{remark}

\begin{exercise}
\label{exe.1.7}Let $G=\left(  V,E\right)  $ be a simple graph.

An edge $e=\left\{  u,v\right\}  $ of $G$ will be called \textbf{odd} if the
number $\deg u+\deg v$ is odd.

Prove that the number of odd edges of $G$ is even.\medskip

[\textbf{Hint:} There are several solutions. One uses
\href{https://en.wikipedia.org/wiki/Modular_arithmetic}{modular arithmetic}
and (in particular) the congruence $m^{2} \equiv m \mod 2$ for every integer
$m$. Other solutions use nothing but common sense.]
\end{exercise}

\begin{exercise}
\label{exe.2.6}Let $G=\left(  V,E\right)  $ be a simple graph. Let $S$ be a
subset of $V$, and let $k=\left\vert S\right\vert $. Prove that
\[
\sum_{v\in S}\deg v\leq k\left(  k-1\right)  +\sum_{v\in V\setminus S}%
\min\left\{  \deg v,k\right\}  .
\]

\end{exercise}

\begin{remark}
\label{rmk.erdos-gallai-easy.eg}Exercise \ref{exe.2.6} has a converse (the
so-called
\textbf{\href{https://en.wikipedia.org/wiki/Erdos-Gallai_theorem}{Erd\"{o}s--Gallai
theorem}}): If $d_{1},d_{2},\ldots,d_{n}$ are $n$ nonnegative integers such
that $d_{1}+d_{2}+\cdots+d_{n}$ is even and such that $d_{1}\geq d_{2}%
\geq\cdots\geq d_{n}$ and such that each $k\in\left\{  1,2,\ldots,n\right\}  $
satisfies%
\[
\sum_{i=1}^{k}d_{i}\leq k\left(  k-1\right)  +\sum_{i=k+1}^{n}\min\left\{
d_{i},k\right\}  ,
\]
then there exists a simple graph with vertex set $\left\{  1,2,\ldots
,n\right\}  $ whose vertices have degrees $d_{1},d_{2},\ldots,d_{n}$.
\end{remark}

\bigskip

\subsection{\label{sec.sg.iso}Graph isomorphism}

Two graphs can be distinct and yet \textquotedblleft the same up to the names
of their vertices\textquotedblright: for instance,
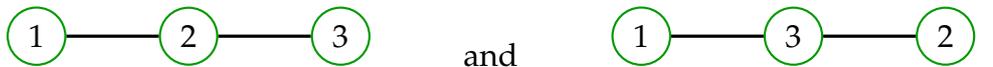
\[%
\begin{tikzpicture}[scale=2]
\begin{scope}[every node/.style={circle,thick,draw=green!60!black}]
\node(1) at (0,0) {$1$};
\node(2) at (1,0) {$2$};
\node(3) at (2,0) {$3$};
\end{scope}
\begin{scope}[every edge/.style={draw=black,very thick}]
\path[-] (1) edge (2) (2) edge (3);
\end{scope}
\end{tikzpicture}%
\ \ \ \ \ \ \ \ \ \ \text{and}\ \ \ \ \ \ \ \ \ \
\begin{tikzpicture}[scale=2]
\begin{scope}[every node/.style={circle,thick,draw=green!60!black}]
\node(1) at (0,0) {$1$};
\node(3) at (1,0) {$3$};
\node(2) at (2,0) {$2$};
\end{scope}
\begin{scope}[every edge/.style={draw=black,very thick}]
\path[-] (1) edge (3) (3) edge (2);
\end{scope}
\end{tikzpicture}%
\ \ .
\]
Let us formalize this:

\begin{definition}
Let $G$ and $H$ be two simple graphs.

\begin{enumerate}
\item[\textbf{(a)}] A \textbf{graph isomorphism} (or \textbf{isomorphism})
from $G$ to $H$ means a bijection $\phi:\operatorname*{V}\left(  G\right)
\rightarrow\operatorname*{V}\left(  H\right)  $ that \textquotedblleft
preserves edges\textquotedblright, i.e., that has the following property: For
any two vertices $u$ and $v$ of $G$, we have%
\[
\left(  uv\in\operatorname*{E}\left(  G\right)  \right)  \Longleftrightarrow
\left(  \phi\left(  u\right)  \phi\left(  v\right)  \in\operatorname*{E}%
\left(  H\right)  \right)  .
\]

\item[\textbf{(b)}] We say that $G$ and $H$ are \textbf{isomorphic} (this is
written $G\cong H$) if there exists a graph isomorphism from $G$ to $H$.
\end{enumerate}
\end{definition}

Here are two examples:

\begin{itemize}
\item The two graphs%
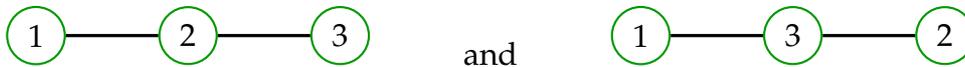
\[%
%
\]
are isomorphic, because the bijection between their vertex sets that sends
$1,2,3$ to $1,3,2$ is an isomorphism. Another isomorphism between the same two
graphs sends $1,2,3$ to $2,3,1$.

\item The two graphs%
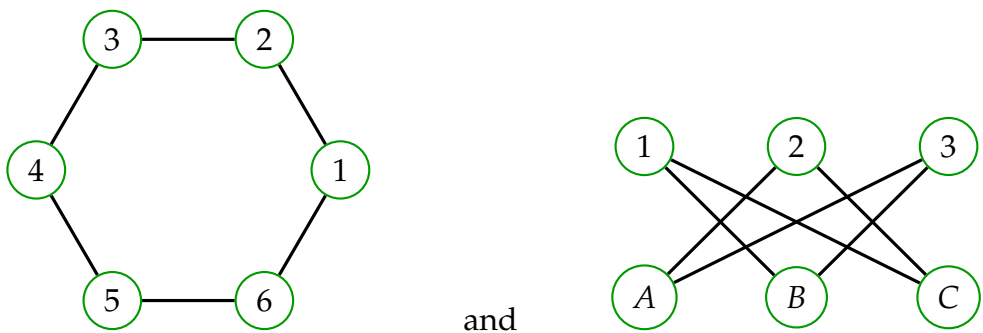
\[%
%
\]
are isomorphic, because the bijection between their vertex sets that sends
$1,2,3,4,5,6$ to $1,B,3,A,2,C$ is an isomorphism.
\end{itemize}

Here are some basic properties of isomorphisms (the proofs are straightforward):

\begin{proposition}
Let $G$ and $H$ be two graphs. The inverse of a graph isomorphism $\phi$ from
$G$ to $H$ is a graph isomorphism from $H$ to $G$.
\end{proposition}

\begin{proposition}
Let $G$, $H$ and $I$ be three graphs. If $\phi$ is a graph isomorphism from
$G$ to $H$, and $\psi$ is a graph isomorphism from $H$ to $I$, then $\psi
\circ\phi$ is a graph isomorphism from $G$ to $I$.
\end{proposition}

As a consequence of these two propositions, it is easy to see that the
relation $\cong$ (on the class of all graphs) is an equivalence relation.

Graph isomorphisms preserve all \textquotedblleft intrinsic\textquotedblright%
\ properties of a graph. For example:

\begin{proposition}
Let $G$ and $H$ be two simple graphs, and $\phi$ a graph isomorphism from $G$
to $H$. Then:

\begin{enumerate}
\item[\textbf{(a)}] For every $v\in\operatorname*{V}\left(  G\right)  $, we
have $\deg_{G}v=\deg_{H}\left(  \phi\left(  v\right)  \right)  $. Here,
$\deg_{G}v$ means the degree of $v$ as a vertex of $G$, whereas $\deg
_{H}\left(  \phi\left(  v\right)  \right)  $ means the degree of $\phi\left(
v\right)  $ as a vertex of $H$.

\item[\textbf{(b)}] We have $\left\vert \operatorname*{E}\left(  H\right)
\right\vert =\left\vert \operatorname*{E}\left(  G\right)  \right\vert $.

\item[\textbf{(c)}] We have $\left\vert \operatorname*{V}\left(  H\right)
\right\vert =\left\vert \operatorname*{V}\left(  G\right)  \right\vert $.
\end{enumerate}
\end{proposition}

One use of graph isomorphisms is to relabel the vertices of a graph. For
example, we can relabel the vertices of an $n$-vertex graph as $1,2,\ldots,n$,
or as any other $n$ distinct objects:

\begin{proposition}
\label{prop.intro.iso.rename} Let $G$ be a simple graph. Let $S$ be a finite
set such that $\left\vert S\right\vert =\left\vert \operatorname*{V}\left(
G\right)  \right\vert $. Then, there exists a simple graph $H$ that is
isomorphic to $G$ and has vertex set $\operatorname*{V}\left(  H\right)  =S$.
\end{proposition}

\begin{proof}
Straightforward.
\end{proof}

\subsection{\label{sec.sg.complete}Some families of graphs}

We will now define some particularly significant families of graphs.

\subsubsection{Complete and empty graphs}

The simplest families of graphs are the complete graphs and the empty graphs:

\begin{definition}
\label{def.sg.complete-empty}Let $V$ be a finite set.

\begin{enumerate}
\item[\textbf{(a)}] The \textbf{complete graph} on $V$ means the simple graph
$\left(  V,\ \mathcal{P}_{2}\left(  V\right)  \right)  $. It is the simple
graph with vertex set $V$ in which every two distinct vertices are adjacent.

If $V=\left\{  1,2,\ldots,n\right\}  $ for some $n\in\mathbb{N}$, then the
complete graph on $V$ is denoted $K_{n}$.

\item[\textbf{(b)}] The \textbf{empty graph} on $V$ means the simple graph
$\left(  V,\ \varnothing\right)  $. It is the simple graph with vertex set $V$
and no edges.
\end{enumerate}
\end{definition}

The following pictures show the complete graph and the empty graph on the set
$\left\{  1,2,3,4,5\right\}  $:%

\[%
\begin{tabular}
[c]{|c|c|}\hline
complete graph & empty graph\\\hline
$\ \
\begin{tikzpicture}
\begin{scope}[every node/.style={circle,thick,draw=green!60!black}]
\node(A) at (0:2) {$1$};
\node(B) at (360/5:2) {$2$};
\node(C) at (2*360/5:2) {$3$};
\node(D) at (3*360/5:2) {$4$};
\node(E) at (4*360/5:2) {$5$};
\end{scope}
\begin{scope}[every edge/.style={draw=black,very thick}]
\path
[-] (A) edge (B) (A) edge (C) (A) edge (D) (A) edge (E) (B) edge (C) (B) edge (D) (B) edge (E) (C) edge (D) (C) edge (E) (D) edge (E);
\end{scope}
\end{tikzpicture}%
\ \ $ & $\ \
\begin{tikzpicture}
\begin{scope}[every node/.style={circle,thick,draw=green!60!black}]
\node(A) at (0:2) {$1$};
\node(B) at (360/5:2) {$2$};
\node(C) at (2*360/5:2) {$3$};
\node(D) at (3*360/5:2) {$4$};
\node(E) at (4*360/5:2) {$5$};
\end{scope}
\end{tikzpicture}%
\ \ $\\\hline
\end{tabular}
\
\]
The complete one is called $K_{5}$.

Here are the complete graphs $K_{0},K_{1},K_{2},K_{3},K_{4}$:%
\[%
\begin{tabular}
[c]{|c|c|c|c|c|}\hline
$K_{0}$ & $K_{1}$ & $K_{2}$ & $K_{3}$ & $K_{4}$\\
\hline
\end{tabular}
\
\]

Note that a simple graph $G$ is isomorphic to the complete graph $K_{n}$ if
and only if it has $n$ vertices and is a complete graph (i.e., every two
distinct vertices are adjacent).\medskip

\textbf{Question:} Given two finite sets $V$ and $W$, what are the
isomorphisms from the complete graph on $V$ to the complete graph on $W$ ?

\textbf{Answer:} If $\left\vert V\right\vert \neq\left\vert W\right\vert $,
then there are none. If $\left\vert V\right\vert =\left\vert W\right\vert $,
then any bijection from $V$ to $W$ is an isomorphism. The same holds for empty graphs.

\subsubsection{Path and cycle graphs}

Next come two families of graphs with fairly simple shapes:

\begin{definition}
\label{def.sg.path}For each $n\in\mathbb{N}$, we define the $n$\textbf{-th
path graph }$P_{n}$ to be the simple graph%
\begin{align*}
&  \left(  \left\{  1,2,\ldots,n\right\}  ,\ \ \left\{  \left\{
i,i+1\right\}  \ \mid\ 1\leq i<n\right\}  \right) \\
&  =\left(  \left\{  1,2,\ldots,n\right\}  ,\ \ \left\{  12,\ 23,\ 34,\ \ldots
,\ \left(  n-1\right)  n\right\}  \right)  .
\end{align*}
This graph has $n$ vertices and $n-1$ edges (unless $n=0$, in which case it
has $0$ edges).
\end{definition}

\begin{definition}
\label{def.sg.cycle}For each $n>1$, we define the $n$\textbf{-th cycle graph
}$C_{n}$ to be the simple graph%
\begin{align*}
&  \left(  \left\{  1,2,\ldots,n\right\}  ,\ \ \left\{  \left\{
i,i+1\right\}  \ \mid\ 1\leq i<n\right\}  \cup\left\{  \left\{  n,1\right\}
\right\}  \right) \\
&  =\left(  \left\{  1,2,\ldots,n\right\}  ,\ \ \left\{  12,\ 23,\ 34,\ \ldots
,\ \left(  n-1\right)  n,\ n1\right\}  \right)  .
\end{align*}
This graph has $n$ vertices and $n$ edges (unless $n=2$, in which case it has
$1$ edge only). (We will later modify the definition of the $2$-nd cycle graph
$C_{2}$ somewhat, in order to force it to have $2$ edges. But we cannot do
this yet, since a simple graph with $2$ vertices cannot have $2$ edges.)
\end{definition}

The following pictures show the path graph $P_{5}$ and the cycle graph $C_{5}$:%

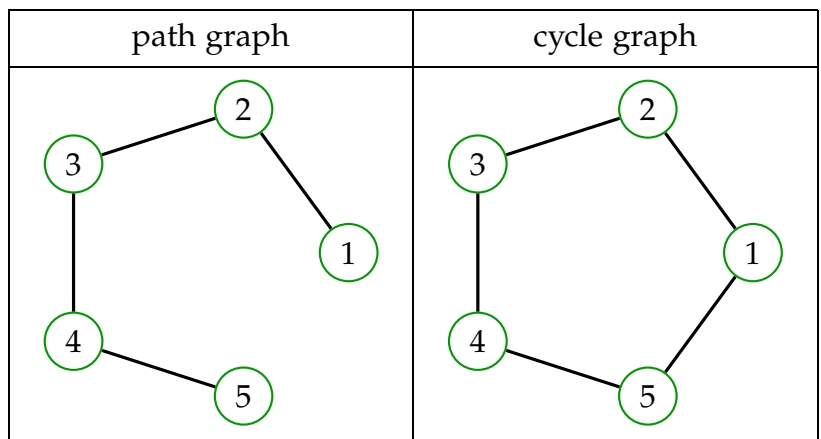
\[%
\begin{tabular}
[c]{|c|c|}\hline
path graph & cycle graph\\\hline
$\ \
\begin{tikzpicture}
\begin{scope}[every node/.style={circle,thick,draw=green!60!black}]
\node(A) at (0:2) {$1$};
\node(B) at (360/5:2) {$2$};
\node(C) at (2*360/5:2) {$3$};
\node(D) at (3*360/5:2) {$4$};
\node(E) at (4*360/5:2) {$5$};
\end{scope}
\begin{scope}[every edge/.style={draw=black,very thick}]
\path[-] (A) edge (B) (B) edge (C) (C) edge (D) (D) edge (E);
\end{scope}
\end{tikzpicture}%
\ \ $ & $\ \
\begin{tikzpicture}
\begin{scope}[every node/.style={circle,thick,draw=green!60!black}]
\node(A) at (0:2) {$1$};
\node(B) at (360/5:2) {$2$};
\node(C) at (2*360/5:2) {$3$};
\node(D) at (3*360/5:2) {$4$};
\node(E) at (4*360/5:2) {$5$};
\end{scope}
\begin{scope}[every edge/.style={draw=black,very thick}]
\path[-] (A) edge (B) (B) edge (C) (C) edge (D) (D) edge (E) (E) edge (A);
\end{scope}
\end{tikzpicture}%
\ \ $\\\hline
\end{tabular}
\]
Of course, it is more common to draw the path graph stretched out
horizontally:
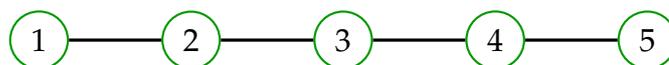
\[%
%
\]

Note that the cycle graph $C_{3}$ is identical with the complete graph $K_{3}%
$. \medskip

\textbf{Question:} What are the graph isomorphisms from $P_{n}$ to itself?

\textbf{Answer:} One such isomorphism is the identity map $\operatorname*{id}%
:\left\{  1,2,\ldots,n\right\}  \rightarrow\left\{  1,2,\ldots,n\right\}  $.
Another is the \textquotedblleft reversal\textquotedblright\ map
\begin{align*}
\left\{  1,2,\ldots,n\right\}   &  \rightarrow\left\{  1,2,\ldots,n\right\}
,\\
i  &  \mapsto n+1-i.
\end{align*}
There are no others. \medskip

\textbf{Question:} What are the graph isomorphisms from $C_{n}$ to itself?

\textbf{Answer:} For any $k\in\mathbb{Z}$, we can define a \textquotedblleft
rotation by $k$ vertices\textquotedblright, which is the map%
\begin{align*}
\left\{  1,2,\ldots,n\right\}   &  \rightarrow\left\{  1,2,\ldots,n\right\}
,\\
i  &  \mapsto\left(  i+k\text{ reduced modulo }n\text{ to an element of
}\left\{  1,2,\ldots,n\right\}  \right)  .
\end{align*}
Thus we get $n$ rotations (one for each $k\in\left\{  1,2,\ldots,n\right\}
$); all of them are graph isomorphisms.

There are also the reflections, which are the maps%
\begin{align*}
\left\{  1,2,\ldots,n\right\}   &  \rightarrow\left\{  1,2,\ldots,n\right\}
,\\
i  &  \mapsto\left(  k-i\text{ reduced modulo }n\text{ to an element of
}\left\{  1,2,\ldots,n\right\}  \right)
\end{align*}
for $k\in\mathbb{Z}$. There are $n$ of them, too, and they are isomorphisms as well.

Altogether we obtain $2n$ isomorphisms (for $n>2$), and there are no others.
(The group they form is the $n$-th dihedral group.)

\subsubsection{\label{subsec.sg.complete.kneser}Kneser graphs}

Here is a more exotic family of graphs:

\begin{example}
If $S$ is a finite set, and if $k\in\mathbb{N}$, then we define the
$k$\textbf{-th Kneser graph of }$S$ to be the simple graph%
\[
K_{S,k}:=\left(  \mathcal{P}_{k}\left(  S\right)  ,\ \ \left\{  IJ\ \mid
\ I,J\in\mathcal{P}_{k}\left(  S\right)  \text{ and }I\cap J=\varnothing
\right\}  \right)  .
\]
The vertices of $K_{S,k}$ are the $k$-element subsets of $S$, and two such
subsets are adjacent if they are disjoint.
\end{example}

The graph $K_{\left\{  1,2,\ldots,5\right\}  ,2}$ is called the
\textbf{Petersen graph}; here is how it looks like:%
\[%
\begin{tikzpicture}
\begin{scope}[every node/.style={circle,thick,draw=green!60!black}]
\node(A) at (0:2) {$\set{1,2}$};
\node(B) at (360/5:2) {$\set{2,3}$};
\node(C) at (2*360/5:2) {$\set{3,4}$};
\node(D) at (3*360/5:2) {$\set{4,5}$};
\node(E) at (4*360/5:2) {$\set{1,5}$};
\node(A2) at (0:5) {$\set{3,5}$};
\node(B2) at (360/5:5) {$\set{1,4}$};
\node(C2) at (2*360/5:5) {$\set{2,5}$};
\node(D2) at (3*360/5:5) {$\set{1,3}$};
\node(E2) at (4*360/5:5) {$\set{2,4}$};
\end{scope}
\begin{scope}[every edge/.style={draw=black,very thick}]
\path
[-] (A) edge (C) (C) edge (E) (E) edge (B) (B) edge (D) (D) edge (A) (A) edge (A2) (B) edge (B2) (C) edge (C2) (D) edge (D2) (E) edge (E2) (A2) edge (B2) (B2) edge (C2) (C2) edge (D2) (D2) edge (E2) (E2) edge (A2);
\end{scope}
\end{tikzpicture}%
\]

\subsection{\label{sec.sg.subgraphs}Subgraphs}

\begin{definition}
\label{def.sg.subgraph}Let $G=\left(  V,E\right)  $ be a simple graph.

\begin{enumerate}
\item[\textbf{(a)}] A \textbf{subgraph} of $G$ means a simple graph of the
form $H=\left(  W,F\right)  $, where $W\subseteq V$ and $F\subseteq E$. In
other words, a subgraph of $G$ means a simple graph whose vertices are
vertices of $G$ and whose edges are edges of $G$.

\item[\textbf{(b)}] Let $S$ be a subset of $V$. The \textbf{induced subgraph
of }$G$\textbf{ on the set }$S$ denotes the subgraph%
\[
\left(  S,\ \ E\cap\mathcal{P}_{2}\left(  S\right)  \right)
\]
of $G$. In other words, it denotes the subgraph of $G$ whose vertices are the
elements of $S$, and whose edges are precisely those edges of $G$ whose both
endpoints belong to $S$. We denote this induced subgraph by $G\left[
S\right]  $.

\item[\textbf{(c)}] An \textbf{induced subgraph} of $G$ means a subgraph of
$G$ that is the induced subgraph of $G$ on $S$ for some $S\subseteq V$.
\end{enumerate}
\end{definition}

Thus, a subgraph of a graph $G$ is obtained by throwing away some vertices and
some edges of $G$ (in such a way, of course, that no edges remain
\textquotedblleft dangling\textquotedblright\ -- i.e., if you throw away a
vertex, then you must throw away all edges that contain this vertex). Such a
subgraph is an induced subgraph if no edges are removed without need -- i.e.,
if you removed only those edges that lost some of their endpoints. Thus,
induced subgraphs can be characterized as follows:

\begin{proposition}
\label{prop.sg.subgraph.induced.crit}Let $H$ be a subgraph of a simple graph
$G$. Then, $H$ is an induced subgraph of $G$ if and only if each edge $uv$ of
$G$ whose endpoints $u$ and $v$ belong to $\operatorname*{V}\left(  H\right)
$ is an edge of $H$.
\end{proposition}

\begin{proof}
This is a matter of understanding the definition.
\end{proof}

\begin{example}
\label{exa.sg.K4-24}The following table shows a graph $G$ and two of its
subgraphs: one that is induced and one that is not.%
\[%
\begin{tabular}
[c]{|c|c|c|}\hline
graph $G$ & induced subgraph & subgraph that is\\
& on the set $\left\{  1,2,3\right\}  $ & not induced\\\hline
$\ \
\begin{tikzpicture}[scale=2]
\begin{scope}[every node/.style={circle,thick,draw=green!60!black}]
\node(1) at (0,0) {$1$};
\node(2) at (0,1) {$2$};
\node(3) at (1,1) {$3$};
\node(4) at (1,0) {$4$};
\end{scope}
\draw[very thick] (1) -- (2) -- (3) -- (1) -- (4) -- (3);
\end{tikzpicture}%
\ \ $ & $\ \
\begin{tikzpicture}[scale=2]
\begin{scope}[every node/.style={circle,thick,draw=green!60!black}]
\node(1) at (0,0) {$1$};
\node(2) at (0,1) {$2$};
\node(3) at (1,1) {$3$};
\end{scope}
\draw[very thick] (1) -- (2) -- (3) -- (1);
\end{tikzpicture}%
\ \ $ & $\ \
\begin{tikzpicture}[scale=2]
\begin{scope}[every node/.style={circle,thick,draw=green!60!black}]
\node(1) at (0,0) {$1$};
\node(2) at (0,1) {$2$};
\node(3) at (1,1) {$3$};
\end{scope}
\draw[very thick] (2) -- (1) -- (3);
\end{tikzpicture}%
\ \ $\\\hline
\end{tabular}
\
\]

\end{example}

\begin{example}
Let $n>1$ be an integer.

\begin{enumerate}
\item[\textbf{(a)}] The path graph $P_{n}$ is a subgraph of the cycle graph
$C_{n}$. It is not an induced subgraph (for $n>2$), because it contains the
two vertices $n$ and $1$ of $C_{n}$ but does not contain the edge $n1$.

\item[\textbf{(b)}] The path graph $P_{n-1}$ is an induced subgraph of $P_{n}%
$. (Namely, it is the induced subgraph of $P_{n}$ on the set $\left\{
1,2,\ldots,n-1\right\}  $.)

\item[\textbf{(c)}] Assume that $n>3$. Is $C_{n-1}$ a subgraph of $C_{n}$ ?
No, because the edge $\left(  n-1\right)  1$ belongs to $C_{n-1}$ but not to
$C_{n}$.
\end{enumerate}
\end{example}

The following is easy:

\begin{proposition}
Let $G$ be a simple graph, and let $H$ be a subgraph of $G$. Assume that $H$
is a complete graph. Then, $H$ is automatically an induced subgraph of $G$.
\end{proposition}

\begin{proof}
This follows from Proposition \ref{prop.sg.subgraph.induced.crit}, since the
completeness of $H$ means that each $2$-element subset $\left\{  u,v\right\}
$ of the vertex set of $H$ is an edge of $H$.
\end{proof}

We note that triangles in a graph can be characterized in terms of complete
subgraphs. Namely, a triangle \textquotedblleft is\textquotedblright\ the same
as a complete subgraph (or, equivalently, induced complete subgraph) with
three vertices:

\begin{remark}
Let $G$ be a simple graph. Let $u,v,w$ be three distinct vertices of $G$. The
following are equivalent:

\begin{enumerate}
\item The set $\left\{  u,v,w\right\}  $ is a triangle of $G$.

\item The induced subgraph of $G$ on $\left\{  u,v,w\right\}  $ is isomorphic
to $K_{3}$.

\item The induced subgraph of $G$ on $\left\{  u,v,w\right\}  $ is isomorphic
to $C_{3}$.
\end{enumerate}
\end{remark}

Thus, instead of saying \textquotedblleft triangle of $G$\textquotedblright,
one often says \textquotedblleft a $K_{3}$ in $G$\textquotedblright\ or
\textquotedblleft a $C_{3}$ in $G$\textquotedblright. Generally,
\textquotedblleft an $H$ in $G$\textquotedblright\ (where $H$ and $G$ are two
graphs) means a subgraph of $G$ that is isomorphic to $H$. (In the case when
$H=K_{3}=C_{3}$, it does not matter whether we require it to be a subgraph or
an induced subgraph, since a complete subgraph has to be induced automatically.)

\begin{example}
Let $G$ be the following simple graph:%
\[%
\begin{tikzpicture}[scale=2.5]
\begin{scope}[every node/.style={circle,thick,draw=green!60!black}]
\node(1) at (0,0) {$1$};
\node(2) at (1,0) {$2$};
\node(3) at (2,0) {$3$};
\node(4) at (3,0) {$4$};
\node(5) at (0,1) {$5$};
\node(6) at (1,1) {$6$};
\node(7) at (2,1) {$7$};
\node(8) at (3,1) {$8$};
\end{scope}
\draw[very thick] (1) -- (2) -- (6) -- (5) -- (1);
\draw[very thick] (2) -- (3) -- (6) -- (7) -- (2);
\draw[very thick] (3) -- (7) -- (8) -- (3) -- (4) -- (8);
\end{tikzpicture}%
\ \ .
\]
It has several subgraphs isomorphic to $C_{4}$, including the three
\textquotedblleft little squares\textquotedblright\ (note that the middle
square produces three different subgraphs isomorphic to $C_{4}$: you can
either keep the edges $23,37,76,62$ or keep the edges $23,36,67,72$ or keep
the edges $26,63,37,72$) as well as some less obvious subgraphs (e.g., the one
with vertices $3,6,7,8$ and edges $36,67,78,83$). But only one of these
subgraphs is an induced subgraph (namely, the leftmost \textquotedblleft
little square\textquotedblright, which is $G\left[  \left\{  1,2,5,6\right\}
\right]  $), since all the others have extraneous edges. The middle square
yields an induced subgraph isomorphic to $K_{4}$, namely $G\left[  \left\{
2,3,6,7\right\}  \right]  $.
\end{example}

\begin{exercise}
\label{exe.1.3}Let $n$ be a positive integer. Let $S$ be a simple graph with
$2n$ vertices. Prove that $S$ has two distinct vertices that have an even
number of common neighbors.
\end{exercise}

\begin{exercise}
\label{exe.1.8}Let $n\geq2$ be an integer. Let $G$ be a simple graph with $n$ vertices.

\begin{enumerate}
\item[\textbf{(a)}] Describe $G$ if the degrees of the vertices of $G$ are
$1,1,\ldots,1,n-1$.

\item[\textbf{(b)}] Let $a$ and $b$ be two positive integers such that $a + b
= n$. Describe $G$ if the degrees of the vertices of $G$ are $1, 1, \ldots, 1,
a, b$.
\end{enumerate}

Here, to \textquotedblleft describe\textquotedblright\ $G$ means to explicitly
determine (with proof) a graph that is isomorphic to $G$.
\end{exercise}

\begin{remark}
The situations in Exercise \ref{exe.1.8} are, in a sense, exceptional.
Typically, the degrees of the vertices of a graph do not uniquely determine
the graph up to isomorphism. For example, the two graphs
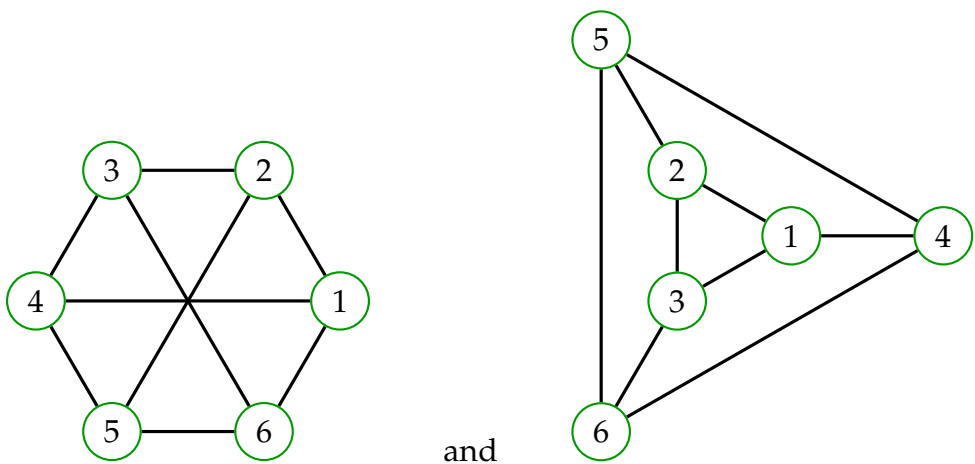
\[
\begin{tikzpicture}
\begin{scope}[every node/.style={circle,thick,draw=green!60!black}]
\node(A) at (0:2) {$1$};
\node(B) at (60:2) {$2$};
\node(C) at (120:2) {$3$};
\node(D) at (180:2) {$4$};
\node(E) at (240:2) {$5$};
\node(F) at (300:2) {$6$} ;
\end{scope}
\begin{scope}[every edge/.style={draw=black,very thick}]
\path[-] (A) edge (B) (B) edge (C) (C) edge (D) (D) edge (E) (E) edge (F) (F) edge (A);
\path[-] (A) edge (D);
\path[-] (B) edge (E);
\path[-] (C) edge (F);
\end{scope}
\end{tikzpicture}
\qquad\text{ and } \qquad\begin{tikzpicture}
\begin{scope}[every node/.style={circle,thick,draw=green!60!black}]
\node(A) at (0:1) {$1$};
\node(B) at (120:1) {$2$};
\node(C) at (240:1) {$3$};
\node(D) at (0:3) {$4$};
\node(E) at (120:3) {$5$};
\node(F) at (240:3) {$6$} ;
\end{scope}
\begin{scope}[every edge/.style={draw=black,very thick}]
\path[-] (A) edge (B) (B) edge (C) (C) edge (A);
\path[-] (D) edge (E) (E) edge (F) (F) edge (D);
\path[-] (A) edge (D);
\path[-] (B) edge (E);
\path[-] (C) edge (F);
\end{scope}
\end{tikzpicture}
\]
are not isomorphic\footnotemark, but have the same degrees (namely, each
vertex of either graph has degree $3$).
\end{remark}

\footnotetext{The easiest way to see this is to observe that the second graph
has a triangle (i.e., three distinct vertices that are mutually adjacent),
while the first graph does not.}

\subsection{\label{sec.sg.djun}Disjoint unions}

Another way of constructing new graphs from old is the disjoint union. The
idea is simple: Taking the disjoint union $G_{1}\sqcup G_{2}\sqcup\cdots\sqcup
G_{k}$ of several simple graphs $G_{1},G_{2},\ldots,G_{k}$ means putting the
graphs alongside each other and treating the result as one big graph. To make
this formally watertight, we have to relabel each vertex $v$ of each graph
$G_{i}$ as the pair $\left(  i,v\right)  $, so that vertices coming from
different graphs appear as different even if they were equal. For example, the
disjoint union $C_{3}\sqcup C_{4}$ of the two cycle graphs $C_{3}$ and $C_{4}$
should not be%
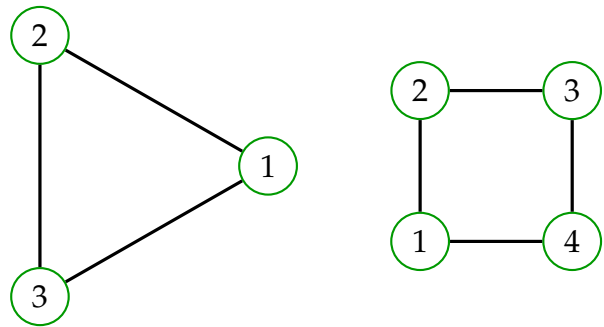
\[%
%
\]
(which makes no sense, because there are two points labelled $1$ in this
picture, but a graph can have only one vertex $1$), but rather should be%
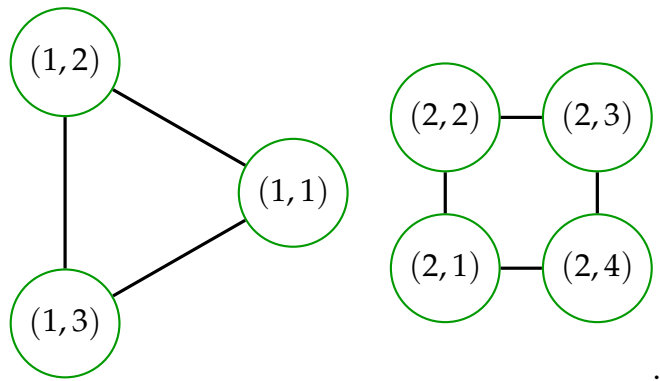
\[%
%
.
\]

So here is the formal definition:

\begin{definition}
Let $G_{1},G_{2},\ldots,G_{k}$ be simple graphs, where $G_{i}=\left(
V_{i},E_{i}\right)  $ for each $i\in\left\{  1,2,\ldots,k\right\}  $. The
\textbf{disjoint union} of these $k$ graphs $G_{1},G_{2},\ldots,G_{k}$ is
defined to be the simple graph $\left(  V,E\right)  $, where%
\begin{align*}
V  &  =\left\{  \left(  i,v\right)  \ \mid\ i\in\left\{  1,2,\ldots,k\right\}
\text{ and }v\in V_{i}\right\}  \ \ \ \ \ \ \ \ \ \ \text{and}\\
E  &  =\left\{  \left\{  \left(  i,v_{1}\right)  ,\left(  i,v_{2}\right)
\right\}  \ \mid\ i\in\left\{  1,2,\ldots,k\right\}  \text{ and }\left\{
v_{1},v_{2}\right\}  \in E_{i}\right\}  .
\end{align*}
This disjoint union is denoted by $G_{1}\sqcup G_{2}\sqcup\cdots\sqcup G_{k}$.
\end{definition}

Note: If $G$ and $H$ are two graphs, then the two graphs $G\sqcup H$ and
$H\sqcup G$ are isomorphic, but not the same graph (unless $G=H$). For
example, $C_{3}\sqcup C_{4}$ has a vertex $\left(  2,4\right)  $, but
$C_{4}\sqcup C_{3}$ does not.

\subsection{\label{sec.sg.walks}Walks and paths}

We now come to the definitions of walks and paths -- two of the most
fundamental features that graphs can have. In particular, Euler's 1736 paper,
where graphs were first studied, is about certain kinds of walks.

\subsubsection{Definitions}

Imagine a graph as a road network, where each vertex is a town and each edge
is a (bidirectional) road. By successively walking along several edges, you
can often get from a town to another even if they are not adjacent. This is
made formal in the concept of a \textquotedblleft walk\textquotedblright:

\begin{definition}
Let $G$ be a simple graph. Then:

\begin{enumerate}
\item[\textbf{(a)}] A \textbf{walk} (in $G$) means a finite sequence $\left(
v_{0},v_{1},\ldots,v_{k}\right)  $ of vertices of $G$ (with $k\geq0$) such
that all of $v_{0}v_{1},\ v_{1}v_{2},\ v_{2}v_{3},\ \ldots,\ v_{k-1}v_{k}$ are
edges of $G$. (The latter condition is vacuously true if $k=0$.)

\item[\textbf{(b)}] If $\mathbf{w}=\left(  v_{0},v_{1},\ldots,v_{k}\right)  $
is a walk in $G$, then:

\begin{itemize}
\item The \textbf{vertices} of $\mathbf{w}$ are defined to be $v_{0}%
,v_{1},\ldots,v_{k}$.

\item The \textbf{edges} of $\mathbf{w}$ are defined to be $v_{0}v_{1}%
,\ v_{1}v_{2},\ v_{2}v_{3},\ \ldots,\ v_{k-1}v_{k}$.

\item The nonnegative integer $k$ is called the \textbf{length} of
$\mathbf{w}$. (This is the number of all edges of $\mathbf{w}$, counted with
multiplicity. It is $1$ smaller than the number of all vertices of
$\mathbf{w}$, counted with multiplicity.)

\item The vertex $v_{0}$ is called the \textbf{starting point} of $\mathbf{w}%
$. We say that $\mathbf{w}$ \textbf{starts} (or \textbf{begins}) at $v_{0}$.

\item The vertex $v_{k}$ is called the \textbf{ending point} of $\mathbf{w}$.
We say that $\mathbf{w}$ \textbf{ends} at $v_{k}$.
\end{itemize}

\item[\textbf{(c)}] A \textbf{path} (in $G$) means a walk (in $G$) whose
vertices are distinct. In other words, a path means a walk $\left(
v_{0},v_{1},\ldots,v_{k}\right)  $ such that $v_{0},v_{1},\ldots,v_{k}$ are distinct.

\item[\textbf{(d)}] Let $p$ and $q$ be two vertices of $G$. A \textbf{walk
from }$p$ \textbf{to }$q$ means a walk that starts at $p$ and ends at $q$. A
\textbf{path from }$p$ \textbf{to }$q$ means a path that starts at $p$ and
ends at $q$.

\item[\textbf{(e)}] We often say \textquotedblleft walk of $G$%
\textquotedblright\ and \textquotedblleft path of $G$\textquotedblright%
\ instead of \textquotedblleft walk in $G$\textquotedblright\ and
\textquotedblleft path in $G$\textquotedblright, respectively.
\end{enumerate}
\end{definition}

\begin{example}
\label{exa.sg.walk.example1}Let $G$ be the graph%
\[
\left(  \left\{  1,2,3,4,5,6\right\}  ,\ \ \left\{
12,\ 23,\ 34,\ 45,\ 56,\ 61,\ 13\right\}  \right)  .
\]
This graph looks as follows:%
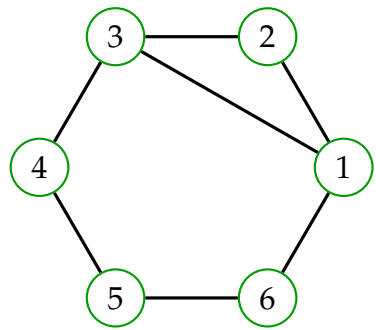
\[%
\begin{tikzpicture}
\begin{scope}[every node/.style={circle,thick,draw=green!60!black}]
\node(A) at (0:2) {$1$};
\node(B) at (60:2) {$2$};
\node(C) at (120:2) {$3$};
\node(D) at (180:2) {$4$};
\node(E) at (240:2) {$5$};
\node(F) at (300:2) {$6$};
\end{scope}
\begin{scope}[every edge/.style={draw=black,very thick}]
\path
[-] (A) edge (B) (B) edge (C) (C) edge (A) (C) edge (D) (D) edge (E) (E) edge (F) (F) edge (A);
\end{scope}
\end{tikzpicture}%
\]
Then:

\begin{itemize}
\item The sequence $\left(  1,3,4,5,6,1,3,2\right)  $ of vertices of $G$ is a
walk in $G$. This walk is a walk from $1$ to $2$. It is not a path. The length
of this walk is $7$.

\item The sequence $\left(  1,2,4,3\right)  $ of vertices of $G$ is not a
walk, since $24$ is not an edge of $G$. Hence, it is not a path either.

\item The sequence $\left(  1,3,2,1\right)  $ is a walk from $1$ to $1$. It
has length $3$. It is not a path.

\item The sequence $\left(  1,2,1\right)  $ is a walk from $1$ to $1$. It has
length $2$. It is not a path.

\item The sequence $\left(  5\right)  $ is a walk from $5$ to $5$. It has
length $0$. It is a path. More generally, each vertex $v$ of $G$ produces a
length-$0$ path $\left(  v\right)  $.

\item The sequence $\left(  5,4\right)  $ is a walk from $5$ to $4$. It has
length $1$. It is a path. More generally, each edge $uv$ of $G$ produces a
length-$1$ path $\left(  u,v\right)  $.
\end{itemize}
\end{example}

Intuitively, we can think of walks and paths as follows:

\begin{itemize}
\item A \textbf{walk} of a graph is a way of walking from one vertex to
another (or to the same vertex) by following a sequence of edges.

\item A \textbf{path} is a walk whose vertices are distinct (i.e., each vertex
appears at most once in the walk).
\end{itemize}

\begin{exercise}
\label{exe.intro.path.edges-dist} Let $G$ be a simple graph. Let $\mathbf{w}$
be a path in $G$. Prove that the edges of $\mathbf{w}$ are distinct. (This may
look obvious when you can point to a picture; but we ask you to give a
rigorous proof!) \medskip

[\textbf{Solution:} This is Exercise 3 on homework set \#1 from my Spring 2017
course; see \href{https://www.cip.ifi.lmu.de/~grinberg/t/17s/}{the course
page} for solutions.]
\end{exercise}

\subsubsection{Composing/concatenating and reversing walks}

Here are some simple things we can do with walks and paths.

First, we can \textquotedblleft splice\textquotedblright\ two walks together
if the ending point of the first is the starting point of the second:

\begin{proposition}
\label{prop.sg.walk-concat}Let $G$ be a simple graph. Let $u$, $v$ and $w$ be
three vertices of $G$. Let $\mathbf{a}=\left(  a_{0},a_{1},\ldots
,a_{k}\right)  $ be a walk from $u$ to $v$. Let $\mathbf{b}=\left(
b_{0},b_{1},\ldots,b_{\ell}\right)  $ be a walk from $v$ to $w$. Then,%
\begin{align*}
\left(  a_{0},a_{1},\ldots,a_{k},b_{1},b_{2},\ldots,b_{\ell}\right)   &
=\left(  a_{0},a_{1},\ldots,a_{k-1},b_{0},b_{1},\ldots,b_{\ell}\right) \\
&  =\left(  a_{0},a_{1},\ldots,a_{k-1},v,b_{1},b_{2},\ldots,b_{\ell}\right)
\end{align*}
is a walk from $u$ to $w$. This walk shall be denoted $\mathbf{a}%
\ast\mathbf{b}$.
\end{proposition}

\begin{proof}
Intuitively clear and straightforward to verify.
\end{proof}

Note that \textquotedblleft splicing\textquotedblright\ does not preserve
path-ness: If $\mathbf{a}$ and $\mathbf{b}$ are two paths, then the walk
$\mathbf{a}\ast\mathbf{b}$ is not necessarily a path. \medskip

Another thing that can be done with walks (and this one does preserve
path-ness) is walking them backwards:

\begin{proposition}
\label{prop.sg.walk-rev}Let $G$ be a simple graph. Let $u$ and $v$ be two
vertices of $G$. Let $\mathbf{a}=\left(  a_{0},a_{1},\ldots,a_{k}\right)  $ be
a walk from $u$ to $v$. Then:

\begin{enumerate}
\item[\textbf{(a)}] The list $\left(  a_{k},a_{k-1},\ldots,a_{0}\right)  $ is
a walk from $v$ to $u$. We denote this walk by $\operatorname*{rev}\mathbf{a}$
and call it the \textbf{reversal} of $\mathbf{a}$.

\item[\textbf{(b)}] If $\mathbf{a}$ is a path, then $\operatorname*{rev}%
\mathbf{a}$ is a path again.
\end{enumerate}
\end{proposition}

\begin{proof}
Intuitively clear and straightforward to verify.
\end{proof}

\subsubsection{Reducing walks to paths}

A path is just a walk without repeated vertices. If you have a walk, you can
turn it into a path by removing \textquotedblleft loops\textquotedblright\ (or
\textquotedblleft digressions\textquotedblright):

\begin{proposition}
\label{prop.sg.walk-to-path-1}Let $G$ be a simple graph. Let $u$ and $v$ be
two vertices of $G$. Let $\mathbf{a}=\left(  a_{0},a_{1},\ldots,a_{k}\right)
$ be a walk from $u$ to $v$. Assume that $\mathbf{a}$ is not a path. Then,
there exists a walk from $u$ to $v$ whose length is smaller than $k$.
\end{proposition}

\begin{proof}
Since $\mathbf{a}$ is not a path, two of its vertices are equal. In other
words, there exist $i<j$ such that $a_{i}=a_{j}$. Consider these $i$ and $j$.
Now, consider the tuple%
\[
\left(  \underbrace{a_{0},a_{1},\ldots,a_{i}}_{\text{the first }i+1\text{
vertices of }\mathbf{a}},\underbrace{a_{j+1},a_{j+2},\ldots,a_{k}}_{\text{the
last }k-j\text{ vertices of }\mathbf{a}}\right)
\]
(this is just $\mathbf{a}$ with the part between $a_{i}$ and $a_{j}$ cut out).
This tuple is a walk from $u$ to $v$, and its length is $\underbrace{i}%
_{<j}+\left(  k-j\right)  <j+\left(  k-j\right)  =k$. So we have found a walk
from $u$ to $v$ whose length is smaller than $k$. This proves the proposition.
\end{proof}

\begin{example}
Consider the walk $\left(  1,3,4,5,6,1,3,2\right)  $ from Example
\ref{exa.sg.walk.example1}. Then, Proposition \ref{prop.sg.walk-to-path-1}
tells us that there is a walk from $1$ to $2$ that has smaller length. You can
find this walk by removing the part between the two $3$'s. You get the walk
$\left(  1,3,2\right)  $. This is actually a path.
\end{example}

\begin{corollary}
[When there is a walk, there is a path]\label{cor.sg.walk-thus-path}Let $G$ be
a simple graph. Let $u$ and $v$ be two vertices of $G$. Assume that there is a
walk from $u$ to $v$ of length $k$ for some $k\in\mathbb{N}$. Then, there is a
path from $u$ to $v$ of length $\leq k$.
\end{corollary}

\begin{proof}
Proposition \ref{prop.sg.walk-to-path-1} says that if there is a walk from $u$
to $v$ that is not a path, then there is a walk from $u$ to $v$ having shorter
length. Apply this repeatedly, until you get a path. (You will eventually get
a path, because the length cannot keep decreasing forever.)
\end{proof}

\subsubsection{\label{subsec.sg.walks.algo}Remark on algorithms}

We take a little break from proving structural theorems in order to address
some important computational questions. As always in these notes, we will only
scratch the surface and content ourselves with simple but not quite optimal algorithms.

Given a simple graph $G$ and two vertices $u$ and $v$ of $G$, we can ask
ourselves the following questions:

\begin{statement}
\textbf{Question 1:} Does $G$ have a walk from $u$ to $v$ ?

\textbf{Question 2:} Does $G$ have a path from $u$ to $v$ ?

\textbf{Question 3:} Find a shortest path from $u$ to $v$ (that is, a path
from $u$ to $v$ having the smallest possible length), or determine that no
such path exists.

\textbf{Question 4:} Given a number $k\in\mathbb{N}$, find a walk from $u$ to
$v$ having length $k$, or determine that no such walk exists.

\textbf{Question 5:} Given a number $k\in\mathbb{N}$, find a path from $u$ to
$v$ having length $k$, or determine that no such path exists.
\end{statement}

Corollary \ref{cor.sg.walk-thus-path} reveals that Questions 1 and 2 are
equivalent (indeed, the existence of a walk from $u$ to $v$ entails the
existence of a path from $u$ to $v$ by Corollary \ref{cor.sg.walk-thus-path},
whereas the converse is obvious). Question 3 is clearly a stronger version of
Question 2 (in the sense that any answer to Question 3 will automatically
answer Question 2 as well). \medskip

With a bit more thought, it is easily seen that Question 4 is a stronger
version of Question 3. Indeed, Corollary \ref{cor.sg.walk-thus-path} shows
that a shortest walk from $u$ to $v$ (if it exists) must also be a shortest
path from $u$ to $v$. However, any path from $u$ to $v$ must have length $\leq
n-1$, where $n$ is the number of vertices of $G$ (since a path of length $k$
has $k+1$ distinct vertices, but $G$ has only $n$ vertices to spare). Hence,
if there is no walk of length $\leq n-1$ from $u$ to $v$, then there is no
path from $u$ to $v$ whatsoever. Thus, if we answer Question 4 for all values
$k\in\left\{  0,1,\ldots,n-1\right\}  $, then we obtain either a shortest path
from $u$ to $v$ (by taking the smallest $k$ for which the answer is positive,
and then picking the resulting walk, which must be a shortest path by what we
previously said), or proof positive that no path from $u$ to $v$ exists (if
the answer for each $k\in\left\{  0,1,\ldots,n-1\right\}  $ is negative).

Thus, answering Question 4 will yield answers to Questions 1, 2 and 3.
\medskip

Let us now outline a way how Question 4 can be answered using a recursive
algorithm. Specifically, we recurse on $k$. The base case ($k=0$) is easy: A
walk from $u$ to $v$ having length $0$ exists if $u=v$ and does not exist
otherwise. The interesting part is the recursion step: Assume that the integer
$k$ is positive, and that we already know how to answer Question 4 for $k-1$
instead of $k$. Now, let us answer it for $k$. To do so, we observe that any
walk from $u$ to $v$ having length $k$ must have the form $\left(
u,\ldots,w,v\right)  $, where the penultimate vertex $w$ is some neighbor of
$v$. Moreover, if we remove the last vertex $v$ from our walk $\left(
u,\ldots,w,v\right)  $, then we obtain a walk $\left(  u,\ldots,w\right)  $ of
length $k-1$. Hence, we can find a walk from $u$ to $v$ having length $k$ as follows:

\begin{itemize}
\item We make a list of all neighbors of $v$. We go through this list in some
arbitrary order.

\item For each neighbor $w$ in this list, we try to find a walk from $u$ to
$w$ having length $k-1$ (this is a matter of answering Question 4 for $k-1$
instead of $k$, so we supposedly already know how to do this). If such a walk
exists, then we simply insert $v$ at its end, and thus obtain a walk from $u$
to $v$ having length $k$. Thus we obtain a positive answer to our question.

\item If we have gone through our whole list of neighbors of $v$ without
finding a walk from $u$ to $v$ having length $k$, then no such walk exists,
and thus we have found a negative answer.
\end{itemize}

This recursive algorithm answers Question 4, and is fast enough to be
practically viable if implemented well. (In the language of complexity theory,
it is a \href{https://en.wikipedia.org/wiki/Time_complexity}{polynomial time
algorithm}\footnote{To be specific: Its running time can be bounded in a
polynomial of $n$ and $k$, where $n$ is the number of vertices of $G$.}.) Much
more efficient algorithms exist, however. In applications, a generalized
version of Question 3 often appears, asking for a path that is shortest not in
the sense of smallest length, but in the sense of smallest \textquotedblleft
weighted length\textquotedblright\ (i.e., different edges contribute
differently to this \textquotedblleft length\textquotedblright). This
generalized question is one of the most fundamental algorithmic problems in
computer science, known as the \textbf{shortest path problem}, and various
algorithms can be found on
\href{https://en.wikipedia.org/wiki/Shortest_path_problem}{its Wikipedia page}
and in algorithm-focussed texts such as \cite[\S 3.5]{Griffi21},
\cite[\S 12.3]{KelTro17}, \cite[\S 1.5]{Even12}, \cite[Chapter 1]{Schrij-ACO}
or (for a royal treatment) \cite[Chapters 6--8]{Schrij-CO1}.\medskip

Question 5 looks superficially similar to Question 4, yet it differs in the
most important way: There is no efficient algorithm known for answering it! In
the language of complexity theory, it is an
\href{https://en.wikipedia.org/wiki/NP-hardness}{NP-hard problem}, which means
that a polynomial-time algorithm for it is not expected to exist (although
this is the kind of negative that appears near-impossible to prove at the
current stage of the discipline). It is still technically a finite problem
(there are only finitely many possible paths in $G$, and thus one can
theoretically try them all), and there is even a polynomial-time algorithm for
any fixed value of $k$ (again, a trivial one: check all the $n^{k+1}$ possible
$\left(  k+1\right)  $-tuples of vertices of $G$ for whether they are paths
from $u$ to $v$), but the complexity of this algorithm grows exponentially in
$k$, which makes it useless in practice.

\subsubsection{The equivalence relation \textquotedblleft
path-connected\textquotedblright}

We can use the concepts of walks and paths to define a certain equivalence
relation on the vertex set $\operatorname*{V}\left(  G\right)  $ of any graph
$G$:

\begin{definition}
\label{def.sg.pc}Let $G$ be a simple graph. We define a binary relation
$\simeq_{G}$ on the set $\operatorname*{V}\left(  G\right)  $ as follows: For
two vertices $u$ and $v$ of $G$, we shall have $u\simeq_{G}v$ if and only if
there exists a walk from $u$ to $v$ in $G$.

This binary relation $\simeq_{G}$ is called \textquotedblleft%
\textbf{path-connectedness}\textquotedblright\ or just \textquotedblleft%
\textbf{connectedness}\textquotedblright. When two vertices $u$ and $v$
satisfy $u\simeq_{G}v$, we say that \textquotedblleft$u$ and $v$ are
\textbf{path-connected}\textquotedblright.
\end{definition}

\begin{proposition}
\label{prop.sg.pc.eqrel}Let $G$ be a simple graph. Then, the relation
$\simeq_{G}$ is an equivalence relation.
\end{proposition}

\begin{proof}
We need to show that $\simeq_{G}$ is symmetric, reflexive and transitive.

\begin{itemize}
\item \textbf{Symmetry:} If $u\simeq_{G}v$, then $v\simeq_{G}u$, because we
can take a walk from $u$ to $v$ and reverse it.

\item \textbf{Reflexivity:} We always have $u\simeq_{G}u$, since the trivial
walk $\left(  u\right)  $ is a walk from $u$ to $u$.

\item \textbf{Transitivity:} If $u\simeq_{G}v$ and $v\simeq_{G}w$, then
$u\simeq_{G}w$, because (as we know from Proposition \ref{prop.sg.walk-concat}%
) we can take a walk $\mathbf{a}$ from $u$ to $v$ and a walk $\mathbf{b}$ from
$v$ to $w$ and combine them to form the walk $\mathbf{a}\ast\mathbf{b}$
defined in Proposition \ref{prop.sg.walk-concat}.
\end{itemize}
\end{proof}

\begin{proposition}
\label{prop.sg.pc.there-is-path}Let $G$ be a simple graph. Let $u$ and $v$ be
two vertices of $G$. Then, $u\simeq_{G}v$ if and only if there exists a path
from $u$ to $v$.
\end{proposition}

\begin{proof}
$\Longleftarrow:$ Clear, since any path is a walk.

$\Longrightarrow:$ This is just saying that if there is a walk from $u$ to
$v$, then there is a path from $u$ to $v$. But this follows from Corollary
\ref{cor.sg.walk-thus-path}.
\end{proof}

\subsubsection{Connected components and connectedness}

The equivalence relation $\simeq_{G}$ introduced in Definition \ref{def.sg.pc}
allows us to define two important concepts:

\begin{definition}
\label{def.sg.component}Let $G$ be a simple graph. The equivalence classes of
the equivalence relation $\simeq_{G}$ are called the \textbf{connected
components} (or, for short, \textbf{components}) of $G$.
\end{definition}

\begin{definition}
\label{def.sg.connected}Let $G$ be a simple graph. We say that $G$ is
\textbf{connected} if $G$ has exactly one component.
\end{definition}

Thus, a simple graph $G$ is connected if and only if it has at least one
component (i.e., it has at least one vertex) and it has at most one component
(i.e., each two of its vertices are path-connected).

\begin{example}
Let $G$ be the graph with vertex set $\left\{  1,2,\ldots,9\right\}  $ and
such that two vertices $i$ and $j$ are adjacent if and only if $\left\vert
i-j\right\vert =3$. What are the components of $G$ ?

The graph $G$ looks like this:%
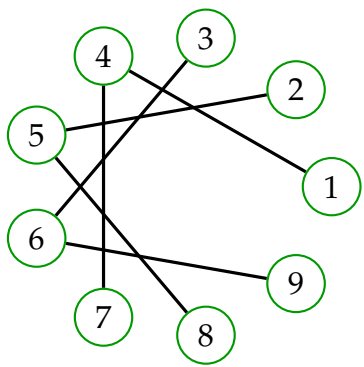
\[%
\begin{tikzpicture}
\begin{scope}[every node/.style={circle,thick,draw=green!60!black}]
\node(A) at (0:2) {$1$};
\node(B) at (360/9:2) {$2$};
\node(C) at (2*360/9:2) {$3$};
\node(D) at (3*360/9:2) {$4$};
\node(E) at (4*360/9:2) {$5$};
\node(F) at (5*360/9:2) {$6$};
\node(G) at (6*360/9:2) {$7$};
\node(H) at (7*360/9:2) {$8$};
\node(I) at (8*360/9:2) {$9$};
\end{scope}
\begin{scope}[every edge/.style={draw=black,very thick}]
\path
[-] (A) edge (D) (D) edge (G) (B) edge (E) (E) edge (H) (C) edge (F) (F) edge (I);
\end{scope}
\end{tikzpicture}%
\ \ .
\]
This looks like a jumbled mess, so you might think that all vertices are
mutually path-connected. But this is not the case, because edges that cross in
a drawing do not necessarily have endpoints in common. Walks can only move
from one edge to another at a common endpoint. Thus, there are much fewer
walks than the picture might suggest. We have $1\simeq_{G}4\simeq_{G}7$ and
$2\simeq_{G}5\simeq_{G}8$ and $3\simeq_{G}6\simeq_{G}9$, but there are no
further $\simeq_{G}$-relations. In fact, two vertices of $G$ are adjacent only
if they are congruent modulo $3$ (as numbers), and therefore you cannot move
from one modulo-$3$ congruence class to another by walking along edges of $G$.
So the components of $G$ are $\left\{  1,4,7\right\}  $ and $\left\{
2,5,8\right\}  $ and $\left\{  3,6,9\right\}  $. The graph $G$ is not connected.
\end{example}

\begin{example}
Let $G$ be the graph with vertex set $\left\{  1,2,\ldots,9\right\}  $ and
such that two vertices $i$ and $j$ are adjacent if and only if $\left\vert
i-j\right\vert =6$. This graph looks like this:%
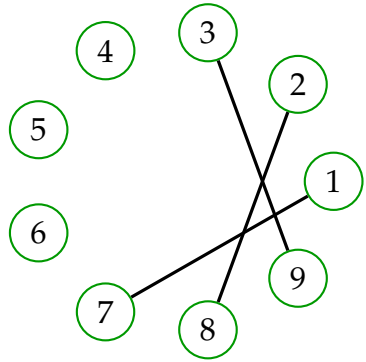
\[%
\begin{tikzpicture}
\begin{scope}[every node/.style={circle,thick,draw=green!60!black}]
\node(A) at (0:2) {$1$};
\node(B) at (360/9:2) {$2$};
\node(C) at (2*360/9:2) {$3$};
\node(D) at (3*360/9:2) {$4$};
\node(E) at (4*360/9:2) {$5$};
\node(F) at (5*360/9:2) {$6$};
\node(G) at (6*360/9:2) {$7$};
\node(H) at (7*360/9:2) {$8$};
\node(I) at (8*360/9:2) {$9$};
\end{scope}
\begin{scope}[every edge/.style={draw=black,very thick}]
\path[-] (A) edge (G) (B) edge (H) (C) edge (I);
\end{scope}
\end{tikzpicture}%
\ \ .
\]
What are the components of $G$ ? They are $\left\{  1,7\right\}  $ and
$\left\{  2,8\right\}  $ and $\left\{  3,9\right\}  $ and $\left\{  4\right\}
$ and $\left\{  5\right\}  $ and $\left\{  6\right\}  $. Note that three of
these six components are singleton sets. The graph $G$ is not connected.
\end{example}

\begin{example}
Let $G$ be the graph with vertex set $\left\{  1,2,\ldots,9\right\}  $ and
such that two vertices $i$ and $j$ are adjacent if and only if $\left\vert
i-j\right\vert =3$ or $\left\vert i-j\right\vert =4$. This graph looks like
this:%
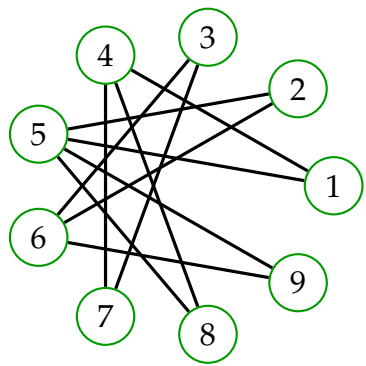
\[%
\begin{tikzpicture}
\begin{scope}[every node/.style={circle,thick,draw=green!60!black}]
\node(A) at (0:2) {$1$};
\node(B) at (360/9:2) {$2$};
\node(C) at (2*360/9:2) {$3$};
\node(D) at (3*360/9:2) {$4$};
\node(E) at (4*360/9:2) {$5$};
\node(F) at (5*360/9:2) {$6$};
\node(G) at (6*360/9:2) {$7$};
\node(H) at (7*360/9:2) {$8$};
\node(I) at (8*360/9:2) {$9$};
\end{scope}
\begin{scope}[every edge/.style={draw=black,very thick}]
\path
[-] (A) edge (D) (D) edge (G) (B) edge (E) (E) edge (H) (C) edge (F) (F) edge (I) (A) edge (E) (B) edge (F) (C) edge (G) (D) edge (H) (E) edge (I);
\end{scope}
\end{tikzpicture}%
\ \ .
\]
We can take a long walk through $G$:%
\[
\left(  1,4,7,3,6,9,5,2,5,8\right)  .
\]
This walk traverses every vertex of $G$; thus, any two vertices of $G$ are
path-connected. Hence, $G$ has only one component, namely $\left\{
1,2,\ldots,9\right\}  $. Thus, $G$ is connected.
\end{example}

\begin{example}
The complete graph on a nonempty set is connected. The complete graph on the
empty set is not connected, since it has $0$ (not $1$) components.
\end{example}

\begin{example}
The empty graph on a finite set $V$ has $\left\vert V\right\vert $ many
components (those are the singleton sets $\left\{  v\right\}  $ for $v\in V$).
Thus, it is connected if and only if $\left\vert V\right\vert =1$.
\end{example}

\begin{exercise}
\label{exe.3.7}Let $k\in\mathbb{N}$. Let $S$ be a finite set.

Recall that the \textbf{Kneser graph} $K_{S,k}$ is the simple graph whose
vertices are the $k$-element subsets of $S$, and whose edges are the unordered
pairs $\left\{  A,B\right\}  $ consisting of two such subsets $A$ and $B$ that
satisfy $A\cap B=\varnothing$.

Prove that this Kneser graph $K_{S,k}$ is connected if $\left\vert
S\right\vert \geq2k+1$. \medskip

[\textbf{Remark:} Can the \textquotedblleft if\textquotedblright\ here be
replaced by an \textquotedblleft if and only if\textquotedblright? Not quite,
because the graph $K_{S,k}$ is also connected if $\left\vert S\right\vert =2$
and $k=1$ (in which case it has two vertices and one edge), or if $\left\vert
S\right\vert =k$ (in which case it has only one vertex), or if $k=0$ (in which
case it has only one vertex). But these are the only \textquotedblleft
exceptions\textquotedblright.]
\end{exercise}

\subsubsection{Induced subgraphs on components}

Recall that if $S$ is a subset of the vertex set of a simple graph $G$, then
$G\left[  S\right]  $ denotes the induced subgraph of $G$ on $S$ (see
Definition \ref{def.sg.subgraph} \textbf{(b)}). The following is not hard to see:

\begin{proposition}
\label{prop.sg.component.ind-connected}Let $G$ be a simple graph. Let $C$ be a
component of $G$.

Then, the graph $G\left[  C\right]  $ (that is, the induced subgraph of $G$ on
the set $C$) is connected.
\end{proposition}

\begin{proof}
We need to show that $G\left[  C\right]  $ is connected. In other words, we
need to show that $G\left[  C\right]  $ has exactly $1$ component.

Clearly, $G\left[  C\right]  $ has at least one vertex (since $C$ is a
component, i.e., an equivalence class of $\simeq_{G}$, but equivalence classes
are always nonempty), thus has at least $1$ component. So we only need to show
that $G\left[  C\right]  $ has no more than $1$ component. In other words, we
need to show that any two vertices of $G\left[  C\right]  $ are path-connected
in $G\left[  C\right]  $.

So let $u$ and $v$ be two vertices of $G\left[  C\right]  $. Then, $u,v\in C$,
and therefore $u\simeq_{G}v$ (since $C$ is a component of $G$). In other
words, there exists a walk $\mathbf{w}=\left(  w_{0},w_{1},\ldots
,w_{k}\right)  $ from $u$ to $v$ in $G$. We shall now prove that this walk
$\mathbf{w}$ is actually a walk of $G\left[  C\right]  $. In other words, we
shall prove that all vertices of $\mathbf{w}$ belong to $C$.

But this is easy: If $w_{i}$ is a vertex of $\mathbf{w}$, then $\left(
w_{0},w_{1},\ldots,w_{i}\right)  $ is a walk from $u$ to $w_{i}$ in $G$, and
therefore we have $u\simeq_{G}w_{i}$, so that $w_{i}$ belongs to the same
component of $G$ as $u$; but that component is $C$. Thus, we have shown that
each vertex $w_{i}$ of $\mathbf{w}$ belongs to $C$. Therefore, $\mathbf{w}$ is
a walk of the graph $G\left[  C\right]  $. Consequently, it shows that
$u\simeq_{G\left[  C\right]  }v$.

We have now proved that $u\simeq_{G\left[  C\right]  }v$ for any two vertices
$u$ and $v$ of $G\left[  C\right]  $. Hence, the relation $\simeq_{G\left[
C\right]  }$ has no more than $1$ equivalence class. In other words, the graph
$G\left[  C\right]  $ has no more than $1$ component. This completes our proof.
\end{proof}

Note that the converse of Proposition \ref{prop.sg.component.ind-connected} is
not true: If $S$ is a subset of the vertex set of a simple graph $G$ such that
$G\left[  S\right]  $ is connected, then it does not follow that $S$ is a
component of $G$.

\begin{noncompile}
In the following proposition, we are using the notation $G\left[  S\right]  $
for the induced subgraph of a simple graph $G$ on a subset $S$ of its vertex set.
\end{noncompile}

\begin{proposition}
\label{prop.sg.component.djun}Let $G$ be a simple graph. Let $C_{1}%
,C_{2},\ldots,C_{k}$ be all components of $G$ (listed without repetition).

Thus, $G$ is isomorphic to the disjoint union $G\left[  C_{1}\right]  \sqcup
G\left[  C_{2}\right]  \sqcup\cdots\sqcup G\left[  C_{k}\right]  $.
\end{proposition}

\begin{proof}
Consider the bijection from $\operatorname*{V}\left(  G\left[  C_{1}\right]
\sqcup G\left[  C_{2}\right]  \sqcup\cdots\sqcup G\left[  C_{k}\right]
\right)  $ to $\operatorname*{V}\left(  G\right)  $ that sends each vertex
$\left(  i,v\right)  $ of $G\left[  C_{1}\right]  \sqcup G\left[
C_{2}\right]  \sqcup\cdots\sqcup G\left[  C_{k}\right]  $ to the vertex $v$ of
$G$. We claim that this bijection is a graph isomorphism. In order to prove
this, we need to check that there are no edges of $G$ that join vertices in
different components. But this is easy: If two vertices in different
components of $G$ were adjacent, then they would be path-connected, and thus
would actually belong to the same component.
\end{proof}

The upshot of these results is that every simple graph can be decomposed into
a disjoint union of its components (or, more precisely, of the induced
subgraphs on its components). Each of these components is a connected graph.
Moreover, this is easily seen to be the only way to decompose the graph into a
disjoint union of connected graphs.

\subsubsection{Some exercises on connectedness}

\begin{exercise}
\label{exe.paths.connectivity-by-split} Let $G$ be a simple graph with
$\operatorname{V}\left(  G\right)  \neq\varnothing$. Show that the following
two statements are equivalent:

\begin{itemize}
\item \textit{Statement 1:} The graph $G$ is connected.

\item \textit{Statement 2:} For every two nonempty subsets $A$ and $B$ of
$\operatorname{V}\left(  G\right)  $ satisfying $A\cap B=\varnothing$ and
$A\cup B=\operatorname{V}\left(  G\right)  $, there exist $a\in A$ and $b\in
B$ such that $ab\in\operatorname{E}\left(  G\right)  $. (In other words:
Whenever we subdivide the vertex set $\operatorname{V}\left(  G\right)  $ of
$G$ into two nonempty subsets, there will be at least one edge of $G$
connecting a vertex in one subset to a vertex in another.)
\end{itemize}

[\textbf{Solution:} This is Exercise 7 on homework set \#1 from my Spring 2017
course; see \href{https://www.cip.ifi.lmu.de/~grinberg/t/17s/}{the course
page} for solutions.]
\end{exercise}

\begin{exercise}
\label{exe.paths.connectivity-GH} Let $V$ be a nonempty finite set. Let $G$
and $H$ be two simple graphs such that $\operatorname{V}\left(  G\right)
=\operatorname{V}\left(  H\right)  =V$. Assume that for each $u\in V$ and
$v\in V$, there exists a path from $u$ to $v$ in $G$ or a path from $u$ to $v$
in $H$. Prove that at least one of the graphs $G$ and $H$ is
connected.\medskip

[\textbf{Solution:} This is Exercise 8 on homework set \#1 from my Spring 2017
course; see \href{https://www.cip.ifi.lmu.de/~grinberg/t/17s/}{the course
page} for solutions.]
\end{exercise}

\begin{exercise}
\label{exe.paths.qedmo4c1} Let $G=\left(  V,E\right)  $ be a simple graph. The
\textbf{complement graph} $\overline{G}$ of $G$ is defined to be the simple
graph $\left(  V,\ \mathcal{P}_{2}\left(  V\right)  \setminus E\right)  $.
(Thus, two distinct vertices $u$ and $v$ in $V$ are adjacent in $\overline{G}$
if and only if they are not adjacent in $G$.)

Prove that at least one of the following two statements holds:

\begin{itemize}
\item \textit{Statement 1:} For each $u\in V$ and $v\in V$, there exists a
path from $u$ to $v$ in $G$ of length $\leq3$.

\item \textit{Statement 2:} For each $u\in V$ and $v\in V$, there exists a
path from $u$ to $v$ in $\overline{G}$ of length $\leq2$.
\end{itemize}

[\textbf{Solution:} This is Exercise 9 on homework set \#1 from my Spring 2017
course; see \href{https://www.cip.ifi.lmu.de/~grinberg/t/17s/}{the course
page} for solutions.]
\end{exercise}

\begin{exercise}
\label{exe.1.4}Let $n\geq2$ be an integer. Let $G$ be a connected simple graph
with $n$ vertices.

\begin{enumerate}
\item[\textbf{(a)}] Describe $G$ if the degrees of the vertices of $G$ are
$1,1,2,2,\ldots,2$ (exactly two $1$'s and $n-2$ many $2$'s).

\item[\textbf{(b)}] Describe $G$ if the degrees of the vertices of $G$ are
$1,1,\ldots,1,n-1$.

\item[\textbf{(c)}] Describe $G$ if the degrees of the vertices of $G$ are
$2,2,\ldots,2$.
\end{enumerate}

Here, to \textquotedblleft describe\textquotedblright\ $G$ means to explicitly
determine (with proof) a graph that is isomorphic to $G$.
\end{exercise}

The following exercise is not explicitly concerned with connectedness and
components, but it might help to think about components to solve it (although
there are solutions that do not use them):

\begin{exercise}
\label{exe.1.6}Let $G$ be a simple graph with $n$ vertices. Assume that each
vertex of $G$ has at least one neighbor.

A \textbf{matching} of $G$ shall mean a set $F$ of edges of $G$ such that no
two edges in $F$ have a vertex in common. Let $m$ be the largest size of a
matching of $G$.

An \textbf{edge cover} of $G$ shall mean a set $F$ of edges of $G$ such that
each vertex of $G$ is contained in at least one edge $e\in F$. Let $c$ be the
smallest size of an edge cover of $G$.

Prove that $c+m=n$.
\end{exercise}

\begin{remark}
Let $G$ be the cycle graph $C_{5}$ shown in Example \ref{exa.dom.C5}. Then,
$\left\{  12,34\right\}  $ is a matching of $G$ of largest possible size
(why?), whereas $\left\{  12,34,25\right\}  $ is an edge cover of $G$ of
smallest possible size (why?). Thus, Exercise \ref{exe.1.6} says that $2+3=5$
here, which is indeed true.
\end{remark}

\subsection{\label{sec.sg.cycs}Closed walks and cycles}

Here are two further kinds of walks:

\begin{definition}
Let $G$ be a simple graph.

\begin{enumerate}
\item[\textbf{(a)}] A \textbf{closed walk} of $G$ means a walk whose first
vertex is identical with its last vertex. In other words, it means a walk
$\left(  w_{0},w_{1},\ldots,w_{k}\right)  $ with $w_{0}=w_{k}$. Sometimes,
closed walks are also known as \textbf{circuits} (but many authors use this
latter word for something slightly different).

\item[\textbf{(b)}] A \textbf{cycle} of $G$ means a closed walk $\left(
w_{0},w_{1},\ldots,w_{k}\right)  $ such that $k\geq3$ and such that the
vertices $w_{0},w_{1},\ldots,w_{k-1}$ are distinct.
\end{enumerate}
\end{definition}

\begin{example}
\label{exa.cycles.1}Let $G$ be the simple graph%
\[
\left(  \left\{  1,2,3,4,5,6\right\}  ,\ \ \left\{
12,\ 23,\ 34,\ 45,\ 56,\ 61,\ 13\right\}  \right)  .
\]
This graph looks as follows (we have already seen it in Example
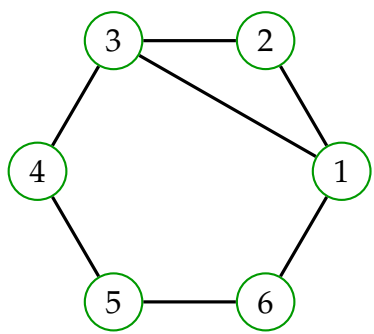
\ref{exa.sg.walk.example1}):%
\[%
\begin{tikzpicture}
\begin{scope}[every node/.style={circle,thick,draw=green!60!black}]
\node(A) at (0:2) {$1$};
\node(B) at (60:2) {$2$};
\node(C) at (120:2) {$3$};
\node(D) at (180:2) {$4$};
\node(E) at (240:2) {$5$};
\node(F) at (300:2) {$6$};
\end{scope}
\begin{scope}[every edge/.style={draw=black,very thick}]
\path
[-] (A) edge (B) (B) edge (C) (C) edge (A) (C) edge (D) (D) edge (E) (E) edge (F) (F) edge (A);
\end{scope}
\end{tikzpicture}%
\]
Then:

\begin{itemize}
\item The sequence $\left(  1,3,2,1,6,5,6,1\right)  $ is a closed walk of $G$.
But it is very much not a cycle.

\item The sequences $\left(  1,2,3,1\right)  $ and $\left(
1,3,4,5,6,1\right)  $ and $\left(  1,2,3,4,5,6,1\right)  $ are cycles of $G$.
You can get further cycles by rotating these sequences (in a proper sense of
this word -- e.g., rotating $\left(  1,2,3,1\right)  $ gives $\left(
2,3,1,2\right)  $ and $\left(  3,1,2,3\right)  $) and by reversing them. Every
cycle of $G$ can be obtained in this way.

\item The sequences $\left(  1\right)  $ and $\left(  1,2,1\right)  $ are
closed walks, but not cycles of $G$ (since they fail the $k\geq3$ condition).

\item The sequence $\left(  1,2,3\right)  $ is a walk, but not a closed walk,
since $1\neq3$.
\end{itemize}
\end{example}

Authors have different opinions about whether $\left(  1,2,3,1\right)  $ and
$\left(  1,3,2,1\right)  $ count as different cycles. Fortunately, this
matters only if you want to count cycles, but not for the existence or
non-existence of cycles.

We have now defined paths (in an arbitrary graph) and also path graphs $P_{n}%
$; we have also defined cycles (in an arbitrary graph) and also cycle graphs
$C_{n}$. Besides their similar names, are they related? The answer is
\textquotedblleft yes\textquotedblright:

\begin{proposition}
Let $G$ be a simple graph.

\begin{enumerate}
\item[\textbf{(a)}] If $\left(  p_{0},p_{1},\ldots,p_{k}\right)  $ is a path
of $G$, then there is a subgraph of $G$ isomorphic to the path graph $P_{k+1}%
$, namely the subgraph $\left(  \left\{  p_{0},p_{1},\ldots,p_{k}\right\}
,\ \left\{  p_{i}p_{i+1}\ \mid\ 0\leq i<k\right\}  \right)  $. (If this
subgraph is actually an induced subgraph of $G$, then the path $\left(
p_{0},p_{1},\ldots,p_{k}\right)  $ is called an \textquotedblleft induced
path\textquotedblright.)

Conversely, any subgraph of $G$ isomorphic to $P_{k+1}$ gives a path of $G$.

\item[\textbf{(b)}] Now, assume that $k\geq3$. If $\left(  c_{0},c_{1}%
,\ldots,c_{k}\right)  $ is a cycle of $G$, then there is a subgraph of $G$
isomorphic to the cycle graph $C_{k}$, namely the subgraph $\left(  \left\{
c_{0},c_{1},\ldots,c_{k}\right\}  ,\ \left\{  c_{i}c_{i+1}\ \mid\ 0\leq
i<k\right\}  \right)  $. (If this subgraph is actually an induced subgraph of
$G$, then the cycle $\left(  c_{0},c_{1},\ldots,c_{k}\right)  $ is called an
\textquotedblleft induced cycle\textquotedblright.)

Conversely, any subgraph of $G$ isomorphic to $C_{k}$ gives a cycle of $G$.
\end{enumerate}
\end{proposition}

\begin{proof}
Straightforward.
\end{proof}

Certain graphs contain cycles; other graphs don't. For instance, the complete
graph $K_{n}$ contains a lot of cycles (when $n\geq3$), whereas the path graph
$P_{n}$ contains none. Let us try to find some criteria for when a graph can
and when it cannot have cycles\footnote{Mantel's theorem already gives such a
criterion for cycles of length $3$ (because a cycle of length $3$ is the same
as a triangle).}:

\begin{proposition}
\label{prop.cyc.btf-walk-cyc}Let $G$ be a simple graph. Let $\mathbf{w}$ be a
walk of $G$ such that no two consecutive edges of $\mathbf{w}$ are identical.
(By \textquotedblleft consecutive edges\textquotedblright, we mean edges of
the form $w_{i-1}w_{i}$ and $w_{i}w_{i+1}$, where $w_{i-1},w_{i},w_{i+1}$ are
three consecutive vertices of $\mathbf{w}$.)

Then, $\mathbf{w}$ either is a path or contains a cycle (i.e., there exists a
cycle of $G$ whose edges are edges of $\mathbf{w}$).
\end{proposition}

\begin{example}
Let $G$ be as in Example \ref{exa.cycles.1}. Then, $\left(
2,1,3,2,1,6\right)  $ is a walk $\mathbf{w}$ of $G$ such that no two
consecutive edges of $\mathbf{w}$ are identical (even though the edge $21$
appears twice in this walk). On the other hand, $\left(  2,1,3,1,6\right)  $
is not such a walk (since its two consecutive edges $13$ and $31$ are identical).
\end{example}

\begin{proof}
[Proof of Proposition \ref{prop.cyc.btf-walk-cyc}.]We assume that $\mathbf{w}$
is not a path. We must then show that $\mathbf{w}$ contains a cycle.

Write $\mathbf{w}$ as $\mathbf{w}=\left(  w_{0},w_{1},\ldots,w_{k}\right)  $.
Since $\mathbf{w}$ is not a path, two of the vertices $w_{0},w_{1}%
,\ldots,w_{k}$ must be equal. In other words, there exists a pair $\left(
i,j\right)  $ of integers $i$ and $j$ with $i<j$ and $w_{i}=w_{j}$. Among all
such pairs, we pick one with \textbf{minimum} difference $j-i$. We shall show
that the walk $\left(  w_{i},w_{i+1},\ldots,w_{j}\right)  $ is a cycle.

First, this walk is clearly a closed walk (since $w_{i}=w_{j}$). It thus
remains to show that $j-i\geq3$ and that the vertices $w_{i},w_{i+1}%
,\ldots,w_{j-1}$ are distinct. The distinctness of $w_{i},w_{i+1}%
,\ldots,w_{j-1}$ follows from the minimality of $j-i$. To show that $j-i\geq
3$, we assume the contrary. Thus, $j-i$ is either $1$ or $2$ (since $i<j$).
But $j-i$ cannot be $1$, since the endpoints of an edge cannot be equal (since
our graph is a simple graph). So $j-i$ must be $2$. Thus, $w_{i}=w_{i+2}$.
Therefore, the two edges $w_{i}w_{i+1}$ and $w_{i+1}w_{i+2}$ are identical.
But this contradicts the fact that no two consecutive edges of $\mathbf{w}$
are identical. Contradiction, qed.
\end{proof}

\begin{corollary}
\label{cor.cyc.btw-clow-cyc}Let $G$ be a simple graph. Assume that $G$ has a
closed walk $\mathbf{w}$ of length $>0$ such that no two consecutive edges of
$\mathbf{w}$ are identical. Then, $G$ has a cycle.
\end{corollary}

\begin{proof}
This follows from Proposition \ref{prop.cyc.btf-walk-cyc}, since $\mathbf{w}$
is not a path.
\end{proof}

\begin{theorem}
\label{thm.cyc.two-paths-cyc}Let $G$ be a simple graph. Let $u$ and $v$ be two
vertices in $G$. Assume that there are two distinct paths from $u$ to $v$.
Then, $G$ has a cycle.
\end{theorem}

\begin{proof}
More generally, we shall prove this theorem with the word \textquotedblleft
path\textquotedblright\ replaced by \textquotedblleft backtrack-free
walk\textquotedblright, where a \textquotedblleft\textbf{backtrack-free
walk}\textquotedblright\ means a walk $\mathbf{w}$ such that no two
consecutive edges of $\mathbf{w}$ are identical. This is a generalization of
the theorem, since every path is a backtrack-free walk (why?).

So we claim the following:

\begin{statement}
\textit{Claim 1:} Let $\mathbf{p}$ and $\mathbf{q}$ be two distinct
backtrack-free walks that start at the same vertex and end at the same vertex.
Then, $G$ has a cycle.
\end{statement}

We shall prove Claim 1 by induction on the length of $\mathbf{p}$:

\textit{Base case:} As the base case, we note that Claim 1 is vacuously true
when the length of $\mathbf{p}$ is $-1$ (since there are no walks of length
$-1$).

\textit{Induction step:} We fix an integer $N$, and we assume that Claim 1 is
proved in the case when the length of $\mathbf{p}$ is $N-1$. We must now show
that it is also true when the length of $\mathbf{p}$ is $N$.

So let $\mathbf{p}=\left(  p_{0},p_{1},\ldots,p_{a}\right)  $ and
$\mathbf{q}=\left(  q_{0},q_{1},\ldots,q_{b}\right)  $ be two distinct
backtrack-free walks that start at the same vertex and end at the same vertex
and satisfy $a=N$. We must find a cycle.

The walks $\mathbf{p}$ and $\mathbf{q}$ are distinct but start at the same
vertex, so they cannot both be trivial\footnote{We say that a walk is
\textbf{trivial} if it has length $0$.}. If one of them is trivial, then the
other is a closed walk (because a trivial walk is a closed walk), and then our
goal follows from Corollary \ref{cor.cyc.btw-clow-cyc} in this case (because
we have a nontrivial closed backtrack-free walk). Hence, from now on, we WLOG
assume that neither of the two walks $\mathbf{p}$ and $\mathbf{q}$ is trivial.
Thus, each of these two walks has a last edge. The last edge of $\mathbf{p}$
is $p_{a-1}p_{a}$, whereas the last edge of $\mathbf{q}$ is $q_{b-1}q_{b}$.

Two cases are possible:

\textit{Case 1:} We have $p_{a-1}p_{a}=q_{b-1}q_{b}$.

\textit{Case 2:} We have $p_{a-1}p_{a}\neq q_{b-1}q_{b}$.

Let us consider Case 1 first. In this case, the last edges $p_{a-1}p_{a}$ and
$q_{b-1}q_{b}$ of the two walks $\mathbf{p}$ and $\mathbf{q}$ are identical,
so the second-to-last vertices of these two walks must also be identical.
Thus, if we remove these last edges from both walks, then we obtain two
shorter backtrack-free walks $\left(  p_{0},p_{1},\ldots,p_{a-1}\right)  $ and
$\left(  q_{0},q_{1},\ldots,q_{b-1}\right)  $ that again start at the same
vertex and end at the same vertex, but the length of the first of them is
$a-1=N-1$. Hence, by the induction hypothesis, we can apply Claim 1 to these
two shorter walks (instead of $\mathbf{p}$ and $\mathbf{q}$), and we conclude
that $G$ has a cycle. So we are done in Case 1.

Let us now consider Case 2. In this case, we combine the two walks
$\mathbf{p}$ and $\mathbf{q}$ (more precisely, $\mathbf{p}$ and the reversal
of $\mathbf{q}$) to obtain the closed walk%
\[
\left(  p_{0},p_{1},\ldots,p_{a-1},p_{a}=q_{b},q_{b-1},\ldots,q_{0}\right)  .
\]
This closed walk is backtrack-free (since $\left(  p_{0},p_{1},\ldots
,p_{a}\right)  $ and $\left(  q_{0},q_{1},\ldots,q_{b}\right)  $ are
backtrack-free, and since $p_{a-1}p_{a}\neq q_{b-1}q_{b}$) and has length $>0$
(since it contains at least the edge $p_{a-1}p_{a}$). Hence, Corollary
\ref{cor.cyc.btw-clow-cyc} entails that $G$ has a cycle.

We have thus found a cycle in both Cases 1 and 2. This completes the induction
step. Thus, we have proved Claim 1. As we said, Theorem
\ref{thm.cyc.two-paths-cyc} follows from it.
\end{proof}

\begin{exercise}
\label{exe.2.1}Let $G$ be a simple graph.

\begin{enumerate}
\item[\textbf{(a)}] Prove that if $G$ has a closed walk of odd length, then
$G$ has a cycle of odd length.

\item[\textbf{(b)}] Is it true that if $G$ has a closed walk of length not
divisible by $3$, then $G$ has a cycle of length not divisible by $3$ ?

\item[\textbf{(c)}] Does the answer to part \textbf{(b)} change if we replace
\textquotedblleft walk\textquotedblright\ by \textquotedblleft
non-backtracking walk\textquotedblright? (A walk $\mathbf{w}$ with edges
$e_{1},e_{2},\ldots,e_{k}$ (in this order) is said to be
\textbf{non-backtracking} if each $i\in\left\{  1,2,\ldots,k-1\right\}  $
satisfies $e_{i}\neq e_{i+1}$.)

\item[\textbf{(d)}] A \textbf{trail} (in a graph) means a walk whose edges are
distinct (but whose vertices are not necessarily distinct). Does the answer to
part \textbf{(b)} change if we replace \textquotedblleft
walk\textquotedblright\ by \textquotedblleft trail\textquotedblright?
\end{enumerate}

(Proofs and counterexamples should be given.)
\end{exercise}

\subsection{\label{sec.sg.lpt}The longest path trick}

Here is another proposition that guarantees the existence of cycles in a graph
under certain circumstances. More importantly, its proof illustrates a useful
tactic in dealing with graphs:

\begin{proposition}
\label{prop.cycle.deg-d}Let $G$ be a simple graph with at least one vertex.
Let $d>1$ be an integer. Assume that each vertex of $G$ has degree $\geq d$.
Then, $G$ has a cycle of length $\geq d+1$.
\end{proposition}

\begin{proof}
Let $\mathbf{p}=\left(  v_{0},v_{1},\ldots,v_{m}\right)  $ be a
\textbf{longest} path of $G$. (Why does $G$ have a longest path? Let's see:
Any path of $G$ has length $\leq\left\vert V\right\vert -1$, since its
vertices have to be distinct. Moreover, $G$ has at least one vertex and thus
has at least one path. A finite nonempty set of integers has a largest
element. Thus, $G$ has a longest path.)

The vertex $v_{0}$ has degree $\geq d$ (by assumption), and thus has $\geq d$
neighbors (since the degree of a vertex is the number of its neighbors).

If all neighbors of $v_{0}$ belonged to the set $\left\{  v_{1},v_{2}%
,\ldots,v_{d-1}\right\}  $\ \ \ \ \footnote{If $d-1>m$, then this set should
be understood to mean $\left\{  v_{1},v_{2},\ldots,v_{m}\right\}  $.}, then
the number of neighbors of $v_{0}$ would be at most $d-1$, which would
contradict the previous sentence. Thus, there exists at least one neighbor $u$
of $v_{0}$ that does \textbf{not} belong to this set $\left\{  v_{1}%
,v_{2},\ldots,v_{d-1}\right\}  $. Consider this $u$. Then, $u\neq v_{0}$
(since a vertex cannot be its own neighbor).

Attaching the vertex $u$ to the front of the path $\mathbf{p}$, we obtain a
walk%
\[
\mathbf{p}^{\prime}:=\left(  u,v_{0},v_{1},\ldots,v_{m}\right)  .
\]
If we had $u\notin\left\{  v_{0},v_{1},\ldots,v_{m}\right\}  $, then this walk
$\mathbf{p}^{\prime}$ would be a path; but this would contradict the fact that
$\mathbf{p}$ is a \textbf{longest} path of $G$. Thus, we must have
$u\in\left\{  v_{0},v_{1},\ldots,v_{m}\right\}  $. In other words, $u=v_{i}$
for some $i\in\left\{  0,1,\ldots,m\right\}  $. Consider this $i$. Since
$u\neq v_{0}$ and $u\notin\left\{  v_{1},v_{2},\ldots,v_{d-1}\right\}  $, we
thus have $i\geq d$. Here is a picture:%
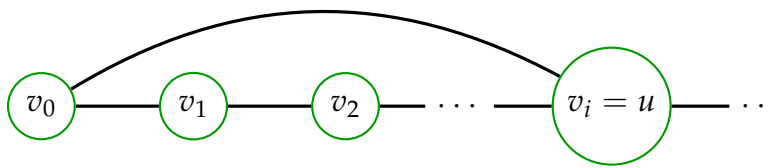
\[%
\begin{tikzpicture}
\begin{scope}[every node/.style={circle,thick,draw=green!60!black}]
\node(A) at (0,0) {$v_0$};
\node(B) at (2,0) {$v_1$};
\node(C) at (4,0) {$v_2$};
\node(E) at (7.5,0) {$v_i=u$};
\end{scope}
\node(D) at (5.5,0) {$\cdots$};
\node(F) at (9.5,0) {$\cdots$};
\begin{scope}[every edge/.style={draw=black,very thick}]
\path
[-] (A) edge (B) (B) edge (C) (C) edge (D) (D) edge (E) (E) edge (F) (A) edge[bend left=30] (E);
\end{scope}
\end{tikzpicture}%
\]

Now, consider the walk%
\[
\mathbf{c}:=\left(  u,v_{0},v_{1},\ldots,v_{i}\right)  .
\]
This is a closed walk (since $u=v_{i}$) and has length $i+1\geq d+1$ (since
$i\geq d$). If we can show that $\mathbf{c}$ is a cycle, then we have thus
found a cycle of length $\geq d+1$, so we will be done.

It thus remains to prove that $\mathbf{c}$ is a cycle. Let us do this. We need
to check that the vertices $u,v_{0},v_{1},\ldots,v_{i-1}$ are distinct, and
that the length of $\mathbf{c}$ is $\geq3$. The latter claim is clear: The
length of $\mathbf{c}$ is $i+1\geq d+1\geq3$ (since $d>1$ and $d\in\mathbb{Z}%
$). The former claim is not much harder: Since $u=v_{i}$, the vertices
$u,v_{0},v_{1},\ldots,v_{i-1}$ are just the vertices $v_{i},v_{0},v_{1}%
,\ldots,v_{i-1}$, and thus are distinct because they are distinct vertices of
the path $\mathbf{p}$. The proof of Proposition \ref{prop.cycle.deg-d} is thus complete.
\end{proof}

\subsection{\label{sec.sg.bridges}Bridges}

One question that will later prove crucial is: What happens to a graph if we
remove a single edge from it? Let us first define a notation for this:

\begin{definition}
\label{def.sg.G-e}Let $G=\left(  V,E\right)  $ be a simple graph. Let $e$ be
an edge of $G$. Then, $G\setminus e$ will mean the graph obtained from $G$ by
removing this edge $e$. In other words,%
\[
G\setminus e:=\left(  V,\ E\setminus\left\{  e\right\}  \right)  .
\]

\end{definition}

Some authors write $G-e$ for $G\setminus e$.

\begin{theorem}
\label{thm.G-e.conn}Let $G$ be a simple graph. Let $e$ be an edge of $G$. Then:

\begin{enumerate}
\item[\textbf{(a)}] If $e$ is an edge of some cycle of $G$, then the
components of $G\setminus e$ are precisely the components of $G$. (Keep in
mind that the components are sets of vertices. It is these sets that we are
talking about here, not the induced subgraphs on these sets.)

\item[\textbf{(b)}] If $e$ appears in no cycle of $G$ (in other words, there
exists no cycle of $G$ such that $e$ is an edge of this cycle), then the graph
$G\setminus e$ has one more component than $G$.
\end{enumerate}
\end{theorem}

\begin{example}
Let $G$ be the graph shown in the following picture:%
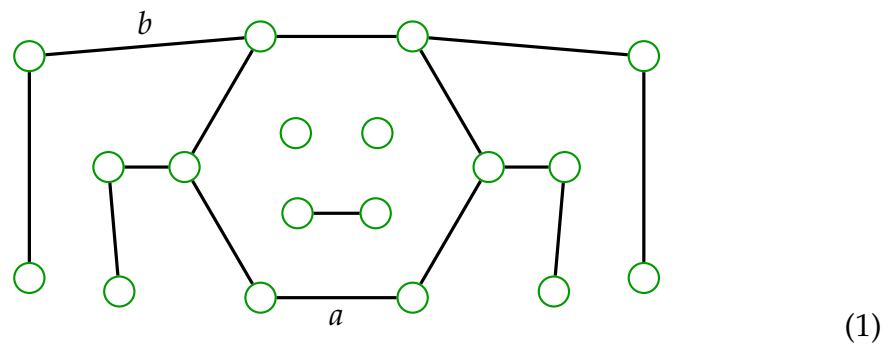
\begin{equation}%

\label{eq.exa.lem.graph.Gwoe-conn.1}%
\end{equation}
(where we have labeled the edges $a$ and $b$ for further reference). This
graph has $4$ components. The edge $a$ is an edge of a cycle of $G$, whereas
the edge $b$ appears in no cycle of $G$. Thus, if we set $e=a$, then Theorem
\ref{thm.G-e.conn} \textbf{(a)} shows that the components of $G\setminus e$
are precisely the components of $G$. This graph $G\setminus e$ for $e=a$ looks
as follows:%
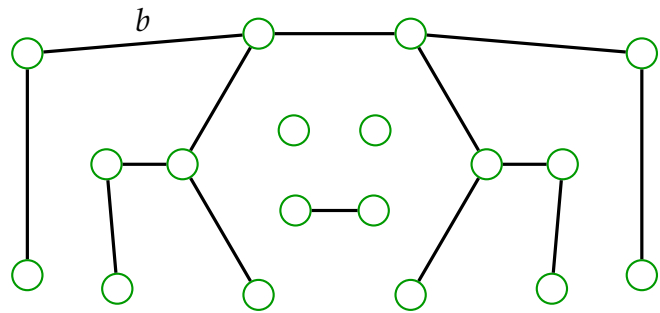
\[%
%
\]
and visibly has $4$ components. On the other hand, if we set $e=b$, then
Theorem \ref{thm.G-e.conn} \textbf{(b)} shows that the graph $G\setminus e$
has one more component than $G$. This graph $G\setminus e$ for $e=b$ looks as
follows:%
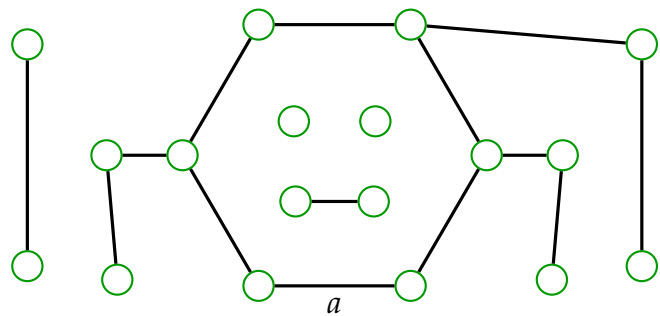
\[%
%
\]
and visibly has $5$ components.
\end{example}

\begin{proof}
[Proof of Theorem \ref{thm.G-e.conn}.]We will only sketch the proof. For
details, see \cite[\S 6.7]{21f6}.

Let $u$ and $v$ be the endpoints of $e$, so that $e=uv$. Note that $\left(
u,v\right)  $ is a path of $G$, and thus we have $u\simeq_{G}v$. \medskip

\textbf{(a)} Assume that $e$ is an edge of some cycle of $G$. Then, if you
remove $e$ from this cycle, then you still have a path from $u$ to $v$ left
(as the remaining edges of the cycle function as a detour), and this path is a
path of $G\setminus e$. Thus, $u\simeq_{G\setminus e}v$.

Now, we must show that the components of $G\setminus e$ are precisely the
components of $G$. This will clearly follow if we can show that the relation
$\simeq_{G\setminus e}$ is precisely the relation $\simeq_{G}$ (because the
components of a graph are the equivalence classes of its $\simeq$ relation).
So let us prove the latter fact.

We must show that two vertices $x$ and $y$ of $G$ satisfy $x\simeq_{G\setminus
e}y$ if and only if they satisfy $x\simeq_{G}y$. The \textquotedblleft only
if\textquotedblright\ part is obvious (since a walk of $G\setminus e$ is
always a walk of $G$). It thus remains to prove the \textquotedblleft
if\textquotedblright\ part. So we assume that $x$ and $y$ are two vertices of
$G$ satisfying $x\simeq_{G}y$, and we want to show that $x\simeq_{G\setminus
e}y$.

From $x\simeq_{G}y$, we conclude that $G$ has a path from $x$ to $y$ (by
Proposition \ref{prop.sg.pc.there-is-path}). If this path does not
use\footnote{We say that a walk $\mathbf{w}$ \textbf{uses} an edge $f$ if $f$
is an edge of $\mathbf{w}$.} the edge $e$, then it is a path from $x$ to $y$
in $G\setminus e$, and thus we have $x\simeq_{G\setminus e}y$, which is what
we wanted to prove. So we WLOG assume that this path does use the edge $e$.
Thus, this path contains the endpoints $u$ and $v$ of this edge $e$. We WLOG
assume that $u$ appears before $v$ on this path (otherwise, just swap $u$ with
$v$). Thus, this path looks as follows:%
\[
\left(  x,\ldots,u,v,\ldots,y\right)  .
\]
If we remove the edge $e=uv$, then this path breaks into two smaller paths%
\[
\left(  x,\ldots,u\right)  \ \ \ \ \ \ \ \ \ \ \text{and}%
\ \ \ \ \ \ \ \ \ \ \left(  v,\ldots,y\right)
\]
(since the edges of a path are distinct, so $e$ appears only once in it). Both
of these two smaller paths are paths of $G\setminus e$. Thus, $x\simeq
_{G\setminus e}u$ and $v\simeq_{G\setminus e}y$. Now, recalling that
$\simeq_{G\setminus e}$ is an equivalence relation, we combine these results
to obtain%
\[
x\simeq_{G\setminus e}u\simeq_{G\setminus e}v\simeq_{G\setminus e}y.
\]
Hence, $x\simeq_{G\setminus e}y$. This completes the proof of Theorem
\ref{thm.G-e.conn} \textbf{(a)}. \medskip

\textbf{(b)} Assume that $e$ appears in no cycle of $G$. We must prove that
the graph $G\setminus e$ has one more component than $G$. To do so, it
suffices to show the following:

\begin{statement}
\textit{Claim 1:} The component of $G$ that contains $u$ and $v$ (this
component exists, since $u\simeq_{G}v$) breaks into two components of
$G\setminus e$ when the edge $e$ is removed.
\end{statement}

\begin{statement}
\textit{Claim 2:} All other components of $G$ remain components of $G\setminus
e$.
\end{statement}

Claim 2 is pretty clear: The components of $G$ that don't contain $u$ and $v$
do not change at all when $e$ is removed (since they contain neither endpoint
of $e$). Thus, they remain components of $G\setminus e$. (Formalizing this is
a nice exercise in formalization; see \cite[\S 6.7]{21f6}.)

It remains to prove Claim 1. We introduce some notations:

\begin{itemize}
\item Let $C$ be the component of $G$ that contains $u$ and $v$.

\item Let $A$ be the component of $G\setminus e$ that contains $u$.

\item Let $B$ be the component of $G\setminus e$ that contains $v$.
\end{itemize}

Then, we must show that $A\cup B=C$ and $A\cap B=\varnothing$.

To see that $A\cap B=\varnothing$, we need to show that $u\simeq_{G\setminus
e}v$ does \textbf{not} hold (since $A$ and $B$ are the equivalence classes of
$u$ and $v$ with respect to the relation $\simeq_{G\setminus e}$). So let us
do this. Assume the contrary. Thus, $u\simeq_{G\setminus e}v$. Hence, there
exists a path from $u$ to $v$ in $G\setminus e$. Since $e=uv$, we can
\textquotedblleft close\textquotedblright\ this path by appending the vertex
$u$ to its end; the result is a cycle of the graph $G$ that contains the edge
$e$. But this contradicts our assumption that no cycle of $G$ contains $e$.
This contradiction shows that our assumption was wrong. Thus, we conclude that
$u\simeq_{G\setminus e}v$ does \textbf{not} hold. Hence, as we said, $A\cap
B=\varnothing$.

It remains to show that $A\cup B=C$. Since $A$ and $B$ are clearly subsets of
$C$ (because each walk of $G\setminus e$ is a walk of $G$, and thus each
component of $G\setminus e$ is a subset of a component of $G$), we have $A\cup
B\subseteq C$, and therefore we only need to show that $C\subseteq A\cup B$.
In other words, we need to show that each $c\in C$ belongs to $A\cup B$.

Let us show this. Let $c\in C$ be a vertex. Then, $c\simeq_{G}u$ (since $C$ is
the component of $G$ containing $u$). Therefore, $G$ has a path $\mathbf{p}$
from $c$ to $u$. Consider this path $\mathbf{p}$. Two cases are possible:

\begin{itemize}
\item \textit{Case 1:} This path $\mathbf{p}$ does not use the edge $e$. In
this case, $\mathbf{p}$ is a path of $G\setminus e$, and thus we obtain
$c\simeq_{G\setminus e}u$. In other words, $c\in A$ (since $A$ is the
component of $G\setminus e$ containing $u$).

\item \textit{Case 2:} This path $\mathbf{p}$ does use the edge $e$. In this
case, the edge $e$ must be the last edge of $\mathbf{p}$ (since the path
$\mathbf{p}$ would otherwise contain the vertex $u$ twice\footnote{Indeed, the
path $\mathbf{p}$ already ends in $u$. If it would contain $e$ anywhere other
than at the very end, then it would thus contain the vertex $u$ twice (since
$u$ is an endpoint of $e$).}; but a path cannot contain a vertex twice), and
the last two vertices of $\mathbf{p}$ must be $v$ and $u$ in this order. Thus,
by removing the last vertex from $\mathbf{p}$, we obtain a path from $c$ to
$v$, and this latter path is a path of $G\setminus e$ (since it no longer
contains $u$ and therefore does not use $e$). This yields $c\simeq_{G\setminus
e}v$. In other words, $c\in B$ (since $B$ is the component of $G\setminus e$
containing $v$).
\end{itemize}

\noindent In either of these two cases, we have shown that $c$ belongs to one
of $A$ and $B$. In other words, $c\in A\cup B$. This is precisely what we
wanted to show. This completes the proof of Theorem \ref{thm.G-e.conn}
\textbf{(b)}.
\end{proof}

We introduce some fairly standard terminology:

\begin{definition}
\label{def.sg.bridges-cuts}Let $e$ be an edge of a simple graph $G$.

\begin{enumerate}
\item[\textbf{(a)}] We say that $e$ is a \textbf{bridge} (of $G$) if $e$
appears in no cycle of $G$.

\item[\textbf{(b)}] We say that $e$ is a \textbf{cut-edge} (of $G$) if the
graph $G\setminus e$ has more components than $G$.
\end{enumerate}
\end{definition}

\begin{corollary}
\label{cor.bridge=cut}Let $e$ be an edge of a simple graph $G$. Then, $e$ is a
bridge if and only if $e$ is a cut-edge.
\end{corollary}

\begin{proof}
Follows from Theorem \ref{thm.G-e.conn}.
\end{proof}

We can also define \textquotedblleft cut-vertices\textquotedblright: A vertex
$v$ of a graph $G$ is said to be a \textbf{cut-vertex} if the graph
$G\setminus v$ (that is, the graph $G$ with the vertex $v$
removed\footnote{When we remove a vertex, we must of course also remove all
edges that contain this vertex.}) has more components than $G$. Unfortunately,
there doesn't seem to be an analogue of Corollary \ref{cor.bridge=cut} for
cut-vertices. Note also that removing a vertex (unlike removing an edge) can
add more than one component to the graph (or it can also subtract $1$
component if this vertex had degree $0$). For example, removing the vertex $0$
from the graph%
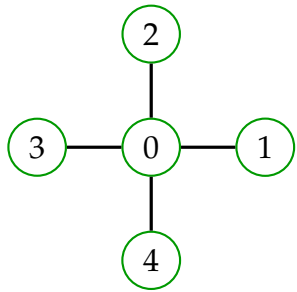
\[%
\begin{tikzpicture}
\begin{scope}[every node/.style={circle,thick,draw=green!60!black}]
\node(1) at (0:1.5) {$1$};
\node(2) at (90:1.5) {$2$};
\node(3) at (180:1.5) {$3$};
\node(4) at (270:1.5) {$4$};
\node(0) at (0,0) {$0$};
\end{scope}
\begin{scope}[every edge/.style={draw=black,very thick}]
\path[-] (0) edge (1) (0) edge (2) (0) edge (3) (0) edge (4);
\end{scope}
\end{tikzpicture}%
\]
results in an empty graph on the set $\left\{  1,2,3,4\right\}  $, so the
number of components has increased from $1$ to $4$.

\subsection{\label{sec.sg.dom}Dominating sets}

\subsubsection{Definition and basic facts}

Here is another concept we can define for a graph:

\begin{definition}
\label{def.sg.dom}Let $G=\left(  V,E\right)  $ be a simple graph.

A subset $U$ of $V$ is said to be \textbf{dominating} (for $G$) if it has the
following property: Each vertex $v\in V\setminus U$ has at least one neighbor
in $U$.

A \textbf{dominating set} for $G$ (or \textbf{dominating set} of $G$) will
mean a subset of $V$ that is dominating.
\end{definition}

\begin{example}
\label{exa.dom.C5} Consider the cycle graph
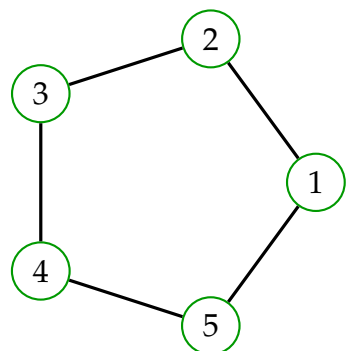
\[
C_{5}=\left(  \left\{  1,2,3,4,5\right\}  ,\ \ \left\{
12,\ 23,\ 34,\ 45,\ 51\right\}  \right)  =%
\begin{tikzpicture}
\begin{scope}[every node/.style={circle,thick,draw=green!60!black}]
\node(A) at (0:2) {$1$};
\node(B) at (360/5:2) {$2$};
\node(C) at (2*360/5:2) {$3$};
\node(D) at (3*360/5:2) {$4$};
\node(E) at (4*360/5:2) {$5$};
\end{scope}
\begin{scope}[every edge/.style={draw=black,very thick}]
\path[-] (A) edge (B) (B) edge (C) (C) edge (D) (D) edge (E) (E) edge (A);
\end{scope}
\end{tikzpicture}%
.
\]
The set $\left\{  1,3\right\}  $ is a dominating set for $C_{5}$, since all
three vertices $2,4,5$ that don't belong to $\left\{  1,3\right\}  $ have
neighbors in $\left\{  1,3\right\}  $. The set $\left\{  1,5\right\}  $ is not
a dominating set for $C_{5}$, since the vertex $3$ has no neighbor in
$\left\{  1,5\right\}  $. There is no dominating set for $C_{5}$ that has size
$0$ or $1$, but there are several of size $2$, and every subset of size
$\geq3$ is dominating.
\end{example}

Here are some more examples:

\begin{itemize}
\item If $G=\left(  V,E\right)  $ is a simple graph, then the whole vertex set
$V$ is always dominating, whereas the empty set $\varnothing$ is dominating
only when $V=\varnothing$.

\item If $G=\left(  V,E\right)  $ is a complete graph, then any nonempty
subset of $V$ is dominating.

\item If $G=\left(  V,E\right)  $ is an empty graph, then only $V$ is dominating.
\end{itemize}

Clearly, the \textquotedblleft denser\textquotedblright\ a graph is (i.e., the
more edges it has), the \textquotedblleft easier\textquotedblright\ it is for
a set to be dominating. Often, a graph is given, and one is interested in
finding a dominating set of the smallest possible size\footnote{Supposedly,
this has applications in mobile networking: For example, you might want to
choose a set of routers in a given network so that each node is either a
router or directly connected (i.e., adjacent) to one.}. As the case of an
empty graph reveals, sometimes the only choice is the whole vertex set.
However, in many cases, we can do better. Namely, we need to require that the
graph has no isolated vertices:

\begin{definition}
Let $G$ be a simple graph. A vertex $v$ of $G$ is said to be \textbf{isolated}
if it has no neighbors (i.e., if $\deg v=0$).
\end{definition}

An isolated vertex has to belong to every dominating set (since otherwise, it
would need a neighbor in that set, but it has no neighbors). Thus, isolated
vertices do not contribute much to the study of dominating sets, other than
inflating their size. Therefore, when we look for dominating sets, we can
restrict ourselves to graphs with no isolated vertices. There, we have the
following result:

\begin{proposition}
\label{prop.dom.V/2}Let $G=\left(  V,E\right)  $ be a simple graph that has no
isolated vertices. Then:

\begin{enumerate}
\item[\textbf{(a)}] There exists a dominating subset of $V$ that has size
$\leq\left\vert V\right\vert /2$.

\item[\textbf{(b)}] There exist two disjoint dominating subsets $A$ and $B$ of
$V$ such that $A\cup B=V$.
\end{enumerate}
\end{proposition}

One proof of this proposition will be given in Exercise \ref{exe.2.4} below
(homework set \#2 exercise 4). Another appears in \cite[\S 3.6]{17s}.

For specific graphs, the bound $\left\vert V\right\vert /2$ in Proposition
\ref{prop.dom.V/2} \textbf{(a)} can often be improved. Here is an example:

\begin{exercise}
\label{exe.2.2}Let $n\geq3$ be an integer. Find a formula for the smallest
size of a dominating set of the cycle graph $C_{n}$. You can use the
\textbf{ceiling function} $x\mapsto\left\lceil x\right\rceil $, which sends a
real number $x$ to the smallest integer that is $\geq x$.
\end{exercise}

\begin{exercise}
\label{exe.2.3}Let $n$ and $k$ be positive integers such that $n\geq k\left(
k+1\right)  $ and $k>1$. Recall (from Subsection
\ref{subsec.sg.complete.kneser}) the Kneser graph $KG_{n,k}$, whose vertices
are the $k$-element subsets of $\left\{  1,2,\ldots,n\right\}  $, and whose
edges are the unordered pairs $\left\{  A,B\right\}  $ of such subsets with
$A\cap B=\varnothing$.

Prove that the minimum size of a dominating set of $KG_{n,k}$ is $k+1$.
\end{exercise}

\begin{exercise}
\label{exe.2.4}Let $G=\left(  V,E\right)  $ be a connected simple graph with
at least two vertices.

The \textbf{distance} $d\left(  v, w \right)  $ between two vertices $v$ and
$w$ of $G$ is defined to be the smallest length of a path from $v$ to $w$. (In
particular, $d\left(  v, v \right)  = 0$ for each $v \in V$.)

Fix a vertex $v \in V$. Define two subsets
\[
A = \left\{  w \in V \mid d\left(  v, w \right)  \text{ is even} \right\}
\qquad\text{ and } \qquad B = \left\{  w \in V \mid d\left(  v, w \right)
\text{ is odd} \right\}
\]
of $V$.

\begin{enumerate}
\item[\textbf{(a)}] Prove that $A$ is dominating.

\item[\textbf{(b)}] Prove that $B$ is dominating.

\item[\textbf{(c)}] Prove that there exists a dominating set of $G$ that has
size $\leq\left|  V \right|  /2$.

\item[\textbf{(d)}] Prove that the claim of part \textbf{(c)} holds even if we
don't assume that $G$ is connected, as long as we assume that each vertex of
$G$ has at least one neighbor. (In other words, prove Proposition
\ref{prop.dom.V/2} \textbf{(a)}.)
\end{enumerate}
\end{exercise}

\subsubsection{The number of dominating sets}

Next, we state a rather surprising recent result about the number of
dominating sets of a graph:

\begin{theorem}
[Brouwer's dominating set theorem]\label{thm.dom.brouwer}Let $G$ be a simple
graph. Then, the number of dominating sets of $G$ is odd.
\end{theorem}

Three proofs of this theorem are given in Brouwer's note \cite{Brouwe09}%
.\footnote{Other proofs can be found in the AoPS thread
\url{https://artofproblemsolving.com/community/c6h358772p1960068} . (This
thread is concerned with a superficially different contest problem, but the
latter problem is quickly revealed to be Theorem \ref{thm.dom.brouwer} in a
number-theoretical disguise.)} Let me show the one I like the most. We first
need a notation:

\begin{definition}
Let $G=\left(  V,E\right)  $ be a simple graph. A \textbf{detached pair} will
mean a pair $\left(  A,B\right)  $ of two disjoint subsets $A$ and $B$ of $V$
such that there exists no edge $ab\in E$ with $a\in A$ and $b\in B$.
\end{definition}

\begin{example}
Consider the cycle graph%
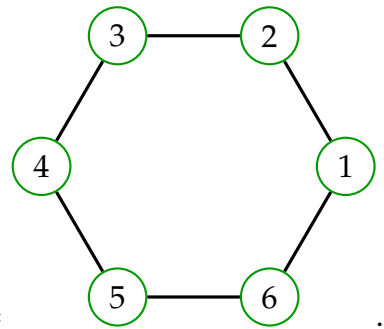
\[
C_{6}=\left(  \left\{  1,2,3,4,5,6\right\}  ,\ \left\{
12,\ 23,\ 34,\ 45,\ 56,\ 61\right\}  \right)  =%
\begin{tikzpicture}
\begin{scope}[every node/.style={circle,thick,draw=green!60!black}]
\node(A) at (0:2) {$1$};
\node(B) at (60:2) {$2$};
\node(C) at (120:2) {$3$};
\node(D) at (180:2) {$4$};
\node(E) at (240:2) {$5$};
\node(F) at (300:2) {$6$};
\end{scope}
\begin{scope}[every edge/.style={draw=black,very thick}]
\path
[-] (A) edge (B) (B) edge (C) (C) edge (D) (D) edge (E) (E) edge (F) (F) edge (A);
\end{scope}
\end{tikzpicture}%
.
\]
Then, $\left(  \left\{  1,2\right\}  ,\left\{  4,5\right\}  \right)  $ is a
detached pair, whereas $\left(  \left\{  1,2\right\}  ,\left\{  3,4\right\}
\right)  $ is not (since $23$ is an edge). Of course, there are many other
detached pairs; in particular, any pair of the form $\left(  \varnothing
,B\right)  $ or $\left(  A,\varnothing\right)  $ is detached.
\end{example}

Let me stress that the word \textquotedblleft pair\textquotedblright\ always
means \textquotedblleft ordered pair\textquotedblright\ unless I say
otherwise. So, if $\left(  A,B\right)  $ is a detached pair, then $\left(
B,A\right)  $ is a different detached pair, unless $A=B=\varnothing$.

Here is an attempt at a proof of Theorem \ref{thm.dom.brouwer}. It is a nice
example of how to apply known results to new graphs to obtain new results. The
only problem is, it shows a result that is a bit at odds with the claim of the theorem...

\begin{proof}
[Proof of Theorem \ref{thm.dom.brouwer}, attempt 1.]Write the graph $G$ as
$\left(  V,E\right)  $.

Recall that $\mathcal{P}\left(  V\right)  $ denotes the set of all subsets of
$V$.

Construct a new graph $H$ with the vertex set $\mathcal{P}\left(  V\right)  $
as follows: Two subsets $A$ and $B$ of $V$ are adjacent as vertices of $H$ if
and only if $\left(  A,B\right)  $ is a detached pair. (Note that if the
original graph $G$ has $n$ vertices, then this graph $H$ has $2^{n}$ vertices.
It is huge!)

I claim that the vertices of $H$ that have odd degree are precisely the
subsets of $V$ that are dominating. In other words:

\begin{statement}
\textit{Claim 1:} Let $A$ be a subset of $V$. Then, the vertex $A$ of $H$ has
odd degree if and only if $A$ is a dominating set of $G$.
\end{statement}

[\textit{Proof of Claim 1:} We let $N\left(  A\right)  $ denote the set of all
vertices of $G$ that have a neighbor in $A$. (This may or may not be disjoint
from $A$.)

The neighbors of $A$ (as a vertex in $H$) are precisely the subsets $B$ of $V$
such that $\left(  A,B\right)  $ is a detached pair (by the definition of
$H$). In other words, they are the subsets $B$ of $V$ that are disjoint from
$A$ and also have no neighbors in $A$ (by the definition of a
\textquotedblleft detached pair\textquotedblright). In other words, they are
the subsets $B$ of $V$ that are disjoint from $A$ and also disjoint from
$N\left(  A\right)  $. In other words, they are the subsets of the set
$V\setminus\left(  A\cup N\left(  A\right)  \right)  $. Hence, the number of
such subsets $B$ is $2^{\left\vert V\setminus\left(  A\cup N\left(  A\right)
\right)  \right\vert }$.

The degree of $A$ (as a vertex of $H$) is the number of neighbors of $A$ in
$H$. Thus, this degree is $2^{\left\vert V\setminus\left(  A\cup N\left(
A\right)  \right)  \right\vert }$ (because we have just shown that the number
of neighbors of $A$ is $2^{\left\vert V\setminus\left(  A\cup N\left(
A\right)  \right)  \right\vert }$). But $2^{k}$ is odd if and only if $k=0$.
Thus, we conclude that the degree of $A$ (as a vertex of $H$) is odd if and
only if $\left\vert V\setminus\left(  A\cup N\left(  A\right)  \right)
\right\vert =0$. The condition $\left\vert V\setminus\left(  A\cup N\left(
A\right)  \right)  \right\vert =0$ can be rewritten as follows:%
\begin{align*}
&  \ \left(  \left\vert V\setminus\left(  A\cup N\left(  A\right)  \right)
\right\vert =0\right) \\
&  \Longleftrightarrow\ \left(  V\setminus\left(  A\cup N\left(  A\right)
\right)  =\varnothing\right) \\
&  \Longleftrightarrow\ \left(  V\subseteq A\cup N\left(  A\right)  \right) \\
&  \Longleftrightarrow\ \left(  V\setminus A\subseteq N\left(  A\right)
\right) \\
&  \Longleftrightarrow\ \left(  \text{each vertex }v\in V\setminus A\text{
belongs to }N\left(  A\right)  \right) \\
&  \Longleftrightarrow\ \left(  \text{each vertex }v\in V\setminus A\text{ has
a neighbor in }A\right) \\
&  \Longleftrightarrow\ \left(  A\text{ is dominating}\right)
\ \ \ \ \ \ \ \ \ \ \left(  \text{by the definition of \textquotedblleft
dominating\textquotedblright}\right)  .
\end{align*}
Thus, what we have just shown is that the degree of $A$ (as a vertex of $H$)
is odd if and only if $A$ is dominating. This proves Claim 1.] \medskip

Claim 1 shows that the vertices of $H$ that have odd degree are precisely the
dominating sets of $G$. But the handshake lemma (Corollary
\ref{cor.deg.odd-deg-even}) tells us that any simple graph has an even number
of vertices of odd degree. Applying this to $H$, we conclude that there is an
even number of dominating sets of $G$.

Huh? We want to show that there is an \textbf{odd} number of dominating sets
of $G$, not an even number! Why did we just get the opposite result?

Puzzle: Find the mistake in our above reasoning! The answer will be revealed
on the next page.
\end{proof}

\pagebreak

So what was the mistake in our reasoning?

The mistake is that our definition of $H$ requires the vertex $\varnothing$ of
$H$ to be adjacent to itself (since $\left(  \varnothing,\varnothing\right)  $
is a detached pair); but a vertex of a simple graph cannot be adjacent to
itself. So we need to tweak the definition of $H$ somewhat:

\begin{proof}
[Correction of the above proof of Theorem \ref{thm.dom.brouwer}.]Define the
graph $H$ as above, but do not try to have $\varnothing$ adjacent to itself.
(This is the only vertex that creates any trouble, because a detached pair
$\left(  A,B\right)  $ cannot satisfy $A=B$ unless both $A$ and $B$ are
$\varnothing$.)

We WLOG assume that $V\neq\varnothing$ (otherwise, the claim is obvious).
Thus, the empty set $\varnothing$ is not dominating.

Our Claim 1 needs to be modified as follows:

\begin{statement}
\textit{Claim 1':} Let $A$ be a subset of $V$. Then, the vertex $A$ of $H$ has
odd degree if and only if $A$ is empty or a dominating set of $G$.
\end{statement}

This can be proved in the same way as we \textquotedblleft
proved\textquotedblright\ Claim 1 above; we just need to treat the
$A=\varnothing$ case separately now (but this case is easy: $\varnothing$ is
adjacent to all other vertices of $H$, and thus has degree $2^{\left\vert
V\right\vert }-1$, which is odd).

So we conclude (using the handshake lemma) that the number of empty or
dominating sets is even. Subtracting $1$ for the empty set, we conclude that
the number of dominating sets is odd (since the empty set is not dominating).
This proves Brouwer's theorem (Theorem \ref{thm.dom.brouwer}).
\end{proof}

There are other ways to prove Brouwer's theorem as well. A particularly nice
one was found by Irene Heinrich and Peter Tittmann in 2017; they gave an
\textquotedblleft explicit\textquotedblright\ formula for the number of
dominating sets that shows that this number is odd (\cite[Theorem 8]%
{HeiTit17}, restated using the language of detached pairs):

\begin{theorem}
[Heinrich--Tittmann formula]Let $G=\left(  V,E\right)  $ be a simple graph
with $n$ vertices. Assume that $n>0$.

Let $\alpha$ be the number of all detached pairs $\left(  A,B\right)  $ such
that both numbers $\left\vert A\right\vert $ and $\left\vert B\right\vert $
are even and positive.

Let $\beta$ be the number of all detached pairs $\left(  A,B\right)  $ such
that both numbers $\left\vert A\right\vert $ and $\left\vert B\right\vert $
are odd.

Then:

\begin{enumerate}
\item[\textbf{(a)}] The numbers $\alpha$ and $\beta$ are even.

\item[\textbf{(b)}] The number of dominating sets of $G$ is $2^{n}%
-1+\alpha-\beta$.
\end{enumerate}
\end{theorem}

Part \textbf{(a)} of this theorem is obvious (recall that if $\left(
A,B\right)  $ is a detached pair, then so is $\left(  B,A\right)  $). Part
\textbf{(b)} is the interesting part. In \cite[\S 3.3--\S 3.4]{17s}, I give a
long but elementary proof.

More recently (\cite{HeiTit18}), Heinrich and Tittmann have refined their
formula to allow counting dominating sets of a given size. Their main result
is the following formula (exercise 5 on homework set \#2):

\begin{exercise}
\label{exe.2.5}Let $G=\left(  V,E\right)  $ be a simple graph with at least
one vertex. Let $n=\left\vert V\right\vert $. A \textbf{detached pair} means a
pair $\left(  A,B\right)  $ of two disjoint subsets $A$ and $B$ of $V$ such
that there exists no edge $ab\in E$ with $a\in A$ and $b\in B$.

Prove the following generalization of the Heinrich--Tittmann formula:
\[
\sum_{\substack{S\text{ is a dominating}\\\text{set of }G}}x^{\left\vert
S\right\vert }=\left(  1+x\right)  ^{n}-1+\sum_{\substack{\left(  A,B\right)
\text{ is a detached pair;}\\A\neq\varnothing;\ B\neq\varnothing}}\left(
-1\right)  ^{\left\vert A\right\vert }x^{\left\vert B\right\vert }.
\]
(Here, both sides are polynomials in a single indeterminate $x$ with
coefficients in $\mathbb{Z}$.) \medskip

[\textbf{Hint:} This is a generalization of the Heinrich--Tittmann formula for
the number of dominating sets. (The latter formula can be obtained fairly
easily by substituting $x=1$ into the above and subsequently cancelling the
addends with $\left\vert A\right\vert \not \equiv \left\vert B\right\vert
\operatorname{mod}2$ against each other.) You are free to copy arguments from
\cite{17s} and change whatever needs to be changed. (Some lemmas can even be
used without any changes -- they can then be cited without proof.)]
\end{exercise}

The following exercise gives a generalization of Theorem \ref{thm.dom.brouwer}
(to recover Theorem \ref{thm.dom.brouwer} from it, set $k=1$):

\begin{exercise}
\label{exe.2017-1.6}Let $k$ be a positive integer. Let $G=\left(  V,E\right)
$ be a simple graph. A subset $U$ of $V$ will be called $k$%
\textbf{-path-dominating} if for every $v\in V$, there exists a path of length
$\leq k$ from $v$ to some element of $U$.

Prove that the number of all $k$-path-dominating subsets of $V$ is
odd.\medskip

[\textbf{Hint:} This is not as substantial a generalization as it may look.
The shortest proof is very short.]\medskip

[\textbf{Solution:} This is Exercise 6 on homework set \#1 from my Spring 2017
course; see \href{https://www.cip.ifi.lmu.de/~grinberg/t/17s/}{the course
page} for solutions.]
\end{exercise}

\subsection{\label{sec.sg.ham}Hamiltonian paths and cycles}

\subsubsection{Basics}

Now to something different. Here is a quick question: Given a simple graph
$G$, when is there a closed \textbf{walk} that contains each vertex of $G$ ?

The answer is easy: When $G$ is connected. Indeed, if a simple graph $G$ is
connected, then we can label its vertices by $v_{1},v_{2},\ldots,v_{n}$
arbitrarily, and we then get a closed walk by composing a walk from $v_{1}$ to
$v_{2}$ with a walk from $v_{2}$ to $v_{3}$ with a walk from $v_{3}$ to
$v_{4}$ and so on, ending with a walk from $v_{n}$ to $v_{1}$. This closed
walk will certainly contain each vertex. Conversely, such a walk cannot exist
if $G$ is not connected.

The question becomes a lot more interesting if we replace \textquotedblleft
closed walk\textquotedblright\ by \textquotedblleft path\textquotedblright\ or
\textquotedblleft cycle\textquotedblright. The resulting objects have a name:

\begin{definition}
Let $G=\left(  V,E\right)  $ be a simple graph.

\begin{enumerate}
\item[\textbf{(a)}] A \textbf{Hamiltonian path} in $G$ means a walk of $G$
that contains each vertex of $G$ exactly once. Obviously, it is a path.

\item[\textbf{(b)}] A \textbf{Hamiltonian cycle} in $G$ means a cycle $\left(
v_{0},v_{1},\ldots,v_{k}\right)  $ of $G$ such that each vertex of $G$ appears
exactly once among $v_{0},v_{1},\ldots,v_{k-1}$.
\end{enumerate}
\end{definition}

Some graphs have Hamiltonian paths; some don't. Having a Hamiltonian cycle is
even stronger than having a Hamiltonian path, because if $\left(  v_{0}%
,v_{1},\ldots,v_{k}\right)  $ is a Hamiltonian cycle of $G$, then $\left(
v_{0},v_{1},\ldots,v_{k-1}\right)  $ is a Hamiltonian path of $G$.

\begin{convention}
In the following, we will abbreviate:

\begin{itemize}
\item \textquotedblleft Hamiltonian path\textquotedblright\ as
\textquotedblleft\textbf{hamp}\textquotedblright;

\item \textquotedblleft Hamiltonian cycle\textquotedblright\ as
\textquotedblleft\textbf{hamc}\textquotedblright.
\end{itemize}
\end{convention}

\newpage

\begin{example}
\label{exa.hamps-hamcs} Which of the following eight graphs have hamps? Which
have hamcs?
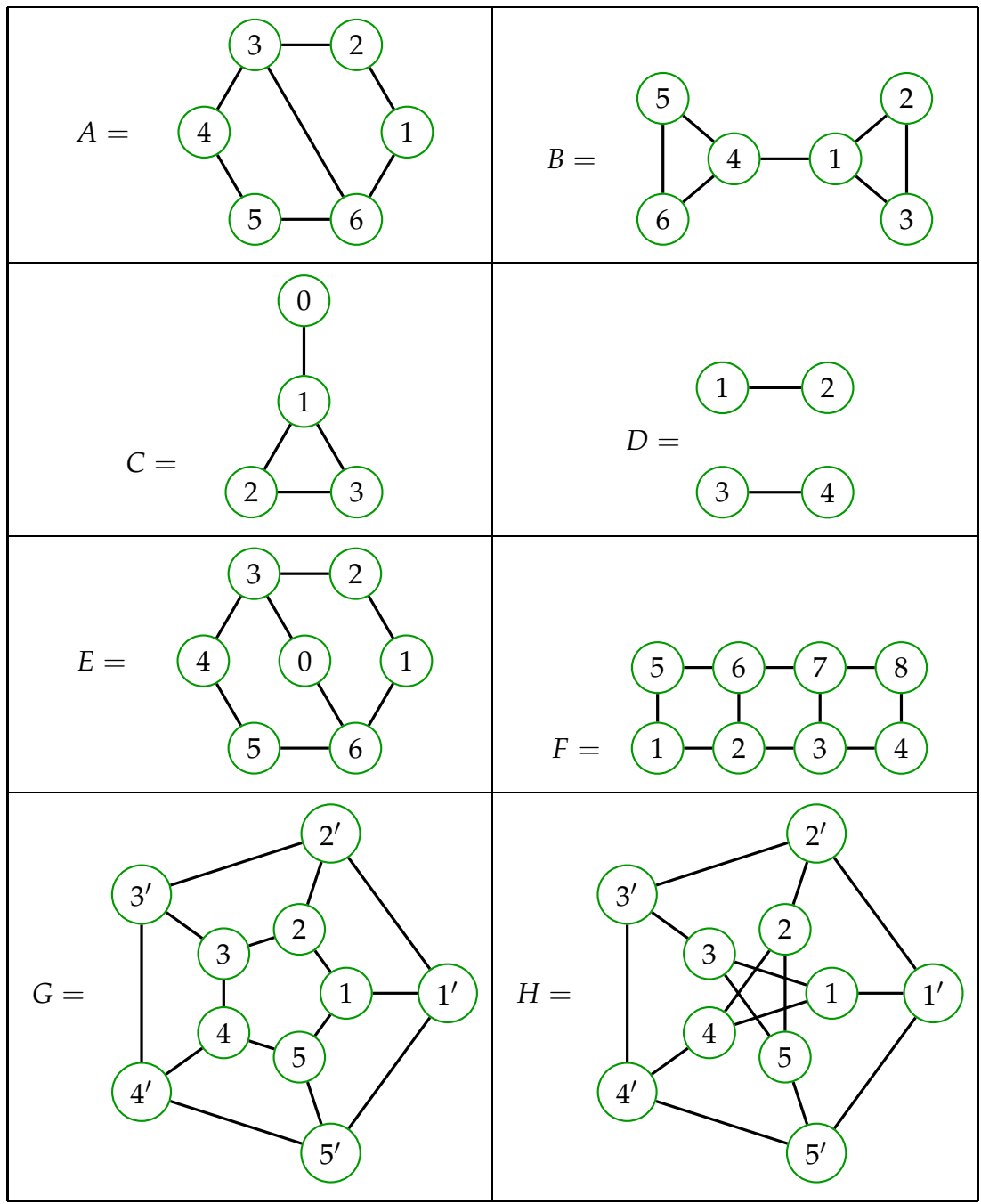
\[%

\]

\textbf{Answers:}

\begin{itemize}
\item The graph $A$ has a hamc $\left(  1,2,3,4,5,6,1\right)  $, and thus a
hamp $\left(  1,2,3,4,5,6\right)  $. (Recall that a graph that has a hamc
always has a hamp, since we can simply remove the last vertex from a hamc to
obtain a hamp.)

\item The graph $B$ has a hamp $\left(  2,3,1,4,5,6\right)  $, but no hamc.
The easiest way to see that $B$ has no hamc is the following: The edge $14$ is
a cut-edge (i.e., removing it renders the graph disconnected), thus a bridge
(i.e., an edge that appears in no cycle); therefore, any cycle must stay
entirely \textquotedblleft on one side\textquotedblright\ of this edge.

\item The graph $C$ has a hamp $\left(  0,1,2,3\right)  $, but no hamc. The
argument for the non-existence of a hamc is the same as for $B$: The edge $01$
is a bridge.

\item The graph $D$ has neither a hamp nor a hamc, because it is not
connected. Only a connected graph can have a hamp.

\item The graph $E$ has a hamp $\left(  0,3,2,1,6,5,4\right)  $, but no hamc
(checking this requires some work, though).

\item The graph $F$ has a hamc $\left(  1,2,3,4,8,7,6,5,1\right)  $, thus also
a hamp.

\item The graph $G$ has a hamc $\left(  1,2,3,4,5,5^{\prime},4^{\prime
},3^{\prime},2^{\prime},1^{\prime},1\right)  $, thus also a hamp.

\item The graph $H$ (which, by the way, is isomorphic to the Petersen graph
from Subsection \ref{subsec.sg.complete.kneser}) has a hamp $\left(
1,3,5,2,4,4^{\prime},3^{\prime},2^{\prime},1^{\prime},5^{\prime}\right)  $,
but no hamc (but this is not obvious! see
\href{https://en.wikipedia.org/wiki/Petersen_graph#Hamiltonian_paths_and_cycles}{the
Wikipedia article} for an argument).
\end{itemize}
\end{example}

In general, finding a hamp or a hamc, or proving that none exists, is a hard
problem. It can always be solved by brute force (i.e., by trying all lists of
distinct vertices and checking if there is a hamp among them, and likewise for
hamcs), but this quickly becomes forbiddingly laborious as the size of the
graph increases. Some faster algorithms exist (in particular, there is one of
running time $O\left(  n^{2}2^{n}\right)  $, where $n$ is the number of
vertices), but no polynomial-time algorithm is known. The problem (both in its
hamp version and in its hamc version) is known to be NP-hard (in the language
of complexity theory). In practice, hamps and hamcs can often be found with
some wit and perseverance; proofs of their non-existence can often be obtained
with some logic and case analysis (see the above example for some sample
arguments). See
\href{https://en.wikipedia.org/wiki/Hamiltonian path problem}{the Wikipedia
page for \textquotedblleft Hamiltonian path problem\textquotedblright\ for
more information}.

The problem of finding hamps is related to the so-called \textquotedblleft
traveling salesman problem\textquotedblright\ (TSP), which asks for a hamp
with \textquotedblleft minimum weight\textquotedblright\ in a weighted graph
(each edge has a number assigned to it, which is called its \textquotedblleft
weight\textquotedblright, and the weight of a hamp is the sum of the weights
of the edges it uses). There is a lot of computer-science literature about
this problem.

\subsubsection{Sufficient criteria: Ore and Dirac}

We shall now show some necessary criteria and some sufficient criteria (but no
necessary-and-sufficient criteria) for the existence of hamps and hamcs. Here
is the most famous sufficient criterion, known as \textbf{Ore's theorem}:

\begin{theorem}
[Ore's theorem]\label{thm.hamc.ore}Let $G=\left(  V,E\right)  $ be a simple
graph with $n$ vertices, where $n\geq3$.

Assume that $\deg x+\deg y\geq n$ for any two non-adjacent distinct vertices
$x$ and $y$.

Then, $G$ has a hamc.
\end{theorem}

There are various proofs of this theorem scattered around; see \cite[Theorem
3.6]{Harju14} or \cite[Theorem 5.3.2]{Guicha16}. We shall give a proof that is
essentially Ore's original proof \cite{Ore60} (somewhat restated as in the
\textquotedblleft Algorithm\textquotedblright\ section on
\href{https://en.wikipedia.org/wiki/Ore's theorem}{the Wikipedia page for
\textquotedblleft Ore's theorem\textquotedblright}):

\begin{proof}
[Proof of Theorem \ref{thm.hamc.ore}.]A \textbf{listing} (of $V$) shall mean a
list of elements of $V$ that contains each element exactly once. It must
clearly be an $n$-tuple.

The \textbf{hamness} of a listing $\left(  v_{1},v_{2},\ldots,v_{n}\right)  $
will mean the number of all $i\in\left\{  1,2,\ldots,n\right\}  $ such that
$v_{i}v_{i+1}\in E$. Here, we set $v_{n+1}=v_{1}$. (Visually, it is best to
represent a listing $\left(  v_{1},v_{2},\ldots,v_{n}\right)  $ by drawing the
vertices $v_{1},v_{2},\ldots,v_{n}$ on a circle in this order. Its hamness
then counts how often two successive vertices on the circle are adjacent in
the graph $G$.) Note that the hamness of a listing $\left(  v_{1},v_{2}%
,\ldots,v_{n}\right)  $ does not change if we cyclically rotate the listing
(i.e., transform it into $\left(  v_{2},v_{3},\ldots,v_{n},v_{1}\right)  $).

Clearly, if we can find a listing $\left(  v_{1},v_{2},\ldots,v_{n}\right)  $
of hamness $\geq n$, then all of $v_{1}v_{2},\ v_{2}v_{3},\ \ldots
,\ v_{n}v_{1}$ are edges of $G$, and thus $\left(  v_{1},v_{2},\ldots
,v_{n},v_{1}\right)  $ is a hamc of $G$. Thus, we need to find a listing of
hamness $\geq n$.

To do so, I will show that if you have a listing of hamness $<n$, then you can
slightly modify it to get a listing of larger hamness. In other words, I will
show the following:

\begin{statement}
\textit{Claim 1:} Let $\left(  v_{1},v_{2},\ldots,v_{n}\right)  $ be a listing
of hamness $k<n$. Then, there exists a listing of hamness larger than $k$.
\end{statement}

[\textit{Proof of Claim 1:} Since the listing $\left(  v_{1},v_{2}%
,\ldots,v_{n}\right)  $ has hamness $k<n$, there exists some $i\in\left\{
1,2,\ldots,n\right\}  $ such that $v_{i}v_{i+1}\notin E$. Pick such an $i$.
Thus, the vertices $v_{i}$ and $v_{i+1}$ of $G$ are non-adjacent (and
distinct, of course). The \textquotedblleft$\deg x+\deg y\geq n$%
\textquotedblright\ assumption of the theorem thus yields $\deg\left(
v_{i}\right)  +\deg\left(  v_{i+1}\right)  \geq n$.

However,%
\begin{align*}
\deg\left(  v_{i}\right)   &  =\left\vert \left\{  w\in V\ \mid\ v_{i}w\in
E\right\}  \right\vert \\
&  =\left\vert \left\{  j\in\left\{  1,2,\ldots,n\right\}  \ \mid\ v_{i}%
v_{j}\in E\right\}  \right\vert \\
&  =\left\vert \left\{  j\in\left\{  1,2,\ldots,n\right\}  \setminus\left\{
i\right\}  \ \mid\ v_{i}v_{j}\in E\right\}  \right\vert
\end{align*}
(because $j=i$ could not satisfy $v_{i}v_{j}\in E$ anyway) and%
\begin{align*}
\deg\left(  v_{i+1}\right)   &  =\left\vert \left\{  w\in V\ \mid\ v_{i+1}w\in
E\right\}  \right\vert \\
&  =\left\vert \left\{  j\in\left\{  1,2,\ldots,n\right\}  \ \mid
\ v_{i+1}v_{j+1}\in E\right\}  \right\vert \\
&  \ \ \ \ \ \ \ \ \ \ \ \ \ \ \ \ \ \ \ \ \left(
\begin{array}
[c]{c}%
\text{since }\left(  v_{2},v_{3},\ldots,v_{n+1}\right)  \text{ is a listing of
}V\\
\text{(because }v_{n+1}=v_{1}\text{)}%
\end{array}
\right) \\
&  =\left\vert \left\{  j\in\left\{  1,2,\ldots,n\right\}  \setminus\left\{
i\right\}  \ \mid\ v_{i+1}v_{j+1}\in E\right\}  \right\vert
\end{align*}
(because $j=i$ could not satisfy $v_{i+1}v_{j+1}\in E$ anyway). In light of
these two equalities, we can rewrite the inequality $\deg\left(  v_{i}\right)
+\deg\left(  v_{i+1}\right)  \geq n$ as%
\begin{align*}
&  \left\vert \left\{  j\in\left\{  1,2,\ldots,n\right\}  \setminus\left\{
i\right\}  \ \mid\ v_{i}v_{j}\in E\right\}  \right\vert \\
&  \ \ \ \ \ \ \ \ \ \ +\left\vert \left\{  j\in\left\{  1,2,\ldots,n\right\}
\setminus\left\{  i\right\}  \ \mid\ v_{i+1}v_{j+1}\in E\right\}  \right\vert
\geq n.
\end{align*}
Thus, the two subsets $\left\{  j\in\left\{  1,2,\ldots,n\right\}
\setminus\left\{  i\right\}  \ \mid\ v_{i}v_{j}\in E\right\}  $ and
\newline$\left\{  j\in\left\{  1,2,\ldots,n\right\}  \setminus\left\{
i\right\}  \ \mid\ v_{i+1}v_{j+1}\in E\right\}  $ of the $\left(  n-1\right)
$-element set $\left\{  1,2,\ldots,n\right\}  \setminus\left\{  i\right\}  $
have total size $\geq n$ (that is, the sum of their sizes is $\geq n$). Hence,
these two subsets must overlap (i.e., have an element in common). In other
words, there exists a $j\in\left\{  1,2,\ldots,n\right\}  \setminus\left\{
i\right\}  $ that satisfies both $v_{i}v_{j}\in E$ and $v_{i+1}v_{j+1}\in E$.
Pick such a $j$.

Now, consider a new listing obtained from the old listing $\left(  v_{1}%
,v_{2},\ldots,v_{n}\right)  $ as follows:

\begin{itemize}
\item First, cyclically rotate the old listing so that it begins with
$v_{i+1}$. Thus, you get the listing $\left(  v_{i+1},v_{i+2},\ldots
,v_{n},v_{1},v_{2},\ldots,v_{i}\right)  $.

\item Then, reverse the part of the listing starting at $v_{i+1}$ and ending
at $v_{j}$. Thus, you get the new listing%
\[
\left(  \underbrace{v_{j},v_{j-1},\ldots,v_{i+1}}_{\substack{\text{This is the
reversed part;}\\\text{it may or may not \textquotedblleft wrap
around\textquotedblright}\\\text{(i.e., contain }\ldots,v_{1},v_{n}%
,\ldots\text{ somewhere).}}},\underbrace{v_{j+1},v_{j+2},\ldots,v_{i}%
}_{\substack{\text{This is the part that}\\\text{was not reversed.}}}\right)
.
\]
This is the new listing we want.
\end{itemize}

I claim that this new listing has hamness larger than $k$. Indeed, rotating
the old listing clearly did not change its hamness. But reversing the part
from $v_{i+1}$ to $v_{j}$ clearly did: After the reversal, the edges
$v_{i}v_{i+1}$ and $v_{j}v_{j+1}$ no longer count towards the hamness (if they
were edges to begin with), but the edges $v_{i}v_{j}$ and $v_{i+1}v_{j+1}$
started counting towards the hamness. This is a good bargain, because it means
that the hamness gained $+2$ from the newly-counted edges $v_{i}v_{j}$ and
$v_{i+1}v_{j+1}$ (which, as we know, both exist), while only losing $0$ or $1$
(since the edge $v_{i}v_{i+1}$ did not exist, whereas the edge $v_{j}v_{j+1}$
may or may not have been lost). Thus, the hamness of the new listing is larger
than the hamness of the old listing either by $1$ or $2$. In other words, it
is larger than $m$ by at least $1$ or $2$.\ This proves Claim 1.] \medskip

Now, we can start with \textbf{any} listing of $V$ and keep modifying it using
Claim 1, increasing its hamness each time, until its hamness becomes $\geq n$.
But once its hamness is $\geq n$, we have found a hamc (as explained above).
Theorem \ref{thm.hamc.ore} is thus proven.
\end{proof}

As a particular case, we get a less general but more memorable criterion
called \textbf{Dirac's criterion}:

\begin{corollary}
[Dirac's theorem]\label{cor.hamc.dirac}Let $G=\left(  V,E\right)  $ be a
simple graph with $n$ vertices, where $n\geq3$.

Assume that $\deg x\geq\dfrac{n}{2}$ for each vertex $x\in V$.

Then, $G$ has a hamc.
\end{corollary}

\begin{proof}
Follows from Ore's theorem, since any two vertices $x$ and $y$ of $G$ satisfy
$\underbrace{\deg x}_{\geq\dfrac{n}{2}}+\underbrace{\deg y}_{\geq\dfrac{n}{2}%
}\geq\dfrac{n}{2}+\dfrac{n}{2}=n$.
\end{proof}

\begin{exercise}
\label{exe.4.1}\ \ 

\begin{enumerate}
\item[\textbf{(a)}] Let $G=\left(  V,E\right)  $ be a simple graph, and let
$u$ and $v$ be two distinct vertices of $G$ that are not adjacent. Let
$n=\left\vert V\right\vert $. Assume that $\deg u+\deg v\geq n$. Let
$G^{\prime}=\left(  V,E\cup\left\{  uv\right\}  \right)  $ be the simple graph
obtained from $G$ by adding a new edge $uv$. Assume that $G^{\prime}$ has a
hamc. Prove that $G$ has a hamc.

\item[\textbf{(b)}] Does this remain true if we replace \textquotedblleft
hamc\textquotedblright\ by \textquotedblleft hamp\textquotedblright?
\end{enumerate}
\end{exercise}

\subsubsection{A necessary criterion}

So much for sufficient criteria. What about necessary criteria?

\begin{proposition}
\label{prop.hamc.nec1}Let $G=\left(  V,E\right)  $ be a simple graph.

For each subset $S$ of $V$, we let $G\setminus S$ be the induced subgraph of
$G$ on the set $V\setminus S$. (In other words, this is the graph obtained
from $G$ by removing all vertices in $S$ and removing all edges that have at
least one endpoint in $S$.)

(For example, if $%
%
$ .)

Also, we let $\operatorname*{conn}\left(  H\right)  $ denote the number of
connected components of a simple graph $H$.

\begin{enumerate}
\item[\textbf{(a)}] If $G$ has a hamc, then every nonempty $S\subseteq V$
satisfies $\operatorname*{conn}\left(  G\setminus S\right)  \leq\left\vert
S\right\vert $.

\item[\textbf{(b)}] If $G$ has a hamp, then every $S\subseteq V$ satisfies
$\operatorname*{conn}\left(  G\setminus S\right)  \leq\left\vert S\right\vert
+1$.
\end{enumerate}
\end{proposition}

For example, part \textbf{(a)} of this proposition shows that the graph $E$
from Example \ref{exa.hamps-hamcs} has no hamc, because if we take $S$ to be
$\left\{  3,6\right\}  $, then $\operatorname*{conn}\left(  G\setminus
S\right)  =3$ whereas $\left\vert S\right\vert =2$. Thus, the proposition can
be used to rule out the existence of hamps and hamcs in some cases.

\begin{proof}
[Proof of Proposition \ref{prop.hamc.nec1}.]\textbf{(a)} Let $S\subseteq V$ be
a nonempty set. If we cut $\left\vert S\right\vert $ many vertices out of a
cycle, then the cycle splits into at most $\left\vert S\right\vert $ paths:%
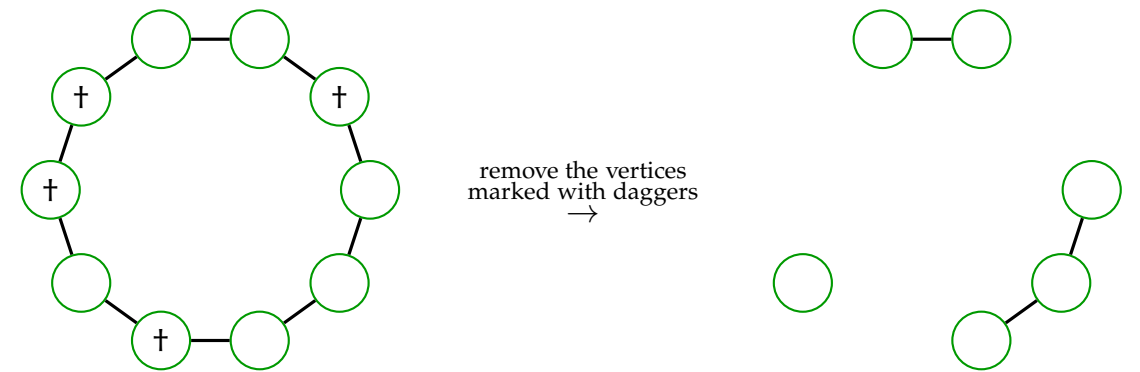
\[%
%
\]
Of course, our graph $G$ itself may not be a cycle, but if it has a hamc, then
the removal of the vertices in $S$ will split the hamc into at most
$\left\vert S\right\vert $ paths (according to the preceding sentence), and
thus the graph $G\setminus S$ will have $\leq\left\vert S\right\vert $ many
components (just using the surviving edges of the hamc alone). Taking into
account all the other edges of $G$ can only decrease the number of components.
\medskip

\textbf{(b)} This is analogous to part \textbf{(a)}.
\end{proof}

This proposition often (but not always) gives a quick way of convincing
yourself that a graph has no hamc or hamp. Alas, its converse is false. Case
in point: The Petersen graph (defined in Subsection
\ref{subsec.sg.complete.kneser}) has no hamc, but it does satisfy the
\textquotedblleft every nonempty $S\subseteq V$ satisfies
$\operatorname*{conn}\left(  G\setminus S\right)  \leq\left\vert S\right\vert
$\textquotedblright\ condition of Proposition \ref{prop.hamc.nec1}
\textbf{(a)}.

\subsubsection{\label{subsec.hamp.hypercube}Hypercubes}

Now, let us move on to a concrete example of a graph that has a hamc.

\begin{definition}
\label{def.hypercube}Let $n\in\mathbb{N}$. The $n$\textbf{-hypercube} $Q_{n}$
(more precisely, the $n$\textbf{-th hypercube graph}) is the simple graph with
vertex set%
\[
\left\{  0,1\right\}  ^{n}=\left\{  \left(  a_{1},a_{2},\ldots,a_{n}\right)
\ \mid\ \text{each }a_{i}\text{ belongs to }\left\{  0,1\right\}  \right\}
\]
and edge set defined as follows: A vertex $\left(  a_{1},a_{2},\ldots
,a_{n}\right)  \in\left\{  0,1\right\}  ^{n}$ is adjacent to a vertex $\left(
b_{1},b_{2},\ldots,b_{n}\right)  \in\left\{  0,1\right\}  ^{n}$ if and only if
there exists \textbf{exactly} one $i\in\left\{  1,2,\ldots,n\right\}  $ such
that $a_{i}\neq b_{i}$. (For example, in $Q_{4}$, the vertex $\left(
0,1,1,0\right)  $ is adjacent to $\left(  0,1,0,0\right)  $.)

The elements of $\left\{  0,1\right\}  ^{n}$ are often called
\textbf{bitstrings} (or \textbf{binary words}), and their entries are called
their \textbf{bits} (or \textbf{letters}). So two bitstrings are adjacent in
$Q_{n}$ if and only if they differ in exactly one bit.

We often write a bitstring $\left(  a_{1},a_{2},\ldots,a_{n}\right)  $ as
$a_{1}a_{2}\cdots a_{n}$. (For example, we write $\left(  0,1,1,0\right)  $ as
$0110$.)
\end{definition}

\newpage

\begin{example}
Here is how the $n$-hypercubes $Q_{n}$ look like for $n=1,2,3$:%
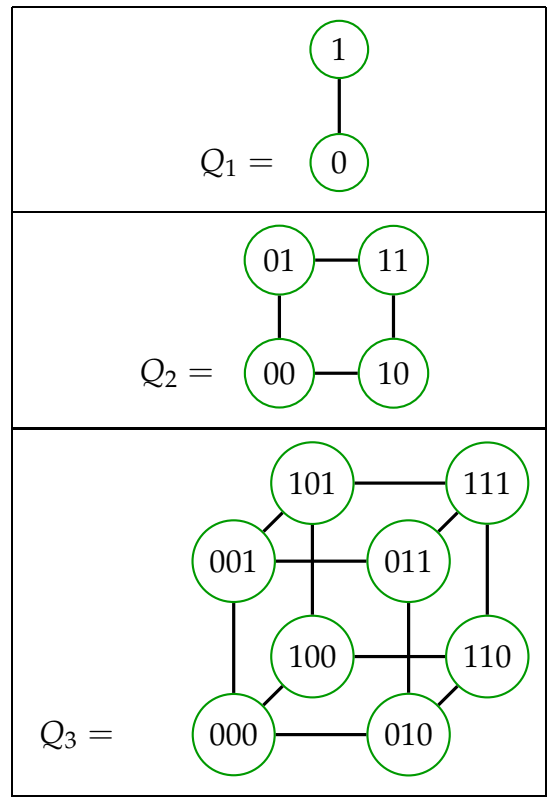
\[
\text{%

}%
\]
This should explain the name \textquotedblleft hypercube\textquotedblright.
The $0$-hypercube $Q_{0}$ is a graph with just one vertex (namely, the empty
bitstring $\left(  {}\right)  $).
\end{example}

\begin{theorem}
[Gray]\label{thm.hamc.gray}Let $n\geq2$. Then, the graph $Q_{n}$ has a hamc.
\end{theorem}

Such hamcs are known as \textbf{Gray codes}. They are circular lists of
bitstrings of length $n$ such that two consecutive bitstrings in the list
always differ in exactly one bit. See
\href{https://en.wikipedia.org/wiki/Gray code}{the Wikipedia article on
\textquotedblleft Gray codes\textquotedblright}\ for applications.

\begin{example}
Here is one hamp of $Q_{3}$ (depicted by drawing its edges while omitting all
other edges of the graph):%
\[%
\begin{tikzpicture}[scale=2.3]
\begin{scope}[every node/.style={circle,thick,draw=green!60!black}]
\node(000) at (0, 0) {$000$};
\node(001) at (0, 1) {$001$};
\node(010) at (1, 0) {$010$};
\node(011) at (1, 1) {$011$};
\node(100) at (0.45, 0.45) {$100$};
\node(101) at (0.45, 1.45) {$101$};
\node(110) at (1.45, 0.45) {$110$};
\node(111) at (1.45, 1.45) {$111$};
\end{scope}
\draw
[blue, very thick] (000) -- (001) -- (011) -- (010) -- (110) -- (111) -- (101) -- (100) -- (000);
\end{tikzpicture}%
\text{ }%
\]

\end{example}

\begin{proof}
[Proof of Theorem \ref{thm.hamc.gray}.]We will show something stronger:

\begin{statement}
\textit{Claim 1:} For each $n\geq1$, the $n$-hypercube $Q_{n}$ has a hamp from
$00\cdots0$ to $100\cdots0$.

(Keep in mind that $00\cdots0$ and $100\cdots0$ are bitstrings, not numbers:%
\[
00\cdots0=\left(  \underbrace{0,0,\ldots,0}_{n\text{ zeroes}}\right)
;\ \ \ \ \ \ \ \ \ \ 100\cdots0=\left(  1,\underbrace{0,0,\ldots,0}_{n-1\text{
zeroes}}\right)  .
\]
)
\end{statement}

[\textit{Proof of Claim 1:} We induct on $n$.

\textit{Induction base:} A look at $Q_{1}$ reveals a hamp from $0$ to $1$.

\textit{Induction step:} Fix $N\geq2$. We assume that Claim 1 holds for
$n=N-1$. In other words, $Q_{N-1}$ has a hamp from $\underbrace{00\cdots
0}_{N-1\text{ zeroes}}$ to $1\underbrace{00\cdots0}_{N-2\text{ zeroes}}$. Let
$\mathbf{p}$ be such a hamp.

By attaching a $0$ to the front of each bitstring (= vertex) in $\mathbf{p}$,
we obtain a path $\mathbf{q}$ from $\underbrace{00\cdots0}_{N\text{ zeroes}}$
to $01\underbrace{00\cdots0}_{N-2\text{ zeroes}}$ in $Q_{N}$.

By attaching a $1$ to the front of each bitstring (= vertex) in $\mathbf{p}$,
we obtain a path $\mathbf{r}$ from $1\underbrace{00\cdots0}_{N-1\text{
zeroes}}$ to $11\underbrace{00\cdots0}_{N-2\text{ zeroes}}$ in $Q_{N}$.

Now, we assemble a hamp from $\underbrace{00\cdots0}_{N\text{ zeroes}}$ to
$1\underbrace{00\cdots0}_{N-1\text{ zeroes}}$ in $Q_{N}$ as follows:

\begin{itemize}
\item Start at $\underbrace{00\cdots0}_{N\text{ zeroes}}$, and follow the path
$\mathbf{q}$ to its end (i.e., to $01\underbrace{00\cdots0}_{N-2\text{
zeroes}}$).

\item Then, move to the adjacent vertex $11\underbrace{00\cdots0}_{N-2\text{
zeroes}}$.

\item Then, follow the path $\mathbf{r}$ backwards, ending up at
$1\underbrace{00\cdots0}_{N-1\text{ zeroes}}$.
\end{itemize}

This shows that Claim 1 holds for $n=N$, too.]

Claim 1 tells us that the $n$-hypercube $Q_{n}$ has a hamp from $00\cdots0$ to
$100\cdots0$. Since its starting point $00\cdots0$ and its ending point
$100\cdots0$ are adjacent, we can turn this hamp into a hamc by appending the
starting point $00\cdots0$ again at the end. This proves Theorem
\ref{thm.hamc.gray}.
\end{proof}

\subsubsection{Cartesian products}

Theorem \ref{thm.hamc.gray} can in fact be generalized. To state the
generalization, we define the \textbf{Cartesian product} of two graphs:

\begin{definition}
\label{def.sg.cartprod}Let $G=\left(  V,E\right)  $ and $H=\left(  W,F\right)
$ be two simple graphs. The \textbf{Cartesian product} $G\times H$ of these
two graphs is defined to be the simple graph $\left(  V\times W,\ E^{\prime
}\cup F^{\prime}\right)  $, where
\begin{align*}
E^{\prime}  &  :=\left\{  \left(  v_{1},w\right)  \left(  v_{2},w\right)
\ \mid\ v_{1}v_{2}\in E\text{ and }w\in W\right\}
\ \ \ \ \ \ \ \ \ \ \text{and}\\
F^{\prime}  &  :=\left\{  \left(  v,w_{1}\right)  \left(  v,w_{2}\right)
\ \mid\ w_{1}w_{2}\in F\text{ and }v\in V\right\}  .
\end{align*}
In other words, it is the graph whose vertices are pairs $\left(  v,w\right)
\in V\times W$ consisting of a vertex of $G$ and a vertex of $H$, and whose
edges are of the forms%
\[
\left(  v_{1},w\right)  \left(  v_{2},w\right)
\ \ \ \ \ \ \ \ \ \ \text{where }v_{1}v_{2}\in E\text{ and }w\in W
\]
and%
\[
\left(  v,w_{1}\right)  \left(  v,w_{2}\right)
\ \ \ \ \ \ \ \ \ \ \text{where }w_{1}w_{2}\in F\text{ and }v\in V.
\]

\end{definition}

For example, the Cartesian product $G\times P_{2}$ of a simple graph $G$ with
the $2$-path graph $P_{2}$ can be constructed by overlaying two copies of $G$
and additionally joining each vertex of the first copy to the corresponding
vertex of the second copy by an edge. (The vertices of the first copy are the
$\left(  v,1\right)  $, whereas the vertices of the second copy are the
$\left(  v,2\right)  $.) For a specific example, here is the $5$-cycle graph
$C_{5}$ and the Cartesian product $C_{5}\times P_{2}$:%
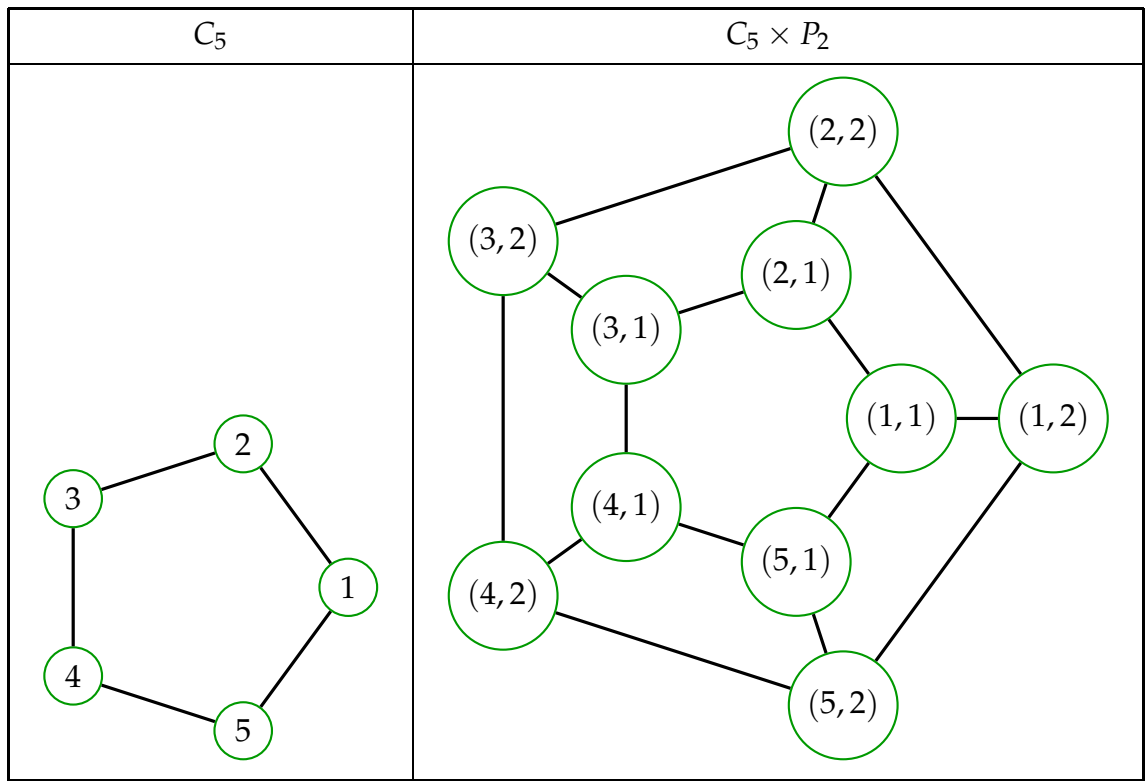
\[%
\begin{tabular}
[c]{|c|c|}\hline
$C_{5}$ & $C_{5}\times P_{2}$\\
\hline
\end{tabular}
\]
As another instance of the above description of $G\times P_{2}$, it is easy to
see the following:

\begin{proposition}
\label{prop.cartes.Qn}We have $Q_{n}\cong Q_{n-1}\times P_{2}$ for each
$n\geq1$. (See Definition \ref{def.hypercube} for the definitions of $Q_{n}$
and $Q_{n-1}$.)
\end{proposition}

\begin{proof}
This is Exercise 1 \textbf{(a)} on homework set \#2 from my Spring 2017
course; see \href{https://www.cip.ifi.lmu.de/~grinberg/t/17s/}{the course
page} for solutions.
\end{proof}

Now, we claim the following:

\begin{theorem}
\label{thm.hamc.cartes}Let $G$ and $H$ be two simple graphs. Assume that each
of the two graphs $G$ and $H$ has a hamp. Then:

\begin{enumerate}
\item[\textbf{(a)}] The Cartesian product $G\times H$ has a hamp.

\item[\textbf{(b)}] Now assume furthermore that at least one of the two
numbers $\left\vert \operatorname*{V}\left(  G\right)  \right\vert $ and
$\left\vert \operatorname*{V}\left(  H\right)  \right\vert $ is even, and that
both numbers $\left\vert \operatorname*{V}\left(  G\right)  \right\vert $ and
$\left\vert \operatorname*{V}\left(  H\right)  \right\vert $ are larger than
$1$. Then, the Cartesian product $G\times H$ has a hamc.
\end{enumerate}
\end{theorem}

\begin{proof}
This is Exercise 1 on homework set \#2 from my Spring 2017 course
(specifically, its parts \textbf{(b)} and \textbf{(c)}). Its solution can be
found on \href{https://www.cip.ifi.lmu.de/~grinberg/t/17s/}{the course page}.
(Specifically, see
\href{https://www.cip.ifi.lmu.de/~grinberg/t/17s/hw2s.pdf}{the solution to
Exercise 1 on homework set \#2 from Spring 2017}.)
\end{proof}

Now, Theorem \ref{thm.hamc.gray} can be reproved (again by inducting on $n$)
using Theorem \ref{thm.hamc.cartes} \textbf{(b)} and Proposition
\ref{prop.cartes.Qn}, since $P_{2}$ has a hamp and since $\left\vert
\operatorname*{V}\left(  P_{2}\right)  \right\vert =2$ is even. (Convince
yourself that this works!)

\subsubsection{Subset graphs}

The $n$-hypercube $Q_{n}$ can be reinterpreted in terms of subsets of
$\left\{  1,2,\ldots,n\right\}  $. Namely: Let $n\in\mathbb{N}$. Let $G_{n}$
be the simple graph whose vertex set is the powerset $\mathcal{P}\left(
\left\{  1,2,\ldots,n\right\}  \right)  $ of $\left\{  1,2,\ldots,n\right\}  $
(that is, the vertices are all $2^{n}$ subsets of $\left\{  1,2,\ldots
,n\right\}  $), and whose edges are determined as follows: Two vertices $S$
and $T$ are adjacent if and only if one of the two sets $S$ and $T$ is
obtained from the other by inserting an extra element (i.e., we have either
$S=T\cup\left\{  s\right\}  $ for some $s\notin T$, or $T=S\cup\left\{
t\right\}  $ for some $t\notin S$). Then, $G_{n}\cong Q_{n}$, since the map%
\begin{align*}
\left\{  0,1\right\}  ^{n}  &  \rightarrow\mathcal{P}\left(  \left\{
1,2,\ldots,n\right\}  \right)  ,\\
\left(  a_{1},a_{2},\ldots,a_{n}\right)   &  \mapsto\left\{  i\in\left\{
1,2,\ldots,n\right\}  \ \mid\ a_{i}=1\right\}
\end{align*}
is a graph isomorphism from $Q_{n}$ to $G_{n}$.

Thus, Theorem \ref{thm.hamc.gray} shows that for each $n\geq2$, the graph
$G_{n}$ has a hamc. In other words, for each $n\geq2$, we can list all subsets
of $\left\{  1,2,\ldots,n\right\}  $ in a circular list in such a way that
each subset on this list is obtained from the previous one by inserting or
removing a single element. For example, for $n=3$, here is such a list:%
\[
\varnothing,\ \ \left\{  1\right\}  ,\ \ \left\{  1,2\right\}  ,\ \ \left\{
2\right\}  ,\ \ \left\{  2,3\right\}  ,\ \ \left\{  1,2,3\right\}
,\ \ \left\{  1,3\right\}  ,\ \ \left\{  3\right\}  .
\]

A long-standing question only resolved a few years ago asked whether the same
can be done with the subsets of $\left\{  1,2,\ldots,n\right\}  $ having size
$\dfrac{n\pm1}{2}$ when $n$ is odd. For example, for $n=3$, we can do it as
follows:%
\[
\left\{  1\right\}  ,\ \ \left\{  1,2\right\}  ,\ \ \left\{  2\right\}
,\ \ \left\{  2,3\right\}  ,\ \ \left\{  3\right\}  ,\ \ \left\{  1,3\right\}
.
\]
In other words, if $n\geq3$ is odd, and if $G_{n}^{\prime}$ is the induced
subgraph of $G_{n}$ on the set of all subsets $J$ of $\left\{  1,2,\ldots
,n\right\}  $ that satisfy $\left\vert J\right\vert \in\left\{  \dfrac{n-1}%
{2},\dfrac{n+1}{2}\right\}  $, then does $G_{n}^{\prime}$ have a hamc?

Since $G_{n}\cong Q_{n}$, we can restate this question equivalently as
follows: If $n\geq3$ is odd, and if $Q_{n}^{\prime}$ is the induced subgraph
of $Q_{n}$ on the set
\[
\left\{  a_{1}a_{2}\cdots a_{n}\in\left\{  0,1\right\}  ^{n}\ \mid\ \sum
_{i=1}^{n}a_{i}\in\left\{  \dfrac{n-1}{2},\dfrac{n+1}{2}\right\}  \right\}  ,
\]
then does $Q_{n}^{\prime}$ have a hamc?

In 2014, Torsten M\"{u}tze proved that the answer is \textquotedblleft
yes\textquotedblright. See \cite{Mutze14} for his truly nontrivial proof,
\cite{MeMeMu22} for an analogue for Kneser graphs, and \cite{Mutze22} for a
recent survey of similar questions. (Cf. also
\href{https://en.wikipedia.org/wiki/Change_ringing}{change ringing}.)

The following exercise provides another generalization of Theorem
\ref{thm.hamc.gray}:

\begin{exercise}
\label{exe.3.4}Let $n$ and $k$ be two integers such that $n>k>0$. Define the
simple graph $Q_{n,k}$ as follows: Its vertices are the bitstrings $\left(
a_{1},a_{2},\ldots,a_{n}\right)  \in\left\{  0,1\right\}  ^{n}$; two such
bitstrings are adjacent if and only if they differ in exactly $k$ bits (in
other words: two vertices $\left(  a_{1},a_{2},\ldots,a_{n}\right)  $ and
$\left(  b_{1},b_{2},\ldots,b_{n}\right)  $ are adjacent if and only if the
number of $i\in\left\{  1,2,\ldots,n\right\}  $ satisfying $a_{i}\neq b_{i}$
equals $k$). (Thus, $Q_{n,1}$ is the $n$-hypercube graph $Q_{n}$.)

\begin{enumerate}
\item[\textbf{(a)}] Does $Q_{n,k}$ have a hamc when $k$ is even? (Recall that
\textquotedblleft hamc\textquotedblright\ is short for \textquotedblleft
Hamiltonian cycle\textquotedblright.)

\item[\textbf{(b)}] Does $Q_{n,k}$ have a hamc when $k$ is odd?
\end{enumerate}

[\textbf{Hint:} One way to approach part \textbf{(b)} is by identifying the
set $\left\{  0, 1 \right\}  $ with the field $\mathbb{F}_{2}$ with two
elements. The bitstrings $\left(  a_{1}, a_{2}, \ldots, a_{n} \right)
\in\left\{  0, 1 \right\}  ^{n}$ thus become the size-$n$ row vectors in the
$\mathbb{F}_{2}$-vector space $\mathbb{F}_{2}^{n}$. Let $e_{1}, e_{2}, \ldots,
e_{n}$ be the standard basis vectors of $\mathbb{F}_{2}^{n}$ (so that $e_{i}$
has a $1$ in its $i$-th position and zeroes everywhere else). Then, two
vectors are adjacent in the $n$-hypercube graph $Q_{n}$ (resp. in the graph
$Q_{n,k}$) if and only if their difference is one of the standard basis
vectors (resp., a sum of $k$ distinct standard basis vectors). Try to use this
to find a graph isomorphism from $Q_{n}$ to a subgraph of $Q_{n, k}$.]
\end{exercise}

The next exercise extends the idea of our proof of Theorem \ref{thm.hamc.gray}:

\begin{exercise}
\label{exe.3.8}Let $n\geq1$. Let $Q_{n}$ be the $n$-hypercube graph, as in
Definition \ref{def.hypercube}. Recall that \textquotedblleft
hamp\textquotedblright\ is short for \textquotedblleft Hamiltonian
path\textquotedblright.

At what vertices can a hamp of $Q_{n}$ end if it starts at the vertex
$00\cdots0$ ? (Find all possibilities, and prove that they are possible and
all other vertices are impossible.)
\end{exercise}

\section{\label{chp.mg}Multigraphs}

\subsection{\label{sec.mg.defs}Definitions}

So far, we have been working with simple graphs. We shall now introduce
several other kinds of graphs, starting with the \textbf{multigraphs}.

\begin{definition}
Let $V$ be a set. Then, $\mathcal{P}_{1,2}\left(  V\right)  $ shall mean the
set of all $1$-element or $2$-element subsets of $V$. In other words,%
\begin{align*}
\mathcal{P}_{1,2}\left(  V\right)  :=  &  \left\{  S\subseteq V\ \mid
\ \left\vert S\right\vert \in\left\{  1,2\right\}  \right\} \\
=  &  \left\{  \left\{  u,v\right\}  \ \mid\ u,v\in V\text{ not necessarily
distinct}\right\}  .
\end{align*}

\end{definition}

For instance,%
\[
\mathcal{P}_{1,2}\left(  \left\{  1,2,3\right\}  \right)  =\left\{  \left\{
1\right\}  ,\ \left\{  2\right\}  ,\ \left\{  3\right\}  ,\ \left\{
1,2\right\}  ,\ \left\{  1,3\right\}  ,\ \left\{  2,3\right\}  \right\}  .
\]

We can now define multigraphs:

\begin{definition}
A \textbf{multigraph} is a triple $\left(  V,E,\varphi\right)  $, where $V$
and $E$ are two finite sets, and $\varphi:E\rightarrow\mathcal{P}_{1,2}\left(
V\right)  $ is a map.
\end{definition}

\begin{example}
\label{exa.mg.exa1}Here is a multigraph:%
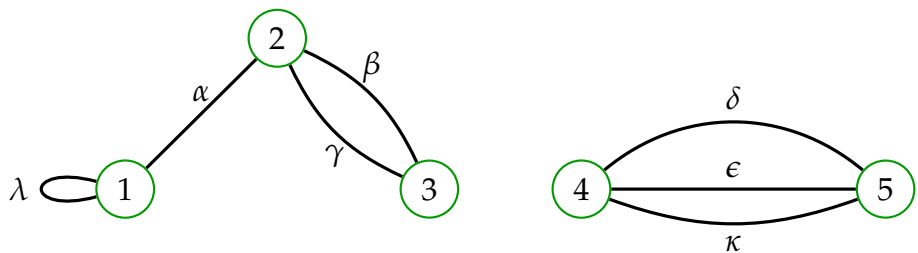
\[%
\begin{tikzpicture}[scale=2]
\begin{scope}[every node/.style={circle,thick,draw=green!60!black}]
\node(1) at (-1,0) {$1$};
\node(2) at (0,1) {$2$};
\node(3) at (1,0) {$3$};
\node(4) at (2,0) {$4$};
\node(5) at (4,0) {$5$};
\end{scope}
\begin{scope}[every edge/.style={draw=black,very thick}, every loop/.style={}]
\path[-] (1) edge[loop left] node[left] {$\lambda$} (1);
\path[-] (1) edge node[above] {$\alpha$} (2);
\path[-] (2) edge[bend left=20] node[above] {$\beta$} (3);
\path[-] (2) edge[bend right=20] node[below] {$\gamma$} (3);
\path[-] (4) edge[bend left=40] node[above] {$\delta$} (5);
\path[-] (4) edge node[above] {$\varepsilon$} (5);
\path[-] (4) edge[bend right=20] node[below] {$\kappa$} (5);
\end{scope}
\end{tikzpicture}%
\]
Formally speaking, this multigraph is the triple $\left(  V,E,\varphi\right)
$, where
\[
V=\left\{  1,2,3,4,5\right\}  ,\ \ \ \ \ \ \ \ \ \ E=\left\{  \alpha
,\beta,\gamma,\delta,\varepsilon,\kappa,\lambda\right\}  ,
\]
and where $\varphi:E\rightarrow\mathcal{P}_{1,2}\left(  V\right)  $ is the map
that sends $\alpha,\beta,\gamma,\delta,\varepsilon,\kappa,\lambda$ to
$\left\{  1,2\right\}  ,\left\{  2,3\right\}  ,\left\{  2,3\right\}  ,\left\{
4,5\right\}  ,\left\{  4,5\right\}  ,\left\{  4,5\right\}  ,\left\{
1\right\}  $, respectively. (Of course, you can write $\left\{  1\right\}  $
as $\left\{  1,1\right\}  $.)
\end{example}

This suggests the following terminology (most of which is a calque of our
previously defined terminology for simple graphs):

\begin{definition}
\label{def.mg.walks-etc}Let $G=\left(  V,E,\varphi\right)  $ be a multigraph. Then:

\begin{enumerate}
\item[\textbf{(a)}] The elements of $V$ are called the \textbf{vertices} of
$G$.

The set $V$ is called the \textbf{vertex set} of $G$, and is denoted
$\operatorname*{V}\left(  G\right)  $.

\item[\textbf{(b)}] The elements of $E$ are called the \textbf{edges} of $G$.

The set $E$ is called the \textbf{edge set} of $G$, and is denoted
$\operatorname*{E}\left(  G\right)  $.

\item[\textbf{(c)}] If $e$ is an edge of $G$, then the elements of
$\varphi\left(  e\right)  $ are called the \textbf{endpoints} of $e$.

\item[\textbf{(d)}] We say that an edge $e$ \textbf{contains} a vertex $v$ if
$v\in\varphi\left(  e\right)  $ (in other words, if $v$ is an endpoint of $e$).

\item[\textbf{(e)}] Two vertices $u$ and $v$ are said to be \textbf{adjacent}
if there exists an edge $e \in E$ whose endpoints are $u$ and $v$.

\item[\textbf{(f)}] Two edges $e$ and $f$ are said to be \textbf{parallel} if
$\varphi\left(  e\right)  =\varphi\left(  f\right)  $. (In the above example,
any two of the edges $\delta,\varepsilon,\kappa$ are parallel.)

\item[\textbf{(g)}] We say that $G$ has \textbf{no parallel edges} if $G$ has
no two distinct edges that are parallel.

\item[\textbf{(h)}] An edge $e$ is called a \textbf{loop} (or
\textbf{self-loop}) if $\varphi\left(  e\right)  $ is a $1$-element set (i.e.,
if $e$ has only one endpoint). (In Example \ref{exa.mg.exa1}, the edge
$\lambda$ is a loop.)

\item[\textbf{(i)}] We say that $G$ is \textbf{loopless} if $G$ has no loops
(among its edges).

\item[\textbf{(j)}] The \textbf{degree} $\deg v$ (also written $\deg_{G}v$) of
a vertex $v$ of $G$ is defined to be the number of edges that contain $v$,
where loops are counted twice. In other words,%
\begin{align*}
\deg v=  &  \deg_{G}v\\
:=  &  \underbrace{\left\vert \left\{  e\in E\ \mid\ v\in\varphi\left(
e\right)  \right\}  \right\vert }_{\substack{\text{this counts all
edges}\\\text{that contain }v}}+\underbrace{\left\vert \left\{  e\in
E\ \mid\ \varphi\left(  e\right)  =\left\{  v\right\}  \right\}  \right\vert
}_{\substack{\text{this counts all loops}\\\text{that contain }v\text{ once
again}}}.
\end{align*}
(Note that, unlike in the case of a simple graph, $\deg v$ is \textbf{not} the
number of neighbors of $v$, unless it happens that $v$ is not contained in any
loops or parallel edges.)

(For example, in Example \ref{exa.mg.exa1}, we have $\deg1 = 3$ and $\deg2 =
3$ and $\deg3 = 2$ and $\deg4 = 3$ and $\deg5 = 3$.)

\item[\textbf{(k)}] A \textbf{walk} in $G$ means a list of the form%
\[
\left(  v_{0},e_{1},v_{1},e_{2},v_{2},\ldots,e_{k},v_{k}\right)
\ \ \ \ \ \ \ \ \ \ \left(  \text{with }k\geq0\right)  ,
\]
where $v_{0},v_{1},\ldots,v_{k}$ are vertices of $G$, where $e_{1}%
,e_{2},\ldots,e_{k}$ are edges of $G$, and where each $i\in\left\{
1,2,\ldots,k\right\}  $ satisfies%
\[
\varphi\left(  e_{i}\right)  =\left\{  v_{i-1},v_{i}\right\}
\]
(that is, the endpoints of each edge $e_{i}$ are $v_{i-1}$ and $v_{i}$). Note
that we have to record both the vertices \textbf{and} the edges in our walk,
since we want the walk to \textquotedblleft know\textquotedblright\ which
edges it traverses. (For instance, in Example \ref{exa.mg.exa1}, the two walks
$\left(  1,\alpha,2,\beta,3\right)  $ and $\left(  1,\alpha,2,\gamma,3\right)
$ are distinct.)

The \textbf{vertices} of a walk $\left(  v_{0},e_{1},v_{1},e_{2},v_{2}%
,\ldots,e_{k},v_{k}\right)  $ are $v_{0},v_{1},\ldots,v_{k}$; the
\textbf{edges} of this walk are $e_{1},e_{2},\ldots,e_{k}$. This walk is said
to \textbf{start} at $v_{0}$ and \textbf{end} at $v_{k}$; it is also said to
be a \textbf{walk from }$v_{0}$ \textbf{to }$v_{k}$. Its \textbf{starting
point} is $v_{0}$, and its \textbf{ending point} is $v_{k}$. Its
\textbf{length} is $k$.

\item[\textbf{(l)}] A \textbf{path} means a walk whose vertices are distinct.

\item[\textbf{(m)}] The notions of \textquotedblleft\textbf{path-connected}%
\textquotedblright\ and \textquotedblleft\textbf{connected}\textquotedblright%
\ and \textquotedblleft\textbf{component}\textquotedblright\ are defined
exactly as for simple graphs. The symbol $\simeq_{G}$ still means
\textquotedblleft path-connected\textquotedblright.

\item[\textbf{(n)}] A \textbf{closed walk} (or \textbf{circuit}) means a walk
$\left(  v_{0},e_{1},v_{1},e_{2},v_{2},\ldots,e_{k},v_{k}\right)  $ with
$v_{k}=v_{0}$.

\item[\textbf{(o)}] A \textbf{cycle} means a closed walk $\left(  v_{0}%
,e_{1},v_{1},e_{2},v_{2},\ldots,e_{k},v_{k}\right)  $ such that

\begin{itemize}
\item the vertices $v_{0},v_{1},\ldots,v_{k-1}$ are distinct;

\item the edges $e_{1},e_{2},\ldots,e_{k}$ are distinct;

\item we have $k\geq1$.
\end{itemize}

(Note that we are not requiring $k\geq3$ any more, as we did for simple
graphs. Thus, in Example \ref{exa.mg.exa1}, both $\left(  2,\beta
,3,\gamma,2\right)  $ and $\left(  1,\lambda,1\right)  $ are cycles, but
$\left(  2,\beta,3,\beta,2\right)  $ is not. The purpose of the
\textquotedblleft$k\geq3$\textquotedblright\ requirement for cycles in simple
graphs was to disallow closed walks such as $\left(  2,\beta,3,\beta,2\right)
$ from being cycles; but they are now excluded by the \textquotedblleft the
edges $e_{1},e_{2},\ldots,e_{k}$ are distinct\textquotedblright\ condition.)

\item[\textbf{(p)}] Hamiltonian paths and cycles are defined as for simple graphs.

\item[\textbf{(q)}] We draw a multigraph by drawing each vertex as a point,
each edge as a curve, and labeling both the vertices and the edges (or not, if
we don't care about what they are). An example of such a drawing appeared in
Example \ref{exa.mg.exa1}.
\end{enumerate}
\end{definition}

So there are two differences between simple graphs and multigraphs:

\begin{enumerate}
\item A multigraph can have loops, whereas a simple graph cannot.

\item In a simple graph, an edge $e$ \textbf{is} a set of two vertices,
whereas in a multigraph, an edge $e$ \textbf{has} a set of two vertices
(possibly two equal ones, if $e$ is a loop) assigned to it by the map
$\varphi$. This not only allows for parallel edges, but also lets us store
some information in the \textquotedblleft identities\textquotedblright\ of the edges.
\end{enumerate}

Nevertheless, the two notions have much in common; thus, they are both called
\textquotedblleft graphs\textquotedblright:

\begin{convention}
The word \textquotedblleft\textbf{graph}\textquotedblright\ means either
\textquotedblleft simple graph\textquotedblright\ or \textquotedblleft
multigraph\textquotedblright. The precise meaning should usually be understood
from the context. (I will try not to use it when it could cause confusion.)
\end{convention}

Fortunately, simple graphs and multigraphs have many properties in common, and
often it is not hard to derive a result about multigraphs from the analogous
result about simple graphs or vice versa. We will soon explore how some of the
properties we have seen in the previous chapter can be adapted to multigraphs.
First, however, let us explain how to convert multigraphs into simple graphs
and vice versa.

\subsection{\label{sec.mg.conv}Conversions}

We can turn each multigraph into a simple graph, but at a cost of losing some information:

\begin{definition}
Let $G=\left(  V,E,\varphi\right)  $ be a multigraph. Then, the
\textbf{underlying simple graph} $G^{\operatorname*{simp}}$ of $G$ means the
simple graph%
\[
\left(  V,\ \left\{  \varphi\left(  e\right)  \ \mid\ e\in E\text{ is not a
loop}\right\}  \right)  .
\]
In other words, it is the simple graph with vertex set $V$ in which two
distinct vertices $u$ and $v$ are adjacent if and only if $u$ and $v$ are
adjacent in $G$. Thus, $G^{\operatorname*{simp}}$ is obtained from $G$ by
removing loops and \textquotedblleft collapsing\textquotedblright\ parallel
edges to a single edge.
\end{definition}

For example, the underlying simple graph of the multigraph $G$ in Example
\ref{exa.mg.exa1} is%
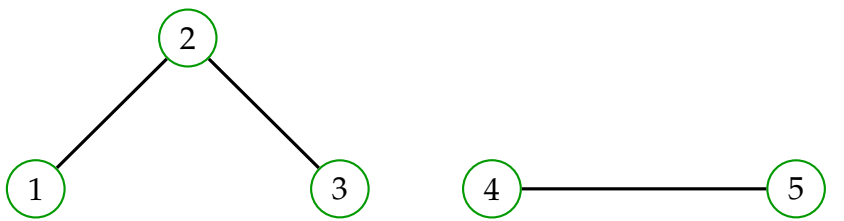
\[%
%
\ \ .
\]

Conversely, each simple graph can be viewed as a multigraph:

\begin{definition}
Let $G=\left(  V,E\right)  $ be a simple graph. Then, the
\textbf{corresponding multigraph} $G^{\operatorname*{mult}}$ is defined to be
the multigraph%
\[
\left(  V,E,\iota\right)  ,
\]
where $\iota:E\rightarrow\mathcal{P}_{1,2}\left(  V\right)  $ is the map
sending each $e\in E$ to $e$ itself.
\end{definition}

\begin{example}
If%
\[%
%
\ \ .
\]

\end{example}

As we said, the \textquotedblleft underlying simple graph\textquotedblright%
\ construction $G\mapsto G^{\operatorname*{simp}}$ destroys information, so it
is irreversible. This being said, the two constructions $G\mapsto
G^{\operatorname*{simp}}$ and $G\mapsto G^{\operatorname*{mult}}$ come fairly
close to undoing one another:\footnote{In the following proposition, we will
use the notion of an \textquotedblleft isomorphism of
multigraphs\textquotedblright. A rigorous definition of this notion is given
in Definition \ref{def.mg.iso} further below (but it is more or less what you
would expect: it is a way to relabel the vertices and the edges of one
multigraph to obtain those of another).}

\begin{proposition}
\ 

\begin{enumerate}
\item[\textbf{(a)}] If $G$ is a simple graph, then $\left(
G^{\operatorname*{mult}}\right)  ^{\operatorname*{simp}}=G$.

\item[\textbf{(b)}] If $G$ is a loopless multigraph that has no parallel
edges, then $\left(  G^{\operatorname*{simp}}\right)  ^{\operatorname*{mult}%
}\cong G$. (This is just an isomorphism, not an equality, since the
\textquotedblleft identities\textquotedblright\ of the edges of $G$ have been
forgotten in $G^{\operatorname*{simp}}$ and cannot be recovered.)

\item[\textbf{(c)}] If $G$ is a multigraph that has loops or (distinct)
parallel edges, then the multigraph $\left(  G^{\operatorname*{simp}}\right)
^{\operatorname*{mult}}$ has fewer edges than $G$ and thus is not isomorphic
to $G$.
\end{enumerate}
\end{proposition}

\begin{proof}
A matter of understanding the definitions.
\end{proof}

We will often identify a simple graph $G$ with the corresponding multigraph
$G^{\operatorname*{mult}}$. This may be dangerous, because we have defined
notions such as adjacency, walks, paths, cycles, etc. both for simple graphs
and for multigraphs; thus, when we identify a simple graph $G$ with the
multigraph $G^{\operatorname*{mult}}$, we are potentially inviting ambiguity
(for example, does \textquotedblleft cycle of $G$\textquotedblright\ mean a
cycle of the simple graph $G$ or of the multigraph $G^{\operatorname*{mult}}$
?). Fortunately, this ambiguity is harmless, because whenever $G$ is a simple
graph, any of the notions we defined for $G$ is equivalent to the
corresponding notion for the multigraph $G^{\operatorname*{mult}}$. For
example, for the notions of a cycle, we have the following:

\begin{proposition}
Let $G$ be a simple graph. Then:

\begin{enumerate}
\item[\textbf{(a)}] If $\left(  v_{0},e_{1},v_{1},e_{2},v_{2},\ldots
,e_{k},v_{k}\right)  $ is a cycle of the multigraph $G^{\operatorname*{mult}}%
$, then $\left(  v_{0},v_{1},\ldots,v_{k}\right)  $ is a cycle of the simple
graph $G$.

\item[\textbf{(b)}] Conversely, if $\left(  v_{0},v_{1},\ldots,v_{k}\right)  $
is a cycle of the simple graph $G$, then $\left(  v_{0},\left\{  v_{0}%
,v_{1}\right\}  ,v_{1},\left\{  v_{1},v_{2}\right\}  ,v_{2},\ldots
,v_{k-1},\left\{  v_{k-1},v_{k}\right\}  ,v_{k}\right)  $ is a cycle of the
multigraph $G^{\operatorname*{mult}}$.
\end{enumerate}
\end{proposition}

\begin{proof}
This is not completely obvious, since our definitions of a cycle of a simple
graph and of a cycle of a multigraph were somewhat different. The proof boils
down to checking the following two statements:

\begin{enumerate}
\item If $\left(  v_{0},v_{1},\ldots,v_{k}\right)  $ is a cycle of the simple
graph $G$, then its edges \newline$\left\{  v_{0},v_{1}\right\}  ,\left\{
v_{1},v_{2}\right\}  ,\ldots,\left\{  v_{k-1},v_{k}\right\}  $ are distinct.

\item If $\left(  v_{0},e_{1},v_{1},e_{2},v_{2},\ldots,e_{k},v_{k}\right)  $
is a cycle of the multigraph $G^{\operatorname*{mult}}$, then $k\geq3$.
\end{enumerate}

Checking statement 2 is easy (we cannot have $k=1$ since
$G^{\operatorname*{mult}}$ has no loops, and we cannot have $k=2$ since this
would lead to $e_{1}=e_{2}$). Statement 1 is also clear, since the
distinctness of the $k$ vertices $v_{0},v_{1},\ldots,v_{k-1}$ forces the
$2$-element sets formed from these $k$ vertices to also be distinct (and since
the edges $\left\{  v_{0},v_{1}\right\}  ,\left\{  v_{1},v_{2}\right\}
,\ldots,\left\{  v_{k-1},v_{k}\right\}  =\left\{  v_{k-1},v_{0}\right\}  $ are
such $2$-element sets).
\end{proof}

For all other notions discussed above, it is even more obvious that there is
no ambiguity.

\subsection{\label{sec.mg.gen}Generalizing from simple graphs to multigraphs}

Now, as promised, we shall revisit the results of Chapter \ref{chp.sg}, and
see which of them also hold for multigraphs instead of simple graphs.

\subsubsection{The Ramsey number $R\left(  3,3\right)  $}

One of the first properties of simple graphs that we proved is the following
(Proposition \ref{prop.simple.R33}):

\begin{proposition}
Let $G$ be a simple graph with $\left|  \operatorname{V}\left(  G \right)
\right|  \geq6$ (that is, $G$ has at least $6$ vertices). Then, at least one
of the following two statements holds:

\begin{itemize}
\item \textit{Statement 1:} There exist three distinct vertices $a$, $b$ and
$c$ of $G$ such that $ab$, $bc$ and $ca$ are edges of $G$.

\item \textit{Statement 2:} There exist three distinct vertices $a$, $b$ and
$c$ of $G$ such that none of $ab$, $bc$ and $ca$ is an edge of $G$.
\end{itemize}
\end{proposition}

This is still true for multigraphs\footnote{Of course, we should understand it
appropriately: i.e., we should read \textquotedblleft$ab$ is an
edge\textquotedblright\ as \textquotedblleft there is an edge with endpoints
$a$ and $b$\textquotedblright.}, because replacing a multigraph $G$ by the
underlying simple graph $G^{\operatorname*{simp}}$ does not change the meaning
of the statement.

\subsubsection{Degrees}

In Definition \ref{def.sg.deg}, we defined the degree of a vertex $v$ in a
simple graph $G=\left(  V,E\right)  $ by%
\begin{align*}
\deg v:=  &  \left(  \text{the number of edges }e\in E\text{ that contain
}v\right) \\
=  &  \left(  \text{the number of neighbors of }v\right) \\
=  &  \left\vert \left\{  u\in V\ \mid\ uv\in E\right\}  \right\vert \\
=  &  \left\vert \left\{  e\in E\ \mid\ v\in e\right\}  \right\vert .
\end{align*}
These equalities \textbf{no longer hold} when $G$ is a multigraph. Parallel
edges correspond to the same neighbor, so the number of neighbors of $v$ is
only a lower bound on $\deg v$. (For instance, in Example \ref{exa.mg.exa1},
the vertex $2$ has degree $3$ but only $2$ neighbors.)

\bigskip

Proposition \ref{prop.degv.lessn} (which says that if $G$ is a simple graph
with $n$ vertices, then any vertex $v$ of $G$ has degree $\deg v\in\left\{
0,1,\ldots,n-1\right\}  $) also \textbf{no longer holds} for multigraphs,
because you can have arbitrarily many edges in a multigraph with just $1$ or
$2$ vertices. (You can even have parallel loops!)

\bigskip

Is Proposition \ref{prop.deg.euler} true for multigraphs? Yes, because we have
said that loops should count twice in the definition of the degree. The proof
needs some tweaking, though. Let me give a slightly different proof; but
first, let me state the claim for multigraphs as a proposition of its own:

\begin{proposition}
[Euler 1736 for multigraphs]\label{prop.mg.sum-deg}Let $G$ be a multigraph.
Then, the sum of the degrees of all vertices of $G$ equals twice the number of
edges of $G$. In other words,%
\[
\sum_{v\in\operatorname*{V}\left(  G\right)  }\deg v=2\cdot\left\vert
\operatorname*{E}\left(  G\right)  \right\vert .
\]

\end{proposition}

\begin{proof}
Write $G$ as $G=\left(  V,E,\varphi\right)  $; thus, $\operatorname*{V}\left(
G\right)  =V$ and $\operatorname*{E}\left(  G\right)  =E$.

For each edge $e$, let us (arbitrarily) choose one endpoint of $e$ and denote
it by $\alpha\left(  e\right)  $. The other endpoint will be called
$\beta\left(  e\right)  $. If $e$ is a loop, then we set $\beta\left(
e\right)  =\alpha\left(  e\right)  $. Then, for each vertex $v$, we have%
\begin{align*}
\deg v  &  =\left(  \text{the number of }e\in E\text{ such that }%
v=\alpha\left(  e\right)  \right) \\
&  \ \ \ \ \ \ \ \ \ \ +\left(  \text{the number of }e\in E\text{ such that
}v=\beta\left(  e\right)  \right)
\end{align*}
(note how loops get counted twice on the right hand side, because if $e\in E$
is a loop, then $v$ is both $\alpha\left(  e\right)  $ and $\beta\left(
e\right)  $ at the same time). Summing up this equality over all $v\in V$, we
obtain%
\begin{align*}
\sum_{v\in V}\deg v  &  =\sum_{v\in V}\left(  \text{the number of }e\in
E\text{ such that }v=\alpha\left(  e\right)  \right) \\
&  \ \ \ \ \ \ \ \ \ \ +\sum_{v\in V}\left(  \text{the number of }e\in E\text{
such that }v=\beta\left(  e\right)  \right)  .
\end{align*}
However,
\[
\sum_{v\in V}\left(  \text{the number of }e\in E\text{ such that }%
v=\alpha\left(  e\right)  \right)  =\left\vert E\right\vert ,
\]
since each edge $e\in E$ is counted in exactly one addend of this sum.
Similarly,
\[
\sum_{v\in V}\left(  \text{the number of }e\in E\text{ such that }%
v=\beta\left(  e\right)  \right)  =\left\vert E\right\vert .
\]
Thus, the above equality becomes%
\begin{align*}
\sum_{v\in V}\deg v  &  =\underbrace{\sum_{v\in V}\left(  \text{the number of
}e\in E\text{ such that }v=\alpha\left(  e\right)  \right)  }_{=\left\vert
E\right\vert }\\
&  \ \ \ \ \ \ \ \ \ \ +\underbrace{\sum_{v\in V}\left(  \text{the number of
}e\in E\text{ such that }v=\beta\left(  e\right)  \right)  }_{=\left\vert
E\right\vert }\\
&  =\left\vert E\right\vert +\left\vert E\right\vert =2\cdot\left\vert
E\right\vert .
\end{align*}
This proves Proposition \ref{prop.mg.sum-deg}.
\end{proof}

This is a good motivation for counting loops twice in the definition of a degree.

\bigskip

The handshake lemma (Corollary \ref{cor.deg.odd-deg-even}) \textbf{still holds
for multigraphs}. In other words, we have the following:

\begin{corollary}
[handshake lemma]\label{cor.mg.odd-deg-even}Let $G$ be a multigraph. Then, the
number of vertices $v$ of $G$ whose degree $\deg v$ is odd is even.
\end{corollary}

\begin{proof}
This follows from Proposition \ref{prop.mg.sum-deg} in the same way as for
simple graphs.
\end{proof}

\bigskip

Proposition \ref{prop.sg.deg.equal-degs} \textbf{fails for multigraphs}. For
example, the multigraph \newline$%
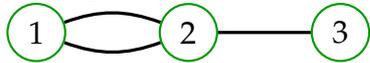
\begin{tikzpicture}[scale=2]
\begin{scope}[every node/.style={circle,thick,draw=green!60!black}]
\node(1) at (-1,0) {$1$};
\node(2) at (0,0) {$2$};
\node(3) at (1,0) {$3$};
\end{scope}
\begin{scope}[every edge/.style={draw=black,very thick}, every loop/.style={}]
\path[-] (1) edge[bend left=20] (2);
\path[-] (1) edge[bend right=20] (2);
\path[-] (2) edge (3);
\end{scope}
\end{tikzpicture}%
$ has three vertices with degrees $1,2,3$. Fortunately, Proposition
\ref{prop.sg.deg.equal-degs} was more of a curiosity than a useful fact.

\bigskip

Mantel's theorem (Theorem \ref{thm.sg.mantel}) also \textbf{fails for
multigraphs}, because we can join two vertices with a lot of parallel edges
and thus satisfy $e>n^{2}/4$ for stupid reasons without ever creating a
triangle. Thus, Turan's theorem (Theorem \ref{thm.sg.turan}) also fails for multigraphs.

\subsubsection{Graph isomorphisms}

Graph isomorphy (and isomorphisms) can still be defined for multigraphs, but
the definition is not the same as for simple graphs. Graph isomorphisms can no
longer be defined merely as bijections between the vertex sets, since we also
need to specify what they do to the edges. Instead, we define them as follows:

\begin{definition}
\label{def.mg.iso}Let $G=\left(  V,E,\varphi\right)  $ and $H=\left(
W,F,\psi\right)  $ be two multigraphs.

\begin{enumerate}
\item[\textbf{(a)}] A \textbf{graph isomorphism} (or \textbf{isomorphism})
from $G$ to $H$ means a \textbf{pair} $\left(  \alpha,\beta\right)  $ of
bijections%
\[
\alpha:V\rightarrow W\ \ \ \ \ \ \ \ \ \ \text{and}\ \ \ \ \ \ \ \ \ \ \beta
:E\rightarrow F
\]
with the property that if $e\in E$, then the endpoints of $\beta\left(
e\right)  $ are the images under $\alpha$ of the endpoints of $e$. (This
property can also be restated as a commutative diagram
\[%
\xymatrix@C=4pc{
E \ar[r]^\beta\ar[d]_{\varphi} & F \ar[d]^\psi\\
\mathcal{P}_{1,2}\left(V\right) \ar[r]_{\mathcal{P}\left(\alpha\right)}
& \mathcal{P}_{1,2}\left(W\right)
}%
\ \ ,
\]
where $\mathcal{P}\left(  \alpha\right)  $ is the map from $\mathcal{P}%
_{1,2}\left(  V\right)  $ to $\mathcal{P}_{1,2}\left(  W\right)  $ that sends
each subset $\left\{  u,v\right\}  \in\mathcal{P}_{1,2}\left(  V\right)  $ to
$\left\{  \alpha\left(  u\right)  ,\alpha\left(  v\right)  \right\}
\in\mathcal{P}_{1,2}\left(  W\right)  $. If you are used to category theory,
this restatement may look more natural to you.)

\item[\textbf{(b)}] We say that $G$ and $H$ are \textbf{isomorphic} (this is
written $G\cong H$) if there exists a graph isomorphism from $G$ to $H$.
\end{enumerate}
\end{definition}

Again, isomorphy of multigraphs is an equivalence relation.

\subsubsection{Complete graphs, paths, cycles}

In Definition \ref{def.sg.complete-empty}, Definition \ref{def.sg.path} and
Definition \ref{def.sg.cycle}, we defined the complete graphs $K_{n}$, the
path graphs $P_{n}$ and the cycle graphs $C_{n}$ as simple graphs. Thus, all
of them can be viewed as multigraphs if one so desires (since each simple
graph $G$ gives rise to a multigraph $G^{\operatorname*{mult}}$).

However, using multigraphs, we can extend our definition of $n$-th cycle
graphs $C_{n}$ to the case $n=1$ and also tweak it in the case $n=2$ to make
it more natural. We do this as follows:

\begin{definition}
\label{def.mg.Cn} We modify the definition of cycle graphs (Definition
\ref{def.sg.cycle}) as follows:

\begin{enumerate}
\item[\textbf{(a)}] We \textbf{redefine} the $2$-nd cycle graph $C_{2}$ to be
the multigraph with two vertices $1$ and $2$ and two parallel edges with
endpoints $1$ and $2$. (We don't care what the edges are, only that there are
two of them and each has endpoints $1$ and $2$.) Thus, it looks as follows: $%
\begin{tikzpicture}[scale=2]
\begin{scope}[every node/.style={circle,thick,draw=green!60!black}]
\node(1) at (-1,0) {$1$};
\node(2) at (0,0) {$2$};
\end{scope}
\begin{scope}[every edge/.style={draw=black,very thick}, every loop/.style={}]
\path[-] (1) edge[bend left=20] (2);
\path[-] (1) edge[bend right=20] (2);
\end{scope}
\end{tikzpicture}%
$\ \ .

\item[\textbf{(b)}] We define the $1$-st cycle graph $C_{1}$ to be the
multigraph with one vertex $1$ and one edge (which is necessarily a loop).
Thus, it looks as follows: $%
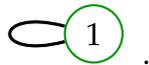
\begin{tikzpicture}[scale=2]
\begin{scope}[every node/.style={circle,thick,draw=green!60!black}]
\node(1) at (-1,0) {$1$};
\end{scope}
\begin{scope}[every edge/.style={draw=black,very thick}, every loop/.style={}]
\path[-] (1) edge[loop left] (1);
\end{scope}
\end{tikzpicture}%
$\ \ .
\end{enumerate}
\end{definition}

This has the effect that the $n$-th cycle graph $C_{n}$ has exactly $n$ edges
for each $n\geq1$ (rather than having $1$ edge for $n=2$, as it did back when
it was a simple graph).

\subsubsection{Induced submultigraphs}

In Definition \ref{def.sg.subgraph}, we defined subgraphs and induced
subgraphs of a simple graph. The corresponding notions for multigraphs are
defined as follows:

\begin{definition}
Let $G=\left(  V,E,\varphi\right)  $ be a multigraph.

\begin{enumerate}
\item[\textbf{(a)}] A \textbf{submultigraph} of $G$ means a multigraph of the
form $H=\left(  W,F,\psi\right)  $, where $W\subseteq V$ and $F\subseteq E$
and $\psi=\varphi\mid_{F}$. In other words, a submultigraph of $G$ means a
multigraph $H$ whose vertices are vertices of $G$ and whose edges are edges of
$G$ and whose edges have the same endpoints in $H$ as they do in $G$.

We often abbreviate \textquotedblleft submultigraph\textquotedblright\ as
\textquotedblleft\textbf{subgraph}\textquotedblright.

\item[\textbf{(b)}] Let $S$ be a subset of $V$. The \textbf{induced
submultigraph of }$G$\textbf{ on the set }$S$ denotes the submultigraph%
\[
\left(  S,\ \ E^{\prime},\ \ \varphi\mid_{E^{\prime}}\right)
\]
of $G$, where%
\[
E^{\prime}:=\left\{  e\in E\ \mid\ \text{all endpoints of }e\text{ belong to
}S\right\}  .
\]
In other words, it denotes the submultigraph of $G$ whose vertices are the
elements of $S$, and whose edges are precisely those edges of $G$ whose both
endpoints belong to $S$. We denote this induced submultigraph by $G\left[
S\right]  $.

\item[\textbf{(c)}] An \textbf{induced submultigraph} of $G$ means a
submultigraph of $G$ that is the induced submultigraph of $G$ on $S$ for some
$S\subseteq V$.
\end{enumerate}

The infix \textquotedblleft multi\textquotedblright\ is often omitted. So we
often speak of \textquotedblleft subgraphs\textquotedblright\ instead of
\textquotedblleft submultigraphs\textquotedblright.
\end{definition}

With these definitions, we can now identify cycles in a multigraph with
subgraphs isomorphic to a cycle graph: A cycle of length $n$ in a multigraph
$G$ is \textquotedblleft the same as\textquotedblright\ a submultigraph of $G$
isomorphic to $C_{n}$. (We leave the details to the reader.)

\subsubsection{Disjoint unions}

In Section \ref{sec.sg.djun}, we defined the disjoint union of two or more
simple graphs. The analogous definition for multigraphs is straightforward and
left to the reader.

\subsubsection{Walks}

We already defined walks, paths, closed walks and cycles for multigraphs back
in Section \ref{sec.mg.defs}. The \textbf{length} of a walk is still defined
to be its number of edges. Now, let's see which of their basic properties
(seen in Section \ref{sec.sg.walks}) still hold for multigraphs.

First of all, the edges of a path are still always distinct. This is just as
easy to prove as for simple graphs.

Next, let us see how two walks can be \textquotedblleft
spliced\textquotedblright\ together:

\begin{proposition}
\label{prop.mg.walk-concat}Let $G$ be a multigraph. Let $u$, $v$ and $w$ be
three vertices of $G$. Let $\mathbf{a}=\left(  a_{0},e_{1},a_{1},\ldots
,e_{k},a_{k}\right)  $ be a walk from $u$ to $v$. Let $\mathbf{b}=\left(
b_{0},f_{1},b_{1},\ldots,f_{\ell},b_{\ell}\right)  $ be a walk from $v$ to
$w$. Then,%
\begin{align*}
&  \left(  a_{0},e_{1},a_{1},\ldots,e_{k},a_{k},f_{1},b_{1},f_{2},b_{2}%
,\ldots,f_{\ell},b_{\ell}\right) \\
&  =\left(  a_{0},e_{1},a_{1},\ldots,a_{k-1},e_{k},b_{0},f_{1},b_{1}%
,\ldots,f_{\ell},b_{\ell}\right) \\
&  =\left(  a_{0},e_{1},a_{1},\ldots,a_{k-1},e_{k},v,f_{1},b_{1}%
,\ldots,f_{\ell},b_{\ell}\right)
\end{align*}
is a walk from $u$ to $w$. This walk shall be denoted $\mathbf{a}%
\ast\mathbf{b}$.
\end{proposition}

Walks can be reversed (i.e., walked in backwards direction):

\begin{proposition}
\label{prop.mg.walk-rev}Let $G$ be a multigraph. Let $u$ and $v$ be two
vertices of $G$. Let $\mathbf{a}=\left(  a_{0},e_{1},a_{1},\ldots,e_{k}%
,a_{k}\right)  $ be a walk from $u$ to $v$. Then:

\begin{enumerate}
\item[\textbf{(a)}] The list $\left(  a_{k},e_{k},a_{k-1},e_{k-1},\ldots
,e_{1},a_{0}\right)  $ is a walk from $v$ to $u$. We denote this walk by
$\operatorname*{rev}\mathbf{a}$ and call it the \textbf{reversal} of
$\mathbf{a}$.

\item[\textbf{(b)}] If $\mathbf{a}$ is a path, then $\operatorname*{rev}%
\mathbf{a}$ is a path again.
\end{enumerate}
\end{proposition}

Walks that are not paths contain smaller walks between the same vertices:

\begin{proposition}
\label{prop.mg.walk-to-path-1}Let $G$ be a multigraph. Let $u$ and $v$ be two
vertices of $G$. Let $\mathbf{a}=\left(  a_{0},e_{1},a_{1},\ldots,e_{k}%
,a_{k}\right)  $ be a walk from $u$ to $v$. Assume that $\mathbf{a}$ is not a
path. Then, there exists a walk from $u$ to $v$ whose length is smaller than
$k$.
\end{proposition}

\begin{corollary}
[When there is a walk, there is a path]\label{cor.mg.walk-thus-path}Let $G$ be
a multigraph. Let $u$ and $v$ be two vertices of $G$. Assume that there is a
walk from $u$ to $v$ of length $k$ for some $k\in\mathbb{N}$. Then, there is a
path from $u$ to $v$ of length $\leq k$.
\end{corollary}

All these results can be proved in the same way as their counterparts for
simple graphs (Proposition \ref{prop.sg.walk-concat}, Proposition
\ref{prop.sg.walk-rev}, Proposition \ref{prop.sg.walk-to-path-1} and Corollary
\ref{cor.sg.walk-thus-path}); the only change needed is to record the edges in
the walk.

Corollary \ref{cor.mg.walk-thus-path} can be strengthened a bit:

\begin{corollary}
[When there is a walk, there is a path]\label{cor.mg.walk-thus-path2}Let $G$
be a multigraph. Let $u$ and $v$ be two vertices of $G$. Assume that there is
a walk $\mathbf{w}$ from $u$ to $v$ of length $k$ for some $k\in\mathbb{N}$.
Then, there is a path $\mathbf{p}$ from $u$ to $v$ of length $\leq k$ that is
contained in $\mathbf{w}$ (meaning that each edge of $\mathbf{p}$ is an edge
of $\mathbf{w}$, and each vertex of $\mathbf{p}$ is a vertex of $\mathbf{w}$).
\end{corollary}

\begin{proof}
Replay the proof of Corollary \ref{cor.mg.walk-thus-path} (and the proofs of
the underlying results) and observe that each time we shorten our walk by
\textquotedblleft removing a loop\textquotedblright, we obtain a walk that is
contained in $\mathbf{w}$.
\end{proof}

Given a multigraph $G$ and two vertices $u$ and $v$ of $G$, we can ask
ourselves the same five Questions 1, 2, 3, 4 and 5 that we asked for a simple
graph $G$ in Subsection \ref{subsec.sg.walks.algo}. The answers we gave in
that subsection still apply without requiring substantial changes; the only
necessary modification is that we now have to keep track of the edges in a
path or walk. (The reader can easily fill in the details here.)

\subsubsection{Path-connectedness}

The relation \textquotedblleft path-connected\textquotedblright\ is defined
for multigraphs just as it is for simple graphs (Definition \ref{def.sg.pc}),
and is still denoted $\simeq_{G}$. It is still an equivalence relation (and
the proof is the same as for simple graphs). The following also holds (with
the same proof as for simple graphs):

\begin{proposition}
Let $G$ be a multigraph. Let $u$ and $v$ be two vertices of $G$. Then,
$u\simeq_{G}v$ if and only if there exists a path from $u$ to $v$.
\end{proposition}

The definitions of \textquotedblleft components\textquotedblright\ and
\textquotedblleft connected\textquotedblright\ for multigraphs are the same as
for simple graphs (Definition \ref{def.sg.component} and Definition
\ref{def.sg.connected}). The following propositions can be proved in the same
way as we proved their analogues for simple graphs (Proposition
\ref{prop.sg.component.ind-connected} and Proposition
\ref{prop.sg.component.djun}):

\begin{proposition}
\label{prop.mg.conn.induced-is-conn}Let $G$ be a multigraph. Let $C$ be a
component of $G$.

Then, the multigraph $G\left[  C\right]  $ (that is, the induced submultigraph
of $G$ on the set $C$) is connected.
\end{proposition}

\begin{proposition}
\label{prop.mg.conn.djun-of-comps}Let $G$ be a multigraph. Let $C_{1}%
,C_{2},\ldots,C_{k}$ be all components of $G$ (listed without repetition).

Thus, $G$ is isomorphic to the disjoint union $G\left[  C_{1}\right]  \sqcup
G\left[  C_{2}\right]  \sqcup\cdots\sqcup G\left[  C_{k}\right]  $.
\end{proposition}

The following proposition is an analogue of Proposition
\ref{prop.cyc.btf-walk-cyc} for multigraphs:

\begin{proposition}
\label{prop.mg.cyc.btf-walk-cyc}Let $G$ be a multigraph. Let $\mathbf{w}$ be a
walk of $G$ such that no two consecutive edges of $\mathbf{w}$ are identical.
(By \textquotedblleft consecutive edges\textquotedblright, we mean edges of
the form $e_{i-1}$ and $e_{i}$, where $e_{1},e_{2},\ldots,e_{k}$ are the edges
of $\mathbf{w}$ from first to last.)

Then, $\mathbf{w}$ either is a path or contains a cycle (i.e., there exists a
cycle of $G$ whose edges are edges of $\mathbf{w}$).
\end{proposition}

\begin{proof}
The proof of this proposition for multigraphs is more or less the same as it
was for simple graphs (i.e., as the proof of Proposition
\ref{prop.cyc.btf-walk-cyc}), with a mild difference in how we prove that the
walk $\left(  w_{i},w_{i+1},\ldots,w_{j}\right)  $ is a cycle (of course, this
walk is no longer $\left(  w_{i},w_{i+1},\ldots,w_{j}\right)  $ now, but
rather $\left(  w_{i},e_{i+1},w_{i+1},\ldots,e_{j},w_{j}\right)  $, because
the edges need to be included).\footnote{Here are some details:
\par
We assume that $\mathbf{w}$ is not a path, and we write the walk $\mathbf{w}$
as $\left(  w_{0},e_{1},w_{1},e_{2},w_{2},\ldots,e_{k},w_{k}\right)  $. Then,
there exists a pair $\left(  i,j\right)  $ of integers $i$ and $j$ with $i<j$
and $w_{i}=w_{j}$. Among all such pairs, we pick one with \textbf{minimum}
difference $j-i$. Then, $\left(  w_{i},e_{i+1},w_{i+1},\ldots,e_{j}%
,w_{j}\right)  $ is a closed walk. We claim that this closed walk is a cycle.
\par
To do so, we need to show that
\par
\begin{enumerate}
\item the vertices $w_{i},w_{i+1},\ldots,w_{j-1}$ are distinct;
\par
\item the edges $e_{i+1},e_{i+2},\ldots,e_{j}$ are distinct;
\par
\item we have $j-i\geq1$.
\end{enumerate}
\par
The first of these claims follows from the minimality of $j-i$. The third
follows from $i<j$. It remains to prove the second claim. In other words, it
remains to prove that the edges $e_{i+1},e_{i+2},\ldots,e_{j}$ are distinct,
i.e., that we have $e_{a}\neq e_{b}$ for any two integers $a$ and $b$
satisfying $i<a<b\leq j$. Let us do this. Let $a$ and $b$ be two integers
satisfying $i<a<b\leq j$. We must show that $e_{a}\neq e_{b}$. We distinguish
two cases: the case $a=b-1$ and the case $a\neq b-1$.
\par
\begin{itemize}
\item If $a=b-1$, then $e_{a}$ and $e_{b}$ are two consecutive edges of
$\mathbf{w}$ and thus distinct (since we assumed that no two consecutive edges
of $\mathbf{w}$ are identical). Thus, $e_{a}\neq e_{b}$ is proved in the case
when $a=b-1$.
\par
\item Now, consider the case when $a\neq b-1$. In this case, we must have
$a<b-1$ (since $a<b$ entails $a\leq b-1$). Also, $i\leq a-1$ (since $i<a$).
Hence, $i\leq a-1<a<b-1\leq j-1$ (since $b\leq j$). Therefore, $a-1$, $a$ and
$b-1$ are three distinct elements of the set $\left\{  i,i+1,\ldots
,j-1\right\}  $. Consequently, $w_{a-1},w_{a},w_{b-1}$ are three distinct
vertices (since the vertices $w_{i},w_{i+1},\ldots,w_{j-1}$ are distinct).
Therefore, $w_{b-1}\notin\left\{  w_{a-1},w_{a}\right\}  =\varphi\left(
e_{a}\right)  $ (since $\mathbf{w}$ is a walk, so that the edge $e_{a}$ has
endpoints $w_{a-1}$ and $w_{a}$). However, $\varphi\left(  e_{b}\right)
=\left\{  w_{b-1},w_{b}\right\}  $ (since $\mathbf{w}$ is a walk, so that the
edge $e_{b}$ has endpoints $w_{b-1}$ and $w_{b}$). Now, comparing $w_{b-1}%
\in\left\{  w_{b-1},w_{b}\right\}  =\varphi\left(  e_{b}\right)  $ with
$w_{b-1}\notin\varphi\left(  e_{a}\right)  $, we see that the sets
$\varphi\left(  e_{b}\right)  $ and $\varphi\left(  e_{a}\right)  $ must be
distinct (since $\varphi\left(  e_{b}\right)  $ contains $w_{b-1}$ but
$\varphi\left(  e_{a}\right)  $ does not). In other words, $\varphi\left(
e_{b}\right)  \neq\varphi\left(  e_{a}\right)  $. Hence, $e_{b}\neq e_{a}$. In
other words, $e_{a}\neq e_{b}$. Thus, $e_{a}\neq e_{b}$ is proved in the case
when $a\neq b-1$.
\end{itemize}
\par
We have now proved $e_{a}\neq e_{b}$ in both cases, so we are done.}
\end{proof}

Just as for simple graphs, we get the following corollary:

\begin{corollary}
Let $G$ be a multigraph. Assume that $G$ has a closed walk $\mathbf{w}$ of
length $>0$ such that no two consecutive edges of $\mathbf{w}$ are identical.
Then, $G$ has a cycle.
\end{corollary}

The analogue of Theorem \ref{thm.cyc.two-paths-cyc} for multigraphs is true as well:

\begin{theorem}
\label{prop.mg.cyc.two-paths-cyc}Let $G$ be a multigraph. Let $u$ and $v$ be
two vertices in $G$. Assume that there are two distinct paths from $u$ to $v$.
Then, $G$ has a cycle.
\end{theorem}

\begin{proof}
For simple graphs, this was proved as Theorem \ref{thm.cyc.two-paths-cyc}
above. The same proof applies to multigraphs, once the obvious changes are
made (e.g., instead of $p_{a-1}p_{a}$ and $q_{b-1}q_{b}$, we need to take the
last edges of the two walks $\mathbf{p}$ and $\mathbf{q}$).
\end{proof}

\bigskip

In contrast, Proposition \ref{prop.cycle.deg-d} is \textbf{false} for
multigraphs. In fact, we can take a multigraph with a single vertex and lots
of loops around it. In that case, its degree can be very large, but it has no
cycles of length $>1$.

\subsubsection{$G\setminus e$, bridges and cut-edges}

Next, we extend the definition of $G\setminus e$ (Definition \ref{def.sg.G-e})
to multigraphs:

\begin{definition}
\label{def.mg.G-e}Let $G=\left(  V,E,\varphi\right)  $ be a multigraph. Let
$e$ be an edge of $G$. Then, $G\setminus e$ will mean the graph obtained from
$G$ by removing this edge $e$. In other words,%
\[
G\setminus e:=\left(  V,\ E\setminus\left\{  e\right\}  ,\ \varphi
\mid_{E\setminus\left\{  e\right\}  }\right)  .
\]

\end{definition}

Some authors write $G-e$ for $G\setminus e$.

The analogue of Theorem \ref{thm.G-e.conn} for multigraphs holds (and can be
proved in the same way as Theorem \ref{thm.G-e.conn}):

\begin{theorem}
\label{thm.mg.G-e.conn}Let $G$ be a multigraph. Let $e$ be an edge of $G$. Then:

\begin{enumerate}
\item[\textbf{(a)}] If $e$ is an edge of some cycle of $G$, then the
components of $G\setminus e$ are precisely the components of $G$. (Keep in
mind that the components are sets of vertices. It is these sets that we are
talking about here, not the induced subgraphs on these sets.)

\item[\textbf{(b)}] If $e$ appears in no cycle of $G$ (in other words, there
exists no cycle of $G$ such that $e$ is an edge of this cycle), then the graph
$G\setminus e$ has one more component than $G$.
\end{enumerate}
\end{theorem}

Note that an edge $e$ that is a loop always is an edge of a cycle (indeed, it
creates a cycle of length $1$), and can never appear on any path; thus,
removing such an edge $e$ obviously does not change the path-connectedness relation.

Defining cut-edges and bridges just as we did for simple graphs (Definition
\ref{def.sg.bridges-cuts}), we equally recover the following corollary:

\begin{corollary}
\label{cor.mg.bridge=cut}Let $e$ be an edge of a multigraph $G$. Then, $e$ is
a bridge if and only if $e$ is a cut-edge.
\end{corollary}

\begin{proof}
Just like the proof of Corollary \ref{cor.bridge=cut}.
\end{proof}

\subsubsection{Dominating sets}

We defined and studied dominating sets in Section \ref{sec.sg.dom}. We could
define dominating sets for multigraphs in the same way as for simple graphs,
but we would not get anything new this way. Indeed, if $G$ is a multigraph,
then the dominating sets of $G$ are precisely the dominating sets of
$G^{\operatorname*{simp}}$. Thus, we can reduce any claims about dominating
sets of multigraphs to analogous claims about simple graphs.

\subsubsection{Hamiltonian paths and cycles}

Hamiltonian paths and Hamiltonian cycles were defined for multigraphs in
Definition \ref{def.mg.walks-etc}. It is easy to see that a multigraph $G$ has
a Hamiltonian path or Hamiltonian cycle if and only if the corresponding
simple graph $G^{\operatorname*{simp}}$ has one, with only one exception: A
multigraph $G$ with a single vertex has a Hamiltonian cycle if it has a loop,
but the corresponding simple graph $G^{\operatorname*{simp}}$ does not have a
Hamiltonian cycle.

This does not mean, however, that everything we proved about Hamiltonian paths
still applies to multigraphs. For instance, neither Ore's theorem (Theorem
\ref{thm.hamc.ore}) nor Dirac's theorem (Corollary \ref{cor.hamc.dirac}) holds
for multigraphs, because we could duplicate edges to make degrees arbitrarily
large, without necessarily creating a hamc.

Proposition \ref{prop.hamc.nec1} still holds for multigraphs, but this is
clear because it can be derived from the corresponding property of
$G^{\operatorname*{simp}}$.

\subsubsection{Exercises}

\begin{exercise}
\label{exe.3.1}Which of the Exercises \ref{exe.1.1}, \ref{exe.1.2},
\ref{exe.1.3}, \ref{exe.1.4}, \ref{exe.1.6}, \ref{exe.1.7} and \ref{exe.1.8}
remain true if \textquotedblleft simple graph\textquotedblright\ is replaced
by \textquotedblleft multigraph\textquotedblright?

(For each exercise that becomes false, provide a counterexample. For each
exercise that remains true, either provide a new solution that works for
multigraphs, or argue that the solution we have seen applies verbatim to
multigraphs, or derive the multigraph case from the simple graph case.)
\end{exercise}

\begin{exercise}
\label{exe.mt1.deg-cycle-2} Let $G$ be a multigraph with at least one edge.
Assume that each vertex of $G$ has even degree. Prove that $G$ has a
cycle.\medskip

[\textbf{Solution:} This is Exercise 4 on midterm \#1 from my Spring 2017
course; see \href{https://www.cip.ifi.lmu.de/~grinberg/t/17s/}{the course
page} for solutions.]
\end{exercise}

\begin{exercise}
\label{exe.3.2}Let $G$ be a multigraph with at least one vertex. Let $d>2$ be
an integer. Assume that $\deg v>2$ for each vertex $v$ of $G$. Prove that $G$
has a cycle whose length is not divisible by $d$.
\end{exercise}

\begin{exercise}
\label{exe.3.10}Let $G$ be a multigraph. Assume that $G$ has exactly two
vertices of odd degree. Prove that these two vertices are path-connected.
\end{exercise}

\begin{exercise}
\label{exe.3.3}Let $G=\left(  V,E,\varphi\right)  $ be a multigraph that has
no loops.

If $e\in E$ is an edge that contains a vertex $v\in V$, then we let $e/v$
denote the endpoint of $e$ distinct from $v$.

\begin{noncompile}
(If $e$ is a loop, then this is understood to mean $v$ itself.)
\end{noncompile}

For each $v \in V$, we define a rational number $q_{v}$ by
\[
q_{v} = \sum_{\substack{e \in E; \\v \in\varphi\left(  e \right)  }}
\dfrac{\deg\left(  e/v \right)  }{\deg v}.
\]
(Note that the denominator $\deg v$ on the right hand side is nonzero whenever
the sum is nonempty!)

(Thus, $q_{v}$ is the average degree of the neighbors of $v$, weighted with
the number of edges that join $v$ to the respective neighbors. If $v$ has no
neighbors, then $q_{v} = 0$.)

Prove that
\[
\sum_{v \in V} q_{v} \geq\sum_{v \in V} \deg v.
\]

(In other words, in a social network, your average friend has, on average,
more friends than you do!) \medskip

[\textbf{Hint:} Any positive reals $x$ and $y$ satisfy $\dfrac{x}{y}+\dfrac
{y}{x}\geq2$. Why, and how does this help?]
\end{exercise}

\begin{exercise}
\label{exe.5.3}Let $F$ be any field. (For instance, $F$ can be $\mathbb{Q}$ or
$\mathbb{R}$ or $\mathbb{C}$.)

Let $G = \left(  V, E, \varphi\right)  $ be a multigraph, where $V = \left\{
1, 2, \ldots, n \right\}  $ for some $n \in\mathbb{N}$.

For each edge $e \in E$, we construct a column vector $\chi_{e} \in F^{n}$
(that is, a column vector with $n$ entries) as follows:

\begin{itemize}
\item If $e$ is a loop, then we let $\chi_{e}$ be the zero vector.

\item Otherwise, we let $u$ and $v$ be the two endpoints of $e$, and we let
$\chi_{e}$ be the column vector that has a $1$ in its $u$-th position, a $-1$
in its $v$-th position, and $0$s in all other positions. (This depends on
which endpoint we call $u$ and which endpoint we call $v$, but we just make
some choice and stick with it. The result will be true no matter how we choose.)
\end{itemize}

Let $M$ be the $n\times\left\vert E\right\vert $-matrix over $F$ whose columns
are the column vectors $\chi_{e}$ for all $e\in E$ (we order them in some way;
the exact order doesn't matter). Prove that
\[
\operatorname*{rank}M=\left\vert V\right\vert -\operatorname*{conn}G,
\]
where $\operatorname*{conn}G$ denotes the number of components of $G$.
\medskip

[\textbf{Example:} Here is an example: Let $G$ be the multigraph
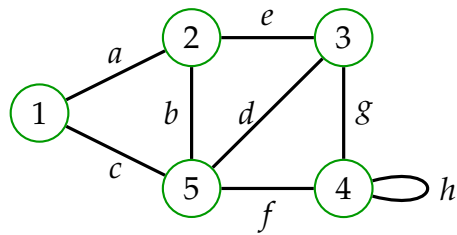
\[
\begin{tikzpicture}[scale=2]
\begin{scope}[every node/.style={circle,thick,draw=green!60!black}]
\node(1) at (-1,0) {$1$};
\node(2) at (0,0.5) {$2$};
\node(3) at (1,0.5) {$3$};
\node(4) at (1,-0.5) {$4$};
\node(5) at (0,-0.5) {$5$};
\end{scope}
\begin{scope}[every edge/.style={draw=black,very thick}, every loop/.style={}]
\path[-] (1) edge node[above] {$a$} (2);
\path[-] (1) edge node[below] {$c$} (5);
\path[-] (2) edge node[left] {$b$} (5);
\path[-] (2) edge node[above] {$e$} (3);
\path[-] (3) edge node[left] {$d$} (5);
\path[-] (4) edge node[below] {$f$} (5);
\path[-] (3) edge node[right] {$g$} (4);
\path[-] (4) edge[loop right] node[right] {$h$} (4);
\end{scope}
\end{tikzpicture}
\]
(so that $n=5$). Then, if we choose the endpoints of $b$ to be $2$ and $5$ in
this order, then we have $\chi_{b}=%
\begin{pmatrix}
0\\
1\\
0\\
0\\
-1
\end{pmatrix}
$. (Choosing them to be $5$ and $2$ instead, we would obtain $\chi_{b}=%
\begin{pmatrix}
0\\
-1\\
0\\
0\\
1
\end{pmatrix}
$.) If we do the same for all edges of $G$ (that is, we choose the smaller
endpoint as $u$ and the larger endpoint as $v$), and if we order the columns
so that they correspond to the edges $a,b,c,d,e,f,g,h$ from left to right,
then the matrix $M$ comes out as follows:
\[
M=%
\begin{pmatrix}
1 & 0 & 1 & 0 & 0 & 0 & 0 & 0\\
-1 & 1 & 0 & 0 & 1 & 0 & 0 & 0\\
0 & 0 & 0 & 1 & -1 & 0 & 1 & 0\\
0 & 0 & 0 & 0 & 0 & 1 & -1 & 0\\
0 & -1 & -1 & -1 & 0 & -1 & 0 & 0
\end{pmatrix}
.
\]
It is easy to see that $\operatorname*{rank}M=4$, which is precisely
$\left\vert V\right\vert -\operatorname*{conn}G$.] \medskip

[\textbf{Remark:} The claim of the exercise can be restated as follows: The
span of the vectors $\chi_{e}$ for all $e\in E$ has dimension $\left\vert
V\right\vert -\operatorname*{conn}G$.

Topologists will recognize the matrix $M$ as (a matrix that represents) the
boundary operator $\partial:C_{1}\left(  G\right)  \rightarrow C_{0}\left(
G\right)  $, where $G$ is viewed as a CW-complex.]
\end{exercise}

\begin{exercise}
\label{exe.hw3.conn} If $G$ is a multigraph, then $\operatorname{conn} G$
shall denote the number of connected components of $G$. (Note that this is $0$
when $G$ has no vertices, and $1$ if $G$ is connected.)

Let $\left(  V,H,\varphi\right)  $ be a multigraph. Let $E$ and $F$ be two
subsets of $H$.

\begin{enumerate}
\item[\textbf{(a)}] Prove that
\begin{align}
&  \operatorname{conn}\left(  V,\ E,\ \varphi\mid_{E}\right)
+\operatorname{conn}\left(  V,\ F,\ \varphi\mid_{F}\right) \nonumber\\
&  \leq\operatorname{conn}\left(  V,\ E\cup F,\ \varphi\mid_{E\cup F}\right)
+\operatorname{conn}\left(  V,\ E\cap F,\ \varphi\mid_{E\cap F}\right)  .
\label{eq.exe.hw3.conn.ineq}%
\end{align}

[\textbf{Hint:} Feel free to restrict yourself to the case of a simple graph;
in this case, $E$ and $F$ are two subsets of $\mathcal{P}_{2}\left(  V\right)
$, and you have to show that
\[
\operatorname{conn}\left(  V,\ E\right)  +\operatorname{conn}\left(
V,\ F\right)  \leq\operatorname{conn}\left(  V,\ E\cup F\right)
+\operatorname{conn}\left(  V,\ E\cap F\right)  .
\]
This isn't any easier than the general case, but saves you the hassle of
carrying the map $\varphi$ around.]

\item[\textbf{(b)}] Give an example where the inequality
\eqref{eq.exe.hw3.conn.ineq} does \textbf{not} become an equality.
\end{enumerate}

[\textbf{Solution:} This is Exercise 3 on homework set \#3 from my Spring 2017
course; see \href{https://www.cip.ifi.lmu.de/~grinberg/t/17s/}{the course
page} for solutions.]
\end{exercise}

\begin{exercise}
\label{exe.3.5}Let $G=\left(  V,E,\varphi\right)  $ be a connected multigraph
with $2m$ edges, where $m\in\mathbb{N}$. A set $\left\{  e,f\right\}  $ of two
distinct edges will be called a \textbf{friendly couple} if $e$ and $f$ have
at least one endpoint in common. Prove that the edge set of $G$ can be
decomposed into $m$ disjoint friendly couples (i.e., there exist $m$ disjoint
friendly couples $\left\{  e_{1},f_{1}\right\}  ,\left\{  e_{2},f_{2}\right\}
,\ldots,\left\{  e_{m},f_{m}\right\}  $ such that $E=\left\{  e_{1}%
,f_{1},e_{2},f_{2},\ldots,e_{m},f_{m}\right\}  $). (\textquotedblleft
Disjoint\textquotedblright\ means \textquotedblleft disjoint as
sets\textquotedblright\ -- i.e., having no edges in common.)

[\textbf{Example:} Here is a graph with an even number of edges:
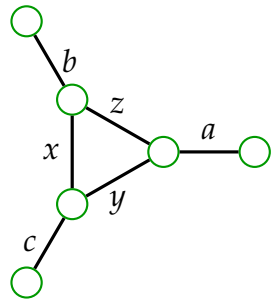
\[
\begin{tikzpicture}
\begin{scope}[every node/.style={circle,thick,draw=green!60!black}]
\node(A) at (0:0.8) {};
\node(B) at (120:0.8) {};
\node(C) at (240:0.8) {};
\node(D) at (0:2) {};
\node(E) at (120:2) {};
\node(F) at (240:2) {} ;
\end{scope}
\begin{scope}[every edge/.style={draw=black,very thick}]
\path[-] (A) edge node[above] {$z$} (B);
\path[-] (B) edge node[left] {$x$} (C);
\path[-] (C) edge node[below] {$y$} (A);
\path[-] (D) edge node[above] {$a$} (A);
\path[-] (E) edge node[right] {$b$} (B);
\path[-] (F) edge node[left] {$c$} (C);
\end{scope}
\end{tikzpicture}
\]
One possible decomposition of its edge set into disjoint friendly couples is
$\left\{  a,y\right\}  ,\left\{  b,z\right\}  ,\left\{  c,x\right\}  $.]
\medskip

[\textbf{Hint:} Induct on $\left\vert E\right\vert $. Pick a vertex $v$ of
degree $>1$ and consider the components of $G\setminus v$.]
\end{exercise}

\begin{exercise}
\label{exe.3.6}Let $n\geq0$. Let $d_{1},d_{2},\ldots,d_{n}$ be $n$ nonnegative
integers such that $d_{1}+d_{2}+\cdots+d_{n}$ is even.

\begin{enumerate}
\item[\textbf{(a)}] Prove that there exists a multigraph $G$ with vertex set
$\left\{  1,2,\ldots,n\right\}  $ such that all $i\in\left\{  1,2,\ldots
,n\right\}  $ satisfy $\deg i=d_{i}$.

\item[\textbf{(b)}] Prove that there exists a loopless multigraph $G$ with
vertex set $\left\{  1,2,\ldots,n\right\}  $ such that all $i\in\left\{
1,2,\ldots,n\right\}  $ satisfy $\deg i=d_{i}$ if and only if each
$i\in\left\{  1,2,\ldots,n\right\}  $ satisfies the inequality
\begin{equation}
\sum_{\substack{j\in\left\{  1,2,\ldots,n\right\}  ;\\j\neq i}}d_{j}\geq
d_{i}. \label{eq.deg-mulg-seq.n-gon}%
\end{equation}

\end{enumerate}

[\textbf{Remark:} The inequality \eqref{eq.deg-mulg-seq.n-gon} is the
\textquotedblleft$n$-gon inequality\textquotedblright: It is equivalent to the
existence of a (possibly degenerate) $n$-gon with sidelengths $d_{1}%
,d_{2},\ldots,d_{n}$.]
\end{exercise}

\begin{exercise}
\label{exe.2.7}Let $G$ be a loopless multigraph. Recall that a \textbf{trail}
(in $G$) means a walk whose edges are distinct (but whose vertices are not
necessarily distinct). Let $u$ and $v$ be two vertices of $G$. As usual,
\textquotedblleft trail from $u$ to $v$\textquotedblright\ means
\textquotedblleft trail that starts at $u$ and ends at $v$\textquotedblright.
Prove that%
\begin{align*}
&  \left(  \text{the number of trails from }u\text{ to }v\text{ in }G\right)
\\
&  \equiv\left(  \text{the number of paths from }u\text{ to }v\text{ in
}G\right)  \operatorname{mod}2.
\end{align*}

[\textbf{Hint:} Try to pair up the non-path trails into pairs. Make sure to
prove that this pairing is well-defined (i.e., each non-path trail
$\mathbf{t}$ has exactly one partner, which is not itself, and that
$\mathbf{t}$ is the designated partner of its partner!).]
\end{exercise}

\begin{exercise}
\label{exe.3.9}Let $G$ be a multigraph such that every vertex of $G$ has even
degree. Let $u$ and $v$ be two distinct vertices of $G$. Prove that the number
of paths from $u$ to $v$ is even. \medskip

[\textbf{Hint:} When you add an edge joining $u$ to $v$, the graph $G$ becomes
a graph with exactly two odd-degree vertices $u$ and $v$, and the claim
becomes ``the number of paths from $u$ to $v$ is odd'' (why?). In this form,
the claim turns out to be easier to prove. Indeed, any path must start with
some edge...

Keep in mind that paths can be replaced by trails, by Exercise \ref{exe.2.7}.]
\end{exercise}

\begin{exercise}
\label{exe.mt1.cyctree} Let $G=\left(  V,E,\varphi\right)  $ be a multigraph
such that $\left\vert E\right\vert >\left\vert V\right\vert $. Prove that $G$
has a cycle of length $\leq\dfrac{2n+2}{3}$, where $n=\left\vert V\right\vert
$. \medskip

[\textbf{Solution:} This is Exercise 8 on midterm \#3 from my Spring 2017
course (except that the simple graph was replaced by a multigraph); see
\href{https://www.cip.ifi.lmu.de/~grinberg/t/17s/}{the course page} for solutions.]
\end{exercise}

\subsection{\label{sec.mg.euler}Eulerian circuits and walks}

\subsubsection{Definitions}

Let us now move on to a new feature of multigraphs, one that we have not yet
studied (even for simple graphs).

Recall that a Hamiltonian path or cycle is a path or cycle that contains all
vertices of the graph. Being a path or cycle, it has to contain each of them
exactly once (except, in the case of a cycle, of its starting point).

What about a walk or closed walk that contains all \textbf{edges} exactly once
instead? These are called \textquotedblleft Eulerian\textquotedblright\ walks
or circuits; here is the formal definition:

\begin{definition}
Let $G$ be a multigraph.

\begin{enumerate}
\item[\textbf{(a)}] A walk of $G$ is said to be \textbf{Eulerian} if each edge
of $G$ appears exactly once in this walk.

(In other words: A walk $\left(  v_{0},e_{1},v_{1},e_{2},v_{2},\ldots
,e_{k},v_{k}\right)  $ of $G$ is said to be \textbf{Eulerian} if for each edge
$e$ of $G$, there exists exactly one $i\in\left\{  1,2,\ldots,k\right\}  $
such that $e=e_{i}$.)

\item[\textbf{(b)}] An \textbf{Eulerian circuit} of $G$ means a circuit (i.e.,
closed walk) of $G$ that is Eulerian. (Strictly speaking, the preceding
sentence is redundant, but we still said it to stress the notion of an
Eulerian circuit.)
\end{enumerate}
\end{definition}

Unlike for Hamiltonian paths and cycles, an Eulerian walk or circuit is
usually not a path or cycle. Also, finding an Eulerian walk in a multigraph
$G$ is not the same as finding an Eulerian walk in the simple graph
$G^{\operatorname*{simp}}$. (Nevertheless, some authors call Eulerian walks
\textquotedblleft Eulerian paths\textquotedblright\ and call Eulerian circuits
\textquotedblleft Eulerian cycles\textquotedblright. This is rather confusing.)

\Needspace{51pc}

\begin{example}
\label{exa.eulerian.exa1}Consider the following multigraphs:%
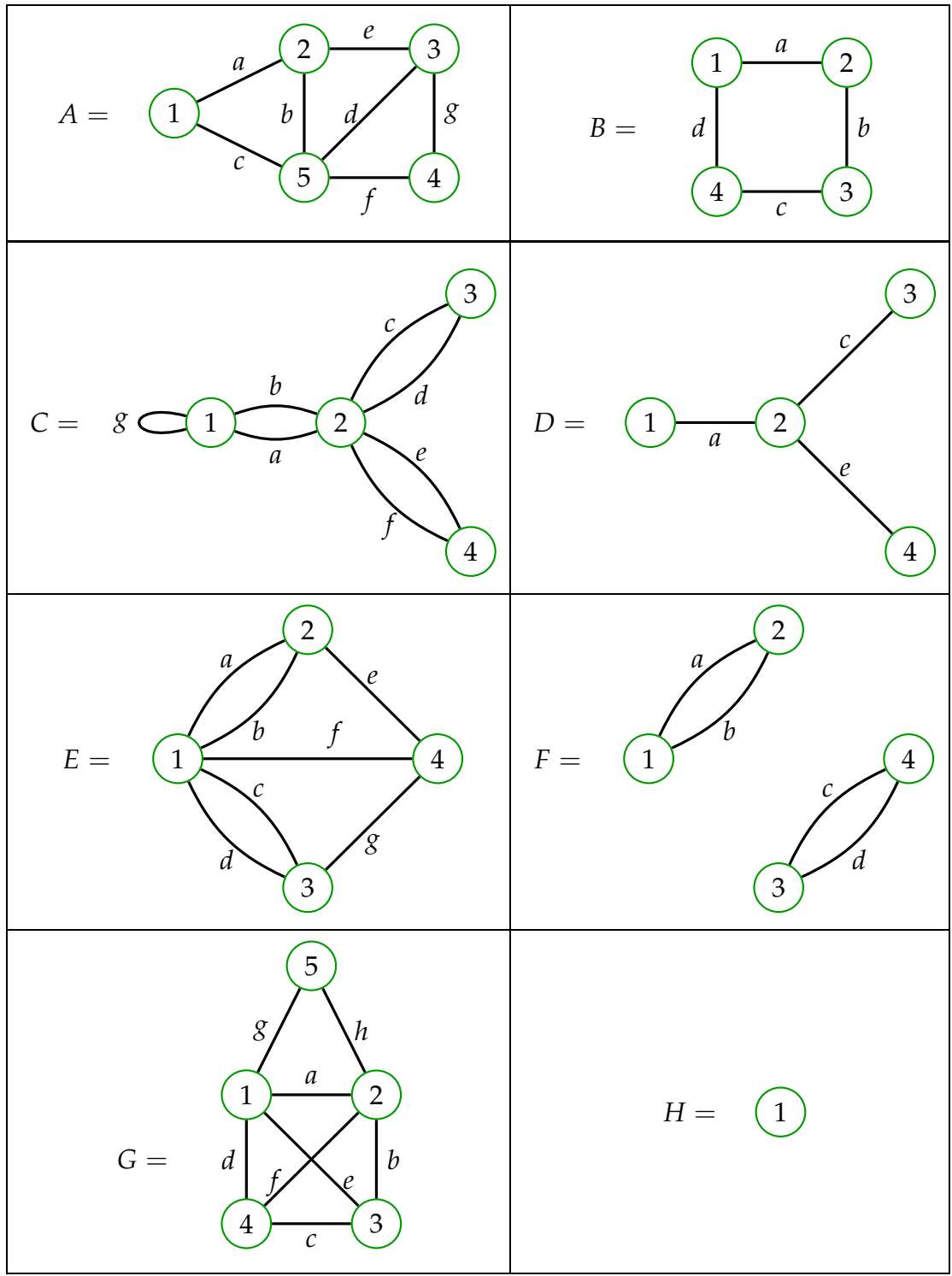
\[%

\ \
\]

\begin{itemize}
\item The multigraph $A$ has an Eulerian walk $\left(
3,d,5,b,2,e,3,g,4,f,5,c,1,a,2\right)  $. But $A$ has no Eulerian circuit. The
easiest way to see this is by observing that $A$ has a vertex of odd degree
(e.g., the vertex $2$). If an Eulerian circuit were to exist, then it would
have to enter this vertex as often as it exited it; but this would mean that
the degree of this vertex would be even (because each edge containing this
vertex would be used exactly once either to enter or to exit it, except for
loops, which would be used twice). So, more generally, any multigraph that has
a vertex of odd degree cannot have an Eulerian circuit.

\item The multigraph $B$ has an Eulerian circuit $\left(
1,a,2,b,3,c,4,d,1\right)  $, and thus of course an Eulerian walk (since any
Eulerian circuit is an Eulerian walk).

\item The multigraph $C$ has an Eulerian circuit $\left(
1,g,1,b,2,c,3,d,2,e,4,f,2,a,1\right)  $.

\item The multigraph $D$ has no Eulerian walk. Indeed, it has four vertices of
odd degree. If $v$ is a vertex of odd degree, then any Eulerian walk has to
either start or end at $v$ (since otherwise, the walk would enter and leave
$v$ equally often, but then the degree of $v$ would be even). But a walk can
only have one starting point and one ending point. This allows for two
vertices of odd degree, but not more than two. So, more generally, any
multigraph that has more than two vertices of odd degree cannot have an
Eulerian walk.

\item The multigraph $E$ has no Eulerian walk. The reason is the same as for
$D$. Note that $E$ is the famous multigraph of bridges in K\"{o}nigsberg, as
studied by Euler in 1736 (see
\href{https://en.wikipedia.org/wiki/Seven_Bridges_of_Konigsberg}{the Wikipedia
page for \textquotedblleft Seven bridges of K\"{o}nigsberg\textquotedblright}
for the backstory).

\item The multigraph $F$ has no Eulerian walk, since it has two components,
each containing at least one edge. (An Eulerian walk would have to contain
both edges $b$ and $c$, but there is no way to walk between them, since they
belong to different components.)

\item The multigraph $G$ has an Eulerian walk, namely $\left(
3,b,2,h,5,g,1,a,2,f,4,d,1,e,3,c,4\right)  $. It has no Eulerian circuit, since
it has two vertices of odd degree.

\item The multigraph $H$ has an Eulerian circuit, namely $\left(  1\right)  $.
\end{itemize}
\end{example}

\begin{remark}
For the pedants: A multigraph can have an Eulerian circuit even if it is not
connected, as long as all its edges belong to the same component (i.e., all
but one components are just singletons with no edges). Here is an example:
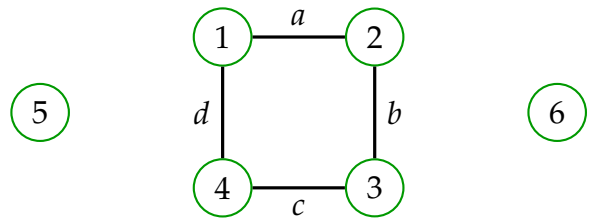
\[%
\begin{tikzpicture}[scale=2]
\begin{scope}[every node/.style={circle,thick,draw=green!60!black}]
\node(1) at (0,0.5) {$1$};
\node(2) at (1,0.5) {$2$};
\node(3) at (1,-0.5) {$3$};
\node(4) at (0,-0.5) {$4$};
\node(5) at (-1.2, 0) {$5$};
\node(6) at (2.2, 0) {$6$};
\end{scope}
\begin{scope}[every edge/.style={draw=black,very thick}, every loop/.style={}]
\path[-] (1) edge node[above] {$a$} (2);
\path[-] (2) edge node[right] {$b$} (3);
\path[-] (3) edge node[below] {$c$} (4);
\path[-] (4) edge node[left] {$d$} (1);
\end{scope}
\end{tikzpicture}%
\]

\end{remark}

\begin{exercise}
\label{exe.eulertrails.Kn} Let $n$ be a positive integer. Recall from
Definition \ref{def.sg.complete-empty} \textbf{(a)} that $K_{n}$ denotes the
complete graph on $n$ vertices. This is the graph with vertex set $V=\left\{
1,2,\ldots,n\right\}  $ and edge set $\mathcal{P}_{2}\left(  V\right)  $ (so
each two distinct vertices are adjacent).

Find Eulerian circuits for the graphs $K_{3}$, $K_{5}$, and $K_{7}$.\medskip

[\textbf{Solution:} This is Exercise 2 on homework set \#2 from my Spring 2017
course; see \href{https://www.cip.ifi.lmu.de/~grinberg/t/17s/}{the course
page} for solutions.]
\end{exercise}

\subsubsection{\label{subsec.mg.euler.eh}The Euler--Hierholzer theorem}

How hard is it to find an Eulerian walk or circuit in a multigraph, or to
check if there is any? Surprisingly, this is a lot easier than the same
questions for Hamiltonian paths or cycles. The second question in particular
is answered (for connected multigraphs) by the \textbf{Euler--Hierholzer
theorem}:

\begin{theorem}
[Euler, Hierholzer]\label{thm.mg.euler}Let $G$ be a connected multigraph. Then:

\begin{enumerate}
\item[\textbf{(a)}] The multigraph $G$ has an Eulerian circuit if and only if
each vertex of $G$ has even degree.

\item[\textbf{(b)}] The multigraph $G$ has an Eulerian walk if and only if all
but at most two vertices of $G$ have even degree.
\end{enumerate}
\end{theorem}

We already proved the \textquotedblleft$\Longrightarrow$\textquotedblright%
\ directions of both parts \textbf{(a)} and \textbf{(b)} in Example
\ref{exa.eulerian.exa1}. It remains to prove the \textquotedblleft%
$\Longleftarrow$\textquotedblright\ directions. Euler did not actually prove
them in his 1736 paper\footnote{although he pretended to do so, in a
beautifully transparent example of \textquotedblleft proof by
trivialization\textquotedblright\ \cite[\S 21]{Euler36}}, but Hierholzer did
in his 1873 paper \cite{Hierho73}. The \textquotedblleft
standard\textquotedblright\ proof can be found in many texts, such as
\cite[Theorem 5.2.2 and Theorem 5.2.3]{Guicha16}. I will sketch a different
proof, which I learnt from \cite[Problem 12.35]{LeLeMe} and which is actually
rather close to Hierholzer's original one \cite{Hierho73}. We begin with the
following definition:

\begin{definition}
\label{def.mg.trail}Let $G$ be a multigraph. A \textbf{trail} of $G$ means a
walk of $G$ whose edges are distinct.
\end{definition}

So a trail can repeat vertices, but cannot repeat edges.

Thus, an Eulerian walk has to be a trail. A trail cannot be longer than an
Eulerian walk. Hence, a reasonable way to try constructing an Eulerian walk is
to start with some trail, and make it progressively longer until it becomes
Eulerian (hopefully).

This suggests the following approach to proving the \textquotedblleft%
$\Longleftarrow$\textquotedblright\ directions of Theorem \ref{thm.mg.euler}:
We pick the longest trail of $G$ and argue that (under the right assumptions)
it has to be Eulerian, since otherwise there would be a way to make it longer.
Of course, we need to find such a way. Here is the first step:

\begin{lemma}
\label{lem.mg.euler.longest-trail-exists}Let $G$ be a multigraph with at least
one vertex. Then, $G$ has a longest trail.
\end{lemma}

\begin{proof}
Clearly, $G$ has at least one trail (e.g., a length-$0$ trail from a vertex to
itself). Moreover, $G$ has only finitely many trails (since each edge of $G$
can only be used once in a trail, and there are only finitely many edges).
Hence, the maximum principle proves the lemma.
\end{proof}

Our goal now is to show that under appropriate conditions, such a longest
trail will be Eulerian. This will require two further lemmas.

First, one more piece of notation: We say that an edge $e$ of a multigraph $G$
\textbf{intersects} a walk $\mathbf{w}$ if at least one endpoint of $e$ is a
vertex of $\mathbf{w}$. Here is how this can look like:%
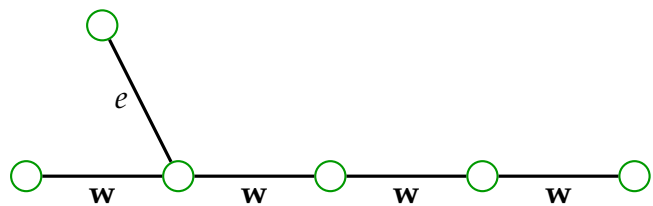
\[%
%
\]
(here, the edges of $\mathbf{w}$ are marked with a \textquotedblleft%
$\mathbf{w}$\textquotedblright\ underneath them) or%
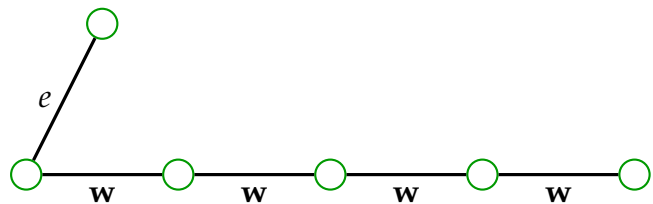
\[%
%
\]
(here, the endpoint of $e$ that is a vertex of $\mathbf{w}$ happens to be the
starting point of $\mathbf{w}$) or%
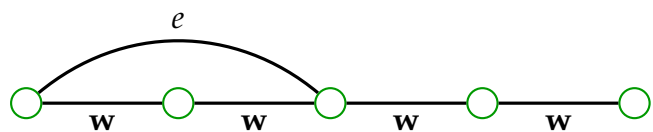
\[%
%
\]
(here, both endpoints of $e$ happen to be vertices of $\mathbf{w}$). Be
careful with such pictures, though: A walk doesn't have to be a path; it can
visit a vertex any number of times!

\begin{lemma}
\label{lem.mg.euler.walk-intersects-edge}Let $G$ be a connected multigraph.
Let $\mathbf{w}$ be a walk of $G$. Assume that there exists an edge of $G$
that is not an edge of $\mathbf{w}$.

Then, there exists an edge of $G$ that is not an edge of $\mathbf{w}$ but
intersects $\mathbf{w}$.
\end{lemma}

\begin{proof}
We assumed that there exists an edge of $G$ that is not an edge of
$\mathbf{w}$. Pick such an edge, and call it $f$.

A \textquotedblleft$\mathbf{w}$-$f$-path\textquotedblright\ will mean a path
from a vertex of $\mathbf{w}$ to an endpoint of $f$. Such a path clearly
exists, since $G$ is connected. Thus, we can pick a \textbf{shortest} such
path. If this shortest path has length $0$, then we are done (since $f$
intersects $\mathbf{w}$ in this case). If not, we consider the first edge of
this path. This first edge cannot be an edge of $\mathbf{w}$, because
otherwise we could remove it from the path and get an even shorter
$\mathbf{w}$-$f$-path. But it clearly intersects $\mathbf{w}$. So we have
found an edge of $G$ that is not an edge of $\mathbf{w}$ but intersects
$\mathbf{w}$. This proves the lemma.
\end{proof}

\begin{lemma}
\label{lem.mg.euler.longest-trail-closed}Let $G$ be a multigraph such that
each vertex of $G$ has even degree. Let $\mathbf{w}$ be a longest trail of
$G$. Then, $\mathbf{w}$ is a closed walk.
\end{lemma}

\begin{proof}
Assume the contrary. Let $u$ be the starting point and $v$ the ending point of
$\mathbf{w}$. Since we assumed that $\mathbf{w}$ is not a closed walk, we thus
have $u\neq v$.

Consider the edges of $\mathbf{w}$ that contain $v$. Such edges are of two
kinds: those by which $\mathbf{w}$ enters $v$ (this means that $v$ comes
immediately after this edge in $\mathbf{w}$), and those by which $\mathbf{w}$
leaves $v$ (this means that $v$ comes immediately before this edge in
$\mathbf{w}$).\ \ \ \ \footnote{Loops whose only endpoint is $v$ count as
both.} Except for the very last edge of $\mathbf{w}$, each edge of the former
kind is immediately followed by an edge of the latter kind; conversely, each
edge of the latter kind is immediately preceded by an edge of the former kind
(since $\mathbf{w}$ starts at the vertex $u$, which is distinct from $v$).
Hence, the walk $\mathbf{w}$ has exactly one more edge entering $v$ than it
has edges leaving $v$. Thus, the number of edges of $\mathbf{w}$ that contain
$v$ (with loops counting twice) is odd. However, the total number of edges of
$G$ that contain $v$ (with loops counting twice) is even (because it is the
degree of $v$, but we assumed that each vertex of $G$ has even degree). So
these two numbers are distinct. Thus, there is at least one edge of $G$ that
contains $v$ but is not an edge of $\mathbf{w}$.

Fix such an edge and call it $f$. Now, append $f$ to the trail $\mathbf{w}$ at
the end. The result will be a trail (since $f$ is not an edge of $\mathbf{w}$)
that is longer than $\mathbf{w}$. But this contradicts the fact that
$\mathbf{w}$ is a longest trail. Thus, the lemma is proved.
\end{proof}

We can now finish the proof of the Euler--Hierholzer theorem:

\begin{proof}
[Proof of Theorem \ref{thm.mg.euler}.]\textbf{(a)} $\Longrightarrow:$ We
proved this back in Example \ref{exa.eulerian.exa1}. \medskip

$\Longleftarrow:$ Assume that each vertex of $G$ has even degree.

By Lemma \ref{lem.mg.euler.longest-trail-exists}, we know that $G$ has a
longest trail. Fix such a longest trail, and call it $\mathbf{w}$. Then, Lemma
\ref{lem.mg.euler.longest-trail-closed} shows that $\mathbf{w}$ is a closed walk.

We claim that $\mathbf{w}$ is Eulerian. Indeed, assume the contrary. Then,
there exists an edge of $G$ that is not an edge of $\mathbf{w}$. Hence, Lemma
\ref{lem.mg.euler.walk-intersects-edge} shows that there exists an edge of $G$
that is not an edge of $\mathbf{w}$ but intersects $\mathbf{w}$. Fix such an
edge, and call it $f$.

Since $f$ intersects $\mathbf{w}$, there exists an endpoint $v$ of $f$ that is
a vertex of $\mathbf{w}$. Consider this $v$. Since $\mathbf{w}$ is a
\textbf{closed} trail, we can WLOG assume that $\mathbf{w}$ starts and ends at
$v$ (since we can otherwise achieve this by
rotating\footnote{\textbf{Rotating} a closed walk $\left(  w_{0},e_{1}%
,w_{1},e_{2},w_{2},\ldots,e_{k},w_{k}\right)  $ means moving its first vertex
and its first edge to the end, i.e., replacing the walk by $\left(
w_{1},e_{2},w_{2},e_{3},w_{3},\ldots,e_{k},w_{k},e_{1},w_{1}\right)  $. This
always results in a closed walk again. For example, if $\left(
1,a,2,b,3,c,1\right)  $ is a closed walk, then we can rotate it to obtain
$\left(  2,b,3,c,1,a,2\right)  $; then, rotating it one more time, we obtain
$\left(  3,c,1,a,2,b,3\right)  $.
\par
Clearly, by rotating a closed walk several times, we can make it start at any
of its vertices. Moreover, if we rotate a closed trail, then we obtain a
closed trail.} $\mathbf{w}$). Then, we can append the edge $f$ to the trail
$\mathbf{w}$. This results in a new trail (since $f$ is not an edge of
$\mathbf{w}$) that is longer than $\mathbf{w}$. And this contradicts the fact
that $\mathbf{w}$ is a longest trail of $G$.

This contradiction proves that $\mathbf{w}$ is Eulerian. Hence, $\mathbf{w}$
is an Eulerian circuit (since $\mathbf{w}$ is a closed walk). Thus, the
\textquotedblleft$\Longleftarrow$\textquotedblright\ direction of Theorem
\ref{thm.mg.euler} \textbf{(a)} is proven. \medskip

\textbf{(b)} $\Longrightarrow:$ Already proved in Example
\ref{exa.eulerian.exa1}.\medskip

$\Longleftarrow:$ Assume that all but at most two vertices of $G$ have even
degree. We must prove that $G$ has an Eulerian walk.

If each vertex of $G$ has even degree, then this follows from Theorem
\ref{thm.mg.euler} \textbf{(a)}, since every Eulerian circuit is an Eulerian
walk. Thus, we WLOG assume that not each vertex of $G$ has even degree. In
other words, the number of vertices of $G$ having odd degree is positive.

The handshake lemma for multigraphs (i.e., Corollary \ref{cor.mg.odd-deg-even}%
) shows that the number of vertices of $G$ having odd degree is even.
Furthermore, this number is at most $2$ (since all but at most two vertices of
$G$ have even degree). So this number is even, positive and at most $2$. Thus,
this number is $2$. In other words, the multigraph $G$ has exactly two
vertices having odd degree. Let $u$ and $v$ be these two vertices.

Add a new edge $e$ that has endpoints $u$ and $v$ to the multigraph $G$ (do
this even if there already is such an edge!\footnote{This is a time to be
grateful for the notion of a multigraph. We could not do this with simple
graphs!}). Let $G^{\prime}$ denote the resulting multigraph. Then, in
$G^{\prime}$, each vertex has even degree (since the newly added edge $e$ has
increased the degrees of $u$ and $v$ by $1$, thus turning them from odd to
even). Moreover, $G^{\prime}$ is still connected (since $G$ was connected, and
the newly added edge $e$ can hardly take that away). Thus, we can apply
Theorem \ref{thm.mg.euler} \textbf{(a)} to $G^{\prime}$ instead of $G$. As a
result, we conclude that $G^{\prime}$ has an Eulerian circuit. Cutting the
newly added edge $e$ out of this Eulerian circuit\footnote{More precisely: We
rotate this circuit until $e$ becomes its last edge, and then we remove this
last edge to obtain a walk.}, we obtain an Eulerian walk of $G$. Hence, $G$
has an Eulerian walk. Thus, the \textquotedblleft$\Longleftarrow
$\textquotedblright\ direction of Theorem \ref{thm.mg.euler} \textbf{(b)} is proven.
\end{proof}

\bigskip

\textbf{Note:} If you look closely at the above proof, you will see hidden in
it an algorithm for \textbf{finding} Eulerian circuits and walks.\footnote{You
might be skeptical about this. After all, in order to apply Lemma
\ref{lem.mg.euler.longest-trail-closed}, we need a longest trail, so you might
wonder how we can find a longest trail to begin with.
\par
Fortunately, we don't need to take Lemma
\ref{lem.mg.euler.longest-trail-closed} this literally. Our above proof of
Lemma \ref{lem.mg.euler.longest-trail-closed} can be used even if $\mathbf{w}$
is \textbf{not} a longest trail. In this case, however, instead of showing
that $\mathbf{w}$ is a closed walk, this proof may show us a way how to make
$\mathbf{w}$ longer. In other words, by following this proof, we may discover
a trail longer than $\mathbf{w}$. In this case, we can replace $\mathbf{w}$ by
this longer trail, and then apply Lemma
\ref{lem.mg.euler.longest-trail-closed} again. We can repeat this over and
over again, until we do end up with a closed walk. (This will eventually
happen, since we know that a trail cannot be longer than the total number of
edges of $G$.)} \medskip

\begin{exercise}
\label{exe.mt1.euler-add} Let $G$ be a connected multigraph. Let $m$ be the
number of vertices of $G$ that have odd degree. Prove that we can add $m/2$
new edges to $G$ in such a way that the resulting multigraph will have an
Eulerian circuit. (It is allowed to add an edge even if there is already an
edge between the same two vertices.) \medskip

[\textbf{Solution:} This exercise is Exercise 6 on midterm \#1 from my Spring
2017 course; see \href{https://www.cip.ifi.lmu.de/~grinberg/t/17s/}{the course
page} for solutions.]
\end{exercise}

\begin{exercise}
\label{exe.mt1.L-hamil} Let $G=\left(  V,E,\varphi\right)  $ be a multigraph.
The \textbf{line graph} $L\left(  G\right)  $ is defined as the simple graph
$\left(  E,F\right)  $, where
\[
F=\left\{  \left\{  e_{1},e_{2}\right\}  \in\mathcal{P}_{2}\left(  {E}\right)
\ \mid\ \varphi\left(  e_{1}\right)  \cap\varphi\left(  e_{2}\right)
\neq\varnothing\right\}  .
\]
(In other words, $L\left(  G\right)  $ is the graph whose \textbf{vertices}
are the \textbf{edges} of $G$, and in which two vertices $e_{1}$ and $e_{2}$
are adjacent if and only if the edges $e_{1}$ and $e_{2}$ of $G$ share a
common endpoint.)

[\textbf{Example:} Here is a multigraph $G$ along with its line graph
$L\left(  G\right)  $:%
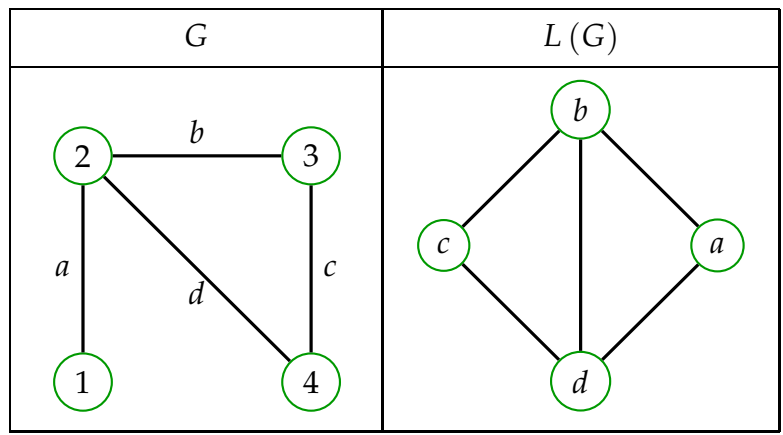
\[%
\begin{tabular}
[c]{|c|c|}\hline
$G$ & $L\left(  G\right)  $\\\hline
$\ \
\begin{tikzpicture}[scale=3]
\begin{scope}[every node/.style={circle,thick,draw=green!60!black}]
\node(1) at (0,0) {$1$};
\node(2) at (0,1) {$2$};
\node(3) at (1,1) {$3$};
\node(4) at (1,0) {$4$};
\end{scope}
\begin{scope}[every edge/.style={draw=black,very thick}, every loop/.style={}]
\path[-] (1) edge node[left] {$a$} (2);
\path[-] (2) edge node[above] {$b$} (3);
\path[-] (2) edge node[below] {$d$} (4);
\path[-] (3) edge node[right] {$c$} (4);
\end{scope}
\end{tikzpicture}%
\ \ $ & \ \ $%
\begin{tikzpicture}[scale=1.8]
\begin{scope}[every node/.style={circle,thick,draw=green!60!black}]
\node(a) at (0:1) {$a$};
\node(b) at (90:1) {$b$};
\node(c) at (180:1) {$c$};
\node(d) at (270:1) {$d$};
\end{scope}
\begin{scope}[every edge/.style={draw=black,very thick}, every loop/.style={}]
\path[-] (a) edge (b) (a) edge (d) (b) edge (d) (b) edge (c) (c) edge (d);
\end{scope}
\end{tikzpicture}%
\ \ $\\\hline
\end{tabular}
\ \ \ .
\]
Note that $L\left(  G\right)  $ does not always determine $G$ uniquely.]

Assume that $\left\vert V\right\vert >1$. Prove the following:

\begin{enumerate}
\item[\textbf{(a)}] If $G$ has a Hamiltonian path, then $L\left(  G\right)  $
has a Hamiltonian path.

\item[\textbf{(b)}] If $G$ has an Eulerian walk, then $L\left(  G\right)  $
has a Hamiltonian path.
\end{enumerate}

[\textbf{Solution:} This exercise is Exercise 2 on midterm \#1 from my Spring
2017 course (generalized from simple graphs to multigraphs); see
\href{https://www.cip.ifi.lmu.de/~grinberg/t/17s/}{the course page} for solutions.]
\end{exercise}

\bigskip

\section{Digraphs and multidigraphs}

\subsection{Definitions}

We have so far seen two concepts of graphs: simple graphs and multigraphs.

For all their differences, these two concepts have one thing in common: The
two endpoints of an edge are equal in rights. Thus, when defining walks, each
edge serves as a \textquotedblleft two-way road\textquotedblright. Hence, such
graphs are good at modelling symmetric relations between things.

We shall now introduce two analogous versions of \textquotedblleft
graphs\textquotedblright\ in which the edges have directions. These versions
are known as \textbf{directed graphs} (short: \textbf{digraphs}). In such
directed graphs, each edge will have a specified starting point (its
\textquotedblleft source\textquotedblright) and a specified ending point (its
\textquotedblleft target\textquotedblright). Correspondingly, we will draw
these edges as arrows, and we will only allow using them in one direction
(viz., from source to target) when we walk down the graph. Here are the
definitions in detail:

\begin{definition}
A \textbf{simple digraph} is a pair $\left(  V,A\right)  $, where $V$ is a
finite set, and where $A$ is a subset of $V\times V$.
\end{definition}

\begin{definition}
Let $D=\left(  V,A\right)  $ be a simple digraph.

\begin{enumerate}
\item[\textbf{(a)}] The set $V$ is called the \textbf{vertex set} of $D$; it
is denoted by $\operatorname*{V}\left(  D\right)  $.

Its elements are called the \textbf{vertices} (or \textbf{nodes}) of $D$.

\item[\textbf{(b)}] The set $A$ is called the \textbf{arc set} of $D$; it is
denoted by $\operatorname*{A}\left(  D\right)  $.

Its elements are called the \textbf{arcs} (or \textbf{directed edges}) of $D$.

When $u$ and $v$ are two elements of $V$, we will occasionally use $uv$ as a
shorthand for the pair $\left(  u,v\right)  $. Note that this means an ordered
pair now!

\item[\textbf{(c)}] If $\left(  u,v\right)  $ is an arc of $D$ (or, more
generally, a pair in $V\times V$), then $u$ is called the \textbf{source} of
this arc, and $v$ is called the \textbf{target} of this arc.

\item[\textbf{(d)}] We draw $D$ as follows: We represent each vertex of $D$ by
a point, and each arc $uv$ by an arrow that goes from the point representing
$u$ to the point representing $v$.

\item[\textbf{(e)}] An arc $\left(  u,v\right)  $ is called a \textbf{loop}
(or \textbf{self-loop}) if $u=v$. (In other words, an arc is a loop if and
only if its source is its target.)
\end{enumerate}
\end{definition}

\begin{example}
For each $n\in\mathbb{N}$, we define the \textbf{divisibility digraph on
}$\left\{  1,2,\ldots,n\right\}  $ to be the simple digraph $\left(
V,A\right)  $, where $V=\left\{  1,2,\ldots,n\right\}  $ and%
\[
A=\left\{  \left(  i,j\right)  \in V\times V\ \mid\ i\text{ divides
}j\right\}  .
\]

For example, for $n=6$, this digraph looks as follows:%
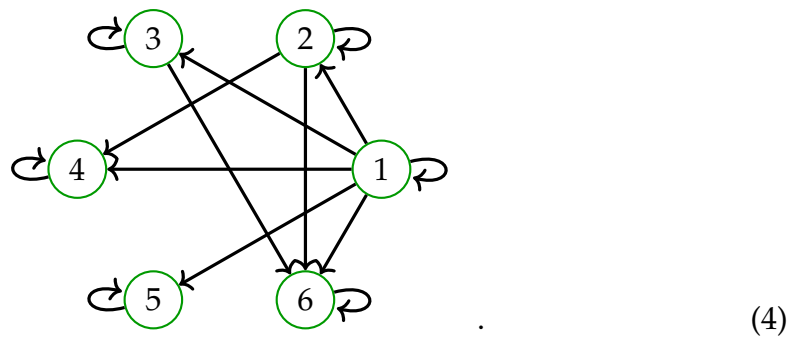
\begin{equation}%
\begin{tikzpicture}
\begin{scope}[every node/.style={circle,thick,draw=green!60!black}]
\node(1) at (0:2) {$1$};
\node(2) at (60:2) {$2$};
\node(3) at (120:2) {$3$};
\node(4) at (180:2) {$4$};
\node(5) at (240:2) {$5$};
\node(6) at (300:2) {$6$};
\end{scope}
\begin{scope}[every edge/.style={draw=black,very thick}]
\path[->] (1) edge (2) edge (3) edge (4) edge (5) edge (6);
\path[->] (2) edge (4) edge (6);
\path[->] (3) edge (6);
\path[->] (1) edge[loop right] (1);
\path[->] (2) edge[loop right] (2);
\path[->] (3) edge[loop left] (3);
\path[->] (4) edge[loop left] (4);
\path[->] (5) edge[loop left] (5);
\path[->] (6) edge[loop right] (6);
\end{scope}
\end{tikzpicture}%
\ \ . \label{eq.exa.digr.divisible6}%
\end{equation}
Note that the divisibility digraph on $\left\{  1,2,\ldots,n\right\}  $ has
$n$ loops $\left(  1,1\right)  ,\ \left(  2,2\right)  ,\ \ldots,\ \left(
n,n\right)  $.
\end{example}

Note that simple digraphs (unlike simple graphs) are allowed to have loops
(i.e., arcs of the form $\left(  v,v\right)  $).

\begin{definition}
A \textbf{multidigraph} is a triple $\left(  V,A,\psi\right)  $, where $V$ and
$A$ are two finite sets, and $\psi:A\rightarrow V\times V$ is a map.
\end{definition}

\begin{definition}
Let $D=\left(  V,A,\psi\right)  $ be a multidigraph.

\begin{enumerate}
\item[\textbf{(a)}] The set $V$ is called the \textbf{vertex set} of $D$; it
is denoted by $\operatorname*{V}\left(  D\right)  $.

Its elements are called the \textbf{vertices} (or \textbf{nodes}) of $D$.

\item[\textbf{(b)}] The set $A$ is called the \textbf{arc set} of $D$; it is
denoted by $\operatorname*{A}\left(  D\right)  $.

Its elements are called the \textbf{arcs} (or \textbf{directed edges}) of $D$.

\item[\textbf{(c)}] If $a$ is an arc of $D$, and if $\psi\left(  a\right)
=\left(  u,v\right)  $, then the vertex $u$ is called the \textbf{source} of
$a$, and the vertex $v$ is called the \textbf{target} of $a$.

\item[\textbf{(d)}] We draw $D$ as follows: We represent each vertex of $D$ by
a point, and each arc $a$ by an arrow that goes from the point representing
$u$ to the point representing $v$, where $\left(  u,v\right)  =\psi\left(
a\right)  $.

\item[\textbf{(e)}] An arc $a$ of $D$ is called a \textbf{loop} (or
\textbf{self-loop}) if its source is its target (i.e., if it satisfies $u=v$,
where $\left(  u,v\right)  =\psi\left(  a\right)  $).
\end{enumerate}
\end{definition}

\begin{example}
Here is a multidigraph:%
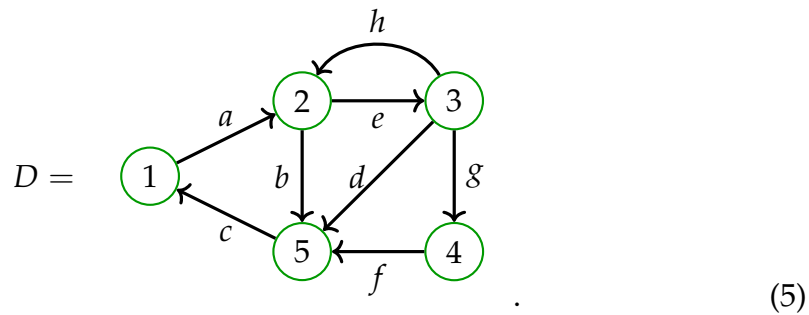
\begin{equation}%
\begin{tikzpicture}[scale=2]
\begin{scope}[every node/.style={circle,thick,draw=green!60!black}]
\node(1) at (-1,0) {$1$};
\node(2) at (0,0) {$2$};
\node(3) at (1,0.5) {$3$};
\node(4) at (1,-0.5) {$4$};
\end{scope}
\node(X) at (-1.7, 0) {$D = $};
\begin{scope}[every edge/.style={draw=black,very thick}, every loop/.style={}]
\path[->] (1) edge[bend left=60] node[above] {$a$} (2);
\path[->] (1) edge[bend right=60] node[below] {$b$} (2);
\path[->] (2) edge node[above] {$c$} (3);
\path[->] (3) edge node[right] {$d$} (4);
\path[->] (4) edge node[below] {$e$} (2);
\end{scope}
\end{tikzpicture}%
\ \ . \label{eq.exa.digr.mdg.2.eq}%
\end{equation}
Formally speaking, this multidigraph is the triple $\left(  V,A,\psi\right)
$, where $V=\left\{  1,2,3,4\right\}  $ and $A=\left\{  a,b,c,d,e\right\}  $
and $\psi\left(  a\right)  =\left(  1,2\right)  $ and $\psi\left(  b\right)
=\left(  1,2\right)  $ and $\psi\left(  c\right)  =\left(  2,3\right)  $ and
$\psi\left(  d\right)  =\left(  3,4\right)  $ and $\psi\left(  e\right)
=\left(  4,2\right)  $.
\end{example}

\begin{example}
Here is another multidigraph:%
\begin{equation}%
\begin{tikzpicture}[scale=2]
\begin{scope}[every node/.style={circle,thick,draw=green!60!black}]
\node(1) at (-1,0) {$1$};
\node(2) at (0,0.5) {$2$};
\node(3) at (1,0.5) {$3$};
\node(4) at (1,-0.5) {$4$};
\node(5) at (0,-0.5) {$5$};
\end{scope}
\node(X) at (-1.7, 0) {$D = $};
\begin{scope}[every edge/.style={draw=black,very thick}, every loop/.style={}]
\path[->] (1) edge node[above] {$a$} (2);
\path[<-] (1) edge node[below] {$c$} (5);
\path[->] (2) edge node[left] {$b$} (5);
\path[->] (2) edge node[below] {$e$} (3);
\path[<-] (2) edge[bend left=60] node[above] {$h$} (3);
\path[->] (3) edge node[left] {$d$} (5);
\path[->] (4) edge node[below] {$f$} (5);
\path[->] (3) edge node[right] {$g$} (4);
\end{scope}
\end{tikzpicture}%
\ \ . \label{eq.exa.digr.mdg.1.eq}%
\end{equation}
Formally speaking, this multidigraph is the triple $\left(  V,A,\psi\right)
$, where $V=\left\{  1,2,3,4,5\right\}  $ and $A=\left\{
a,b,c,d,e,f,g,h\right\}  $ and $\psi\left(  a\right)  =\left(  1,2\right)  $
and $\psi\left(  b\right)  =\left(  2,5\right)  $ and so on.
\end{example}

Thus, simple digraphs and multidigraphs are analogues of simple graphs and
multigraphs, respectively, in which the edges have been replaced by arcs
(\textquotedblleft edges endowed with a direction\textquotedblright). The
analogy is perfect but for the fact that simple graphs forbid loops but simple
digraphs allow loops (but different authors have different opinions on this).

\begin{convention}
The word \textquotedblleft\textbf{digraph}\textquotedblright\ means either
\textquotedblleft simple digraph\textquotedblright\ or \textquotedblleft
multidigraph\textquotedblright, depending on the context.
\end{convention}

The word \textquotedblleft digraph\textquotedblright\ was originally a
shorthand for \textquotedblleft\textbf{directed graph}\textquotedblright, but
by now it is a technical term that is perfectly understood by everyone in the
subject. (It is also understood by linguists, but in a rather different way.)

\begin{definition}
A digraph $D$ is said to be \textbf{loopless} if it has no loops.
\end{definition}

\subsection{Outdegrees and indegrees}

What can we do with digraphs? Many of the things we have done with graphs can
be modified to work with digraphs (although not all their properties will
still hold). For example, the notion of the degree of a vertex in a graph has
the following two counterpart notions for digraphs:

\begin{definition}
Let $D$ be a digraph with vertex set $V$. (This can be either a simple digraph
or a multidigraph.) Let $v\in V$ be any vertex. Then:

\begin{enumerate}
\item[\textbf{(a)}] The \textbf{outdegree} of $v$ denotes the number of arcs
of $D$ whose source is $v$. This outdegree is denoted $\deg^{+}v$.

\item[\textbf{(b)}] The \textbf{indegree} of $v$ denotes the number of arcs of
$D$ whose target is $v$. This indegree is denoted $\deg^{-}v$.
\end{enumerate}
\end{definition}

\begin{example}
In the divisibility digraph on $\left\{  1,2,3,4,5,6\right\}  $ (see
(\ref{eq.exa.digr.divisible6}) for a drawing), we have%
\begin{align*}
\deg^{+}1  &  =6,\ \ \ \ \ \ \ \ \ \ \deg^{-}1=1,\ \ \ \ \ \ \ \ \ \ \deg
^{+}2=3,\ \ \ \ \ \ \ \ \ \ \deg^{-}2=2,\\
\deg^{+}3  &  =2,\ \ \ \ \ \ \ \ \ \ \deg^{-}3=2,\ \ \ \ \ \ \ \ \ \ \deg
^{+}4=1,\ \ \ \ \ \ \ \ \ \ \deg^{-}4=3,\\
\deg^{+}5  &  =1,\ \ \ \ \ \ \ \ \ \ \deg^{-}5=2,\ \ \ \ \ \ \ \ \ \ \deg
^{+}6=1,\ \ \ \ \ \ \ \ \ \ \deg^{-}6=4.
\end{align*}

\end{example}

Recall Euler's result (Proposition \ref{prop.mg.sum-deg}) saying that in a
graph, the sum of all degrees is twice the number of edges. Here is an
analogue of this result for digraphs:

\begin{proposition}
[diEuler]\label{prop.mdg.sum-deg}Let $D$ be a digraph with vertex set $V$ and
arc set $A$. Then,%
\[
\sum_{v\in V}\deg^{+}v=\sum_{v\in V}\deg^{-}v=\left\vert A\right\vert .
\]

\end{proposition}

\begin{proof}
By the definition of an outdegree, we have%
\[
\deg^{+}v=\left(  \text{the number of arcs of }D\text{ whose source is
}v\right)
\]
for each $v\in V$. Thus,%
\begin{align*}
\sum_{v\in V}\deg^{+}v  &  =\sum_{v\in V}\left(  \text{the number of arcs of
}D\text{ whose source is }v\right) \\
&  =\left(  \text{the number of all arcs of }D\right) \\
&  \ \ \ \ \ \ \ \ \ \ \ \ \ \ \ \ \ \ \ \ \left(
\begin{array}
[c]{c}%
\text{since each arc of }D\text{ has exactly one source,}\\
\text{and thus is counted exactly once in the sum}%
\end{array}
\right) \\
&  =\left\vert A\right\vert .
\end{align*}
Similarly, $\sum_{v\in V}\deg^{-}v=\left\vert A\right\vert $.
\end{proof}

(\textquotedblleft diEuler\textquotedblright\ is not a real mathematician; I
just gave that moniker to Proposition \ref{prop.mdg.sum-deg} in order to
stress its analogy with Euler's 1736 result.)

\subsection{Subdigraphs}

Just as we defined subgraphs of a multigraph, we can define subdigraphs (or
\textquotedblleft submultidigraphs\textquotedblright, to be very precise) of a digraph:

\begin{definition}
\label{def.mdg.submdg}Let $D=\left(  V,A,\psi\right)  $ be a multidigraph.

\begin{enumerate}
\item[\textbf{(a)}] A \textbf{submultidigraph} (or, for short,
\textbf{subdigraph}) of $D$ means a multidigraph of the form $E=\left(
W,B,\chi\right)  $, where $W\subseteq V$ and $B\subseteq A$ and $\chi=\psi
\mid_{B}$. In other words, a submultidigraph of $D$ means a multidigraph $E$
whose vertices are vertices of $D$ and whose arcs are arcs of $D$ and whose
arcs have the same sources and targets in $E$ as they have in $D$.

\item[\textbf{(b)}] Let $S$ be a subset of $V$. The \textbf{induced subdigraph
of }$D$\textbf{ on the set }$S$ denotes the subdigraph%
\[
\left(  S,\ \ A^{\prime},\ \ \psi\mid_{A^{\prime}}\right)
\]
of $D$, where%
\[
A^{\prime}:=\left\{  a\in A\ \mid\ \text{both the source and the target of
}a\text{ belong to }S\right\}  .
\]
In other words, it denotes the subdigraph of $D$ whose vertices are the
elements of $S$, and whose arcs are precisely those arcs of $D$ whose sources
and targets both belong to $S$. We denote this induced subdigraph by $D\left[
S\right]  $.

\item[\textbf{(c)}] An \textbf{induced subdigraph} of $D$ means a subdigraph
of $D$ that is the induced subdigraph of $D$ on $S$ for some $S\subseteq V$.
\end{enumerate}
\end{definition}

\subsection{Conversions}

\subsubsection{Multidigraphs to multigraphs}

Any multidigraph $D$ can be turned into an (undirected) multigraph $G$ by
\textquotedblleft removing the arrowheads\textquotedblright\ (aka
\textquotedblleft forgetting the directions of the arcs\textquotedblright):

\begin{definition}
\label{def.mdg.Dund} Let $D$ be a multidigraph. Then, $D^{\operatorname*{und}%
}$ will denote the multigraph obtained from $D$ by replacing each arc with an
edge whose endpoints are the source and the target of this arc. Formally, this
is defined as follows: If $D=\left(  V,A,\psi\right)  $, then
$D^{\operatorname*{und}}=\left(  V,A,\varphi\right)  $, where the map
$\varphi:A\rightarrow\mathcal{P}_{1,2}\left(  V\right)  $ sends each arc $a\in
A$ to the set of the entries of $\psi\left(  a\right)  $ (that is, to the set
consisting of the source of $a$ and the target of $a$).

We call $D^{\operatorname*{und}}$ the \textbf{underlying undirected graph} of
$D$.
\end{definition}

For example, if $D$ is the multidigraph from (\ref{eq.exa.digr.mdg.1.eq}),
then $D^{\operatorname*{und}}$ is the following multigraph:%
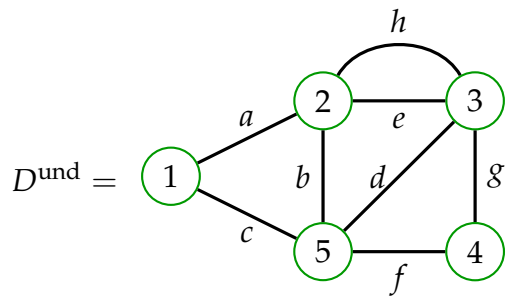
\[%
\begin{tikzpicture}[scale=2]
\begin{scope}[every node/.style={circle,thick,draw=green!60!black}]
\node(1) at (-1,0) {$1$};
\node(2) at (0,0.5) {$2$};
\node(3) at (1,0.5) {$3$};
\node(4) at (1,-0.5) {$4$};
\node(5) at (0,-0.5) {$5$};
\end{scope}
\node(X) at (-1.7, 0) {$D^{\operatorname{und}} = $};
\begin{scope}[every edge/.style={draw=black,very thick}, every loop/.style={}]
\path[-] (1) edge node[above] {$a$} (2);
\path[-] (1) edge node[below] {$c$} (5);
\path[-] (2) edge node[left] {$b$} (5);
\path[-] (2) edge node[below] {$e$} (3);
\path[-] (2) edge[bend left=60] node[above] {$h$} (3);
\path[-] (3) edge node[left] {$d$} (5);
\path[-] (4) edge node[below] {$f$} (5);
\path[-] (3) edge node[right] {$g$} (4);
\end{scope}
\end{tikzpicture}%
\ \ .
\]

\subsubsection{Multigraphs to multidigraphs}

We have just seen how to turn any multidigraph $D$ into a multigraph
$D^{\operatorname*{und}}$ by forgetting the directions of the arcs.

Conversely, we can turn a multigraph $G$ into a multidigraph
$G^{\operatorname*{bidir}}$ by \textquotedblleft duplicating\textquotedblright%
\ each edge (more precisely: turning each edge into two arcs with opposite
orientations). Here is a formal definition:

\begin{definition}
\label{def.mg.bidir} Let $G=\left(  V,E,\varphi\right)  $ be a multigraph. For
each edge $e\in E$, let us choose one of the endpoints of $e$ and call it
$s_{e}$; the other endpoint will then be called $t_{e}$. (If $e$ is a loop,
then we understand $t_{e}$ to mean $s_{e}$.)

We then define $G^{\operatorname*{bidir}}$ to be the multidigraph $\left(
V,\ \ E\times\left\{  1,2\right\}  ,\ \ \psi\right)  $, where the map
$\psi:E\times\left\{  1,2\right\}  \rightarrow V\times V$ is defined as
follows: For each edge $e\in E$, we set
\[
\psi\left(  e,1\right)  =\left(  s_{e},t_{e}\right)
\ \ \ \ \ \ \ \ \ \ \text{and}\ \ \ \ \ \ \ \ \ \ \psi\left(  e,2\right)
=\left(  t_{e},s_{e}\right)  .
\]
We call $G^{\operatorname*{bidir}}$ the \textbf{bidirectionalized
multidigraph} of $G$.
\end{definition}

Note that the map $\psi$ depends on our choice of $s_{e}$'s (that is, it
depends on which endpoint of an edge $e$ we choose to be $s_{e}$). This makes
the definition of $G^{\operatorname*{bidir}}$ non-canonical\footnote{If you
want to eradicate this dependency, you can rename the arcs $\left(
e,1\right)  $ and $\left(  e,2\right)  $ of $G^{\operatorname*{bidir}}$ as
$\left(  e,s_{e}\right)  $ and $\left(  e,t_{e}\right)  $, except in the case
when $e$ is a loop, in which case they remain $\left(  e,1\right)  $ and
$\left(  e,2\right)  $. This way, $\psi$ no longer depends on any choices. But
the independence might not be worth the awkwardness of having loops and
non-loops treated differently.}. Fortunately, all choices of $s_{e}$'s will
lead to mutually isomorphic multidigraphs $G^{\operatorname*{bidir}}$. (The
notion of \textbf{isomorphism} for multidigraphs is exactly the one that you expect.)

\begin{example}
If%
\[%
%
\ \ .
\]
(Here, for example, we have chosen $s_{a}$ to be $2$, so that $t_{a}=3$ and
$\psi\left(  a,1\right)  =\left(  2,3\right)  $ and $\psi\left(  a,2\right)
=\left(  3,2\right)  $.) Yes, even the loops of $G$ are duplicated in
$G^{\operatorname*{bidir}}$ !
\end{example}

The operation that assigns a multidigraph $G^{\operatorname*{bidir}}$ to a
multigraph $G$ is injective -- i.e., the original graph $G$ can be uniquely
reconstructed from $G^{\operatorname*{bidir}}$. This is in stark difference to
the operation $D\mapsto D^{\operatorname*{und}}$, which destroys information
(the directions of the arcs). Note that the multigraph $\left(
G^{\operatorname*{bidir}}\right)  ^{\operatorname*{und}}$ is not isomorphic to
$G$, since each edge of $G$ is doubled in $\left(  G^{\operatorname*{bidir}%
}\right)  ^{\operatorname*{und}}$.

\subsubsection{Simple digraphs to multidigraphs}

Next, we introduce another operation: one that turns simple digraphs into
multidigraphs. This is very similar to the operation $G\mapsto
G^{\operatorname*{mult}}$ that turns simple graphs into multigraphs, so we
will even use the same notation for it. Its definition is as follows:

\begin{definition}
\label{def.sdg.Dmult}Let $D=\left(  V,A\right)  $ be a simple digraph. Then,
the \textbf{corresponding multidigraph} $D^{\operatorname*{mult}}$ is defined
to be the multidigraph%
\[
\left(  V,A,\iota\right)  ,
\]
where $\iota:A\rightarrow V\times V$ is the map sending each $a\in A$ to $a$ itself.
\end{definition}

\begin{example}
If%
\[%
%
\ \ .
\]

\end{example}

\subsubsection{Multidigraphs to simple digraphs}

There is also an operation $D\mapsto D^{\operatorname*{simp}}$ that turns
multidigraphs into simple digraphs:\footnote{I will use a notation that I
probably should have introduced before: If $u$ and $v$ are two vertices of a
digraph, then an \textquotedblleft\textbf{arc from }$u$ \textbf{to }%
$v$\textquotedblright\ means an arc with source $u$ and target $v$.}

\begin{definition}
Let $D=\left(  V,A,\psi\right)  $ be a multidigraph. Then, the
\textbf{underlying simple digraph} $D^{\operatorname*{simp}}$ of $D$ means the
simple digraph%
\[
\left(  V,\ \left\{  \psi\left(  a\right)  \ \mid\ a\in A\right\}  \right)  .
\]
In other words, it is the simple digraph with vertex set $V$ in which there is
an arc from $u$ to $v$ if there exists an arc from $u$ to $v$ in $D$. Thus,
$D^{\operatorname*{simp}}$ is obtained from $D$ by \textquotedblleft
collapsing\textquotedblright\ parallel arcs (i.e., arcs having the same source
and the same target) to a single arc.
\end{definition}

\begin{example}
If%
\[%
%
\ \ .
\]
Note that the arcs $c$ and $d$ have not been \textquotedblleft
collapsed\textquotedblright\ into one arc, since they do not have the same
source and the same target. Likewise, the loop $g$ has been preserved (unlike
for undirected graphs).
\end{example}

\subsubsection{Multidigraphs as a big tent}

We have now established (in this section and also in Section \ref{sec.mg.conv}%
) several operations for converting between different types of graphs. Let us
bring them together in a single picture:\footnote{The enthusiastic reader can
also define conversions between simple graphs and simple digraphs, which will
complete this picture to a \textquotedblleft square\textquotedblright; but we
shall not have a use for such conversions.}%
\[%
\begin{tikzpicture}[scale=6]
\begin{scope}[every node/.style={circle,thick,draw=green!60!black}]
\node(sg) at (0,0) {simple graphs};
\node(mg) at (1,0) {multigraphs};
\node(sdg) at (0,-1) {simple digraphs};
\node(mdg) at (1,-1) {multidigraphs};
\end{scope}
\draw[->, very thick] (sg) edge[bend left=20] node[above] {$G \mapsto
G^{\operatorname{mult}}$} (mg);
\draw[->, very thick] (mg) edge[bend left=20] node[below] {$G \mapsto
G^{\operatorname{simp}}$} (sg);
\draw[->, very thick] (sdg) edge[bend left=20] node[above] {$D \mapsto
D^{\operatorname{mult}}$} (mdg);
\draw[->, very thick] (mdg) edge[bend left=20] node[below] {$D \mapsto
D^{\operatorname{simp}}$} (sdg);
\draw[->, very thick] (mg) edge[bend left=20] node[right] {$G \mapsto
G^{\operatorname{bidir}}$} (mdg);
\draw[->, very thick] (mdg) edge[bend left=20] node[left] {$D \mapsto
D^{\operatorname{und}}$} (mg);
\end{tikzpicture}%
\ \ .
\]

A takeaway from this all is that multidigraphs are the \textquotedblleft most
general\textquotedblright\ notion of graphs we have introduced so far. Indeed,
using the operations we have seen so far, we can convert every notion of
graphs into a multidigraph:

\begin{itemize}
\item Each simple graph becomes a multigraph via the $G\mapsto
G^{\operatorname*{mult}}$ operation.

\item Each multigraph, in turn, becomes a multidigraph via the $G\mapsto
G^{\operatorname*{bidir}}$ operation.

\item Each simple digraph becomes a multidigraph via the $D\mapsto
D^{\operatorname*{mult}}$ operation.
\end{itemize}

Since all three of these operations are injective (i.e., lose no information),
we thus can encode each of our four notions of graphs as a multidigraph.
Consequently, any theorem about multidigraphs can be specialized to the other
three types of graphs. This doesn't mean that any theorem on any other type of
graphs can be generalized to multidigraphs, though (e.g., Mantel's theorem
holds only for simple graphs) -- but when it can, we will try to state it at
the most general level possible, to avoid doing the same work twice.

\subsection{\label{sec.dg.walks}Walks, paths, closed walks, cycles}

\subsubsection{Definitions}

Let us now define various kinds of walks for simple digraphs and for multidigraphs.

For simple digraphs, we imitate the definitions from Sections
\ref{sec.sg.walks} and \ref{sec.sg.cycs} as best as we can, making sure to
require all arcs to be traversed in the correct direction:

\begin{definition}
Let $D$ be a simple digraph. Then:

\begin{enumerate}
\item[\textbf{(a)}] A \textbf{walk} (in $D$) means a finite sequence $\left(
v_{0},v_{1},\ldots,v_{k}\right)  $ of vertices of $D$ (with $k\geq0$) such
that all of the pairs $v_{0}v_{1},\ v_{1}v_{2},\ v_{2}v_{3},\ \ldots
,\ v_{k-1}v_{k}$ are arcs of $D$. (The latter condition is vacuously true if
$k=0$.)

\item[\textbf{(b)}] If $\mathbf{w}=\left(  v_{0},v_{1},\ldots,v_{k}\right)  $
is a walk in $D$, then:

\begin{itemize}
\item The \textbf{vertices} of $\mathbf{w}$ are defined to be $v_{0}%
,v_{1},\ldots,v_{k}$.

\item The \textbf{arcs} of $\mathbf{w}$ are defined to be the pairs
$v_{0}v_{1},\ v_{1}v_{2},\ v_{2}v_{3},\ \ldots,\ v_{k-1}v_{k}$.

\item The nonnegative integer $k$ is called the \textbf{length} of
$\mathbf{w}$. (This is the number of all arcs of $\mathbf{w}$, counted with
multiplicity. It is $1$ smaller than the number of all vertices of
$\mathbf{w}$, counted with multiplicity.)

\item The vertex $v_{0}$ is called the \textbf{starting point} of $\mathbf{w}%
$. We say that $\mathbf{w}$ \textbf{starts} (or \textbf{begins}) at $v_{0}$.

\item The vertex $v_{k}$ is called the \textbf{ending point} of $\mathbf{w}$.
We say that $\mathbf{w}$ \textbf{ends} at $v_{k}$.
\end{itemize}

\item[\textbf{(c)}] A \textbf{path} (in $D$) means a walk (in $D$) whose
vertices are distinct. In other words, a path means a walk $\left(
v_{0},v_{1},\ldots,v_{k}\right)  $ such that $v_{0},v_{1},\ldots,v_{k}$ are distinct.

\item[\textbf{(d)}] Let $p$ and $q$ be two vertices of $D$. A \textbf{walk
from }$p$ \textbf{to }$q$ means a walk that starts at $p$ and ends at $q$. A
\textbf{path from }$p$ \textbf{to }$q$ means a path that starts at $p$ and
ends at $q$.

\item[\textbf{(e)}] A \textbf{closed walk} of $D$ means a walk whose first
vertex is identical with its last vertex. In other words, it means a walk
$\left(  w_{0},w_{1},\ldots,w_{k}\right)  $ with $w_{0}=w_{k}$. Sometimes,
closed walks are also known as \textbf{circuits} (but many authors use this
latter word for something slightly different).

\item[\textbf{(f)}] A \textbf{cycle} of $D$ means a closed walk $\left(
w_{0},w_{1},\ldots,w_{k}\right)  $ such that $k\geq1$ and such that the
vertices $w_{0},w_{1},\ldots,w_{k-1}$ are distinct.
\end{enumerate}
\end{definition}

Note that we replaced the condition $k\geq3$ by $k\geq1$ in the definition of
a cycle, since simple digraphs can have loops. Fortunately, with the arcs
being directed, we no longer have to worry about the same arc being traversed
back and forth, so we need no extra condition to rule this out.

\begin{example}
\label{exa.digr.sdg.walks1}Consider the simple digraph%
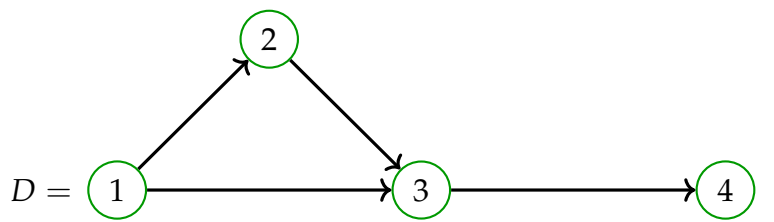
\[%
%
\ \ .
\]
Then, $\left(  1,2,3,4\right)  $ and $\left(  1,3,4\right)  $ are two walks of
$D$, and these walks are paths. But $\left(  2,3,1\right)  $ is not a walk
(since you cannot use the arc $13$ to get from $3$ to $1$). This digraph $D$
has no cycles, and its only closed walks have length $0$.
\end{example}

\begin{example}
\label{exa.digr.sdg.walks2}Consider the simple digraph%
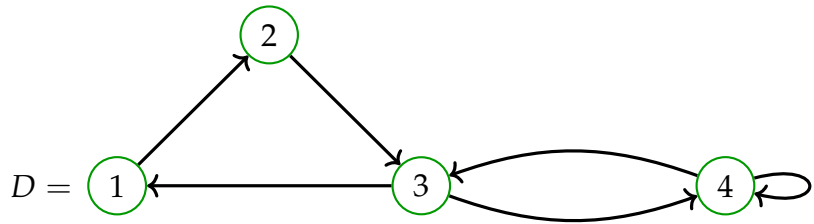
\[%
%
\ \ .
\]
Then, $\left(  1,2,3,1\right)  $ and $\left(  3,4,3\right)  $ and $\left(
4,4\right)  $ are cycles of $D$. Moreover, $\left(  1,2,3,4,3,1\right)  $ is a
closed walk but not a cycle.
\end{example}

Now let's define the same concepts for multidigraphs, by modifying the
analogous definitions for multigraphs we saw in Definition
\ref{def.mg.walks-etc}:

\begin{definition}
\label{def.mdg.walks-etc}Let $D=\left(  V,A,\psi\right)  $ be a multidigraph. Then:

\begin{enumerate}
\item[\textbf{(a)}] A \textbf{walk} in $D$ means a list of the form%
\[
\left(  v_{0},a_{1},v_{1},a_{2},v_{2},\ldots,a_{k},v_{k}\right)
\ \ \ \ \ \ \ \ \ \ \left(  \text{with }k\geq0\right)  ,
\]
where $v_{0},v_{1},\ldots,v_{k}$ are vertices of $D$, where $a_{1}%
,a_{2},\ldots,a_{k}$ are arcs of $D$, and where each $i\in\left\{
1,2,\ldots,k\right\}  $ satisfies%
\[
\psi\left(  a_{i}\right)  =\left(  v_{i-1},v_{i}\right)
\]
(that is, each arc $a_{i}$ has source $v_{i-1}$ and target $v_{i}$). Note that
we have to record both the vertices \textbf{and} the arcs in our walk, since
we want the walk to \textquotedblleft know\textquotedblright\ which arcs it traverses.

The \textbf{vertices} of a walk $\left(  v_{0},a_{1},v_{1},a_{2},v_{2}%
,\ldots,a_{k},v_{k}\right)  $ are $v_{0},v_{1},\ldots,v_{k}$; the
\textbf{arcs} of this walk are $a_{1},a_{2},\ldots,a_{k}$. This walk is said
to \textbf{start} at $v_{0}$ and \textbf{end} at $v_{k}$; it is also said to
be a \textbf{walk from }$v_{0}$ \textbf{to }$v_{k}$. Its \textbf{starting
point} is $v_{0}$, and its \textbf{ending point} is $v_{k}$. Its
\textbf{length} is $k$.

\item[\textbf{(b)}] A \textbf{path} means a walk whose vertices are distinct.

\item[\textbf{(c)}] A \textbf{closed walk} (or \textbf{circuit}) means a walk
$\left(  v_{0},a_{1},v_{1},a_{2},v_{2},\ldots,a_{k},v_{k}\right)  $ with
$v_{k}=v_{0}$.

\item[\textbf{(d)}] A \textbf{cycle} means a closed walk $\left(  v_{0}%
,a_{1},v_{1},a_{2},v_{2},\ldots,a_{k},v_{k}\right)  $ such that

\begin{itemize}
\item the vertices $v_{0},v_{1},\ldots,v_{k-1}$ are distinct;

\item we have $k\geq1$.
\end{itemize}

(This automatically implies that the arcs $a_{1},a_{2},\ldots,a_{k}$ are
distinct, since each arc $a_{i}$ has source $v_{i-1}$.)
\end{enumerate}
\end{definition}

\begin{example}
\label{exa.digr.mdg.cycles1}Consider the multidigraph%
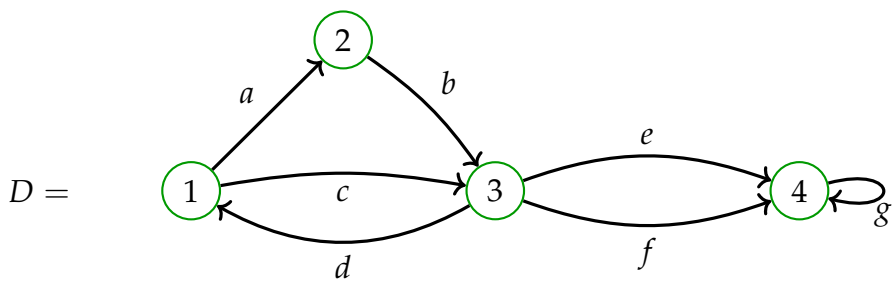
\[%
\begin{tikzpicture}[scale=2]
\begin{scope}[every node/.style={circle,thick,draw=green!60!black}]
\node(1) at (-1,0) {$1$};
\node(2) at (0,1) {$2$};
\node(3) at (1,0) {$3$};
\node(4) at (3,0) {$4$};
\end{scope}
\node(X) at (-2, 0) {$D=$};
\begin{scope}[every edge/.style={draw=black,very thick}, every loop/.style={}]
\path[->] (1) edge node[above left] {$a$} (2);
\path[->] (2) edge[bend left=10] node[above right] {$b$} (3);
\path[->] (1) edge[bend left=10] node[below] {$c$} (3);
\path[->] (3) edge[bend left=30] node[below] {$d$} (1);
\path[->] (3) edge[bend left=20] node[above] {$e$} (4);
\path[->] (3) edge[bend right=20] node[below] {$f$} (4);
\path[->] (4) edge[loop right] node[below] {$g$} (4);
\end{scope}
\end{tikzpicture}%
\ \ .
\]
Then, $\left(  1,a,2,b,3,d,1\right)  $ and $\left(  3,d,1,c,3\right)  $ and
$\left(  4,g,4\right)  $ are three cycles of $D$, whereas $\left(
3,d,1,a,2,b,3,d,1,c,3\right)  $ is a circuit but not a cycle.
\end{example}

\subsubsection{Basic properties}

Now, let us see which properties of walks, paths, closed walks and cycles
remain valid for digraphs.

In Proposition \ref{prop.sg.walk-concat}, we saw how two walks in a simple
graph could be combined (\textquotedblleft spliced together\textquotedblright)
if the ending point of the first is the starting point of the second. In
Proposition \ref{prop.mg.walk-concat}, we generalized this to multigraphs. The
same holds for multidigraphs:

\begin{proposition}
\label{prop.mdg.walk-concat}Let $D$ be a multidigraph. Let $u$, $v$ and $w$ be
three vertices of $D$. Let $\mathbf{a}=\left(  a_{0},e_{1},a_{1},\ldots
,e_{k},a_{k}\right)  $ be a walk from $u$ to $v$. Let $\mathbf{b}=\left(
b_{0},f_{1},b_{1},\ldots,f_{\ell},b_{\ell}\right)  $ be a walk from $v$ to
$w$. Then,%
\begin{align*}
&  \left(  a_{0},e_{1},a_{1},\ldots,e_{k},a_{k},f_{1},b_{1},f_{2},b_{2}%
,\ldots,f_{\ell},b_{\ell}\right) \\
&  =\left(  a_{0},e_{1},a_{1},\ldots,a_{k-1},e_{k},b_{0},f_{1},b_{1}%
,\ldots,f_{\ell},b_{\ell}\right) \\
&  =\left(  a_{0},e_{1},a_{1},\ldots,a_{k-1},e_{k},v,f_{1},b_{1}%
,\ldots,f_{\ell},b_{\ell}\right)
\end{align*}
is a walk from $u$ to $w$. This walk shall be denoted $\mathbf{a}%
\ast\mathbf{b}$.
\end{proposition}

\begin{proof}
The same (trivial) argument as for undirected graphs works here.
\end{proof}

However, unlike for undirected graphs, we can no longer reverse walks or paths
in digraphs. Thus, it often happens that there is a walk from $u$ to $v$, but
no walk from $v$ to $u$. \medskip

Reducing a walk to a path (as we did in Proposition
\ref{prop.sg.walk-to-path-1} for simple graphs and in Proposition
\ref{prop.mg.walk-to-path-1} for multigraphs) still works for multidigraphs:

\begin{proposition}
\label{prop.mdg.walk-to-path-1}Let $D$ be a multidigraph. Let $u$ and $v$ be
two vertices of $D$. Let $\mathbf{a}$ be a walk from $u$ to $v$. Let $k$ be
the length of $\mathbf{a}$. Assume that $\mathbf{a}$ is not a path. Then,
there exists a walk from $u$ to $v$ whose length is smaller than $k$.
\end{proposition}

\begin{corollary}
[When there is a walk, there is a path]\label{cor.mdg.walk-thus-path}Let $D$
be a multidigraph. Let $u$ and $v$ be two vertices of $D$. Assume that there
is a walk from $u$ to $v$ of length $k$ for some $k\in\mathbb{N}$. Then, there
is a path from $u$ to $v$ of length $\leq k$.
\end{corollary}

The proofs of these facts are the same as for multigraphs (see Proposition
\ref{prop.mg.walk-to-path-1} and Corollary \ref{cor.mg.walk-thus-path},
respectively). \medskip

The following proposition is an analogue of Proposition
\ref{prop.cyc.btf-walk-cyc} for multidigraphs:

\begin{proposition}
\label{prop.mdg.cyc.btf-walk-cyc}Let $D$ be a multidigraph. Let $\mathbf{w}$
be a walk of $D$. Then, $\mathbf{w}$ either is a path or contains a cycle
(i.e., there exists a cycle of $D$ whose arcs are arcs of $\mathbf{w}$).
\end{proposition}

\begin{proof}
This follows by the same argument as Proposition \ref{prop.cyc.btf-walk-cyc}.
\end{proof}

\subsubsection{\label{subsec.dg.walks.algo}Remark on algorithms}

Given a multidigraph $D$ and two vertices $u$ and $v$ of $D$, we can pose the
same five algorithmic questions (Questions 1, 2, 3, 4 and 5) that we posed for
a simple graph $G$ in Subsection \ref{subsec.sg.walks.algo}. As with
multigraphs, the same answers that we gave back then are still valid in our
new setting, as long as we replace \textquotedblleft neighbors of
$v$\textquotedblright\ by \textquotedblleft in-neighbors of $v$%
\textquotedblright\ (that is, vertices $w$ such that $D$ has an arc from $w$
to $v$), and as long as we keep track of the arcs in our paths or walks.

\subsubsection{Exercises}

\begin{exercise}
\label{exe.4.2}Let $D$ be a multidigraph with at least one vertex. Prove the following:

\begin{enumerate}
\item[\textbf{(a)}] If each vertex $v$ of $D$ satisfies $\deg^{+} v > 0$, then
$D$ has a cycle.

\item[\textbf{(b)}] If each vertex $v$ of $D$ satisfies $\deg^{+}v=\deg
^{-}v=1$, then each vertex of $D$ belongs to exactly one cycle of $D$. Here,
two cycles are considered to be identical if one can be obtained from the
other by cyclic rotation.
\end{enumerate}
\end{exercise}

\begin{exercise}
\label{exe.8.8}Let $p$ be a prime number. Let $\left(  a_{1},a_{2}%
,a_{3},\ldots\right)  $ be a sequence of integers that is periodic with period
$p$ (that is, that satisfies $a_{i}=a_{i+p}$ for each $i>0$). Assume that
$a_{1}+a_{2}+\cdots+a_{p}$ is not divisible by $p$. Prove that there exists an
$i\in\left\{  1,2,\ldots,p\right\}  $ such that none of the $p$ numbers
\[
a_{i},\ a_{i}+a_{i+1},\ a_{i}+a_{i+1}+a_{i+2},\ \ldots,\ a_{i}+a_{i+1}%
+\cdots+a_{i+p-1}%
\]
(that is, of the $p$ sums $a_{i}+a_{i+1}+\cdots+a_{j}$ for $i\leq j<i+p$) is
divisible by $p$. \medskip

[\textbf{Remark:} This would be false if $p$ was not prime. For instance, for
$p=4$, the sequence $\left(  0,2,2,2,0,2,2,2,\ldots\right)  $ would be a
counterexample.] \medskip

[\textbf{Hint:} Use Exercise \ref{exe.4.2} \textbf{(a)}. What is the digraph,
and why does it have a cycle?]
\end{exercise}

\begin{exercise}
\label{exe.5.1}Let $D=\left(  V,A,\psi\right)  $ be a multidigraph.

For two vertices $u$ and $v$ of $D$, we shall write $u\overset{\ast
}{\rightarrow}v$ if there exists a path from $u$ to $v$.

A \textbf{root} of $D$ means a vertex $u\in V$ such that each vertex $v\in V$
satisfies $u\overset{\ast}{\rightarrow}v$.

A \textbf{common ancestor} of two vertices $u$ and $v$ means a vertex $w\in V$
such that $w\overset{\ast}{\rightarrow}u$ and $w\overset{\ast}{\rightarrow}v$.

Assume that $D$ has at least one vertex. Prove that $D$ has a root if and only
if every two vertices in $D$ have a common ancestor.
\end{exercise}

The following exercise is both a directed analogue and a generalization of
Mantel's theorem (Theorem \ref{thm.sg.mantel}):

\begin{exercise}
\label{exe.4.6}Let $D$ be a simple digraph with $n$ vertices and $a$ arcs.
Assume that $D$ has no loops, and that we have $a>n^{2}/2$. Prove the following:

\begin{enumerate}
\item[\textbf{(a)}] The digraph $D$ has a cycle of length $3$.

\item[\textbf{(b)}] We define an \textbf{enhanced }$3$\textbf{-cycle} to be a
triple $\left(  u,v,w\right)  $ of distinct vertices of $D$ such that all four
pairs $\left(  u,v\right)  $, $\left(  v,w\right)  $, $\left(  w,u\right)  $
and $\left(  u,w\right)  $ are arcs of $D$. Then, the digraph $D$ has an
enhanced $3$-cycle.
\end{enumerate}
\end{exercise}

\begin{exercise}
\label{exe.6.5}Let $D=\left(  V,A\right)  $ be a simple digraph that has no cycles.

If $\mathbf{v} = \left(  v_{1}, v_{2}, \ldots, v_{n} \right)  $ is a list of
vertices of $D$ (not necessarily a walk!), then a \textbf{back-cut} of
$\mathbf{v}$ shall mean an arc $a \in A$ whose source is $v_{i}$ and whose
target is $v_{j}$ for some $i, j \in\left\{  1, 2, \ldots, n \right\}  $
satisfying $i > j$. (Colloquially speaking, a back-cut of $\mathbf{v}$ is an
arc of $D$ that leads from some vertex of $\mathbf{v}$ to some earlier vertex
of $\mathbf{v}$.)

A list $\mathbf{v}=\left(  v_{1},v_{2},\ldots,v_{n}\right)  $ of vertices of
$D$ is said to be a \textbf{toposort}\footnotemark\ of $D$ if it contains each
vertex of $D$ exactly once and has no back-cuts.

Prove the following:

\begin{enumerate}
\item[\textbf{(a)}] The digraph $D$ has at least one toposort.

\item[\textbf{(b)}] If $D$ has only one toposort, then this toposort is a
Hamiltonian path of $D$.

Here, a \textbf{Hamiltonian path} in $D$ means a walk of $D$ that contains
each vertex of $D$ exactly once.
\end{enumerate}

[\textbf{Example:} For example, the digraph
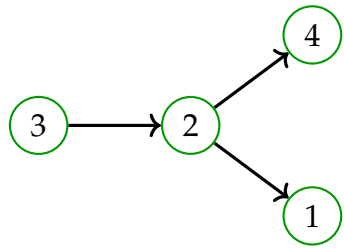
\[
\begin{tikzpicture}[scale=2]
\begin{scope}[every node/.style={circle,thick,draw=green!60!black}]
\node(1) at (-1,0) {$3$};
\node(2) at (0,0) {$2$};
\node(3) at (0.8,-0.6) {$1$};
\node(4) at (0.8,0.6) {$4$};
\end{scope}
\begin{scope}[every edge/.style={draw=black,very thick}, every loop/.style={}]
\path[->] (1) edge (2) (2) edge (3) edge (4);
\end{scope}
\end{tikzpicture}
\]
has two toposorts: $\left(  3,2,1,4\right)  $ and $\left(  3,2,4,1\right)  $.]
\end{exercise}

\footnotetext{This is short for \textquotedblleft topological
sorting\textquotedblright. I don't know where this name comes from.}

\begin{exercise}
\label{exe.9.8}Let $n$ be a positive integer. Let $D$ be a digraph that has no
cycles of length $\leq2$. Assume that $D$ has at least $2^{n-1}$ vertices.
Prove that $D$ has an induced subdigraph that has $n$ vertices and has no cycles.
\end{exercise}

\subsubsection{The adjacency matrix}

A simple way to find the number of walks from a given vertex to a given vertex
in a multidigraph is provided by matrix algebra:

\begin{theorem}
\label{thm.adjmat.walks}Let $D=\left(  V,A,\psi\right)  $ be a multidigraph,
where $V=\left\{  1,2,\ldots,n\right\}  $ for some $n\in\mathbb{N}$.

If $M$ is any matrix, and if $i$ and $j$ are two positive integers, then
$M_{i,j}$ shall denote the $\left(  i,j\right)  $-th entry of $M$ (that is,
the entry of $M$ in the $i$-th row and the $j$-th column).

Let $C$ be the $n\times n$-matrix (with real entries) defined by
\begin{align*}
C_{i,j}  &  =\left(  \text{the number of all arcs $a\in A$ with source $i$ and
target $j$}\right) \\
&  \ \ \ \ \ \ \ \ \ \ \ \ \ \ \ \ \ \ \ \ \text{for all $i,j\in V$.}%
\end{align*}

Let $k\in\mathbb{N}$, and let $i,j\in V$. Then, $\left(  C^{k}\right)  _{i,j}$
equals the number of all walks of $D$ having starting point $i$, ending point
$j$ and length $k$.
\end{theorem}

\begin{remark}
The matrix $C$ in Theorem \ref{thm.adjmat.walks} is known as the
\textbf{adjacency matrix} of $D$. For example, if the multidigraph is
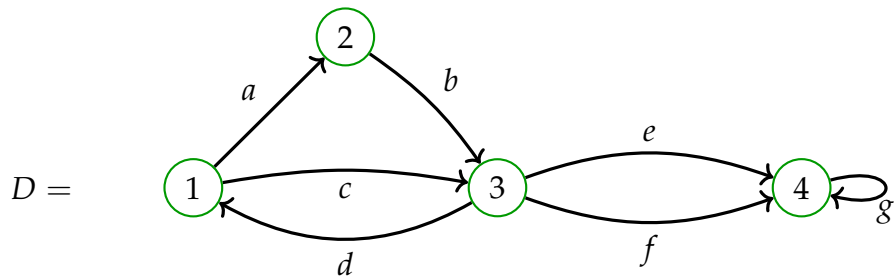
\[
\begin{tikzpicture}[scale=2]
\begin{scope}[every node/.style={circle,thick,draw=green!60!black}]
\node(1) at (-1,0) {$1$};
\node(2) at (0,1) {$2$};
\node(3) at (1,0) {$3$};
\node(4) at (3,0) {$4$};
\end{scope}
\node(X) at (-2, 0) {$D=$};
\begin{scope}[every edge/.style={draw=black,very thick}, every loop/.style={}]
\path[->] (1) edge node[above left] {$a$} (2);
\path[->] (2) edge[bend left=10] node[above right] {$b$} (3);
\path[->] (1) edge[bend left=10] node[below] {$c$} (3);
\path[->] (3) edge[bend left=30] node[below] {$d$} (1);
\path[->] (3) edge[bend left=20] node[above] {$e$} (4);
\path[->] (3) edge[bend right=20] node[below] {$f$} (4);
\path[->] (4) edge[loop right] node[below] {$g$} (4);
\end{scope}
\end{tikzpicture}
\]
then its adjacency matrix is
\[
C=%
\begin{pmatrix}
0 & 1 & 1 & 0\\
0 & 0 & 1 & 0\\
1 & 0 & 0 & 2\\
0 & 0 & 0 & 1
\end{pmatrix}
,
\]
and thus Theorem \ref{thm.adjmat.walks} yields (among other things) that the
$\left(  1,3\right)  $-rd entry $\left(  C^{k}\right)  _{1,3}$ of its $k$-th
power $C^{k}$ equals the number of all walks of $D$ having starting point $1$,
ending point $3$ and length $k$.

The adjacency matrix of a multidigraph $D$ determines $D$ up to the identities
of the arcs, and thus is often used as a convenient way to encode a multidigraph.
\end{remark}

\begin{proof}
[Proof of Theorem \ref{thm.adjmat.walks}.]Forget that we fixed $i$, $j$ and
$k$. We want to prove the following claim:

\begin{statement}
\textit{Claim 1:} Let $i\in V$ and $j\in V$ and $k\in\mathbb{N}$. Then,%
\[
\left(  C^{k}\right)  _{i,j}=\left(  \text{the number of walks from }i\text{
to }j\text{ that have length }k\right)  .
\]

\end{statement}

Before we prove this claim, let us recall that $C$ is the adjacency matrix of
$D$. Thus, for each $i\in V$ and $j\in V$, we have%
\[
C_{i,j}=\left(  \text{the number of all arcs $a\in A$ with source $i$ and
target $j$}\right)
\]
(by the definition of the adjacency matrix). In other words, for each $i\in V$
and $j\in V$, we have%
\[
C_{i,j}=\left(  \text{the number of arcs from }i\text{ to }j\right)  ,
\]
where we agree that an \textquotedblleft arc from $i$ to $j$\textquotedblright%
\ means an arc $a\in A$ with source $i$ and target $j$.

Renaming $i$ as $w$ in this statement, we obtain the following: For each $w\in
V$ and $j\in V$, we have%
\begin{equation}
C_{w,j}=\left(  \text{the number of arcs from }w\text{ to }j\right)  .
\label{sol.graph.count-walks.Aij=}%
\end{equation}

Let us also recall that any two $n\times n$-matrices $M$ and $N$ satisfy
\begin{equation}
\left(  MN\right)  _{i,j}=\sum_{w=1}^{n}M_{i,w}N_{w,j}
\label{sol.graph.count-walks.mmult}%
\end{equation}
for any $i\in V$ and $j\in V$. (Indeed, this is just the rule for how matrices
are multiplied.)

We can now prove Claim 1:

[\textit{Proof of Claim 1:} We shall prove Claim 1 by induction on $k$:

\textit{Induction base:} We shall first prove Claim 1 for $k=0$.

Indeed, let $i\in V$ and $j\in V$. The $0$-th power of any $n\times n$-matrix
is defined to be the $n\times n$ identity matrix $I_{n}$; thus, $C^{0}=I_{n}$.
Hence,
\begin{equation}
\left(  C^{0}\right)  _{i,j}=\left(  I_{n}\right)  _{i,j}=%
\begin{cases}
1, & \text{if }i=j;\\
0, & \text{if }i\neq j
\end{cases}
\label{pf.thm.adjmat.walks.c1.pf.1}%
\end{equation}
(by the definition of the identity matrix).

On the other hand, how many walks from $i$ to $j$ have length $0$ ? A walk
that has length $0$ must consist of a single vertex, which is simultaneously
the starting point and the ending point of this walk. Thus, a walk from $i$ to
$j$ that has length $0$ exists only when $i=j$, and in this case there is
exactly one such walk (namely, the walk $\left(  i\right)  $). Hence,%
\[
\left(  \text{the number of walks from }i\text{ to }j\text{ that have length
}0\right)  =%
\begin{cases}
1, & \text{if }i=j;\\
0, & \text{if }i\neq j.
\end{cases}
\]
Comparing this with (\ref{pf.thm.adjmat.walks.c1.pf.1}), we conclude that%
\begin{equation}
\left(  C^{0}\right)  _{i,j}=\left(  \text{the number of walks from }i\text{
to }j\text{ that have length }0\right)  .
\label{sol.graph.count-walks.goal.IB.1}%
\end{equation}

Now, forget that we fixed $i$ and $j$. We thus have proven
(\ref{sol.graph.count-walks.goal.IB.1}) for any $i\in V$ and $j\in V$. In
other words, Claim 1 holds for $k=0$. Thus, the induction base is complete.

\textit{Induction step:} Let $g$ be a positive integer. Assume that Claim 1
holds for $k=g-1$. We must show that Claim 1 holds for $k=g$ as well.

We have assumed that Claim 1 holds for $k=g-1$. In other words, for any $i\in
V$ and $j\in V$, we have%
\[
\left(  C^{g-1}\right)  _{i,j}=\left(  \text{the number of walks from }i\text{
to }j\text{ that have length }g-1\right)  .
\]
Renaming $j$ as $w$ in this statement, we obtain the following: For any $i\in
V$ and $w\in V$, we have%
\begin{equation}
\left(  C^{g-1}\right)  _{i,w}=\left(  \text{the number of walks from }i\text{
to }w\text{ that have length }g-1\right)  .
\label{sol.graph.count-walks.goal.IH}%
\end{equation}

Each walk from $i$ to $j$ that has length $g$ has the form
\[
\mathbf{w}=\left(  v_{0},a_{1},v_{1},a_{2},v_{2},\ldots,a_{g-1},v_{g-1}%
,a_{g},v_{g}\right)
\]
for some vertices $v_{0},v_{1},\ldots,v_{g}$ of $D$ and some arcs $a_{1}%
,a_{2},\ldots,a_{g}$ of $D$ satisfying $v_{0}=i$ and $v_{g}=j$ and $\left(
\psi\left(  a_{h}\right)  =\left(  v_{h-1},v_{h}\right)  \text{ for all }%
h\in\left\{  1,2,\ldots,g\right\}  \right)  $. Thus, each such walk
$\mathbf{w}$ can be constructed by the following algorithm:

\begin{itemize}
\item First, we choose a vertex $w$ of $D$ to serve as the vertex $v_{g-1}$
(that is, as the penultimate vertex of the walk $\mathbf{w}$). This vertex $w$
must belong to $V$.

\item Now, we choose the vertices $v_{0},v_{1},\ldots,v_{g-1}$ (that is, all
vertices of our walk except for the last one) and the arcs $a_{1},a_{2}%
,\ldots,a_{g-1}$ (that is, all arcs of our walk except for the last one) in
such a way that $v_{g-1}=w$. This is tantamount to choosing a walk $\left(
v_{0},a_{1},v_{1},a_{2},v_{2},\ldots,a_{g-1},v_{g-1}\right)  $ from $i$ to $w$
that has length $g-1$. This choice can be made in $\left(  C^{g-1}\right)
_{i,w}$ many ways (because (\ref{sol.graph.count-walks.goal.IH}) shows that
the number of walks from $i$ to $w$ that have length $g-1$ is $\left(
C^{g-1}\right)  _{i,w}$).

\item We have now determined all but the last vertex and all but the last arc
of our walk $\mathbf{w}$. We set the last vertex $v_{g}$ of our walk to be
$j$. (This is the only possible option, since our walk $\mathbf{w}$ has to be
a walk from $i$ to $j$.)

\item We choose the last arc $a_{g}$ of our walk $\mathbf{w}$. This arc
$a_{g}$ must have source $v_{g-1}$ and target $v_{g}$; in other words, it must
have source $w$ and target $j$ (since $v_{g-1}=w$ and $v_{g}=j$). In other
words, it must be an arc from $w$ to $j$. Thus, it can be chosen in $C_{w,j}$
many ways (because (\ref{sol.graph.count-walks.Aij=}) shows that the number of
arcs from $w$ to $j$ is $C_{w,j}$).
\end{itemize}

Conversely, of course, this algorithm always constructs a walk from $i$ to $j$
that has length $g$, and different choices in the algorithm lead to distinct
walks. Thus, the total number of walks from $i$ to $j$ that have length $g$
equals the total number of choices in the algorithm. But the latter number is
$\sum_{w\in V}\left(  C^{g-1}\right)  _{i,w}C_{w,j}$ (since the algorithm
first chooses a $w\in V$, then involves a step with $\left(  C^{g-1}\right)
_{i,w}$ choices, and then involves a step with $C_{w,j}$ choices). Hence, the
total number of walks from $i$ to $j$ that have length $g$ is $\sum_{w\in
V}\left(  C^{g-1}\right)  _{i,w}C_{w,j}$. In other words,%
\[
\left(  \text{the number of walks from }i\text{ to }j\text{ that have length
}g\right)  =\sum_{w\in V}\left(  C^{g-1}\right)  _{i,w}C_{w,j}.
\]
Comparing this with%
\begin{align*}
\left(  \underbrace{C^{g}}_{=C^{g-1}C}\right)  _{i,j}  &  =\left(
C^{g-1}C\right)  _{i,j}=\sum_{w=1}^{n}\left(  C^{g-1}\right)  _{i,w}C_{w,j}\\
&  \ \ \ \ \ \ \ \ \ \ \ \ \ \ \ \ \ \ \ \ \left(  \text{by
(\ref{sol.graph.count-walks.mmult}) (applied to }M=C^{g-1}\text{ and
}N=C\text{)}\right) \\
&  =\sum_{w\in V}\left(  C^{g-1}\right)  _{i,w}C_{w,j}%
\ \ \ \ \ \ \ \ \ \ \left(  \text{since }\left\{  1,2,\ldots,n\right\}
=V\right)  ,
\end{align*}
we obtain%
\begin{equation}
\left(  C^{g}\right)  _{i,j}=\left(  \text{the number of walks from }i\text{
to }j\text{ that have length }g\right)  .
\label{sol.graph.count-walks.goal.IG}%
\end{equation}

Now, forget that we fixed $i$ and $j$. We thus have proven
(\ref{sol.graph.count-walks.goal.IG}) for any $i\in V$ and $j\in V$. In other
words, Claim 1 holds for $k=g$. Thus, the induction step is complete. Hence,
Claim 1 is proven by induction.]

Theorem \ref{thm.adjmat.walks} follows immediately from Claim 1.
\end{proof}

\begin{exercise}
\label{exe.4.4}Let $E$ be the following multidigraph:
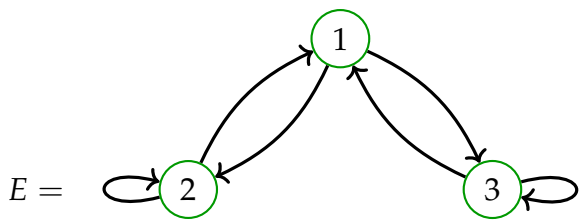
\[
\begin{tikzpicture}[scale=2]
\begin{scope}[every node/.style={circle,thick,draw=green!60!black}]
\node(1) at (0,1) {$1$};
\node(2) at (-1,0) {$2$};
\node(3) at (1,0) {$3$};
\end{scope}
\node(X) at (-2, 0) {$E=$};
\begin{scope}[every edge/.style={draw=black,very thick}, every loop/.style={}]
\path[->] (1) edge[bend left=20] (2);
\path[->] (2) edge[bend left=20] (1);
\path[->] (1) edge[bend left=20] (3);
\path[->] (3) edge[bend left=20] (1);
\path[->] (2) edge[loop left] (2);
\path[->] (3) edge[loop right] (3);
\end{scope}
\end{tikzpicture}
\]
Let $k\in\mathbb{N}$. Compute the number of walks from $1$ to $1$ having
length $k$.
\end{exercise}

\subsubsection{Walks and bidirectionalization}

Let us take a look at what bidirectionalization (i.e., the operation $G\mapsto
G^{\operatorname*{bidir}}$ that sends a multigraph $G$ to the multidigraph
$G^{\operatorname*{bidir}}$) does to walks, paths, closed walks and cycles:

\begin{proposition}
\label{prop.mdg.walks-bidir}Let $G$ be a multigraph. Then:

\begin{enumerate}
\item[\textbf{(a)}] The walks of $G$ are \textquotedblleft more or less the
same as\textquotedblright\ the walks of the multidigraph
$G^{\operatorname*{bidir}}$. More precisely, each walk of $G$ gives rise to a
walk of $G^{\operatorname*{bidir}}$ (with the same starting point and the same
ending point), and conversely, each walk of $G^{\operatorname*{bidir}}$ gives
rise to a walk of $G$. If $G$ has no loops, then this is a one-to-one
correspondence (i.e., a bijection) between the walks of $G$ and the walks of
$G^{\operatorname*{bidir}}$.

\item[\textbf{(b)}] The paths of $G$ are \textquotedblleft more or less the
same as\textquotedblright\ the paths of the multidigraph
$G^{\operatorname*{bidir}}$. This is always a one-to-one correspondence, since
paths cannot contain loops.

\item[\textbf{(c)}] The closed walks of $G$ are \textquotedblleft more or less
the same as\textquotedblright\ the closed walks of the multidigraph
$G^{\operatorname*{bidir}}$.

\item[\textbf{(d)}] The cycles of $G$ are not quite the same as the cycles of
$G^{\operatorname*{bidir}}$. In fact, if $e$ is an edge of $G$ with two
distinct endpoints $u$ and $v$, then $\left(  u,e,v,e,u\right)  $ is not a
cycle of $G$, but either $\left(  u,\left(  e,1\right)  ,v,\left(  e,2\right)
,u\right)  $ or $\left(  u,\left(  e,2\right)  ,v,\left(  e,1\right)
,u\right)  $ is a cycle of $G^{\operatorname*{bidir}}$ (this is best seen on a
picture: $G$ has the edge $%
\begin{tikzpicture}[scale=2]
\begin{scope}[every node/.style={circle,thick,draw=green!60!black}]
\node(u) at (0,0) {$u$};
\node(v) at (1,0) {$v$};
\end{scope}
\begin{scope}[every edge/.style={draw=black,very thick}, every loop/.style={}]
\path[-] (u) edge node[above] {$e$} (v);
\end{scope}
\end{tikzpicture}%
$ whereas $G^{\operatorname*{bidir}}$ has the arc-pair $%
\begin{tikzpicture}[scale=2]
\begin{scope}[every node/.style={circle,thick,draw=green!60!black}]
\node(u) at (0,0) {$u$};
\node(v) at (1,0) {$v$};
\end{scope}
\begin{scope}[every edge/.style={draw=black,very thick}, every loop/.style={}]
\path[->] (u) edge[bend left=20] node[above] {$\tup{e,1}$} (v);
\path[->] (v) edge[bend left=20] node[below] {$\tup{e,2}$} (u);
\end{scope}
\end{tikzpicture}%
$), so $G^{\operatorname*{bidir}}$ usually has more cycles than $G$ has. But
it is true that each cycle of $G$ gives rise to a cycle of
$G^{\operatorname*{bidir}}$.
\end{enumerate}
\end{proposition}

\begin{proof}
This is all straightforward. \medskip

\textbf{(a)} To transform a walk $\mathbf{w}$ of $G$ into a walk
$\mathbf{w}^{\prime}$ of $G^{\operatorname*{bidir}}$ (with the same vertices),
you must replace each edge $e$ of $\mathbf{w}$ by one of the two arcs $\left(
e,1\right)  $ or $\left(  e,2\right)  $ of $G^{\operatorname*{bidir}}$.
(Unless $e$ is a loop, only one of these two arcs will do the job, since its
direction must match the direction in which $\mathbf{w}$ uses the edge $e$.
When $e$ is a loop, however, you get to choose between $\left(  e,1\right)  $
and $\left(  e,2\right)  $.)

Conversely, to transform a walk of $G^{\operatorname*{bidir}}$ into a walk of
$G$, you must replace each arc $\left(  e,i\right)  $ (with $i\in\left\{
1,2\right\}  $) by the corresponding edge $e$ of $G$. \medskip

Parts \textbf{(b)} and \textbf{(c)} easily follow from part \textbf{(a)}
(indeed, the correspondences between the paths and between the closed walks of
$G$ and of $G^{\operatorname*{bidir}}$ are just restrictions of the
correspondence between the walks explained in part \textbf{(a)}). Part
\textbf{(d)} is just as easy.
\end{proof}

\subsection{Connectedness strong and weak}

\subsubsection{Strong connectivity and strong components}

We defined the \textquotedblleft path-connected\textquotedblright\ relation
for undirected graphs using the existence of walks (see Definition
\ref{def.sg.pc}). For a digraph, however, the relations \textquotedblleft
there is a walk from $u$ to $v$\textquotedblright\ and \textquotedblleft there
is a walk from $v$ to $u$\textquotedblright\ are (in general) distinct and
non-symmetric, so I prefer not to give them a symmetric-looking symbol such as
$\simeq_{D}$. Instead, we define \textbf{strong path-connectedness} to mean
the existence of \textbf{both} walks:

\begin{definition}
\label{def.mdg.spc}Let $D$ be a multidigraph. We define a binary relation
$\simeq_{D}$ on the set $\operatorname*{V}\left(  D\right)  $ as follows: For
two vertices $u$ and $v$ of $D$, we shall have $u\simeq_{D}v$ if and only if
there exists a walk from $u$ to $v$ in $D$ and there exists a walk from $v$ to
$u$ in $D$.

This binary relation $\simeq_{D}$ is called \textquotedblleft\textbf{strong
path-connectedness}\textquotedblright. When two vertices $u$ and $v$ satisfy
$u\simeq_{D}v$, we say that \textquotedblleft$u$ and $v$ are \textbf{strongly
path-connected}\textquotedblright.
\end{definition}

\begin{example}
Let $D$ be as in Example \ref{exa.digr.mdg.cycles1}. Then, $1\simeq_{D}2$,
because there exists a walk from $1$ to $2$ in $D$ (for instance, $\left(
1,a,2\right)  $) and there also exists a walk from $2$ to $1$ in $D$ (for
instance, $\left(  2,b,3,d,1\right)  $). However, we don't have $3\simeq_{D}%
4$. Indeed, while there exists a walk from $3$ to $4$ in $D$, there exists no
walk from $4$ to $3$ in $D$.
\end{example}

\begin{proposition}
Let $D$ be a multidigraph. Then, the relation $\simeq_{D}$ is an equivalence relation.
\end{proposition}

\begin{proof}
Easy, like for simple graphs.
\end{proof}

Again, we can replace \textquotedblleft walk\textquotedblright\ by
\textquotedblleft path\textquotedblright\ in the definition of the relation
$\simeq_{D}$:

\begin{proposition}
\label{prop.mdg.spc.there-is-path}Let $D$ be a multidigraph. Let $u$ and $v$
be two vertices of $D$. Then, $u\simeq_{D}v$ if and only if there exist a path
from $u$ to $v$ and a path from $v$ to $u$.
\end{proposition}

\begin{proof}
Easy, like for simple graphs.
\end{proof}

\begin{definition}
\label{def.mdg.scomp}Let $D$ be a multidigraph. The equivalence classes of the
equivalence relation $\simeq_{D}$ are called the \textbf{strong components} of
$D$.
\end{definition}

\begin{definition}
\label{def.mdg.sconn}Let $D$ be a multidigraph. We say that $D$ is
\textbf{strongly connected} if $D$ has exactly one strong component.
\end{definition}

Thus, a multidigraph $D$ is strongly connected if and only if it has at least
one vertex and there is a path from any vertex to any vertex.

\begin{example}
If $D$ is as in Example \ref{exa.digr.mdg.cycles1}, then the strong components
of $D$ are $\left\{  1,2,3\right\}  $ and $\left\{  4\right\}  $.
\end{example}

In analogy to Proposition \ref{prop.sg.component.ind-connected}, strong
components have the following property:

\begin{proposition}
\label{prop.mdg.strong-component.ind-connected}Let $D$ be a multidigraph. Let
$C$ be a strong component of $D$.

Then, the multidigraph $D\left[  C\right]  $ (that is, the induced subdigraph
of $D$ on the set $C$) is strongly connected.
\end{proposition}

\begin{proof}
Just as in the above proof of Proposition
\ref{prop.sg.component.ind-connected}, we can see that $D\left[  C\right]  $
has at least $1$ strong component. It remains to show that $D\left[  C\right]
$ has no more than $1$ component. In other words, it remains to show that any
two vertices of $D\left[  C\right]  $ are strongly path-connected in $D\left[
C\right]  $.

So let $u$ and $v$ be two vertices of $D\left[  C\right]  $. Then, $u,v\in C$,
and therefore $u\simeq_{D}v$ (since $C$ is a component of $D$). In other
words, there exists a walk $\mathbf{u}=\left(  u_{0},u_{1},\ldots
,u_{k}\right)  $ from $u$ to $v$ in $D$, and there exists a walk
$\mathbf{v}=\left(  v_{0},v_{1},\ldots,v_{\ell}\right)  $ from $v$ to $u$ in
$D$. We shall now prove that these walks $\mathbf{u}$ and $\mathbf{v}$
actually are walks of $D\left[  C\right]  $.

In fact: If $u_{i}$ is any vertex of $\mathbf{u}$, then

\begin{itemize}
\item there is a walk from $u$ to $u_{i}$ in $D$ (namely, the first part
$\left(  u_{0},u_{1},\ldots,u_{i}\right)  $ of $\mathbf{u}$), and

\item there is a walk from $u_{i}$ to $u$ in $D$ (namely, the second part
$\left(  u_{i},u_{i+1},\ldots,u_{k}\right)  $ of $\mathbf{u}$, spliced with
the walk $\mathbf{v}$),
\end{itemize}

\noindent and therefore we have $u\simeq_{D}u_{i}$, so that $u_{i}$ belongs to
the same strong component of $D$ as $u$; but that strong component is $C$.
Thus, we have shown that each vertex $u_{i}$ of $\mathbf{u}$ belongs to $C$.
Therefore, $\mathbf{u}$ is a walk of the digraph $D\left[  C\right]  $.
Likewise, $\mathbf{v}$ is a walk of $D\left[  C\right]  $ as well. Hence, the
digraph $D\left[  C\right]  $ has a walk from $u$ to $v$ (namely, $\mathbf{u}%
$) and a walk from $v$ to $u$ (namely, $\mathbf{v}$). In other words,
$u\simeq_{D\left[  C\right]  }v$.

We have now proved that $u\simeq_{D\left[  C\right]  }v$ for any two vertices
$u$ and $v$ of $D\left[  C\right]  $. Hence, the relation $\simeq_{D\left[
C\right]  }$ has no more than $1$ equivalence class. In other words, the
digraph $D\left[  C\right]  $ has no more than $1$ strong component. This
completes our proof.
\end{proof}

\subsubsection{Weak connectivity and weak components}

In comparison, here is a weaker notion of connected components and connectedness:

\begin{definition}
\label{def.mdg.wcomp-wconn}Let $D$ be a multidigraph. Consider its underlying
undirected multigraph $D^{\operatorname*{und}}$. The components of this
undirected multigraph $D^{\operatorname*{und}}$ (that is, the equivalence
classes of the equivalence relation $\simeq_{D^{\operatorname*{und}}}$) are
called the \textbf{weak components} of $D$. We say that $D$ is \textbf{weakly
connected} if $D$ has exactly one weak component (i.e., if
$D^{\operatorname*{und}}$ is connected).
\end{definition}

\begin{example}
\label{exa.mdg.scomp-wcomp}Let $D$ be the following simple digraph:%
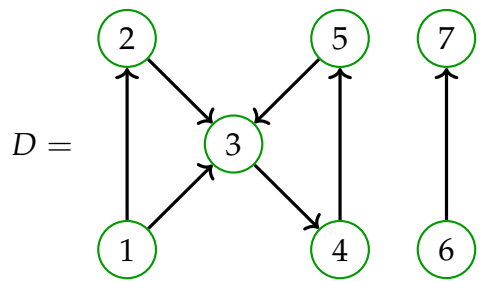
\[%
\begin{tikzpicture}[scale=1.4]
\begin{scope}[every node/.style={circle,thick,draw=green!60!black}]
\node(1) at (-1,-1) {$1$};
\node(2) at (-1,1) {$2$};
\node(3) at (0,0) {$3$};
\node(4) at (1,-1) {$4$};
\node(5) at (1,1) {$5$};
\node(6) at (2,-1) {$6$};
\node(7) at (2,1) {$7$};
\end{scope}
\node(X) at (-1.8, 0) {$D=$};
\begin{scope}[every edge/.style={draw=black,very thick}, every loop/.style={}]
\path[->] (1) edge (2) (2) edge (3) (1) edge (3);
\path[->] (4) edge (5) (5) edge (3) (3) edge (4);
\path[->] (6) edge (7);
\end{scope}
\end{tikzpicture}%
\ \ .
\]
We treat $D$ as a multidigraph (namely, $D^{\operatorname*{mult}}$).

The weak components of $D$ are $\left\{  1,2,3,4,5\right\}  $ and $\left\{
6,7\right\}  $.

The strong components of $D$ are $\left\{  1\right\}  $, $\left\{  2\right\}
$, $\left\{  3,4,5\right\}  $, $\left\{  6\right\}  $ and $\left\{  7\right\}
$. (Indeed, for example, we have $1\not \simeq _{D}2\not \simeq _{D}3$ but
$3\simeq_{D}4\simeq_{D}5$.)

So $D$ is neither strongly nor weakly connected, but has more strong than weak components.
\end{example}

\begin{example}
The digraph from Example \ref{exa.digr.sdg.walks1} is weakly connected, but
not at all strongly connected (indeed, each of its strong components has size
$1$). The digraph from Example \ref{exa.digr.sdg.walks2}, on the other hand,
is strongly connected.
\end{example}

\begin{proposition}
\label{prop.mdg.spc-is-wpc}Any strongly connected digraph is weakly connected.
\end{proposition}

\begin{proof}
Let $D$ be a multidigraph. Then, any walk of $D$ is (or, more precisely, gives
rise to) a walk of $D^{\operatorname*{und}}$. Hence, if two vertices $u$ and
$v$ of $D$ are strongly path-connected in $D$, then they are path-connected in
$D^{\operatorname*{und}}$. Therefore, if $D$ is strongly connected, then
$D^{\operatorname*{und}}$ is connected, but this means that $D$ is weakly connected.
\end{proof}

\begin{exercise}
\label{exe.4.7}Let $D$ be a multidigraph. Prove that the strong components of
$D$ are the weak components of $D$ if and only if each arc of $D$ is contained
in at least one cycle.
\end{exercise}

\subsubsection{Sink components}

Let us come back to strong components for a little sidenote. Not all strong
components are alike; some have special properties:

\begin{definition}
\label{def.mdg.sinkcomp}Let $D$ be a multidigraph. A strong component $C$ of
$D$ is said to be a \textbf{sink component} if it has the following property:
If $a$ is an arc of $D$ whose source lies in $C$, then the target of $a$ also
lies in $C$.
\end{definition}

Informally, a sink component is a strong component \textquotedblleft with no
way out\textquotedblright\ (i.e., no arcs leading out of it).

\begin{example}
If $D$ is the digraph from Example \ref{exa.mdg.scomp-wcomp}, then the sink
components of $D$ are $\left\{  3,4,5\right\}  $ and $\left\{  7\right\}  $.
\end{example}

\begin{example}
Let $D$ be the following simple digraph:%
\[%
\begin{tikzpicture}[scale=1.4]
\begin{scope}[every node/.style={circle,thick,draw=green!60!black}]
\node(1) at (-2, 1) {$1$};
\node(2) at (0, 2) {$2$};
\node(3) at (-1, -0.5) {$3$};
\node(4) at (3, 1) {$4$};
\node(5) at (-4.5, -0.5) {$5$};
\node(6) at (-4, 1) {$6$};
\node(7) at (-3, -0.5) {$7$};
\node(8) at (3, -0.5) {$8$};
\node(9) at (1, -0.5) {$9$};
\end{scope}
\node(X) at (-5.4, 0) {$D=$};
\begin{scope}[every edge/.style={draw=black,very thick}, every loop/.style={}]
\path[->] (1) edge[bend left=20] (3) (3) edge[bend left=20] (1);
\path[->] (4) edge[bend left=20] (8) (8) edge[bend left=20] (4);
\path[->] (6) edge (7) (7) edge[bend left=20] (5) (5) edge (6);
\path[->] (2) edge[loop above] (2);
\path
[->] (2) edge (1) (1) edge (6) (1) edge[bend right=10] (7) (2) edge (4) (8) edge (9) (3) edge (9);
\end{scope}
\end{tikzpicture}%
\ \ .
\]
We treat $D$ as a multidigraph (namely, $D^{\operatorname*{mult}}$). Its
strong components are $\left\{  1,3\right\}  ,\ \left\{  2\right\}
,\ \left\{  4,8\right\}  ,\ \left\{  5,6,7\right\}  ,\ \left\{  9\right\}  $.
Among these five strong components, the only sink components are $\left\{
5,6,7\right\}  $ and $\left\{  9\right\}  $. For example, $\left\{
4,8\right\}  $ is not a sink component, since the arc $89$ has source in
$\left\{  4,8\right\}  $ but target not in $\left\{  4,8\right\}  $.
\end{example}

\begin{theorem}
\label{thm.mdg.sinkcomp}Let $D$ be a multidigraph with at least one vertex.
Then, $D$ has at least one sink component.
\end{theorem}

\begin{proof}
Write the multidigraph $D$ as $D=\left(  V,A,\psi\right)  $. For each vertex
$v\in V$, define $\operatorname*{Desc}v$ to be the set of all vertices $w\in
V$ such that $D$ has a walk from $v$ to $w$. (These vertices are sometimes
called the \textquotedblleft descendants\textquotedblright\ of $v$; thus the
notation $\operatorname*{Desc}v$.) Note that the existence of a walk from $v$
to $w$ is easily testable algorithmically (see Subsection
\ref{subsec.dg.walks.algo}).

The following is easy to see:

\begin{statement}
\textit{Claim 1:} Let $p,q\in V$ be two vertices such that $q\in
\operatorname*{Desc}p$. Then, $\operatorname*{Desc}q\subseteq
\operatorname*{Desc}p$.
\end{statement}

\begin{proof}
[Proof of Claim 1.]Let $w\in\operatorname*{Desc}q$. We must show that
$w\in\operatorname*{Desc}p$.

We have $q\in\operatorname*{Desc}p$; in other words, $D$ has a walk from $p$
to $q$ (by the definition of $\operatorname*{Desc}p$). Likewise, from
$w\in\operatorname*{Desc}q$, we conclude that $D$ has a walk from $q$ to $w$.
Let us denote these two walks by $\mathbf{a}$ and $\mathbf{b}$. Splicing these
two walks together (using Proposition \ref{prop.mdg.walk-concat}, applied to
$u=p$ and $v=q$), we obtain a walk $\mathbf{a}\ast\mathbf{b}$ from $p$ to $w$.
Hence, $D$ has a walk from $p$ to $w$. In other words, we have $w\in
\operatorname*{Desc}p$ (by the definition of $\operatorname*{Desc}p$). This
proves Claim 1.
\end{proof}

Now, recall that the set $V$ is nonempty (since $D$ has at least one vertex)
and finite. Among all the vertices $v\in V$, we pick one vertex $v$ for which
the size $\left\vert \operatorname*{Desc}v\right\vert $ is minimal. Thus,%
\begin{equation}
\left\vert \operatorname*{Desc}v\right\vert \leq\left\vert
\operatorname*{Desc}w\right\vert \ \ \ \ \ \ \ \ \ \ \text{for each }w\in V.
\label{pf.thm.mdg.sinkcomp.min}%
\end{equation}

Now, let $C$ be the strong component of $D$ that contains $v$. We claim that
$C$ is a sink component.

To prove this, we consider an arbitrary arc $a$ whose source lies in $C$. Let
$p$ be this source, and let $q$ be the target of $a$. Thus, $\left(
p,a,q\right)  $ is a walk from $p$ to $q$. Hence, $D$ has a walk from $p$ to
$q$; in other words, $q\in\operatorname*{Desc}p$ (by the definition of
$\operatorname*{Desc}p$). Thus, Claim 1 yields
\[
\operatorname*{Desc}q\subseteq\operatorname*{Desc}p.
\]

On the other hand, $p$ is the source of $a$ and thus lies in $C$ (since we
said that the source of $a$ lies in $C$). In other words, $p$ lies in the same
strong component of $D$ as $v$ (since $C$ is the strong component of $D$ that
contains $v$). In other words, $p\simeq_{D}v$ (since the strong components of
$D$ are the equivalence classes of $\simeq_{D}$). In other words, there exists
a walk from $p$ to $v$ in $D$ and there exists a walk from $v$ to $p$ in $D$
(by the definition of $\simeq_{D}$). In particular, $D$ has a walk from $v$ to
$p$; in other words, $p\in\operatorname*{Desc}v$. Thus, Claim 1 (applied to
$v$ and $p$ instead of $p$ and $q$) yields $\operatorname*{Desc}%
p\subseteq\operatorname*{Desc}v$. Hence,
\[
\left\vert \operatorname*{Desc}p\right\vert \leq\left\vert
\operatorname*{Desc}v\right\vert \leq\left\vert \operatorname*{Desc}%
q\right\vert \ \ \ \ \ \ \ \ \ \ \left(  \text{by
(\ref{pf.thm.mdg.sinkcomp.min}), applied to }w=q\right)  .
\]

However, if two finite sets $X$ and $Y$ satisfy $X\subseteq Y$ and $\left\vert
Y\right\vert \leq\left\vert X\right\vert $, then $X=Y$ (because otherwise, $X$
would be a proper subset of $Y$ and thus satisfy $\left\vert X\right\vert
<\left\vert Y\right\vert $). Applying this fundamental fact to
$X=\operatorname*{Desc}q$ and $Y=\operatorname*{Desc}p$, we conclude that
$\operatorname*{Desc}q=\operatorname*{Desc}p$ (since $\operatorname*{Desc}%
q\subseteq\operatorname*{Desc}p$ and $\left\vert \operatorname*{Desc}%
p\right\vert \leq\left\vert \operatorname*{Desc}q\right\vert $).

But $D$ has a walk from $p$ to $p$ (namely, the trivial walk $\left(
p\right)  $); in other words, $p\in\operatorname*{Desc}p$. In other words,
$p\in\operatorname*{Desc}q$ (since $\operatorname*{Desc}q=\operatorname*{Desc}%
p$). Hence, $D$ has a walk from $q$ to $p$. Since $D$ also has a walk from $p$
to $q$, we thus conclude that there exist both a walk from $p$ to $q$ and a
walk from $q$ to $p$ in $D$. In other words, $p\simeq_{D}q$ (by the definition
of $\simeq_{D}$). Hence, $p$ and $q$ lie in the same strong component of $D$
(since the strong components of $D$ are the equivalence classes of $\simeq
_{D}$). Thus, from $p\in C$, we obtain $q\in C$ (since $C$ is a strong
component). In other words, the target of $a$ lies in $C$ (since $q$ is the
target of $a$).

Forget that we fixed $a$. We thus have shown that if $a$ is an arc of $D$
whose source lies in $C$, then the target of $a$ also lies in $C$. In other
words, $C$ is a sink component of $D$. Hence, $D$ has a sink component, qed.
\end{proof}

\subsubsection{Exercises}

\begin{exercise}
\label{exe.mt2.eclectic-cycle.di} Let $D=\left(  V,E,\psi\right)  $ be a multidigraph.

Let $A$, $B$ and $C$ be three subsets of $V$ such that the induced subdigraphs
$D\left[  A\right]  $, $D\left[  B\right]  $ and $D\left[  C\right]  $ are
strongly connected.

A cycle of $D$ will be called \textbf{eclectic} if it contains at least one
arc of $D\left[  A\right]  $, at least one arc of $D\left[  B\right]  $ and at
least one arc of $D\left[  C\right]  $ (although these three arcs are not
required to be distinct).

Prove the following:

\begin{enumerate}
\item[\textbf{(a)}] If the sets $B\cap C$, $C\cap A$ and $A\cap B$ are
nonempty, but $A\cap B\cap C$ is empty, then $D$ has an eclectic cycle.

\item[\textbf{(b)}] If the induced subdigraphs $D\left[  B\cap C\right]  $,
$D\left[  C\cap A\right]  $ and $D\left[  A\cap B\right]  $ are strongly
connected, but the induced subdigraph $D\left[  A\cap B\cap C\right]  $ is not
strongly connected, then $D$ has an eclectic cycle.
\end{enumerate}

[\textbf{Note:} Keep in mind that the multidigraph with $0$ vertices does not
count as strongly connected.] \medskip

[\textbf{Solution:} This is a generalization of Exercise 7 on midterm \#2 from
my Spring 2017 course; see
\href{https://www.cip.ifi.lmu.de/~grinberg/t/17s/}{the course page} for solutions.]
\end{exercise}

\subsection{Eulerian walks and circuits}

We have studied Eulerian walks and circuits for (undirected) multigraphs in
Section \ref{sec.mg.euler}. Let us now define analogous concepts for multidigraphs:

\begin{definition}
\label{def.digr.eulerian}Let $D$ be a multidigraph.

\begin{enumerate}
\item[\textbf{(a)}] A walk of $D$ is said to be \textbf{Eulerian} if each arc
of $D$ appears exactly once in this walk.

(In other words: A walk $\left(  v_{0},a_{1},v_{1},a_{2},v_{2},\ldots
,a_{k},v_{k}\right)  $ of $D$ is said to be \textbf{Eulerian} if for each arc
$a$ of $D$, there exists exactly one $i\in\left\{  1,2,\ldots,k\right\}  $
such that $a=a_{i}$.)

\item[\textbf{(b)}] An \textbf{Eulerian circuit} of $D$ means a circuit (i.e.,
closed walk) of $D$ that is Eulerian.
\end{enumerate}
\end{definition}

The Euler--Hierholzer theorem gives a necessary and sufficient criterion for a
multigraph to have an Eulerian circuit or walk. For multidigraphs, there is an
analogous result:

\begin{theorem}
[diEuler, diHierholzer]\label{thm.digr.euler-hier}Let $D$ be a weakly
connected multidigraph. Then:

\begin{enumerate}
\item[\textbf{(a)}] The multidigraph $D$ has an Eulerian circuit if and only
if each vertex $v$ of $D$ satisfies $\deg^{+}v=\deg^{-}v$.

\item[\textbf{(b)}] The multidigraph $D$ has an Eulerian walk if and only if
all but two vertices $v$ of $D$ satisfy $\deg^{+}v=\deg^{-}v$, and the
remaining two vertices $v$ satisfy $\left\vert \deg^{+}v-\deg^{-}v\right\vert
\leq1$.
\end{enumerate}
\end{theorem}

\begin{exercise}
\label{exe.4.3}Prove Theorem \ref{thm.digr.euler-hier}.
\end{exercise}

Incidentally, the \textquotedblleft each vertex $v$ of $D$ satisfies $\deg
^{+}v=\deg^{-}v$\textquotedblright\ condition has a name:

\begin{definition}
\label{def.mdg.balanced}A multidigraph $D$ is said to be \textbf{balanced} if
each vertex $v$ of $D$ satisfies $\deg^{+}v=\deg^{-}v$.
\end{definition}

So balancedness is necessary and sufficient for the existence of an Eulerian
circuit in a weakly connected multidigraph.

The following proposition is obvious:

\begin{proposition}
\label{prop.Gdir-balanced}Let $G$ be a multigraph. Then, the multidigraph
$G^{\operatorname*{bidir}}$ is balanced.
\end{proposition}

\begin{proof}
The definition of $G^{\operatorname*{bidir}}$ yields that each vertex $v$ of
$G^{\operatorname*{bidir}}$ satisfies $\deg^{+}v=\deg v$ and $\deg^{-}v=\deg
v$, where $\deg v$ denotes the degree of $v$ as a vertex of $G$. Hence, each
vertex $v$ of $G^{\operatorname*{bidir}}$ satisfies $\deg^{+}v=\deg v=\deg
^{-}v$. In other words, $G^{\operatorname*{bidir}}$ is balanced.
\end{proof}

Combining this proposition with Theorem \ref{thm.digr.euler-hier}
\textbf{(a)}, we can obtain a curious fact about undirected(!) multigraphs:

\begin{theorem}
\label{thm.digr.euler-bidir}Let $G$ be a connected multigraph. Then, the
multidigraph $G^{\operatorname*{bidir}}$ has an Eulerian circuit. In other
words, there is a circuit of $G$ that contains each edge \textbf{exactly
twice}, and uses it once in each direction.
\end{theorem}

\begin{proof}
The multidigraph $G^{\operatorname*{bidir}}$ is balanced (by Proposition
\ref{prop.Gdir-balanced}) and weakly connected (this follows easily from the
connectedness of $G$). Hence, Theorem \ref{thm.digr.euler-hier} \textbf{(a)}
can be applied to $D=G^{\operatorname*{bidir}}$. Thus,
$G^{\operatorname*{bidir}}$ has an Eulerian circuit. Reinterpreting this
circuit as a circuit of $G$, we obtain a circuit of $G$ that contains each
edge \textbf{exactly twice}, and uses it once in each direction. This proves
Theorem \ref{thm.digr.euler-bidir}.
\end{proof}

\begin{example}
Let $G$ be the multigraph%
\[%
\begin{tikzpicture}[scale=2]
\begin{scope}[every node/.style={circle,thick,draw=green!60!black}]
\node(0) at (-2.5,0) {$0$};
\node(1) at (-0.8,0) {$1$};
\node(2) at (0,1) {$2$};
\node(3) at (0.8,0) {$3$};
\node(4) at (2.5,0) {$4$};
\end{scope}
\begin{scope}[every edge/.style={draw=black,very thick}, every loop/.style={}]
\path[-] (0) edge node[above] {$x$} (1);
\path[-] (1) edge node[above left] {$a$} (2);
\path[-] (2) edge node[above right] {$b$} (3);
\path[-] (1) edge node[below] {$c$} (3);
\path[-] (3) edge node[above] {$y$} (4);
\end{scope}
\end{tikzpicture}%
\ \ .
\]
This graph $G$ has no Eulerian circuit, nor even an Eulerian walk (since the
vertices $0,1,3,4$ have odd degrees). But Theorem \ref{thm.digr.euler-bidir}
shows that it has a circuit that contains each edge \textbf{exactly twice},
and uses it once in each direction. And indeed, for example,
\[
\left(  0,x,1,a,2,b,3,y,4,y,3,b,2,a,1,c,3,c,1,x,0\right)
\]
is such a circuit.
\end{example}

The following exercise provides yet another analogue of the Euler--Hierholzer
theorem (Theorem \ref{thm.mg.euler} \textbf{(a)} to be specific), this time
not for multidigraphs but for multigraphs with colored edges:

\begin{exercise}
Let $G$ be a connected multigraph. Assume that each edge of $G$ has been
colored either red or blue. For any vertex $v$ of $G$, we let $\deg^{r}v$
denote the number of red edges containing $v$ (where red loops are counted
twice), and we let $\deg^{b}v$ denote the number of blue edges containing $v$
(where blue loops are counted twice). A walk $\mathbf{w}$ of $G$ is said to be
\textbf{alternating} if no two consecutive edges of $\mathbf{w}$ have the same
color (i.e., if its edges alternate between red and blue, starting with either
red or blue). An alternating walk $\mathbf{w}$ of $G$ is said to be
\textbf{fully alternating} if the first and the last edges of $\mathbf{w}$
have different colors (or else if $\mathbf{w}$ has no edges).

Prove that the multigraph $G$ has a fully alternating Eulerian circuit if and
only if each vertex $v$ of $G$ satisfies $\deg^{r}v=\deg^{b}v$.
\end{exercise}

An example of the situation described in this exercise is shown in the
following picture:%
\[%
\begin{tikzpicture}[scale=2]
\begin{scope}[every node/.style={circle,thick,draw=green!60!black}]
\node(a) at (2,0.5) {$1$};
\node(b) at (1,0) {$2$};
\node(c) at (1,1) {$3$};
\node(d) at (0,1) {$4$};
\node(e) at (0,0) {$5$};
\end{scope}
\begin{scope}[every edge/.style={draw=black,very thick}, every loop/.style={}]
\path[-, draw=blue] (a) edge[blue] node[midway, below] {$a$} (b);
\path[-, red] (b) edge[red, bend left=20] node[midway, left] {$b$} (c);
\path[-, blue] (c) edge[blue, bend left=20] node[midway, right] {$c$} (b);
\path[-, red] (a) edge[red] node[midway, above] {$d$} (c);
\path[-, blue] (b) edge[blue] node[midway, below] {$e$} (e);
\path[-, red] (e) edge[red] node[midway, right] {$f$} (d);
\path[-, blue] (c) edge[blue] node[midway, above] {$g$} (d);
\path[-, red] (b) edge[red, loop below] node {$h$} (b);
\end{scope}
\end{tikzpicture}%
\]
(where the edges $a,c,e,g$ are colored blue while the edges $b,d,f,h$ are
colored red). In this example, $\left(
1,a,2,b,3,g,4,f,5,e,2,h,2,c,3,d,1\right)  $ is a fully alternating Eulerian circuit.

\subsection{Hamiltonian cycles and paths}

We can define Hamiltonian paths and cycles for simple digraphs in the same way
as we defined them for simple graphs:

\begin{definition}
\label{def.sdg.hamp-hamc}Let $D=\left(  V,A\right)  $ be a simple digraph.

\begin{enumerate}
\item[\textbf{(a)}] A \textbf{Hamiltonian path} in $D$ means a walk of $D$
that contains each vertex of $D$ exactly once. Obviously, it is a path.

\item[\textbf{(b)}] A \textbf{Hamiltonian cycle} in $D$ means a cycle $\left(
v_{0},v_{1},\ldots,v_{k}\right)  $ of $D$ such that each vertex of $D$ appears
exactly once among $v_{0},v_{1},\ldots,v_{k-1}$.
\end{enumerate}
\end{definition}

\begin{convention}
In the following, we will abbreviate:

\begin{itemize}
\item \textquotedblleft Hamiltonian path\textquotedblright\ as
\textquotedblleft\textbf{hamp}\textquotedblright;

\item \textquotedblleft Hamiltonian cycle\textquotedblright\ as
\textquotedblleft\textbf{hamc}\textquotedblright.
\end{itemize}
\end{convention}

Hamps and hamcs for digraphs have some properties analogous to those for
simple graphs. In particular, there is an analogue of Ore's theorem:

\begin{theorem}
[Meyniel]Let $D=\left(  V,A\right)  $ be a strongly connected loopless simple
digraph with $n$ vertices. Assume that for each pair $\left(  u,v\right)  \in
V\times V$ of two vertices $u$ and $v$ satisfying $u\neq v$ and $\left(
u,v\right)  \notin A$ and $\left(  v,u\right)  \notin A$, we have $\deg u+\deg
v\geq2n-1$. Here, $\deg w$ means $\deg^{+}w+\deg^{-}w$. Then, $D$ has a hamc.
\end{theorem}

This is known as \textbf{Meyniel's theorem}. For its proof (which is far more
complicated than the proof of Ore's theorem), see \cite{BonTho} or
\cite[\S 10.3, Theorem 7]{Berge91} or \cite{BerLiu98} (which shows an even
more general result). Note that the \textquotedblleft strongly
connected\textquotedblright\ condition is needed.

\subsection{The reverse and complement digraphs}

We take a break from studying hamps (Hamiltonian paths) in order to introduce
two more operations on simple digraphs.

\begin{definition}
\label{def.sdg.rev-comp}Let $D=\left(  V,A\right)  $ be a simple digraph. Then:

\begin{enumerate}
\item[\textbf{(a)}] The elements of $\left(  V\times V\right)  \setminus A$
will be called the \textbf{non-arcs} of $D$.

\item[\textbf{(b)}] The \textbf{reversal} of a pair $\left(  i,j\right)  \in
V\times V$ means the pair $\left(  j,i\right)  $.

\item[\textbf{(c)}] We define $D^{\operatorname*{rev}}$ as the simple digraph
$\left(  V,A^{\operatorname*{rev}}\right)  $, where
\[
A^{\operatorname*{rev}}=\left\{  \left(  j,i\right)  \ \mid\ \left(
i,j\right)  \in A\right\}  .
\]
Thus, $D^{\operatorname*{rev}}$ is the digraph obtained from $D$ by reversing
each arc (i.e., swapping its source and its target). This is called the
\textbf{reversal} of $D$.

\item[\textbf{(d)}] We define $\overline{D}$ as the simple digraph $\left(
V,\ \ \left(  V\times V\right)  \setminus A\right)  $. This is the digraph
that has the same vertices as $D$, but whose arcs are precisely the non-arcs
of $D$. This digraph $\overline{D}$ is called the \textbf{complement} of $D$.
\end{enumerate}
\end{definition}

\begin{example}
Let
\[
D=%
%
\ \ .
\]

\end{example}

\begin{convention}
In the following, the symbol \# means \textquotedblleft
number\textquotedblright. For example,%
\[
\left(  \text{\# of subsets of }\left\{  1,2,3\right\}  \right)  =8.
\]

\end{convention}

We shall now try to count hamps in simple digraphs\footnote{See
\cite{17s-lec7} for a more detailed treatment of this topic.}. As a warmup,
here is a particularly simple case:

\begin{proposition}
\label{prop.hamp.123n}Let $D$ be the simple digraph $\left(  V,A\right)  $,
where%
\[
V=\left\{  1,2,\ldots,n\right\}  \ \ \ \ \ \ \ \ \ \ \text{for some }%
n\in\mathbb{N},
\]
and where%
\[
A=\left\{  \left(  i,j\right)  \ \mid\ i<j\right\}  .
\]
Then, $\left(  \text{\# of hamps of }D\right)  =1$.
\end{proposition}

\begin{proof}
It is easy to see that the only hamp of $D$ is $\left(  1,2,\ldots,n\right)  $.
\end{proof}

The following is easy, too:

\begin{proposition}
\label{prop.tour.berge-triv}Let $D$ be a simple digraph. Then,%
\[
\left(  \text{\# of hamps of }D^{\operatorname*{rev}}\right)  =\left(
\text{\# of hamps of }D\right)  .
\]

\end{proposition}

\begin{proof}
The hamps of $D^{\operatorname*{rev}}$ are obtained from the hamps of $D$ by
walking backwards.
\end{proof}

So far, so boring. What about this:

\begin{theorem}
[Berge's theorem]\label{thm.hamp.Dbar}Let $D$ be a simple digraph. Then,%
\[
\left(  \text{\# of hamps of }\overline{D}\right)  \equiv\left(  \text{\# of
hamps of }D\right)  \operatorname{mod}2.
\]

\end{theorem}

This is much less obvious or even expected. We first give an example:

\begin{example}
\label{exa.hamp.comp}Let $D$ be the following digraph:%
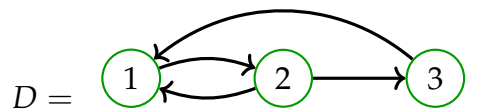
\[
D=\ \
\begin{tikzpicture}
\begin{scope}[every node/.style={circle,thick,draw=green!60!black}]
\node(A) at (0,0) {$1$};
\node(B) at (2,0) {$2$};
\node(C) at (4,0) {$3$};
\end{scope}
\begin{scope}[every edge/.style={draw=black,very thick}]
\path[->] (B) edge (C);
\path[->] (A) edge[bend left=20] (B);
\path[->] (B) edge[bend left=20] (A);
\path[->] (C) edge[bend right=40] (A);
\end{scope}
\end{tikzpicture}%
\ \ .
\]
This digraph has $3$ hamps: $\left(  1,2,3\right)  $ and $\left(
2,3,1\right)  $ and $\left(  3,1,2\right)  $.

Its complement $\overline{D}$ looks as follows:%
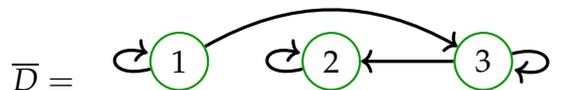
\[
\overline{D}=\ \
\begin{tikzpicture}
\begin{scope}[every node/.style={circle,thick,draw=green!60!black}]
\node(A) at (0,0) {$1$};
\node(B) at (2,0) {$2$};
\node(C) at (4,0) {$3$};
\end{scope}
\begin{scope}[every edge/.style={draw=black,very thick}]
\path[->] (A) edge[loop left] (A);
\path[->] (B) edge[loop left] (B);
\path[->] (C) edge[loop right] (C);
\path[->] (C) edge (B);
\path[->] (A) edge[bend left=30] (C);
\end{scope}
\end{tikzpicture}%
\ \ .
\]
It has only $1$ hamp: $\left(  1,3,2\right)  $.

Thus, in this case, Theorem \ref{thm.hamp.Dbar} says that $1\equiv
3\operatorname{mod}2$.
\end{example}

\begin{proof}
[Proof of Theorem \ref{thm.hamp.Dbar}.](This is an outline; see \cite[proof of
Theorem 1.3.6]{17s-lec7} for more details.)

Write the simple digraph $D$ as $D=\left(  V,A\right)  $, and assume WLOG that
$V\neq\varnothing$. Set $n=\left\vert V\right\vert $.

A $V$\textbf{-listing} will mean a list of elements of $V$ that contains each
element of $V$ exactly once. (Thus, each $V$-listing is an $n$-tuple, and
there are $n!$ many $V$-listings.) Note that a $V$-listing is the same as a
hamp of the \textquotedblleft complete\textquotedblright\ digraph $\left(
V,V\times V\right)  $. Any hamp of $D$ or of $\overline{D}$ is therefore a
$V$-listing, but not every $V$-listing is a hamp of $D$ or $\overline{D}$.

If $\sigma=\left(  \sigma_{1},\sigma_{2},\ldots,\sigma_{n}\right)  $ is a
$V$-listing, then we define a set%
\[
P\left(  \sigma\right)  :=\left\{  \sigma_{1}\sigma_{2},\ \sigma_{2}\sigma
_{3},\ \ldots,\ \sigma_{n-1}\sigma_{n}\right\}  .
\]
We call this set $P\left(  \sigma\right)  $ the \textbf{arc set} of $\sigma$.
When we regard $\sigma$ as a hamp of $\left(  V,V\times V\right)  $, this set
$P\left(  \sigma\right)  $ is just the set of all arcs of $\sigma$. Note that
this is an $\left(  n-1\right)  $-element set. We make a few easy observations
(prove them!):

\begin{statement}
\textit{Observation 1:} We can reconstruct a $V$-listing $\sigma$ from its arc
set $P\left(  \sigma\right)  $. In other words, the map $\sigma\mapsto
P\left(  \sigma\right)  $ is injective.
\end{statement}

\begin{statement}
\textit{Observation 2:} Let $\sigma$ be a $V$-listing. Then, $\sigma$ is a
hamp of $D$ if and only if $P\left(  \sigma\right)  \subseteq A$.
\end{statement}

\begin{statement}
\textit{Observation 3:} Let $\sigma$ be a $V$-listing. Then, $\sigma$ is a
hamp of $\overline{D}$ if and only if $P\left(  \sigma\right)  \subseteq
\left(  V\times V\right)  \setminus A$.
\end{statement}

Now, let $N$ be the \# of pairs $\left(  \sigma,B\right)  $, where $\sigma$ is
a $V$-listing and $B$ is a subset of $A$ satisfying $B\subseteq P\left(
\sigma\right)  $. Thus,%
\[
N=\sum_{\sigma\text{ is a }V\text{-listing}}N_{\sigma},
\]
where%
\[
N_{\sigma}=\left(  \text{\# of subsets }B\text{ of }A\text{ satisfying
}B\subseteq P\left(  \sigma\right)  \right)  .
\]
But we also have%
\[
N=\sum_{B\text{ is a subset of }A}N^{B},
\]
where%
\[
N^{B}=\left(  \text{\# of }V\text{-listings }\sigma\text{ satisfying
}B\subseteq P\left(  \sigma\right)  \right)  .
\]

Let us now relate these two sums to hamps. We begin with $\sum\limits_{\sigma
\text{ is a }V\text{-listing}}N_{\sigma}$.

We shall use the \textbf{Iverson bracket notation}: i.e., the notation
$\left[  \mathcal{A}\right]  $ for the truth value of a statement
$\mathcal{A}$. This truth value is defined to be the number $1$ if
$\mathcal{A}$ is true, and $0$ if $\mathcal{A}$ is false. For instance,%
\[
\left[  2+2=4\right]  =1\ \ \ \ \ \ \ \ \ \ \text{and}%
\ \ \ \ \ \ \ \ \ \ \left[  2+2=5\right]  =0.
\]

For any $V$-listing $\sigma$, we have%
\begin{align*}
N_{\sigma}  &  =\left(  \text{\# of subsets }B\text{ of }A\text{ satisfying
}B\subseteq P\left(  \sigma\right)  \right) \\
&  =\left(  \text{\# of subsets }B\text{ of }A\cap P\left(  \sigma\right)
\right) \\
&  =2^{\left\vert A\cap P\left(  \sigma\right)  \right\vert }\\
&  \equiv\left[  \left\vert A\cap P\left(  \sigma\right)  \right\vert
=0\right]  \ \ \ \ \ \ \ \ \ \ \left(  \text{since }2^{m}\equiv\left[
m=0\right]  \operatorname{mod}2\text{ for each }m\in\mathbb{N}\right) \\
&  =\left[  A\cap P\left(  \sigma\right)  =\varnothing\right]
\ \ \ \ \ \ \ \ \ \ \left(
\begin{array}
[c]{c}%
\text{since equivalent statements have the}\\
\text{same truth value}%
\end{array}
\right) \\
&  =\left[  P\left(  \sigma\right)  \subseteq\left(  V\times V\right)
\setminus A\right]  \ \ \ \ \ \ \ \ \ \ \left(  \text{since }P\left(
\sigma\right)  \text{ is always a subset of }V\times V\right) \\
&  =\left[  \sigma\text{ is a hamp of }\overline{D}\right]  \operatorname{mod}%
2\ \ \ \ \ \ \ \ \ \ \left(  \text{by Observation 3}\right)  .
\end{align*}
So%
\begin{align*}
N  &  =\sum_{\sigma\text{ is a }V\text{-listing}}\underbrace{N_{\sigma}%
}_{\equiv\left[  \sigma\text{ is a hamp of }\overline{D}\right]
\operatorname{mod}2}\\
&  \equiv\sum_{\sigma\text{ is a }V\text{-listing}}\left[  \sigma\text{ is a
hamp of }\overline{D}\right] \\
&  =\left(  \text{\# of }V\text{-listings }\sigma\text{ that are hamps of
}\overline{D}\right) \\
&  \ \ \ \ \ \ \ \ \ \ \ \ \ \ \ \ \ \ \ \ \left(
\begin{array}
[c]{c}%
\text{because }\sum_{\sigma\text{ is a }V\text{-listing}}\left[  \sigma\text{
is a hamp of }\overline{D}\right]  \text{ is a sum}\\
\text{of several }1\text{'s and several }0\text{'s, and the }1\text{'s in
this}\\
\text{sum correspond precisely to}\\
\text{the }V\text{-listings }\sigma\text{ that are hamps of }\overline{D}%
\end{array}
\right) \\
&  =\left(  \text{\# of hamps of }\overline{D}\right)  \operatorname{mod}2.
\end{align*}

What about the other expression for $N$ ? Recall that%
\[
N=\sum_{B\text{ is a subset of }A}N^{B},
\]
where%
\[
N^{B}=\left(  \text{\# of }V\text{-listings }\sigma\text{ satisfying
}B\subseteq P\left(  \sigma\right)  \right)  .
\]
We want to prove that this sum equals $\left(  \text{\# of hamps of }D\right)
$, at least modulo $2$.

So let $B$ be a subset of $A$. We want to know $N^{B}\operatorname{mod}2$. In
other words, we want to know when $N^{B}$ is odd.

Let us first assume that $N^{B}$ is odd, and see what follows from this.

Since $N^{B}$ is odd, we have $N^{B}>0$. Thus, there exists \textbf{at least
one} $V$-listing $\sigma$ satisfying $B\subseteq P\left(  \sigma\right)  $. We
shall now draw some conclusions from this.

First, a definition: A \textbf{path cover} of $V$ means a set of paths in the
\textquotedblleft complete\textquotedblright\ digraph $\left(  V,V\times
V\right)  $ such that each vertex $v\in V$ is contained in exactly one of
these paths. The \textbf{set of arcs} of such a path cover is simply the set
of all arcs of all its paths. For example, if $V=\left\{
1,2,3,4,5,6,7\right\}  $, then
\[
\left\{  \left(  1,3,5\right)  ,\ \ \left(  2\right)  ,\ \ \left(  6\right)
,\ \ \left(  7,4\right)  \right\}
\]
is a path cover of $V$, and its set of arcs is $\left\{  13,\ 35,\ 74\right\}
$.

Now, ponder the following: If we remove an arc $v_{i}v_{i+1}$ from a path
$\left(  v_{1},v_{2},\ldots,v_{k}\right)  $, then this path breaks up into two
paths $\left(  v_{1},v_{2},\ldots,v_{i}\right)  $ and $\left(  v_{i+1}%
,v_{i+2},\ldots,v_{k}\right)  $. Thus, if we remove some arcs from the arc set
$P\left(  \sigma\right)  $ of a $V$-listing $\sigma$, then we obtain the set
of arcs of a path cover of $V$. (For instance, removing the arcs $52$, $26$
and $67$ from the arc set $P\left(  \sigma\right)  $ of the $V$-listing
$\sigma=\left(  1,3,5,2,6,7,4\right)  $ yields precisely the path cover
$\left\{  \left(  1,3,5\right)  ,\ \ \left(  2\right)  ,\ \ \left(  6\right)
,\ \ \left(  7,4\right)  \right\}  $ that we just showed as an example.)

Now, recall that there exists \textbf{at least one} $V$-listing $\sigma$
satisfying $B\subseteq P\left(  \sigma\right)  $. Hence, $B$ is obtained by
removing some arcs from the arc set $P\left(  \sigma\right)  $ of this
$V$-listing $\sigma$. Therefore, $B$ is the set of arcs of a path cover of $V$
(by the claim of the preceding paragraph). Let us say that this path cover
consists of exactly $r$ paths. Then,
\[
\left(  \text{\# of }V\text{-listings }\sigma\text{ satisfying }B\subseteq
P\left(  \sigma\right)  \right)  =r!,
\]
because any such $V$-listing $\sigma$ can be constructed by concatenating the
$r$ paths in our path cover in some order (and there are $r!$ possible orders).

Thus, $N^{B}=\left(  \text{\# of }V\text{-listings }\sigma\text{ satisfying
}B\subseteq P\left(  \sigma\right)  \right)  =r!$. But we have assumed that
$N^{B}$ is odd. So $r!$ is odd. Since $r$ is positive (because $V\neq
\varnothing$, so our path cover must contain at least one path), this entails
that $r=1$. So our path cover is just a single path; this path is a path of
$D$ (since its set of arcs $B$ is a subset of $A$) and therefore is a hamp of
$D$ (since it constitutes a path cover of $V$ all by itself). If we denote it
by $\sigma$, then we have $B=P\left(  \sigma\right)  $ (since $B$ is the set
of arcs of the path cover that consists of $\sigma$ alone).

Forget our assumption that $N^{B}$ is odd. We have thus shown that if $N^{B}$
is odd, then $B=P\left(  \sigma\right)  $ for some hamp $\sigma$ of $D$.

Conversely, it is easy to see that if $B=P\left(  \sigma\right)  $ for some
hamp $\sigma$ of $D$, then $N^{B}$ is odd (and actually equals $1$).

Combining these two results, we see that $N^{B}$ is odd \textbf{if and only
if} $B=P\left(  \sigma\right)  $ for some hamp $\sigma$ of $D$. Therefore,%
\[
\left[  N^{B}\text{ is odd}\right]  =\left[  B=P\left(  \sigma\right)  \text{
for some hamp }\sigma\text{ of }D\right]  .
\]
However,%
\begin{align*}
N^{B}  &  \equiv\left[  N^{B}\text{ is odd}\right]
\ \ \ \ \ \ \ \ \ \ \left(  \text{since }m\equiv\left[  m\text{ is
odd}\right]  \operatorname{mod}2\text{ for any }m\in\mathbb{Z}\right) \\
&  =\left[  B=P\left(  \sigma\right)  \text{ for some hamp }\sigma\text{ of
}D\right]  \operatorname{mod}2.
\end{align*}

We have proved this congruence for every subset $B$ of $A$. Thus,%
\begin{align*}
N  &  =\sum_{B\text{ is a subset of }A}\underbrace{N^{B}}_{\equiv\left[
B=P\left(  \sigma\right)  \text{ for some hamp }\sigma\text{ of }D\right]
\operatorname{mod}2}\\
&  \equiv\sum_{B\text{ is a subset of }A}\left[  B=P\left(  \sigma\right)
\text{ for some hamp }\sigma\text{ of }D\right] \\
&  =\left(  \text{\# of subsets }B\text{ of }A\text{ such that }B=P\left(
\sigma\right)  \text{ for some hamp }\sigma\text{ of }D\right) \\
&  =\left(  \text{\# of sets of the form }P\left(  \sigma\right)  \text{ for
some hamp }\sigma\text{ of }D\right) \\
&  \ \ \ \ \ \ \ \ \ \ \ \ \ \ \ \ \ \ \ \ \left(
\begin{array}
[c]{c}%
\text{because each set of the form }P\left(  \sigma\right)  \text{ for some}\\
\text{hamp }\sigma\text{ of }D\text{ is a subset of }A\text{ (by Observation
2)}%
\end{array}
\right) \\
&  =\left(  \text{\# of hamps of }D\right)  \operatorname{mod}2
\end{align*}
(indeed, Observation 1 shows that different hamps $\sigma$ have different sets
$P\left(  \sigma\right)  $, so counting the sets $P\left(  \sigma\right)  $
for all hamps $\sigma$ is equivalent to counting the hamps $\sigma$ themselves).

Now we have proved that $N\equiv\left(  \text{\# of hamps of }\overline
{D}\right)  \operatorname{mod}2$ and \newline$N\equiv\left(  \text{\# of hamps
of }D\right)  \operatorname{mod}2$. Comparing these two congruences, we obtain%
\[
\left(  \text{\# of hamps of }\overline{D}\right)  \equiv\left(  \text{\# of
hamps of }D\right)  \operatorname{mod}2.
\]
This proves Berge's theorem.
\end{proof}

\subsection{Tournaments}

\subsubsection{Definition}

We now introduce a special class of simple digraphs.

\begin{definition}
\label{def.tour.def}A \textbf{tournament} is defined to be a loopless simple
digraph $D$ that satisfies the

\begin{itemize}
\item \textbf{Tournament axiom:} For any two distinct vertices $u$ and $v$ of
$D$, \textbf{exactly} one of $\left(  u,v\right)  $ and $\left(  v,u\right)  $
is an arc of $D$.
\end{itemize}
\end{definition}

\begin{example}
\label{exa.tour.exas1}\ \ 

\begin{enumerate}
\item[\textbf{(a)}] The following two digraphs are tournaments:%
\[%
\begin{tabular}
[c]{|c|c|}\hline
\ \ $%
\begin{tikzpicture}
\begin{scope}[every node/.style={circle,thick,draw=green!60!black}]
\node(A) at (0:1.5) {$1$};
\node(B) at (120:1.5) {$2$};
\node(C) at (240:1.5) {$3$};
\end{scope}
\begin{scope}[every edge/.style={draw=black,very thick}]
\path[->] (A) edge (B) (B) edge (C) (A) edge (C);
\end{scope}
\end{tikzpicture}%
\ \ $ & \ \ $%
\begin{tikzpicture}
\begin{scope}[every node/.style={circle,thick,draw=green!60!black}]
\node(A) at (0:1.5) {$1$};
\node(B) at (120:1.5) {$2$};
\node(C) at (240:1.5) {$3$};
\end{scope}
\begin{scope}[every edge/.style={draw=black,very thick}]
\path[->] (A) edge (B) (B) edge (C) (C) edge (A);
\end{scope}
\end{tikzpicture}%
\ \ $\\\hline
\end{tabular}
\ \ .
\]

\item[\textbf{(b)}] None of the following three digraphs is a tournament:%
\[%
\begin{tabular}
[c]{|c|c|c|}\hline
\ \ $%
\begin{tikzpicture}
\begin{scope}[every node/.style={circle,thick,draw=green!60!black}]
\node(A) at (0:1.5) {$1$};
\node(B) at (120:1.5) {$2$};
\node(C) at (240:1.5) {$3$};
\end{scope}
\begin{scope}[every edge/.style={draw=black,very thick}]
\path[->] (A) edge (B) (B) edge (C);
\end{scope}
\end{tikzpicture}%
\ \ $ & \ \ $%
\begin{tikzpicture}
\begin{scope}[every node/.style={circle,thick,draw=green!60!black}]
\node(A) at (0:1.5) {$1$};
\node(B) at (120:1.5) {$2$};
\node(C) at (240:1.5) {$3$};
\end{scope}
\begin{scope}[every edge/.style={draw=black,very thick}]
\path[->] (B) edge (C);
\path[->] (A) edge[bend left=20] (B);
\path[->] (B) edge[bend left=20] (A);
\path[->] (C) edge[bend left=10] (A);
\end{scope}
\end{tikzpicture}%
\ \ $ & \ \ $%
\begin{tikzpicture}
\begin{scope}[every node/.style={circle,thick,draw=green!60!black}]
\node(A) at (0:1.5) {$1$};
\node(B) at (120:1.5) {$2$};
\node(C) at (240:1.5) {$3$};
\end{scope}
\begin{scope}[every edge/.style={draw=black,very thick}]
\path[->] (A) edge (B) (B) edge (C) (A) edge (C);
\path[->] (A) edge[loop right] (A);
\end{scope}
\end{tikzpicture}%
\ \ $\\\hline
\end{tabular}
\ \ ,
\]
because

\begin{itemize}
\item the first of these three digraphs violates the tournament axiom for
$u=1$ and $v=3$;

\item the second of these three digraphs violates the tournament axiom for
$u=1$ and $v=2$;

\item the third of these three digraphs is not a tournament, since it is not loopless.
\end{itemize}

\item[\textbf{(c)}] The digraph $D$ in Proposition \ref{prop.hamp.123n} always
is a tournament.
\end{enumerate}
\end{example}

\begin{example}
\label{exe.tour.5-v-t}Here is a tournament with $5$ vertices:%
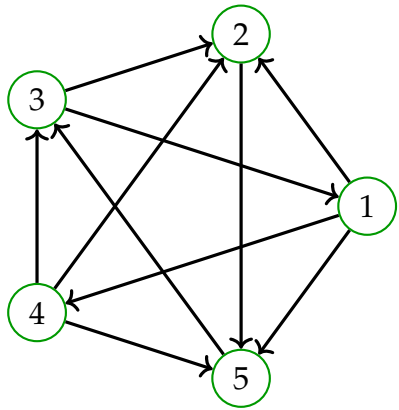
\[%
%
\ \ .
\]

\end{example}

A tournament can also be viewed as a complete graph, whose each edge has been
given a direction.

Using Definition \ref{def.sdg.rev-comp}, we can restate the definition of a
tournament as follows:

\begin{proposition}
\label{prop.tour.rev-bar-def}Let $D=\left(  V,A\right)  $ be a loopless simple
digraph. Then, $D$ is a tournament if and only if the non-loop arcs of
$\overline{D}$ are precisely the arcs of $D^{\operatorname*{rev}}$.
\end{proposition}

\begin{proof}
Easy consequence of definitions.
\end{proof}

\begin{exercise}
\label{exe.tour.kingchicken}Let $D$ be a tournament with at least one vertex.

We say that a vertex $u$ of $D$ \textbf{directly owns} a vertex $w$ of $D$ if
$\left(  u,w\right)  $ is an arc of $D$.

We say that a vertex $u$ of $D$ \textbf{indirectly owns} a vertex $w$ of $D$
if there exists a vertex $v$ of $D$ such that both $\left(  u,v\right)  $ and
$\left(  v,w\right)  $ are arcs of $D$.

Prove that $D$ has a vertex that (directly or indirectly) owns all other
vertices. \medskip

[\textbf{Solution:} This exercise appears in \cite[Exercise 6.3.1]{20f}
(restated in the language of players and matches) and in \cite[Theorem
1]{Maurer80} (restated in the language of chickens and pecking orders). It
originates in a study of pecking orders by Landau \cite{Landau53}.]
\end{exercise}

\subsubsection{The R\'{e}dei theorems}

Which tournaments have hamps? The answer is surprisingly simple:\footnote{Here
we agree to consider the empty list $\left(  {}\right)  $ to be a hamp of the
digraph $\left(  \varnothing,\varnothing\right)  $.}

\begin{theorem}
[Easy R\'{e}dei theorem]\label{thm.tour.redei-easy}A tournament always has at
least one hamp.
\end{theorem}

Even better, and perhaps even more surprisingly:

\begin{theorem}
[Hard R\'{e}dei theorem]\label{thm.tour.redei-hard}Let $D$ be a tournament.
Then,%
\[
\left(  \text{\# of hamps of }D\right)  \text{ is odd.}%
\]

\end{theorem}

Our goal now is to prove these two theorems. Clearly, the Easy R\'{e}dei
Theorem follows from the Hard one, since an odd number cannot be $0$. Thus, it
will suffice to prove the Hard one.

The proof of the hard R\'{e}dei theorem will rely on the following crucial lemma:

\begin{lemma}
\label{lem.tour.redei-lem}Let $D=\left(  V,A\right)  $ be a tournament, and
let $vw\in A$ be an arc of $D$.

Let $D^{\prime}$ be the digraph obtained from $D$ by reversing the arc $vw$.
In other words, let%
\[
D^{\prime}:=\left(  V,\ \ \left(  A\setminus\left\{  vw\right\}  \right)
\cup\left\{  wv\right\}  \right)  .
\]
Then, $D^{\prime}$ is again a tournament, and satisfies%
\[
\left(  \text{\# of hamps of }D\right)  \equiv\left(  \text{\# of hamps of
}D^{\prime}\right)  \operatorname{mod}2.
\]

\end{lemma}

Here is a visualization of the setup of Lemma \ref{lem.tour.redei-lem}:
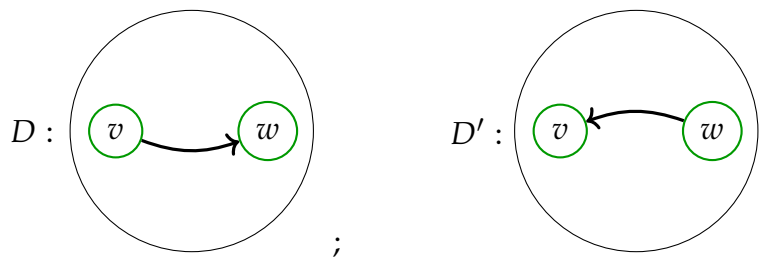
\[%
\begin{tikzpicture}
\draw(1,0) circle (1.6);
\begin{scope}[every node/.style={circle,thick,draw=green!60!black}]
\node(v) at (0,0) {$v$};
\node(w) at (2,0) {$w$};
\end{scope}
\begin{scope}[every edge/.style={draw=black,very thick}]
\path[->] (v) edge[bend right=20] (w);
\end{scope}
\node(D) at (-1.1, 0) {$D : $};
\end{tikzpicture}%
\ \ ;\ \ \ \ \ \ \ \ \ \
\begin{tikzpicture}
\draw(1,0) circle (1.6);
\begin{scope}[every node/.style={circle,thick,draw=green!60!black}]
\node(v) at (0,0) {$v$};
\node(w) at (2,0) {$w$};
\end{scope}
\begin{scope}[every edge/.style={draw=black,very thick}]
\path[->] (w) edge[bend right=20] (v);
\end{scope}
\node(D) at (-1.1, 0) {$D' : $};
\end{tikzpicture}%
\ \ .
\]
(Here, we are only showing the arcs joining $v$ with $w$, since $D$ and
$D^{\prime}$ agree in all other arcs.)

\begin{proof}
[Proof of Lemma \ref{lem.tour.redei-lem}.](This is an outline; see \cite[proof
of Lemma 1.6.2]{17s-lec7} for more details.)

First of all, $D^{\prime}$ is clearly a tournament. It remains to prove the congruence.

We introduce two more digraphs: Let
\begin{align*}
D_{0}  &  :=\left(  \text{the digraph }D\text{ with the arc }vw\text{
removed}\right)  \ \ \ \ \ \ \ \ \ \ \text{and}\\
D_{2}  &  :=\left(  \text{the digraph }D\text{ with the arc }wv\text{
added}\right)  .
\end{align*}
Note that these are not tournaments any more. Here is a comparative
illustration of all four digraphs $D$, $D^{\prime}$, $D_{0}$ and $D_{2}$
(again showing only the arcs joining $v$ with $w$, since there are no
differences in the other arcs):%
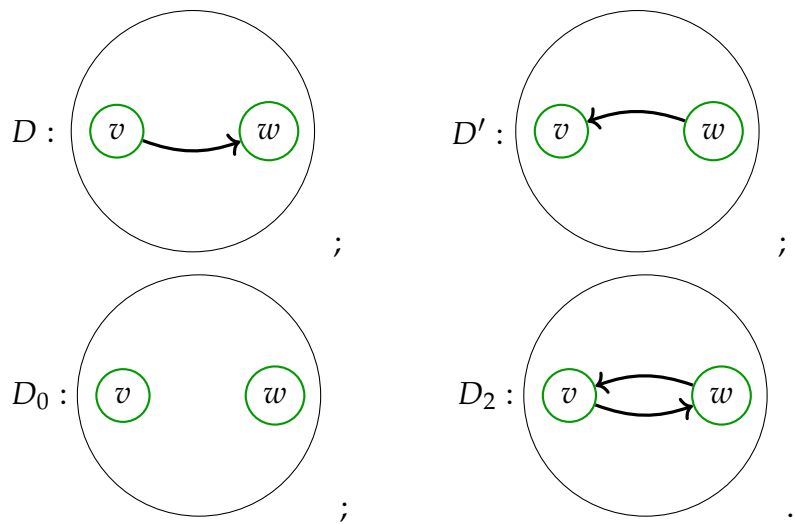
\begin{align*}
&
%
\ \ .
\end{align*}

The digraph $D_{0}$ is $D^{\prime}$ with the arc $wv$ removed. Therefore, a
hamp of $D_{0}$ is the same as a hamp of $D^{\prime}$ that does not use the
arc $wv$. Hence,%
\begin{align*}
&  \left(  \text{\# of hamps of }D_{0}\right) \\
&  =\left(  \text{\# of hamps of }D^{\prime}\text{ that do not use the arc
}wv\right) \\
&  =\left(  \text{\# of hamps of }D^{\prime}\right)  -\left(  \text{\# of
hamps of }D^{\prime}\text{ that use the arc }wv\right)  .
\end{align*}
Similarly, since $D$ is $D_{2}$ with the arc $wv$ removed, we have%
\begin{align*}
&  \left(  \text{\# of hamps of }D\right) \\
&  =\left(  \text{\# of hamps of }D_{2}\right)  -\left(  \text{\# of hamps of
}D_{2}\text{ that use the arc }wv\right) \\
&  =\left(  \text{\# of hamps of }D_{2}\right)  -\left(  \text{\# of hamps of
}D^{\prime}\text{ that use the arc }wv\right)
\end{align*}
(the last equality is because a hamp of $D_{2}$ that uses the arc $wv$ cannot
use the arc $vw$, and therefore is automatically a hamp of $D^{\prime}$ as
well, and of course the converse is obviously true).

However, from the previously proved equality%
\begin{align*}
&  \left(  \text{\# of hamps of }D_{0}\right) \\
&  =\left(  \text{\# of hamps of }D^{\prime}\right)  -\left(  \text{\# of
hamps of }D^{\prime}\text{ that use the arc }wv\right)  ,
\end{align*}
we obtain%
\begin{align*}
&  \left(  \text{\# of hamps of }D^{\prime}\right) \\
&  =\left(  \text{\# of hamps of }D_{0}\right)  +\left(  \text{\# of hamps of
}D^{\prime}\text{ that use the arc }wv\right) \\
&  \equiv\left(  \text{\# of hamps of }D_{0}\right)  -\left(  \text{\# of
hamps of }D^{\prime}\text{ that use the arc }wv\right)  \operatorname{mod}2
\end{align*}
(since $x+y\equiv x-y\operatorname{mod}2$ for any integers $x$ and $y$). Thus,
if we can show that%
\[
\left(  \text{\# of hamps of }D_{2}\right)  \equiv\left(  \text{\# of hamps of
}D_{0}\right)  \operatorname{mod}2,
\]
then we will be able to conclude that%
\begin{align*}
&  \left(  \text{\# of hamps of }D\right) \\
&  =\underbrace{\left(  \text{\# of hamps of }D_{2}\right)  }_{\equiv\left(
\text{\# of hamps of }D_{0}\right)  \operatorname{mod}2}-\left(  \text{\# of
hamps of }D^{\prime}\text{ that use the arc }wv\right) \\
&  \equiv\left(  \text{\# of hamps of }D_{0}\right)  -\left(  \text{\# of
hamps of }D^{\prime}\text{ that use the arc }wv\right) \\
&  \equiv\left(  \text{\# of hamps of }D^{\prime}\right)  \operatorname{mod}2,
\end{align*}
and the proof of the lemma will be complete.

So let us show this. Recall that $D$ is a tournament. Thus, the non-loop arcs
of $\overline{D}$ are precisely the arcs of $D^{\operatorname*{rev}}$ (by
Proposition \ref{prop.tour.rev-bar-def}). Hence, the non-loop arcs of
$\overline{D_{0}}$ are precisely the arcs of $D_{2}^{\operatorname*{rev}}$
(since $\overline{D_{0}}$ is just $\overline{D}$ with the extra arc $vw$
added, and since $D_{2}^{\operatorname*{rev}}$ is just $D^{\operatorname*{rev}%
}$ with the extra arc $vw$ added). Therefore, the digraphs $\overline{D_{0}}$
and $D_{2}^{\operatorname*{rev}}$ are equal \textquotedblleft up to
loops\textquotedblright\ (i.e., they have the same vertices and the same
non-loop arcs). Since loops don't matter for hamps, these two digraphs thus
have the same of hamps. Hence,%
\[
\left(  \text{\# of hamps in }\overline{D_{0}}\right)  =\left(  \text{\# of
hamps in }D_{2}^{\operatorname*{rev}}\right)  =\left(  \text{\# of hamps in
}D_{2}\right)
\]
(by Proposition \ref{prop.tour.berge-triv}), and therefore%
\[
\left(  \text{\# of hamps in }D_{2}\right)  =\left(  \text{\# of hamps in
}\overline{D_{0}}\right)  \equiv\left(  \text{\# of hamps in }D_{0}\right)
\operatorname{mod}2
\]
(by Theorem \ref{thm.hamp.Dbar}). As explained above, this completes the proof
of Lemma \ref{lem.tour.redei-lem}.
\end{proof}

Now, the Hard R\'{e}dei theorem has become easy:

\begin{proof}
[Proof of Theorem \ref{thm.tour.redei-hard}.](This is an outline; see
\cite[proof of Theorem 1.6.1]{17s-lec7} for more details.)

We need to prove that the \# of hamps of $D$ is odd. Lemma
\ref{lem.tour.redei-lem} tells us that the parity of this \# does not change
when we reverse a single arc of $D$. Thus, of course, if we reverse several
arcs of $D$, then this parity does not change either. However, we can WLOG
assume that the vertices of $D$ are $1,2,\ldots,n$ for some $n\in\mathbb{N}$,
and then, by reversing the appropriate arcs, we can ensure that the arcs of
$D$ are%
\begin{align*}
&  12,\ 13,\ 14,\ \ldots,\ 1n,\\
&  \ \ \ \ \ \ \text{\ }23,\ 24,\ \ldots,\ 2n,\\
&  \ \ \ \ \ \ \ \ \ \ \cdots,\\
&  \ \ \ \ \ \ \ \ \ \ \ \ \ \ \ \ \ \ \ \ \text{\ \ \ \ \ }\left(
n-1\right)  n
\end{align*}
(i.e., each arc of $D$ has the form $ij$ with $i<j$). But at this point, the
tournament $D$ has only one hamp: namely, $\left(  1,2,\ldots,n\right)  $. So
$\left(  \text{\# of hamps of }D\right)  =1$ is odd at this point. Since the
parity of the \# of hamps of $D$ has not changed as we reversed our arcs, we
thus conclude that it has always been odd. This proves the Hard R\'{e}dei
theorem (Theorem \ref{thm.tour.redei-hard}).
\end{proof}

As we already mentioned, the Easy R\'{e}dei theorem follows from the Hard
R\'{e}dei theorem. But it also has a short self-contained proof (\cite[Theorem
1.4.9]{17s-lec7}).

\begin{remark}
Theorem \ref{thm.tour.redei-hard} shows that the \# of hamps in a tournament
is an odd positive integer. Can it be any odd positive integer, or are certain
odd positive integers impossible?

Surprisingly, $7$ and $21$ are impossible. All other odd numbers between $1$
and $80555$ are possible. For higher numbers, the answer is not known so far.
See MathOverflow question \#232751 (\cite{MO232751}) for more details.
\end{remark}

\subsubsection{Hamiltonian cycles in tournaments}

By the Easy R\'{e}dei theorem, every tournament has a hamp. But of course, not
every tournament has a hamc\footnote{Recall that \textquotedblleft%
\textbf{hamc}\textquotedblright\ is our shorthand for \textquotedblleft
Hamiltonian cycle\textquotedblright.}. One obstruction is clear:

\begin{proposition}
If a digraph $D$ has a hamc, then $D$ is strongly connected.
\end{proposition}

In general, this is only a necessary criterion for a hamc, not a sufficient
one. Not every strongly connected digraph has a hamc. However, it turns out
that for tournaments, it is also sufficient, as long as the tournament has
enough vertices:

\begin{theorem}
[Camion's theorem]If a tournament $D$ is strongly connected and has at least
two vertices, then $D$ has a hamc.
\end{theorem}

\begin{proof}
[Proof sketch.]A detailed proof can be found in \cite[Theorem 1.5.5]%
{17s-lec7}; here is just a very rough sketch.

Let $D=\left(  V,A\right)  $ be a strongly connected tournament with at least
two vertices.\footnote{By the way, a tournament with exactly two vertices
cannot be strongly connected (as it has only $1$ arc). Thus, by requiring $D$
to have at least two vertices, we have actually guaranteed that $D$ has at
least three vertices.} We must show that $D$ has a hamc.

It is easy to see that $D$ has a cycle. Let $\mathbf{c}=\left(  v_{1}%
,v_{2},\ldots,v_{k},v_{1}\right)  $ be a cycle of maximum length. We shall
show that $\mathbf{c}$ is a hamc.

Let $C$ be the set $\left\{  v_{1},v_{2},\ldots,v_{k}\right\}  $ of all
vertices of this cycle $\mathbf{c}$.

A vertex $w\in V\setminus C$ will be called a \textbf{to-vertex} if there
exists an arc from some $v_{i}$ to $w$.

A vertex $w\in V\setminus C$ will be called a \textbf{from-vertex} if there
exists an arc from $w$ to some $v_{i}$.

Since $D$ is a tournament, each vertex in $V\setminus C$ is a to-vertex or a
from-vertex. In theory, a vertex could be both (having an arc from some
$v_{i}$ and also an arc to some other $v_{j}$). However, this does not
actually happen. To see why, argue as follows:

\begin{itemize}
\item If a to-vertex $w$ has an arc from some $v_{i}$, then it must also have
an arc from $v_{i+1}$\ \ \ \ \footnote{Here, indices are periodic modulo $k$,
so that $v_{k+1}$ means $v_{1}$.} (because otherwise there would be an arc
from $w$ to $v_{i+1}$, and then we could make our cycle $\mathbf{c}$ longer by
interjecting $w$ between $v_{i}$ and $v_{i+1}$; but this would contradict the
fact that $\mathbf{c}$ is a cycle of maximum length).

\item Iterating this argument, we see that if a to-vertex $w$ has an arc from
some $v_{i}$, then it must also have an arc from $v_{i+1}$, an arc from
$v_{i+2}$, an arc from $v_{i+3}$, and so on; i.e., it must have an arc from
each vertex of $\mathbf{c}$. Consequently, $w$ cannot be a from-vertex. This
shows that a to-vertex cannot be a from-vertex.
\end{itemize}

Let $F$ be the set of all from-vertices, and let $T$ be the set of all
to-vertices. Then, as we have just shown, $F$ and $T$ are disjoint. Moreover,
$F\cup T=V\setminus C$. Since a to-vertex cannot be a from-vertex, we
furthermore conclude that any to-vertex has an arc from each vertex of
$\mathbf{c}$ (otherwise, it would be a from-vertex), and that any from-vertex
has an arc to each vertex of $\mathbf{c}$ (otherwise, it would be a to-vertex).

Next, we argue that there cannot be an arc from a to-vertex $t$ to a
from-vertex $f$. Indeed, if there was such an arc, then we could make the
cycle $\mathbf{c}$ longer by interjecting $t$ and $f$ between (say) $v_{1}$
and $v_{2}$.

In total, we now know that every vertex of $D$ belongs to one of the three
disjoint sets $C$, $F$ and $T$, and furthermore there is no arc from $T$ to
$F$, no arc from $T$ to $C$, and no arc from $C$ to $F$. Thus, there exists no
walk from a vertex in $T$ to a vertex in $C$ (because there is no way out of
$T$). This would contradict the fact that $D$ is strongly connected, unless
the set $T$ is empty. Hence, $T$ must be empty. Similarly, $F$ must be empty.
Since $F\cup T=V\setminus C$, this entails that $V\setminus C$ is empty, so
that $V=C$. In other words, each vertex of $D$ is on our cycle $\mathbf{c}$.
Therefore, $\mathbf{c}$ is a hamc. This proves Camion's theorem.
\end{proof}

\subsubsection{Application of tournaments to the Vandermonde determinant}

To wrap up the topic of tournaments, let me briefly discuss a curious
application of their theory: a combinatorial proof of the Vandermonde
determinant formula. See \cite{17s-lec8} for the many details I'll be omitting.

Recall the Vandermonde determinant formula:

\Needspace{15pc}

\begin{theorem}
[Vandermonde determinant formula]Let $x_{1},x_{2},\ldots,x_{n}$ be $n$ numbers
(or, more generally, elements of a commutative ring). Consider the $n\times
n$-matrix%
\[
V:=\left(
\begin{array}
[c]{ccccc}%
1 & 1 & 1 & \cdots & 1\\
x_{1} & x_{2} & x_{3} & \cdots & x_{n}\\
x_{1}^{2} & x_{2}^{2} & x_{3}^{2} & \cdots & x_{n}^{2}\\
\vdots & \vdots & \vdots & \ddots & \vdots\\
x_{1}^{n-1} & x_{2}^{n-1} & x_{3}^{n-1} & \cdots & x_{n}^{n-1}%
\end{array}
\right)  =\left(  x_{j}^{i-1}\right)  _{1\leq i\leq n,\ 1\leq j\leq n}.
\]
Then, its determinant is
\[
\det V=\prod_{1\leq i<j\leq n}\left(  x_{j}-x_{i}\right)  .
\]

\end{theorem}

There are many simple proofs of this theorem (e.g., a few on
\href{https://proofwiki.org/wiki/Value_of_Vandermonde_Determinant/Formulation_1}{its
ProofWiki page}, which works with the transpose matrix). I will now outline a
combinatorial one, using tournaments. This proof goes back to Ira Gessel's
1979 paper \cite{Gessel79}.

First, how do $\det V$ and $\prod\limits_{1\leq i<j\leq n}\left(  x_{i}%
-x_{j}\right)  $ relate to tournaments?

As a warmup, let's assume that we have some number $y_{\left(  i,j\right)  }$
given for each pair $\left(  i,j\right)  $ of integers, and let's expand the
product%
\[
\left(  y_{\left(  1,2\right)  }+y_{\left(  2,1\right)  }\right)  \left(
y_{\left(  1,3\right)  }+y_{\left(  3,1\right)  }\right)  \left(  y_{\left(
2,3\right)  }+y_{\left(  3,2\right)  }\right)  .
\]
The result is a sum of $8$ products, one for each way to pluck an addend out
of each of the three little sums:%
\begin{align*}
&  \left(  y_{\left(  1,2\right)  }+y_{\left(  2,1\right)  }\right)  \left(
y_{\left(  1,3\right)  }+y_{\left(  3,1\right)  }\right)  \left(  y_{\left(
2,3\right)  }+y_{\left(  3,2\right)  }\right) \\
&  =y_{\left(  1,2\right)  }y_{\left(  1,3\right)  }y_{\left(  2,3\right)
}+y_{\left(  1,2\right)  }y_{\left(  1,3\right)  }y_{\left(  3,2\right)
}+y_{\left(  1,2\right)  }y_{\left(  3,1\right)  }y_{\left(  2,3\right)
}+y_{\left(  1,2\right)  }y_{\left(  3,1\right)  }y_{\left(  3,2\right)  }\\
&  \ \ \ \ \ \ \ \ \ \ +y_{\left(  2,1\right)  }y_{\left(  1,3\right)
}y_{\left(  2,3\right)  }+y_{\left(  2,1\right)  }y_{\left(  1,3\right)
}y_{\left(  3,2\right)  }+y_{\left(  2,1\right)  }y_{\left(  3,1\right)
}y_{\left(  2,3\right)  }+y_{\left(  2,1\right)  }y_{\left(  3,1\right)
}y_{\left(  3,2\right)  }.
\end{align*}
Note that each of the $8$ products obtained has the form $y_{a}y_{b}y_{c}$, where

\begin{itemize}
\item $a$ is one of the pairs $\left(  1,2\right)  $ and $\left(  2,1\right)
$,

\item $b$ is one of the pairs $\left(  1,3\right)  $ and $\left(  3,1\right)
$, and

\item $c$ is one of the pairs $\left(  2,3\right)  $ and $\left(  3,2\right)
$.
\end{itemize}

\noindent We can view these pairs $a$, $b$ and $c$ as the arcs of a tournament
with vertex set $\left\{  1,2,3\right\}  $. Thus, our above expansion can be
rewritten more compactly as follows:%
\begin{align*}
&  \left(  y_{\left(  1,2\right)  }+y_{\left(  2,1\right)  }\right)  \left(
y_{\left(  1,3\right)  }+y_{\left(  3,1\right)  }\right)  \left(  y_{\left(
2,3\right)  }+y_{\left(  3,2\right)  }\right) \\
&  =\sum_{\substack{D\text{ is a tournament}\\\text{with vertex set }\left\{
1,2,3\right\}  }}\ \ \prod_{\left(  i,j\right)  \text{ is an arc of }%
D}y_{\left(  i,j\right)  }.
\end{align*}
For reference, here are all the $8$ tournaments with vertex set $\left\{
1,2,3\right\}  $:%
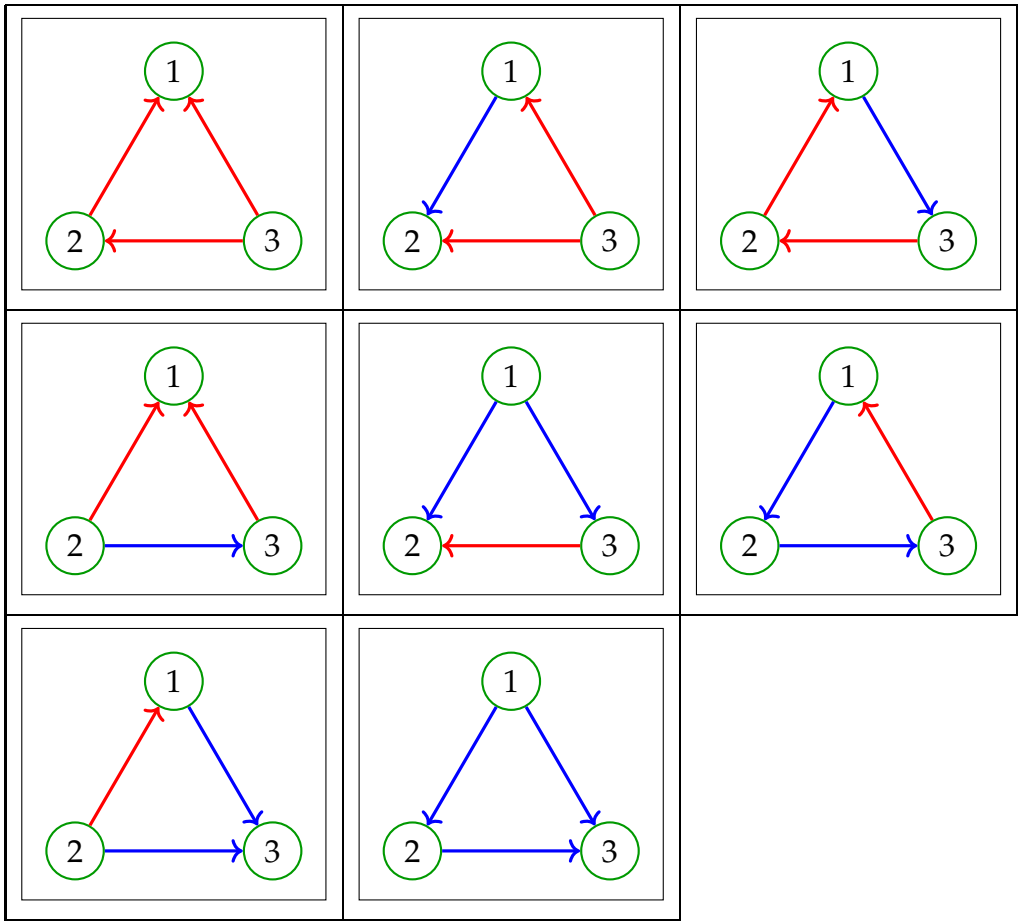
\[%

\ \ \ .
\]
Here, for convenience, we are drawing an arc $ij$ in blue if $i<j$ and in red otherwise.

This expansion can be generalized: We have%
\[
\prod_{1\leq i<j\leq n}\left(  y_{\left(  i,j\right)  }+y_{\left(  j,i\right)
}\right)  =\sum_{\substack{D\text{ is a tournament}\\\text{with vertex set
}\left\{  1,2,\ldots,n\right\}  }}\ \ \prod_{\left(  i,j\right)  \text{ is an
arc of }D}y_{\left(  i,j\right)  }.
\]
Substituting $y_{\left(  i,j\right)  }=%
\begin{cases}
x_{j}, & \text{if }i<j;\\
-x_{j}, & \text{if }i\geq j
\end{cases}
$ in this equality, we obtain%
\begin{align*}
\prod_{1\leq i<j\leq n}\left(  x_{j}-x_{i}\right)   &  =\sum
_{\substack{D\text{ is a tournament}\\\text{with vertex set }\left\{
1,2,\ldots,n\right\}  }}\ \ \underbrace{\prod_{\left(  i,j\right)  \text{ is
an arc of }D}%
\begin{cases}
x_{j}, & \text{if }i<j;\\
-x_{j}, & \text{if }i\geq j
\end{cases}
}_{\substack{=\left(  -1\right)  ^{\left(  \text{\# of red arcs of }D\right)
}\prod\limits_{j=1}^{n}x_{j}^{\deg^{-}j}\\\text{(where }\deg^{-}j\text{ means
the indegree of }j\text{ in }D\text{,}\\\text{and where the \textquotedblleft
red arcs\textquotedblright\ are the arcs }ij\text{ with }i>j\text{)}}}\\
&  =\sum_{\substack{D\text{ is a tournament}\\\text{with vertex set }\left\{
1,2,\ldots,n\right\}  }}\left(  -1\right)  ^{\left(  \text{\# of red arcs of
}D\right)  }\prod\limits_{j=1}^{n}x_{j}^{\deg^{-}j}.
\end{align*}
We shall refer to this sum as the \textquotedblleft big sum\textquotedblright.

On the other hand, if we let $S_{n}$ be the group of permutations of $\left\{
1,2,\ldots,n\right\}  $, and if we denote the sign of a permutation $\sigma$
by $\operatorname*{sign}\sigma$, then we have%
\[
\det V=\det\left(  V^{T}\right)  =\sum_{\sigma\in S_{n}}\operatorname*{sign}%
\sigma\cdot\prod_{j=1}^{n}x_{j}^{\sigma\left(  j\right)  -1}%
\]
(by the definition of a determinant). We shall refer to this sum as the
\textquotedblleft small sum\textquotedblright.

Our goal is to prove that the big sum equals the small sum. To prove this, we
must verify the following:

\begin{enumerate}
\item Each addend of the small sum is an addend of the big sum. Indeed, for
each permutation $\sigma\in S_{n}$, there is a certain tournament $T_{\sigma}$
that has
\[
\left(  -1\right)  ^{\left(  \text{\# of red arcs of }T_{\sigma}\right)
}\prod\limits_{j=1}^{n}x_{j}^{\deg^{-}j}=\operatorname*{sign}\sigma\cdot
\prod_{j=1}^{n}x_{j}^{\sigma\left(  j\right)  -1}.
\]
Can you find this $T_{\sigma}$ ?

\item All the addends of the big sum that are \textbf{not} addends of the
small sum cancel each other out. Why?

The basic idea is to argue that if a tournament $D$ appears in the big sum but
not in the small sum, then $D$ has a $3$-cycle (i.e., a cycle of length $3$).
When we reverse such a $3$-cycle (i.e., we reverse each of its arcs), the
indegrees of all vertices are preserved, but the sign $\left(  -1\right)
^{\left(  \text{\# of red arcs of }D\right)  }$ is flipped (since three arcs
change their orientation).

This suffices to show that for each addend that appears in the big sum but not
in the small sum, there is another addend with the same magnitude but with
opposite sign. Unfortunately, this in itself does not suffice to ensure that
all these addends cancel out; for example, the sum $1+1+1+\left(  -1\right)  $
has the same property but does not equal $0$. We need to show that the \# of
addends with positive sign (i.e., with $\left(  -1\right)  ^{\left(  \text{\#
of red arcs of }D\right)  }=1$) and a given magnitude equals the \# of addends
with negative sign (i.e., with $\left(  -1\right)  ^{\left(  \text{\# of red
arcs of }D\right)  }=-1$) and the same magnitude.

One way to achieve this would be by constructing a bijection (aka
\textquotedblleft perfect matching\textquotedblright)\ between the
\textquotedblleft positive\textquotedblright\ and the \textquotedblleft
negative\textquotedblright\ addends. This is tricky here: We would have to
decide \textbf{which} $3$-cycle to reverse (as there are usually many of
them), and this has to be done in a bijective way (so that two
\textquotedblleft positive\textquotedblright\ addends don't get assigned the
same \textquotedblleft negative\textquotedblright\ partner).

A less direct, but easier way is the following: Fix a positive integer $k$,
and consider only the tournaments with exactly $k$ many $3$-cycles. For each
such tournament, we can reverse any of its $k$ many $3$-cycles. It can be
shown (nice exercise!) that reversing the arcs of a $3$-cycle does not change
the \# of all $3$-cycles; thus, we don't accidentally change our $k$ in the
process. Thus, we find a \textquotedblleft$k$-to-$k$\textquotedblright%
\ correspondence between the \textquotedblleft positive\textquotedblright%
\ addends of a given magnitude and the \textquotedblleft
negative\textquotedblright\ addends of the same magnitude. As one can easily
see, this entails that the former and the latter are equinumerous, and thus
really cancel out. The addends that remain are exactly those in the small sum.
\end{enumerate}

As already mentioned, this is only a rough summary of the proof; the details
can be found in \cite{17s-lec8}.

\subsection{Exercises on tournaments}

There is, of course, much more to say about tournaments. See \cite{Moon13} for
a selection of topics. Let us merely hint at some possible directions by
giving a few exercises.

The next three exercises use the notion of a \textquotedblleft$3$%
-cycle\textquotedblright:

\begin{definition}
A $3$\textbf{-cycle} in a tournament $D=\left(  V,A\right)  $ means a triple
$\left(  u,v,w\right)  $ of vertices in $V$ such that all three pairs $\left(
u,v\right)  $, $\left(  v,w\right)  $ and $\left(  w,u\right)  $ belong to $A$.
\end{definition}

For example, the tournament shown in Example \ref{exe.tour.5-v-t} has the nine
different $3$-cycles%
\begin{align*}
&  \left(  1,4,3\right)  ,\ \ \ \ \ \ \ \ \ \ \left(  1,5,3\right)
,\ \ \ \ \ \ \ \ \ \ \left(  2,5,3\right)  ,\ \ \ \ \ \ \ \ \ \ \left(
3,1,4\right)  ,\ \ \ \ \ \ \ \ \ \ \left(  3,1,5\right)  ,\\
&  \left(  3,2,5\right)  ,\ \ \ \ \ \ \ \ \ \ \left(  4,3,1\right)
,\ \ \ \ \ \ \ \ \ \ \left(  5,3,1\right)  ,\ \ \ \ \ \ \ \ \ \ \left(
5,3,2\right)  .
\end{align*}
(Yes, we are counting a $3$-cycle $\left(  u,v,w\right)  $ as being distinct
from $\left(  v,w,u\right)  $ and $\left(  w,u,v\right)  $.)

\begin{exercise}
\label{exe.tournament.num-3-cycs}Let $D=\left(  V,A\right)  $ be a tournament.
Set $n=\left\vert V\right\vert $ and $m=\sum_{v\in V}\dbinom{\deg^{-}\left(
v\right)  }{2}$.

\begin{enumerate}
\item[\textbf{(a)}] Show that $m=\sum_{v\in V}\dbinom{\deg^{+}\left(
v\right)  }{2}$.

\item[\textbf{(b)}] Show that the number of $3$-cycles in $D$ is $3\left(
\dbinom{n}{3}-m\right)  $.\medskip
\end{enumerate}

[\textbf{Solution:} This is Exercise 5 on homework set \#2 from my Spring 2017
course; see \href{https://www.cip.ifi.lmu.de/~grinberg/t/17s/}{the course
page} for solutions.]
\end{exercise}

The next exercise uses the notation $\deg_{D}^{-}v$ for the indegree of a
vertex $v$ in a digraph $D$. (We usually denote this by $\deg^{-}v$, but
sometimes it is important to stress the dependence on $D$, since $v$ can be a
vertex of two different digraphs.)

\begin{exercise}
\label{exe.tournament.3-cyc-rev}If a tournament $D$ has a $3$-cycle $\left(
u,v,w\right)  $, then we can define a new tournament $D_{u,v,w}^{\prime}$ as
follows: The vertices of $D_{u,v,w}^{\prime}$ shall be the same as those of
$D$. The arcs of $D_{u,v,w}^{\prime}$ shall be the same as those of $D$,
except that the three arcs $\left(  u,v\right)  $, $\left(  v,w\right)  $ and
$\left(  w,u\right)  $ are replaced by the three new arcs $\left(  v,u\right)
$, $\left(  w,v\right)  $ and $\left(  u,w\right)  $. (Visually speaking,
$D_{u,v,w}^{\prime}$ is obtained from $D$ by turning the arrows on the arcs
$\left(  u,v\right)  $, $\left(  v,w\right)  $ and $\left(  w,u\right)  $
around.) We say that the new tournament $D_{u,v,w}^{\prime}$ is obtained from
the old tournament $D$ by a $3$\textbf{-cycle reversal operation}.

Now, let $V$ be a finite set, and let $E$ and $F$ be two tournaments with
vertex set $V$. Prove that $F$ can be obtained from $E$ by a sequence of
$3$-cycle reversal operations if and only if each $v\in V$ satisfies $\deg
_{E}^{-}\left(  v\right)  =\deg_{F}^{-}\left(  v\right)  $. (Note that a
sequence may be empty, which allows handling the case $E=F$ even if $E$ has no
$3$-cycles to reverse.)\medskip

[\textbf{Solution:} This is Exercise 6 on homework set \#2 from my Spring 2017
course; see \href{https://www.cip.ifi.lmu.de/~grinberg/t/17s/}{the course
page} for solutions.]
\end{exercise}

\begin{exercise}
\label{exe.tournament.2-path-rev}A tournament $D=\left(  V,A\right)  $ is
called \textbf{transitive} if it has no $3$-cycles.

If a tournament $D=\left(  V,A\right)  $ has three distinct vertices $u$, $v$
and $w$ satisfying $\left(  u,v\right)  \in A$ and $\left(  v,w\right)  \in
A$, then we can define a new tournament $D_{u,v,w}^{\prime\prime}$ as follows:
The vertices of $D_{u,v,w}^{\prime\prime}$ shall be the same as those of $D$.
The arcs of $D_{u,v,w}^{\prime\prime}$ shall be the same as those of $D$,
except that the two arcs $\left(  u,v\right)  $ and $\left(  v,w\right)  $ are
replaced by the two new arcs $\left(  v,u\right)  $ and $\left(  w,v\right)
$. We say that the new tournament $D_{u,v,w}^{\prime\prime}$ is obtained from
the old tournament $D$ by a $2$\textbf{-path reversal operation}.

Let $D$ be any tournament. Prove that there is a sequence of $2$-path reversal
operations that transforms $D$ into a transitive tournament.\medskip

[\textbf{Solution:} This is Exercise 7 on homework set \#2 from my Spring 2017
course; see \href{https://www.cip.ifi.lmu.de/~grinberg/t/17s/}{the course
page} for solutions.]
\end{exercise}

\section{Trees and arborescences}

Trees are particularly nice graphs. Among other things, they can be
characterized as

\begin{itemize}
\item the minimal connected graphs on a given set of vertices, or

\item the maximal acyclic (= having no cycles) graphs on a given set of
vertices, or

\item in many other ways.
\end{itemize}

Arborescences are their closest analogue for digraphs.

In this chapter, we will discuss the theory of trees and some of their
applications. Further applications are usually covered in courses in
theoretical computer science, but their notion of a tree is somewhat different
from ours.

\subsection{Some general properties of components and cycles}

\subsubsection{Backtrack-free walks revisited}

Before we start with trees, let us recall and prove some more facts about
general multigraphs. Recall the notion of a \textquotedblleft backtrack-free
walk\textquotedblright\ that already had a brief appearance in the proof of
Theorem \ref{thm.cyc.two-paths-cyc}:

\begin{definition}
\label{def.btf-walk}Let $G$ be a multigraph. A \textbf{backtrack-free walk} of
$G$ means a walk $\mathbf{w}$ such that no two consecutive edges of
$\mathbf{w}$ are identical.
\end{definition}

Here are a few properties of this notion:

\begin{proposition}
\label{prop.btf-walk.cyc}Let $G$ be a multigraph. Let $\mathbf{w}$ be a
backtrack-free walk of $G$. Then, $\mathbf{w}$ either is a path or contains a cycle.
\end{proposition}

\begin{proof}
We have already proved this for simple graphs (in Proposition
\ref{prop.cyc.btf-walk-cyc}). More or less the same argument works for
multigraphs. (\textquotedblleft More or less\textquotedblright\ because the
definition of a cycle in a multigraph is slightly different from that in a
simple graph; but the proof is easy to adapt.)
\end{proof}

\begin{theorem}
\label{thm.btf-walk.two-cyc}Let $G$ be a multigraph. Let $u$ and $v$ be two
vertices of $G$. Assume that there are two distinct backtrack-free walks from
$u$ to $v$ in $G$. Then, $G$ has a cycle.
\end{theorem}

\begin{proof}
We have already proved this for simple graphs (Claim 1 in the proof of Theorem
\ref{thm.cyc.two-paths-cyc}). More or less the same argument works for multigraphs.
\end{proof}

\subsubsection{Counting components}

Next, we shall derive a few properties of the number of components of a graph.
Again, we have already done most of the hard work, and we can now derive
corollaries. First, we give this number a name:

\begin{definition}
Let $G$ be a multigraph. Then, $\operatorname*{conn}G$ means the number of
components of $G$. (Some authors also call this number $b_{0}\left(  G\right)
$. This notation comes from algebraic topology, where it stands for the $0$-th
Betti number. This makes sense, because we can regard a multigraph $G$ as a
topological space. But we won't need this.)
\end{definition}

So a multigraph $G$ satisfies $\operatorname*{conn}G=1$ if and only if $G$ is
connected. Moreover, $\operatorname*{conn}G=0$ if and only if $G$ has no vertices.

Let us next recall Definition \ref{def.mg.G-e} and Theorem
\ref{thm.mg.G-e.conn} (which is an analogue of Theorem \ref{thm.G-e.conn} and
can be proved in more or less the same way). As a consequence of the latter
theorem, we obtain the following:

\begin{corollary}
\label{cor.conn.+1}Let $G$ be a multigraph. Let $e$ be an edge of $G$. Then:

\begin{enumerate}
\item[\textbf{(a)}] If $e$ is an edge of some cycle of $G$, then
$\operatorname*{conn}\left(  G\setminus e\right)  =\operatorname*{conn}G$.

\item[\textbf{(b)}] If $e$ appears in no cycle of $G$, then
$\operatorname*{conn}\left(  G\setminus e\right)  =\operatorname*{conn}G+1$.

\item[\textbf{(c)}] In either case, we have $\operatorname*{conn}\left(
G\setminus e\right)  \leq\operatorname*{conn}G+1$.
\end{enumerate}
\end{corollary}

\begin{proof}
Part \textbf{(a)} follows from Theorem \ref{thm.mg.G-e.conn} \textbf{(a)}.
Part \textbf{(b)} follows from Theorem \ref{thm.mg.G-e.conn} \textbf{(b)}.
Part \textbf{(c)} follows by combining parts \textbf{(a)} and \textbf{(b)}.
\end{proof}

\begin{corollary}
\label{cor.conn.geq}Let $G=\left(  V,E,\varphi\right)  $ be a multigraph.
Then, $\operatorname*{conn}G\geq\left\vert V\right\vert -\left\vert
E\right\vert $.
\end{corollary}

\begin{proof}
We induct on $\left\vert E\right\vert $:

\textit{Base case:} If $\left\vert E\right\vert =0$, then
$\operatorname*{conn}G=\left\vert V\right\vert $ (since $\left\vert
E\right\vert =0$ means that the graph $G$ has no edges, and thus no two
distinct vertices are path-connected); but this rewrites as
$\operatorname*{conn}G=\left\vert V\right\vert -\left\vert E\right\vert $
(since $\left\vert E\right\vert =0$). Thus, Corollary \ref{cor.conn.geq} is
proved for $\left\vert E\right\vert =0$.

\textit{Induction step:} Let $k\in\mathbb{N}$. Assume (as the induction
hypothesis) that Corollary \ref{cor.conn.geq} holds for $\left\vert
E\right\vert =k$. We must now show that it also holds for $\left\vert
E\right\vert =k+1$.

So let us consider a multigraph $G=\left(  V,E,\varphi\right)  $ with
$\left\vert E\right\vert =k+1$. Thus, $\left\vert E\right\vert -1=k$. Pick any
edge $e\in E$ (such an edge exists, since $\left\vert E\right\vert
=k+1\geq1>0$). Then, the multigraph $G\setminus e$ has edge set $E\setminus
\left\{  e\right\}  $ and therefore has $\left\vert E\setminus\left\{
e\right\}  \right\vert =\left\vert E\right\vert -1=k$ many edges. Hence, by
the induction hypothesis, we have
\[
\operatorname*{conn}\left(  G\setminus e\right)  \geq\left\vert V\right\vert
-\left\vert E\setminus\left\{  e\right\}  \right\vert
\]
(since $G\setminus e$ is a multigraph with vertex set $V$ and edge set
$E\setminus\left\{  e\right\}  $). However, Corollary \ref{cor.conn.+1}
\textbf{(c)} yields $\operatorname*{conn}\left(  G\setminus e\right)
\leq\operatorname*{conn}G+1$. Thus,%
\begin{align*}
\operatorname*{conn}G  &  \geq\underbrace{\operatorname*{conn}\left(
G\setminus e\right)  }_{\geq\left\vert V\right\vert -\left\vert E\setminus
\left\{  e\right\}  \right\vert }-\,1\geq\left\vert V\right\vert
-\underbrace{\left\vert E\setminus\left\{  e\right\}  \right\vert
}_{=\left\vert E\right\vert -1}-\,1\\
&  =\left\vert V\right\vert -\left(  \left\vert E\right\vert -1\right)
-1=\left\vert V\right\vert -\left\vert E\right\vert .
\end{align*}
This completes the induction step. Thus, Corollary \ref{cor.conn.geq} is proven.
\end{proof}

\begin{corollary}
\label{cor.conn.eq}Let $G=\left(  V,E,\varphi\right)  $ be a multigraph that
has no cycles. Then, $\operatorname*{conn}G=\left\vert V\right\vert
-\left\vert E\right\vert $.
\end{corollary}

\begin{proof}
Replay the proof of Corollary \ref{cor.conn.geq}, with just a few changes:
Instead of applying Corollary \ref{cor.conn.+1} \textbf{(c)}, apply Corollary
\ref{cor.conn.+1} \textbf{(b)} (this is allowed because $G$ has no cycles and
thus $e$ appears in no cycle of $G$). The induction hypothesis can be used
because when $G$ has no cycles, $G\setminus e$ has no cycles either. All
$\leq$ and $\geq$ signs in the above proof now can be replaced by $=$ signs
(since Corollary \ref{cor.conn.+1} \textbf{(b)} claims an equality, not an
inequality). The result is therefore $\operatorname*{conn}G=\left\vert
V\right\vert -\left\vert E\right\vert $.
\end{proof}

\begin{corollary}
\label{cor.conn.geq+1}Let $G=\left(  V,E,\varphi\right)  $ be a multigraph
that has at least one cycle. Then, $\operatorname*{conn}G\geq\left\vert
V\right\vert -\left\vert E\right\vert +1$.
\end{corollary}

\begin{proof}
Pick an edge $e\in E$ that belongs to some cycle (such an edge exists, since
$G$ has at least one cycle). Then, Corollary \ref{cor.conn.+1} \textbf{(a)}
yields $\operatorname*{conn}\left(  G\setminus e\right)  =\operatorname*{conn}%
G$. However, Corollary \ref{cor.conn.geq} (applied to $G\setminus e$ and
$E\setminus\left\{  e\right\}  $ instead of $G$ and $E$) yields
\[
\operatorname*{conn}\left(  G\setminus e\right)  \geq\left\vert V\right\vert
-\underbrace{\left\vert E\setminus\left\{  e\right\}  \right\vert
}_{=\left\vert E\right\vert -1}=\left\vert V\right\vert -\left(  \left\vert
E\right\vert -1\right)  =\left\vert V\right\vert -\left\vert E\right\vert +1.
\]
Since $\operatorname*{conn}\left(  G\setminus e\right)  =\operatorname*{conn}%
G$, this rewrites as $\operatorname*{conn}G\geq\left\vert V\right\vert
-\left\vert E\right\vert +1$.
\end{proof}

We summarize what we have proved into one convenient theorem:

\begin{theorem}
\label{thm.conn.conn-V-E}Let $G=\left(  V,E,\varphi\right)  $ be a multigraph. Then:

\begin{enumerate}
\item[\textbf{(a)}] We always have $\operatorname*{conn}G\geq\left\vert
V\right\vert -\left\vert E\right\vert $.

\item[\textbf{(b)}] We have $\operatorname*{conn}G=\left\vert V\right\vert
-\left\vert E\right\vert $ if and only if $G$ has no cycles.
\end{enumerate}
\end{theorem}

\begin{proof}
\textbf{(a)} This is Corollary \ref{cor.conn.geq}. \medskip

\textbf{(b)} $\Longleftarrow:$ This is Corollary \ref{cor.conn.eq}.

$\Longrightarrow:$ Assume that $\operatorname*{conn}G=\left\vert V\right\vert
-\left\vert E\right\vert $. If $G$ had any cycles, then Corollary
\ref{cor.conn.geq+1} would yield $\operatorname*{conn}G\geq\left\vert
V\right\vert -\left\vert E\right\vert +1>\left\vert V\right\vert -\left\vert
E\right\vert $, which would contradict $\operatorname*{conn}G=\left\vert
V\right\vert -\left\vert E\right\vert $. So $G$ has no cycles. This proves the
\textquotedblleft$\Longrightarrow$\textquotedblright\ direction of Theorem
\ref{thm.conn.conn-V-E}.
\end{proof}

\begin{remark}
\label{rmk.conn.cyclomatic-num}Let $G=\left(  V,E,\varphi\right)  $ be a
multigraph. The number%
\[
\operatorname*{conn}G-\left(  \left\vert V\right\vert -\left\vert E\right\vert
\right)
\]
is known as \href{https://en.wikipedia.org/wiki/Circuit_rank}{the
\textbf{circuit rank} or the \textbf{cyclomatic number}} of $G$. By Theorem
\ref{thm.conn.conn-V-E} \textbf{(a)}, this number is always nonnegative; by
Theorem \ref{thm.conn.conn-V-E} \textbf{(b)}, it equals $0$ if and only if $G$
has no cycles.

One might optimistically hope that this number counts the cycles of $G$. But
this is not true (no matter whether we count reversals and cyclic rotations of
a cycle as being distinct or as being equal). For example, a complete graph
$K_{n}$ (with $n\geq1$) has cyclomatic number $1-\left(  n-\dbinom{n}%
{2}\right)  $, but it usually has many more cycles than this.

However, the cyclomatic number $\operatorname*{conn}G-\left(  \left\vert
V\right\vert -\left\vert E\right\vert \right)  $ does have a meaning, albeit a
subtler one: It is the dimension of a certain vector space that, in some way,
consists of cycles. See \cite[\S 2.1]{Berge91} for details.
\end{remark}

\subsection{Forests and trees}

\subsubsection{Definitions}

We now introduce two of the heroes of this chapter:

\begin{definition}
\label{def.trees.forest}A \textbf{forest} is a multigraph with no cycles.

(In particular, a forest therefore cannot contain two distinct parallel edges.
It also cannot contain loops.)
\end{definition}

\begin{definition}
\label{def.trees.tree}A \textbf{tree} is a connected forest.
\end{definition}

\Needspace{25pc}

\begin{example}
\label{exa.trees.exa1}Consider the following multigraphs:%
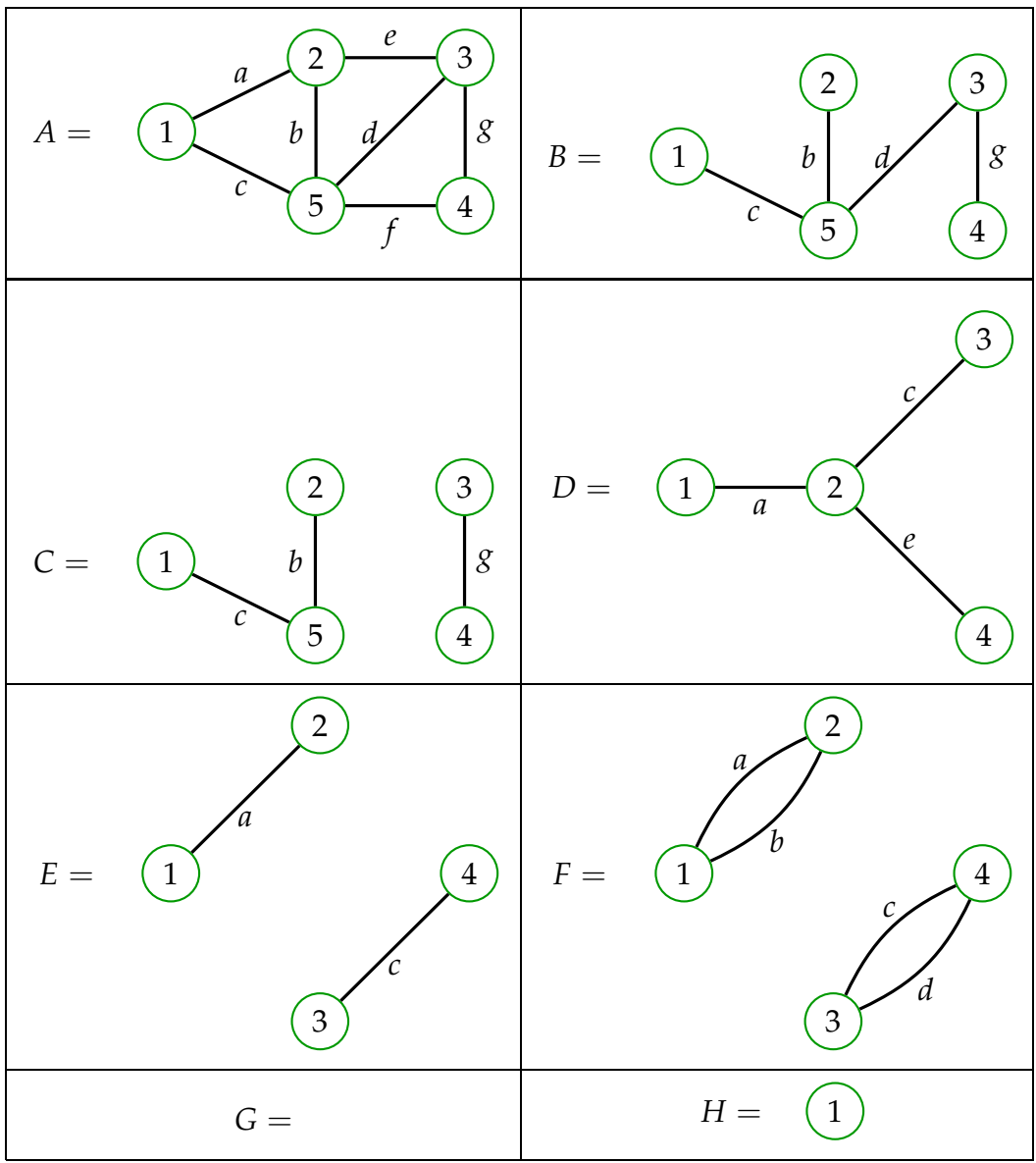
\[%

\ \ \ \ \ \ .
\]
(Yes, $G$ is an empty graph with no vertices.) Which of them are forests, and
which are trees?

\begin{itemize}
\item The graph $A$ is not a forest, since it has a cycle (actually, several
cycles). Thus, $A$ is not a tree either.

\item The graph $B$ is a tree.

\item The graph $C$ is a forest, but not a tree, since it is not connected.

\item The graph $D$ is a tree.

\item The graph $E$ is a forest, but not a tree.

\item The graph $F$ is not a forest, since it has cycles.

\item The graph $G$ (which has no vertices and no edges) is a forest, but not
a tree, since it is not connected (recall: a graph is connected if it has $1$
component; but $G$ has $0$ components).

\item The graph $H$ is a tree.
\end{itemize}
\end{example}

\subsubsection{The tree equivalence theorem}

Trees can be described in many ways:

\begin{theorem}
[The tree equivalence theorem]\label{thm.trees.T1-8}Let $G=\left(
V,E,\varphi\right)  $ be a multigraph. Then, the following eight statements
are equivalent:

\begin{itemize}
\item \textbf{Statement T1:} The multigraph $G$ is a tree.

\item \textbf{Statement T2:} The multigraph $G$ has no loops, and we have
$V\neq\varnothing$, and for each $u\in V$ and $v\in V$, there is a
\textbf{unique} path from $u$ to $v$.

\item \textbf{Statement T3:} We have $V\neq\varnothing$, and for each $u\in V$
and $v\in V$, there is a \textbf{unique} backtrack-free walk from $u$ to $v$.

\item \textbf{Statement T4:} The multigraph $G$ is connected, and we have
$\left\vert E\right\vert =\left\vert V\right\vert -1$.

\item \textbf{Statement T5:} The multigraph $G$ is connected, and we have
$\left\vert E\right\vert <\left\vert V\right\vert $.

\item \textbf{Statement T6:} We have $V\neq\varnothing$, and the graph $G$ is
a forest, but adding any new edge to $G$ creates a cycle.

\item \textbf{Statement T7:} The multigraph $G$ is connected, but removing any
edge from $G$ yields a disconnected (i.e., non-connected) graph.

\item \textbf{Statement T8:} The multigraph $G$ is a forest, and we have
$\left\vert E\right\vert \geq\left\vert V\right\vert -1$ and $V\neq
\varnothing$.
\end{itemize}
\end{theorem}

\begin{proof}
We shall prove the following implications:%
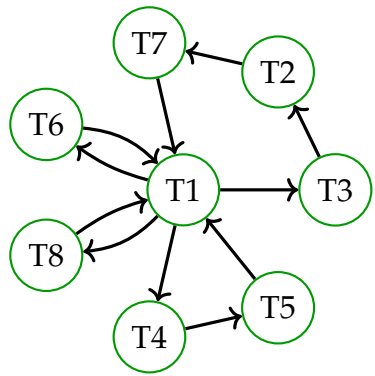
\[%
\begin{tikzpicture}
\begin{scope}[every node/.style={circle,thick,draw=green!60!black}]
\node(3) at (0:2) {T3};;
\node(2) at (360/7:2) {T2};
\node(7) at (2*360/7:2) {T7};
\node(6) at (3*360/7:2) {T6};
\node(8) at (4*360/7:2) {T8};
\node(4) at (5*360/7:2) {T4};
\node(5) at (6*360/7:2) {T5};
\node(1) at (0,0) {T1};
\end{scope}
\begin{scope}[every edge/.style={draw=black,very thick}]
\path[->] (1) edge (3) (3) edge (2) (2) edge (7) (7) edge (1);
\path[->] (1) edge[bend left=10] (6) (6) edge[bend left=20] (1);
\path[->] (1) edge[bend left=20] (8) (8) edge[bend left=10] (1);
\path[->] (1) edge (4) (4) edge (5) (5) edge (1);
\end{scope}
\end{tikzpicture}%
\ \ .
\]
In this digraph, an arc from T$i$ to T$j$ stands for the implication
T$i\Longrightarrow$T$j$. Since this digraph is strongly connected (i.e., you
can travel from Statement T$i$ to Statement T$j$ along its arcs for any
$i,j$), this will prove the theorem. So let us prove the implications.
\medskip

\textit{Proof of T1}$\Longrightarrow$\textit{T3:} Assume that Statement T1
holds. Thus, $G$ is a tree. Therefore, $G$ is connected, so that
$V\neq\varnothing$. We must prove that for each $u\in V$ and $v\in V$, there
is a \textbf{unique} backtrack-free walk from $u$ to $v$. The existence of
such a walk is clear (since $G$ is connected, so there is a path from $u$ to
$v$). Thus, we only need to show that it is unique. But this is easy: If there
were two distinct backtrack-free walks from $u$ to $v$ (for some $u\in V$ and
$v\in V$), then Theorem \ref{thm.btf-walk.two-cyc} would show that $G$ has a
cycle, and thus $G$ could not be a forest, let alone a tree. Thus, the
backtrack-free walk from $u$ to $v$ is unique. So we have proved Statement T3.
The implication T1$\Longrightarrow$T3 is thus proved. \medskip

\textit{Proof of T3}$\Longrightarrow$\textit{T2:} Assume that Statement T3
holds. We must prove that Statement T2 holds. First, $G$ has no loops, because
if there was a loop $e$ with endpoint $u$, then the two walks $\left(
u\right)  $ and $\left(  u,e,u\right)  $ would be two distinct backtrack-free
walks from $u$ to $u$. It remains to prove that for each each $u\in V$ and
$v\in V$, there is a \textbf{unique} path from $u$ to $v$. However, the
existence of a walk from $u$ to $v$ always implies the existence of a path
from $u$ to $v$ (by Corollary \ref{cor.mg.walk-thus-path}). Moreover, the
uniqueness of a backtrack-free walk from $u$ to $v$ implies the uniqueness of
a path from $u$ to $v$ (since any path is a backtrack-free walk). Thus,
Statement T2 follows from Statement T3. \medskip

\textit{Proof of T2}$\Longrightarrow$\textit{T7:} Assume that Statement T2
holds. Then, $G$ is connected. Now, let us remove any edge $e$ from $G$. Let
$u$ and $v$ be the endpoints of $e$. Then, $u\neq v$ (since $G$ has no loops).
There cannot be a path from $u$ to $v$ in the graph $G\setminus e$ (because if
there was such a path, then it would also be a path from $u$ to $v$ in the
graph $G$, and this path would be distinct from the path $\left(
u,e,v\right)  $; thus, the graph $G$ would have at least two paths from $u$ to
$v$; but this would contradict the uniqueness part of Statement T2). Hence,
the graph $G\setminus e$ is disconnected. So we have shown that $G$ is
connected, but removing any edge from $G$ yields a disconnected graph. In
other words, Statement T7 holds. \medskip

\textit{Proof of T7}$\Longrightarrow$\textit{T1:} Assume that Statement T7
holds. We must show that $G$ is a tree. Since $G$ is connected (by Statement
T7), it suffices to show that $G$ is a forest, i.e., that $G$ has no cycles.
However, if $G$ had any cycle, then we could pick any edge $e$ of this cycle,
and then we would know that $G\setminus e$ is still connected (since Corollary
\ref{cor.conn.+1} \textbf{(a)} would yield $\operatorname*{conn}\left(
G\setminus e\right)  =\operatorname*{conn}G=1$), and this would contradict
Statement T7. Thus, $G$ has no cycles, hence is a forest. This proves
Statement T1. \medskip

\textit{Proof of T1}$\Longrightarrow$\textit{T6:} Assume that Statement T1
holds. Thus, $G$ is a tree. We must show that adding any new edge to $G$
creates a cycle (since all other parts of Statement T6 are clear).

Indeed, let us add a new edge $f$ to $G$. Let $u$ and $v$ be the endpoints of
$f$. The graph $G$ is connected, so there is already a path from $u$ to $v$ in
$G$. Combining this path with the edge $f$, we obtain a cycle. Thus, the graph
obtained from $G$ by adding the new edge $f$ has a cycle. This completes our
proof that Statement T6 holds. \medskip

\textit{Proof of T6}$\Longrightarrow$\textit{T1:} Assume that Statement T6
holds. Thus, $G$ is a forest. We must only show that $G$ is connected.

Assume the contrary. Thus, there exist two vertices $u$ and $v$ of $G$ that
are not path-connected in $G$. Hence, adding a new edge $f$ with endpoints $u$
and $v$ to the graph $G$ cannot create a new cycle (because any such cycle
would have to contain $f$ (otherwise, it would already be a cycle of $G$, but
$G$ has no cycles), and then we could remove $f$ from it to obtain a path from
$u$ to $v$ in $G$; but such a path cannot exist, since $u$ and $v$ are not
path-connected in $G$). This contradicts Statement T6.

So we have shown that $G$ is connected, and thus $G$ is a tree. This proves
Statement T1. \medskip

\textit{Proof of T1}$\Longrightarrow$\textit{T8:} Assume that Statement T1
holds. So $G$ is a tree. Clearly, $G$ is then a forest. We must show that
$\left\vert E\right\vert \geq\left\vert V\right\vert -1$.

Theorem \ref{thm.conn.conn-V-E} \textbf{(a)} yields $\operatorname*{conn}%
G\geq\left\vert V\right\vert -\left\vert E\right\vert $. But we have
$\operatorname*{conn}G=1$ because $G$ is connected. Thus,
$1=\operatorname*{conn}G\geq\left\vert V\right\vert -\left\vert E\right\vert
$. In other words, $\left\vert E\right\vert \geq\left\vert V\right\vert -1$.
This proves Statement T8. \medskip

\textit{Proof of T8}$\Longrightarrow$\textit{T1:} Assume that Statement T8
holds. Thus, $G$ is a forest. We must only show that $G$ is connected.
However, $G$ is a forest, and thus has no cycles. Hence, Theorem
\ref{thm.conn.conn-V-E} \textbf{(b)} yields $\operatorname*{conn}G=\left\vert
V\right\vert -\left\vert E\right\vert \leq1$ (since Statement 8 yields
$\left\vert E\right\vert \geq\left\vert V\right\vert -1$). On the other hand,
$\operatorname*{conn}G\geq1$ (since $V\neq\varnothing$). Combining these two
inequalities, we obtain $\operatorname*{conn}G=1$. In other words, $G$ is
connected. This yields Statement T1 (since $G$ is a forest). \medskip

\textit{Proof of T1}$\Longrightarrow$\textit{T4:} Assume that Statement T1
holds. Then, $G$ is a tree, hence a connected forest. Therefore, $G$ has no
cycles (by the definition of a forest). Theorem \ref{thm.conn.conn-V-E}
\textbf{(b)} therefore yields $\operatorname*{conn}G=\left\vert V\right\vert
-\left\vert E\right\vert $. Thus, $\left\vert V\right\vert -\left\vert
E\right\vert =\operatorname*{conn}G=1$ (since $G$ is connected), so that
$\left\vert E\right\vert =\left\vert V\right\vert -1$. Thus, Statement T4 is
proved. \medskip

\textit{Proof of T4}$\Longrightarrow$\textit{T5:} The implication
T4$\Longrightarrow$T5 is obvious. \medskip

\textit{Proof of T5}$\Longrightarrow$\textit{T1:} Assume that Statement T5
holds. Thus, the multigraph $G$ is connected, and we have $\left\vert
E\right\vert <\left\vert V\right\vert $. Thus, $\left\vert E\right\vert
\leq\left\vert V\right\vert -1$. In other words, $1\leq\left\vert V\right\vert
-\left\vert E\right\vert $. Since $G$ is connected, we have
$\operatorname*{conn}G=1\leq\left\vert V\right\vert -\left\vert E\right\vert
$. However, Theorem \ref{thm.conn.conn-V-E} \textbf{(a)} yields
$\operatorname*{conn}G\geq\left\vert V\right\vert -\left\vert E\right\vert $.
Combining these two inequalities, we obtain $\operatorname*{conn}G=\left\vert
V\right\vert -\left\vert E\right\vert $. Thus, Theorem \ref{thm.conn.conn-V-E}
\textbf{(b)} shows that $G$ has no cycles. In other words, $G$ is a forest.
Hence, $G$ is a tree (since $G$ is connected). This proves Statement T1.
\medskip

We have now proved all necessary implications to conclude that all eight
statements T1, T2, $\ldots$, T8 are equivalent. Theorem \ref{thm.trees.T1-8}
is thus proved.
\end{proof}

We also observe the following connection between trees and forests:

\begin{proposition}
\label{prop.tree.forest=conn=trees}Let $G$ be a multigraph, and let
$C_{1},C_{2},\ldots,C_{k}$ be its components. Then, $G$ is a forest if and
only if all the induced subgraphs $G\left[  C_{1}\right]  ,G\left[
C_{2}\right]  ,\ldots,G\left[  C_{k}\right]  $ are trees.
\end{proposition}

\begin{proof}
$\Longrightarrow:$ Assume that $G$ is a forest. Thus, $G$ has no cycles.
Hence, the induced subgraphs $G\left[  C_{1}\right]  ,G\left[  C_{2}\right]
,\ldots,G\left[  C_{k}\right]  $ have no cycles either (since a cycle in any
of them would be a cycle of $G$); in other words, they are forests. But they
are furthermore connected (since the induced subgraph on a component is always
connected\footnote{This is Proposition \ref{prop.mg.conn.induced-is-conn}.}).
Hence, they are connected forests, i.e., trees. \medskip

$\Longleftarrow:$ Assume that the induced subgraphs $G\left[  C_{1}\right]
,G\left[  C_{2}\right]  ,\ldots,G\left[  C_{k}\right]  $ are trees. Hence,
none of them has a cycle. Thus, $G$ has no cycles either (since a cycle of $G$
would have to be fully contained in one of these induced
subgraphs\footnote{Indeed, if it wasn't, then it would contain vertices from
different components. But this is impossible, since there are no walks between
vertices in different components.}). In other words, $G$ is a forest.
\end{proof}

\subsubsection{Summary}

Let us briefly summarize some properties of trees:

If $T=\left(  V,E,\varphi\right)  $ is a tree, then...

\begin{itemize}
\item $T$ is a connected forest. (This is how trees were defined.) Thus, $T$
has no cycles. (This is how forests were defined.)

\item we have $\left\vert E\right\vert =\left\vert V\right\vert -1$. (This
follows from the implication T1$\Longrightarrow$T4 in Theorem
\ref{thm.trees.T1-8}.)

\item adding any new edge to $T$ creates a cycle. (This follows from the
implication T1$\Longrightarrow$T6 in Theorem \ref{thm.trees.T1-8}.)

\item removing any edge from $T$ yields a disconnected (i.e., non-connected)
graph. (This follows from the implication T1$\Longrightarrow$T7 in Theorem
\ref{thm.trees.T1-8}.)

\item for each $u\in V$ and $v\in V$, there is a \textbf{unique}
backtrack-free walk from $u$ to $v$. (This follows from the implication
T1$\Longrightarrow$T3 in Theorem \ref{thm.trees.T1-8}.) Moreover, this
backtrack-free walk is a path (since any walk from $u$ to $v$ contains a path
from $u$ to $v$).
\end{itemize}

\begin{remark}
Computer scientists use some notions of \textquotedblleft
trees\textquotedblright\ that are similar to ours, but not quite the same. In
particular, their trees often have \textbf{roots} (i.e., one vertex is chosen
to be called \textquotedblleft the root\textquotedblright\ of the tree), which
leads to a parent/child relationship on each edge (namely: the endpoint closer
to the root is called the \textquotedblleft parent\textquotedblright\ of the
endpoint further away from the root). Often, they also impose a total order on
the children of each given vertex. With these extra data, a tree can be used
for addressing objects, since each vertex has a unique \textquotedblleft path
description\textquotedblright\ from the root leading to it (e.g.,
\textquotedblleft the second child of the fourth child of the
root\textquotedblright). But this all is going too far afield for us here; we
are mainly interested in trees as graphs, and won't impose any extra structure
unless we need it for something.
\end{remark}

\begin{exercise}
\label{exe.5.2}Let $G$ be a multigraph that has no loops. Assume that there
exists a vertex $u$ of $G$ such that
\[
\text{for each vertex $v$ of $G$, there is a \textbf{unique} path from $u$ to
$v$ in $G$.}%
\]
Prove that $G$ is a tree. \medskip

[\textbf{Remark:} Pay attention to the quantifiers used here: $\exists
u\forall v$. This differs from the $\forall u\forall v$ in Statement T2 of the
tree equivalence theorem (Theorem \ref{thm.trees.T1-8}).]
\end{exercise}

\subsection{Leaves}

Continuing with our faux-botanical terminology, we define leaves in a tree:

\begin{definition}
Let $T$ be a tree. A vertex of $T$ is said to be a \textbf{leaf} if its degree
is $1$.
\end{definition}

For example, the tree
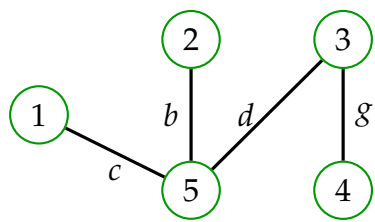
\[%
%
\]
has three leaves: $1$, $2$ and $4$.

How to find a tree with as many leaves as possible (for a given number of
vertices)? For any $n\geq3$, the simple graph%
\[
\left(  \left\{  0,1,\ldots,n-1\right\}  ,\ \left\{  0i\ \mid\ i>0\right\}
\right)
\]
is a tree (when considered as a multigraph), and has $n-1$ leaves (namely, all
of $1,2,\ldots,n-1$). This tree is called an $n$\textbf{-star graph}, as it
looks as follows:%
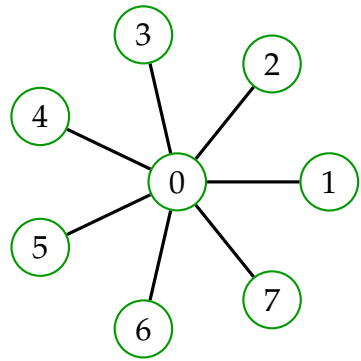
\[%
%
\ \ \ \ \ \ \ \ \ \ \left(  \text{for }n=8\right)  .
\]
It is easy to see that no tree with $n\geq3$ vertices can have more than $n-1$
leaves, so the $n$-star graph is optimal in this sense. Note that for $n=2$,
the $n$-star graph has $2$ leaves, not $1$.

How to find a tree with as few leaves as possible? For any $n\geq2$, the
$n$-path graph%
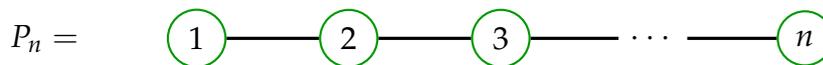
\[%
%
\]
is a tree with only $2$ leaves (viz., the vertices $1$ and $n$). Can we find a
tree with fewer leaves? For $n=1$, yes, because the $1$-path graph $P_{1}$
(this is simply the graph with $1$ vertex and no edges) has no leaves at all.
However, for $n\geq2$, the $n$-path graph is the best we can do:

\begin{theorem}
\label{thm.tree.leaves.2}Let $T$ be a tree with at least $2$ vertices. Then:

\begin{enumerate}
\item[\textbf{(a)}] The tree $T$ has at least $2$ leaves.

\item[\textbf{(b)}] Let $v$ be a vertex of $T$. Then, there exist two distinct
leaves $p$ and $q$ of $T$ such that $v$ lies on the path from $p$ to $q$.
\end{enumerate}
\end{theorem}

Note that I'm saying \textquotedblleft the path\textquotedblright\ rather than
\textquotedblleft a path\textquotedblright\ here. This is allowed, because in
a tree, for any two vertices $p$ and $q$, there is a \textbf{unique} path from
$p$ to $q$. This follows from Statement T2 in the tree equivalence theorem
(Theorem \ref{thm.trees.T1-8}).

\begin{proof}
[Proof of Theorem \ref{thm.tree.leaves.2}.]\textbf{(b)} We apply a variant of
the \textquotedblleft longest path trick\textquotedblright: Among all paths
that contain the vertex $v$, let $\mathbf{w}$ be a longest one. Let $p$ be the
starting point of $\mathbf{w}$, and let $q$ be the ending point of
$\mathbf{w}$. We shall show that $p$ and $q$ are two distinct leaves.

[Here is a picture of $\mathbf{w}$, for what it's worth:%
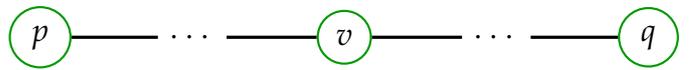
\[%
\begin{tikzpicture}[scale=2]
\begin{scope}[every node/.style={circle,thick,draw=green!60!black}]
\node(1) at (0,0) {$p$};
\node(3) at (2,0) {$v$};
\node(5) at (4,0) {$q$};
\end{scope}
\node(4) at (3,0) {$\cdots$};
\node(2) at (1,0) {$\cdots$};
\begin{scope}[every edge/.style={draw=black,very thick}, every loop/.style={}]
\path[-] (1) edge (2) (2) edge (3) (3) edge (4) (4) edge (5);
\end{scope}
\end{tikzpicture}%
\ \ .
\]
Of course, the tree $T$ can have other edges as well, not just those of
$\mathbf{w}$.]

First, we observe that $T$ is connected (since $T$ is a tree), and has at
least one vertex $u$ distinct from $v$ (since $T$ has at least $2$ vertices).
Hence, $T$ has a path $\mathbf{r}$ that connects $v$ to $u$. This path
$\mathbf{r}$ must contain at least one edge (since $u\neq v$). Thus, we have
found a path $\mathbf{r}$ of $T$ that contains $v$ and contains at least one
edge. Hence, the path $\mathbf{w}$ must contain at least one edge as well
(since $\mathbf{w}$ is a longest path that contains $v$, and thus cannot be
shorter than $\mathbf{r}$). Since $\mathbf{w}$ is a path from $p$ to $q$, we
thus conclude that $p\neq q$ (because if a path contains at least one edge,
then its starting point is distinct from its ending point).

Now, assume (for the sake of contradiction) that $p$ is not a leaf. Then,
$\deg p\neq1$. The path $\mathbf{w}$ already contains one edge that contains
$p$ (namely, the first edge of $\mathbf{w}$). Since $\deg p\neq1$, there must
be another edge $f$ of $T$ that contains $\mathbf{w}$. Consider this $f$. Let
$p^{\prime}$ be its endpoint distinct from $p$ (if $f$ is a loop, then we set
$p^{\prime}=p$). Appending this edge $f$ (and its endpoint) to the beginning
of the path $\mathbf{w}$, we obtain a backtrack-free walk
\[
\left(  p^{\prime},f,\underbrace{p,\ldots,v,\ldots,q}_{\text{This is
}\mathbf{w}}\right)
\]
(this is backtrack-free since $f$ is not the first edge of $\mathbf{w}$).
According to Proposition \ref{prop.btf-walk.cyc}, this backtrack-free walk
either is a path or contains a cycle. Since $T$ has no cycle (because $T$ is a
forest), we thus conclude that this backtrack-free walk is a path. It is
furthermore a path that contains $v$ and is longer than $\mathbf{w}$ (longer
by $1$, in fact). But this contradicts the fact that $\mathbf{w}$ is a longest
path that contains $v$. This contradiction shows that our assumption (that $p$
is not a leaf) was wrong.

Hence, $p$ is a leaf. A similar argument shows that $q$ is a leaf (here, we
need to append the new edge at the end of $\mathbf{w}$ rather than at the
beginning). Thus, $p$ and $q$ are two distinct leaves of $T$ (distinct because
$p\neq q$) such that $v$ lies on the path from $p$ to $q$ (since $v$ lies on
the path $\mathbf{w}$, which is a path from $p$ to $q$). This proves Theorem
\ref{thm.tree.leaves.2} \textbf{(b)}. \medskip

\textbf{(a)} Pick any vertex $v$ of $T$. Then, Theorem \ref{thm.tree.leaves.2}
\textbf{(b)} shows that there exist two distinct leaves $p$ and $q$ of $T$
such that $v$ lies on the path from $p$ to $q$. Thus, in particular, there
exist two distinct leaves $p$ and $q$ of $T$. In other words, $T$ has at least
two leaves. This proves Theorem \ref{thm.tree.leaves.2} \textbf{(a)}.

[\textit{Remark:} Another way to prove part \textbf{(a)} is to write the tree
$T$ as $T=\left(  V,E,\varphi\right)  $, and recall the handshake lemma, which
yields%
\begin{align*}
\sum_{v\in V}\deg v  &  =2\cdot\left\vert E\right\vert =2\cdot\left(
\left\vert V\right\vert -1\right)  \ \ \ \ \ \ \ \ \ \ \left(  \text{since
}\left\vert E\right\vert =\left\vert V\right\vert -1\text{ in a tree}\right)
\\
&  =2\cdot\left\vert V\right\vert -2.
\end{align*}
Since each $v\in V$ satisfies $\deg v\geq1$ (why?), this equality entails that
at least two vertices $v\in V$ must satisfy $\deg v\leq1$ (since otherwise,
the sum $\sum_{v\in V}\deg v$ would be $\geq2\cdot\left\vert V\right\vert
-1$), and therefore these two vertices are leaves.]
\end{proof}

\bigskip

Leaves are particularly helpful for performing induction on trees. The formal
reason for this is the following theorem:

\begin{theorem}
[induction principle for trees]\label{thm.tree.leaf-ind.1}Let $T$ be a tree
with at least $2$ vertices. Let $v$ be a leaf of $T$. Let $T\setminus v$ be
the multigraph obtained from $T$ by removing $v$ and all edges that contain
$v$ (note that there is only one such edge, since $v$ is a leaf). Then,
$T\setminus v$ is again a tree.
\end{theorem}

Here is an example of a tree $T$ and of the smaller tree $T\setminus v$
obtained by removing a leaf $v$ (namely, $v=3$):%
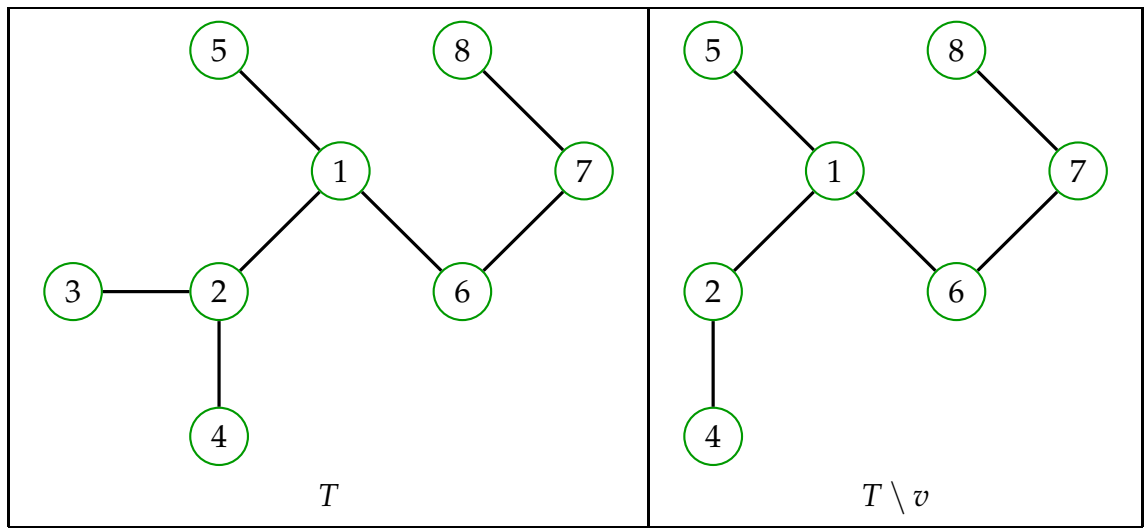
\[%

\]

\begin{proof}
[Proof of Theorem \ref{thm.tree.leaf-ind.1}.]Write $T$ as $T=\left(
V,E,\varphi\right)  $. Thus, $T\setminus v$ is the induced subgraph $T\left[
V\setminus\left\{  v\right\}  \right]  $.

The graph $T$ is a tree, thus a forest; hence, it has no cycles. Thus, the
graph $T\setminus v$ has no cycles either. Hence, it is a forest.

Furthermore, this forest $T\setminus v$ has at least $1$ vertex (since $T$ has
at least $2$ vertices).

We shall now show that any two vertices $p$ and $q$ of $T\setminus v$ are
path-connected in $T\setminus v$.

Indeed, let $p$ and $q$ be two vertices of $T\setminus v$. Then, $p$ and $q$
are path-connected in $T$ (since $T$ is connected). Hence, there exists a path
$\mathbf{w}$ from $p$ to $q$ in $T$. Consider this path $\mathbf{w}$. Note
that $v$ is neither the starting point nor the ending point of this path
$\mathbf{w}$ (since $p$ and $q$ are vertices of $T\setminus v$, and thus
distinct from $v$). Hence, if $v$ was a vertex of $\mathbf{w}$, then
$\mathbf{w}$ would contain \textbf{two distinct} edges that contain $v$
(namely, the edge just before $v$ and the edge just after $v$). But this is
impossible, since there is only one edge available that contains $v$ (because
$v$ is a leaf). Thus, $v$ cannot be a vertex of $\mathbf{w}$. Hence, the path
$\mathbf{w}$ does not use the vertex $v$, and thus is a path in the graph
$T\setminus v$ as well. So the vertices $p$ and $q$ are path-connected in
$T\setminus v$.

We have now shown that any two vertices $p$ and $q$ of $T\setminus v$ are
path-connected in $T\setminus v$. This shows that $T\setminus v$ is connected
(since $T\setminus v$ has at least $1$ vertex). Hence, $T\setminus v$ is a
tree (since $T\setminus v$ is a forest).
\end{proof}

Theorem \ref{thm.tree.leaf-ind.1} has a converse as well:

\begin{theorem}
\label{thm.tree.leaf-ind.2}Let $G$ be a multigraph. Let $v$ be a vertex of $G$
such that $\deg v=1$ and such that $G\setminus v$ is a tree. (Here,
$G\setminus v$ means the multigraph obtained from $G$ by removing the vertex
$v$ and all edges that contain $v$.) Then, $G$ is a tree.
\end{theorem}

\begin{proof}
Left to the reader. (The main step is to show that a cycle of $G$ cannot
contain $v$.)
\end{proof}

Theorem \ref{thm.tree.leaf-ind.1} helps prove many properties of trees by
induction on the number of vertices. In the induction step, remove a leaf $v$
and apply the induction hypothesis to $T\setminus v$.

The following exercise is essentially a generalization of Theorem
\ref{thm.tree.leaves.2} \textbf{(a)}:

\begin{exercise}
\label{exe.5.7}Let $T$ be a tree. Let $w$ be any vertex of $T$. Prove that $T$
has at least $\deg w$ many leaves.
\end{exercise}

\begin{exercise}
Let $T$ be a tree that has at least two vertices. A \textbf{branch vertex} of
$T$ means a vertex $v$ of $T$ such that $\deg v\geq3$.

\begin{enumerate}
\item[\textbf{(a)}] Prove that if $T$ has $k$ leaves (for some $k\in
\mathbb{N}$), then $T$ has at most $k-2$ branch vertices.

\item[\textbf{(b)}] Prove that $T$ has two distinct leaves $u$ and $v$ such
that the path from $u$ to $v$ contains at most one branch vertex.
\end{enumerate}
\end{exercise}

\begin{exercise}
\label{exe.5.8c}A dominating set of a multigraph $G$ is defined to be a
dominating set of its underlying simple graph $G^{\operatorname{simp}}$.

Let $G$ be a forest. Prove that
\[
\sum_{D\text{ is a dominating set of }G}\left(  -1\right)  ^{\left\vert
D\right\vert }=\pm1.
\]

\end{exercise}

\begin{exercise}
\label{exe.hw3.whamilton} Let $T$ be a tree having more than $1$ vertex. Let
$L$ be the set of leaves of $T$. Prove that it is possible to add $\left\vert
L\right\vert -1$ new edges to $T$ in such a way that the resulting multigraph
has a Hamiltonian cycle. \medskip

[\textbf{Solution:} This is Exercise 4 on homework set \#3 from my Spring 2017
course; see \href{https://www.cip.ifi.lmu.de/~grinberg/t/17s/}{the course
page} for solutions.]
\end{exercise}

\subsection{Spanning trees}

\subsubsection{Spanning subgraphs}

We now proceed to a crucial application of trees. First we define a concept
that makes sense for any multigraphs:

\begin{definition}
A \textbf{spanning subgraph} of a multigraph $G=\left(  V,E,\varphi\right)  $
means a multigraph of the form $\left(  V,F,\varphi\mid_{F}\right)  $, where
$F$ is a subset of $E$.

In other words, it means a submultigraph of $G$ with the same vertex set as
$G$.

In other words, it means a multigraph obtained from $G$ by removing some
edges, but leaving all vertices undisturbed.
\end{definition}

Compare this to the notion of an induced subgraph:

\begin{itemize}
\item To build an \textbf{induced} subgraph, we throw away some vertices but
keep all the edges that we can keep. (As usual in mathematics, the words
\textquotedblleft some vertices\textquotedblright\ include \textquotedblleft
no vertices\textquotedblright\ and \textquotedblleft all
vertices\textquotedblright.)

\item In contrast, to build a \textbf{spanning} subgraph, we keep all vertices
but throw away some edges.
\end{itemize}

\subsubsection{Spanning trees}

Spanning subgraphs are particularly useful when they are trees:

\begin{definition}
\label{def.spt}A \textbf{spanning tree} of a multigraph $G$ means a spanning
subgraph of $G$ that is a tree.
\end{definition}

\begin{example}
Let $G$ be the following multigraph:%
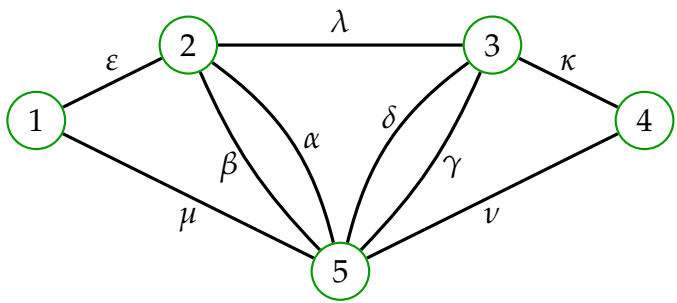
\[%
%
\ \ .
\]
Here is a spanning tree of $G$:%
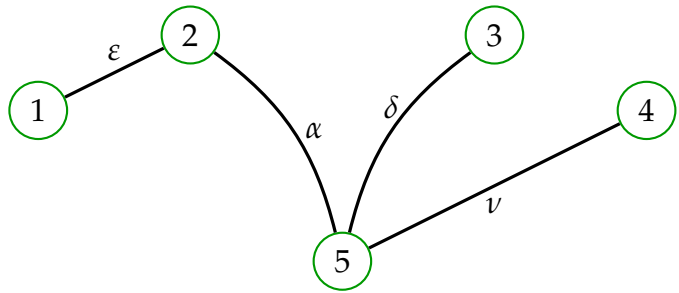
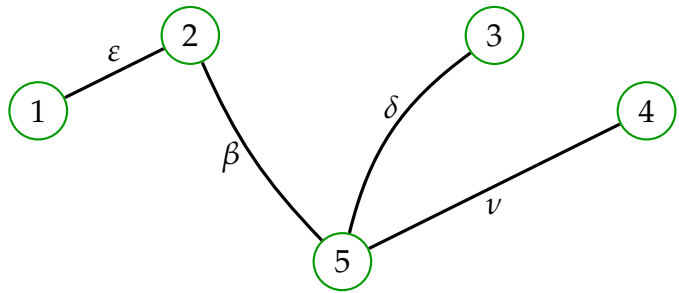
\[%
%
\ \ .
\]
Here is another:
\[%
%
\ \ .
\]
(Yes, this is a different one, because $\alpha\neq\beta$.) And here is yet
another spanning tree of $G$:%
\[%
%
\ \ .
\]

\end{example}

\begin{example}
\label{exa.spt.Cn}Let $n$ be a positive integer. Consider the cycle graph
$C_{n}$. (We defined this graph $C_{n}$ in Definition \ref{def.sg.cycle} for
all $n\geq2$, but we later redefined $C_{2}$ and defined $C_{1}$ in Definition
\ref{def.mg.Cn}. Here, we are using the latter modified definition.)

The graph $C_{n}$ has exactly $n$ spanning trees. Indeed, any graph obtained
from $C_{n}$ by removing a single edge is a spanning tree of $C_{n}$.
\end{example}

\begin{proof}
A tree with $n$ vertices must have exactly $n-1$ edges (by the implication
T1$\Longrightarrow$T4 in Theorem \ref{thm.trees.T1-8}). Thus, a spanning
subgraph of $C_{n}$ can be a tree only if it has $n-1$ edges, i.e., only if it
is obtained from $C_{n}$ by removing a single edge (since $C_{n}$ has $n$
edges in total). Thus, $C_{n}$ has at most $n$ spanning trees (since $C_{n}$
has $n$ edges that can be removed). It remains to check that any subgraph
obtained from $C_{n}$ by removing a single edge is indeed a spanning tree. But
this is easy, since all such subgraphs are isomorphic to the path graph
$P_{n}$. This proves Example \ref{exa.spt.Cn}.
\end{proof}

\begin{exercise}
\label{exe.hw3.countST-ac}Fix $m\geq1$. Let $G$ be the simple graph with
$3m+2$ vertices
\[
a,b,x_{1},x_{2},\ldots,x_{m},y_{1},y_{2},\ldots,y_{m},z_{1},z_{2},\ldots,z_{m}%
\]
and the following $3m+3$ edges:
\begin{align*}
&  ax_{1},\ \ ay_{1},\ \ az_{1},\\
&  x_{i}x_{i+1},\ \ y_{i}y_{i+1},\ \ z_{i}z_{i+1}\qquad\text{ for all }%
i\in\left\{  1,2,\ldots,m-1\right\}  ,\\
&  x_{m}b,\ \ y_{m}b,\ \ z_{m}b.
\end{align*}
(Thus, the graph consists of two vertices $a$ and $b$ connected by three
paths, each of length $m+1$, with no overlaps between the paths except for
their starting and ending points. Here is a picture for $m=3$:
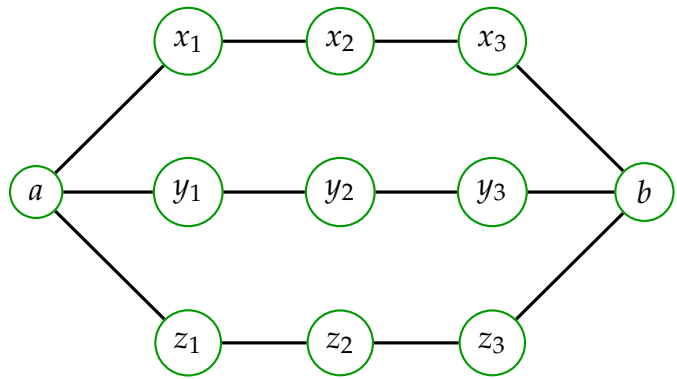
\[%
\begin{tikzpicture}
\begin{scope}[every node/.style={circle,thick,draw=green!60!black}]
\node(A) at (0,0) {$a$};
\node(B) at (2,0) {$y_1$};
\node(C) at (4,0) {$y_2$};
\node(D) at (6,0) {$y_3$};
\node(B2) at (2,2) {$x_1$};
\node(C2) at (4,2) {$x_2$};
\node(D2) at (6,2) {$x_3$};
\node(B3) at (2,-2) {$z_1$};
\node(C3) at (4,-2) {$z_2$};
\node(D3) at (6,-2) {$z_3$};
\node(E) at (8,0) {$b$};
\end{scope}
\begin{scope}[every edge/.style={draw=black,very thick}]
\path[-] (A) edge (B) (B) edge (C) (C) edge (D) (D) edge (E);
\path[-] (A) edge (B2) (B2) edge (C2) (C2) edge (D2) (D2) edge (E);
\path[-] (A) edge (B3) (B3) edge (C3) (C3) edge (D3) (D3) edge (E);
\end{scope}
\end{tikzpicture}%
\]
) Compute the number of spanning trees of $G$.

[To argue why your number is correct, a sketch of the argument in 1-2
sentences should be enough; a fully rigorous proof is not required.]\medskip

[\textbf{Solution:} This is Exercise 2 \textbf{(c)} on homework set \#3 from
my Spring 2017 course; see
\href{https://www.cip.ifi.lmu.de/~grinberg/t/17s/}{the course page} for solutions.]
\end{exercise}

\subsubsection{Spanning forests}

A spanning tree of a graph $G$ can be regarded as a minimum \textquotedblleft
backbone\textquotedblright\ of $G$ -- that is, a way to keep $G$ connected
using as few edges as possible. Of course, if $G$ is not connected, then this
is not possible at all, so $G$ has no spanning trees in this case. The best
one can hope for is a spanning subgraph that keeps each component of $G$
connected using as few edges as possible. This is known as a \textquotedblleft
spanning forest\textquotedblright:

\begin{definition}
\label{def.spf}A \textbf{spanning forest} of a multigraph $G$ means a spanning
subgraph $H$ of $G$ that is a forest and satisfies $\operatorname*{conn}%
H=\operatorname*{conn}G$.
\end{definition}

When $G$ is a connected multigraph, a spanning forest of $G$ means the same as
a spanning tree of $G$.

\begin{exercise}
Here is an example of a graph and one of its spanning forests:%
\[%
\begin{tabular}
[c]{|c|c|}\hline
graph $G$ & spanning forest of $G$\\\hline
\ \ $%
\begin{tikzpicture}[scale=1.8]
\begin{scope}[every node/.style={circle,thick,draw=green!60!black}]
\node(1) at (0,0) {$1$};
\node(2) at (0,1) {$2$};
\node(3) at (1,1) {$3$};
\node(4) at (1,0) {$4$};
\node(5) at (2,0) {$5$};
\node(6) at (2,1) {$6$};
\node(7) at (3,1) {$7$};
\node(8) at (3,0) {$8$};
\end{scope}
\draw[very thick] (5) -- (4) -- (3) -- (2) -- (1) -- (4);
\draw[very thick] (6) -- (7) -- (8) -- (6);
\end{tikzpicture}%
\ \ $ & \ \ $%
\begin{tikzpicture}[scale=1.8]
\begin{scope}[every node/.style={circle,thick,draw=green!60!black}]
\node(1) at (0,0) {$1$};
\node(2) at (0,1) {$2$};
\node(3) at (1,1) {$3$};
\node(4) at (1,0) {$4$};
\node(5) at (2,0) {$5$};
\node(6) at (2,1) {$6$};
\node(7) at (3,1) {$7$};
\node(8) at (3,0) {$8$};
\end{scope}
\draw[very thick] (5) -- (4) -- (3) -- (2);
\draw[very thick] (1) -- (4);
\draw[very thick] (8) -- (6) -- (7);
\end{tikzpicture}%
\ \ $\\\hline
\end{tabular}
\]

\end{exercise}

\subsubsection{Existence and construction of a spanning tree}

The following theorem is crucial, which is why we will outline four different proofs:

\begin{theorem}
\label{thm.spt.exists}Each connected multigraph $G$ has at least one spanning tree.
\end{theorem}

\begin{proof}
[First proof.]Let $G$ be a connected multigraph. We want to construct a
spanning tree of $G$. We try to achieve this by removing edges from $G$ one by
one, until $G$ becomes a tree. When doing so, we must be careful not to
disconnect the graph (i.e., not to destroy its connectedness). According to
Theorem \ref{thm.mg.G-e.conn}, this can be achieved by making sure that we
never remove a bridge (i.e., an edge that appears in no cycle). Thus, we keep
removing non-bridges (i.e., edges that are not bridges) as long as we can
(i.e., until we end up with a graph in which every edge is a bridge).

So here is the algorithm: We start with $G$, and we successively remove
non-bridges one by one until we no longer have any non-bridges
left\footnote{\textbf{Warning:} We cannot remove several non-bridges at once!
We have to remove them one by one. Indeed, if $e$ and $f$ are two non-bridges
of $G$, then there is no guarantee that $f$ remains a non-bridge in
$G\setminus e$. So we cannot remove both $e$ and $f$ simultaneously; we have
to remove one of them and check whether the other is still a non-bridge.}.
This procedure cannot go on forever, since $G$ has only finitely many edges.
Thus, after finitely many steps, we will end up with a graph that has no
non-bridges any more. This resulting graph therefore has no cycles (since any
cycle would have at least one edge, and this edge would be a non-bridge), but
is still connected (since $G$ was connected, and we never lost connectedness
as we removed only non-bridges). Thus, this resulting graph is a tree. Since
it is also a spanning subgraph of $G$ (by construction), it is therefore a
spanning tree of $G$. This proves Theorem \ref{thm.spt.exists}.
\end{proof}

\begin{proof}
[Second proof (sketched).]In the above first proof, we constructed a spanning
tree of $G$ by starting with $G$ and successively removing edges until we got
a tree. Now let us take the opposite strategy: Start with an empty graph on
the same vertex set as $G$, and successively add edges (from $G$) until we get
a connected graph.

Here are some details: We start with a graph $L$ that has the same vertex set
as $G$, but has no edges. Now, we inspect all edges $e$ of $G$ one by one (in
some order). For each such edge $e$, we add it to $L$, but only if it does not
create a cycle in $L$; otherwise, we discard this edge. Notice that adding an
edge $e$ with endpoints $u$ and $v$ to $L$ creates a cycle if and only if $u$
and $v$ lie in the same component of $L$ (before we add $e$). Thus, we only
add an edge to $L$ if its endpoints lie in different components of $L$;
otherwise, we discard it. This way, at the end of the procedure, our graph $L$
will still have no cycles (since we never create any cycles). In other words,
it will be a forest.

Let me denote this forest by $H$. (Thus, $H$ is the $L$ at the end of the
procedure.) I claim that this forest $H$ is a spanning tree of $G$. Why? Since
we know that $H$ is a forest, we only need to show that $H$ is connected.
Assume the contrary. Thus, there is at least one edge $e$ of $G$ whose
endpoints lie in different components of $H$ (why?). This edge $e$ is
therefore not an edge of $H$. Therefore, at some point during our construction
of $H$, we must have discarded this edge $e$ (instead of adding it to $L$). As
we know, this means that the endpoints of $e$ used to lie in the same
component of $L$ at the point at which we discarded $e$. But this entails that
these two endpoints lie in the same component of $L$ at the end of the
procedure as well (because the graph $L$ never loses any edges during the
procedure, so that any two vertices that used to lie in the same component of
$L$ at some point will still lie in the same component of $L$ ever after). In
other words, the endpoints of $e$ lie in the same component of $H$. This
contradicts our assumption that the endpoints of $e$ lie in different
components of $H$. This contradiction completes our proof that $H$ is
connected. Hence, $H$ is a spanning tree of $G$, and we have proved Theorem
\ref{thm.spt.exists} again.
\end{proof}

\begin{proof}
[Third proof.]This proof takes yet another approach to constructing a spanning
tree of $G$: We choose an arbitrary vertex $r$ of $G$, and then progressively
\textquotedblleft spread a rumor\textquotedblright\ from $r$. The rumor starts
at vertex $r$. On day $0$, only $r$ has heard the rumor. Every day, every
vertex that knows the rumor spreads it to all its neighbors (i.e., all
vertices adjacent to it). Since $G$ is connected, the rumor will eventually
spread to every vertex of $G$. Now, each vertex $v$ (other than $r$) remembers
which other vertex $v^{\prime}$ it has first heard the rumor from (if it heard
it from several vertices at the same time, it just picks one of them), and
picks some edge $e_{v}$ that has endpoints $v$ and $v^{\prime}$ (such an edge
must exist, since $v$ must have heard the rumor from a neighbor). The edges
$e_{v}$ for all $v\in V\setminus\left\{  r\right\}  $ (where $V$ is the vertex
set of $G$) then form a spanning tree of $G$ (that is, the graph with vertex
set $V$ and edge set $\left\{  e_{v}\ \mid\ v\in V\setminus\left\{  r\right\}
\right\}  $ is a spanning tree). Why?

Intuitively, this is quite convincing: This graph cannot have cycles (because
that would require a time loop) and must be connected (because for any vertex
$v$, we can trace back the path of the rumor from $r$ to $v$ by following the
edges $e_{v}$ backwards). To obtain a rigorous proof, we formalize this
construction mathematically:

Write $G$ as $G=\left(  V,E,\varphi\right)  $. Choose any vertex $r$ of $G$.

We shall recursively construct a sequence of subgraphs%
\[
\left(  V_{0},E_{0},\varphi_{0}\right)  ,\ \ \ \ \left(  V_{1},E_{1}%
,\varphi_{1}\right)  ,\ \ \ \ \left(  V_{2},E_{2},\varphi_{2}\right)
,\ \ \ \ \ldots
\]
of $G$. The idea behind these subgraphs is that for each $i\in\mathbb{N}$, the
set $V_{i}$ will consist of all vertices $v$ that have heard the rumor by day
$i$, and the set $E_{i}$ will consist of the corresponding edges $e_{v}$. The
map $\varphi_{i}$ will be the restriction of $\varphi$ to $E_{i}$, of course.

Here is the exact construction of this sequence of subgraphs:

\begin{itemize}
\item \textit{Recursion base:} Set $V_{0}:=\left\{  r\right\}  $ and
$E_{0}:=\varnothing$. Let $\varphi_{0}$ be the restriction of $\varphi$ to the
(empty) set $E_{0}$.

\item \textit{Recursion step:} Let $i\in\mathbb{N}$. Assume that the subgraph
$\left(  V_{i},E_{i},\varphi_{i}\right)  $ of $G$ has already been defined.
Now, we set%
\[
V_{i+1}:=V_{i}\cup\left\{  v\in V\ \mid\ v\text{ is adjacent to some vertex in
}V_{i}\right\}  .
\]
For each $v\in V_{i+1}\setminus V_{i}$, we choose \textbf{one} edge $e_{v}$
that joins\footnote{We say that an edge \textbf{joins} a vertex $p$ to a
vertex $q$ if the endpoints of this edge are $p$ and $q$.} $v$ to a vertex in
$V_{i}$ (such an edge exists, since $v\in V_{i+1}$; if there are several, we
just choose a random one). Set%
\[
E_{i+1}:=E_{i}\cup\left\{  e_{v}\ \mid\ v\in V_{i+1}\setminus V_{i}\right\}
.
\]
Finally, we let $\varphi_{i+1}$ be the restriction of the map $\varphi$ to the
set $E_{i+1}$. This is a map from $E_{i+1}$ to $\mathcal{P}_{1,2}\left(
V_{i+1}\right)  $ (because any edge $e_{v}$ with $v\in V_{i+1}\setminus V_{i}$
has one endpoint $v$ in $V_{i+1}\setminus V_{i}\subseteq V_{i+1}$ and the
other endpoint in $V_{i}\subseteq V_{i+1}$). Thus, $\left(  V_{i+1}%
,E_{i+1},\varphi_{i+1}\right)  $ is a well-defined subgraph of $G$.
\end{itemize}

This construction yields that $\left(  V_{i},E_{i},\varphi_{i}\right)  $ is a
subgraph of $\left(  V_{i+1},E_{i+1},\varphi_{i+1}\right)  $ for each
$i\in\mathbb{N}$. Hence, $V_{0}\subseteq V_{1}\subseteq V_{2}\subseteq\cdots$,
so that $\left\vert V_{0}\right\vert \leq\left\vert V_{1}\right\vert
\leq\left\vert V_{2}\right\vert \leq\cdots$. Since a sequence of integers
bounded from above cannot keep increasing forever (and the sizes $\left\vert
V_{i}\right\vert $ are bounded from above by $\left\vert V\right\vert $, since
each $V_{i}$ is a subset of $V$), we thus see that there exists some
$i\in\mathbb{N}$ such that $\left\vert V_{i}\right\vert =\left\vert
V_{i+1}\right\vert $. Consider this $i$. From $\left\vert V_{i}\right\vert
=\left\vert V_{i+1}\right\vert $, we obtain $V_{i}=V_{i+1}$ (since
$V_{i}\subseteq V_{i+1}$).

In our colloquial model above, $V_{i}=V_{i+1}$ means that no new vertices
learn the rumor on day $i+1$; it is reasonable to expect that at this point,
every vertex has heard the rumor. In other words, we claim that $V_{i}=V$. A
rigorous proof of this can be easily given using the fact that $G$ is
connected\footnote{Here is the \textit{proof} in detail: We must show that
$V_{i}=V$. Assume the contrary. Thus, there exists a vertex $u\in V\setminus
V_{i}$. Consider this $u$. The path from $r$ to $u$ starts at a vertex in
$V_{i}$ (since $r\in V_{0}\subseteq V_{i}$) and ends at a vertex in
$V\setminus V_{i}$ (since $u\in V\setminus V_{i}$). Thus, it must cross over
from $V_{i}$ into $V\setminus V_{i}$ at some point. Therefore, there exists an
edge with one endpoint in $V_{i}$ and the other endpoint in $V\setminus V_{i}%
$. Let $v$ and $w$ be these two endpoints, so that $v\in V_{i}$ and $w\in
V\setminus V_{i}$. Then, $w$ is adjacent to some vertex in $V_{i}$ (namely, to
$v$), and therefore belongs to $V_{i+1}$ (by the definition of $V_{i+1}$).
Hence, $w\in V_{i+1}=V_{i}$. But this contradicts $w\notin V\setminus V_{i}$.
This contradiction shows that our assumption was wrong, qed.}.

Now, we claim that the subgraph $\left(  V_{i},E_{i},\varphi_{i}\right)  $ is
a spanning tree of $G$. To see this, we must show that this subgraph is a
forest and is connected (since $V_{i}=V$ already shows that it is a spanning
subgraph). Before we do this, let us give an example:

\begin{example}
Let $G$ be the following multigraph:%
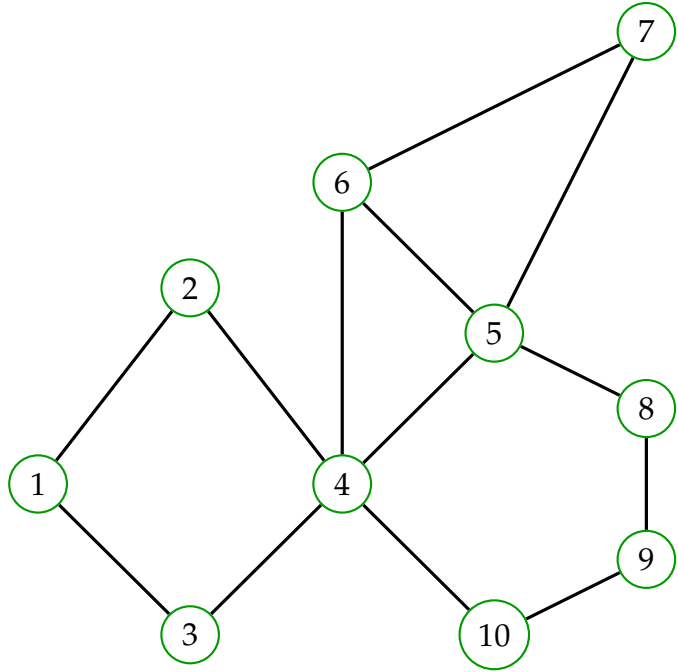
\[%
%
\ \ .
\]
Set $r=3$. Then, the above construction yields%
\begin{align*}
V_{0}  &  =\left\{  3\right\}  ,\\
V_{1}  &  =\left\{  3,1,4\right\}  ,\\
V_{2}  &  =\left\{  3,1,4,2,5,6,10\right\}  ,\\
V_{3}  &  =\left\{  3,1,4,2,5,6,10,8,9,7\right\}  =V,
\end{align*}
so that $V_{k}=V$ for all $k\geq3$. Thus, we can take $i=3$. Here is an image
of the $V_{k}$ as progressively growing circles:%
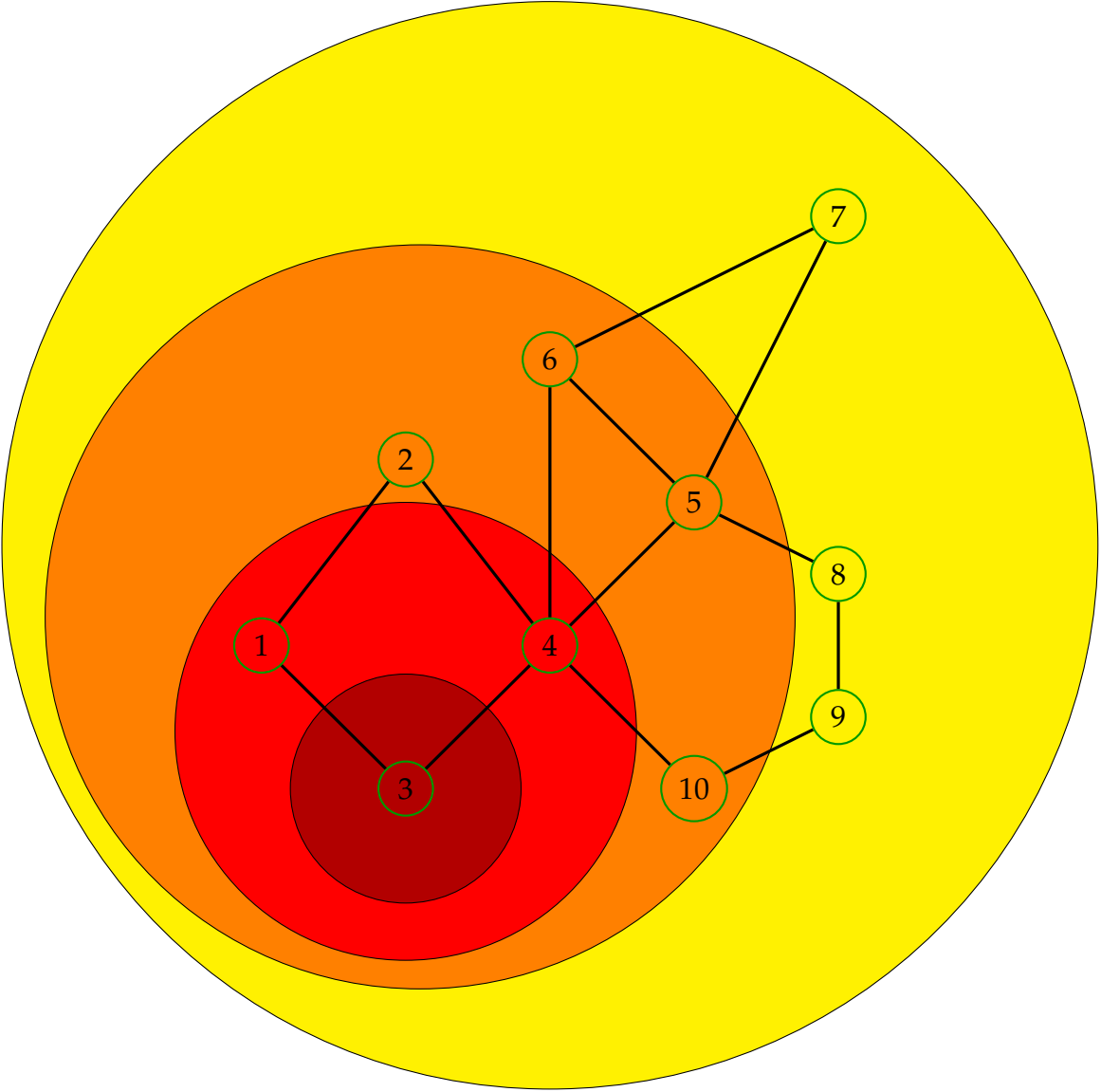
\[%
%
\ \ .
\]
(The dark-red inner circle is $V_{0}$; the red circle is $V_{1}$; the orange
circle is $V_{2}$; the yellow circle is $V_{3}=V_{4}=V_{5}=\cdots=V$.)
Finally, the edges $e_{v}$ can be chosen to be the following (we are painting
them red for clarity):%
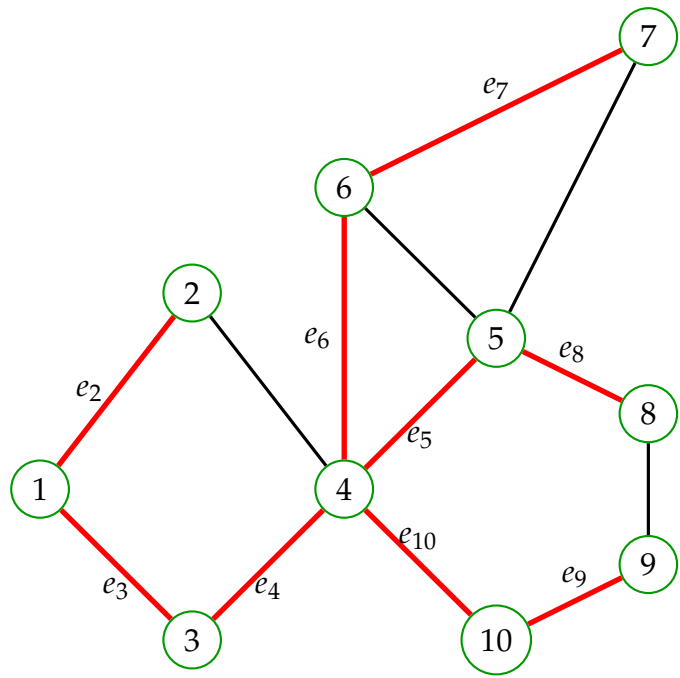
\[%
%
\ \ .
\]
(Here, we have made two choices: We chose $e_{2}$ to be the edge joining $2$
with $1$ rather than the edge joining $2$ with $4$, and we chose $e_{7}$ to be
the edge joining $7$ with $6$ rather than $7$ with $5$. The other options
would have been equally fine.)
\end{example}

We now return to the general proof. Let us first show the following:

\begin{statement}
\textit{Claim 1:} Let $j\in\mathbb{N}$. Each vertex of the graph $\left(
V_{j},E_{j},\varphi_{j}\right)  $ is path-connected to $r$ in this graph.
\end{statement}

[\textit{Proof of Claim 1:} We induct on $j$:

\textit{Base case:} For $j=0$, Claim 1 is obvious, since $V_{0}=\left\{
r\right\}  $ (so the only vertex of the graph in question is $r$ itself).

\textit{Induction step:} Fix some positive integer $k$. Assume (as the
induction hypothesis) that Claim 1 holds for $j=k-1$. That is, each vertex of
the graph $\left(  V_{k-1},E_{k-1},\varphi_{k-1}\right)  $ is path-connected
to $r$ in this graph.

Now, let $v$ be a vertex of the graph $\left(  V_{k},E_{k},\varphi_{k}\right)
$. We must show that $v$ is path-connected to $r$ in this graph. If $v\in
V_{k-1}$, then this follows from the induction hypothesis (since $\left(
V_{k-1},E_{k-1},\varphi_{k-1}\right)  $ is a subgraph of $\left(  V_{k}%
,E_{k},\varphi_{k}\right)  $). Thus, we WLOG assume that $v\notin V_{k-1}$
from now on. Hence, $v\in V_{k}\setminus V_{k-1}$. According to the recursive
definition of $E_{k}$, this entails that there is an edge $e_{v}\in E_{k}$
that joins $v$ to some vertex $u\in V_{k-1}$. Consider this latter vertex $u$.
Then, $v$ is path-connected to $u$ in the graph $\left(  V_{k},E_{k}%
,\varphi_{k}\right)  $ (since the edge $e_{v}$ provides a length-$1$ path from
$v$ to $u$). However, $u$ is path-connected to $r$ in the graph $\left(
V_{k-1},E_{k-1},\varphi_{k-1}\right)  $ (by the induction hypothesis, since
$u\in V_{k-1}$), hence also in the graph $\left(  V_{k},E_{k},\varphi
_{k}\right)  $ (since $\left(  V_{k-1},E_{k-1},\varphi_{k-1}\right)  $ is a
subgraph of $\left(  V_{k},E_{k},\varphi_{k}\right)  $). Since the relation
\textquotedblleft path-connected\textquotedblright\ is transitive, we conclude
from the previous two sentences that $v$ is path-connected to $r$ in the graph
$\left(  V_{k},E_{k},\varphi_{k}\right)  $.

So we have shown that each vertex $v$ of the graph $\left(  V_{k}%
,E_{k},\varphi_{k}\right)  $ is path-connected to $r$ in the graph $\left(
V_{k},E_{k},\varphi_{k}\right)  $. In other words, Claim 1 holds for $j=k$.
This completes the induction step, and Claim 1 is proved.] \medskip

Claim 1 (applied to $j=i$) shows that each vertex of the graph $\left(
V_{i},E_{i},\varphi_{i}\right)  $ is path-connected to $r$ in this graph.
Since the relation \textquotedblleft path-connected\textquotedblright\ is an
equivalence relation, this entails that any two vertices of this graph are
path-connected. Thus, the graph $\left(  V_{i},E_{i},\varphi_{i}\right)  $ is
connected (since it has at least one vertex). It remains to prove that this
graph $\left(  V_{i},E_{i},\varphi_{i}\right)  $ is a forest.

Again, we do this using an auxiliary claim:

\begin{statement}
\textit{Claim 2:} Let $j\in\mathbb{N}$. Then, the graph $\left(  V_{j}%
,E_{j},\varphi_{j}\right)  $ has no cycles.
\end{statement}

[\textit{Proof of Claim 2:} We induct on $j$:

\textit{Base case:} The graph $\left(  V_{0},E_{0},\varphi_{0}\right)  $ has
no edges (because $E_{0}=\varnothing$) and thus no cycles. Thus, Claim 2 holds
for $j=0$.

\textit{Induction step:} Fix some positive integer $k$. Assume (as the
induction hypothesis) that Claim 2 holds for $j=k-1$. That is, the graph
$\left(  V_{k-1},E_{k-1},\varphi_{k-1}\right)  $ has no cycles.

Now, let $\mathbf{c}$ be a cycle of the graph $\left(  V_{k},E_{k},\varphi
_{k}\right)  $. Then, $\mathbf{c}$ must use at least one edge from
$E_{k}\setminus E_{k-1}$ (since otherwise, $\mathbf{c}$ would be a cycle of
the graph $\left(  V_{k-1},E_{k-1},\varphi_{k-1}\right)  $, but this is
impossible, since $\left(  V_{k-1},E_{k-1},\varphi_{k-1}\right)  $ has no
cycles). However, each edge from $E_{k}\setminus E_{k-1}$ has the form $e_{v}$
for some $v\in V_{k}\setminus V_{k-1}$ (because of how $E_{k}$ was defined).
Thus, $\mathbf{c}$ must have an edge of this form. Consider the corresponding
vertex $v\in V_{k}\setminus V_{k-1}$. The cycle $\mathbf{c}$ contains the edge
$e_{v}$ and therefore also contains its endpoint $v$. However, (again by the
definition of $E_{k}$) the edge $e_{v}$ is the \textbf{only} edge in $E_{k}$
that contains the vertex $v$. Since the edge $e_{v}$ is not a loop (because it
joins the vertex $v\in V_{k}\setminus V_{k-1}$ with a vertex in $V_{k-1}$), we
thus conclude that the vertex $v$ has degree $1$ in the graph $\left(
V_{k},E_{k},\varphi_{k}\right)  $. Thus, the vertex $v$ cannot be contained in
any cycle of $\left(  V_{k},E_{k},\varphi_{k}\right)  $ (because a cycle
cannot contain a vertex of degree $1$). This contradicts the fact that the
cycle $\mathbf{c}$ contains $v$.

Forget that we fixed $\mathbf{c}$. We thus have obtained a contradiction for
each cycle $\mathbf{c}$ of the graph $\left(  V_{k},E_{k},\varphi_{k}\right)
$. Hence, the graph $\left(  V_{k},E_{k},\varphi_{k}\right)  $ has no cycles.
In other words, Claim 2 holds for $j=k$. This completes the induction step,
and Claim 2 is proved.] \medskip

Applying Claim 2 to $j=i$, we see that the graph $\left(  V_{i},E_{i}%
,\varphi_{i}\right)  $ has no cycles. In other words, this graph is a forest.
Since it is connected, it is therefore a tree. Since it is a spanning subgraph
of $G$, we thus conclude that it is a spanning tree of $G$. Hence, we have
constructed a spanning tree of $G$. \medskip

We note an important property of this construction:

\begin{statement}
\textit{Claim 3:} For each $k\in\mathbb{N}$, we have%
\[
V_{k}=\left\{  v\in V\ \mid\ d\left(  r,v\right)  \leq k\right\}  ,
\]
where $d\left(  r,v\right)  $ means the length of a shortest path from $r$ to
$v$.
\end{statement}

This is easily proved by induction on $k$. Thus, the spanning tree $\left(
V_{i},E_{i},\varphi_{i}\right)  $ we have constructed has the following
property: For each $v\in V$, the path from $r$ to $v$ in this spanning tree is
a shortest path from $r$ to $v$ in $G$. For this reason, this spanning tree is
called a \textbf{breadth-first search (\textquotedblleft BFS\textquotedblright%
) tree}. Note that the choice of root $r$ is important here: It is usually not
true that the path from an arbitrary vertex $u$ to an arbitrary vertex $v$
along our spanning tree is a shortest path in $G$. No spanning tree of $G$ has
this property, unless $G$ itself is \textquotedblleft more or less a
tree\textquotedblright\ (more precisely, unless $G^{\operatorname*{simp}}$ is
a tree)!
\end{proof}

\begin{proof}
[Fourth proof of Theorem \ref{thm.spt.exists} (sketched).]We imagine a snake
that slithers along the edges of $G$, trying to eventually bite each vertex.
It starts at some vertex $r$, which it immediately bites. Any time the snake
enters a vertex $v$, it makes the following step:

\begin{itemize}
\item If some neighbor of $v$ has not been bitten yet, then the snake picks
such a neighbor $w$ as well as some edge $f$ that joins $w$ with $v$; the
snake then moves to $w$ along the edge $f$, bites the vertex $w$ and marks the
edge $f$.

\item If not, then the snake marks the vertex $v$ as fully digested and
backtracks (along the marked edges) to the last vertex it has visited but not
fully digested yet. (To \textquotedblleft\textbf{backtrack}\textquotedblright%
\ along a marked edge means to use it in the direction opposite to how it was
originally used. In other words, if a marked edge $f$ is initially used to go
from $v$ to $w$, then backtracking along it means walking it from $w$ to $v$.)
\end{itemize}

Once backtracking is no longer possible (because there are no more vertices
left that are not fully digested), the procedure is finished. I claim that the
marked edges at that moment are the edges of a spanning tree of $G$.

I won't prove this claim in detail, but I will give some hints. First,
however, an example:

\begin{example}
\label{exa.spt.exists.4th-proof.1}Let $G$ be the following connected
multigraph:%
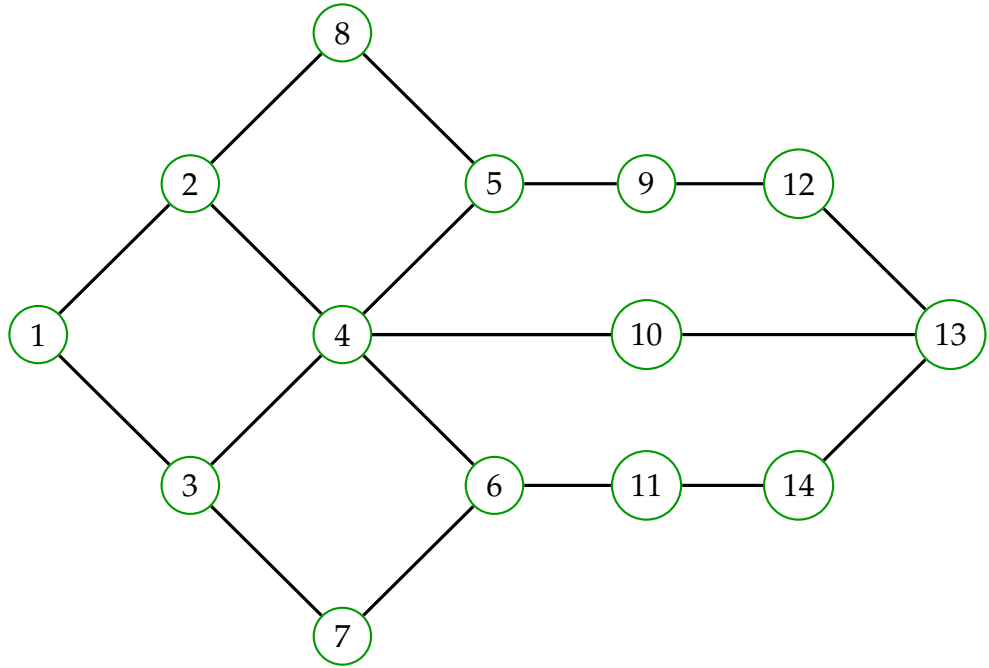
\[%
\begin{tikzpicture}[scale=2]
\begin{scope}[every node/.style={circle,thick,draw=green!60!black}]
\node(1) at (-1,0) {$1$};
\node(2) at (0,1) {$2$};
\node(3) at (0,-1) {$3$};
\node(4) at (1,0) {$4$};
\node(7) at (1,-2) {$7$};
\node(8) at (1, 2) {$8$};
\node(5) at (2, 1) {$5$};
\node(6) at (2, -1) {$6$};
\node(9) at (3, 1) {$9$};
\node(10) at (3, 0) {$10$};
\node(11) at (3, -1) {$11$};
\node(12) at (4, 1) {$12$};
\node(13) at (5, 0) {$13$};
\node(14) at (4, -1) {$14$};
\end{scope}
\begin{scope}[every edge/.style={draw=black,very thick}, every loop/.style={}]
\path[-] (1) edge (2) edge (3);
\path[-] (2) edge (4) edge (8);
\path[-] (3) edge (4) edge (7);
\path[-] (8) edge (5);
\path[-] (4) edge (5) edge (6) edge (10);
\path[-] (7) edge (6);
\path[-] (5) edge (9);
\path[-] (4) edge (10);
\path[-] (6) edge (11);
\path[-] (9) edge (12);
\path[-] (10) edge (13);
\path[-] (11) edge (14);
\path[-] (13) edge (12) edge (14);
\end{scope}
\end{tikzpicture}%
\ \ .
\]
Let our snake start its journey at $r=3$. It bites this vertex. Then, let's
say that it picks the vertex $1$ as its next victim (it could just as well go
to $4$ or $7$; the snake has many choices, but we follow one possible trip).
Thus, it next arrives at vertex $1$, bites it and marks the edge that brought
it to this vertex. As its next destination, it necessarily picks the vertex
$2$ (since vertex $3$ has already been bitten). It moves to vertex $2$, bites
it and marks the edge. Next, let's say that it picks the vertex $4$ (the other
option would be $8$). It thus moves to $4$, bites it and marks the edge.
Proceeding likewise, it then moves to $5$ (the other options are $6$ and $10$;
the vertices $2$ and $3$ do not qualify since they are already bitten), bites
$5$ and marks an edge. From there, let's say it moves to $8$, bites $8$ and
marks an edge. Now, there is no longer an unbitten neighbor of $8$ to move to.
Thus, the snake marks the vertex $8$ as fully digested and backtracks to the
last vertex not fully digested -- which, at this point, is $5$. From this
vertex $5$, it moves on to $9$ (this is the only option, since $4$ and $8$
have already been bitten). And so on. Here is one possible outcome of this
journey (there are a few more decisions that the snake can make here, so you
may get a different one):%
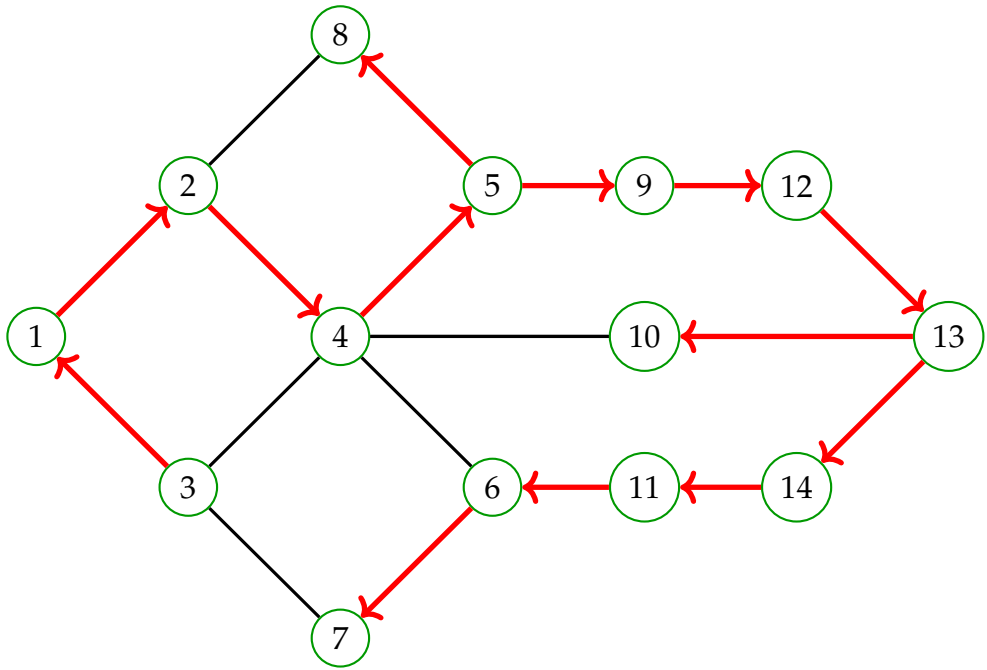
\[%
%
\ \ .
\]
Here, the marked edges are drawn in bold red ink, and endowed with an arrow
that represents the direction in which they were first used (e.g., the edge
joining $2$ with $4$ has an arrow towards $4$ because it was first used to get
from $2$ to $4$).
\end{example}

Now, as promised, let me outline a proof of the above claim (that the marked
edges form a spanning tree of $G$). To wit, argue the following four
observations (ideally in this order):

\begin{enumerate}
\item After each step, the marked edges are precisely the edges along which
the snake has moved so far.

\item After each step, the network of bitten vertices and marked edges is a tree.

\item After enough steps, each bitten vertex is fully digested.

\item At that point, the network of bitten vertices and marked edges is a
spanning tree (since each neighbor of a fully digested vertex is bitten, thus
fully digested by observation 3).
\end{enumerate}

\noindent Details are left to the reader.

The result is that Theorem \ref{thm.spt.exists} is proved once again. However,
more comes out of the above construction if you know where to look. The
spanning tree $T$ of $G$ whose edges are the edges marked by the snake is
called a \textbf{depth-first search (\textquotedblleft DFS\textquotedblright)
tree}. It has the following extra property: If $u$ and $v$ are two adjacent
vertices of $G$, then either $u$ lies on the path from $r$ to $v$ in $T$, or
$v$ lies on the path from $r$ to $u$ in $T$. (This called a \textquotedblleft
lineal spanning tree\textquotedblright. See \cite[\S 6.1]{BenWil06} for details.)
\end{proof}

\subsubsection{Applications}

Spanning trees have lots of applications:

\begin{itemize}
\item A spanning tree of a graph can be viewed as a kind of \textquotedblleft
backbone\textquotedblright\ of the graph, which in particular provides
\textquotedblleft canonical\textquotedblright\ paths between any two vertices.
This is useful, e.g., for networking applications where having a choice
between different paths would be problematic (see, e.g.,
\href{https://en.wikipedia.org/wiki/Spanning_Tree_Protocol}{the Spanning Tree
Protocol}).

\item A $w$-minimum spanning tree (see Exercise \ref{exe.5.6} = Homework set
\#5 exercise 6) solves a global version of the cheapest-path problem. It can
also be used for detecting clusters.

\item Depth-first search (the algorithm used in our fourth proof of Theorem
\ref{thm.spt.exists}) can also be used as a way to traverse all vertices of a
given graph and return back to the starting point. In particular, this
provides an algorithmic way to solve mazes (since a maze can be modeled as a
graph, where the vertices correspond to \textquotedblleft
rooms\textquotedblright\ and the edges correspond to \textquotedblleft
doors\textquotedblright). This appears to have been the original motivation
for Tr\'{e}maux to invent depth-first search back in the 19th century.
\end{itemize}

Here is a more theoretical application of spanning trees:

\begin{definition}
A vertex $v$ of a connected multigraph $G$ is said to be a \textbf{cut-vertex}
if the graph $G\setminus v$ is disconnected. (Recall that $G\setminus v$ is
the multigraph obtained from $G$ by removing the vertex $v$ and all edges that
contain $v$.)
\end{definition}

\begin{proposition}
Let $G$ be a connected multigraph with $\geq2$ vertices. Then, there are at
least $2$ vertices of $G$ that are \textbf{not} cut-vertices.
\end{proposition}

\begin{proof}
Pick a spanning tree $T$ of $G$ (we know from Theorem \ref{thm.spt.exists}
that such a spanning tree exists). Then, $T$ has at least $2$ leaves (by
Theorem \ref{thm.tree.leaves.2} \textbf{(a)}). But each leaf of $T$ is a
non-cut-vertex of $G$ (why?).
\end{proof}

\begin{remark}
It is not true that conversely, any non-leaf of $T$ is a cut-vertex of $G$. So
we cannot get any lower bound on the number of cut-vertices. And this is not
surprising: Lots of graphs (e.g., the complete graph $K_{n}$ for $n\geq2$)
have no cut-vertices at all. These graphs are said to be \textbf{2-connected},
and their properties have been amply studied (see, e.g., \cite[\S 4.2]{West01}
for an introduction).
\end{remark}

\subsubsection{Exercises}

\begin{exercise}
\label{exe.5.0}Let $G$ be a connected multigraph. Let $T_{1}$ and $T_{2}$ be
two spanning trees of $G$.

Prove the following:\footnotemark

\begin{enumerate}
\item[\textbf{(a)}] For any $e\in\operatorname*{E}\left(  T_{1}\right)
\setminus\operatorname*{E}\left(  T_{2}\right)  $, there exists an
$f\in\operatorname*{E}\left(  T_{2}\right)  \setminus\operatorname*{E}\left(
T_{1}\right)  $ with the property that replacing $e$ by $f$ in $T_{1}$ (that
is, removing the edge $e$ from $T_{1}$ and adding the edge $f$) results in a
spanning tree of $G$.

\item[\textbf{(b)}] For any $f\in\operatorname*{E}\left(  T_{2}\right)
\setminus\operatorname*{E}\left(  T_{1}\right)  $, there exists an
$e\in\operatorname*{E}\left(  T_{1}\right)  \setminus\operatorname*{E}\left(
T_{2}\right)  $ with the property that replacing $e$ by $f$ in $T_{1}$ (that
is, removing the edge $e$ from $T_{1}$ and adding the edge $f$) results in a
spanning tree of $G$.
\end{enumerate}

[\textbf{Hint:} The two parts look very similar, but (to my knowledge) their
proofs are not.]
\end{exercise}

\footnotetext{Recall that $\operatorname*{E}\left(  H\right)  $ denotes the
edge set of any graph $H$.}

\begin{exercise}
Let $G$ be a connected multigraph. Let $\mathcal{S}$ be the simple graph whose
vertices are the spanning trees of $G$, and whose edges are defined as
follows: Two spanning trees $T_{1}$ and $T_{2}$ of $G$ are adjacent (as
vertices of $\mathcal{S}$) if and only if $T_{2}$ can be obtained from $T_{1}$
by removing an edge and adding another (i.e., if and only if there exist an
edge $e_{1}$ of $T_{1}$ and an edge $e_{2}$ of $T_{2}$ such that $e_{2} \neq
e_{1}$ and $T_{2} \setminus e_{2} = T_{1} \setminus e_{1}$).

Prove that the simple graph $\mathcal{S}$ is itself connected. (In simpler
language: Prove that any spanning tree of $G$ can be transformed into any
other spanning tree of $G$ by a sequence of legal \textquotedblleft remove an
edge and add another\textquotedblright\ operations, where such an operation is
called \textbf{legal} if its result is a spanning tree of $G$.) \medskip

[\textbf{Example:} If $G$ is the multigraph
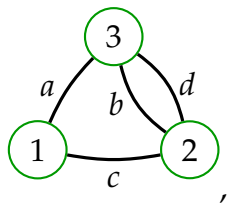
\[

,
\]
then the graph $\mathcal{S}$ looks as follows:
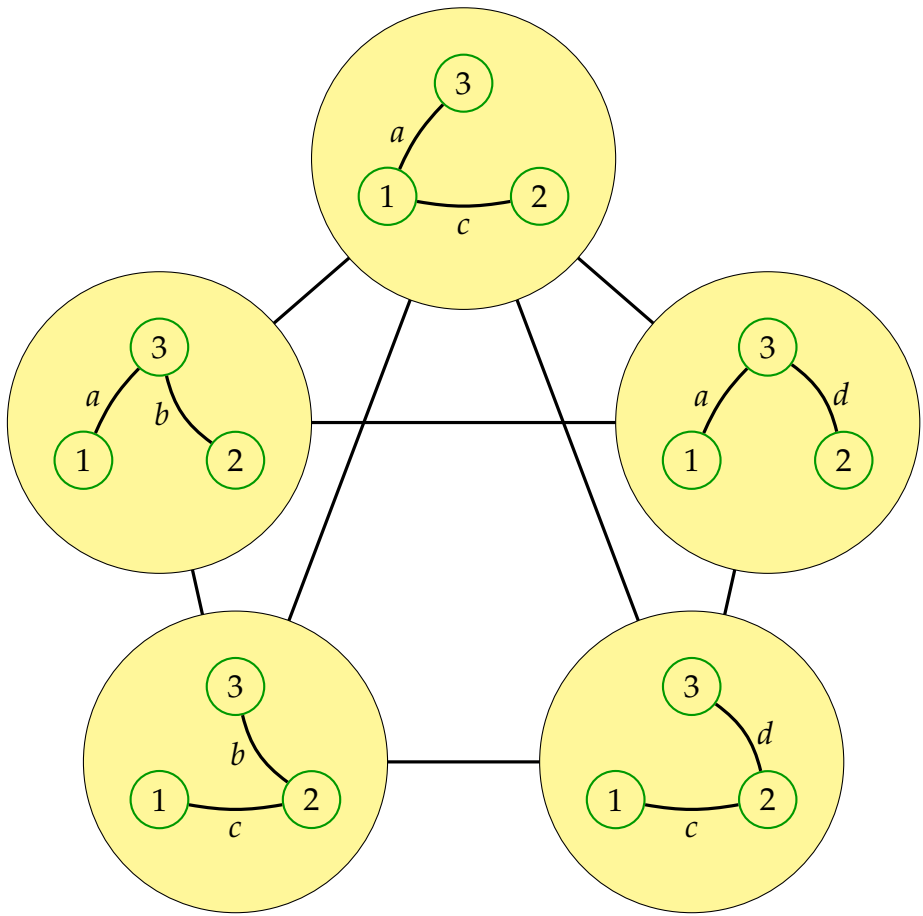
\[
%
\]
]
\end{exercise}

\begin{exercise}
\label{exe.5.6}Let $G=\left(  V,E,\varphi\right)  $ be a connected multigraph.
Let $w:E\rightarrow\mathbb{R}$ be a map that assigns a real number $w\left(
e\right)  $ to each edge $e$. We shall call this real number $w\left(
e\right)  $ the \textbf{weight} of the edge $e$.

If $H = \left(  W, F, \varphi\mid_{F} \right)  $ is a subgraph of $G$, then
the \textbf{weight} $w\left(  H \right)  $ of $H$ is defined to be $\sum_{f
\in F} w\left(  f \right)  $ (that is, the sum of the weights of all edges of
$H$).

A \textbf{$w$-minimum spanning tree} of $G$ means a spanning tree of $G$ that
has the smallest weight among all spanning trees of $G$.

In our first proof of Theorem \ref{thm.spt.exists}, we have seen a way to
construct a spanning tree of $G$ by successively removing non-bridges until
only bridges remain. (A \textbf{non-bridge} means an edge that is not a bridge.)

Now, let us perform this algorithm, but taking care to choose a non-bridge of
largest weight (among all non-bridges) at each step. Prove that the result
will be a $w$-minimum spanning tree.
\end{exercise}

\begin{exercise}
\label{exe.6.1}Let $G$ be a connected multigraph with an even number of
vertices. Prove that there exists a spanning subgraph $H$ of $G$ such that
each vertex of $H$ has odd degree (in $H$). \medskip

\textbf{[Hint:} One way to solve this begins by reducing the problem to the
case when $G$ is a tree.]
\end{exercise}

\subsubsection{Existence and construction of a spanning forest}

So we have learnt that connected graphs have spanning trees. What do
disconnected graphs have?

\begin{corollary}
Each multigraph has a spanning forest.
\end{corollary}

\begin{proof}
Apply Theorem \ref{thm.spt.exists} to each component of the multigraph. Then,
combine the resulting spanning trees into a spanning forest.
\end{proof}

\subsection{Centers of graphs and trees}

\subsubsection{Distances}

Given a graph, we can define a \textquotedblleft distance\textquotedblright%
\ between any two of its vertices, simply by counting edges on the shortest
path from one to the other:

\begin{definition}
\label{def.mg.dist}Let $G$ be a multigraph.

For any two vertices $u$ and $v$ of $G$, we define the \textbf{distance}
between $u$ and $v$ to be the smallest length of a path from $u$ to $v$. If no
such path exists, then this distance is defined to be $\infty$.

The distance between $u$ and $v$ is denoted by $d\left(  u,v\right)  $ or by
$d_{G}\left(  u,v\right)  $ when the graph $G$ is not clear from the context.
\end{definition}

\begin{example}
If $G$ is the multigraph from Example \ref{exa.spt.exists.4th-proof.1}, then%
\[
d_{G}\left(  1,9\right)  =4,\ \ \ \ \ \ \ \ \ \ d_{G}\left(  4,13\right)
=2,\ \ \ \ \ \ \ \ \ \ d_{G}\left(  4,4\right)  =0.
\]

\end{example}

\begin{remark}
Distances in a multigraph satisfy the rules that you would expect a distance
function to satisfy:

\begin{enumerate}
\item[\textbf{(a)}] We have $d\left(  u,u\right)  =0$ for any vertex $u$.

\item[\textbf{(b)}] We have $d\left(  u,v\right)  =d\left(  v,u\right)  $ for
any vertices $u$ and $v$.

\item[\textbf{(c)}] We have $d\left(  u,v\right)  +d\left(  v,w\right)  \geq
d\left(  u,w\right)  $ for any vertices $u$, $v$ and $w$. (Here, we understand
that $\infty\geq m$ and $\infty+m=\infty$ for any $m\in\mathbb{N}$.)
\end{enumerate}

Also:

\begin{enumerate}
\item[\textbf{(d)}] The distances $d\left(  u,v\right)  $ do not change if we
replace \textquotedblleft path\textquotedblright\ by \textquotedblleft
walk\textquotedblright\ in the definition of the distance.

\item[\textbf{(e)}] If $V$ is the vertex set of our multigraph, then $d\left(
u,v\right)  \leq\left\vert V\right\vert -1$ for any vertices $u$ and $v$.
\end{enumerate}
\end{remark}

\begin{proof}
Part \textbf{(d)} follows from Corollary \ref{cor.mg.walk-thus-path}. The
proofs of \textbf{(a)}, \textbf{(b)} and \textbf{(c)} are then straightforward
(the proof of \textbf{(c)} relies on part \textbf{(d)}, because splicing two
paths generally only yields a walk, not a path). Finally, in order to prove
part \textbf{(e)}, observe that any path of our multigraph has length
$\leq\left\vert V\right\vert -1$ (since its vertices are distinct).
\end{proof}

We note that the definition of a distance becomes simpler if our multigraph is
a tree: Namely, if $T$ is a tree, then the distance $d\left(  u,v\right)  $
between two vertices $u$ and $v$ is the length of the \textbf{only} path from
$u$ to $v$ in $T$. Thus, in a tree, we do not have to worry whether a given
path is the shortest.

We also notice that if $G$ is a multigraph, and if $u$ and $v$ are two
vertices of $G$, then the distance $d_{G}\left(  u,v\right)  $ in $G$ equals
the distance $d_{G^{\operatorname*{simp}}}\left(  u,v\right)  $ in the simple
graph $G^{\operatorname*{simp}}$. (The reason for this is that any path of $G$
can be converted into a path of $G^{\operatorname*{simp}}$ having the same
length, and vice versa. Of course, this is not a one-to-one correspondence,
but it suffices for our purposes.) Thus, when studying distances on a
multigraph, we can WLOG restrict ourselves to simple graphs. \medskip

The following few exercises give some curious properties of distances in
various kinds of graphs.

\begin{exercise}
\label{exe.mt1.d+d+d} Let $a$, $b$ and $c$ be three vertices of a connected
multigraph $G=\left(  V,E,\varphi\right)  $. Prove that $d\left(  b,c\right)
+d\left(  c,a\right)  +d\left(  a,b\right)  \leq2\left\vert V\right\vert -2$.
\medskip

[\textbf{Solution:} This is Exercise 7 on midterm \#1 from my Spring 2017
course, except that the simple graph has been replaced by a multigraph (but
this makes no serious difference); see
\href{https://www.cip.ifi.lmu.de/~grinberg/t/17s/}{the course page} for solutions.]
\end{exercise}

\begin{exercise}
\label{exe.hw3.d+d+d.directed} Let $a$, $b$ and $c$ be three vertices of a
strongly connected multidigraph $D=\left(  V,A,\psi\right)  $ such that
$\left\vert V\right\vert \geq4$. For any two vertices $u$ and $v$ of $D$, we
define the distance $d\left(  u,v\right)  $ to be the smallest length of a
path from $u$ to $v$. (This definition is the obvious analogue of Definition
\ref{def.mg.dist} for digraphs.)

\begin{enumerate}
\item[\textbf{(a)}] Prove that $d\left(  b,c\right)  +d\left(  c,a\right)
+d\left(  a,b\right)  \leq3\left\vert V\right\vert -4$.

\item[\textbf{(b)}] For each $n\geq5$, construct an example in which
$\left\vert V\right\vert =n$ and $d\left(  b,c\right)  +d\left(  c,a\right)
+d\left(  a,b\right)  =3\left\vert V\right\vert -4$. (No proof is required for
the example.)
\end{enumerate}

[\textbf{Solution:} This is Exercise 5 on homework set \#3 from my Spring 2017
course, except that the simple digraph has been replaced by a multidigraph
(but this makes no serious difference); see
\href{https://www.cip.ifi.lmu.de/~grinberg/t/17s/}{the course page} for solutions.]
\end{exercise}

\begin{exercise}
\label{exe.mt2.tropigrass} Let $G$ be a tree. Let $x$, $y$, $z$ and $w$ be
four vertices of $G$.

Show that the two largest ones among the three numbers
\[
d\left(  x,y\right)  +d\left(  z,w\right)  ,\qquad d\left(  x,z\right)
+d\left(  y,w\right)  \qquad\text{and}\qquad d\left(  x,w\right)  +d\left(
y,z\right)
\]
are equal. \medskip

[\textbf{Solution:} This is Exercise 6 on midterm \#2 from my Spring 2017
course; see \href{https://www.cip.ifi.lmu.de/~grinberg/t/17s/}{the course
page} for solutions.]
\end{exercise}

\begin{exercise}
\label{exe.mt3.tropigrass2} Let $G$ be a connected multigraph. Let $x$, $y$,
$z$ and $w$ be four vertices of $G$.

Assume that the two largest ones among the three numbers
\[
d\left(  x,y\right)  +d\left(  z,w\right)  ,\qquad d\left(  x,z\right)
+d\left(  y,w\right)  \qquad\text{and}\qquad d\left(  x,w\right)  +d\left(
y,z\right)
\]
are \textbf{not} equal.

Prove that $G$ has a cycle of length $\leq d\left(  x,z\right)  +d\left(
y,w\right)  +d\left(  x,w\right)  +d\left(  y,z\right)  $. \medskip

[\textbf{Hint:} This is a strengthening of Exercise \ref{exe.mt2.tropigrass}.
Try deriving it by applying the latter exercise to a strategically chosen
subgraph of $G$.] \medskip

[\textbf{Solution:} This is Exercise 1 on midterm \#3 from my Spring 2017
course; see \href{https://www.cip.ifi.lmu.de/~grinberg/t/17s/}{the course
page} for solutions.]
\end{exercise}

\subsubsection{Eccentricity and centers}

We can now define \textquotedblleft eccentricities\textquotedblright:

\begin{definition}
Let $v$ be a vertex of a multigraph $G=\left(  V,E,\varphi\right)  $. The
\textbf{eccentricity} of $v$ (with respect to $G$) is defined to be the number%
\[
\max\left\{  d\left(  v,u\right)  \ \mid\ u\in V\right\}  \in\mathbb{N}%
\cup\left\{  \infty\right\}  .
\]
This eccentricity is denoted by $\operatorname*{ecc}v$ or $\operatorname*{ecc}%
\nolimits_{G}v$.
\end{definition}

Intuitively, the eccentricity of a vertex $v$ is \textquotedblleft how far you
can get away from $v$\textquotedblright\ on $G$.

\begin{definition}
Let $G=\left(  V,E,\varphi\right)  $ be a multigraph. Then, a \textbf{center}
of $G$ means a vertex of $G$ whose eccentricity is minimum (among all vertices).
\end{definition}

(Some authors have a slightly different definition of a \textquotedblleft
center\textquotedblright: They define the \textbf{center} of $G$ to be the
\textbf{set} of all vertices of $G$ whose eccentricity is minimum. That is,
what they call \textquotedblleft center\textquotedblright\ is the set of what
we call \textquotedblleft centers\textquotedblright.)

\begin{example}
\label{exa.eccent.eccent1}Let $G$ be the following multigraph:%
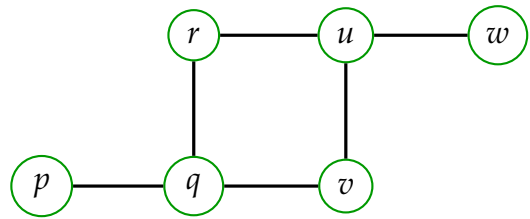
\[%
%
\ \ .
\]
Then, the eccentricities of its vertices are as follows (we are just labeling
each vertex with its eccentricity):%
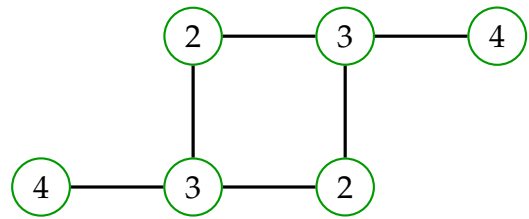
\[%
%
\ \ .
\]
Thus, the centers of $G$ are the vertices $r$ and $v$.
\end{example}

\begin{example}
Let $G$ be a complete graph $K_{n}$ (with $n$ vertices). Then, each vertex of
$G$ has the same eccentricity (which is $1$ if $n\geq2$ and $0$ if $n=1$), and
thus each vertex of $G$ is a center of $G$.
\end{example}

\begin{example}
\label{exa.eccent.eccent3}Let $G$ be a graph with more than one component.
Then, each vertex $v$ of $G$ has eccentricity $\infty$ (because there exists
at least one vertex $u$ that lies in a different component of $G$ than $v$,
and thus this vertex $u$ satisfies $d\left(  v,u\right)  =\infty$). Hence,
each vertex of $G$ is a center of $G$.
\end{example}

\subsubsection{The centers of a tree}

As we see from Example \ref{exa.eccent.eccent3}, eccentricity and centers are
not very useful notions when the graph is disconnected. Even for a connected
graph, Example \ref{exa.eccent.eccent1} shows that the centers do not
necessarily form a connected subgraph. However, in a tree, they behave a lot better:

\begin{theorem}
\label{thm.eccent.tree-main}Let $T$ be a tree. Then:

\begin{enumerate}
\item[\textbf{(a)}] The tree $T$ has either $1$ or $2$ centers.

\item[\textbf{(b)}] If $T$ has $2$ centers, then these $2$ centers are adjacent.

\item[\textbf{(c)}] Moreover, these centers can be found by the following algorithm:

If $T$ has more than $2$ vertices, then we remove all leaves from $T$
(simultaneously). What remains is again a tree. If that tree still has more
than $2$ vertices, we remove all leaves from it (simultaneously). The result
is again a tree. If that tree still has more than $2$ vertices, we remove all
leaves from it (simultaneously), and continue doing so until we are left with
a tree that has only $1$ or $2$ vertices. These vertices are the centers of
$T$.
\end{enumerate}
\end{theorem}

To prove Theorem \ref{thm.eccent.tree-main}, we first study how a tree is
affected when all its leaves are removed:

\begin{lemma}
\label{lem.eccent.tree-TL}Let $T=\left(  V,E,\varphi\right)  $ be a tree with
more than $2$ vertices.

Let $L$ be the set of all leaves of $T$.

Let $T\setminus L$ be the induced submultigraph of $T$ on the set $V\setminus
L$. (Thus, $T\setminus L$ is obtained from $T$ by removing all the vertices in
$L$ and all edges that contain a vertex in $L$.)

Then:

\begin{enumerate}
\item[\textbf{(a)}] The multigraph $T\setminus L$ is a tree.

\item[\textbf{(b)}] For any $u\in V\setminus L$ and $v\in V\setminus L$, we
have
\[
\left\{  \text{paths of }T\text{ from }u\text{ to }v\right\}  =\left\{
\text{paths of }T\setminus L\text{ from }u\text{ to }v\right\}
\]
(that is, the paths of $T$ from $u$ to $v$ are precisely the paths of
$T\setminus L$ from $u$ to $v$).

\item[\textbf{(c)}] For any $u\in V\setminus L$ and $v\in V\setminus L$, we
have $d_{T}\left(  u,v\right)  =d_{T\setminus L}\left(  u,v\right)  $.

\item[\textbf{(d)}] Each vertex $v\in V\setminus L$ satisfies
$\operatorname*{ecc}\nolimits_{T}v=\operatorname*{ecc}\nolimits_{T\setminus
L}v+1$.

\item[\textbf{(e)}] Each leaf $v\in L$ satisfies $\operatorname*{ecc}%
\nolimits_{T}v=\operatorname*{ecc}\nolimits_{T}w+1$, where $w$ is the unique
neighbor of $v$ in $T$. (A \textbf{neighbor} of $v$ means a vertex that is
adjacent to $v$.)

\item[\textbf{(f)}] The centers of $T$ are precisely the centers of
$T\setminus L$.
\end{enumerate}
\end{lemma}

\begin{example}
Let $T$ be the following tree:%
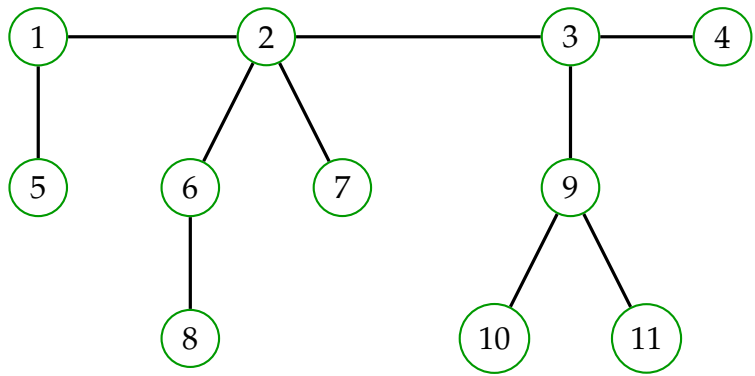
\[%
%
\ \ .
\]
Then, the set $L$ from Lemma \ref{lem.eccent.tree-TL} is $\left\{
4,5,7,8,10,11\right\}  $, and the tree $T\setminus L$ looks as follows:%
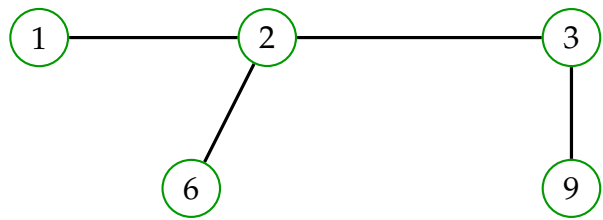
\[%
%
\ \ .
\]

\end{example}

\begin{proof}
[Proof of Lemma \ref{lem.eccent.tree-TL}.]First, we notice that $T$ is a
forest (since $T$ is a tree), and thus has no cycles. In particular, $T$
therefore has no loops and no parallel edges. Also, for any two vertices $u$
and $v$ of $T$, there is a unique path from $u$ to $v$ in $T$.

Next, we introduce some terminology: If $\mathbf{p}$ is a path of some
multigraph, then an \textbf{intermediate vertex} of $\mathbf{p}$ shall mean a
vertex of $\mathbf{p}$ that is neither the starting point nor the ending point
of $\mathbf{p}$. In other words, if $\mathbf{p}=\left(  p_{0},e_{1}%
,p_{1},e_{2},p_{2},\ldots,e_{k},p_{k}\right)  $ is a path of some multigraph,
then the intermediate vertices of $\mathbf{p}$ are $p_{1},p_{2},\ldots
,p_{k-1}$. Clearly, any intermediate vertex of a path $\mathbf{p}$ must have
degree $\geq2$ (since the path $\mathbf{p}$ enters it along some edge, and
leaves it along another). Hence, if $\mathbf{p}$ is a path of $T$, then
\begin{equation}
\text{any intermediate vertex of }\mathbf{p}\text{ must belong to }V\setminus
L \label{pf.lem.eccent.tree-TL.0}%
\end{equation}
(because it must have degree $\geq2$, thus cannot be a leaf of $T$; but this
means that it cannot belong to $L$; therefore, it must belong to $V\setminus
L$). \medskip

\textbf{(b)} Let $u\in V\setminus L$ and $v\in V\setminus L$. Let $\mathbf{p}$
be a path of $T$ from $u$ to $v$. We shall show that $\mathbf{p}$ is a path of
$T\setminus L$ as well.

Indeed, let us first check that all vertices of $\mathbf{p}$ belong to
$V\setminus L$. This is clear for the vertices $u$ and $v$ (since $u\in
V\setminus L$ and $v\in V\setminus L$); but it also holds for every
intermediate vertex of $\mathbf{p}$ (by (\ref{pf.lem.eccent.tree-TL.0})).
Thus, it does indeed hold for all vertices of $\mathbf{p}$.

We have thus shown that all vertices of $\mathbf{p}$ belong to $V\setminus L$.
Hence, $\mathbf{p}$ is a path of $T\setminus L$ (since $T\setminus L$ is the
induced submultigraph of $T$ on the set $V\setminus L$).

Forget that we fixed $\mathbf{p}$. We have thus shown that every path
$\mathbf{p}$ of $T$ from $u$ to $v$ is also a path of $T\setminus L$. Hence,%
\[
\left\{  \text{paths of }T\text{ from }u\text{ to }v\right\}  \subseteq
\left\{  \text{paths of }T\setminus L\text{ from }u\text{ to }v\right\}  .
\]
Conversely, we have%
\[
\left\{  \text{paths of }T\setminus L\text{ from }u\text{ to }v\right\}
\subseteq\left\{  \text{paths of }T\text{ from }u\text{ to }v\right\}  ,
\]
since every path of $T\setminus L$ is a path of $T$ (because $T\setminus L$ is
a submultigraph of $T$). Combining these two facts, we obtain%
\[
\left\{  \text{paths of }T\text{ from }u\text{ to }v\right\}  =\left\{
\text{paths of }T\setminus L\text{ from }u\text{ to }v\right\}  .
\]
This proves Lemma \ref{lem.eccent.tree-TL} \textbf{(b)}. \medskip

\textbf{(c)} This follows from Lemma \ref{lem.eccent.tree-TL} \textbf{(b)},
since the distance $d_{G}\left(  u,v\right)  $ of two vertices $u$ and $v$ in
a graph $G$ is defined to be the smallest length of a path from $u$ to $v$.
\medskip

\textbf{(a)} The graph $T$ is a tree, thus a forest. Hence, its submultigraph
$T\setminus L$ is a forest as well (since any cycle of $T\setminus L$ would be
a cycle of $T$). It thus remains to show that $T\setminus L$ is connected.

First, it is easy to see that $T\setminus L$ has at least one
vertex\footnote{\textit{Proof.} We assumed that $T$ has more than $2$
vertices. In other words, there exist three distinct vertices $u,v,w$ of $T$.
Consider these $u,v,w$. If all three distances $d_{T}\left(  u,v\right)  $,
$d_{T}\left(  v,w\right)  $ and $d_{T}\left(  w,u\right)  $ were equal to $1$,
then $T$ would have a cycle (of the form $\left(  u,\ast,v,\ast,w,\ast
,u\right)  $, where each asterisk stands for some edge); but this would
contradict the fact that $T$ has no cycles. Thus, not all of these three
distances are equal to $1$. Hence, at least one of them is $\neq1$. WLOG
assume that $d_{T}\left(  u,v\right)  \neq1$ (otherwise, we permute $u,v,w$).
Hence, the path from $u$ to $v$ has more than one edge (indeed, it must have
at least one edge, since $u$ and $v$ are distinct). Therefore, this path has
at least one intermediate vertex. This intermediate vertex then must belong to
$V\setminus L$ (by (\ref{pf.lem.eccent.tree-TL.0})). Hence, it is a vertex of
the subgraph $T\setminus L$. This shows that $T\setminus L$ has at least one
vertex.}. It remains to show that any two vertices of $T\setminus L$ are path-connected.

Let $u$ and $v$ be two vertices of $T\setminus L$. Then, $u\in V\setminus L$
and $v\in V\setminus L$. Hence, Lemma \ref{lem.eccent.tree-TL} \textbf{(b)}
yields%
\[
\left\{  \text{paths of }T\text{ from }u\text{ to }v\right\}  =\left\{
\text{paths of }T\setminus L\text{ from }u\text{ to }v\right\}  .
\]
Thus, $\left\{  \text{paths of }T\setminus L\text{ from }u\text{ to
}v\right\}  =\left\{  \text{paths of }T\text{ from }u\text{ to }v\right\}
\neq\varnothing$ (since there exists a path of $T$ from $u$ to $v$ (because
$T$ is connected)). In other words, there exists a path of $T\setminus L$ from
$u$ to $v$. In other words, $u$ and $v$ are path-connected in $T\setminus L$.

We have now shown that any two vertices $u$ and $v$ of $T\setminus L$ are
path-connected in $T\setminus L$. This entails that $T\setminus L$ is
connected (since $T\setminus L$ has at least one vertex). This proves Lemma
\ref{lem.eccent.tree-TL} \textbf{(a)}. \medskip

\textbf{(d)} If $u$ and $v$ are two vertices of $T\setminus L$, then the two
distances $d_{T}\left(  u,v\right)  $ and $d_{T\setminus L}\left(  u,v\right)
$ are equal (by Lemma \ref{lem.eccent.tree-TL} \textbf{(c)}); thus, we shall
denote both distances by $d\left(  u,v\right)  $ (since there is no confusion
to be afraid of).

Let $v\in V\setminus L$. We must show that $\operatorname*{ecc}\nolimits_{T}%
v=\operatorname*{ecc}\nolimits_{T\setminus L}v+1$.

Let $u$ be a vertex of $T\setminus L$ such that $d\left(  v,u\right)  $ is
maximum. Thus, $\operatorname*{ecc}\nolimits_{T\setminus L}v=d\left(
v,u\right)  $ (by the definition of $\operatorname*{ecc}\nolimits_{T\setminus
L}v$). However, $u$ is a vertex of $T\setminus L$, and thus does not belong to
$L$. Hence, $u$ is not a leaf of $T$ (since $L$ is the set of all leaves of
$T$). Hence, $u$ has degree $\geq2$ in $T$ (since a vertex in a tree with more
than $1$ vertex cannot have degree $0$).

Now, consider the path $\mathbf{p}$ from $v$ to $u$ in the tree $T$. This path
$\mathbf{p}$ has length $d\left(  v,u\right)  $. Since $u$ has degree $\geq2$,
there exist at least two edges of $T$ that contain $u$. Hence, in particular,
there exists at least one edge $f$ that contains $u$ and is distinct from the
last edge of $\mathbf{p}$\ \ \ \ \footnote{If the path $\mathbf{p}$ has no
edges, then $f$ can be any edge that contains $u$.}. Consider this edge $f$.
Let $w$ be the endpoint of $f$ other than $u$. Appending $f$ and $w$ to the
end of the path $\mathbf{p}$, we obtain a walk from $v$ to $w$. This walk is
backtrack-free (since $f$ is distinct from the last edge of $\mathbf{p}$) and
thus must be a path (by Proposition \ref{prop.btf-walk.cyc}, since $T$ has no
cycles). This path has length $d\left(  v,u\right)  +1$ (since it was obtained
by appending an edge to the path $\mathbf{p}$, which has length $d\left(
v,u\right)  $). Hence, $d\left(  v,w\right)  =d\left(  v,u\right)  +1$ (since
there is only one path from $v$ to $w$, and we just showed that this path has
length $d\left(  v,u\right)  +1$). But the definition of eccentricity yields%
\begin{equation}
\operatorname*{ecc}\nolimits_{T}v\geq d\left(  v,w\right)
=\underbrace{d\left(  v,u\right)  }_{=\operatorname*{ecc}\nolimits_{T\setminus
L}v}+\,1=\operatorname*{ecc}\nolimits_{T\setminus L}v+1.
\label{pf.lem.eccent.tree-TL.d.geq}%
\end{equation}

On the other hand, let $x$ be a vertex of $T$ such that $d\left(  v,x\right)
$ is maximum. Thus, $\operatorname*{ecc}\nolimits_{T}v=d\left(  v,x\right)  $
(by the definition of $\operatorname*{ecc}\nolimits_{T}v$). The path from $v$
to $x$ has length $\geq1$ (since otherwise, we would have $x=v$ and therefore
$d\left(  v,x\right)  =d\left(  v,v\right)  =0$, which would easily contradict
the maximality of $d\left(  v,x\right)  $). Thus, it has a second-to-last
vertex. Let $y$ be this second-to-last vertex. Then, the path from $v$ to $y$
is simply the path from $v$ to $x$ with its last edge removed. Consequently,
$d\left(  v,y\right)  =d\left(  v,x\right)  -1$. However, it is easy to see
that $y\in V\setminus L$\ \ \ \ \footnote{\textit{Proof.} Assume the contrary.
Thus, $y\notin V\setminus L$. Hence, $y\neq v$ (since $y\notin V\setminus L$
but $v\in V\setminus L$).
\par
However, $y$ is the second-to-last vertex of the path from $v$ to $x$.
Therefore, $y$ is either the starting point $v$ of this path, or an
intermediate vertex of this path. Since $y\neq v$, we thus conclude that $y$
is an intermediate vertex of this path. Hence, by
(\ref{pf.lem.eccent.tree-TL.0}), we see that $y$ must belong to $V\setminus
L$. But this contradicts $y\notin V\setminus L$. This contradiction shows that
our assumption was false, qed.}. In other words, $y$ is a vertex of
$T\setminus L$. Thus, the definition of eccentricity yields%
\[
\operatorname*{ecc}\nolimits_{T\setminus L}v\geq d\left(  v,y\right)
=\underbrace{d\left(  v,x\right)  }_{=\operatorname*{ecc}\nolimits_{T}%
v}-\,1=\operatorname*{ecc}\nolimits_{T}v-1,
\]
so that $\operatorname*{ecc}\nolimits_{T}v\leq\operatorname*{ecc}%
\nolimits_{T\setminus L}v+1$. Combining this with
(\ref{pf.lem.eccent.tree-TL.d.geq}), we obtain $\operatorname*{ecc}%
\nolimits_{T}v=\operatorname*{ecc}\nolimits_{T\setminus L}v+1$. This proves
Lemma \ref{lem.eccent.tree-TL} \textbf{(d)}. \medskip

\textbf{(e)} If $u$ and $v$ are two vertices of $T\setminus L$, then the two
distances $d_{T}\left(  u,v\right)  $ and $d_{T\setminus L}\left(  u,v\right)
$ are equal (by Lemma \ref{lem.eccent.tree-TL} \textbf{(c)}); thus, we shall
denote both distances by $d\left(  u,v\right)  $ (since there is no confusion
to be afraid of).

Let $v\in L$ be a leaf. Let $w$ be the unique neighbor of $v$ in $T$. We must
prove that $\operatorname*{ecc}\nolimits_{T}v=\operatorname*{ecc}%
\nolimits_{T}w+1$.

We first claim that%
\begin{equation}
d\left(  v,u\right)  =d\left(  w,u\right)  +1\ \ \ \ \ \ \ \ \ \ \text{for
each }u\in V\setminus\left\{  v\right\}  . \label{pf.lem.eccent.tree-TL.e.d}%
\end{equation}

[\textit{Proof of (\ref{pf.lem.eccent.tree-TL.e.d}):} We have $\deg v=1$
(since $v$ is a leaf). In other words, there is a unique edge of $T$ that
contains $v$. Let $e$ be this edge. The endpoints of $e$ are $v$ and $w$
(since $w$ is the unique neighbor of $v$). Thus, $v\neq w$ (since $T$ has no
loops) and $d\left(  v,w\right)  =1$.

Now, let $u\in V\setminus\left\{  v\right\}  $. Then, the path from $v$ to $u$
in $T$ must have length $\geq1$ (since $u\neq v$), and therefore must begin
with the edge $e$ (since $e$ is the only edge that contains $v$). If we remove
this edge $e$ from this path, we thus obtain a path from $w$ to $u$. As a
consequence, the path from $v$ to $u$ is longer by exactly $1$ edge than the
path from $w$ to $u$. In other words, we have $d\left(  v,u\right)  =d\left(
w,u\right)  +1$. This proves (\ref{pf.lem.eccent.tree-TL.e.d}).] \medskip

Now, the definition of eccentricity yields%
\begin{equation}
\operatorname*{ecc}\nolimits_{T}v=\max\left\{  d\left(  v,u\right)
\ \mid\ u\in V\right\}  . \label{pf.lem.eccent.tree-TL.e.eccTv}%
\end{equation}
This maximum is clearly \textbf{not} attained for $u=v$ (since $d\left(
v,v\right)  =0$ is smaller than $d\left(  v,w\right)  =1$). Thus, this maximum
does not change if we remove $v$ from its indexing set $V$. Hence,
(\ref{pf.lem.eccent.tree-TL.e.eccTv}) rewrites as%
\begin{align}
\operatorname*{ecc}\nolimits_{T}v  &  =\max\left\{  \underbrace{d\left(
v,u\right)  }_{\substack{=d\left(  w,u\right)  +1\\\text{(by
(\ref{pf.lem.eccent.tree-TL.e.d}))}}}\ \mid\ u\in V\setminus\left\{
v\right\}  \right\} \nonumber\\
&  =\max\left\{  d\left(  w,u\right)  +1\ \mid\ u\in V\setminus\left\{
v\right\}  \right\} \nonumber\\
&  =\max\left\{  d\left(  w,u\right)  \ \mid\ u\in V\setminus\left\{
v\right\}  \right\}  +1. \label{pf.lem.eccent.tree-TL.e.eccTv2}%
\end{align}

On the other hand, the definition of eccentricity yields%
\begin{equation}
\operatorname*{ecc}\nolimits_{T}w=\max\left\{  d\left(  w,u\right)
\ \mid\ u\in V\right\}  . \label{pf.lem.eccent.tree-TL.e.eccTw}%
\end{equation}
We shall now show that this maximum does not change if we remove $v$ from its
indexing set $V$. In other words, we shall show that%
\begin{equation}
\max\left\{  d\left(  w,u\right)  \ \mid\ u\in V\right\}  =\max\left\{
d\left(  w,u\right)  \ \mid\ u\in V\setminus\left\{  v\right\}  \right\}  .
\label{pf.lem.eccent.tree-TL.e.eccTw2}%
\end{equation}

[\textit{Proof of (\ref{pf.lem.eccent.tree-TL.e.eccTw2}):} Assume that
(\ref{pf.lem.eccent.tree-TL.e.eccTw2}) is false. Then, the maximum
$\max\left\{  d\left(  w,u\right)  \ \mid\ u\in V\right\}  $ is attained
\textbf{only} at $u=v$. In other words, we have%
\begin{equation}
d\left(  w,v\right)  >d\left(  w,u\right)  \ \ \ \ \ \ \ \ \ \ \text{for all
}u\in V\setminus\left\{  v\right\}  .
\label{pf.lem.eccent.tree-TL.e.eccTw2.pf.1}%
\end{equation}
However, the tree $T$ has more than $2$ vertices. Thus, it has a vertex $u$
that is distinct from both $v$ and $w$. Consider this $u$. Thus, $u\in
V\setminus\left\{  v\right\}  $, so that
(\ref{pf.lem.eccent.tree-TL.e.eccTw2.pf.1}) yields $d\left(  w,v\right)
>d\left(  w,u\right)  $. In view of $d\left(  w,v\right)  =d\left(
v,w\right)  =1$, this rewrites as $1>d\left(  w,u\right)  $, so that $d\left(
w,u\right)  <1$. Therefore, $w=u$. But this contradicts the fact that $w$ is
distinct from $u$. This contradiction shows that our assumption was false, and
thus (\ref{pf.lem.eccent.tree-TL.e.eccTw2}) is proved.] \medskip

Now, (\ref{pf.lem.eccent.tree-TL.e.eccTv2}) becomes%
\begin{align*}
\operatorname*{ecc}\nolimits_{T}v  &  =\underbrace{\max\left\{  d\left(
w,u\right)  \ \mid\ u\in V\setminus\left\{  v\right\}  \right\}
}_{\substack{=\max\left\{  d\left(  w,u\right)  \ \mid\ u\in V\right\}
\\\text{(by (\ref{pf.lem.eccent.tree-TL.e.eccTw2}))}}}+\,1\\
&  =\underbrace{\max\left\{  d\left(  w,u\right)  \ \mid\ u\in V\right\}
}_{\substack{=\operatorname*{ecc}\nolimits_{T}w\\\text{(by
(\ref{pf.lem.eccent.tree-TL.e.eccTw}))}}}+\,1=\operatorname*{ecc}%
\nolimits_{T}w+1.
\end{align*}
This proves Lemma \ref{lem.eccent.tree-TL} \textbf{(e)}.\medskip

\textbf{(f)} Lemma \ref{lem.eccent.tree-TL} \textbf{(e)} shows that any vertex
$v\in L$ has a higher eccentricity than its unique neighbor. Thus, a vertex
$v$ of $T$ that minimizes $\operatorname*{ecc}\nolimits_{T}v$ cannot belong to
$L$. In other words, a vertex $v$ of $T$ that minimizes $\operatorname*{ecc}%
\nolimits_{T}v$ must belong to $V\setminus L$.

However, the centers of $T$ are defined to be the vertices of $T$ that
minimize $\operatorname*{ecc}\nolimits_{T}v$. As we just proved, these
vertices must belong to $V\setminus L$. Thus, the centers of $T$ can also be
characterized as the vertices $v\in V\setminus L$ that minimize
$\operatorname*{ecc}\nolimits_{T}v$. However, a vertex $v\in V\setminus L$
minimizes $\operatorname*{ecc}\nolimits_{T}v$ if and only if it minimizes
$\operatorname*{ecc}\nolimits_{T\setminus L}v$ (because Lemma
\ref{lem.eccent.tree-TL} \textbf{(d)} yields $\operatorname*{ecc}%
\nolimits_{T}v=\operatorname*{ecc}\nolimits_{T\setminus L}v+1$ for any such
vertex $v$). Thus, we conclude that the centers of $T$ can be characterized as
the vertices $v\in V\setminus L$ that minimize $\operatorname*{ecc}%
\nolimits_{T\setminus L}v$. But this is precisely the definition of the
centers of $T\setminus L$. As a consequence, we see that the centers of $T$
are precisely the centers of $T\setminus L$. This proves Lemma
\ref{lem.eccent.tree-TL} \textbf{(f)}.
\end{proof}

\begin{proof}
[Proof of Theorem \ref{thm.eccent.tree-main}.]We shall prove parts
\textbf{(a)} and \textbf{(b)} of Theorem \ref{thm.eccent.tree-main} by strong
induction on $\left\vert \operatorname*{V}\left(  T\right)  \right\vert $:

\textit{Induction step:} Consider a tree $T$. Assume that parts \textbf{(a)}
and \textbf{(b)} of Theorem \ref{thm.eccent.tree-main} are true for any tree
with fewer than $\left\vert \operatorname*{V}\left(  T\right)  \right\vert $
many vertices. We must now prove these parts for our tree $T$.

If $\left\vert \operatorname*{V}\left(  T\right)  \right\vert \leq2$, then
both parts are obvious. Hence, WLOG assume that $\left\vert \operatorname*{V}%
\left(  T\right)  \right\vert >2$. Thus, the tree $T$ has more than $2$
vertices. Let $L$ be the set of all leaves of $T$. Note that $\left\vert
L\right\vert \geq2$ (since Theorem \ref{thm.tree.leaves.2} \textbf{(a)} shows
that any tree with at least $2$ vertices has at least $2$ leaves). Define the
multigraph $T\setminus L$ as in Lemma \ref{lem.eccent.tree-TL}. Then, Lemma
\ref{lem.eccent.tree-TL} \textbf{(f)} shows that the centers of $T$ are
precisely the centers of $T\setminus L$.

However, Lemma \ref{lem.eccent.tree-TL} \textbf{(a)} yields that $T\setminus
L$ is again a tree. This tree has fewer vertices than $T$ (since $\left\vert
L\right\vert \geq2>0$). Hence, by the induction hypothesis, both parts
\textbf{(a)} and \textbf{(b)} of Theorem \ref{thm.eccent.tree-main} are true
for the tree $T\setminus L$ instead of $T$. In other words, the tree
$T\setminus L$ has either $1$ or $2$ centers, and if it has $2$ centers, then
these $2$ centers are adjacent. Since the centers of $T$ are precisely the
centers of $T\setminus L$, we can rewrite this as follows: The tree $T$ has
either $1$ or $2$ centers, and if it has $2$ centers, then these $2$ centers
are adjacent. In other words, parts \textbf{(a)} and \textbf{(b)} of Theorem
\ref{thm.eccent.tree-main} hold for our tree $T$. This completes the induction
step. Thus, parts \textbf{(a)} and \textbf{(b)} of Theorem
\ref{thm.eccent.tree-main} are proved. \medskip

\textbf{(c)} This follows from Lemma \ref{lem.eccent.tree-TL} \textbf{(f)}.
Indeed, if $T$ has at most $2$ vertices, then all vertices of $T$ are centers
of $T$ (this is trivial to check). If not, then each \textquotedblleft
leaf-removal\textquotedblright\ step of our algorithm leaves the set of
centers of $T$ unchanged (by Lemma \ref{lem.eccent.tree-TL} \textbf{(f)}), and
thus the centers of the original tree $T$ are precisely the centers of the
tree that remains at the end of the algorithm. But the latter tree has at most
$2$ vertices, and thus its centers are precisely its vertices. So the centers
of $T$ are precisely the vertices that remain at the end of the algorithm.
Theorem \ref{thm.eccent.tree-main} \textbf{(c)} is proven.
\end{proof}

The following exercise shows another approach to the centers of a tree:

\begin{exercise}
\label{exe.6.2}Let $T$ be a tree. Let $\mathbf{p}=\left(  p_{0},\ast
,p_{1},\ast,p_{2},\ldots,\ast,p_{m}\right)  $ be a longest path of $T$. (We
write asterisks for the edges since we don't need to name them.)

Prove the following:

\begin{enumerate}
\item[\textbf{(a)}] If $m$ is even, then the only center of $T$ is $p_{m/2}$.

\item[\textbf{(b)}] If $m$ is odd, then the two centers of $T$ are $p_{\left(
m-1\right)  /2}$ and $p_{\left(  m+1\right)  /2}$.
\end{enumerate}
\end{exercise}

\begin{remark}
Exercise \ref{exe.6.2} is a result by Arthur Cayley from 1875 (see
\cite{Cayley75}). It shows once again that each tree has exactly one center or
two adjacent centers, and also shows that any two longest paths of a tree have
a common vertex.
\end{remark}

The notion of a \textbf{centroid} of a tree is a relative of the notion of a
center. We briefly discuss it in the following exercise:

\begin{exercise}
\label{exe.6.3}Let $T$ be a tree. For any vertex $v$ of $T$, we let $c_{v}$
denote the size of the largest component of the graph $T\setminus v$. (Recall
that $T\setminus v$ is the graph obtained from $T$ by removing the vertex $v$
and all edges that contain $v$. Note that a component (according to our
definition) is a set of vertices; thus, its size is the number of vertices in it.)

The vertices $v$ of $T$ that minimize the number $c_{v}$ are called the
\textbf{centroids} of $T$.

\begin{enumerate}
\item[\textbf{(a)}] Prove that $T$ has no more than two centroids, and
furthermore, if $T$ has two centroids, then these two centroids are adjacent.

\item[\textbf{(b)}] Find a tree $T$ such that the centroid(s) of $T$ are
distinct from the center(s) of $T$.
\end{enumerate}

[\textbf{Example:} Here is an example of a tree $T$, where each vertex $v$ is
labelled with the corresponding number $c_{v}$:%
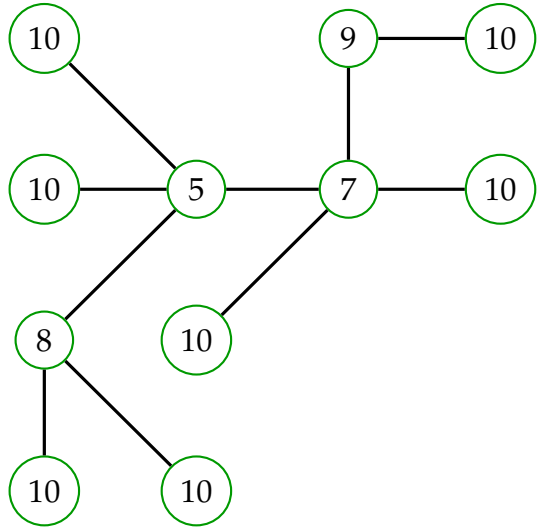
\[
\begin{tikzpicture}[scale=2]
\begin{scope}[every node/.style={circle,thick,draw=green!60!black}]
\node(20) at (2,0) {$10$};
\node(10) at (1,0) {$10$};
\node(11) at (1,1) {$8$};
\node(21) at (2,1) {$10$};
\node(12) at (1,2) {$10$};
\node(13) at (1,3) {$10$};
\node(22) at (2,2) {$5$};
\node(32) at (3,2) {$7$};
\node(33) at (3,3) {$9$};
\node(43) at (4,3) {$10$};
\node(42) at (4,2) {$10$};
\end{scope}
\begin{scope}[every edge/.style={draw=black,very thick}, every loop/.style={}]
\path[-] (11) edge (20) edge (10) edge (22);
\path[-] (22) edge (13) edge (12) edge (32);
\path[-] (32) edge (21) edge (42) edge (33);
\path[-] (43) edge (33);
\end{scope}
\end{tikzpicture}\ \ .
\]
Thus, the vertex labelled $5$ is the only centroid of this tree $T$.]
\end{exercise}

Note the analogy between Exercise \ref{exe.6.3} \textbf{(a)} and Theorem
\ref{thm.eccent.tree-main} \textbf{(a)} and \textbf{(b)}.

\subsection{Arborescences}

\subsubsection{Definitions}

Enough about undirected graphs.

What would be a directed analogue of a tree? I.e., what kind of digraphs play
the same role among digraphs that trees do among undirected graphs?

Trees are graphs that are connected and have no cycles. This suggests two
directed versions:

\begin{itemize}
\item We can study digraphs that are strongly connected and have no cycles.
Unfortunately, there is not much to study: Any such digraph has only $1$
vertex and no arcs. (Make sure you understand why!)

\item We can drop the connectedness requirement. Digraphs that have no cycles
are called \textbf{acyclic}, and more typically they are called \textbf{dags}
(short for \textquotedblleft directed acyclic graphs\textquotedblright).

However, these dags aren't quite like trees. For example, a tree always has
fewer edges than vertices, but a dag can have more arcs than
vertices.\footnote{For example, here is a dag with $4$ vertices and $5$ arcs:%
\[%
%
\ \ .
\]
}
\end{itemize}

Here is a more convincing analogue of trees for digraphs:\footnote{We recall
that we defined a multigraph $D^{\operatorname*{und}}$ for every multidigraph
$D$ (in Definition \ref{def.mdg.Dund}). Roughly speaking, this multigraph
$D^{\operatorname*{und}}$ is obtained by \textquotedblleft forgetting the
directions\textquotedblright\ of the arcs of $D$. Parallel arcs are not merged
into one. For example, $%
%
$ .}

\begin{definition}
\label{def.arbor.arbor-from}Let $D$ be a multidigraph. Let $r$ be a vertex of
$D$.

\begin{enumerate}
\item[\textbf{(a)}] We say that $r$ is a \textbf{from-root} (or, short,
\textbf{root}) of $D$ if for each vertex $v$ of $D$, the digraph $D$ has a
path from $r$ to $v$.

\item[\textbf{(b)}] We say that $D$ is an \textbf{arborescence rooted from
}$r$ if $r$ is a from-root of $D$ and the undirected multigraph
$D^{\operatorname*{und}}$ has no cycles. (Recall that $D^{\operatorname*{und}%
}$ is the multigraph obtained from $D$ by turning each arc into an undirected
edge. Parallel arcs are not merged into one!)
\end{enumerate}
\end{definition}

Of course, there are analogous notions of a \textquotedblleft
to-root\textquotedblright\ and an \textquotedblleft arborescence rooted
towards $r$\textquotedblright, but these are just the same notions that we
just defined with all arrows reversed. So we need not study them separately;
we can just take any property of \textquotedblleft rooted
from\textquotedblright\ and reverse all arcs to make it into a property of
\textquotedblleft rooted to\textquotedblright.

\begin{example}
The multidigraph%
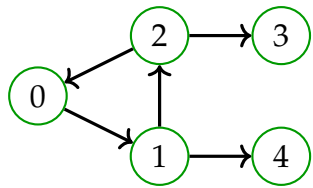
\[%
%
\]
has three from-roots (namely, $0$, $1$ and $2$). It is not an arborescence
rooted from any of them, because turning each arc into an undirected edge
yields a graph with a cycle.

If we reverse the arc from $0$ to $1$, then we obtain a multidigraph%
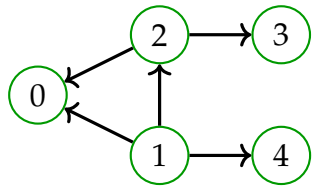
\[%
%
\]
which has only one from-root (namely, $1$) and is still not an arborescence
(for the same reason as before).
\end{example}

\begin{example}
Consider the following multidigraph:%
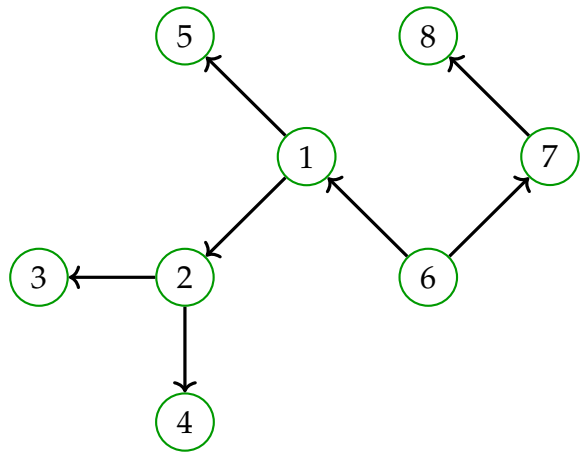
\[%
%
\ \ .
\]
This is an arborescence rooted from $6$. Indeed, it has paths from $6$ to all
vertices, and turning each arc into an undirected edge yields a tree.

If we reverse the arc from $1$ to $2$, we obtain a multidigraph%
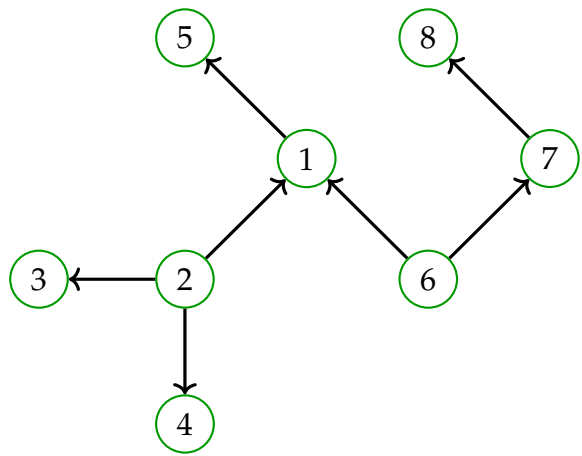
\[%
%
\ \ ,
\]
which is \textbf{not} an arborescence, because it has no from-root anymore.
\end{example}

\subsubsection{Arborescences vs. trees: statement}

The above examples suggest that an arborescence rooted from $r$ is basically
the same as a tree, whose all edges have been \textquotedblleft oriented away
from $r$\textquotedblright. More precisely:

\begin{theorem}
\label{thm.arbor.vs-tree}Let $D$ be a multidigraph, and let $r$ be a vertex of
$D$. Then, the following two statements are equivalent:

\begin{itemize}
\item \textbf{Statement C1:} The multidigraph $D$ is an arborescence rooted
from $r$.

\item \textbf{Statement C2:} The undirected multigraph $D^{\operatorname*{und}%
}$ is a tree, and each arc of $D$ is \textquotedblleft oriented away from
$r$\textquotedblright\ (this means the following: the source of this arc lies
on the unique path between $r$ and the target of this arc on
$D^{\operatorname*{und}}$).
\end{itemize}
\end{theorem}

This is an easy theorem to believe, but an annoyingly hard one to formally
prove in full detail! We shall prove this theorem later.

\subsubsection{The arborescence equivalence theorem}

First, let us show another bunch of equivalent criteria for arborescences,
imitating the tree equivalence theorem (Theorem \ref{thm.trees.T1-8}):

\begin{theorem}
[The arborescence equivalence theorem]\label{thm.arbor.eq-A}Let $D=\left(
V,A,\psi\right)  $ be a multidigraph with a from-root $r$. Then, the following
six statements are equivalent:

\begin{itemize}
\item \textbf{Statement A1:} The multidigraph $D$ is an arborescence rooted
from $r$.

\item \textbf{Statement A2:} We have $\left\vert A\right\vert =\left\vert
V\right\vert -1$.

\item \textbf{Statement A3:} The multigraph $D^{\operatorname*{und}}$ is a tree.

\item \textbf{Statement A4:} For each vertex $v\in V$, the multidigraph $D$
has a unique walk from $r$ to $v$.

\item \textbf{Statement A5:} If we remove any arc from $D$, then the vertex
$r$ will no longer be a from-root of the resulting multidigraph.

\item \textbf{Statement A6:} We have $\deg^{-}r=0$, and each $v\in
V\setminus\left\{  r\right\}  $ satisfies $\deg^{-}v=1$.
\end{itemize}
\end{theorem}

\begin{proof}
We will prove the implications A1$\Longrightarrow$A4$\Longrightarrow
$A5$\Longrightarrow$A6$\Longrightarrow$A2$\Longrightarrow$A3$\Longrightarrow
$A1. Since these implications form a cycle that includes all six statements,
this will entail that all six statements are equivalent. \medskip

Before we prove these implications, we introduce a notation: If $a$ is any arc
of $D$, then $D\setminus a$ shall denote the multidigraph obtained from $D$ by
removing this arc $a$. (Formally, this means that $D\setminus a:=\left(
V,\ A\setminus\left\{  a\right\}  ,\ \psi\mid_{A\setminus\left\{  a\right\}
}\right)  $.)

We now come to the proofs of the promised implications. \medskip

\textit{Proof of the implication A1}$\Longrightarrow$\textit{A4:} Assume that
Statement A1 holds. Thus, $D$ is an arborescence rooted from $r$. In other
words, $r$ is a from-root of $D$ and the undirected multigraph
$D^{\operatorname*{und}}$ has no cycles.

We must show that for each vertex $v\in V$, the multidigraph $D$ has a unique
walk from $r$ to $v$. The existence of such a walk is clear (because $r$ is a
from-root of $D$). It is the uniqueness that we need to prove.

Assume the contrary. Thus, there exists a vertex $v\in V$ such that two
distinct walks $\mathbf{u}$ and $\mathbf{v}$ from $r$ to $v$ exist. However,
the multigraph $D$ has no loops (since any loop of $D$ would be a loop of
$D^{\operatorname*{und}}$, and thus create a cycle of $D^{\operatorname*{und}%
}$, but we know that $D^{\operatorname*{und}}$ has no cycles). Hence, any walk
of $D$ is automatically a backtrack-free walk of $D^{\operatorname*{und}}$
(indeed, it is backtrack-free because the only way two consecutive arcs of a
walk in a \textbf{digraph} can be equal is if they are loops). Therefore, the
two walks $\mathbf{u}$ and $\mathbf{v}$ of $D$ are two backtrack-free walks of
$D^{\operatorname*{und}}$. Thus, there are two distinct backtrack-free walks
from $r$ to $v$ in $D^{\operatorname*{und}}$ (namely, $\mathbf{u}$ and
$\mathbf{v}$). Theorem \ref{thm.btf-walk.two-cyc} thus lets us conclude that
$D^{\operatorname*{und}}$ has a cycle. But this contradicts the fact that
$D^{\operatorname*{und}}$ has no cycles.

This contradiction shows that our assumption was wrong. Hence, we have proved
that for each vertex $v\in V$, the multidigraph $D$ has a unique walk from $r$
to $v$. In other words, Statement A4 holds. \medskip

\textit{Proof of the implication A4}$\Longrightarrow$\textit{A5:} Assume that
Statement A4 holds.

Let now $a$ be any arc of $D$. We shall show that $r$ is not a from-root of
the multidigraph $D\setminus a$.

Indeed, let $s$ be the source and $t$ the target of the arc $a$. We shall show
that the digraph $D\setminus a$ has no path from $r$ to $t$.

Indeed, assume the contrary. Thus, $D\setminus a$ has some path $\mathbf{p}$
from $r$ to $t$. This path does not use the arc $a$ (since it is a path of
$D\setminus a$).

On the other hand, we have assumed that Statement A4 holds. Applying this
statement to $v=s$, we conclude that the multidigraph $D$ has a unique walk
from $r$ to $s$. Let $\left(  v_{0},a_{1},v_{1},a_{2},v_{2},\ldots,a_{k}%
,v_{k}\right)  $ be this walk. By appending the arc $a$ and the vertex $t$ to
its end, we extend it to a longer walk%
\[
\left(  v_{0},a_{1},v_{1},a_{2},v_{2},\ldots,a_{k},v_{k},a,t\right)  ,
\]
which is a walk from $r$ to $t$. We denote this walk by $\mathbf{q}$.

We have now found two walks from $r$ to $t$ in the digraph $D$: namely, the
path $\mathbf{p}$ and the walk $\mathbf{q}$. These two walks are distinct
(since $\mathbf{q}$ uses the arc $a$, but $\mathbf{p}$ does not). However,
Statement A4 (applied to $v=t$) yields that the multidigraph $D$ has a
\textbf{unique} walk from $r$ to $t$. This contradicts the fact that we just
have found two distinct such walks.

This contradiction shows that our assumption was false. Hence, the digraph
$D\setminus a$ has no path from $r$ to $t$. Thus, $r$ is not a from-root of
$D\setminus a$.

Forget that we fixed $a$. We have now proved that if $a$ is any arc of $D$,
then $r$ is not a from-root of $D\setminus a$. In other words, if we remove
any arc from $D$, then the vertex $r$ will no longer be a from-root of the
resulting multidigraph. Thus, Statement A5 holds. \medskip

\textit{Proof of the implication A5}$\Longrightarrow$\textit{A6:} Assume that
Statement A5 holds. We must prove that Statement A6 holds. In other words, we
must prove that $\deg^{-}r=0$, and that each $v\in V\setminus\left\{
r\right\}  $ satisfies $\deg^{-}v=1$.

Let us first prove that $\deg^{-}r=0$. Indeed, assume the contrary. Thus,
$\deg^{-}r\neq0$, so that there exists an arc $a$ with target $r$. We shall
show that $r$ is a from-root of $D\setminus a$.

The arc $a$ has target $r$. Thus, a path that starts at $r$ cannot use this
arc $a$ (because this arc would lead it back to $r$, but a path is not allowed
to revisit any vertex), and therefore must be a path of $D\setminus a$. Thus
we have shown that any path of $D$ that starts at $r$ is also a path of
$D\setminus a$. However, for each vertex $v$ of $D$, the digraph $D$ has a
path from $r$ to $v$ (since $r$ is a from-root of $D$). This path is also a
path of $D\setminus a$ (since any path of $D$ that starts at $r$ is also a
path of $D\setminus a$). Thus, for each vertex $v$ of $D\setminus a$, the
digraph $D\setminus a$ has a path from $r$ to $v$. In other words, $r$ is a
from-root of $D\setminus a$. However, we have assumed that Statement A5 holds.
Thus, in particular, if we remove the arc $a$ from $D$, then the vertex $r$
will no longer be a from-root of the resulting multidigraph. In other words,
$r$ is not a from-root of $D\setminus a$. But this contradicts the fact that
$r$ is a from-root of $D\setminus a$.

This contradiction shows that our assumption was false. Hence, $\deg^{-}r=0$
is proved.

Now, let $v\in V\setminus\left\{  r\right\}  $ be arbitrary. We must show that
$\deg^{-}v=1$.

Indeed, assume the contrary. Thus, $\deg^{-}v\neq1$. Using the fact that $r$
is a from-root of $D$, it is thus easy to see that $\deg^{-}v\geq
2$\ \ \ \ \footnote{\textit{Proof.} Since $r$ is a from-root of $D$, we know
that the digraph $D$ has a path from $r$ to $v$. Since $v\neq r$ (because
$v\in V\setminus\left\{  r\right\}  $), this path must have at least one arc.
The last arc of this path is clearly an arc with target $v$. Thus, there
exists at least one arc with target $v$. In other words, $\deg^{-}v\geq1$.
Combining this with $\deg^{-}v\neq1$, we obtain $\deg^{-}v>1$. In other words,
$\deg^{-}v\geq2$.}. Hence, there exist two distinct arcs $a$ and $b$ with
target $v$. Consider these arcs $a$ and $b$.

We are in one of the following three cases:

\textit{Case 1:} The digraph $D\setminus a$ has a path from $r$ to $v$.

\textit{Case 2:} The digraph $D\setminus b$ has a path from $r$ to $v$.

\textit{Case 3:} Neither the digraph $D\setminus a$ nor the digraph
$D\setminus b$ has a path from $r$ to $v$.

Let us first consider Case 1. In this case, the digraph $D\setminus a$ has a
path from $r$ to $v$. Let $\mathbf{p}$ be such a path.

We have assumed that Statement A5 holds. Thus, in particular, if we remove the
arc $a$ from $D$, then the vertex $r$ will no longer be a from-root of the
resulting multidigraph. In other words, $r$ is not a from-root of $D\setminus
a$. In other words, there exists a vertex $w\in V$ such that the digraph
$D\setminus a$ has no path from $r$ to $w$ (by the definition of a
\textquotedblleft from-root\textquotedblright). Consider this vertex $w$.

The digraph $D$ has a path $\mathbf{q}$ from $r$ to $w$ (since $r$ is a
from-root of $D$). Consider this path $\mathbf{q}$. If the path $\mathbf{q}$
did not use the arc $a$, then it would be a path of $D\setminus a$ as well,
but this would contradict the fact that $D\setminus a$ has no path from $r$ to
$w$. Thus, the path $\mathbf{q}$ must use the arc $a$.

Consider the part of $\mathbf{q}$ that comes after the arc $a$. This part must
be a path from $v$ to $w$ (since the arc $a$ has target $v$, whereas the path
$\mathbf{q}$ has ending point $w$). Let us denote this path by $\mathbf{q}%
^{\prime}$. Thus, the path $\mathbf{q}^{\prime}$ does not use the arc $a$
(since it was defined as the part of $\mathbf{q}$ that comes after $a$).
Hence, $\mathbf{q}^{\prime}$ is a path of $D\setminus a$.

Now, we know that the digraph $D\setminus a$ has a path $\mathbf{p}$ from $r$
to $v$ as well as a path $\mathbf{q}^{\prime}$ from $v$ to $w$. Splicing these
paths together, we obtain a walk $\mathbf{p}\ast\mathbf{q}^{\prime}$ from $r$
to $w$. So we know that $D\setminus a$ has a walk from $r$ to $w$. According
to Corollary \ref{cor.mg.walk-thus-path}, we thus conclude that $D\setminus a$
has a path from $r$ to $w$. This contradicts the fact that $D\setminus a$ has
no path from $r$ to $w$.

We have thus obtained a contradiction in Case 1.

The same argument (but with the roles of $a$ and $b$ interchanged) results in
a contradiction in Case 2.

Let us finally consider Case 3. In this case, neither the digraph $D\setminus
a$ nor the digraph $D\setminus b$ has a path from $r$ to $v$. However, the
digraph $D$ has a path $\mathbf{p}$ from $r$ to $v$ (since $r$ is a from-root
of $D$). Consider this path $\mathbf{p}$. If this path $\mathbf{p}$ did not
use the arc $a$, then it would be a path of $D\setminus a$, but this would
contradict our assumption that the digraph $D\setminus a$ has no path from $r$
to $v$. Thus, this path $\mathbf{p}$ must use the arc $a$. For a similar
reason, it must also use the arc $b$. However, the two arcs $a$ and $b$ have
the same target (viz., $v$) and thus cannot both appear in the same path
(since a path cannot visit a vertex more than once). This contradicts the fact
that the path $\mathbf{p}$ uses both arcs $a$ and $b$. Hence, we have found a
contradiction in Case 3.

We have now found contradictions in all three Cases 1, 2 and 3. This
contradiction shows that our assumption was false. Hence, $\deg^{-}v=1$ is proved.

We have now proved that each $v\in V\setminus\left\{  r\right\}  $ satisfies
$\deg^{-}v=1$. Since we have also shown that $\deg^{-}r=0$, we thus have
proved Statement A6. \medskip

\textit{Proof of the implication A6}$\Longrightarrow$\textit{A2:} Assume that
Statement A6 holds. We must prove that Statement A2 holds. However,
Proposition \ref{prop.mdg.sum-deg} yields%
\begin{align*}
\left\vert A\right\vert  &  =\sum_{v\in V}\deg^{-}v=\underbrace{\deg^{-}%
r}_{\substack{=0\\\text{(by Statement A6)}}}+\sum_{v\in V\setminus\left\{
r\right\}  }\underbrace{\deg^{-}v}_{\substack{=1\\\text{(by Statement A6)}}}\\
&  =0+\sum_{v\in V\setminus\left\{  r\right\}  }1=\sum_{v\in V\setminus
\left\{  r\right\}  }1=\left\vert V\setminus\left\{  r\right\}  \right\vert
=\left\vert V\right\vert -1.
\end{align*}
Hence, Statement A2 holds. \medskip

\textit{Proof of the implication A2}$\Longrightarrow$\textit{A3:} Assume that
Statement A2 holds. We must prove that Statement A3 holds.

For each $v\in V$, the digraph $D$ has a path from $r$ to $v$ (since $r$ is a
from-root of $D$). Thus, for each $v\in V$, the graph $D^{\operatorname*{und}%
}$ has a path from $r$ to $v$ (since any path of $D$ is a path of
$D^{\operatorname*{und}}$). Therefore, any two vertices $u$ and $v$ of
$D^{\operatorname*{und}}$ are path-connected in $D^{\operatorname*{und}}$
(because we can get from $u$ to $v$ via $r$, according to the previous
sentence). Therefore, the graph $D^{\operatorname*{und}}$ is connected (since
it has at least one vertex\footnote{This is because $r\in V$.}). Moreover, its
number of edges is $\left\vert A\right\vert =\left\vert V\right\vert -1$ (by
Statement A2). Therefore, the multigraph $D^{\operatorname*{und}}$ satisfies
the Statement T4 of the tree equivalence theorem (Theorem \ref{thm.trees.T1-8}%
). Consequently, it satisfies Statement T1 of that theorem as well. In other
words, it is a tree. This proves Statement A3. \medskip

\textit{Proof of the implication A3}$\Longrightarrow$\textit{A1:} Assume that
Statement A3 holds. We must prove that Statement A1 holds.

The multigraph $D^{\operatorname*{und}}$ is a tree (by Statement A3), and thus
is a forest; hence, it has no cycles. Since we also know that $r$ is a
from-root of $D$, we thus conclude that $D$ is an arborescence rooted from $r$
(by the definition of an arborescence). In other words, Statement A1 is
satisfied. \medskip

We have now proved all six implications in the chain \newline
A1$\Longrightarrow$A4$\Longrightarrow$A5$\Longrightarrow$A6$\Longrightarrow
$A2$\Longrightarrow$A3$\Longrightarrow$A1. Thus, all six statements A1, A2,
$\ldots$, A6 are equivalent. This proves Theorem \ref{thm.arbor.eq-A}.
\end{proof}

\begin{exercise}
\label{exe.6.4}Let $D=\left(  V,A,\psi\right)  $ be a multidigraph that has no
cycles\footnotemark. Let $r\in V$ be some vertex of $D$. Prove the following:

\begin{enumerate}
\item[\textbf{(a)}] If $\deg^{-} u > 0$ holds for all $u \in V \setminus
\left\{  r \right\}  $, then $r$ is a from-root of $D$.

\item[\textbf{(b)}] If $\deg^{-}u=1$ holds for all $u\in V\setminus\left\{
r\right\}  $, then $D$ is an arborescence rooted from $r$.
\end{enumerate}
\end{exercise}

\footnotetext{Recall that cycles in a digraph have to be directed cycles --
i.e., each arc is traversed from its source to its target.}

\subsection{\label{sec.arbor.vs-tree}Arborescences vs. trees}

Our next goal is to prove Theorem \ref{thm.arbor.vs-tree}, which connects
arborescences with trees.

To prove it formally, we introduce a few notations regarding trees. First, we
recall the notion of a distance (Definition \ref{def.mg.dist}). We claim the
following simple property of distances in trees:

\begin{proposition}
\label{prop.tree.edge-d}Let $T=\left(  V,E,\varphi\right)  $ be a tree. Let
$r\in V$ be a vertex of $T$. Let $e$ be an edge of $T$, and let $u$ and $v$ be
its two endpoints. Then, the distances $d\left(  r,u\right)  $ and $d\left(
r,v\right)  $ differ by exactly $1$ (that is, we have either $d\left(
r,u\right)  =d\left(  r,v\right)  +1$ or $d\left(  r,v\right)  =d\left(
r,u\right)  +1$).
\end{proposition}

\begin{proof}
We recall that since $T$ is a tree, the distance $d\left(  p,q\right)  $
between two vertices $p$ and $q$ of $T$ is simply the length of the path from
$p$ to $q$. (This path is unique, since $T$ is a tree.)

Let $\mathbf{p}$ be the path from $r$ to $u$. Then, we are in one of the
following two cases:

\textit{Case 1:} The edge $e$ is an edge of $\mathbf{p}$.

\textit{Case 2:} The edge $e$ is not an edge of $\mathbf{p}$.

Consider Case 1. In this case, $e$ must be the \textbf{last} edge of
$\mathbf{p}$ (since otherwise, $\mathbf{p}$ would visit $u$ more than once,
but $\mathbf{p}$ cannot do this, since $\mathbf{p}$ is a path). Thus, if we
remove this last edge $e$ (and the vertex $u$) from $\mathbf{p}$, then we
obtain a path from $r$ to $v$. This path is exactly one edge shorter than
$\mathbf{p}$. Thus, $d\left(  r,v\right)  =d\left(  r,u\right)  -1$, so that
$d\left(  r,u\right)  =d\left(  r,v\right)  +1$. So we are done in Case 1.

Now, consider Case 2. In this case, the edge $e$ is not an edge of
$\mathbf{p}$. Thus, we can append $e$ and $v$ to the end of the path
$\mathbf{p}$, and the result will be a backtrack-free walk $\mathbf{p}%
^{\prime}$. However, a backtrack-free walk in a tree is always a path (since
otherwise, it would contain a cycle\footnote{by Proposition
\ref{prop.btf-walk.cyc}}, but a tree has no cycles). Thus, $\mathbf{p}%
^{\prime}$ is a path from $r$ to $v$, and it is exactly one edge longer than
$\mathbf{p}$ (by its construction). Therefore, $d\left(  r,v\right)  =d\left(
r,u\right)  +1$. So we are done in Case 2.

Now, we are done in both cases, so that Proposition \ref{prop.tree.edge-d} is proven.
\end{proof}

\begin{definition}
\label{def.tree.par-child-e}Let $T=\left(  V,E,\varphi\right)  $ be a tree.
Let $r\in V$ be a vertex of $T$. Let $e$ be an edge of $T$. By Proposition
\ref{prop.tree.edge-d}, the distances from the two endpoints of $e$ to the
vertex $r$ differ by exactly $1$. So one of them is smaller than the other.

\begin{enumerate}
\item[\textbf{(a)}] We define the $r$\textbf{-parent} of $e$ to be the
endpoint of $e$ whose distance to $r$ is the smallest. We denote this endpoint
by $e^{-r}$.

\item[\textbf{(b)}] We define the $r$\textbf{-child} of $e$ to be the endpoint
of $e$ whose distance to $r$ is the largest. We denote this endpoint by
$e^{+r}$.
\end{enumerate}

Thus, by Proposition \ref{prop.tree.edge-d}, we have%
\[
d\left(  r,e^{+r}\right)  =d\left(  r,e^{-r}\right)  +1.
\]

\end{definition}

\begin{example}
\label{exa.tree.par-child-e.1}Here is a tree $T$, a vertex $r$, an edge $e$
and its $r$-parent $e^{-r}$ and its $r$-child $e^{+r}$:%
\[%
\begin{tikzpicture}[scale=2]
\begin{scope}[every node/.style={circle,thick,draw=green!60!black}]
\node(1) at (0,2) {$e^{+r}$};
\node(5) at (0,1) {$\phantom{e^{-r}}$};
\node(2) at (1.5,2) {$e^{-r}$};
\node(6) at (1,1) {$\phantom{e^{-r}}$};
\node(7) at (2,1) {$\phantom{e^{-r}}$};
\node(8) at (1,0) {$r\vphantom{e^{-r}}$};
\end{scope}
\begin{scope}[every edge/.style={draw=black,very thick}, every loop/.style={}]
\path[-] (1) edge (5) edge node[above] {$e$}
(2) (6) edge (2) edge (8) (7) edge (2);
\end{scope}
\end{tikzpicture}%
\]

\end{example}

\begin{definition}
\label{def.tree.rto}Let $T=\left(  V,E,\varphi\right)  $ be a tree. Let $r\in
V$ be a vertex of $T$. Then, we define a multidigraph $T^{r\rightarrow}$ by%
\[
T^{r\rightarrow}:=\left(  V,E,\psi\right)  ,
\]
where $\psi:E\rightarrow V\times V$ is the map that sends each edge $e\in E$
to the pair $\left(  e^{-r},e^{+r}\right)  $. Colloquially speaking, this
means that $T^{r\rightarrow}$ is the multidigraph obtained from $T$ by turning
each edge $e$ into an arc from its $r$-parent $e^{-r}$ to its $r$-child
$e^{+r}$. This is what we mean when we speak of \textquotedblleft orienting
each edge of $T$ away from $r$\textquotedblright\ in Theorem
\ref{thm.arbor.vs-tree}.
\end{definition}

\begin{example}
If $T$ is the tree from Example \ref{exa.tree.par-child-e.1}, then
$T^{r\rightarrow}$ is the following multidigraph:%
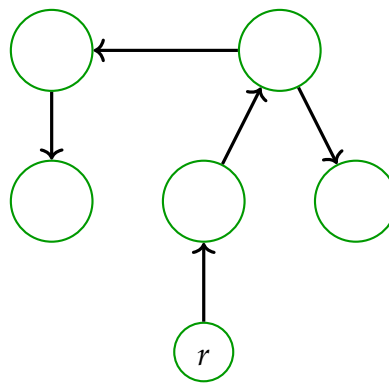
\[%
\begin{tikzpicture}[scale=2]
\begin{scope}[every node/.style={circle,thick,draw=green!60!black}]
\node(1) at (0,2) {$\phantom{e^{-r}}$};
\node(5) at (0,1) {$\phantom{e^{-r}}$};
\node(2) at (1.5,2) {$\phantom{e^{-r}}$};
\node(6) at (1,1) {$\phantom{e^{-r}}$};
\node(7) at (2,1) {$\phantom{e^{-r}}$};
\node(8) at (1,0) {$r\vphantom{e^{-r}}$};
\end{scope}
\begin{scope}[every edge/.style={draw=black,very thick}, every loop/.style={}]
\path[<-] (5) edge (1) (1) edge (2) (2) edge (6) (6) edge (8) (7) edge (2);
\end{scope}
\end{tikzpicture}%
\]

\end{example}

Now, Theorem \ref{thm.arbor.vs-tree} can be rewritten as follows:

\begin{theorem}
\label{thm.arbor.vs-tree.proper}Let $D$ be a multidigraph, and let $r$ be a
vertex of $D$. Then, the following two statements are equivalent:

\begin{itemize}
\item \textbf{Statement C1:} The multidigraph $D$ is an arborescence rooted
from $r$.

\item \textbf{Statement C2:} The undirected multigraph $D^{\operatorname*{und}%
}$ is a tree, and we have $D=\left(  D^{\operatorname*{und}}\right)
^{r\rightarrow}$. (This is a honest equality, not just some isomorphism.)
\end{itemize}
\end{theorem}

The proof of this theorem is best organized by splitting into two lemmas:

\begin{lemma}
\label{lem.arbor-vs-tree.1}Let $T=\left(  V,E,\varphi\right)  $ be a tree. Let
$r\in V$ be a vertex of $T$. Then, the multidigraph $T^{r\rightarrow}$ is an
arborescence rooted from $r$.
\end{lemma}

\begin{proof}
The idea is to show that if $\mathbf{p}$ is a path from $r$ to some vertex $v$
in the tree $T$, then $\mathbf{p}$ is also a path in the digraph
$T^{r\rightarrow}$, because all the edges of $\mathbf{p}$ have been
\textquotedblleft oriented correctly\textquotedblright\ (i.e., their
orientation matches how they are used in $\mathbf{p}$).

Here are the details: Clearly, $\left(  T^{r\rightarrow}\right)
^{\operatorname*{und}}=T$. Hence, the graph $\left(  T^{r\rightarrow}\right)
^{\operatorname*{und}}$ is a tree and hence has no cycles. Thus, it suffices
to prove that $r$ is a from-root of $T^{r\rightarrow}$. In other words, we
must prove that
\begin{equation}
T^{r\rightarrow}\text{ has a path from }r\text{ to }v
\label{pf.lem.arbor-vs-tree.1.goal}%
\end{equation}
for each $v\in V$.

We shall prove (\ref{pf.lem.arbor-vs-tree.1.goal}) by induction on $d\left(
r,v\right)  $ (where $d$ means the distance on the tree $T$):

\textit{Base case:} If $v\in V$ satisfies $d\left(  r,v\right)  =0$, then
$v=r$, and thus $T^{r\rightarrow}$ has a path from $r$ to $v$ (namely, the
trivial path $\left(  r\right)  $). Thus, (\ref{pf.lem.arbor-vs-tree.1.goal})
is proved for $d\left(  r,v\right)  =0$.

\textit{Induction step:} Let $k\in\mathbb{N}$. Assume (as the induction
hypothesis) that (\ref{pf.lem.arbor-vs-tree.1.goal}) holds for each $v\in V$
satisfying $d\left(  r,v\right)  =k$. We must now prove the same for each
$v\in V$ satisfying $d\left(  r,v\right)  =k+1$.

So let $v\in V$ satisfy $d\left(  r,v\right)  =k+1$. Then, the path of $T$
from $r$ to $v$ has length $k+1$. Let $\mathbf{p}$ be this path, let $e$ be
its last edge, and let $u$ be its second-to-last vertex (so that its last edge
$e$ has endpoints $u$ and $v$). Then, by removing the last edge $e$ from the
path $\mathbf{p}$, we obtain a path from $r$ to $u$ that is one edge shorter
than $\mathbf{p}$. Hence, $d\left(  r,u\right)  =d\left(  r,v\right)
-1<d\left(  r,v\right)  $. Consequently, the edge $e$ has $r$-parent $u$ and
$r$-child $v$ (by Definition \ref{def.tree.par-child-e}). In other words,
$e^{-r}=u$ and $e^{+r}=v$. Therefore, in the digraph $T^{r\rightarrow}$, the
edge $e$ is an arc from $u$ to $v$ (by Definition \ref{def.tree.rto}).
Moreover, we have $d\left(  r,u\right)  =d\left(  r,v\right)  -1=k$ (since
$d\left(  r,v\right)  =k+1$); therefore, the induction hypothesis tells us
that (\ref{pf.lem.arbor-vs-tree.1.goal}) holds for $u$ instead of $v$. In
other words, $T^{r\rightarrow}$ has a path from $r$ to $u$. Attaching the arc
$e$ and the vertex $v$ to this path, we obtain a walk of $T^{r\rightarrow}$
from $r$ to $v$ (since $e$ is an arc from $u$ to $v$ in $T^{r\rightarrow}$).
Thus, the digraph $T^{r\rightarrow}$ has a walk from $r$ to $v$, therefore
also a path from $r$ to $v$. Hence, (\ref{pf.lem.arbor-vs-tree.1.goal}) holds
for our $v$. This completes the induction step.

Thus, (\ref{pf.lem.arbor-vs-tree.1.goal}) is proved by induction. As we
explained above, this yields Lemma \ref{lem.arbor-vs-tree.1}.
\end{proof}

\begin{lemma}
\label{lem.arbor-vs-tree.2}Let $D=\left(  V,A,\psi\right)  $ be an
arborescence rooted from $r$ (for some $r\in V$). Let $a\in A$ be an arc of
$D$. Let $s$ be the source of $a$, and let $t$ be the target of $a$. Then:

\begin{enumerate}
\item[\textbf{(a)}] We have $d\left(  r,s\right)  <d\left(  r,t\right)  $,
where $d$ means distance on the tree $D^{\operatorname*{und}}$.

\item[\textbf{(b)}] In the multidigraph $\left(  D^{\operatorname*{und}%
}\right)  ^{r\rightarrow}$, the arc $a$ has source $s$ and target $t$.
\end{enumerate}
\end{lemma}

\begin{proof}
\textbf{(a)} The vertex $r$ is a from-root of $D$ (since $D$ is an
arborescence rooted from $r$). Thus, $D$ has a path from $r$ to $t$. Let
$\mathbf{p}$ be this path. Note that $\deg^{-}t\geq1$, since $t$ is the target
of at least one arc (namely, of $a$).

The digraph $D$ is an arborescence rooted from $r$, and thus satisfies
Statement A6 in the arborescence equivalence theorem (Theorem
\ref{thm.arbor.eq-A}). In other words, we have%
\[
\deg^{-}r=0\ \ \ \ \ \ \ \ \ \ \text{and}\ \ \ \ \ \ \ \ \ \ \deg^{-}v=1\text{
for each }v\in V\setminus\left\{  r\right\}  .
\]
In particular, this entails $\deg^{-}v\leq1$ for each $v\in V$. Applying this
to $v=t$, we obtain $\deg^{-}t\leq1$. Hence, the arc $a$ is the \textbf{only}
arc whose target is $t$.

We have $t\neq r$ (since $\deg^{-}r=0$ but $\deg^{-}t\geq1>0$). Thus, the path
$\mathbf{p}$ from $r$ to $t$ has at least one arc. Its last arc is therefore
an arc whose target is $t$. Hence, this last arc is $a$ (since $a$ is the
\textbf{only} arc whose target is $t$).

If we remove this last arc from the path $\mathbf{p}$, then we obtain a path
$\mathbf{p}^{\prime}$ from $r$ to $s$ (since $s$ is the source of $a$).

However, each path of $D$ is a path of $D^{\operatorname*{und}}$. Thus, in
particular, $\mathbf{p}$ is a path of $D^{\operatorname*{und}}$ from $r$ to
$t$, while $\mathbf{p}^{\prime}$ is a path of $D^{\operatorname*{und}}$ from
$r$ to $s$. Since $\mathbf{p}^{\prime}$ is exactly one edge shorter than
$\mathbf{p}$, we thus obtain $d\left(  r,s\right)  =d\left(  r,t\right)
-1<d\left(  r,t\right)  $. This proves Lemma \ref{lem.arbor-vs-tree.2}
\textbf{(a)}. \medskip

\textbf{(b)} The arc $a$ of the digraph $D$ has source $s$ and target $t$.
Hence, the edge $a$ of the tree $D^{\operatorname*{und}}$ has endpoints $s$
and $t$. Since $d\left(  r,s\right)  <d\left(  r,t\right)  $ (by part
\textbf{(a)}), this entails that its $r$-parent is $s$ and its $r$-child is
$t$ (by Definition \ref{def.tree.par-child-e}). Thus, in the digraph $\left(
D^{\operatorname*{und}}\right)  ^{r\rightarrow}$, this edge $a$ becomes an arc
with source $s$ and target $t$ (by Definition \ref{def.tree.rto}). This proves
Lemma \ref{lem.arbor-vs-tree.2} \textbf{(b)}.
\end{proof}

\begin{proof}
[Proof of Theorem \ref{thm.arbor.vs-tree.proper}.]If $\left(  V,A,\psi\right)
$ is a multidigraph, then we shall refer to the map $\psi:A\rightarrow V\times
V$ (which determines the source and the target of each arc) as the
\textquotedblleft psi-map\textquotedblright\ of this multidigraph.

Write the multidigraph $D$ as $D=\left(  V,A,\psi\right)  $. We shall now
prove the implications C1$\Longrightarrow$C2 and C2$\Longrightarrow$C1
separately: \medskip

\textit{Proof of the implication C1}$\Longrightarrow$\textit{C2:} Assume that
Statement C1 holds. That is, $D$ is an arborescence rooted from $r$. We must
prove Statement C2. In other words, we must prove that the undirected
multigraph $D^{\operatorname*{und}}$ is a tree, and that $D=\left(
D^{\operatorname*{und}}\right)  ^{r\rightarrow}$.

It is clear (by the definition of an arborescence) that
$D^{\operatorname*{und}}$ is a tree. It thus remains to prove that $D=\left(
D^{\operatorname*{und}}\right)  ^{r\rightarrow}$.

The multidigraphs $D$ and $\left(  D^{\operatorname*{und}}\right)
^{r\rightarrow}$ have the same set of vertices (namely, $V$) and the same set
of arcs (namely, $A$); we therefore just need to show that their psi-maps are
the same. In other words, we need to show that $\psi^{\prime}=\psi$, where
$\psi^{\prime}$ is the psi-map of $\left(  D^{\operatorname*{und}}\right)
^{r\rightarrow}$.

Let $a\in A$ be arbitrary. Let $\psi\left(  a\right)  =\left(  s,t\right)  $.
Thus, the arc $a$ of $D$ has source $s$ and target $t$. Lemma
\ref{lem.arbor-vs-tree.2} \textbf{(b)} therefore shows that in the
multidigraph $\left(  D^{\operatorname*{und}}\right)  ^{r\rightarrow}$, the
arc $a$ has source $s$ and target $t$ as well. In other words, $\psi^{\prime
}\left(  a\right)  =\left(  s,t\right)  $ (since $\psi^{\prime}$ is the
psi-map of this multidigraph). Hence, $\psi^{\prime}\left(  a\right)  =\left(
s,t\right)  =\psi\left(  a\right)  $.

Forget that we fixed $a$. We thus have shown that $\psi^{\prime}\left(
a\right)  =\psi\left(  a\right)  $ for each $a\in A$. In other words,
$\psi^{\prime}=\psi$. As explained above, this completes the proof of
Statement C2. \medskip

\textit{Proof of the implication C2}$\Longrightarrow$\textit{C1:} Assume that
Statement C2 holds. Thus, the undirected multigraph $D^{\operatorname*{und}}$
is a tree, and we have $D=\left(  D^{\operatorname*{und}}\right)
^{r\rightarrow}$. Hence, Lemma \ref{lem.arbor-vs-tree.1} (applied to
$T=D^{\operatorname*{und}}$) yields that the multidigraph $\left(
D^{\operatorname*{und}}\right)  ^{r\rightarrow}$ is an arborescence rooted
from $r$. In other words, $D$ is an arborescence rooted from $r$ (since
$D=\left(  D^{\operatorname*{und}}\right)  ^{r\rightarrow}$). This shows that
Statement C1 holds. \medskip

Having now proved both implications C1$\Longrightarrow$C2 and
C2$\Longrightarrow$C1, we conclude that Statements C1 and C2 are equivalent.
Thus, Theorem \ref{thm.arbor.vs-tree.proper} is proved.
\end{proof}

Oof.

Let's get one more consequence out of this. First, let us show that an
arborescence can have only one root:

\begin{proposition}
\label{prop.arbor.o1root}Let $D$ be an arborescence rooted from $r$. Then, $r$
is the \textbf{only} root of $D$.
\end{proposition}

\begin{proof}
[Proof of Proposition \ref{prop.arbor.o1root}.]Assume the contrary. Thus, $D$
has another root $s$ distinct from $r$. Hence, $D$ has a path from $r$ to $s$
(since $r$ is a root) as well as a path from $s$ to $r$ (since $s$ is a root).
Combining these paths gives a circuit of length $>0$. However, a circuit of
length $>0$ in a digraph must always contain a cycle (since Proposition
\ref{prop.mdg.cyc.btf-walk-cyc} shows that it either is a path or contains a
cycle; but it clearly cannot be a path). Hence, $D$ has a cycle. Therefore,
$D^{\operatorname*{und}}$ also has a cycle (since any cycle of $D$ is a cycle
of $D^{\operatorname*{und}}$). However, $D^{\operatorname*{und}}$ has no
cycles (since $D$ is an arborescence rooted from $r$). The preceding two
sentences contradict each other. This shows that the assumption was wrong, and
Proposition \ref{prop.arbor.o1root} is proven.
\end{proof}

\begin{definition}
\label{def.arbor.gen-arbor}A multidigraph $D$ is said to be an
\textbf{arborescence} if there exists a vertex $r$ of $D$ such that $D$ is an
arborescence rooted from $r$. In this case, this $r$ is uniquely determined as
the only root of $D$ (by Proposition \ref{prop.arbor.o1root}).
\end{definition}

\begin{theorem}
\label{thm.arbor-vs-tree.bij}There are two mutually inverse maps%
\begin{align*}
\left\{  \text{pairs }\left(  T,r\right)  \text{ of a tree }T\text{ and a
vertex }r\text{ of }T\right\}   &  \rightarrow\left\{  \text{arborescences}%
\right\}  ,\\
\left(  T,r\right)   &  \mapsto T^{r\rightarrow}%
\end{align*}
and%
\begin{align*}
\left\{  \text{arborescences}\right\}   &  \rightarrow\left\{  \text{pairs
}\left(  T,r\right)  \text{ of a tree }T\text{ and a vertex }r\text{ of
}T\right\}  ,\\
D  &  \mapsto\left(  D^{\operatorname*{und}},\sqrt{D}\right)  ,
\end{align*}
where $\sqrt{D}$ denotes the root of $D$.
\end{theorem}

\begin{proof}
The map%
\begin{align*}
\left\{  \text{pairs }\left(  T,r\right)  \text{ of a tree }T\text{ and a
vertex }r\text{ of }T\right\}   &  \rightarrow\left\{  \text{arborescences}%
\right\}  ,\\
\left(  T,r\right)   &  \mapsto T^{r\rightarrow}%
\end{align*}
is well-defined because of Lemma \ref{lem.arbor-vs-tree.1}. The map
\begin{align*}
\left\{  \text{arborescences}\right\}   &  \rightarrow\left\{  \text{pairs
}\left(  T,r\right)  \text{ of a tree }T\text{ and a vertex }r\text{ of
}T\right\}  ,\\
D  &  \mapsto\left(  D^{\operatorname*{und}},\sqrt{D}\right)  ,
\end{align*}
is well-defined because if $D$ is an arborescence, then
$D^{\operatorname*{und}}$ is a tree. In order to show that these two maps are
mutually inverse, we must check the following two statements:

\begin{enumerate}
\item Each arborescence $D$ satisfies $\left(  D^{\operatorname*{und}}\right)
^{r\rightarrow}=D$, where $r$ is the root of $D$;

\item Each pair $\left(  T,r\right)  $ of a tree $T$ and a vertex $r$ of $T$
satisfies $\left(  T^{r\rightarrow}\right)  ^{\operatorname*{und}}=T$ and
$\sqrt{\left(  T^{r\rightarrow}\right)  ^{\operatorname*{und}}}=r$.
\end{enumerate}

However, Statement 1 follows from Theorem \ref{thm.arbor.vs-tree.proper}
(specifically, from the implication C1$\Longrightarrow$C2 in Theorem
\ref{thm.arbor.vs-tree.proper}). Statement 2 follows from Lemma
\ref{lem.arbor-vs-tree.1} (more precisely, the $\left(  T^{r\rightarrow
}\right)  ^{\operatorname*{und}}=T$ part of Statement 2 is obvious, whereas
the $\sqrt{\left(  T^{r\rightarrow}\right)  ^{\operatorname*{und}}}=r$ part
follows from Lemma \ref{lem.arbor-vs-tree.1}). Thus, Theorem
\ref{thm.arbor-vs-tree.bij} is proved.
\end{proof}

Theorem \ref{thm.arbor-vs-tree.bij} formalizes the idea that an arborescence
is \textquotedblleft just a tree with a chosen vertex\textquotedblright. For
this reason, arborescences are sometimes called \textquotedblleft oriented
trees\textquotedblright, but this name is also shared with a more general
notion, which is why I avoid it.

\begin{exercise}
\label{exe.5.4}Let $G=\left(  V,E,\varphi\right)  $ be a connected multigraph
such that $\left\vert E\right\vert \geq\left\vert V\right\vert $. Show that
there exists an injective map $f:V\rightarrow E$ such that for each vertex
$v\in V$, the edge $f\left(  v\right)  $ contains $v$.

(In other words, show that we can assign to each vertex an edge that contains
this vertex in such a way that no edge is assigned twice.)
\end{exercise}

\subsection{Spanning arborescences}

In analogy to spanning subgraphs of a multigraph, we can define spanning
subdigraphs of a multidigraph:

\begin{definition}
A \textbf{spanning subdigraph} of a multidigraph $D=\left(  V,A,\psi\right)  $
means a multidigraph of the form $\left(  V,B,\psi\mid_{B}\right)  $, where
$B$ is a subset of $A$.

In other words, it means a submultidigraph of $D$ with the same vertex set as
$D$.

In other words, it means a multidigraph obtained from $D$ by removing some
arcs, but leaving all vertices untouched.
\end{definition}

\begin{definition}
Let $D$ be a multidigraph. Let $r$ be a vertex of $D$. A \textbf{spanning
arborescence of }$D$ \textbf{rooted from }$r$ means a spanning subdigraph of
$D$ that is an arborescence rooted from $r$.
\end{definition}

\begin{example}
Let $D=\left(  V,A,\psi\right)  $ be the following multidigraph:
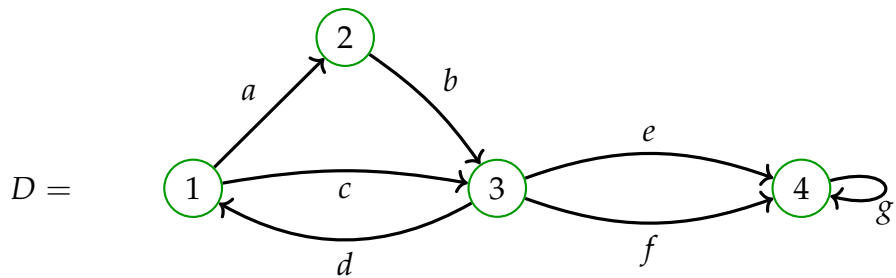
\[%
%
\ \ .
\]

Is there a spanning arborescence of $D$ rooted from $1$ ? Yes, for instance,%
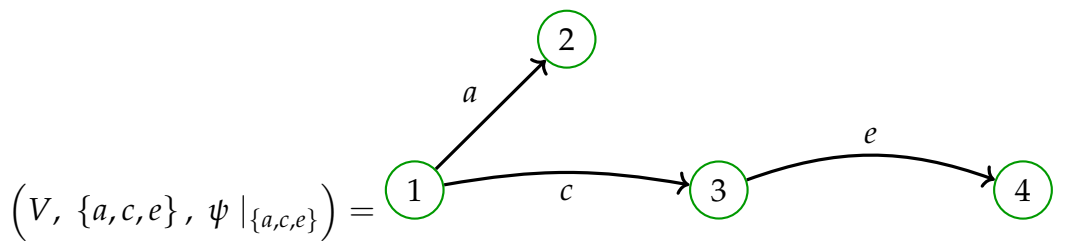
\[
\left(  V,\ \left\{  a,c,e\right\}  ,\ \psi\mid_{\left\{  a,c,e\right\}
}\right)  =%
%
\ \ .
\]
By abuse of notation, we shall refer to this spanning arborescence simply as
$\left\{  a,c,e\right\}  $ (since a spanning subdigraph of $D$ is uniquely
determined by its arc set). Another spanning arborescence of $D$ rooted from
$1$ is $\left\{  a,b,e\right\}  $. Yet another is $\left\{  a,b,f\right\}  $.
A non-example is $\left\{  a,d,f\right\}  $ (indeed, this is an arborescence
rooted from $3$, not from $1$).

Is there a spanning arborescence of $D$ rooted from $2$ ? Yes, for example
$\left\{  b,d,f\right\}  $.

Is there a spanning arborescence of $D$ rooted from $4$ ? No, since $4$ is not
a from-root of $D$.
\end{example}

This illustrates a first obstruction to the existence of spanning
arborescences: Namely, a digraph $D$ can have a spanning arborescence rooted
from $r$ only if $r$ is a from-root. This necessary criterion is also sufficient:

\begin{theorem}
\label{thm.spanning-arbor.exists}Let $D$ be a multidigraph. Let $r$ be a
from-root of $D$. Then, $D$ has a spanning arborescence rooted from $r$.
\end{theorem}

\begin{proof}
This is an analogue of the \textquotedblleft every connected multigraph has a
spanning tree\textquotedblright\ theorem (Theorem \ref{thm.spt.exists}) that
we proved in 4 ways. At least the first proof easily adapts to the directed case:

Remove arcs from $D$ one by one, but in such a way that the \textquotedblleft
rootness of $r$\textquotedblright\ (that is, the property that $r$ is a root
of our multidigraph) is preserved. So we can only remove an arc if $r$ remains
a root afterwards.

Clearly, this removing process will eventually come to an end, since $D$ has
only finitely many arcs. Let $D^{\prime}$ be the multidigraph obtained at the
end of this process. Then, $r$ is still a root of $D^{\prime}$, but we cannot
remove any more arcs from $D^{\prime}$ without breaking the rootness of $r$.
That is, if we remove any arc from $D^{\prime}$, then the vertex $r$ will no
longer be a from-root of the resulting multidigraph. This means that
$D^{\prime}$ satisfies Statement A5 from the arborescence equivalence theorem
(Theorem \ref{thm.arbor.eq-A}). Thus, $D^{\prime}$ satisfies Statement A1 as
well (since all six statements A1,\ A2, $\ldots$, A6 are equivalent). In other
words, $D^{\prime}$ is an arborescence rooted from $r$. Since $D^{\prime}$ is
a spanning subdigraph of $D$, we thus conclude that $D$ has a spanning
arborescence rooted from $r$ (namely, $D^{\prime}$). This proves Theorem
\ref{thm.spanning-arbor.exists}.
\end{proof}

\begin{question}
Can the other three proofs of Theorem \ref{thm.spt.exists} be adapted to
Theorem \ref{thm.spanning-arbor.exists}, too?
\end{question}

\begin{example}
\label{exa.spanning-arbor.Cn}Let $n$ be a positive integer. The $n$%
\textbf{-cycle digraph} $\overrightarrow{C}_{n}$ is defined to be the simple
digraph with vertices $1,2,\ldots,n$ and arcs $12,\ 23,\ 34,\ \ldots,\ \left(
n-1\right)  n,\ n1$. (Here is how it looks for $n=5$:
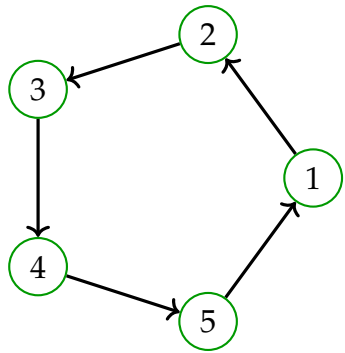
\[%
\begin{tikzpicture}
\begin{scope}[every node/.style={circle,thick,draw=green!60!black}]
\node(A) at (0:2) {$1$};
\node(B) at (360/5:2) {$2$};
\node(C) at (2*360/5:2) {$3$};
\node(D) at (3*360/5:2) {$4$};
\node(E) at (4*360/5:2) {$5$};
\end{scope}
\begin{scope}[every edge/.style={draw=black,very thick}]
\path[->] (A) edge (B) (B) edge (C) (C) edge (D) (D) edge (E) (E) edge (A);
\end{scope}
\end{tikzpicture}%
\]
)

Note that this digraph $\overrightarrow{C}_{n}$ is a directed analogue of the
cycle graph $C_{n}$. As we recall from Example \ref{exa.spt.Cn}, the cycle
graph $C_{n}$ has $n$ spanning trees.

In contrast, the digraph $\overrightarrow{C}_{n}$ has only one spanning
arborescence rooted from $1$. This spanning arborescence is the subdigraph of
$\overrightarrow{C}_{n}$ obtained by removing the arc $n1$.
\end{example}

\begin{proof}
If we remove the arc $n1$ from $\overrightarrow{C}_{n}$, then we obtain the
simple digraph $E$ with vertices $1,2,\ldots,n$ and arcs $12,\ 23,\ \ldots
,\ \left(  n-1\right)  n$. This digraph $E$ is easily seen to be an
arborescence rooted from $1$ (indeed, $1$ is a from-root of $E$, and the
underlying undirected graph $E^{\operatorname*{und}}=P_{n}$ has no cycles).
Thus, $E$ is a spanning arborescence of $\overrightarrow{C}_{n}$ rooted from
$1$.

We shall now prove that it is the only such arborescence. Indeed, let $F$ be
any spanning arborescence of $\overrightarrow{C}_{n}$ rooted from $1$. Then,
$1$ is a from-root of $F$. Hence, for each vertex $v\in\left\{  2,3,\ldots
,n\right\}  $, the digraph $F$ must have a path from $1$ to $v$, and thus must
contain an arc with target $v$ (namely, the last arc of this path). This arc
must be $\left(  v-1,v\right)  $ (since this is the only arc of
$\overrightarrow{C}_{n}$ with target $v$). Thus, for each vertex $v\in\left\{
2,3,\ldots,n\right\}  $, the digraph $F$ must contain the arc $\left(
v-1,v\right)  $. In other words, the digraph $F$ must contain all $n-1$ arcs
$12,\ 23,\ \ldots,\ \left(  n-1\right)  n$. If $F$ were to also contain the
remaining arc $n1$ of $\overrightarrow{C}_{n}$, then the underlying undirected
graph $F^{\operatorname*{und}}=C_{n}$ would contain a cycle, which would
contradict $F$ being an arborescence. Hence, $F$ cannot contain the arc $n1$.
Thus, $F$ contains the $n-1$ arcs $12,\ 23,\ \ldots,\ \left(  n-1\right)  n$
and no others. In other words, $F=E$. This shows that any spanning
arborescence of $\overrightarrow{C}_{n}$ rooted from $1$ must be $E$. In other
words, $E$ is the only spanning arborescence of $\overrightarrow{C}_{n}$
rooted from $1$. This completes the proof of Example
\ref{exa.spanning-arbor.Cn}.
\end{proof}

\subsection{The BEST theorem: statement}

We now come to something much more surprising.

Recall that a multidigraph $D=\left(  V,A,\psi\right)  $ is \textbf{balanced}
if and only if each vertex $v$ satisfies $\deg^{-}v=\deg^{+}v$. This is
necessary for the existence of an Eulerian circuit. If $D$ is weakly
connected, this is also sufficient (by Theorem \ref{thm.digr.euler-hier}
\textbf{(a)}).

Surprisingly, there is a formula for the number of these Eulerian circuits:

\begin{theorem}
[The BEST theorem]\label{thm.BEST.from}Let $D=\left(  V,A,\psi\right)  $ be a
balanced multidigraph such that each vertex has indegree $>0$. Fix an arc $a$
of $D$, and let $r$ be its target. Let $\tau\left(  D,r\right)  $ be the
number of spanning arborescences of $D$ rooted from $r$. Let $\varepsilon
\left(  D,a\right)  $ be the number of Eulerian circuits of $D$ whose last arc
is $a$. Then,%
\[
\varepsilon\left(  D,a\right)  =\tau\left(  D,r\right)  \cdot\prod_{u\in
V}\left(  \deg^{-}u-1\right)  !.
\]

\end{theorem}

The \textquotedblleft BEST\textquotedblright\ in the name of this theorem is
an abbreviation for de Bruijn, van Aardenne--Ehrenfest, Smith and Tutte, who
discovered it in the middle of the 20th century\footnote{More precisely, van
Aardenne--Ehrenfest and de Bruijn discovered it in 1951 (see \cite[\S 6]%
{VanEhr51}) generalizing an earlier result of Smith and Tutte.}. \footnote{We
note that the number of Eulerian circuits of $D$ whose last arc is $a$ is
precisely the number of all Eulerian circuits of $D$ counted up to rotation.
Indeed, each Eulerian circuit of $D$ contains the arc $a$ exactly once, and
thus can be rotated in a unique way to end with $a$.}

To prove this theorem, we shall restate it in terms of \textquotedblleft
arborescences to\textquotedblright\ (as opposed to \textquotedblleft
arborescences from\textquotedblright). Mathematically speaking, this
restatement isn't really necessary (the argument is the same in both cases up
to reversing the directions of all arcs), but it helps make the proof more
intuitive, since it lets us build our Eulerian circuits by moving forwards
rather than backwards.

\subsection{Arborescences rooted to $r$}

Here is the formal definition of \textquotedblleft arborescences
to\textquotedblright:

\begin{definition}
\label{def.arbor.arbor-to}Let $D$ be a multidigraph. Let $r$ be a vertex of
$D$.

\begin{enumerate}
\item[\textbf{(a)}] We say that $r$ is a \textbf{to-root} of $D$ if for each
vertex $v$ of $D$, the digraph $D$ has a path from $v$ to $r$.

\item[\textbf{(b)}] We say that $D$ is an \textbf{arborescence rooted to }$r$
if $r$ is a to-root of $D$ and the undirected multigraph
$D^{\operatorname*{und}}$ has no cycles.
\end{enumerate}
\end{definition}

Clearly, Definition \ref{def.arbor.arbor-from} and Definition
\ref{def.arbor.arbor-to} differ only in the direction of the arcs. In other
words, if we reverse each arc of our digraph (turning its source into its
target and vice versa), then a from-root becomes a to-root, and an
arborescence rooted from $r$ becomes an arborescence rooted to $r$, and vice
versa. Thus, every property that we have proved for arborescences rooted from
$r$ can be translated into the language of arborescences rooted to $r$ by
reversing all arcs.

\begin{fineprint}
If you want to see this stated more rigorously, here is a formal definition of
\textquotedblleft reversing each arc\textquotedblright:

\begin{definition}
Let $D=\left(  V,A,\psi\right)  $ be a multidigraph. Then,
$D^{\operatorname*{rev}}$ shall denote the multidigraph $\left(  V,A,\tau
\circ\psi\right)  $, where $\tau:V\times V\rightarrow V\times V$ is the map
that sends each pair $\left(  s,t\right)  $ to $\left(  t,s\right)  $. Thus,
if an arc $a$ of $D$ has source $s$ and target $t$, then it is also an arc of
$D^{\operatorname*{rev}}$, but in this digraph $D^{\operatorname*{rev}}$ it
has source $t$ and target $s$.

The multidigraph $D^{\operatorname*{rev}}$ is called the \textbf{reversal} of
the multidigraph $D$; we say that it is obtained from $D$ by \textquotedblleft
reversing each arc\textquotedblright.
\end{definition}

This notion of \textquotedblleft reversing each arc\textquotedblright\ allows
us to reverse walks in digraphs: If $\mathbf{w}$ is a walk from a vertex $s$
to $t$ in some multidigraph $D$, then its reversal $\operatorname{rev}%
\mathbf{w}$ (obtained by reading $\mathbf{w}$ backwards) is a walk from $t$ to
$s$ in the multidigraph $D^{\operatorname*{rev}}$. The same holds if we
replace the word \textquotedblleft walk\textquotedblright\ by
\textquotedblleft path\textquotedblright. Thus, we easily obtain the following:

\begin{proposition}
Let $D$ be a multidigraph. Let $r$ be a vertex of $D$. Then:

\begin{enumerate}
\item[\textbf{(a)}] The vertex $r$ is a to-root of $D$ if and only if $r$ is a
from-root of $D^{\operatorname*{rev}}$.

\item[\textbf{(b)}] The digraph $D$ is an arborescence rooted to $r$ if and
only if $D^{\operatorname*{rev}}$ is an arborescence rooted from $r$.
\end{enumerate}
\end{proposition}

\begin{proof}
Completely straightforward unpacking of the definitions.
\end{proof}
\end{fineprint}

Note that when we reverse each arc in a digraph $D$, the outdegrees of its
vertices become their indegrees and vice versa. Hence, a balanced digraph $D$
remains balanced when this happens. In particular, the BEST theorem (Theorem
\ref{thm.BEST.from}) thus gets translated as follows:

\begin{theorem}
[The BEST' theorem]\label{thm.BEST.to}Let $D=\left(  V,A,\psi\right)  $ be a
balanced multidigraph such that each vertex has outdegree $>0$. Fix an arc $a$
of $D$, and let $r$ be its source. Let $\tau\left(  D,r\right)  $ be the
number of spanning arborescences of $D$ rooted to $r$. Let $\varepsilon\left(
D,a\right)  $ be the number of Eulerian circuits of $D$ whose first arc is
$a$. Then,%
\[
\varepsilon\left(  D,a\right)  =\tau\left(  D,r\right)  \cdot\prod_{u\in
V}\left(  \deg^{+}u-1\right)  !.
\]

\end{theorem}

We will soon prove Theorem \ref{thm.BEST.to}, and then derive Theorem
\ref{thm.BEST.from} from it by reversing the arcs.

First, however, let us state the analogue of the Arborescence Equivalence
Theorem (Theorem \ref{thm.arbor.eq-A}) for \textquotedblleft arborescences
rooted to $r$\textquotedblright\ (as opposed to \textquotedblleft
arborescences rooted from $r$\textquotedblright):

\begin{theorem}
[The dual arborescence equivalence theorem]\label{thm.arbor.eq-A-dual}Let
$D=\left(  V,A,\psi\right)  $ be a multidigraph with a to-root $r$. Then, the
following six statements are equivalent:

\begin{itemize}
\item \textbf{Statement A'1:} The multidigraph $D$ is an arborescence rooted
to $r$.

\item \textbf{Statement A'2:} We have $\left\vert A\right\vert =\left\vert
V\right\vert -1$.

\item \textbf{Statement A'3:} The multigraph $D^{\operatorname*{und}}$ is a tree.

\item \textbf{Statement A'4:} For each vertex $v\in V$, the multidigraph $D$
has a unique walk from $v$ to $r$.

\item \textbf{Statement A'5:} If we remove any arc from $D$, then the vertex
$r$ will no longer be a to-root of the resulting multidigraph.

\item \textbf{Statement A'6:} We have $\deg^{+}r=0$, and each $v\in
V\setminus\left\{  r\right\}  $ satisfies $\deg^{+}v=1$.
\end{itemize}
\end{theorem}

\begin{proof}
Upon reversing all arcs of $D$, this turns into the original Arborescence
Equivalence Theorem (Theorem \ref{thm.arbor.eq-A}).
\end{proof}

\subsection{The BEST theorem: proof}

We now come to the proof of the BEST theorem (Theorem \ref{thm.BEST.from}). As
we said, we proceed by proving Theorem \ref{thm.BEST.to} first. We first
outline the idea of the proof; then we will give the details.

\begin{proof}
[Proof idea for Theorem \ref{thm.BEST.to}.]An $a$\textbf{-Eulerian circuit}
shall mean an Eulerian circuit of $D$ whose first arc is $a$.

Let $\mathbf{e}$ be an $a$-Eulerian circuit. Its first arc is $a$; therefore,
its first and last vertex is $r$.

Being an Eulerian circuit, $\mathbf{e}$ must contain each arc of $D$ and
therefore contain each vertex of $D$ (since each vertex has outdegree $>0$).
For each vertex $u\neq r$, we let $e\left(  u\right)  $ be the \textbf{last
exit} of $\mathbf{e}$ from $u$, that is, the last arc of $\mathbf{e}$ that has
source $u$. Let $\operatorname*{Exit}\mathbf{e}$ be the set of these last
exits $e\left(  u\right)  $ for all vertices $u\neq r$. Then, we claim:

\begin{statement}
\textit{Claim 1:} This set $\operatorname*{Exit}\mathbf{e}$ (or, more
precisely, the spanning subdigraph $\left(  V,\ \operatorname*{Exit}%
\mathbf{e},\ \psi\mid_{\operatorname*{Exit}\mathbf{e}}\right)  $) is a
spanning arborescence of $D$ rooted to $r$.
\end{statement}

Let's assume for the moment that Claim 1 is proven. Thus, given any
$a$-Eulerian circuit $\mathbf{e}$, we have constructed a spanning arborescence
of $D$ rooted to $r$.

How many $a$-Eulerian circuits $\mathbf{e}$ lead to a given arborescence in
this way? The answer is rather nice:

\begin{statement}
\textit{Claim 2:} For each spanning arborescence $\left(  V,B,\psi\mid
_{B}\right)  $ of $D$ rooted to $r$, there are exactly $\prod_{u\in V}\left(
\deg^{+}u-1\right)  !$ many $a$-Eulerian circuits $\mathbf{e}$ such that
$\operatorname*{Exit}\mathbf{e}=B$.
\end{statement}

Let us again assume that this is proven. Combining Claim 1 with Claim 2, we
obtain a $\prod_{u\in V}\left(  \deg^{+}u-1\right)  !$-to-$1$ correspondence
between the $a$-Eulerian circuits and the spanning arborescences of $D$ rooted
to $r$. Thus, the number of the former is $\prod_{u\in V}\left(  \deg
^{+}u-1\right)  !$ times the number of the latter. But this is precisely the
claim of Theorem \ref{thm.BEST.to}. Hence, in order to prove Theorem
\ref{thm.BEST.to}, it remains to prove Claim 1 and Claim 2.
\end{proof}

Here is the complete proof:

\begin{proof}
[Proof of Theorem \ref{thm.BEST.to}.]Some notations first:

An \textbf{outgoing arc} from a vertex $u$ will mean an arc whose source is
$u$. An \textbf{incoming arc} into a vertex $u$ will mean an arc whose target
is $u$.

An $a$\textbf{-Eulerian circuit} shall mean an Eulerian circuit of $D$ whose
first arc is $a$.

A \textbf{sparb} shall mean a spanning arborescence of $D$ rooted to $r$.

A spanning subdigraph of $D$ always has the form $\left(  V,B,\psi\mid
_{B}\right)  $ for some subset $B$ of $A$. Thus, it is uniquely determined by
its arc set $B$.

Hence, from now on, we shall identify a spanning subdigraph $\left(
V,B,\psi\mid_{B}\right)  $ of $D$ with its arc set $B$. Conversely, any subset
$B$ of $A$ will be identified with the corresponding spanning subdigraph
$\left(  V,B,\psi\mid_{B}\right)  $ of $D$. Thus, for instance, when we say
that a subset $B$ of $A$ \textquotedblleft is a sparb\textquotedblright, we
shall actually mean that the corresponding spanning subdigraph $\left(
V,B,\psi\mid_{B}\right)  $ is a sparb.

For each $a$-Eulerian circuit $\mathbf{e}$, we define a subset
$\operatorname*{Exit}\mathbf{e}$ of $A$ as follows:

Let $\mathbf{e}$ be an $a$-Eulerian circuit. Its first arc is $a$; thus, its
first and last vertex is $r$. Being an Eulerian circuit, $\mathbf{e}$ must
contain each arc of $D$ and therefore also contain each vertex of $D$ (since
each vertex of $D$ has outdegree $>0$). For each vertex $u\in V\setminus
\left\{  r\right\}  $, we let $e\left(  u\right)  $ be the \textbf{last exit}
of $\mathbf{e}$ from $u$; this means the last arc of $\mathbf{e}$ that has
source $u$. We let $\operatorname*{Exit}\mathbf{e}$ be the set of these last
exits $e\left(  u\right)  $ for all $u\in V\setminus\left\{  r\right\}  $.
Thus, we have defined a subset $\operatorname*{Exit}\mathbf{e}$ of $A$ for
each $a$-Eulerian circuit $\mathbf{e}$.

\begin{example}
\label{exa.pf.thm.BEST.to.1}Here is an example of this construction: Let $D$
be the multidigraph%
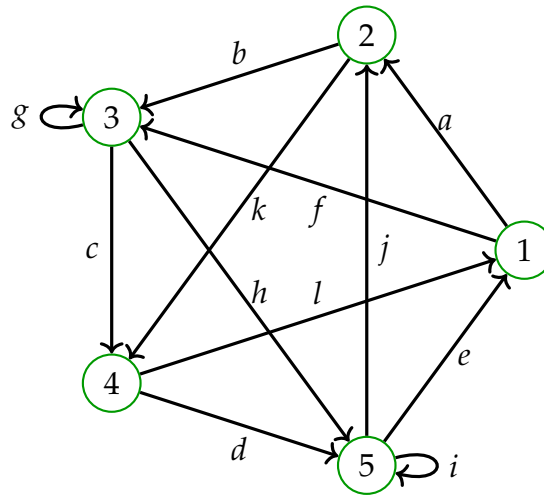
\[%
\begin{tikzpicture}[scale=1.5]
\begin{scope}[every node/.style={circle,thick,draw=green!60!black}]
\node(1) at (0:2) {$1$};
\node(2) at (360/5:2) {$2$};
\node(3) at (2*360/5:2) {$3$};
\node(4) at (3*360/5:2) {$4$};
\node(5) at (4*360/5:2) {$5$};
\end{scope}
\begin{scope}[every edge/.style={draw=black,very thick}, every loop/.style={}]
\path[->] (1) edge node[above] {$a$} (2) (2) edge node[above] {$b$} (3);
\path[->] (3) edge node[left] {$c$} (4) (4) edge node[below] {$d$} (5);
\path[->] (5) edge node[right] {$e$} (1) (1) edge node[below] {$f$} (3);
\path[->] (3) edge[loop left] node[left]{$g$} (3) (3) edge node[right] {$h$}
(5);
\path[->] (5) edge[loop right] node[right]{$i$} (5);
\path[->] (5) edge node[right] {$j$} (2) (2) edge node[right] {$k$} (4);
\path[->] (4) edge node[above] {$l$} (1);
\end{scope}
\end{tikzpicture}%
\]
with $r=1$, and let $\mathbf{e}$ be the $a$-Eulerian circuit%
\[
\left(  1,a,2,b,3,c,4,d,5,e,1,f,3,g,3,h,5,i,5,j,2,k,4,l,1\right)
\]
(we have deliberately named the arcs in such a way that they appear on an
Eulerian circuit in alphabetic order). Then,%
\[
e\left(  2\right)  =k,\ \ \ \ \ \ \ \ \ \ e\left(  3\right)
=h,\ \ \ \ \ \ \ \ \ \ e\left(  4\right)  =l,\ \ \ \ \ \ \ \ \ \ e\left(
5\right)  =j,
\]
so that $\operatorname*{Exit}\mathbf{e}=\left\{  k,h,l,j\right\}  $. Here is
$\operatorname*{Exit}\mathbf{e}$ as a spanning subdigraph:%
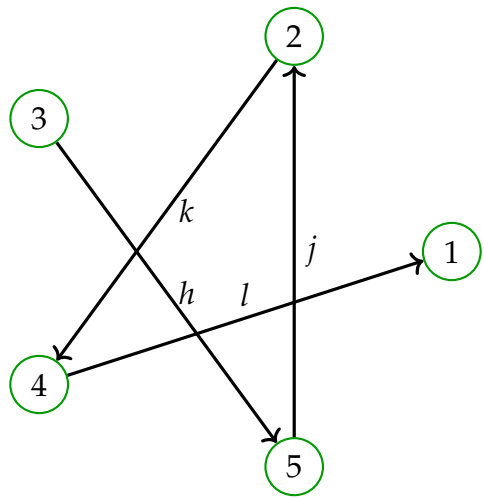
\[%
\begin{tikzpicture}[scale=1.5]
\begin{scope}[every node/.style={circle,thick,draw=green!60!black}]
\node(1) at (0:2) {$1$};
\node(2) at (360/5:2) {$2$};
\node(3) at (2*360/5:2) {$3$};
\node(4) at (3*360/5:2) {$4$};
\node(5) at (4*360/5:2) {$5$};
\end{scope}
\begin{scope}[every edge/.style={draw=black,very thick}, every loop/.style={}]
\path[->] (3) edge node[right] {$h$} (5);
\path[->] (5) edge node[right] {$j$} (2) (2) edge node[right] {$k$} (4);
\path[->] (4) edge node[above] {$l$} (1);
\end{scope}
\end{tikzpicture}%
\]

\end{example}

Now, we claim the following:

\begin{statement}
\textit{Claim 1:} Let $\mathbf{e}$ be an $a$-Eulerian circuit. Then, the set
$\operatorname*{Exit}\mathbf{e}$ is a sparb.
\end{statement}

\begin{statement}
\textit{Claim 2:} For each sparb $B$ (regarded as a subset of $A$), there are
exactly $\prod_{u\in V}\left(  \deg^{+}u-1\right)  !$ many $a$-Eulerian
circuits $\mathbf{e}$ such that $\operatorname*{Exit}\mathbf{e}=B$.
\end{statement}

[\textit{Proof of Claim 1:} The set $\operatorname*{Exit}\mathbf{e}$ contains
exactly one outgoing arc (namely, $e\left(  u\right)  $) from each vertex
$u\in V\setminus\left\{  r\right\}  $, and no outgoing arc from $r$. Thus,
$\left\vert \operatorname*{Exit}\mathbf{e}\right\vert =\left\vert V\right\vert
-1$.

Let us number the arcs of $\mathbf{e}$ as $a_{1},a_{2},\ldots,a_{m}$, in the
order in which they appear in $\mathbf{e}$. (Thus, $a_{1}=a$, since the first
arc of $\mathbf{e}$ is $a$.)

Recall that the arcs in $\operatorname*{Exit}\mathbf{e}$ are the arcs
$e\left(  u\right)  $ for all $u\in V\setminus\left\{  r\right\}  $ (defined
as above -- i.e., the arc $e\left(  u\right)  $ is the last exit of
$\mathbf{e}$ from $u$). We shall refer to these arcs as the \textbf{last-exit
arcs}.

For each $u\in V\setminus\left\{  r\right\}  $, we let $j\left(  u\right)  $
be the unique number $i\in\left\{  1,2,\ldots,m\right\}  $ such that $e\left(
u\right)  =a_{i}$. (This $i$ indeed exists and is unique, since each arc of
$D$ appears exactly once on $\mathbf{e}$.) Thus, $j\left(  u\right)  $ tells
us how late in the Eulerian circuit $\mathbf{e}$ the arc $e\left(  u\right)  $
appears. Since $e\left(  u\right)  $ is the last exit of $\mathbf{e}$ from
$u$, the Eulerian circuit $\mathbf{e}$ never visits the vertex $u$ again after this.

Thus, if a last-exit arc $e\left(  u\right)  $ has target $v\neq r$, then%
\begin{equation}
j\left(  u\right)  <j\left(  v\right)  \label{pf.thm.BEST.to.c1.pf.ineq}%
\end{equation}
(because the arc $e\left(  u\right)  $ leads the circuit $\mathbf{e}$ into the
vertex $v$, which the circuit then has to exit at least once; therefore, the
corresponding last-exit arc $e\left(  v\right)  $ has to appear later in
$\mathbf{e}$ than the arc $e\left(  u\right)  $).

We shall now show that $r$ is a to-root of $\operatorname*{Exit}\mathbf{e}$
(that is, of the spanning subdigraph $\left(  V,\operatorname*{Exit}%
\mathbf{e},\psi\mid_{\operatorname*{Exit}\mathbf{e}}\right)  $). To this
purpose, we must show that for each vertex $v\in V$, there is a path from $v$
to $r$ in the digraph $\left(  V,\operatorname*{Exit}\mathbf{e},\psi
\mid_{\operatorname*{Exit}\mathbf{e}}\right)  $.

Indeed, let $v\in V$ be any vertex. We must find a path from $v$ to $r$ in the
digraph $\left(  V,\operatorname*{Exit}\mathbf{e},\psi\mid
_{\operatorname*{Exit}\mathbf{e}}\right)  $. It will suffice to find a walk
from $v$ to $r$ in this digraph (by Corollary \ref{cor.mdg.walk-thus-path}).
In other words, we must find a way to walk from $v$ to $r$ in $D$ using
last-exit arcs only.

So we start walking at $v$. If $v=r$, then we are already done. Otherwise, we
have $v\in V\setminus\left\{  r\right\}  $, so that the arc $e\left(
v\right)  $ and the number $j\left(  v\right)  $ are well-defined. We thus
take the arc $e\left(  v\right)  $. This brings us to a vertex $v^{\prime}$
(namely, the target of $e\left(  v\right)  $) that satisfies $j\left(
v\right)  <j\left(  v^{\prime}\right)  $ (by (\ref{pf.thm.BEST.to.c1.pf.ineq}%
)). If this vertex $v^{\prime}$ is $r$, then we are done. If not, then
$e\left(  v^{\prime}\right)  $ and $j\left(  v^{\prime}\right)  $ are
well-defined, so we continue our walk by taking the arc $e\left(  v^{\prime
}\right)  $. This brings us to a further vertex $v^{\prime\prime}$ (namely,
the target of $e\left(  v^{\prime}\right)  $) that satisfies $j\left(
v^{\prime}\right)  <j\left(  v^{\prime\prime}\right)  $ (by
(\ref{pf.thm.BEST.to.c1.pf.ineq})). If this vertex $v^{\prime\prime}$ is $r$,
then we are done. Otherwise, we proceed as before. We thus construct a walk
\[
\left(  v,e\left(  v\right)  ,v^{\prime},e\left(  v^{\prime}\right)
,v^{\prime\prime},e\left(  v^{\prime\prime}\right)  ,\ldots\right)
\]
that either goes on indefinitely or stops at the vertex $r$.

However, this walking process cannot go on forever (since the chain of
inequalities $j\left(  v\right)  <j\left(  v^{\prime}\right)  <j\left(
v^{\prime\prime}\right)  <\cdots$ would force the numbers $j\left(  v\right)
,j\left(  v^{\prime}\right)  ,j\left(  v^{\prime\prime}\right)  ,\ldots$ to be
all distinct, but there are only $m$ distinct numbers in $\left\{
1,2,\ldots,m\right\}  $). Thus, it must stop at the vertex $r$. So we have
found a walk from $v$ to $r$ using last-exit arcs only. Thus,
$\operatorname*{Exit}\mathbf{e}$ has a walk from $v$ to $r$. Hence,
$\operatorname*{Exit}\mathbf{e}$ has a path from $v$ to $r$.

Forget that we fixed $v$. We thus have shown that for each vertex $v\in V$,
there is a path from $v$ to $r$ in the digraph $\left(  V,\operatorname*{Exit}%
\mathbf{e},\psi\mid_{\operatorname*{Exit}\mathbf{e}}\right)  $. In other
words, $r$ is a to-root of $\operatorname*{Exit}\mathbf{e}$. Hence, we
conclude (using the implication A'2$\Longrightarrow$A'1 in Theorem
\ref{thm.arbor.eq-A-dual}) that $\operatorname*{Exit}\mathbf{e}$ is an
arborescence rooted to $r$ (since $\left\vert \operatorname*{Exit}%
\mathbf{e}\right\vert =\left\vert V\right\vert -1$). Therefore,
$\operatorname*{Exit}\mathbf{e}$ is a sparb. This proves Claim 1.] \medskip

[\textit{Proof of Claim 2:} Let $B$ be a sparb. (As before, $B$ is a set of
arcs, and we identify it with the spanning subdigraph $\left(  V,B,\psi
\mid_{B}\right)  $.)

We must prove that there are exactly $\prod\limits_{u\in V}\left(  \deg
^{+}u-1\right)  !$ many $a$-Eulerian circuits $\mathbf{e}$ such that
$\operatorname*{Exit}\mathbf{e}=B$.

We shall refer to the arcs in $B$ as the $B$\textbf{-arcs}. Recall that $B$ is
an arborescence rooted to $r$ (since $B$ is a sparb). Hence, by the
implication A'1$\Longrightarrow$A'6 in Theorem \ref{thm.arbor.eq-A-dual}, we
see that the outdegrees of its vertices satisfy%
\[
\deg_{B}^{+}r=0,\ \ \ \ \ \ \ \ \ \ \text{and}\ \ \ \ \ \ \ \ \ \ \deg_{B}%
^{+}v=1\text{ for all }v\in V\setminus\left\{  r\right\}
\]
(where $\deg_{B}^{+}v$ means the outdegree of a vertex in the digraph $\left(
V,B,\psi\mid_{B}\right)  $). In other words, there is no $B$-arc with source
$r$; however, for each vertex $u\in V\setminus\left\{  r\right\}  $, there is
exactly one $B$-arc with source $u$.

Now, we are trying to count the $a$-Eulerian circuits $\mathbf{e}$ such that
$\operatorname*{Exit}\mathbf{e}=B$.

Let us try to construct such an $a$-Eulerian circuit $\mathbf{e}$ as follows:

A turtle wants to walk through the digraph $D$ using each arc of $D$ at most
once. It starts its walk by heading out from the vertex $r$ along the arc $a$.
From that point on, it proceeds in the usual way you would walk on a digraph:
Each time it reaches a vertex, it chooses an arbitrary arc leading out of this
vertex, observing the following two rules:

\begin{enumerate}
\item It never uses an arc that it has already used before.

\item It never uses a $B$-arc unless it has to (i.e., unless this $B$-arc is
the only outgoing arc from its current position that is still unused).
\end{enumerate}

Clearly, the turtle will eventually get stuck at some vertex (with no more
arcs left to continue walking along), since $D$ has only finitely many arcs.

Let $\mathbf{w}$ be the total walk that the turtle has traced by the time it
got stuck. Thus, $\mathbf{w}$ is a trail (i.e., a walk that uses no arc more
than once) that starts with the vertex $r$ and the arc $a$.

We will soon see that $\mathbf{w}$ is an $a$-Eulerian circuit satisfying
$\operatorname*{Exit}\mathbf{w}=B$. First, however, let us see an example:

\begin{example}
\label{exa.pf.thm.BEST.to.turt}Let $D$ be the multidigraph%
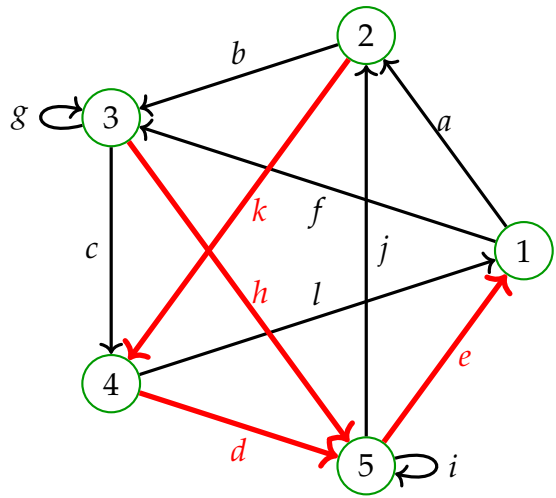
\[%
%
\ \ ,
\]
and let $r=1$ and $a=a$ (we called it $a$ on purpose). Let $B$ be the set
$\left\{  d,e,h,k\right\}  $, regarded as a spanning subdigraph of $D$. (The
arcs of $B$ are drawn bold and in red in the above picture.)

The turtle starts at $r=1$ and walks along the arc $a$. This leads it to the
vertex $2$. It now must choose between the arcs $b$ and $k$, but since it must
not use the $B$-arc $k$ unless it has to, it is actually forced to take the
arc $b$ next. This brings it to the vertex $3$. It now has to choose between
the arcs $c$, $g$ and $h$, but again the arc $h$ is disallowed because it is
not yet time to use a $B$-arc. Let us say that it takes the arc $g$. This
brings it back to the vertex $3$. Next, the turtle must walk along $c$ (since
$g$ is already used, while the $B$-arc still must wait until it is the only
option). This brings it to the vertex $4$. Its next step is to take the arc
$l$ to the vertex $1$. From there, it follows the arc $f$ to the vertex $3$.
Now, it can finally take the $B$-arc $h$, since all the other outgoing arcs
from $3$ have already been used. This brings it to the vertex $5$. Now it has
a choice between the arcs $e$, $i$ and $j$, but the arc $e$ is disallowed
because it is a $B$-arc. Let us say it decides to use the arc $j$. This brings
it to the vertex $2$. From there, it takes the $B$-arc $k$ to the vertex $4$
(since it has no other options). From there, it continues along the $B$-arc
$d$ to the vertex $5$. Now, it has to traverse the loop $i$, and then leave
$5$ along the $B$-arc $e$ to come back to $1$. At this point, the turtle is
stuck, since it has nowhere left to go. The walk $\mathbf{w}$ we obtained is
thus%
\[
\mathbf{w}=\left(  1,a,2,b,3,g,3,c,4,l,1,f,3,h,5,j,2,k,4,d,5,i,5,e,1\right)
.
\]
(Of course, other choices would have led to other walks.)
\end{example}

Returning to the general case, let us analyze the walk $\mathbf{w}$ traversed
by the turtle.

\begin{itemize}
\item First, we claim that $\mathbf{w}$ is a closed walk (i.e., ends at $r$).

[\textit{Proof:} Assume the contrary. Let $u$ be the ending point of
$\mathbf{w}$. Thus, $u$ is the vertex at which the turtle gets stuck.
Moreover, $u\neq r$ (since we just assumed that $\mathbf{w}$ is not a closed
walk). Hence, the walk $\mathbf{w}$ enters the vertex $u$ more often than it
leaves it (since it ends but does not start at $u$). In other words, the
turtle has entered the vertex $u$ more often than it has left it. However,
since $D$ is balanced, we have $\deg^{-}u=\deg^{+}u$. The turtle has entered
the vertex $u$ at most $\deg^{-}u$ times (because it cannot use an arc twice,
but there are only $\deg^{-}u$ many arcs with target $u$). Thus, it has left
the vertex $u$ \textbf{less} than $\deg^{-}u$ times (because it has entered
the vertex $u$ more often than it has left it). Since $\deg^{-}u=\deg^{+}u$,
this means that the turtle has left the vertex $u$ less than $\deg^{+}u$
times. Thus, by the time the turtle has gotten stuck at $u$, there is at least
one outgoing arc from $u$ that has not been used by the turtle. Therefore, the
turtle is not actually stuck at $u$. This is a contradiction. Thus, our
assumption was wrong, so we have proved that $\mathbf{w}$ is a closed walk.]
\end{itemize}

In other words, $\mathbf{w}$ is a circuit. We shall next show that
$\mathbf{w}$ is an Eulerian circuit.

To do so, we introduce one more piece of notation: A vertex $u$ of $D$ will be
called \textbf{exhausted} if the turtle has used each outgoing arc from $u$
(that is, if each outgoing arc from $u$ is used in the circuit $\mathbf{w}$).

Since $\mathbf{w}$ is a circuit, the ending point of $\mathbf{w}$ is its
starting point, i.e., the vertex $r$. Thus, the turtle must have gotten stuck
at $r$. Hence, the vertex $r$ is exhausted.

\begin{itemize}
\item We shall now show that \textbf{all} vertices of $D$ are exhausted.

[\textit{Proof:} Assume the contrary. Thus, there exists a vertex $u$ of $D$
that is not exhausted. Consider this $u$. But $B$ is a sparb, thus an
arborescence rooted to $r$. Hence, $r$ is a to-root of $B$. Therefore, there
exists a path $\mathbf{p}=\left(  p_{0},b_{1},p_{1},b_{2},p_{2},\ldots
,b_{k},p_{k}\right)  $ from $u$ to $r$ in $B$. Consider this path. Thus, we
have $p_{0}=u$ and $p_{k}=r$, and all the arcs $b_{1},b_{2},\ldots,b_{k}$
belong to $B$.

There exists at least one $i\in\left\{  0,1,\ldots,k\right\}  $ such that the
vertex $p_{i}$ is exhausted (for instance, $i=k$ qualifies, since $p_{k}=r$ is
exhausted). Consider the \textbf{smallest} such $i$. Then, $p_{i}\neq p_{0}$
(since $p_{i}$ is exhausted, but $p_{0}=u$ is not). Hence, $i\neq0$, so that
$i\geq1$. Therefore, $p_{i-1}$ exists. Moreover, the vertex $p_{i-1}$ is not
exhausted (since $i$ was defined to be the \textbf{smallest} element of
$\left\{  0,1,\ldots,k\right\}  $ such that $p_{i}$ is exhausted).

The arc $b_{i}$ has source $p_{i-1}$ and target $p_{i}$. Thus, it is an
outgoing arc from $p_{i-1}$ and incoming arc into $p_{i}$. Furthermore, it
belongs to $B$ (since all the arcs $b_{1},b_{2},\ldots,b_{k}$ belong to $B$).

The digraph $D$ is balanced; thus, $\deg^{+}\left(  p_{i}\right)  =\deg
^{-}\left(  p_{i}\right)  $.

The vertex $p_{i}$ is exhausted. In other words, the turtle has used each
outgoing arc from $p_{i}$ (by the definition of \textquotedblleft
exhausted\textquotedblright). Since the turtle never reuses an arc, this
entails that the turtle has used exactly $\deg^{+}\left(  p_{i}\right)  $ many
outgoing arcs from $p_{i}$ (since $\deg^{+}\left(  p_{i}\right)  $ is the
total number of outgoing arcs from $p_{i}$ in $D$). In other words, it has
used exactly $\deg^{-}\left(  p_{i}\right)  $ many outgoing arcs from $p_{i}$
(since $\deg^{+}\left(  p_{i}\right)  =\deg^{-}\left(  p_{i}\right)  $).

However, the turtle's trajectory is a closed walk (in fact, it is the walk
$\mathbf{w}$, which is closed). Thus, it must enter the vertex $p_{i}$ as
often as it leaves this vertex. In other words, the number of incoming arcs
into $p_{i}$ used by the turtle must equal the number of outgoing arcs from
$p_{i}$ used by the turtle. Since we just found (in the preceding paragraph)
that the latter number is $\deg^{-}\left(  p_{i}\right)  $, we thus conclude
that the former number is $\deg^{-}\left(  p_{i}\right)  $ as well. In other
words, the turtle must have used exactly $\deg^{-}\left(  p_{i}\right)  $ many
incoming arcs into $p_{i}$. Since $\deg^{-}\left(  p_{i}\right)  $ is the
total number of incoming arcs into $p_{i}$ in $D$, we thus conclude that the
turtle must have used all incoming arcs into $p_{i}$ (since the turtle never
reuses an arc).

Hence, in particular, the turtle must have used the arc $b_{i}$ (since $b_{i}$
is an incoming arc into $p_{i}$). This arc $b_{i}$ is an outgoing arc from
$p_{i-1}$. But $b_{i}$ is a $B$-arc, and thus our turtle uses this arc only as
a last resort (i.e., after using all other outgoing arcs from $p_{i-1}$).
Hence, we conclude that the turtle must have used all outgoing arcs from
$p_{i-1}$ (since it has used $b_{i}$). In other words, $p_{i-1}$ is exhausted.
But this contradicts the fact that $p_{i-1}$ is not exhausted! This shows that
our assumption was wrong, and our proof is finished.\footnote{For the sake of
diversity, let me sketch a \textit{second proof} of the same claim (i.e., that
all vertices in $D$ are exhausted)\textit{:}
\par
Assume the contrary. Thus, there exists a non-exhausted vertex $u$ of $D$.
Consider this $u$. Then, $u\neq r$ (since $r$ is exhausted but $u$ is not).
Since $u$ is not exhausted, there is at least one outgoing arc from $u$ that
the turtle has not used. Hence, the turtle has not used the $B$-arc outgoing
from $u$ (since the turtle never uses a $B$-arc before it has to). Let $f$ be
this $B$-arc, and let $u^{\prime}$ be its target. Thus, the turtle has not
used all incoming arcs of $u^{\prime}$ (because it has not used the arc $f$).
As a consequence, it has not used all outgoing arcs from $u^{\prime}$ either
(because the turtle has left $u^{\prime}$ as often as it has entered
$u^{\prime}$, but the balancedness of $D$ entails that $\deg^{-}\left(
u^{\prime}\right)  =\deg^{+}\left(  u^{\prime}\right)  $). In other words, the
vertex $u^{\prime}$ is non-exhausted.
\par
Thus, by starting at the non-exhausted vertex $u$ and taking the $B$-arc
outgoing from $u$, we have arrived at a further non-exhausted vertex
$u^{\prime}$. Applying the same argument to $u^{\prime}$ instead of $u$, we
can take a further $B$-arc and arrive at a further non-exhausted vertex
$u^{\prime\prime}$. Continuing like this, we obtain an infinite sequence
$\left(  u,u^{\prime},u^{\prime\prime},\ldots\right)  $ of non-exhausted
vertices such that any vertex in this sequence is reached from the previous
one by traveling along a $B$-arc. Clearly, this sequence must have two equal
vertices (since $D$ has only finitely many vertices). For example, let's say
that $u^{\prime\prime}=u^{\prime\prime\prime\prime\prime}$. Then, if we
consider only the part of the sequence between $u^{\prime\prime}$ and
$u^{\prime\prime\prime\prime\prime}$, then we obtain a closed walk%
\[
\left(  u^{\prime\prime},\ast,u^{\prime\prime\prime},\ast,u^{\prime
\prime\prime\prime},\ast,u^{\prime\prime\prime\prime\prime}\right)  ,
\]
where each asterisk stands for some $B$-arc (not the same one, of course).
This is a closed walk of the digraph $\left(  V,B,\psi\mid_{B}\right)  $.
Since this closed walk has length $>0$, it cannot be a path; therefore, it
contains a cycle (by Proposition \ref{prop.mdg.cyc.btf-walk-cyc}). Thus, we
have found a cycle of the digraph $\left(  V,B,\psi\mid_{B}\right)  $.
However, the digraph $\left(  V,B,\psi\mid_{B}\right)  $ is an arborescence,
and thus has no cycles (because if $D$ is an arborescence, then any cycle of
$D$ would be a cycle of $D^{\operatorname*{und}}$; but the multigraph
$D^{\operatorname*{und}}$ has no cycles by the definition of an arborescence).
The previous two sentences contradict each other. This shows that our
assumption was wrong, and our proof is finished.}.]
\end{itemize}

Thus, we have shown that all vertices of $D$ are exhausted. In other words,
the turtle has used all arcs of $D$. In other words, the trail $\mathbf{w}$
contains all arcs of $D$. Since $\mathbf{w}$ is a trail and a closed walk,
this entails that $\mathbf{w}$ is an Eulerian circuit of $D$. Since
$\mathbf{w}$ starts with $r$ and $a$, this shows further that $\mathbf{w}$ is
an $a$-Eulerian circuit. Since the turtle only used $B$-arcs as a last resort
(and it used each $B$-arc eventually, because $\mathbf{w}$ is Eulerian), we
have $\operatorname*{Exit}\mathbf{w}=B$.

Thus, the turtle's walk has produced an $a$-Eulerian circuit $\mathbf{e}$
satisfying $\operatorname*{Exit}\mathbf{e}=B$ (namely, the walk $\mathbf{w}$).
However, this circuit depends on some decisions the turtle made during its
walk. Namely, every time the turtle was at some vertex $u\in V$, it had to
decide which arc to take next; this arc had to be an unused arc with source
$u$, subject to the conditions that

\begin{enumerate}
\item if $u\neq r$, then the $B$-arc\footnote{We say \textquotedblleft%
\textbf{the} $B$-arc\textquotedblright, because there is exactly one $B$-arc
with source $u$.} has to be used last;

\item if $u=r$, then the arc $a$ has to be used first.
\end{enumerate}

Let us count how many options the turtle has had in total. To make the
argument clearer, we modify the procedure somewhat: Instead of deciding ad-hoc
which arc to take, the turtle should now make all these decisions before
embarking on its journey. To do so, it chooses, for each vertex $u\in V$, a
total order on the set of all arcs with source $u$, such that

\begin{enumerate}
\item if $u\neq r$, then the $B$-arc comes last in this order, and

\item if $u=r$, then the arc $a$ comes first in this order.
\end{enumerate}

Note that this total order can be chosen in $\left(  \deg^{+}u-1\right)  !$
many ways (since there are $\deg^{+}u$ arcs with source $u$, and we can freely
choose their order except that one of them has a fixed position). Thus, in
total, there are $\prod_{u\in V}\left(  \deg^{+}u-1\right)  !$ many options
for how the turtle can choose all these orders. Once these orders have been
chosen, the turtle then uses them to decide which arcs to walk along: Namely,
the first time it visits the vertex $u$, it leaves it along the first arc
(according to its chosen order); the second time, it uses the second arc; the
third time, the third arc; and so on.

So the turtle has $\prod_{u\in V}\left(  \deg^{+}u-1\right)  !$ many options,
and each of these options leads to a different $a$-Eulerian circuit
$\mathbf{e}$ (because the total orders chosen by the turtle are reflected in
$\mathbf{e}$: they are precisely the orders in which the respective arcs
appear in $\mathbf{e}$). Moreover, each $a$-Eulerian circuit $\mathbf{e}$
satisfying $\operatorname*{Exit}\mathbf{e}=B$ comes from one of these
options\footnote{\textit{Proof.} Let $\mathbf{e}$ be an $a$-Eulerian circuit
satisfying $\operatorname*{Exit}\mathbf{e}=B$. We must show that, by choosing
the appropriate total orders ahead of its journey, the turtle will trace this
exact circuit $\mathbf{e}$.
\par
First, let me explain what the \textquotedblleft appropriate total
orders\textquotedblright\ are: They are the orders dictated by $\mathbf{e}$.
That is, for each vertex $u\in V$, the turtle must choose the total order on
the set $\left\{  \text{all arcs with source }u\right\}  $ that coincides with
the order in which these arcs appear on $\mathbf{e}$. This choice is
legitimate, because the arc $a$ is the first arc of $\mathbf{e}$ (so it will
certainly come first in its order), and because each $B$-arc appears in
$\mathbf{e}$ after all other arcs from the same source have appeared (so it
will come last in its total order).
\par
Now, let me explain why the choice of these total orders will force the turtle
to trace the circuit $\mathbf{e}$. Indeed, assume the contrary. Let
$\mathbf{f}$ be the circuit traced by the turtle. Thus, by our assumption,
$\mathbf{f}$ is distinct from $\mathbf{e}$. Hence, $\mathbf{f}$ must diverge
from $\mathbf{e}$ at some point. Consider the first point at which this
happens. Thus, at this point, $\mathbf{f}$ leaves some vertex $u\in V$ along
an arc $e_{1}$, whereas $\mathbf{e}$ instead leaves it along a different arc
$e_{2}$. Since both $\mathbf{f}$ and $\mathbf{e}$ are Eulerian circuits, we
know that both $e_{1}$ and $e_{2}$ must appear in both $\mathbf{f}$ and
$\mathbf{e}$, and moreover, the circuits $\mathbf{f}$ and $\mathbf{e}$ agree
until this point of divergence. Hence, $e_{1}$ appears before $e_{2}$ in
$\mathbf{f}$ but appears after $e_{2}$ in $\mathbf{e}$. Therefore, the order
in which the arcs with source $u$ appear in $\mathbf{f}$ differs from the
order in which they appear in $\mathbf{e}$ (since both $e_{1}$ and $e_{2}$ are
arcs with source $u$). However, this is absurd, since the turtle's walk
$\mathbf{f}$ was constructed in such a way that it takes the arcs with source
$u$ in the same order as $\mathbf{e}$ does (because the total orders were
appropriately chosen). So we found a contradiction, and our proof is
complete.}.

Therefore, the total number of $a$-Eulerian circuits $\mathbf{e}$ satisfying
$\operatorname*{Exit}\mathbf{e}=B$ is the total number of options, which is
$\prod\limits_{u\in V}\left(  \deg^{+}u-1\right)  !$ as we know. This proves
Claim 2.] \medskip

With Claims 1 and 2 proved, we are almost done. The map%
\begin{align*}
\left\{  a\text{-Eulerian circuits of }D\right\}   &  \rightarrow\left\{
\text{sparbs}\right\}  ,\\
\mathbf{e}  &  \mapsto\operatorname*{Exit}\mathbf{e}%
\end{align*}
is well-defined (by Claim 1). Furthermore, Claim 2 shows that this map is a
$\prod\limits_{u\in V}\left(  \deg^{+}u-1\right)  !$-to-$1$
correspondence\footnote{An $m$\textbf{-to-}$1$\textbf{ correspondence} (where
$m$ is a nonnegative integer) means a map $f:X\rightarrow Y$ between two sets
such that each element of $Y$ has exactly $m$ preimages under $f$.} (i.e.,
each sparb $B$ has exactly \newline$\prod\limits_{u\in V}\left(  \deg
^{+}u-1\right)  !$ many preimages under this map). Thus, by the multijection
principle\footnote{The \textbf{multijection principle} is a basic counting
principle that says the following: Let $X$ and $Y$ be two finite sets, and let
$m\in\mathbb{N}$. Let $f:X\rightarrow Y$ be an $m$-to-$1$ correspondence
(i.e., a map such that each element of $Y$ has exactly $m$ preimages under
$f$). Then, $\left\vert X\right\vert =m\cdot\left\vert Y\right\vert $.
\par
For example, $n$ (intact) sheep have $4n$ legs in total, since the map that
sends each leg to its sheep is a $4$-to-$1$ correspondence.}, we conclude
that\footnote{The symbol \textquotedblleft$\#$\textquotedblright\ means
\textquotedblleft number\textquotedblright.}%
\[
\left(  \#\text{ of }a\text{-Eulerian circuits of }D\right)  =\left(
\prod\limits_{u\in V}\left(  \deg^{+}u-1\right)  !\right)  \cdot\left(
\#\text{ of sparbs}\right)  .
\]
Since $\varepsilon\left(  D,a\right)  =\left(  \#\text{ of }a\text{-Eulerian
circuits of }D\right)  $ and $\tau\left(  D,r\right)  =\left(  \#\text{ of
sparbs}\right)  $, we can rewrite this as follows:%
\[
\varepsilon\left(  D,a\right)  =\left(  \prod\limits_{u\in V}\left(  \deg
^{+}u-1\right)  !\right)  \cdot\tau\left(  D,r\right)  =\tau\left(
D,r\right)  \cdot\prod\limits_{u\in V}\left(  \deg^{+}u-1\right)  !.
\]
This proves Theorem \ref{thm.BEST.to}.
\end{proof}

\begin{proof}
[Proof of Theorem \ref{thm.BEST.from}.]As we already mentioned, Theorem
\ref{thm.BEST.from} follows from Theorem \ref{thm.BEST.to} by reversing each
arc (i.e., by applying Theorem \ref{thm.BEST.to} to the digraph
$D^{\operatorname*{rev}}$ instead of $D$).
\end{proof}

The above proof of Theorem \ref{thm.BEST.to} originates from \cite[\S 6]%
{VanEhr51}; I have learned it from \cite[Theorem 10.2]{Stanley-AC} and
\cite[Corollary 4.10]{HLMPPW13} (where it appears in a mildly modified form).

\subsection{A corollary about spanning arborescences}

Before we actually use the BEST (or BEST') theorem to count the Eulerian
circuits on any digraph, let us mention a neat corollary for the number of
spanning arborescences:

\begin{corollary}
\label{cor.BEST.tau=tau}Let $D=\left(  V,A,\psi\right)  $ be a balanced
multidigraph. For each vertex $r\in V$, let $\tau\left(  D,r\right)  $ be the
number of spanning arborescences of $D$ rooted to $r$. Then, $\tau\left(
D,r\right)  $ does not depend on $r$.
\end{corollary}

\begin{proof}
[Proof of Corollary \ref{cor.BEST.tau=tau}.]WLOG assume that $\left\vert
V\right\vert >1$ (else, the claim is obvious). If there is a vertex $v\in V$
with $\deg^{+}v=0$, then this vertex $v$ satisfies $\deg^{-}v=0$ as well
(since the balancedness of $D$ entails $\deg^{-}v=\deg^{+}v=0$), and therefore
$D$ has no spanning arborescences at all (since any spanning arborescence
would have an arc with source or target $v$). Thus, we WLOG assume that
$\deg^{+}v>0$ for all $v\in V$. In other words, each vertex has outdegree $>0$.

Let $r$ and $s$ be two vertices of $D$. We must prove that $\tau\left(
D,r\right)  =\tau\left(  D,s\right)  $.

Pick an arc $a$ with source $r$. (This exists, since $\deg^{+}r>0$.) Pick an
arc $b$ with source $s$. (This exists, since $\deg^{+}s>0$.)

Applying the BEST' theorem (Theorem \ref{thm.BEST.to}), we get
\begin{align*}
\varepsilon\left(  D,a\right)   &  =\tau\left(  D,r\right)  \cdot\prod_{u\in
V}\left(  \deg^{+}u-1\right)  !\ \ \ \ \ \ \ \ \ \ \text{and similarly}\\
\varepsilon\left(  D,b\right)   &  =\tau\left(  D,s\right)  \cdot\prod_{u\in
V}\left(  \deg^{+}u-1\right)  !.
\end{align*}
However, $\varepsilon\left(  D,a\right)  =\varepsilon\left(  D,b\right)  $,
since counting Eulerian circuits that start with $a$ is equivalent to counting
Eulerian circuits that start with $b$ (because an Eulerian circuit can be
rotated uniquely to start with any given arc). Thus, we obtain%
\[
\tau\left(  D,r\right)  \cdot\prod_{u\in V}\left(  \deg^{+}u-1\right)
!=\varepsilon\left(  D,a\right)  =\varepsilon\left(  D,b\right)  =\tau\left(
D,s\right)  \cdot\prod_{u\in V}\left(  \deg^{+}u-1\right)  !.
\]
Cancelling the (nonzero!) number $\prod_{u\in V}\left(  \deg^{+}u-1\right)  !$
from this equality, we obtain $\tau\left(  D,r\right)  =\tau\left(
D,s\right)  $. This proves Corollary \ref{cor.BEST.tau=tau}.
\end{proof}

\subsection{Spanning arborescences vs. spanning trees}

The BEST theorem (Theorem \ref{thm.BEST.to} or Theorem \ref{thm.BEST.from})
connects the $\#$ of Eulerian circuits in a digraph with the $\#$ of spanning
arborescences of the same digraph. Now let us try to find a way to compute the latter.

For example, let us try to do this for digraphs of the form
$G^{\operatorname*{bidir}}$ where $G$ is a multigraph. I claim that the
spanning arborescences of $G^{\operatorname*{bidir}}$ rooted to a given vertex
$r$ are just the spanning trees of $G$ in disguise:

\begin{proposition}
\label{prop.sparb.vs-sptree}Let $G=\left(  V,E,\varphi\right)  $ be a
multigraph. Fix a vertex $r\in V$. Recall that the arcs of
$G^{\operatorname*{bidir}}$ are the pairs $\left(  e,i\right)  \in
E\times\left\{  1,2\right\}  $. Identify each spanning tree of $G$ with its
edge set, and each spanning arborescence of $G^{\operatorname*{bidir}}$ with
its arc set.

If $B$ is a spanning arborescence of $G^{\operatorname*{bidir}}$ rooted to
$r$, then we set%
\[
\overline{B}:=\left\{  e\ \mid\ \left(  e,i\right)  \in B\right\}  .
\]
(Recall that we are identifying spanning arborescences with their arc sets, so
that \textquotedblleft$\left(  e,i\right)  \in B$\textquotedblright\ means
\textquotedblleft$\left(  e,i\right)  $ is an arc of $B$\textquotedblright.)

Then:

\begin{enumerate}
\item[\textbf{(a)}] If $B$ is a spanning arborescence of
$G^{\operatorname*{bidir}}$ rooted to $r$, then $\overline{B}$ is a spanning
tree of $G$.

\item[\textbf{(b)}] The map%
\begin{align*}
\left\{  \text{spanning arborescences of }G^{\operatorname*{bidir}}\text{
rooted to }r\right\}   &  \rightarrow\left\{  \text{spanning trees of
}G\right\}  ,\\
B  &  \mapsto\overline{B}%
\end{align*}
is a bijection.
\end{enumerate}
\end{proposition}

\Needspace{45pc}

\begin{example}
Here is a multigraph $G$ (on the left) with the corresponding multidigraph
$G^{\operatorname*{bidir}}$ (on the right):%
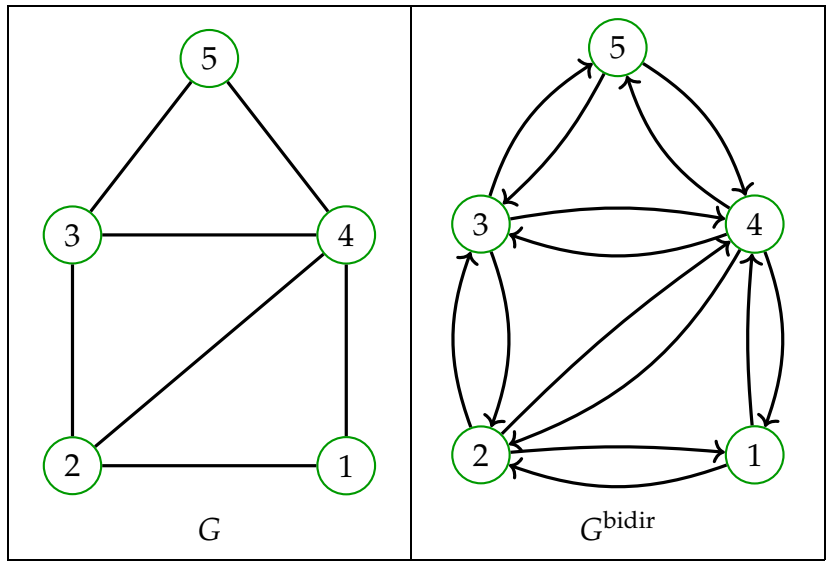
\[%

\ \
\]
Here is a spanning arborescence $B$ of $G^{\operatorname*{bidir}}$ rooted to
$1$, and the corresponding spanning tree $\overline{B}$ of $G$:%
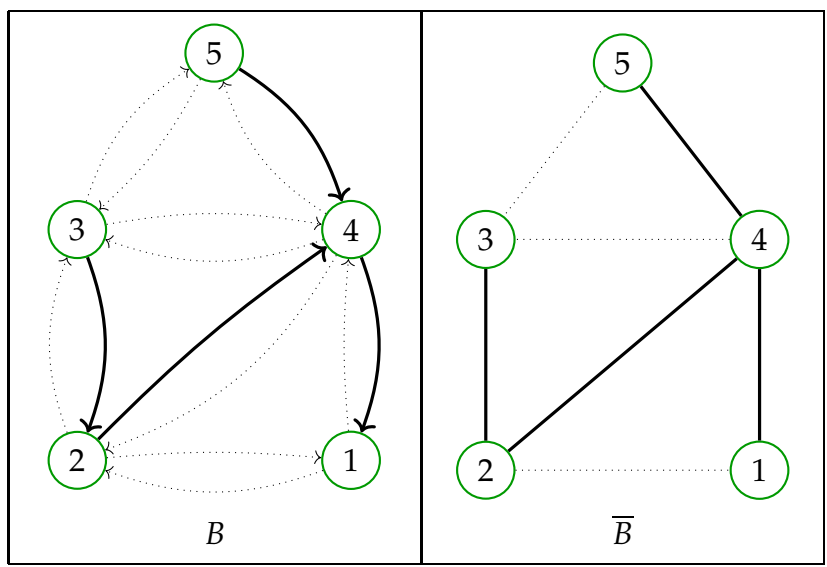
\[%
%
\ \
\]
(here, the arcs of $G^{\operatorname*{bidir}}$ that \textbf{don't} belong to
$B$, as well as the edges of $G$ that \textbf{don't} belong to $\overline{B}$,
have been drawn as dotted arrows). It is fairly easy to see how $B$ can be
reconstructed from $\overline{B}$: You just need to replace each edge of
$\overline{B}$ by the appropriately directed arc (namely, the one that is
\textquotedblleft directed towards $1$\textquotedblright).
\end{example}

\begin{proof}
[Proof of Proposition \ref{prop.sparb.vs-sptree}.]This is an exercise in
yak-shaving (and we have, in fact, shaved a very similar yak in Section
\ref{sec.arbor.vs-tree}; the only difference is that we are no longer dealing
with trees in isolation, but rather with spanning trees of $G$). \medskip

\textbf{(a)} Let $B$ be a spanning arborescence of $G^{\operatorname*{bidir}}$
rooted to $r$. Then, $B^{\operatorname*{und}}$ is a tree (by the implication
A'1$\Longrightarrow$A'3 in Theorem \ref{thm.arbor.eq-A-dual}). However, it is
easy to see that $B^{\operatorname*{und}}\cong\overline{B}$ as multigraphs
(indeed, each vertex $v$ of $B^{\operatorname*{und}}$ corresponds to the same
vertex $v$ of $\overline{B}$, whereas any edge $\left(  e,i\right)  $ of
$B^{\operatorname*{und}}$ corresponds to the edge $e$ of $\overline{B}%
$)\ \ \ \ \footnote{Here we need to use the fact that for each edge $e$ of
$\overline{B}$, \textbf{exactly one} of the two pairs $\left(  e,1\right)  $
and $\left(  e,2\right)  $ is an edge of $B^{\operatorname*{und}}$. But this
is easy to check: At least one of the two pairs $\left(  e,1\right)  $ and
$\left(  e,2\right)  $ must be an arc of $B$ (since $e$ is an edge of
$\overline{B}$). In other words, at least one of the two pairs $\left(
e,1\right)  $ and $\left(  e,2\right)  $ must be an edge of
$B^{\operatorname*{und}}$. But both of these pairs cannot be edges of
$B^{\operatorname*{und}}$ at the same time (since this would create a cycle,
but $B^{\operatorname*{und}}$ is a tree and thus has no cycles). Hence,
exactly one of these pairs is an edge of $B^{\operatorname*{und}}$, qed.}.
Thus, $\overline{B}$ is a tree (since $B^{\operatorname*{und}}$ is a
tree)\footnote{Alternatively, you can prove this as follows: The vertex $r$ is
a to-root of $B$ (since $B$ is an arborescence rooted to $r$). Thus, for each
$v\in V$, there is a path from $v$ to $r$ in $B$. By \textquotedblleft
projecting\textquotedblright\ this path onto $\overline{B}$ (that is,
replacing each arc $\left(  e,i\right)  $ of this path by the corresponding
edge $e$ of $\overline{B}$), we obtain a path from $v$ to $r$ in $\overline
{B}$. This shows that the multigraph $\overline{B}$ is connected. Furthermore,
the definition of $\overline{B}$ shows that $\left\vert \overline
{B}\right\vert \leq\left\vert B\right\vert =\left\vert V\right\vert -1$ (by
Statement A'2 in Theorem \ref{thm.arbor.eq-A-dual}, since $B$ is an
arborescence rooted to $r$). Hence, $\left\vert \overline{B}\right\vert
<\left\vert V\right\vert $. Thus, we can apply the implication
T5$\Longrightarrow$T1 of the Tree Equivalence Theorem (Theorem
\ref{thm.trees.T1-8}) to conclude that $\overline{B}$ is a tree.}, therefore a
spanning tree of $G$ (since $\overline{B}$ is clearly a spanning subgraph of
$G$). This proves Proposition \ref{prop.sparb.vs-sptree} \textbf{(a)}.
\medskip

\textbf{(b)} We must prove that this map is surjective and injective.

\textit{Surjectivity:} Let $T$ be a spanning tree of $G$. Then, the
multidigraph $T^{r\rightarrow}$ (defined in Definition \ref{def.tree.rto}) is
an arborescence rooted from $r$ (by Lemma \ref{lem.arbor-vs-tree.1}).
Reversing each arc in this arborescence $T^{r\rightarrow}$, we obtain a new
multidigraph $T^{r\leftarrow}$, which is thus an arborescence rooted to $r$.
Unfortunately, $T^{r\leftarrow}$ is not a subdigraph of
$G^{\operatorname*{bidir}}$, for a rather stupid reason: The arcs of
$T^{r\leftarrow}$ are elements of $E$, whereas the arcs of
$G^{\operatorname*{bidir}}$ are pairs of the form $\left(  e,i\right)  $ with
$e\in E$ and $i\in\left\{  1,2\right\}  $.

Fortunately, this is easily fixed: For each arc $e$ of $T^{r\leftarrow}$, we
let $e^{\prime}$ be the arc $\left(  e,i\right)  $ of
$G^{\operatorname*{bidir}}$ that has the same source as $e$ (and thus the same
target as $e$). This is uniquely determined, since the arcs $\left(
e,1\right)  $ and $\left(  e,2\right)  $ of $G^{\operatorname*{bidir}}$ have
different sources\footnote{\textit{Proof.} The edge $e$ of $T$ is not a loop
(because $T$ is a tree, but a tree cannot have any loops). Hence, its two
endpoints are distinct. Thus, the arcs $\left(  e,1\right)  $ and $\left(
e,2\right)  $ of $G^{\operatorname*{bidir}}$ have different sources (since
their sources are the two endpoints of $e$).}. If we replace each arc $e$ of
$T^{r\leftarrow}$ by the corresponding arc $e^{\prime}$ of
$G^{\operatorname*{bidir}}$, then we obtain a spanning subdigraph $S$ of
$G^{\operatorname*{bidir}}$ that is an arborescence rooted to $r$ (since
$T^{r\leftarrow}$ is an arborescence rooted to $r$, and we have only replaced
its arcs by equivalent ones with the same sources and the same targets). In
other words, we obtain a spanning arborescence $S$ of
$G^{\operatorname*{bidir}}$ rooted to $r$. It is easy to see that
$\overline{S}=T$. Hence, the map%
\begin{align*}
\left\{  \text{spanning arborescences of }G^{\operatorname*{bidir}}\text{
rooted to }r\right\}   &  \rightarrow\left\{  \text{spanning trees of
}G\right\}  ,\\
B  &  \mapsto\overline{B}%
\end{align*}
sends $S$ to $T$. This shows that $T$ is a value of this map. Since we have
proved this for every spanning tree $T$ of $G$, we have thus shown that this
map is surjective.

\textit{Injectivity:} The main idea is that, in order to recover a spanning
arborescence $B$ back from the corresponding spanning tree $\overline{B}$, we
just need to \textquotedblleft orient the edges of the tree towards
$r$\textquotedblright. Here are the (annoyingly long) details:

Let $B$ and $C$ be two sparbs\footnote{Henceforth, \textquotedblleft
sparb\textquotedblright\ is short for \textquotedblleft spanning arborescence
of $G^{\operatorname*{bidir}}$ rooted to $r$\textquotedblright.} such that
$\overline{B}=\overline{C}$. We must show that $B=C$.

Assume the contrary. Thus, $B\neq C$. Let $T$ be the tree $\overline
{B}=\overline{C}$. Thus, each edge $e$ of $T$ corresponds to either an arc
$\left(  e,1\right)  $ or an arc $\left(  e,2\right)  $ in $B$ (since
$T=\overline{B}$), and likewise for $C$. Conversely, each arc $\left(
e,i\right)  $ of $B$ or of $C$ corresponds to an edge $e$ of $T$. Hence, from
$B\neq C$, we see that there must exist an edge $e$ of $T$ such that

\begin{itemize}
\item \textbf{either} we have $\left(  e,1\right)  \in B$ and $\left(
e,2\right)  \in C$,

\item \textbf{or} we have $\left(  e,1\right)  \in C$ and $\left(  e,2\right)
\in B$.
\end{itemize}

Consider this edge $e$. We WLOG assume that $\left(  e,1\right)  \in B$ and
$\left(  e,2\right)  \in C$ (else, we can just swap $B$ with $C$). Let the arc
$\left(  e,1\right)  $ of $G^{\operatorname*{bidir}}$ have source $s$ and
target $t$, so that $\left(  e,2\right)  $ has source $t$ and target $s$. The
edge $e$ thus has endpoints $s$ and $t$.

Since $B$ is an arborescence rooted to $r$, the vertex $r$ is a to-root of
$B$. Hence, there exists a path $\mathbf{p}$ from $s$ to $r$ in $B$. This path
$\mathbf{p}$ must begin with the arc $\left(  e,1\right)  $%
\ \ \ \ \footnote{\textit{Proof.} Since $r$ is a to-root of $B$, we know that
there exists a path from $t$ to $r$ in $B$. Let $\mathbf{t}$ be this path.
Extending this path $\mathbf{t}$ by the vertex $s$ and the arc $\left(
e,1\right)  $ (which we both insert at the start of $\mathbf{t}$), we obtain a
walk $\mathbf{t}^{\prime}$ from $s$ to $r$ in $B$. (So, if $\mathbf{t}=\left(
t,\ \ldots,\ r\right)  $, then $\mathbf{t}^{\prime}=\left(  s,\ \left(
e,1\right)  ,\ t,\ \ldots,\ r\right)  $.)
\par
However, $B$ is an arborescence rooted to $r$. Thus, Statement A'4 in the Dual
Arborescence Equivalence Theorem (Theorem \ref{thm.arbor.eq-A-dual}) shows
that for each vertex $v\in V$, the digraph $B$ has a unique walk from $v$ to
$r$. Hence, in particular, $B$ has a unique walk from $s$ to $r$. Thus,
$\mathbf{p}=\mathbf{t}^{\prime}$ (since both $\mathbf{p}$ and $\mathbf{t}%
^{\prime}$ are walks from $s$ to $r$ in $B$). Since $\mathbf{t}^{\prime}$
begins with the arc $\left(  e,1\right)  $, we thus conclude that $\mathbf{p}$
begins with the arc $\left(  e,1\right)  $.}. Projecting this path
$\mathbf{p}$ down onto $T$, we obtain a path $\overline{\mathbf{p}}$ from $s$
to $r$ in $T$. (By the word \textquotedblleft projecting\textquotedblright, we
mean replacing each arc $\left(  e,i\right)  $ by the corresponding edge $e$.
Clearly, doing this to a path in $B$ yields a path in $T$, because
$T=\overline{B}$.) Since the path $\mathbf{p}$ begins with the arc $\left(
e,1\right)  $, the \textquotedblleft projected\textquotedblright\ path
$\overline{\mathbf{p}}$ begins with the edge $e$. Thus, in the tree $T$, the
path from $s$ to $r$ begins with the edge $e$ (because this path must be the
path $\overline{\mathbf{p}}$). As a consequence, $t$ must be the second vertex
of this path (since the edge $e$ has endpoints $s$ and $t$), so that removing
the first edge from this path yields the path from $t$ to $r$. Thus, $d\left(
t,r\right)  =d\left(  s,r\right)  -1$, where $d$ denotes distance on the tree
$T$. Hence, $d\left(  t,r\right)  <d\left(  s,r\right)  $.

A similar argument (but with the roles of $B$ and $C$ swapped, as well as the
roles of $s$ and $t$ swapped, and the roles of $\left(  e,1\right)  $ and
$\left(  e,2\right)  $ swapped) shows that $d\left(  s,r\right)  <d\left(
t,r\right)  $. But this contradicts $d\left(  t,r\right)  <d\left(
s,r\right)  $.

This contradiction shows that our assumption was false. Thus, we have proved
that $B=C$.

Forget that we fixed $B$ and $C$. We thus have shown that if $B$ and $C$ are
two sparbs such that $\overline{B}=\overline{C}$, then $B=C$. In other words,
our map
\begin{align*}
\left\{  \text{spanning arborescences of }G^{\operatorname*{bidir}}\text{
rooted to }r\right\}   &  \rightarrow\left\{  \text{spanning trees of
}G\right\}  ,\\
B  &  \mapsto\overline{B}%
\end{align*}
is injective.

We have now shown that this map is both surjective and injective. Hence, it is
a bijection. This proves Proposition \ref{prop.sparb.vs-sptree} \textbf{(b)}.
\end{proof}

\subsection{The matrix-tree theorem}

\subsubsection{Introduction}

So counting spanning trees in a multigraph is a particular case of counting
spanning arborescences (rooted to a given vertex) in a multidigraph. But how
do we do either? Let us begin with some simple examples:

\begin{example}
There is only one spanning tree of the complete graph $K_{1}$:%
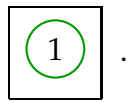
\[%
\begin{tabular}
[c]{|c|}\hline
$%
\begin{tikzpicture}[scale=1.8]
\begin{scope}[every node/.style={circle,thick,draw=green!60!black}]
\node(1) at (0,0) {$1$};
\end{scope}
\begin{scope}[every edge/.style={draw=black,very thick}, every loop/.style={}]
\end{scope}
\end{tikzpicture}%
$\\\hline
\end{tabular}
\ \ .
\]

There is only one spanning tree of the complete graph $K_{2}$:%
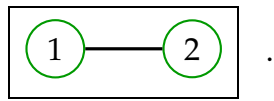
\[%
\begin{tabular}
[c]{|c|}\hline
$%
\begin{tikzpicture}[scale=1.8]
\begin{scope}[every node/.style={circle,thick,draw=green!60!black}]
\node(1) at (0,0) {$1$};
\node(2) at (1,0) {$2$};
\end{scope}
\begin{scope}[every edge/.style={draw=black,very thick}, every loop/.style={}]
\path[-] (1) edge (2);
\end{scope}
\end{tikzpicture}%
$\\\hline
\end{tabular}
\ \ \ .
\]

There are $3$ spanning trees of the complete graph $K_{3}$:%
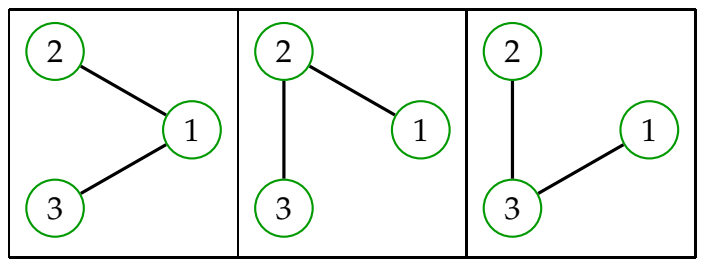
\[%
\begin{tabular}
[c]{|c|c|c|}\hline
$%
\begin{tikzpicture}[scale=1.2]
\begin{scope}[every node/.style={circle,thick,draw=green!60!black}]
\node(1) at (0:1) {$1$};
\node(2) at (120:1) {$2$};
\node(3) at (240:1) {$3$};
\end{scope}
\begin{scope}[every edge/.style={draw=black,very thick}, every loop/.style={}]
\path[-] (1) edge (2) edge (3);
\end{scope}
\end{tikzpicture}%
$ & $%
\begin{tikzpicture}[scale=1.2]
\begin{scope}[every node/.style={circle,thick,draw=green!60!black}]
\node(1) at (0:1) {$1$};
\node(2) at (120:1) {$2$};
\node(3) at (240:1) {$3$};
\end{scope}
\begin{scope}[every edge/.style={draw=black,very thick}, every loop/.style={}]
\path[-] (2) edge (1) edge (3);
\end{scope}
\end{tikzpicture}%
$ & $%
\begin{tikzpicture}[scale=1.2]
\begin{scope}[every node/.style={circle,thick,draw=green!60!black}]
\node(1) at (0:1) {$1$};
\node(2) at (120:1) {$2$};
\node(3) at (240:1) {$3$};
\end{scope}
\begin{scope}[every edge/.style={draw=black,very thick}, every loop/.style={}]
\path[-] (3) edge (1) edge (2);
\end{scope}
\end{tikzpicture}%
$\\\hline
\end{tabular}
\ \ .
\]
(They are all isomorphic, but still distinct.)

\newpage

There are $16$ spanning trees of the complete graph $K_{4}$:%
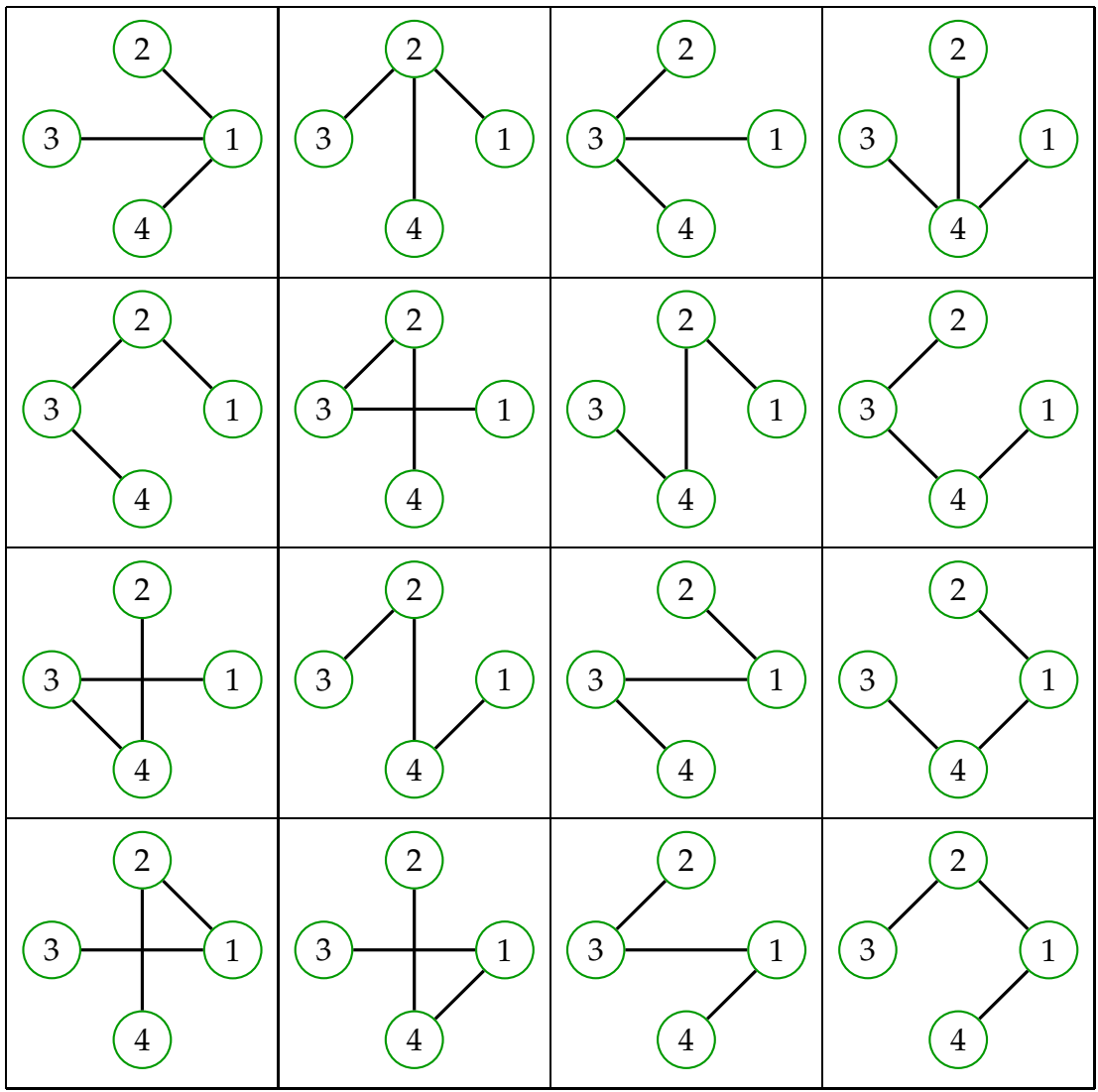
\[%
%
\ \ .
\]
(There are only two non-isomorphic ones among them.)
\end{example}

This example suggests that the \# of spanning trees of a complete graph
$K_{n}$ is $n^{n-2}$.

This is indeed true, and we will prove this later. First, however, let us
discuss the more general problem of counting spanning arborescences of an
arbitrary digraph $D$.

\subsubsection{Notations}

We will use the following convenient shorthand notation:

\begin{definition}
\label{def.iverson} We will use the \textbf{Iverson bracket notation}: If
$\mathcal{A}$ is any logical statement, then we set%
\[
\left[  \mathcal{A}\right]  :=%
\begin{cases}
1, & \text{if }\mathcal{A}\text{ is true;}\\
0, & \text{if }\mathcal{A}\text{ is false.}%
\end{cases}
\]

\end{definition}

For example, $\left[  K_{2}\text{ is a tree}\right]  =1$ whereas $\left[
K_{3}\text{ is a tree}\right]  =0$.

\begin{definition}
Let $M$ be a matrix. Let $i$ and $j$ be two integers. Then,%
\begin{align*}
&  M_{i,j}\text{ will mean the entry of }M\text{ in row }i\text{ and column
}j\text{;}\\
&  M_{\sim i,\sim j}\text{ will mean the matrix }M\text{ with row }i\text{
removed and column }j\text{ removed.}%
\end{align*}

For example,
\[
\left(
\begin{array}
[c]{ccc}%
a & b & c\\
d & e & f\\
g & h & i
\end{array}
\right)  _{2,3}=f\ \ \ \ \ \ \ \ \ \ \text{and}\ \ \ \ \ \ \ \ \ \ \left(
\begin{array}
[c]{ccc}%
a & b & c\\
d & e & f\\
g & h & i
\end{array}
\right)  _{\sim2,\sim3}=\left(
\begin{array}
[c]{cc}%
a & b\\
g & h
\end{array}
\right)  .
\]

\end{definition}

\subsubsection{The Laplacian of a multidigraph}

We shall now assign a matrix to (more or less) any
multidigraph:\footnote{Recall that the symbol \textquotedblleft$\#$%
\textquotedblright\ means \textquotedblleft number\textquotedblright.}

\begin{definition}
\label{def.MTT.L}Let $D=\left(  V,A,\psi\right)  $ be a multidigraph. Assume
that $V=\left\{  1,2,\ldots,n\right\}  $ for some $n\in\mathbb{N}$.

For any $i,j\in V$, we let $a_{i,j}$ be the $\#$ of arcs of $D$ that have
source $i$ and target $j$.

The \textbf{Laplacian} of $D$ is defined to be the $n\times n$-matrix
$L\in\mathbb{Z}^{n\times n}$ whose entries are given by%
\[
L_{i,j}=\left(  \deg^{+}i\right)  \cdot\underbrace{\left[  i=j\right]
}_{\substack{\text{This is also}\\\text{known as }\delta_{i,j}}}-\,a_{i,j}%
\ \ \ \ \ \ \ \ \ \ \text{for all }i,j\in V.
\]
In other words, it is the matrix%
\[
L=\left(
\begin{array}
[c]{cccc}%
\deg^{+}1-a_{1,1} & -a_{1,2} & \cdots & -a_{1,n}\\
-a_{2,1} & \deg^{+}2-a_{2,2} & \cdots & -a_{2,n}\\
\vdots & \vdots & \ddots & \vdots\\
-a_{n,1} & -a_{n,2} & \cdots & \deg^{+}n-a_{n,n}%
\end{array}
\right)  .
\]

\end{definition}

\Needspace{20pc}

\begin{example}
Let $D$ be the digraph%
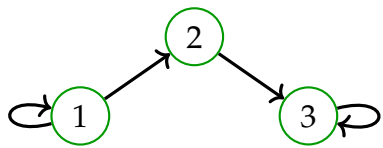
\[%
\begin{tikzpicture}[scale=1.5]
\begin{scope}[every node/.style={circle,thick,draw=green!60!black}]
\node(1) at (-1,0) {$1$};
\node(2) at (0,0.7) {$2$};
\node(3) at (1,0) {$3$};
\end{scope}
\begin{scope}[every edge/.style={draw=black,very thick}, every loop/.style={}]
\path
[->] (1) edge[loop left] (1) edge (2) (2) edge (3) (3) edge[loop right] (3);\end{scope}
\end{tikzpicture}%
\ \ .
\]
Then, its Laplacian is%
\[
\left(
\begin{array}
[c]{ccc}%
2-1 & -1 & -0\\
-0 & 1-0 & -1\\
-0 & -0 & 1-1
\end{array}
\right)  =\left(
\begin{array}
[c]{ccc}%
1 & -1 & 0\\
0 & 1 & -1\\
0 & 0 & 0
\end{array}
\right)  .
\]

\end{example}

One thing we notice from this example is that loops do not matter at all to
the Laplacian $L$. Indeed, a loop with source $i$ and target $i$ counts once
in $\deg^{+}i$ and once in $a_{i,i}$, but these contributions cancel out.

Here is a simple property of Laplacians:

\begin{proposition}
\label{prop.MTT.detL=0}Let $D=\left(  V,A,\psi\right)  $ be a multidigraph.
Assume that $V=\left\{  1,2,\ldots,n\right\}  $ for some positive integer $n$.

Then, the Laplacian $L$ of $D$ is singular; i.e., we have $\det L=0$.
\end{proposition}

\begin{proof}
The sum of all columns of $L$ is the zero vector, because for each $i\in V$ we
have%
\begin{align*}
\sum_{j=1}^{n}L_{i,j}  &  =\sum_{j=1}^{n}\left(  \left(  \deg^{+}i\right)
\cdot\left[  i=j\right]  -a_{i,j}\right)  \ \ \ \ \ \ \ \ \ \ \left(  \text{by
the definition of }L\right) \\
&  =\underbrace{\sum_{j=1}^{n}\left(  \deg^{+}i\right)  \cdot\left[
i=j\right]  }_{\substack{=\deg^{+}i\\\text{(since only the addend}\\\text{for
}j=i\text{ can be nonzero)}}}-\underbrace{\sum_{j=1}^{n}a_{i,j}}%
_{\substack{=\deg^{+}i\\\text{(since this is counting}\\\text{all arcs with
source }i\text{)}}}\\
&  =\deg^{+}i-\deg^{+}i=0.
\end{align*}
In other words, we have $Le=0$ for the vector $e:=\left(  1,1,\ldots,1\right)
^{T}$. Thus, this vector $e$ lies in the kernel (aka nullspace) of $L$, and so
$L$ is singular.

(Note that we used the positivity of $n$ here! If $n=0$, then $e$ is the zero
vector, because a vector with $0$ entries is automatically the zero vector.)
\end{proof}

\subsubsection{The Matrix-Tree Theorem: statement}

Proposition \ref{prop.MTT.detL=0} shows that the determinant of the Laplacian
of a digraph is not very interesting. It is common, however, that when a
matrix has determinant $0$, its largest nonzero minors (= determinants of
submatrices) often carry some interesting information; they are
\textquotedblleft the closest the matrix has\textquotedblright\ to a nonzero
determinant. In the case of the Laplacian, they turn out to count spanning arborescences:

\begin{theorem}
[Matrix-Tree Theorem]\label{thm.MTT.MTT}Let $D=\left(  V,A,\psi\right)  $ be a
multidigraph. Assume that $V=\left\{  1,2,\ldots,n\right\}  $ for some
positive integer $n$.

Let $L$ be the Laplacian of $D$. Let $r$ be a vertex of $D$. Then,%
\[
\left(  \#\text{ of spanning arborescences of }D\text{ rooted to }r\right)
=\det\left(  L_{\sim r,\sim r}\right)  .
\]

\end{theorem}

Before we prove this, some remarks:

\begin{itemize}
\item The determinant $\det\left(  L_{\sim r,\sim r}\right)  $ is the $\left(
r,r\right)  $-th entry of
\href{https://en.wikipedia.org/wiki/Adjugate_matrix}{the adjugate matrix} of
$L$.

\item The $V=\left\{  1,2,\ldots,n\right\}  $ assumption is a typical
\textquotedblleft WLOG assumption\textquotedblright: If you have an arbitrary
digraph $D$, you can always rename its vertices as $1,2,\ldots,n$, and then
this assumption will be satisfied. Thus, Theorem \ref{thm.MTT.MTT} helps you
count the spanning arborescences of any digraph. That said, you can also drop
the $V=\left\{  1,2,\ldots,n\right\}  $ assumption from Theorem
\ref{thm.MTT.MTT} if you are okay with matrices whose rows and columns are
indexed not by numbers but by elements of an arbitrary finite
set\footnote{Such matrices are perfectly fine, just somewhat unusual and hard
to write down (which row do you put on top?). See
\url{https://mathoverflow.net/questions/317105} for details.}.
\end{itemize}

\subsubsection{\label{subsec.MTT.Kn}Application: Counting the spanning trees
of $K_{n}$}

Now, let us use the Matrix-Tree Theorem to count the spanning trees of $K_{n}%
$. This should provide some intuition for the theorem before we come to its proof.

We fix a positive integer $n$. Let $L$ be the Laplacian of the multidigraph
$K_{n}^{\operatorname*{bidir}}$ (where $K_{n}$, as we recall, is the complete
graph on the set $\left\{  1,2,\ldots,n\right\}  $). Then, each vertex of
$K_{n}^{\operatorname*{bidir}}$ has outdegree $n-1$, and thus we have%
\begin{equation}
L=\left(
\begin{array}
[c]{cccc}%
n-1 & -1 & \cdots & -1\\
-1 & n-1 & \cdots & -1\\
\vdots & \vdots & \ddots & \vdots\\
-1 & -1 & \cdots & n-1
\end{array}
\right)  \label{eq.MTT.Kn.L=}%
\end{equation}
(this is the $n\times n$-matrix whose diagonal entries are $n-1$ and whose
off-diagonal entries are $-1$). By Proposition \ref{prop.sparb.vs-sptree}
\textbf{(b)} (applied to $G=K_{n}$ and $r=1$), there is a bijection between
$\left\{  \text{spanning arborescences of }K_{n}^{\operatorname*{bidir}}\text{
rooted to }1\right\}  $ and $\left\{  \text{spanning trees of }K_{n}\right\}
$. Hence, by the bijection principle, we have%
\begin{align*}
&  \left(  \#\text{ of spanning trees of }K_{n}\right) \\
&  =\left(  \#\text{ of spanning arborescences of }K_{n}%
^{\operatorname*{bidir}}\text{ rooted to }1\right) \\
&  =\det\left(  L_{\sim1,\sim1}\right)  \ \ \ \ \ \ \ \ \ \ \left(  \text{by
Theorem \ref{thm.MTT.MTT}, applied to }D=K_{n}^{\operatorname*{bidir}}\text{
and }r=1\right) \\
&  =\det\underbrace{\left(
\begin{array}
[c]{cccc}%
n-1 & -1 & \cdots & -1\\
-1 & n-1 & \cdots & -1\\
\vdots & \vdots & \ddots & \vdots\\
-1 & -1 & \cdots & n-1
\end{array}
\right)  }_{\text{an }\left(  n-1\right)  \times\left(  n-1\right)
\text{-matrix}}\ \ \ \ \ \ \ \ \ \ \left(  \text{by (\ref{eq.MTT.Kn.L=}%
)}\right)  .
\end{align*}
How do we compute this determinant? Here are three ways:\footnote{All three
ways require $n\geq2$, but this is okay since the $n=1$ case is trivial.}

\begin{itemize}
\item The most elementary approach is using row transformations:%
\begin{align*}
&  \det\left(
\begin{array}
[c]{cccc}%
n-1 & -1 & \cdots & -1\\
-1 & n-1 & \cdots & -1\\
\vdots & \vdots & \ddots & \vdots\\
-1 & -1 & \cdots & n-1
\end{array}
\right) \\
&  =\det\left(
\begin{array}
[c]{cccccc}%
n-1 & -1 & -1 & -1 & \cdots & -1\\
-n & n & 0 & 0 & \cdots & 0\\
-n & 0 & n & 0 & \cdots & 0\\
-n & 0 & 0 & n & \cdots & 0\\
\vdots & \vdots & \vdots & \vdots & \ddots & \vdots\\
-n & 0 & 0 & 0 & \cdots & n
\end{array}
\right)  \ \ \ \ \ \ \ \ \ \ \left(
\begin{array}
[c]{c}%
\text{here, we have}\\
\text{subtracted the }1\text{st row}\\
\text{from each other row}%
\end{array}
\right) \\
&  =n^{n-2}\det\left(
\begin{array}
[c]{cccccc}%
n-1 & -1 & -1 & -1 & \cdots & -1\\
-1 & 1 & 0 & 0 & \cdots & 0\\
-1 & 0 & 1 & 0 & \cdots & 0\\
-1 & 0 & 0 & 1 & \cdots & 0\\
\vdots & \vdots & \vdots & \vdots & \ddots & \vdots\\
-1 & 0 & 0 & 0 & \cdots & 1
\end{array}
\right)  \ \ \ \ \ \ \ \ \ \ \left(
\begin{array}
[c]{c}%
\text{here, we have}\\
\text{factored out}\\
\text{an }n\text{ from each}\\
\text{row except for}\\
\text{the first row}%
\end{array}
\right) \\
&  =n^{n-2}\underbrace{\det\left(
\begin{array}
[c]{cccccc}%
1 & 0 & 0 & 0 & \cdots & 0\\
-1 & 1 & 0 & 0 & \cdots & 0\\
-1 & 0 & 1 & 0 & \cdots & 0\\
-1 & 0 & 0 & 1 & \cdots & 0\\
\vdots & \vdots & \vdots & \vdots & \ddots & \vdots\\
-1 & 0 & 0 & 0 & \cdots & 1
\end{array}
\right)  }_{\substack{=1\\\text{(since the matrix is triangular}\\\text{with
diagonal entries }1,1,\ldots,1\text{)}}}\ \ \ \ \ \ \ \ \ \ \left(
\begin{array}
[c]{c}%
\text{here, we have added the }2\text{nd,}\\
3\text{rd, etc. rows to the }1\text{st row}%
\end{array}
\right) \\
&  =n^{n-2}.
\end{align*}

\item \href{https://en.wikipedia.org/wiki/Matrix_determinant_lemma}{The
so-called \textbf{matrix determinant lemma}}\footnote{also known as (and
probably much better titled) the \textquotedblleft rank-one perturbation
formula for the determinant\textquotedblright} says that for any $m\times
m$-matrix $A\in\mathbb{R}^{m\times m}$, any column vector $u\in\mathbb{R}%
^{m\times1}$ and any row vector $v\in\mathbb{R}^{1\times m}$, we have%
\[
\det\left(  A+uv\right)  =\det A+v\left(  \operatorname*{adj}A\right)  u.
\]
This helps us compute our determinant, since%
\begin{align*}
&  \left(
\begin{array}
[c]{cccc}%
n-1 & -1 & \cdots & -1\\
-1 & n-1 & \cdots & -1\\
\vdots & \vdots & \ddots & \vdots\\
-1 & -1 & \cdots & n-1
\end{array}
\right) \\
&  =\underbrace{\left(
\begin{array}
[c]{cccc}%
n & 0 & \cdots & 0\\
0 & n & \cdots & 0\\
\vdots & \vdots & \ddots & \vdots\\
0 & 0 & \cdots & n
\end{array}
\right)  }_{=A}+\underbrace{\left(
\begin{array}
[c]{c}%
-1\\
-1\\
\vdots\\
-1
\end{array}
\right)  }_{=u}\underbrace{\left(
\begin{array}
[c]{cccc}%
1 & 1 & \cdots & 1
\end{array}
\right)  }_{=v}.
\end{align*}

\item Here is an approach that is heavier on linear algebra (specifically,
eigenvectors and eigenvalues\footnote{See \cite[Chapter 4]{Treil17} for a
refresher.}):

Let $\left(  e_{1},e_{2},\ldots,e_{n-1}\right)  $ be the standard basis of the
$\mathbb{R}$-vector space $\mathbb{R}^{n-1}$ (so that $e_{i}$ is the column
vector with its $i$-th coordinate equal to $1$ and all its other coordinates
equal to $0$). Then, we can find the following $n-1$ eigenvectors of our
$\left(  n-1\right)  \times\left(  n-1\right)  $-matrix $\left(
\begin{array}
[c]{cccc}%
n-1 & -1 & \cdots & -1\\
-1 & n-1 & \cdots & -1\\
\vdots & \vdots & \ddots & \vdots\\
-1 & -1 & \cdots & n-1
\end{array}
\right)  $:

\begin{itemize}
\item the $n-2$ eigenvectors $e_{1}-e_{i}$ for all $i\in\left\{
2,3,\ldots,n-1\right\}  $, each of them with eigenvalue $n$ (check this!);

\item the eigenvector $e_{1}+e_{2}+\cdots+e_{n-1}$ with eigenvalue $1$ (check this!).
\end{itemize}

Since these $n-1$ eigenvectors are linearly independent (check this!), they
form a basis of $\mathbb{R}^{n-1}$. Hence, our matrix is similar to the
diagonal matrix with diagonal entries $\underbrace{n,n,\ldots,n}_{n-2\text{
times}},1$ (by \cite[Chapter 4, Theorem 2.1]{Treil17}), and therefore has
determinant $\underbrace{nn\cdots n}_{n-2\text{ times}}1=n^{n-2}$.
\end{itemize}

There are other ways as well. Either way, the result we obtain is $n^{n-2}$.
Thus, we have proved (relying on the Matrix-Tree Theorem, which we haven't yet
proved) the following result:

\begin{theorem}
[Cayley's formula]\label{thm.cayley.cayley}Let $n$ be a positive integer.
Then, the $\#$ of spanning trees of the complete graph $K_{n}$ is $n^{n-2}$.
\end{theorem}

In other words:

\begin{corollary}
\label{cor.cayley.n-trees}Let $n$ be a positive integer. Then, the $\#$ of
simple graphs with vertex set $\left\{  1,2,\ldots,n\right\}  $ that are trees
is $n^{n-2}$.
\end{corollary}

\begin{proof}
This is just Theorem \ref{thm.cayley.cayley}, since the simple graphs with
vertex set $\left\{  1,2,\ldots,n\right\}  $ that are trees are precisely the
spanning trees of $K_{n}$.
\end{proof}

Theorem \ref{thm.cayley.cayley} is widely known as \textbf{Cayley's formula},
despite its earlier discovery by Borchardt in 1860 (\cite{Borcha60}, giving a
proof quite similar to ours). There are many other ways to prove it. I can
particularly recommend the two combinatorial proofs given in \cite[\S 2.4 and
\S 2.5]{Galvin}, as well as Joyal's proof sketched in \cite{Leinst19}. Most
textbooks on enumerative combinatorics give one proof or another;
\cite[Appendix to Chapter 9]{Stanley-AC} gives three. Cayley's formula also
appears in Aigner's and Ziegler's best-of compilation of mathematical proofs
\cite[Chapter 33]{AigZie} with four different proofs. Note that some of the
sources use a matrix-tree theorem for \textbf{undirected} graphs; this is a
particular case of our matrix-tree theorem.\footnote{One more \textbf{remark:}
In Corollary \ref{cor.cayley.n-trees}, we have counted the trees with $n$
vertices (i.e., simple graphs with vertex set $\left\{  1,2,\ldots,n\right\}
$ that are trees). It sounds equally natural to count the \textquotedblleft
unlabelled trees with $n$ vertices\textquotedblright, i.e., the equivalence
classes of such trees up to isomorphism. Unfortunately, this is one of those
\textquotedblleft messy numbers\textquotedblright\ with no good expression:
the best formula known is recursive. There is also an asymptotic formula
(\textquotedblleft Otter's formula\textquotedblright, \cite{Otter48}): the
number of equivalence classes of $n$-vertex trees (up to isomorphism) is%
\[
\approx\beta\dfrac{\alpha^{n}}{n^{5/2}}\ \ \ \ \ \ \ \ \ \ \text{with }%
\alpha\approx2.955\text{ and }\beta\approx0.5349.
\]
}

However, in order to complete our proof, we still need to prove the
Matrix-Tree Theorem.

\subsubsection{Preparations for the proof}

In order to prepare for the proof of the Matrix-Tree Theorem, we state a
simple lemma (yet another criterion for a digraph to be an arborescence):

\begin{lemma}
\label{lem.MTT.0}Let $D=\left(  V,A,\psi\right)  $ be a multidigraph. Let $r$
be a vertex of $D$. Assume that $D$ has no cycles. Assume moreover that $D$
has no arcs with source $r$. Assume furthermore that each vertex $v\in
V\setminus\left\{  r\right\}  $ has outdegree $1$. Then, the digraph $D$ is an
arborescence rooted to $r$.
\end{lemma}

This lemma is precisely Exercise \ref{exe.4.6} \textbf{(b)}, at least after
reversing all arcs. But let us give a self-contained proof here:

\begin{proof}
[Proof of Lemma \ref{lem.MTT.0}.]Let $u$ be any vertex of $D$. Clearly, the
digraph $D$ has at least one path that starts at $u$ (for example, the trivial
path $\left(  u\right)  $). Among all such paths, let $\mathbf{p}=\left(
v_{0},a_{1},v_{1},a_{2},v_{2},\ldots,a_{k},v_{k}\right)  $ be a
\textbf{longest} one.\footnote{Such a longest path exists, since $D$ has at
least one path that starts at $u$, and since a path of $D$ cannot have length
larger than $\left\vert V\right\vert -1$.} Thus, $v_{0}=u$.

We shall show that $v_{k}=r$. Indeed, assume the contrary. Thus, $v_{k}\neq
r$, so that $v_{k}\in V\setminus\left\{  r\right\}  $. Hence, the vertex
$v_{k}$ has outdegree $1$ (since we assumed that each vertex $v\in
V\setminus\left\{  r\right\}  $ has outdegree $1$). Thus, there exists an arc
$b$ of $D$ that has source $v_{k}$. Consider this arc $b$, and let $w$ be its
target. Thus, appending the arc $b$ and the vertex $w$ to the end of the path
$\mathbf{p}$, we obtain a walk%
\[
\mathbf{w}=\left(  v_{0},a_{1},v_{1},a_{2},v_{2},\ldots,a_{k},v_{k}%
,b,w\right)
\]
of $D$ that starts at $u$ (since $v_{0}=u$). Proposition
\ref{prop.mdg.cyc.btf-walk-cyc} shows that this walk $\mathbf{w}$ either is a
path or contains a cycle. Hence, $\mathbf{w}$ is a path (since $D$ has no
cycles). Thus, $\mathbf{w}$ is a path of $D$ that starts at $u$. Since
$\mathbf{w}$ is longer than $\mathbf{p}$ (namely, longer by $1$), this shows
that $\mathbf{p}$ is not a longest path of $D$ that starts at $u$. But this
contradicts the very definition of $\mathbf{p}$.

This contradiction shows that our assumption was false. Hence, $v_{k}=r$.
Thus, $\mathbf{p}$ is a path from $u$ to $r$ (since $v_{0}=u$ and $v_{k}=r$).
Therefore, the digraph $D$ has a path from $u$ to $r$ (namely, $\mathbf{p}$).

Forget that we fixed $u$. We thus have shown that for each vertex $u$ of $D$,
the digraph $D$ has a path from $u$ to $r$. In other words, $r$ is a to-root
of $D$. Furthermore, we have $\deg^{+}r=0$ (since $D$ has no arcs with source
$r$), and each $v\in V\setminus\left\{  r\right\}  $ satisfies $\deg^{+}v=1$
(since we have assumed that each vertex $v\in V\setminus\left\{  r\right\}  $
has outdegree $1$). In other words, the digraph $D$ satisfies Statement A'6
from the dual arborescence equivalence theorem (Theorem
\ref{thm.arbor.eq-A-dual}). Therefore, it satisfies Statement A'1 from that
theorem as well (since all six statements A'1, A'2, $\ldots$, A'6 are
equivalent). In other words, $D$ is an arborescence rooted to $r$. This proves
Lemma \ref{lem.MTT.0}.
\end{proof}

\subsubsection{The Matrix-Tree Theorem: proof}

We shall now prove the Matrix-Tree Theorem (Theorem \ref{thm.MTT.MTT}), guided
by the following battle plan:

\begin{enumerate}
\item First, we will prove it in the case when each vertex $v\in
V\setminus\left\{  r\right\}  $ has outdegree $1$. In this case, after
removing all arcs with source $r$ from $D$ (these arcs do not matter, since
neither the submatrix $L_{\sim r,\sim r}$ nor the spanning arborescences
rooted to $r$ depend on them), we have essentially two options (subcases):
either $D$ is itself an arborescence or $D$ has a cycle.

\item Then, we will prove the matrix-tree theorem in the slightly more general
case when each $v\in V\setminus\left\{  r\right\}  $ has outdegree $\leq1$.
This is easy, since a vertex $v\in V\setminus\left\{  r\right\}  $ having
outdegree $0$ trivializes the theorem.

\item Finally, we will prove the theorem in the general case. This is done by
strong induction on the number of arcs of $D$. Every time you have a vertex
$v\in V\setminus\left\{  r\right\}  $ with outdegree $>1$, you can pick such a
vertex and color the outgoing arcs from it red and blue in such a way that
each color is used at least once. Then, you can consider the subdigraph of $D$
obtained by removing all blue arcs (call it $D^{\operatorname*{red}}$) and the
subdigraph of $D$ obtained by removing all red arcs (call it
$D^{\operatorname*{blue}}$). You can then apply the induction hypothesis to
$D^{\operatorname*{red}}$ and to $D^{\operatorname*{blue}}$ (since each of
these two subdigraphs has fewer arcs than $D$), and add the results together.
The good news is that both the \# of spanning arborescences rooted to $r$ and
the determinant $\det\left(  L_{\sim r,\sim r}\right)  $ \textquotedblleft
behave additively\textquotedblright\ (we will soon see what this means).
\end{enumerate}

So let us begin with Step 1. We first study a very special case:

\begin{lemma}
\label{lem.MTT.1}Let $D=\left(  V,A,\psi\right)  $ be a multidigraph. Let $r$
be a vertex of $D$. Assume that $D$ has no cycles. Assume moreover that $D$
has no arcs with source $r$. Assume furthermore that each vertex $v\in
V\setminus\left\{  r\right\}  $ has outdegree $1$. Then:

\begin{enumerate}
\item[\textbf{(a)}] The digraph $D$ has a unique spanning arborescence rooted
to $r$.

\item[\textbf{(b)}] Assume that $V=\left\{  1,2,\ldots,n\right\}  $ for some
$n\in\mathbb{N}$. Let $L$ be the Laplacian of $D$. Then, $\det\left(  L_{\sim
r,\sim r}\right)  =1$.
\end{enumerate}
\end{lemma}

\begin{proof}
\textbf{(a)} Lemma \ref{lem.MTT.0} shows that the digraph $D$ itself is an
arborescence rooted to $r$.

As a consequence, $D$ itself is a spanning arborescence of $D$ rooted to $r$.

Therefore, $\left\vert A\right\vert =\left\vert V\right\vert -1$ (by Statement
A'2 in the Dual Arborescence Equivalence Theorem (Theorem
\ref{thm.arbor.eq-A-dual})\footnote{or by the fact that $\left\vert
A\right\vert $ is the sum of the outdegrees of all vertices of $D$}). Hence,
$D$ has no spanning arborescences other than itself (because the condition
$\left\vert A\right\vert =\left\vert V\right\vert -1$ would get destroyed as
soon as we remove an arc). So the only spanning arborescence of $D$ rooted to
$r$ is $D$ itself. This proves Lemma \ref{lem.MTT.1} \textbf{(a)}. \medskip

\textbf{(b)} We WLOG assume that $r=n$ (otherwise, we can rename the vertices
$r,r+1,r+2,\ldots,n$ as $n,r,r+1,\ldots,n-1$, so that the matrix $L_{\sim
r,\sim r}$ becomes $L_{\sim n,\sim n}$). Thus, $V\setminus\left\{  r\right\}
=V\setminus\left\{  n\right\}  =\left\{  1,2,\ldots,n-1\right\}  $ (since
$V=\left\{  1,2,\ldots,n\right\}  $). Hence, each vertex $i\in\left\{
1,2,\ldots,n-1\right\}  $ has outdegree $1$ (since each vertex $v\in
V\setminus\left\{  r\right\}  $ has outdegree $1$, and we can apply this to
$v=i$ because $i\in\left\{  1,2,\ldots,n-1\right\}  =V\setminus\left\{
r\right\}  $). In other words,
\begin{equation}
\deg^{+}i=1\ \ \ \ \ \ \ \ \ \ \text{for each }i\in\left\{  1,2,\ldots
,n-1\right\}  . \label{pf.lem.MTT.1.b.1}%
\end{equation}
Moreover, $D$ has no cycles and thus no loops. Thus, for each $i\in\left\{
1,2,\ldots,n-1\right\}  $, the number $a_{i,i}$ defined in Definition
\ref{def.MTT.L} is%
\begin{equation}
a_{i,i}=0 \label{pf.lem.MTT.1.b.2}%
\end{equation}
(since $a_{i,i}$ counts the loops with source $i$).

Let $D^{\prime}$ be the digraph $D$ with a loop added at each vertex -- i.e.,
the multidigraph obtained from $D$ by adding $n$ extra arcs $\ell_{1},\ell
_{2},\ldots,\ell_{n}$ and letting each arc $\ell_{i}$ have source $i$ and
target $i$.

Let $S_{n-1}$ denote the group of permutations of the set
\[
\left\{  1,2,\ldots,n-1\right\}  =\underbrace{\left\{  1,2,\ldots,n\right\}
}_{=V}\setminus\left\{  \underbrace{n}_{=r}\right\}  =V\setminus\left\{
r\right\}  .
\]

Now, from $r=n$, we have%
\begin{equation}
\det\left(  L_{\sim r,\sim r}\right)  =\det\left(  L_{\sim n,\sim n}\right)
=\sum_{\sigma\in S_{n-1}}\operatorname*{sign}\sigma\cdot\prod_{i=1}%
^{n-1}L_{i,\sigma\left(  i\right)  } \label{pf.lem.MTT.1.b.det}%
\end{equation}
(by the Leibniz formula for the determinant). We shall now study the addends
in the sum on the right hand side of this equality. Specifically, we will show
that the only addend whose product $\prod_{i=1}^{n-1}L_{i,\sigma\left(
i\right)  }$ is nonzero is the addend for $\sigma=\operatorname*{id}$.

Indeed, let $\sigma\in S_{n-1}$ be a permutation such that the product
$\prod_{i=1}^{n-1}L_{i,\sigma\left(  i\right)  }$ is nonzero. We shall prove
that $\sigma=\operatorname*{id}$.

Consider an arbitrary $v\in\left\{  1,2,\ldots,n-1\right\}  $. Then,
$L_{v,\sigma\left(  v\right)  }\neq0$ (because $L_{v,\sigma\left(  v\right)
}$ is a factor in the product $\prod_{i=1}^{n-1}L_{i,\sigma\left(  i\right)
}$, which is nonzero). However, the definition of $L$ yields $L_{v,\sigma
\left(  v\right)  }=\left(  \deg^{+}v\right)  \cdot\left[  v=\sigma\left(
v\right)  \right]  -a_{v,\sigma\left(  v\right)  }$. Thus,
\[
\left(  \deg^{+}v\right)  \cdot\left[  v=\sigma\left(  v\right)  \right]
-a_{v,\sigma\left(  v\right)  }=L_{v,\sigma\left(  v\right)  }\neq0.
\]
Hence, at least one of the numbers $\left[  v=\sigma\left(  v\right)  \right]
$ and $a_{v,\sigma\left(  v\right)  }$ is nonzero. In other words, we have
$v=\sigma\left(  v\right)  $ (this is what it means for $\left[
v=\sigma\left(  v\right)  \right]  $ to be nonzero) or the digraph $D$ has an
arc with source $v$ and target $\sigma\left(  v\right)  $ (because this is
what it means for $a_{v,\sigma\left(  v\right)  }$ to be nonzero). In either
case, the digraph $D^{\prime}$ has an arc with source $v$ and target
$\sigma\left(  v\right)  $ (because if $v=\sigma\left(  v\right)  $, then one
of the loops we added to $D$ does the trick). We can apply the same argument
to $\sigma\left(  v\right)  $ instead of $v$, and obtain an arc with source
$\sigma\left(  v\right)  $ and target $\sigma\left(  \sigma\left(  v\right)
\right)  $. Similarly, we obtain an arc with source $\sigma\left(
\sigma\left(  v\right)  \right)  $ and target $\sigma\left(  \sigma\left(
\sigma\left(  v\right)  \right)  \right)  $. We can continue this reasoning
indefinitely. By continuing it for $n$ steps, we obtain a walk%
\[
\left(  v,\ast,\sigma\left(  v\right)  ,\ast,\sigma^{2}\left(  v\right)
,\ast,\sigma^{3}\left(  v\right)  ,\ldots,\ast,\sigma^{n}\left(  v\right)
\right)
\]
in the digraph $D^{\prime}$, where each asterisk means an arc (we don't care
about what these arcs are, so we are not giving them names). This walk cannot
be a path (since it has $n+1$ vertices, but $D^{\prime}$ has only $n$
vertices); thus, it must contain a cycle (by Proposition
\ref{prop.mdg.cyc.btf-walk-cyc}). All arcs of this cycle must be loops
(because otherwise, we could remove the loops from this cycle and obtain a
cycle of $D$, but we know that $D$ has no cycles). In particular, its first
arc is a loop. Thus, our above walk $\left(  v,\ast,\sigma\left(  v\right)
,\ast,\sigma^{2}\left(  v\right)  ,\ast,\sigma^{3}\left(  v\right)
,\ldots,\ast,\sigma^{n}\left(  v\right)  \right)  $ contains a loop (since the
arcs of the cycle come from this walk). In other words, we have $\sigma
^{i}\left(  v\right)  =\sigma^{i+1}\left(  v\right)  $ for some $i\in\left\{
0,1,\ldots,n-1\right\}  $. Since $\sigma$ is injective, we can apply
$\sigma^{-i}$ to both sides of this equality, and conclude that $v=\sigma
\left(  v\right)  $. In other words, $\sigma\left(  v\right)  =v$.

Forget that we fixed $v$. We thus have shown that $\sigma\left(  v\right)  =v$
for each $v\in\left\{  1,2,\ldots,n-1\right\}  $. In other words,
$\sigma=\operatorname*{id}$.

Forget that we fixed $\sigma$. We thus have proved that $\sigma
=\operatorname*{id}$ for each permutation $\sigma\in S_{n-1}$ for which the
product $\prod_{i=1}^{n-1}L_{i,\sigma\left(  i\right)  }$ is nonzero. In other
words, the only permutation $\sigma\in S_{n-1}$ for which the product
$\prod_{i=1}^{n-1}L_{i,\sigma\left(  i\right)  }$ is nonzero is the
permutation $\operatorname*{id}$.

Thus, the only nonzero addend on the right hand side of
(\ref{pf.lem.MTT.1.b.det}) is the addend corresponding to $\sigma
=\operatorname*{id}$. Hence, (\ref{pf.lem.MTT.1.b.det}) can be simplified as
follows:%
\[
\det\left(  L_{\sim n,\sim n}\right)  =\underbrace{\operatorname*{sign}\left(
\operatorname*{id}\right)  }_{=1}\cdot\prod_{i=1}^{n-1}L_{i,\operatorname*{id}%
\left(  i\right)  }=\prod_{i=1}^{n-1}L_{i,\operatorname*{id}\left(  i\right)
}.
\]
Since each $i\in\left\{  1,2,\ldots,n-1\right\}  $ satisfies%
\begin{align*}
L_{i,\operatorname*{id}\left(  i\right)  }  &  =L_{i,i}=\underbrace{\left(
\deg^{+}i\right)  }_{\substack{=1\\\text{(by (\ref{pf.lem.MTT.1.b.1}))}}%
}\cdot\underbrace{\left[  i=i\right]  }_{=1}-\underbrace{a_{i,i}%
}_{\substack{=0\\\text{(by (\ref{pf.lem.MTT.1.b.2}))}}%
}\ \ \ \ \ \ \ \ \ \ \left(  \text{by the definition of }L\right) \\
&  =1\cdot1-0=1,
\end{align*}
this can be simplified to $\det\left(  L_{\sim n,\sim n}\right)  =\prod
_{i=1}^{n-1}1=1$. This proves Lemma \ref{lem.MTT.1} \textbf{(b)}.
\end{proof}

Next, we drop the \textquotedblleft no cycles\textquotedblright\ condition:

\begin{lemma}
\label{lem.MTT.2}Let $D=\left(  V,A,\psi\right)  $ be a multidigraph. Let $r$
be a vertex of $D$. Assume that each vertex $v\in V\setminus\left\{
r\right\}  $ has outdegree $1$. Then, the MTT holds for these $D$ and $r$.
(Here and in the following, \textquotedblleft\textbf{MTT}\textquotedblright%
\ is short for \textquotedblleft Matrix-Tree Theorem\textquotedblright, i.e.,
for Theorem \ref{thm.MTT.MTT}.)
\end{lemma}

\begin{proof}
First of all, we note that an arc with source $r$ cannot appear in any
spanning arborescence of $D$ rooted to $r$ (since any such arborescence
satisfies $\deg^{+}r=0$, according to Statement A'6 in the Dual Arborescence
Equivalence Theorem (Theorem \ref{thm.arbor.eq-A-dual})). Furthermore, the
arcs with source $r$ do not affect the matrix $L_{\sim r,\sim r}$, since they
only appear in the $r$-th row of the matrix $L$ (but this $r$-th row is
removed in $L_{\sim r,\sim r}$).

Hence, any arc with source $r$ can be removed from $D$ without disturbing
anything we currently care about. Thus, we WLOG assume that $D$ has no arcs
with source $r$ (else, we can just remove them from $D$).

We WLOG assume that $r=n$ (otherwise, we can rename the vertices
$r,r+1,r+2,\ldots,n$ as $n,r,r+1,\ldots,n-1$, so that the matrix $L_{\sim
r,\sim r}$ becomes $L_{\sim n,\sim n}$).

We are in one of the following two cases:

\textit{Case 1:} The digraph $D$ has a cycle.

\textit{Case 2:} The digraph $D$ has no cycles.

Consider Case 1. In this case, $D$ has a cycle $\mathbf{v}=\left(  v_{1}%
,\ast,v_{2},\ast,\ldots,\ast,v_{m}\right)  $ (where we again are putting
asterisks in place of the arcs). This cycle cannot contain $r$ (since $D$ has
no arcs with source $r$). Thus, all its vertices $v_{1},v_{2},\ldots,v_{m}$
belong to $V\setminus\left\{  r\right\}  $. Hence, for each $i\in\left\{
1,2,\ldots,m-1\right\}  $, the vertex $v_{i}$ has outdegree $1$ (since we
assumed that each vertex $v\in V\setminus\left\{  r\right\}  $ has outdegree
$1$). Consequently, for each $i\in\left\{  1,2,\ldots,m-1\right\}  $, the only
arc of $D$ that has source $v_{i}$ is the arc that follows $v_{i}$ on the
cycle $\mathbf{v}$. Therefore, in the matrix $L$, the $v_{i}$-th row has a $1$
in the $v_{i}$-th position (because $\deg^{+}\left(  v_{i}\right)  =1$), a
$-1$ in the $v_{i+1}$-th position (since the arc that follows $v_{i}$ on the
cycle $\mathbf{v}$ has source $v_{i}$ and target $v_{i+1}$), and $0$s in all
other positions. Since $r=n$, the same must then be true for the matrix
$L_{\sim r,\sim r}$: That is, the $v_{i}$-th row of the matrix $L_{\sim r,\sim
r}$ has a $1$ in the $v_{i}$-th position, a $-1$ in the $v_{i+1}$-th position,
and $0$s in all other positions. Thus, the sum of the $v_{1}$-th, $v_{2}$-th,
$\ldots$, $v_{m-1}$-th rows of $L_{\sim r,\sim r}$ is the zero vector (since
the $1$s and the $-1$s just cancel out)\footnote{Namely, the $-1$ in the
$v_{i+1}$-th position of the $v_{i}$-th row gets cancelled by the $1$ in the
$v_{i+1}$-th position of the $v_{i+1}$-th row. (We are using the fact that
$v_{m}=v_{1}$ here.)}.\footnote{Let me illustrate this on a representative
example: Assume that the numbers $v_{1},v_{2},\ldots,v_{m-1},v_{m}$ are
$1,2,\ldots,m-1,1$ (respectively). Then, the first $m-1$ rows of $L$ look as
follows:%
\[%
\begin{array}
[c]{cccccccccc}%
1 & -1 &  &  &  &  &  &  &  & \\
& 1 & -1 &  &  &  &  &  &  & \\
&  & 1 & -1 &  &  &  &  &  & \\
&  &  & \ddots & \ddots &  &  &  &  & \\
&  &  &  & 1 & -1 &  &  &  & \\
-1 &  &  &  &  & 1 &  &  &  &
\end{array}
\]
(where all the missing entries are zeroes). Thus, the sum of these $m-1$ rows
is the zero vector. The same is therefore true of the matrix $L_{\sim r, \sim
r}$ (since the first $m-1$ rows of the latter matrix are just the first $m-1$
rows of $L$, with their $r$-th entries removed).
\par
The general case is essentially the same as this example; the only difference
is that the relevant rows are in other positions.}

So we have found a nonempty set of rows of $L_{\sim r,\sim r}$ whose sum is
the zero vector. This yields that the matrix $L_{\sim r,\sim r}$ is singular
(by basic properties of determinants\footnote{Specifically, we are using the
following fact: \textquotedblleft Let $M$ be a square matrix. If there is a
certain nonempty set of rows of $M$ whose sum is the zero vector, then the
matrix $M$ is singular.\textquotedblright.
\par
To prove this fact, we let $S$ be this nonempty set. Choose one row from this
set, and call it the \textbf{chosen row}. Now, add all the other rows from $S$
to this one chosen row. This operation does not change the determinant of $M$
(since the determinant of a matrix is unchanged when we add one row to
another), but the resulting matrix has a zero row (namely, the chosen row) and
thus has determinant $0$. Hence, the original matrix $M$ must have had
determinant $0$ as well. In other words, $M$ was singular, qed.}), so its
determinant is $\det\left(  L_{\sim r,\sim r}\right)  =0$. On the other hand,
the digraph $D$ has no spanning arborescence (because, in order to get a
spanning arborescence of $D$, we would have to remove at least one arc of our
cycle $\mathbf{v}$ (since an arborescence cannot have a cycle); but then, the
source of this arc would have outdegree $0$, and thus we could no longer find
a path from this source to $r$, so we would not obtain a spanning
arborescence). In other words,%
\[
\left(  \text{\# of spanning arborescences of }D\text{ rooted to }r\right)
=0.
\]
Comparing this with $\det\left(  L_{\sim r,\sim r}\right)  =0$, we conclude
that the MTT holds in this case (since it claims that $0=0$). Thus, Case 1 is done.

Next, we consider Case 2. In this case, $D$ has no cycles. Then, $\det\left(
L_{\sim r,\sim r}\right)  =1$ (by Lemma \ref{lem.MTT.1} \textbf{(b)}) and
\[
\left(  \text{\# of spanning arborescences of }D\text{ rooted to }r\right)
=1\ \ \ \ \ \ \ \ \ \ \left(  \text{by Lemma \ref{lem.MTT.1} \textbf{(a)}%
}\right)  .
\]
Thus, the MTT boils down to $1=1$, which is again true.

So Lemma \ref{lem.MTT.2} is proved.
\end{proof}

Next, we venture into a mildly greater generality, corresponding to Step 2 in
our above battle plan:

\begin{lemma}
\label{lem.MTT.3}Let $D=\left(  V,A,\psi\right)  $ be a multidigraph. Let $r$
be a vertex of $D$. Assume that each vertex $v\in V\setminus\left\{
r\right\}  $ has outdegree $\leq1$. Then, the MTT (= Matrix-Tree Theorem)
holds for these $D$ and $r$.
\end{lemma}

We give two proofs for this lemma: a short one using Lemma \ref{lem.MTT.2},
and a longer one that avoids Lemma \ref{lem.MTT.2}.

\begin{proof}
[First proof of Lemma \ref{lem.MTT.3}.]If each vertex $v\in V\setminus\left\{
r\right\}  $ has outdegree $1$, then this is true by Lemma \ref{lem.MTT.2}.

Thus, we WLOG assume that this is not the case. Hence, some vertex $v\in
V\setminus\left\{  r\right\}  $ has outdegree $\neq1$. Consider this $v$. The
outdegree of $v$ is $\neq1$, but also $\leq1$ (by the hypothesis of the
lemma). Hence, this outdegree must be $0$. That is, there is no arc with
source $v$.

WLOG assume that $r=n$ (otherwise, rename the vertices $r,r+1,r+2,\ldots,n$ as
$n,r,r+1,\ldots,n-1$, so that the matrix $L_{\sim r,\sim r}$ becomes $L_{\sim
n,\sim n}$).

We have $v\neq r$. Hence, the digraph $D$ has no path from $v$ to $r$ (since
any such path would include an arc with source $v$, but there is no arc with
source $v$).

Therefore, $D$ has no spanning arborescence rooted to $r$ (because any such
spanning arborescence would have to have a path from $v$ to $r$). In other
words,%
\[
\left(  \text{\# of spanning arborescences of }D\text{ rooted to }r\right)
=0.
\]
Also, $\det\left(  L_{\sim r,\sim r}\right)  =0$ (since the $v$-th row of the
matrix $L_{\sim r,\sim r}$ is $0$ (because there is no arc with source $v$)).
So the MTT boils down to $0=0$ again, and thus Lemma \ref{lem.MTT.3} is proved.
\end{proof}

\begin{fineprint}

\begin{proof}
[Second proof of Lemma \ref{lem.MTT.3}.]Just as in our above proof of Lemma
\ref{lem.MTT.2}, we can WLOG assume

\begin{itemize}
\item that $D$ has no arcs with source $r$, and

\item that $r=n$.
\end{itemize}

\noindent Let us make these two assumptions. We distinguish between two cases:

\textit{Case 1:} The vertex $r$ is a to-root of $D$.

\textit{Case 2:} The vertex $r$ is not a to-root of $D$.

Consider Case 1. In this case, $r$ is a to-root of $D$. Each vertex $v\in
V\setminus\left\{  r\right\}  $ satisfies $\deg^{+}v\leq1$ (since we have
assumed that each vertex $v\in V\setminus\left\{  r\right\}  $ has outdegree
$\leq1$) and therefore $\deg^{+}v=1$ (because otherwise, we would have
$\deg^{+}v=0$, and thus there would be no path from $v$ to $r$, since any such
path would have to start with an arc with source $v$; but this would
contradict the fact that $r$ is a to-root of $D$). Furthermore, $\deg^{+}r=0$
(since $D$ has no arcs with source $r$). Hence, the statement A'6 of Theorem
\ref{thm.arbor.eq-A-dual} holds. Thus, the implication A'6$\Longrightarrow$A'3
of Theorem \ref{thm.arbor.eq-A-dual} shows that the multigraph
$D^{\operatorname*{und}}$ is a tree. Hence, $D^{\operatorname*{und}}$ has no
cycles. Thus, the directed graph $D$ has no cycles either. Therefore, Lemma
\ref{lem.MTT.1} \textbf{(a)} shows that $D$ has a unique spanning arborescence
rooted to $r$. In other words,%
\[
\left(  \text{\# of spanning arborescences of }D\text{ rooted to }r\right)
=1.
\]
Comparing this with $\det\left(  L_{\sim r,\sim r}\right)  =1$ (which follows
from Lemma \ref{lem.MTT.1} \textbf{(b)}), we obtain%
\[
\left(  \#\text{ of spanning arborescences of }D\text{ rooted to }r\right)
=\det\left(  L_{\sim r,\sim r}\right)  .
\]
Thus, the MTT is proved in Case 1.

Let us now consider Case 2. In this case, the vertex $r$ is not a to-root of
$D$. Hence, there exists a vertex $v\in V$ such that $D$ has no path from $v$
to $r$. In other words, the set%
\[
S:=\left\{  v\in V\ \mid\ \text{there exists no path from }v\text{ to
}r\right\}
\]
is nonempty. Consider this set $S$. Clearly, $r\notin S$ (since there is
obviously a path from $r$ to $r$). In other words, $n\notin S$ (since $r=n$).
Thus, $S\subseteq\left\{  1,2,\ldots,n\right\}  \setminus\left\{  n\right\}
=\left\{  1,2,\ldots,n-1\right\}  $.

Recall that $r=n$. Hence, if we permute the vertices in $\left\{
1,2,\ldots,n-1\right\}  $ (that is, we rename them as $\sigma\left(  1\right)
,\sigma\left(  2\right)  ,\ldots,\sigma\left(  n-1\right)  $ for some
permutation $\sigma\in S_{n-1}$), then the $\#$ of spanning arborescences of
$D$ rooted to $r$ remains unchanged, and the determinant $\det\left(  L_{\sim
r,\sim r}\right)  $ remains unchanged as well (since the matrix $L_{\sim
r,\sim r}$ undergoes a permutation of its rows and the same permutation of its
columns, which amounts to conjugating this matrix by a permutation matrix).
Thus, the claim of the MTT does not change if we permute the vertices in
$\left\{  1,2,\ldots,n-1\right\}  $.

Hence, we can WLOG assume that $S=\left\{  1,2,\ldots,k\right\}  $ for some
nonnegative integer $k\leq n-1$ (indeed, we can achieve this by permuting the
vertices in $\left\{  1,2,\ldots,n-1\right\}  $, since $S\subseteq\left\{
1,2,\ldots,n-1\right\}  $). Assume this. Then, the set $\left\{
1,2,\ldots,k\right\}  =S$ is nonempty; thus, $k>0$.

The following property of $S$ will be crucial:

\begin{statement}
\textit{Claim 0:} If the source of an arc $a$ of $D$ belongs to $S$, then the
target of $a$ also belongs to $S$.
\end{statement}

\begin{proof}
[Proof of Claim 0.]Let $a$ be an arc of $D$ whose source belongs to $S$. We
must show that the target of $a$ also belongs to $S$.

Let $s$ and $t$ be the source and the target of $a$. Then, $s\in S$ (since the
source of $a$ belongs to $S$). In other words, there exists no path from $s$
to $r$ (by the definition of $S$).

If there was a path $\mathbf{p}$ from $t$ to $r$, then there would also be a
walk from $s$ to $r$ (namely, the walk $\left(  s,a,\underbrace{t,\ldots
,r}_{\text{the path }\mathbf{p}}\right)  $), and therefore there would also be
a path from $s$ to $r$ (by Corollary \ref{cor.mdg.walk-thus-path}); but this
would contradict the previous sentence. Hence, there exists no path from $t$
to $r$. In other words, $t\in S$ (by the definition of $S$). In other words,
the target of $a$ belongs to $S$ (since this target is $t$). This proves Claim 0.
\end{proof}

Now consider the induced subdigraph $D\left[  S\right]  $ of $D$ on the set
$S=\left\{  1,2,\ldots,k\right\}  $. Let $L_{D\left[  S\right]  }\in
\mathbb{Z}^{k\times k}$ be the Laplacian of this subdigraph $D\left[
S\right]  $. We claim that the Laplacian $L$ of $D$ can be written in
block-matrix notation as follows:%
\begin{equation}
L=\left(
\begin{array}
[c]{cc}%
L_{D\left[  S\right]  } & 0_{k\times\left(  n-k\right)  }\\
P & Q
\end{array}
\right)  \label{pf.lem.MTT.3.Lblock}%
\end{equation}
for some $\left(  n-k\right)  \times k$-matrix $P$ and some $\left(
n-k\right)  \times\left(  n-k\right)  $-matrix $Q$ (where $0_{k\times\left(
n-k\right)  }$ denotes the zero matrix with $k$ rows and $n-k$ columns). To
prove this, we need to show the following two claims:

\begin{statement}
\textit{Claim 1:} The top-left $k\times k$ block of $L$ is $L_{D\left[
S\right]  }$.
\end{statement}

\begin{statement}
\textit{Claim 2:} The top-right $k\times\left(  n-k\right)  $ block of $L$ is
the zero matrix $0_{k\times\left(  n-k\right)  }$.
\end{statement}

\begin{proof}
[Proof of Claim 1.]We need to show that $L_{i,j}=\left(  L_{D\left[  S\right]
}\right)  _{i,j}$ for all $i,j\in\left\{  1,2,\ldots,k\right\}  $.

So let $i,j\in\left\{  1,2,\ldots,k\right\}  $ be arbitrary. Then, $i$ and $j$
are vertices of both digraphs $D$ and $D\left[  S\right]  $.

By the definition of the Laplacian, we have%
\begin{equation}
L_{i,j}=\left(  \deg^{+}i\right)  \cdot\left[  i=j\right]  -a_{i,j}
\label{pf.lem.MTT.3.c1.pf.1}%
\end{equation}
and%
\begin{equation}
\left(  L_{D\left[  S\right]  }\right)  _{i,j}=\left(  \deg^{+}i\right)
\cdot\left[  i=j\right]  -a_{i,j}, \label{pf.lem.MTT.3.c1.pf.2}%
\end{equation}
where $a_{i,j}$ denotes the \# of arcs with source $i$ and target $j$.
However, in theory, the outdegree $\deg^{+}i$ in (\ref{pf.lem.MTT.3.c1.pf.1})
must be distinguished from the outdegree $\deg^{+}i$ in
(\ref{pf.lem.MTT.3.c1.pf.2}): Indeed, the former is defined with respect to
the digraph $D$, while the latter is defined with respect to the subdigraph
$D\left[  S\right]  $. Thus, the former outdegree counts all the arcs of $D$
with source $i$, whereas the latter outdegree counts only those arcs whose
target belongs to $S$. However, in our case, this distinction is unnecessary:
Indeed, any arc of $D$ whose source is $i$ must also have its target belong to
$S$ (by Claim 0, since its source is $i\in\left\{  1,2,\ldots,k\right\}  =S$),
and thus must be an arc of the induced subdigraph $D\left[  S\right]  $ as
well (since both its source and its target belong to $S$). Therefore, the
outdegree $\deg^{+}i$ in (\ref{pf.lem.MTT.3.c1.pf.1}) and the outdegree
$\deg^{+}i$ in (\ref{pf.lem.MTT.3.c1.pf.2}) are equal. (Of course, the
$a_{i,j}$'s in both equations are also equal, since $D\left[  S\right]  $
contains all arcs from $i$ to $j$.) Thus, the right hand sides of the
equalities (\ref{pf.lem.MTT.3.c1.pf.1}) and (\ref{pf.lem.MTT.3.c1.pf.2}) are
equal. Therefore, their left hand sides are also equal. In other words,
$L_{i,j}=\left(  L_{D\left[  S\right]  }\right)  _{i,j}$. This proves Claim 1.
\end{proof}

\begin{proof}
[Proof of Claim 2.]We must prove that $L_{i,j}=0$ for all $i\in\left\{
1,2,\ldots,k\right\}  $ and $j\in\left\{  k+1,k+2,\ldots,n\right\}  $.

So let $i\in\left\{  1,2,\ldots,k\right\}  $ and $j\in\left\{  k+1,k+2,\ldots
,n\right\}  $ be arbitrary. Then, $i\in S$ (since $i\in\left\{  1,2,\ldots
,k\right\}  =S$) and $j\notin S$ (since $j\in\left\{  k+1,k+2,\ldots
,n\right\}  =\left\{  1,2,\ldots,n\right\}  \setminus\underbrace{\left\{
1,2,\ldots,k\right\}  }_{=S}=\left\{  1,2,\ldots,n\right\}  \setminus S$).
Thus, $i\neq j$, so that $\left[  i=j\right]  =0$.

If some arc $a$ of $D$ has source $i$, then the target of $a$ must belong to
$S$ (by Claim 0, since the source of $a$ is $i\in S$), and therefore the
target of $a$ cannot be $j$ (since $j\notin S$). In other words, $D$ has no
arc with source $i$ and target $j$.

Let $a_{i,j}$ be the \# of arcs of $D$ having source $i$ and target $j$. Then,
$a_{i,j}=0$ (since $D$ has no arc with source $i$ and target $j$). But the
definition of the Laplacian $L$ shows that%
\[
L_{i,j}=\left(  \deg^{+}i\right)  \cdot\underbrace{\left[  i=j\right]  }%
_{=0}-\underbrace{a_{i,j}}_{=0}=0.
\]
This proves Claim 2.
\end{proof}

Combining Claim 1 with Claim 2, we see that the $n\times n$-matrix $L$ can be
written in the form (\ref{pf.lem.MTT.3.Lblock}) for some $\left(  n-k\right)
\times k$-matrix $P$ and some $\left(  n-k\right)  \times\left(  n-k\right)
$-matrix $Q$. Thus, its submatrix $L_{\sim n,\sim n}$ (which is obtained from
$L$ by removing the last row and the last column) can be written as%
\[
L_{\sim n,\sim n}=\left(
\begin{array}
[c]{cc}%
L_{D\left[  S\right]  } & 0_{k\times\left(  n-k\right)  }\\
P & Q
\end{array}
\right)  _{\sim n,\sim n}=\left(
\begin{array}
[c]{cc}%
L_{D\left[  S\right]  } & 0_{k\times\left(  n-k-1\right)  }\\
P^{\prime} & Q^{\prime}%
\end{array}
\right)  ,
\]
where $P^{\prime}$ and $Q^{\prime}$ are certain submatrices of $P$ and $Q$
(namely, $P^{\prime}$ is $P$ without its last row, whereas $Q^{\prime}$ is $Q$
without its last row and its last column).\footnote{We have used $k\leq n-1$
here (which is the reason why we aren't cutting into the $L_{D\left[
S\right]  }$ part of $L$ when we are removing its $n$-th row and column).}
This is clearly a block-lower-triangular matrix with diagonal blocks
$L_{D\left[  S\right]  }$ and $Q^{\prime}$. But the determinant of a
block-lower-triangular matrix is known to equal the product of the
determinants of its diagonal blocks\footnote{This is a standard result about
determinants; see, e.g., \cite[Exercise 6.30]{detnotes}.}. Hence, we obtain%
\begin{equation}
\det\left(  L_{\sim n,\sim n}\right)  =\det\left(  L_{D\left[  S\right]
}\right)  \cdot\det\left(  Q^{\prime}\right)  . \label{pf.lem.MTT.3.det*det}%
\end{equation}

But $D\left[  S\right]  $ is a digraph with vertex set $S=\left\{
1,2,\ldots,k\right\}  $, and the number $k$ is a positive integer (since
$k>0$). Hence, Proposition \ref{prop.MTT.detL=0} (applied to $D\left[
S\right]  $ and $k$ instead of $D$ and $n$) yields that the Laplacian
$L_{D\left[  S\right]  }$ of $D\left[  S\right]  $ is singular; i.e., we have
$\det\left(  L_{D\left[  S\right]  }\right)  =0$. Thus,
(\ref{pf.lem.MTT.3.det*det}) becomes%
\[
\det\left(  L_{\sim n,\sim n}\right)  =\underbrace{\det\left(  L_{D\left[
S\right]  }\right)  }_{=0}\cdot\det\left(  Q^{\prime}\right)  =0.
\]
In other words, $\det\left(  L_{\sim r,\sim r}\right)  =0$ (since $r=n$).

On the other hand, $r$ is not a to-root of $D$. Therefore, $r$ is not a
to-root of any spanning subdigraph of $D$ (since a to-root of a spanning
subdigraph of $D$ must always be a to-root of $D$\ \ \ \ \footnote{This is
clear from the definition of a to-root: If $s$ is a to-root of a spanning
subdigraph of $D$, then for each vertex $v$ of $D$, there must be a path from
$v$ to $s$ in this subdigraph, and therefore also a path from $v$ to $s$ in
$D$.}). Hence, the digraph $D$ has no spanning arborescence rooted to $r$
(because $r$ would be a to-root of any such arborescence). In other words,
\[
\left(  \#\text{ of spanning arborescences of }D\text{ rooted to }r\right)
=0.
\]
Comparing this with $\det\left(  L_{\sim r,\sim r}\right)  =0$, we obtain%
\[
\left(  \#\text{ of spanning arborescences of }D\text{ rooted to }r\right)
=\det\left(  L_{\sim r,\sim r}\right)  .
\]
Thus, the MTT is proved in Case 2.

We have now proved the MTT in both Cases 1 and 2. So the MTT holds for our
$D$. This proves Lemma \ref{lem.MTT.3} again.
\end{proof}
\end{fineprint}

We are now ready to prove the MTT in the general case (Step 3 of our battle plan):

\begin{proof}
[Proof of Theorem \ref{thm.MTT.MTT}.]First, we introduce a notation:

\begin{statement}
Let $M$ and $N$ be two $n\times n$-matrices that agree in all but one row.
That is, there exists some $j\in\left\{  1,2,\ldots,n\right\}  $ such that for
each $i\neq j$, we have%
\[
\left(  \text{the }i\text{-th row of }M\right)  =\left(  \text{the }i\text{-th
row of }N\right)  .
\]
Then, we write $M\overset{j}{\equiv}N$, and we let $M\overset{j}{+}N$ be the
$n\times n$-matrix that is obtained from $M$ by adding the $j$-th row of $N$
to the $j$-th row of $M$ (while leaving all remaining rows unchanged).
\end{statement}

For example, if $M=\left(
\begin{array}
[c]{ccc}%
a & b & c\\
d & e & f\\
g & h & i
\end{array}
\right)  $ and $N=\left(
\begin{array}
[c]{ccc}%
a & b & c\\
d^{\prime} & e^{\prime} & f^{\prime}\\
g & h & i
\end{array}
\right)  $, then $M\overset{2}{\equiv}N$ and%
\[
M\overset{2}{+}N=\left(
\begin{array}
[c]{ccc}%
a & b & c\\
d+d^{\prime} & e+e^{\prime} & f+f^{\prime}\\
g & h & i
\end{array}
\right)  .
\]

A well-known property of determinants (the \textbf{multilinearity of the
determinant}) says that if $M$ and $N$ are two $n\times n$-matrices and
$j\in\left\{  1,2,\ldots,n\right\}  $ is a number such that
$M\overset{j}{\equiv}N$, then
\[
\det\left(  M\overset{j}{+}N\right)  =\det M+\det N.
\]

Now, let us prove the MTT. We proceed by strong induction on the \# of arcs of
$D$.

\textit{Induction step:} Let $m\in\mathbb{N}$. Assume (as the induction
hypothesis) that the MTT holds for all digraphs $D$ that have $<m$ arcs. We
must now prove it for our digraph $D$ with $m$ arcs.

WLOG assume that $r=n$ (otherwise, rename the vertices $r,r+1,r+2,\ldots,n$ as
$n,r,r+1,\ldots,n-1$).

If each vertex $v\in V\setminus\left\{  r\right\}  $ has outdegree $\leq1$,
then the MTT holds by Lemma \ref{lem.MTT.3}. Thus, we WLOG assume that some
vertex $v\in V\setminus\left\{  r\right\}  $ has outdegree $>1$. Pick such a
vertex $v$. We color each arc with source $v$ either red or blue, making sure
that at least one arc is red and at least one arc is blue. (We can do this,
since $v$ has outdegree $>1$.) All arcs that do not have source $v$ remain uncolored.

Now, let $D^{\operatorname*{red}}$ be the subdigraph obtained from $D$ by
removing all blue arcs. Then, $D^{\operatorname*{red}}$ has fewer arcs than
$D$. In other words, $D^{\operatorname*{red}}$ has $<m$ arcs. Hence, the
induction hypothesis yields that the MTT holds for $D^{\operatorname*{red}}$.
That is, we have%
\[
\left(  \text{\# of spanning arborescences of }D^{\operatorname*{red}}\text{
rooted to }r\right)  =\det\left(  L_{\sim r,\sim r}^{\operatorname*{red}%
}\right)  ,
\]
where $L^{\operatorname*{red}}$ means the Laplacian of $D^{\operatorname*{red}%
}$.

Likewise, let $D^{\operatorname*{blue}}$ be the subdigraph obtained from $D$
by removing all red arcs. Then, $D^{\operatorname*{blue}}$ has fewer arcs than
$D$. Hence, the induction hypothesis yields that the MTT holds for
$D^{\operatorname*{blue}}$. That is,%
\[
\left(  \text{\# of spanning arborescences of }D^{\operatorname*{blue}}\text{
rooted to }r\right)  =\det\left(  L_{\sim r,\sim r}^{\operatorname*{blue}%
}\right)  ,
\]
where $L^{\operatorname*{blue}}$ means the Laplacian of
$D^{\operatorname*{blue}}$.

\begin{example}
Let $D$ be the multidigraph%
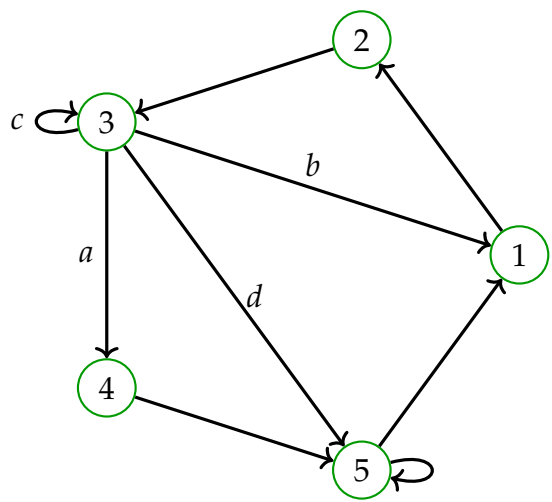
\[%
%
\]
with $r=1$. Its Laplacian is%
\[
L=\left(
\begin{array}
[c]{ccccc}%
1 & -1 & 0 & 0 & 0\\
0 & 1 & -1 & 0 & 0\\
-1 & 0 & 3 & -1 & -1\\
0 & 0 & 0 & 1 & -1\\
-1 & 0 & 0 & 0 & 1
\end{array}
\right)  .
\]
Let us pick $v=3$ (this is a vertex with outdegree $>1$), and let us color the
arcs $a$ and $c$ red and the arcs $b$ and $d$ blue (various other options are
possible). Then, $D^{\operatorname*{red}}$ and $D^{\operatorname*{blue}}$ look
as follows (along with their Laplacians $L^{\operatorname*{red}}$ and
$L^{\operatorname*{blue}}$): \newpage%
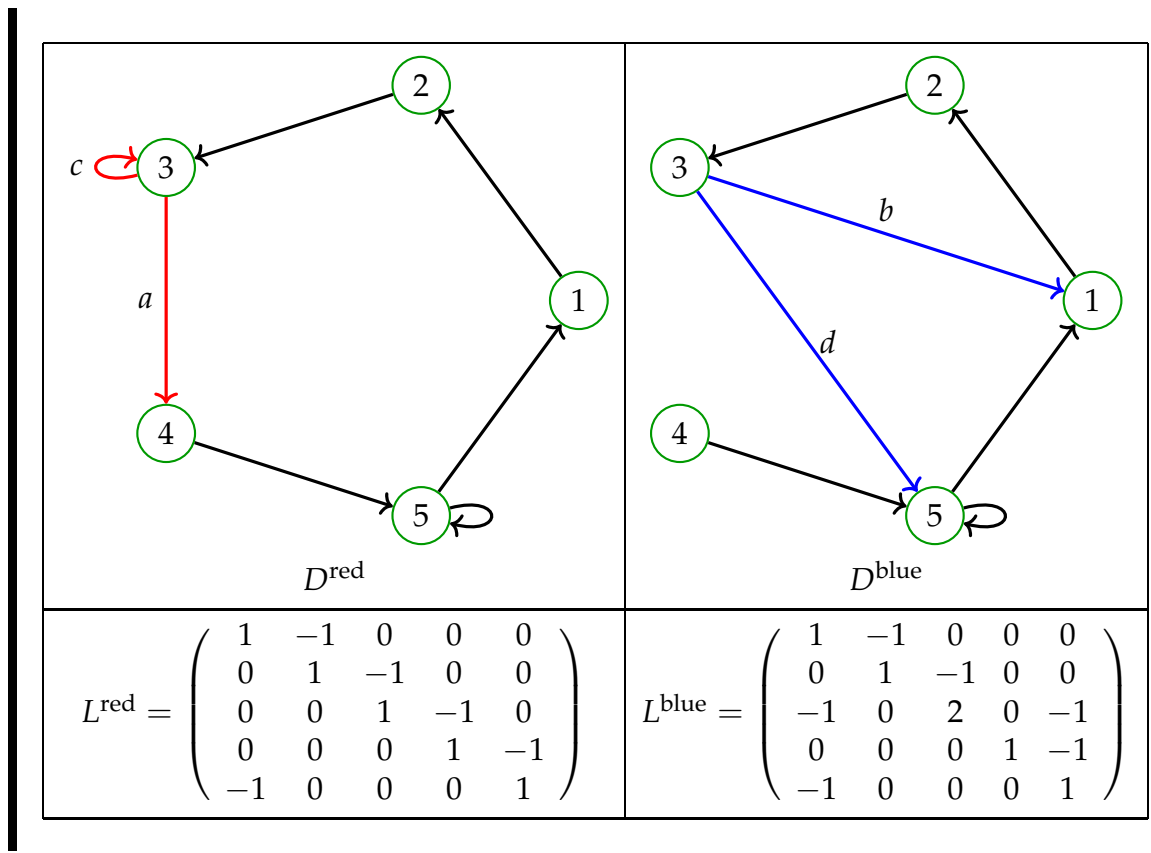
\[%
%
\]

\end{example}

Now, the digraphs $D$, $D^{\operatorname*{blue}}$ and $D^{\operatorname*{red}%
}$ differ only in the arcs with source $v$, and as far as the latter arcs are
concerned, the arcs of $D$ are divided between $D^{\operatorname*{blue}}$ and
$D^{\operatorname*{red}}$. Hence, by the definition of the Laplacian, we have%
\[
L^{\operatorname*{red}}\overset{v}{\equiv}L^{\operatorname*{blue}%
}\ \ \ \ \ \ \ \ \ \ \text{and}\ \ \ \ \ \ \ \ \ \ L^{\operatorname*{red}%
}\overset{v}{+}L^{\operatorname*{blue}}=L.
\]
Thus,%
\[
L_{\sim r,\sim r}^{\operatorname*{red}}\overset{v}{\equiv}L_{\sim r,\sim
r}^{\operatorname*{blue}}\ \ \ \ \ \ \ \ \ \ \text{and}%
\ \ \ \ \ \ \ \ \ \ L_{\sim r,\sim r}^{\operatorname*{red}}\overset{v}{+}%
L_{\sim r,\sim r}^{\operatorname*{blue}}=L_{\sim r,\sim r}%
\]
(here, we have used the fact that $r=n$ and $v\neq r$, so that when we remove
the $r$-th row and the $r$-th column of the matrix $L$, the $v$-th row remains
the $v$-th row). Hence,%
\[
\det\left(  \underbrace{L_{\sim r,\sim r}}_{=L_{\sim r,\sim r}%
^{\operatorname*{red}}\overset{v}{+}L_{\sim r,\sim r}^{\operatorname*{blue}}%
}\right)  =\det\left(  L_{\sim r,\sim r}^{\operatorname*{red}}\overset{v}{+}%
L_{\sim r,\sim r}^{\operatorname*{blue}}\right)  =\det\left(  L_{\sim r,\sim
r}^{\operatorname*{red}}\right)  +\det\left(  L_{\sim r,\sim r}%
^{\operatorname*{blue}}\right)
\]
(by the multilinearity of the determinant).

However, a similar equality holds for the \# of spanning arborescences:
namely, we have%
\begin{align*}
&  \left(  \text{\# of spanning arborescences of }D\text{ rooted to }r\right)
\\
&  =\left(  \text{\# of spanning arborescences of }D^{\operatorname*{red}%
}\text{ rooted to }r\right) \\
&  \ \ \ \ \ \ \ \ \ \ +\left(  \text{\# of spanning arborescences of
}D^{\operatorname*{blue}}\text{ rooted to }r\right)  .
\end{align*}
Here is why: Recall that an arborescence rooted to $r$ must satisfy $\deg
^{+}v=1$ (by Statement A'6 in the Dual Arborescence Equivalence Theorem
(Theorem \ref{thm.arbor.eq-A-dual}), since $v\in V\setminus\left\{  r\right\}
$). In other words, an arborescence rooted to $r$ must contain exactly one arc
with source $v$. In particular, a spanning arborescence of $D$ rooted to $r$
must contain either a red arc or a blue arc, but not both at the same time. In
the former case, it is a spanning arborescence of $D^{\operatorname*{red}}$;
in the latter, it is a spanning arborescence of $D^{\operatorname*{blue}}$.
Conversely, any spanning arborescence of $D^{\operatorname*{red}}$ or of
$D^{\operatorname*{blue}}$ rooted to $r$ is automatically a spanning
arborescence of $D$ rooted to $r$. Thus,
\begin{align*}
&  \left(  \text{\# of spanning arborescences of }D\text{ rooted to }r\right)
\\
&  =\underbrace{\left(  \text{\# of spanning arborescences of }%
D^{\operatorname*{red}}\text{ rooted to }r\right)  }_{\substack{=\det\left(
L_{\sim r,\sim r}^{\operatorname*{red}}\right)  \\\text{(as we saw above)}}}\\
&  \ \ \ \ \ \ \ \ \ \ +\underbrace{\left(  \text{\# of spanning arborescences
of }D^{\operatorname*{blue}}\text{ rooted to }r\right)  }_{\substack{=\det
\left(  L_{\sim r,\sim r}^{\operatorname*{blue}}\right)  \\\text{(as we saw
above)}}}\\
&  =\det\left(  L_{\sim r,\sim r}^{\operatorname*{red}}\right)  +\det\left(
L_{\sim r,\sim r}^{\operatorname*{blue}}\right)  =\det\left(  L_{\sim r,\sim
r}\right)
\end{align*}
(since we proved that $\det\left(  L_{\sim r,\sim r}\right)  =\det\left(
L_{\sim r,\sim r}^{\operatorname*{red}}\right)  +\det\left(  L_{\sim r,\sim
r}^{\operatorname*{blue}}\right)  $). That is, the MTT holds for our digraph
$D$ and its vertex $r$. This completes the induction step, and thus the MTT
(Theorem \ref{thm.MTT.MTT}) is proved.
\end{proof}

Our above proof of Theorem \ref{thm.MTT.MTT} has followed \cite[Theorem
10.4]{Stanley-AC}. Other proofs can be found across the literature, e.g., in
\cite[Theorem 7]{VanEhr51}, in \cite[Theorem 2.8]{Margol10}, in \cite[Theorem
1]{DeLeen19} and in \cite[Theorem 2.5.3]{Holzer22}. (Some of these sources
prove more general versions of the theorem. Confusingly, each source uses
different notations and works in a slightly different setup, although most of
them quickly reveal themselves to be equivalent upon some introspection.)

\subsubsection{Further exercises on the Laplacian}

\begin{exercise}
\label{exe.7.4}Let $G=\left(  V,E,\varphi\right)  $ be a multigraph. Let $L$
be the Laplacian of the digraph $G^{\operatorname*{bidir}}$. Prove that $L$ is
positive semidefinite. \medskip

[\textbf{Hint:} Write $L$ as $N^{T} N$, where $N$ or $N^{T}$ is some matrix
you have seen before.

Note that the statement is not true if we replace $G^{\operatorname*{bidir}}$
by an arbitrary digraph $D$.]
\end{exercise}

The following two exercises stand at the beginning of the theory of
\textbf{chip-firing} and related dynamical systems on a digraph (see
\cite{CorPer18}, \cite{Klivan19} and \cite{JoyMel17} for much more). While the
Laplacian is not mentioned in them directly, it is implicitly involved in the
definition of a \textquotedblleft donation\textquotedblright\ (how?).

\begin{exercise}
\label{exe.8.1}Let $D=\left(  V,A,\psi\right)  $ be a strongly connected multidigraph.

A \textbf{wealth distribution} on $D$ shall mean a family $\left(
k_{v}\right)  _{v\in V}$ of integers (one for each vertex $v\in V$). If
$k=\left(  k_{v}\right)  _{v\in V}$ is a wealth distribution, then we refer to
each value $k_{v}$ as the \textbf{wealth} of the vertex $v$, and we define the
\textbf{total wealth} of $k$ to be the sum $\sum_{v\in V}k_{v}$. We say that a
vertex $v$ is \textbf{in debt} in a given wealth distribution $k=\left(
k_{v}\right)  _{v\in V}$ if its wealth $k_{v}$ is negative.

For any vertices $v$ and $w$, we let $a_{v,w}$ denote the number of arcs that
have source $v$ and $w$.

A \textbf{donation} is an operation that transforms a wealth distribution as
follows: We choose a vertex $v$, and we decrease its wealth by its outdegree
$\deg^{+}v$, and then increase the wealth of each vertex $w\in V$ (including
$v$ itself) by $a_{v,w}$. (You can think of $v$ as donating a unit of wealth
for each arc that has source $v$. This unit flows to the target to this arc.
Note that a donation does not change the total wealth.)

Let $k$ be a wealth distribution on $D$ whose total wealth is larger than
$\left\vert A\right\vert -\left\vert V\right\vert $. Prove that by an
appropriately chosen finite sequence of donations, we can ensure that no
vertex is in debt. \medskip

[\textbf{Example:} For instance, consider the digraph
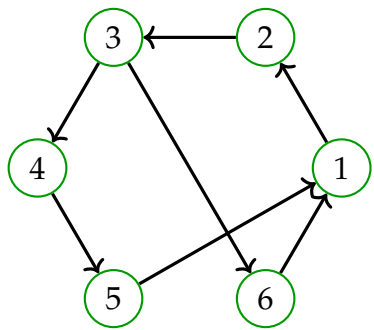
\[
\begin{tikzpicture}[scale=1]
\begin{scope}[every node/.style={circle,thick,draw=green!60!black}]
\node(A) at (0:2) {$1$};
\node(B) at (60:2) {$2$};
\node(C) at (120:2) {$3$};
\node(D) at (180:2) {$4$};
\node(E) at (240:2) {$5$};
\node(F) at (300:2) {$6$};
\end{scope}
\begin{scope}[every edge/.style={draw=black,very thick}]
\path[->] (A) edge (B) (B) edge (C) (D);
\path[->] (C) edge (D) edge (F) (D) edge (E) (E) edge (A);
\path[->] (F) edge (A);
\end{scope}
\end{tikzpicture}
\]
with wealth distribution $\left(  k_{1},k_{2},k_{3},k_{4},k_{5},k_{6}\right)
=\left(  -1,-1,1,2,0,1\right)  $. The vertices $1$ and $2$ are in debt here,
but it is possible to get all vertices out of debt by having the vertices
$4,5,6,1$ donate in some order (the order clearly does not matter for the result\footnotemark).

Note that vertices are allowed to donate multiple times (although in the above
example, this was unnecessary).] \medskip

[\textbf{Hint:} A donation will be called \textbf{safe} if its donor $v$ (that
is, the vertex chosen to lose wealth) satisfies $k_{v}\geq\deg^{+}v$, where
$k$ is the wealth distribution just before this donation. Start by showing
that if the total wealth is larger than $\left\vert A\right\vert -\left\vert
V\right\vert $, then at least one vertex $v$ has wealth $\geq\deg^{+}v$ (and
thus can make a safe donation). Next, show that for any given wealth
distribution $k$, there are only finitely many wealth distributions that can
be obtained from $k$ by a sequence of safe donations. Finally, for any vertex
$v$, find a rational quantity that increases every time that a donor distinct
from $v$ makes a donation. Conclude that in a sufficiently long sequence of
safe donations, every vertex must appear as a donor. But a donor of a safe
donation must be out of debt just before its safe donation, and will never go
back into debt.]
\end{exercise}

\footnotetext{Depending on the order, some vertices will go into debt in the
process, but this is okay as long as they ultimately end up debt-free.}

\begin{exercise}
\label{exe.8.2}We continue with the setting and terminology of Exercise
\ref{exe.8.1}.

A \textbf{clawback} is an operation that transforms a wealth distribution as
follows: We choose a vertex $v$, and we increase its wealth by its outdegree
$\deg^{+}v$, and then decrease the wealth of each vertex $w\in V$ (including
$v$ itself) by $a_{v,w}$. (Thus, a clawback is the inverse of a donation.)

Let $k$ be a wealth distribution on $D$ whose total wealth is larger than
$\left\vert A\right\vert -\left\vert V\right\vert $. Prove that by an
appropriately chosen finite sequence of clawbacks, we can ensure that no
vertex is in debt. \medskip

[\textbf{Remark:} Note that we are still assuming $D$ to be strongly
connected. Otherwise, the truth of the claim is not guaranteed. For instance,
for the digraph
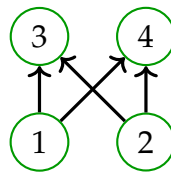
\[
\begin{tikzpicture}[scale=1.4]
\begin{scope}[every node/.style={circle,thick,draw=green!60!black}]
\node(1) at (0, 0) {$1$};
\node(2) at (1, 0) {$2$};
\node(3) at (0, 1) {$3$};
\node(4) at (1, 1) {$4$};
\end{scope}
\begin{scope}[every edge/.style={draw=black,very thick}, every loop/.style={}]
\path[->] (1) edge (3) edge (4) (2) edge (3) edge (4);
\end{scope}
\end{tikzpicture}
\]
with wealth distribution $\left(  k_{1},k_{2},k_{3},k_{4}\right)  =\left(
0,0,-1,2\right)  $, no sequence of donations and clawbacks will result in
every vertex being out of debt (since the wealth difference $k_{4}-k_{3}$ is
preserved under any donation or clawback, but this difference is too large to
come from a debt-free distribution with total weight $1$). ] \medskip

[\textbf{Hint:} Show that any donation is equivalent to an appropriately
chosen composition of clawbacks. Something we know about the Laplacian may
come useful here.]
\end{exercise}

\subsubsection{Application: Counting Eulerian circuits of $K_{n}%
^{\operatorname*{bidir}}$}

Here is one more consequence of the MTT:

\begin{proposition}
Let $n$ be a positive integer. Pick any arc $a$ of the multidigraph
$K_{n}^{\operatorname*{bidir}}$. Then, the \# of Eulerian circuits of
$K_{n}^{\operatorname*{bidir}}$ whose first arc is $a$ is $n^{n-2}\cdot\left(
n-2\right)  !^{n}$.
\end{proposition}

\begin{proof}
Let $r$ be the source of the arc $a$. The digraph $K_{n}%
^{\operatorname*{bidir}}$ is balanced, and each of its vertices has outdegree
$n-1$. By the BEST' theorem (Theorem \ref{thm.BEST.to}), we have%
\begin{align*}
&  \left(  \text{\# of Eulerian circuits of }K_{n}^{\operatorname*{bidir}%
}\text{ whose first arc is }a\right) \\
&  =\underbrace{\left(  \text{\# of spanning arborescences of }K_{n}%
^{\operatorname*{bidir}}\text{ rooted to }r\right)  }_{\substack{=n^{n-2}%
\\\text{(as we saw in Subsection \ref{subsec.MTT.Kn} in the case when
}r=1\text{,}\\\text{and can similarly prove for arbitrary }r\text{)}}%
}\cdot\prod_{u=1}^{n}\left(  \underbrace{\deg^{+}u}_{=n-1}-\,1\right)  !\\
&  =n^{n-2}\cdot\prod_{u=1}^{n}\left(  n-2\right)  !=n^{n-2}\cdot\left(
n-2\right)  !^{n},
\end{align*}
qed.
\end{proof}

In comparison, there is no good formula known for the \# of Eulerian circuits
of the undirected graph $K_{n}$. For $n$ even, this \# is $0$ of course (since
$K_{n}$ has vertices of odd degree in this case). For $n$ odd, the \# grows
very fast, but little else is known about it (see
\url{https://oeis.org/A135388} for some known values, and see Exercise
\ref{exe.7.3} for a divisibility property).

\begin{exercise}
\label{exe.7.2}Let $n$ be a positive integer. Let $N=\left\{  1,2,\ldots
,n\right\}  $. A map $f:N\rightarrow N$ is said to be $n$\textbf{-potent} if
each $i\in N$ satisfies $f^{n-1}\left(  i\right)  =n$. (As usual, $f^{k}$
denotes the $k$-fold composition $f\circ f\circ\cdots\circ f$.)

Prove that the \# of $n$-potent maps $f:N\rightarrow N$ is $n^{n-2}$. \medskip

[\textbf{Hint:} What do these $n$-potent maps have to do with trees?]
\end{exercise}

\begin{exercise}
\label{exe.7.3}Let $n=2m+1>2$ be an odd integer. Let $e$ be an edge of the
(undirected) complete graph $K_{n}$. Prove that the \# of Eulerian circuits of
$K_{n}$ that start with $e$ is a multiple of $\left(  m-1\right)  !^{n}$.
\medskip

[\textbf{Hint:} Argue that each Eulerian circuit of $K_{n}$ is an Eulerian
circuit of a unique balanced tournament. Here, a \textquotedblleft balanced
tournament\textquotedblright\ means a balanced digraph obtained from $K_{n}$
by orienting each edge.]
\end{exercise}

\subsection{The undirected Matrix-Tree Theorem}

\subsubsection{The theorem}

The Matrix-Tree Theorem becomes simpler if we apply it to a digraph of the
form $G^{\operatorname*{bidir}}$:

\begin{theorem}
[undirected Matrix-Tree Theorem]\label{thm.MTT.undir}Let $G=\left(
V,E,\varphi\right)  $ be a multigraph. Assume that $V=\left\{  1,2,\ldots
,n\right\}  $ for some positive integer $n$.

Let $L$ be the Laplacian of the digraph $G^{\operatorname*{bidir}}$.
Explicitly, this is the $n\times n$-matrix $L\in\mathbb{Z}^{n\times n}$ whose
entries are given by%
\[
L_{i,j}=\left(  \deg i\right)  \cdot\left[  i=j\right]  -a_{i,j},
\]
where $a_{i,j}$ is the $\#$ of edges of $G$ that have endpoints $i$ and $j$
(with loops counting twice). Then:

\begin{enumerate}
\item[\textbf{(a)}] For any vertex $r$ of $G$, we have%
\[
\left(  \text{\# of spanning trees of }G\right)  =\det\left(  L_{\sim r,\sim
r}\right)  .
\]

\item[\textbf{(b)}] Let $t$ be an indeterminate. Expand the determinant
$\det\left(  tI_{n}+L\right)  $ (here, $I_{n}$ denotes the $n\times n$
identity matrix) as a polynomial in $t$:%
\[
\det\left(  tI_{n}+L\right)  =c_{n}t^{n}+c_{n-1}t^{n-1}+\cdots+c_{1}%
t^{1}+c_{0}t^{0},
\]
where $c_{0},c_{1},\ldots,c_{n}$ are numbers. (Note that this is the
characteristic polynomial of $L$ up to substituting $-t$ for $t$ and
multiplying by a power of $-1$. Some of its coefficients are $c_{n}=1$ and
$c_{n-1}=\operatorname*{Tr}L$ and $c_{0}=\det L$.) Then,%
\[
\left(  \text{\# of spanning trees of }G\right)  =\dfrac{1}{n}c_{1}.
\]

\item[\textbf{(c)}] Let $\lambda_{1},\lambda_{2},\ldots,\lambda_{n}$ be the
eigenvalues of $L$, listed in such a way that $\lambda_{n}=0$ (we know that
$0$ is an eigenvalue of $L$, since $L$ is singular). Then,%
\[
\left(  \text{\# of spanning trees of }G\right)  =\dfrac{1}{n}\cdot\lambda
_{1}\lambda_{2}\cdots\lambda_{n-1}.
\]

\end{enumerate}
\end{theorem}

\begin{proof}
\textbf{(a)} Let $r$ be a vertex of $G$. Then, Proposition
\ref{prop.sparb.vs-sptree} \textbf{(b)} shows that there is a bijection%
\[
\left\{  \text{spanning arborescences of }G^{\operatorname*{bidir}}\text{
rooted to }r\right\}  \rightarrow\left\{  \text{spanning trees of }G\right\}
.
\]
Hence, by the bijection principle, we have%
\begin{align*}
&  \left(  \text{\# of spanning trees of }G\right) \\
&  =\left(  \text{\# of spanning arborescences of }G^{\operatorname*{bidir}%
}\text{ rooted to }r\right) \\
&  =\det\left(  L_{\sim r,\sim r}\right)  \ \ \ \ \ \ \ \ \ \ \left(  \text{by
the Matrix-Tree Theorem (Theorem \ref{thm.MTT.MTT})}\right)  .
\end{align*}
This proves Theorem \ref{thm.MTT.undir} \textbf{(a)}. \medskip

\textbf{(b)} We claim that
\begin{equation}
c_{1}=\sum_{r=1}^{n}\det\left(  L_{\sim r,\sim r}\right)  .
\label{pf.thm.MTT.undir.b.1}%
\end{equation}
Note that this is a purely linear-algebraic result, and has nothing to do with
the fact that $L$ is the Laplacian of a digraph; it holds just as well if $L$
is replaced by any square matrix.

Once (\ref{pf.thm.MTT.undir.b.1}) is proved, Theorem \ref{thm.MTT.undir}
\textbf{(b)} will easily follow, because (\ref{pf.thm.MTT.undir.b.1}) entails%
\begin{align*}
\dfrac{1}{n}c_{1}  &  =\dfrac{1}{n}\sum_{r=1}^{n}\underbrace{\det\left(
L_{\sim r,\sim r}\right)  }_{\substack{=\left(  \text{\# of spanning trees of
}G\right)  \\\text{(by Theorem \ref{thm.MTT.undir} \textbf{(a)})}}}=\dfrac
{1}{n}\underbrace{\sum_{r=1}^{n}\left(  \text{\# of spanning trees of
}G\right)  }_{=n\cdot\left(  \text{\# of spanning trees of }G\right)  }\\
&  =\dfrac{1}{n}\cdot n\cdot\left(  \text{\# of spanning trees of }G\right)
=\left(  \text{\# of spanning trees of }G\right)  .
\end{align*}
Thus, it remains to prove (\ref{pf.thm.MTT.undir.b.1}).

A rigorous proof of (\ref{pf.thm.MTT.undir.b.1}) can be found in
\cite[Proposition 6.4.29]{21s} or in
\url{https://math.stackexchange.com/a/3989575/} (both of these references
actually describe all coefficients $c_{0},c_{1},\ldots,c_{n}$ of the
polynomial $\det\left(  tI_{n}+L\right)  $, not just the $t^{1}$-coefficient
$c_{1}$). We shall merely outline the proof of (\ref{pf.thm.MTT.undir.b.1}) on
a convenient example. We want to compute $c_{1}$. In other words, we want to
compute the coefficient of $t^{1}$ in the polynomial $\det\left(
tI_{n}+L\right)  $ (since $c_{1}$ is defined to be this very coefficient). Let
us say that $n=4$, so that $L$ has the form%
\[
L=\left(
\begin{array}
[c]{cccc}%
a & b & c & d\\
a^{\prime} & b^{\prime} & c^{\prime} & d^{\prime}\\
a^{\prime\prime} & b^{\prime\prime} & c^{\prime\prime} & d^{\prime\prime}\\
a^{\prime\prime\prime} & b^{\prime\prime\prime} & c^{\prime\prime\prime} &
d^{\prime\prime\prime}%
\end{array}
\right)  .
\]
Thus,
\[
\det\left(  tI_{n}+L\right)  =\det\left(
\begin{array}
[c]{cccc}%
t+a & b & c & d\\
a^{\prime} & t+b^{\prime} & c^{\prime} & d^{\prime}\\
a^{\prime\prime} & b^{\prime\prime} & t+c^{\prime\prime} & d^{\prime\prime}\\
a^{\prime\prime\prime} & b^{\prime\prime\prime} & c^{\prime\prime\prime} &
t+d^{\prime\prime\prime}%
\end{array}
\right)  .
\]
Imagine expanding the right hand side (using the Leibniz formula) and
expanding the resulting products further. For instance, the product
\[
\left(  t+a\right)  \left(  t+b^{\prime}\right)  d^{\prime\prime}%
c^{\prime\prime\prime}%
\]
becomes $ttd^{\prime\prime}c^{\prime\prime\prime}+tb^{\prime}d^{\prime\prime
}c^{\prime\prime\prime}+atd^{\prime\prime}c^{\prime\prime\prime}+ab^{\prime
}d^{\prime\prime}c^{\prime\prime\prime}$. In the huge sum that results, we are
interested in those addends that contain exactly one $t$, because it is
precisely these addends that contribute to the coefficient of $t^{1}$ in the
polynomial $\det\left(  tI_{n}+L\right)  $. Where do these addends come from?
To pick up exactly one $t$ from a product like $\left(  t+a\right)  \left(
t+b^{\prime}\right)  d^{\prime\prime}c^{\prime\prime\prime}$, we need to have
at least one diagonal entry in our product (for example, we cannot pick up any
$t$ from the product $cd^{\prime}b^{\prime\prime}a^{\prime\prime\prime}$), and
we need to pick out the $t$ from this diagonal entry (rather than, e.g., the
$a$ or $b^{\prime}$ or $c^{\prime\prime}$ or $d^{\prime\prime\prime}$). If we
pick the $r$-th diagonal entry, then the rest of the product is part of the
expansion of $\det\left(  L_{\sim r,\sim r}\right)  $ (since we must not pick
any further $t$s and thus can pretend that they are not there in the first
place). Thus, the total $t^{1}$-coefficient in $\det\left(  tI_{n}+L\right)  $
will be $\sum_{r=1}^{n}\det\left(  L_{\sim r,\sim r}\right)  $. This proves
(\ref{pf.thm.MTT.undir.b.1}), and thus the proof of Theorem
\ref{thm.MTT.undir} \textbf{(b)} is complete. \medskip

\textbf{(c)} Consider the polynomial $\det\left(  tI_{n}+L\right)  $
introduced in part \textbf{(b)}, and in particular its $t^{1}$-coefficient
$c_{1}$.

It is known that the characteristic polynomial $\det\left(  tI_{n}-L\right)  $
of $L$ is a monic polynomial of degree $n$, and that its roots are the
eigenvalues $\lambda_{1},\lambda_{2},\ldots,\lambda_{n}$ of $L$. Hence, it can
be factored as follows:%
\[
\det\left(  tI_{n}-L\right)  =\left(  t-\lambda_{1}\right)  \left(
t-\lambda_{2}\right)  \cdots\left(  t-\lambda_{n}\right)  .
\]
Substituting $-t$ for $t$ on both sides of this equality, we obtain%
\[
\det\left(  -tI_{n}-L\right)  =\left(  -t-\lambda_{1}\right)  \left(
-t-\lambda_{2}\right)  \cdots\left(  -t-\lambda_{n}\right)  .
\]
Multiplying both sides of this equality by $\left(  -1\right)  ^{n}$, we find%
\begin{align*}
\det\left(  tI_{n}+L\right)   &  =\left(  t+\lambda_{1}\right)  \left(
t+\lambda_{2}\right)  \cdots\left(  t+\lambda_{n}\right) \\
&  =\left(  t+\lambda_{1}\right)  \left(  t+\lambda_{2}\right)  \cdots\left(
t+\lambda_{n-1}\right)  t\ \ \ \ \ \ \ \ \ \ \left(  \text{since }\lambda
_{n}=0\right)  .
\end{align*}
Hence, the $t^{1}$-coefficient of the polynomial $\det\left(  tI_{n}+L\right)
$ is $\lambda_{1}\lambda_{2}\cdots\lambda_{n-1}$ (since this is clearly the
$t^{1}$-coefficient on the right hand side). Since we defined $c_{1}$ to be
the $t^{1}$-coefficient of the polynomial $\det\left(  tI_{n}+L\right)  $, we
thus conclude that $c_{1}=\lambda_{1}\lambda_{2}\cdots\lambda_{n-1}$. However,
Theorem \ref{thm.MTT.undir} \textbf{(b)} yields%
\[
\left(  \text{\# of spanning trees of }G\right)  =\dfrac{1}{n}%
\underbrace{c_{1}}_{=\lambda_{1}\lambda_{2}\cdots\lambda_{n-1}}=\dfrac{1}%
{n}\cdot\lambda_{1}\lambda_{2}\cdots\lambda_{n-1}.
\]
This proves Theorem \ref{thm.MTT.undir} \textbf{(c)}.
\end{proof}

\subsubsection{Application: counting spanning trees of $K_{n,m}$}

Laplacians of digraphs often have computable eigenvalues, so Theorem
\ref{thm.MTT.undir} \textbf{(c)} is actually pretty useful. A striking example
of a \# of spanning trees (specifically, of the $n$-hypercube graph $Q_{n}$,
which we already met in Subsection \ref{subsec.hamp.hypercube}) that can be
counted using eigenvalues will appear in Exercise \ref{exe.7.5}.

Here, however, let us give a simpler example, in which Theorem
\ref{thm.MTT.undir} \textbf{(a)} suffices:

\begin{exercise}
\label{exe.MTT.Knm}Let $n$ and $m$ be two positive integers. Let $K_{n,m}$ be
the simple graph with $n+m$ vertices
\[
1,2,\ldots,n\ \ \ \ \ \ \ \ \ \ \text{and}\ \ \ \ \ \ \ \ \ \ -1,-2,\ldots
,-m,
\]
where two vertices $i$ and $j$ are adjacent if and only if they have opposite
signs (i.e., each positive vertex is adjacent to each negative vertex, but no
two vertices of the same sign are adjacent).

[For example, here is how $K_{5,2}$ looks like:
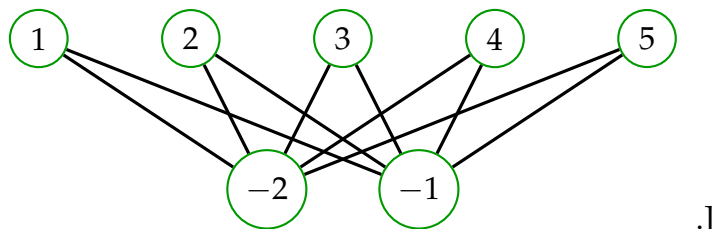
\[%
%
\text{\ \ .]}%
\]

How many spanning trees does $K_{n,m}$ have?
\end{exercise}

\begin{proof}
[Solution.]If we rename the negative vertices $-1,-2,\ldots,-m$ as
$n+1,n+2,\ldots,n+m$, then the Laplacian $L$ of the digraph $K_{n,m}%
^{\operatorname*{bidir}}$ can be written in block-matrix notation as follows:%
\[
L=\left(
\begin{array}
[c]{cc}%
A & B\\
C & D
\end{array}
\right)  ,
\]
where

\begin{itemize}
\item $A$ is a diagonal $n\times n$-matrix whose all diagonal entries are
equal to $m$ (since there are no edges between positive vertices, and since
each positive vertex has degree $m$);

\item $B$ is an $n\times m$-matrix whose all entries equal $-1$;

\item $C$ is an $m\times n$-matrix whose all entries equal $-1$;

\item $D$ is a diagonal $m\times m$-matrix whose all diagonal entries are
equal to $n$.
\end{itemize}

For instance, if $n=3$ and $m=2$, then%
\[
L=\left(
\begin{array}
[c]{ccccc}%
2 & 0 & 0 & -1 & -1\\
0 & 2 & 0 & -1 & -1\\
0 & 0 & 2 & -1 & -1\\
-1 & -1 & -1 & 3 & 0\\
-1 & -1 & -1 & 0 & 3
\end{array}
\right)  .
\]

Theorem \ref{thm.MTT.undir} \textbf{(a)} yields%
\[
\left(  \text{\# of spanning trees of }K_{n,m}\right)  =\det\left(  L_{\sim
r,\sim r}\right)  \ \ \ \ \ \ \ \ \ \ \text{for any vertex }r\text{ of
}K_{n,m};
\]
thus, we need to compute $\det\left(  L_{\sim r,\sim r}\right)  $ for some
vertex $r$. We let $r=1$. Then, the submatrix $L_{\sim r,\sim r}=L_{\sim
1,\sim1}$ of $L$ again can be written in block-matrix notation as follows:%
\begin{equation}
L_{\sim r,\sim r}=\left(
\begin{array}
[c]{cc}%
\widetilde{A} & \widetilde{B}\\
\widetilde{C} & D
\end{array}
\right)  , \label{sol.MTT.Knm.block2}%
\end{equation}
where

\begin{itemize}
\item $\widetilde{A}$ is a diagonal $\left(  n-1\right)  \times\left(
n-1\right)  $-matrix, whose all diagonal entries are equal to $m$;

\item $\widetilde{B}$ is an $\left(  n-1\right)  \times m$-matrix whose all
entries equal $-1$;

\item $\widetilde{C}$ is an $m\times\left(  n-1\right)  $-matrix whose all
entries equal $-1$;

\item $D$ is a diagonal $m\times m$-matrix whose all diagonal entries are
equal to $n$.
\end{itemize}

Fortunately, determinants of block matrices are often not hard to compute, at
least when some of the blocks are invertible. For example, the Schur
complement provides
\href{https://en.wikipedia.org/wiki/Block_matrix#Block_matrix_determinant}{a
neat formula}. Our life here is even easier, since $\widetilde{A}$ and $D$ are
multiples of identity matrices: namely, $\widetilde{A}=mI_{n-1}$ and
$D=nI_{m}$. We perform a \textquotedblleft blockwise row
transformation\textquotedblright\ on the block matrix $L_{\sim r,\sim
r}=\left(
\begin{array}
[c]{cc}%
\widetilde{A} & \widetilde{B}\\
\widetilde{C} & D
\end{array}
\right)  $, specifically subtracting the $\widetilde{C}\widetilde{A}^{-1}%
$-multiple of the first \textquotedblleft block row\textquotedblright%
\ $\left(
\begin{array}
[c]{cc}%
\widetilde{A} & \widetilde{B}%
\end{array}
\right)  $ from the second \textquotedblleft block row\textquotedblright%
\ $\left(
\begin{array}
[c]{cc}%
\widetilde{C} & D
\end{array}
\right)  $ (yes, this is legitimate -- it's the same as left-multiplying by
the block matrix $\left(
\begin{array}
[c]{cc}%
I_{n-1} & 0\\
-\widetilde{C}\widetilde{A}^{-1} & I_{m}%
\end{array}
\right)  $, which has determinant $1$ because it is lower-triangular). As a
result, we obtain%
\begin{align*}
\det\left(
\begin{array}
[c]{cc}%
\widetilde{A} & \widetilde{B}\\
\widetilde{C} & D
\end{array}
\right)   &  =\det\left(
\begin{array}
[c]{cc}%
\widetilde{A} & \widetilde{B}\\
\widetilde{C}-\widetilde{C}\widetilde{A}^{-1}\widetilde{A} & D-\widetilde{C}%
\widetilde{A}^{-1}\widetilde{B}%
\end{array}
\right) \\
&  =\det\left(
\begin{array}
[c]{cc}%
\widetilde{A} & \widetilde{B}\\
0 & D-\widetilde{C}\widetilde{A}^{-1}\widetilde{B}%
\end{array}
\right)  .
\end{align*}
The matrix on the right is \textquotedblleft block-upper
triangular\textquotedblright, so its determinant factors as
follows:\footnote{We are using the fact that if a matrix is block-triangular
(with all diagonal blocks being square matrices), then its determinant is the
product of the determinants of its diagonal blocks. See, e.g.,
\url{https://math.stackexchange.com/a/1221066/} or \cite[Exercise
6.29]{detnotes} for a proof of this fact.}%
\[
\det\left(
\begin{array}
[c]{cc}%
\widetilde{A} & \widetilde{B}\\
0 & D-\widetilde{C}\widetilde{A}^{-1}\widetilde{B}%
\end{array}
\right)  =\det\widetilde{A}\cdot\det\left(  D-\widetilde{C}\widetilde{A}%
^{-1}\widetilde{B}\right)  .
\]
Of course, $\det\widetilde{A}=m^{n-1}$, since $\widetilde{A}$ is a diagonal
matrix with $m,m,\ldots,m$ on the diagonal. Computing $\det\left(
D-\widetilde{C}\widetilde{A}^{-1}\widetilde{B}\right)  $ is a bit more
complicated, but still doable: The matrix $\widetilde{A}^{-1}$ is a diagonal
matrix with $m^{-1},m^{-1},\ldots,m^{-1}$ on the diagonal; thus, its role in
the product $\widetilde{C}\widetilde{A}^{-1}\widetilde{B}$ is merely to
multiply everything by $m^{-1}$. Hence, $\widetilde{C}\widetilde{A}%
^{-1}\widetilde{B}=m^{-1}\widetilde{C}\widetilde{B}$. Since all entries of
$\widetilde{C}$ and $\widetilde{B}$ are $-1$'s, we see that all entries of
$\widetilde{C}\widetilde{B}$ are $\left(  n-1\right)  $'s. Putting all of this
together, we see that $D-\widetilde{C}\widetilde{A}^{-1}\widetilde{B}$ is the
$m\times m$-matrix whose all diagonal entries are equal to $n-m^{-1}\left(
n-1\right)  $ and whose all off-diagonal entries are equal to $-m^{-1}\left(
n-1\right)  $. We have already computed the determinant of a matrix much like
this back in our proof of Cayley's Formula (Subsection \ref{subsec.MTT.Kn});
let us deal with the general case: \Needspace{15pc}

\begin{proposition}
\label{prop.sol.MTT.Knm.det0}Let $n\in\mathbb{N}$. Let $x$ and $a$ be two
numbers. Then,%
\[
\det\underbrace{\left(
\begin{array}
[c]{cccccc}%
x & a & a & \cdots & a & a\\
a & x & a & \cdots & a & a\\
a & a & x & \cdots & a & a\\
\vdots & \vdots & \vdots & \ddots & \vdots & \vdots\\
a & a & a & \cdots & x & a\\
a & a & a & \cdots & a & x
\end{array}
\right)  }_{\substack{\text{the }n\times n\text{-matrix}\\\text{whose diagonal
entries are }x\\\text{and whose off-diagonal entries are }a}}=\left(
x+\left(  n-1\right)  a\right)  \left(  x-a\right)  ^{n-1}.
\]

\end{proposition}

Proposition \ref{prop.sol.MTT.Knm.det0} can be proved using similar reasoning
as the determinant in Subsection \ref{subsec.MTT.Kn}; we will say more about
it later. For now, let us apply it to $m$, $n-m^{-1}\left(  n-1\right)  $ and
$-m^{-1}\left(  n-1\right)  $ instead of $n$, $x$ and $a$, to obtain%
\begin{align*}
\det\left(  D-\widetilde{C}\widetilde{A}^{-1}\widetilde{B}\right)   &
=\underbrace{\left(  \left(  n-m^{-1}\left(  n-1\right)  \right)  +\left(
m-1\right)  \left(  -m^{-1}\left(  n-1\right)  \right)  \right)  }_{=1}\\
&  \ \ \ \ \ \ \ \ \ \ \cdot\left(  \underbrace{\left(  n-m^{-1}\left(
n-1\right)  \right)  -\left(  -m^{-1}\left(  n-1\right)  \right)  }%
_{=n}\right)  ^{m-1}\\
&  =n^{m-1}.
\end{align*}

Now, it is time to combine everything we know. Theorem \ref{thm.MTT.undir}
\textbf{(a)} yields%
\begin{align*}
\left(  \text{\# of spanning trees of }K_{n,m}\right)   &  =\det\left(
L_{\sim r,\sim r}\right) \\
&  =\det\left(
\begin{array}
[c]{cc}%
\widetilde{A} & \widetilde{B}\\
\widetilde{C} & D
\end{array}
\right)  \ \ \ \ \ \ \ \ \ \ \left(  \text{by (\ref{sol.MTT.Knm.block2}%
)}\right) \\
&  =\det\left(
\begin{array}
[c]{cc}%
\widetilde{A} & \widetilde{B}\\
0 & D-\widetilde{C}\widetilde{A}^{-1}\widetilde{B}%
\end{array}
\right) \\
&  =\underbrace{\det\widetilde{A}}_{=m^{n-1}}\cdot\underbrace{\det\left(
D-\widetilde{C}\widetilde{A}^{-1}\widetilde{B}\right)  }_{=n^{m-1}}\\
&  =m^{n-1}\cdot n^{m-1}.
\end{align*}

\end{proof}

Thus, we have obtained the following:

\begin{theorem}
Let $n$ and $m$ be two positive integers. Let $K_{n,m}$ be the simple graph
with $n+m$ vertices
\[
1,2,\ldots,n\ \ \ \ \ \ \ \ \ \ \text{and}\ \ \ \ \ \ \ \ \ \ -1,-2,\ldots
,-m,
\]
where two vertices $i$ and $j$ are adjacent if and only if they have opposite
signs. Then,
\[
\left(  \text{\# of spanning trees of }K_{n,m}\right)  =m^{n-1}\cdot n^{m-1}.
\]

\end{theorem}

See \cite{AbuSbe88} for a combinatorial proof of this theorem.

\begin{exercise}
\label{exe.6.6}Let $n$ be a positive integer. Let $K_{n,2}$ be the simple
graph with vertex set $\left\{  1,2,\ldots,n\right\}  \cup\left\{
-1,-2\right\}  $ such that two vertices of $K_{n,2}$ are adjacent if and only
if they have opposite signs (i.e., each positive vertex is adjacent to each
negative vertex, but no two vertices of the same sign are adjacent). We regard
$K_{n,2}$ as a multigraph in the usual way.

\begin{enumerate}
\item[\textbf{(a)}] Without using the matrix-tree theorem, prove that the
number of spanning trees of $K_{n, 2}$ is $n \cdot2^{n-1}$.

\item[\textbf{(b)}] Let $K_{n,2}^{\prime}$ be the graph obtained by adding a
new edge $\left\{  -1,-2\right\}  $ to $K_{n,2}$. How many spanning trees does
$K_{n,2}^{\prime}$ have?
\end{enumerate}

[\textbf{Example:} Here is the graph $K_{n, 2}$ for $n = 5$:
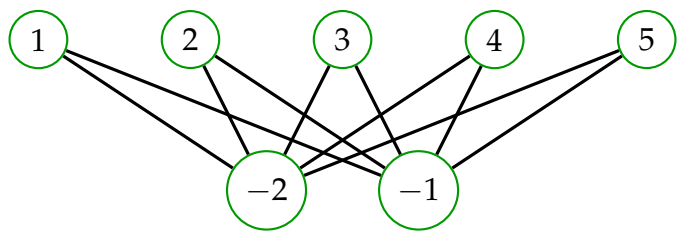
\[
\begin{tikzpicture}[scale=2]
\begin{scope}[every node/.style={circle,thick,draw=green!60!black}]
\node(-1) at (2.5,-1) {$-1$};
\node(-2) at (1.5,-1) {$-2$};
\node(1) at (0,0) {$1$};
\node(2) at (1,0) {$2$};
\node(3) at (2,0) {$3$};
\node(4) at (3,0) {$4$};
\node(5) at (4,0) {$5$};
\end{scope}
\begin{scope}[every edge/.style={draw=black,very thick}, every loop/.style={}]
\path[-] (-1) edge (1) edge (2) edge (3) edge (4) edge (5);
\path[-] (-2) edge (1) edge (2) edge (3) edge (4) edge (5);
\end{scope}
\end{tikzpicture}
\]
And here is the corresponding graph $K^{\prime}_{n, 2}$:
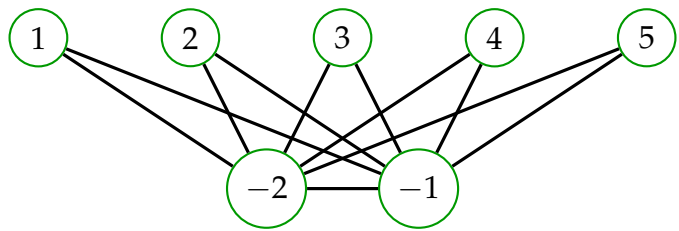
\[
\begin{tikzpicture}[scale=2]
\begin{scope}[every node/.style={circle,thick,draw=green!60!black}]
\node(-1) at (2.5,-1) {$-1$};
\node(-2) at (1.5,-1) {$-2$};
\node(1) at (0,0) {$1$};
\node(2) at (1,0) {$2$};
\node(3) at (2,0) {$3$};
\node(4) at (3,0) {$4$};
\node(5) at (4,0) {$5$};
\end{scope}
\begin{scope}[every edge/.style={draw=black,very thick}, every loop/.style={}]
\path[-] (-1) edge (1) edge (2) edge (3) edge (4) edge (5);
\path[-] (-2) edge (-1) edge (1) edge (2) edge (3) edge (4) edge (5);
\end{scope}
\end{tikzpicture}
\]
]
\end{exercise}

\begin{exercise}
\label{exe.7.11}Let $n$ be a positive integer. Let $A$ be the $\left(
n-1\right)  \times\left(  n-1\right)  $-matrix
\[
\left(
\begin{array}
[c]{ccccc}%
2 & -1 & 0 & \cdots & 0\\
-1 & 2 & -1 & \cdots & 0\\
0 & -1 & 2 & \cdots & 0\\
\vdots & \vdots & \vdots & \ddots & \vdots\\
0 & 0 & 0 & \cdots & 2
\end{array}
\right)  ,
\]
whose $\left(  i,j\right)  $-th entry is%
\[
A_{i,j}:=%
\begin{cases}
2, & \text{if }i=j;\\
-1, & \text{if }\left\vert i-j\right\vert =1;\\
0, & \text{otherwise}%
\end{cases}
\ \ \ \ \ \ \ \ \ \ \text{for all }i,j\in\left\{  1,2,\ldots,n-1\right\}  .
\]
Prove that $\det A=n$. \medskip

[\textbf{Hint:} Recall Example \ref{exa.spt.Cn}.]
\end{exercise}

\begin{exercise}
Let $G$ be a multigraph with an even number of vertices. Assume that each
vertex of $G$ has an even degree. Prove that $G$ has an even number of
spanning trees. \medskip

[\textbf{Hint:} First show a lemma: Let $n\in\mathbb{N}$ be odd. Let $B$ be an
$n\times n$-matrix with integer entries. Assume that $B$ is symmetric (i.e.,
satisfies $B^{T}=B$), and that all diagonal entries of $B$ are even. Then,
$\det B$ is even.]
\end{exercise}

\begin{exercise}
\label{exe.7.5}Let $n$ be a positive integer. Let $Q_{n}$ be the $n$-hypercube
graph (as defined in Definition \ref{def.hypercube}). Recall that its vertex
set is the set $V:=\left\{  0,1\right\}  ^{n}$ of length-$n$ bitstrings, and
that two vertices are adjacent if and only if they differ in exactly one bit.
Our goal is to compute the \# of spanning trees of $Q_{n}$.

Let $D$ be the digraph $Q_{n}^{\operatorname*{bidir}}$. Let $L$ be the
Laplacian of $D$. We regard $L$ as a $V\times V$-matrix (i.e., as a
$2^{n}\times2^{n}$-matrix whose rows and columns are indexed by bitstrings in
$V$).

We shall use the notation $a_{i}$ for the $i$-th entry of a bitstring $a$.
Thus, each bitstring $a \in V$ has the form $a = \left(  a_{1}, a_{2}, \ldots,
a_{n} \right)  $. (We shall avoid the shorthand notation $a_{1} a_{2} \cdots
a_{n}$ here, as it could be mistaken for an actual product.)

For any two bitstrings $a,b\in V$, we define the number $\left\langle
a,b\right\rangle $ to be the integer $a_{1}b_{1}+a_{2}b_{2}+\cdots+a_{n}b_{n}$.

\begin{enumerate}
\item[\textbf{(a)}] Prove that every bitstring $a\in V$ satisfies
\[
\sum_{b\in V}\left(  -1\right)  ^{\left\langle a,b\right\rangle }=%
\begin{cases}
2^{n}, & \text{ if $a=\mathbf{0}$;}\\
0, & \text{ otherwise.}%
\end{cases}
\]
Here, $\mathbf{0}$ denotes the bitstring $\left(  0,0,\ldots,0\right)  \in V$.
\end{enumerate}

Now, define a further $V\times V$-matrix $G$ by requiring that its $\left(
a,b\right)  $-th entry is
\[
G_{a,b}=\left(  -1\right)  ^{\left\langle a,b\right\rangle }\qquad\text{for
any $a,b\in V$}.
\]
Furthermore, define a diagonal $V\times V$-matrix $D$ by requiring that its
$\left(  a,a\right)  $-th entry is
\begin{align*}
D_{a,a}  &  =2\cdot\left(  \#\text{ of $i\in\left\{  1,2,\ldots,n\right\}  $
such that $a_{i}=1$}\right) \\
&  =2\cdot\left(  \text{the number of $1$s in $a$}\right)  \qquad\text{for any
$a\in V$}%
\end{align*}
(and its off-diagonal entries are $0$).

Prove the following:

\begin{enumerate}
\item[\textbf{(b)}] We have $G^{2} = 2^{n} \cdot I$, where $I$ is the identity
$V \times V$-matrix.

\item[\textbf{(c)}] We have $GLG^{-1} = D$.

\item[\textbf{(d)}] The eigenvalues of $L$ are $2k$ for all $k \in\left\{  0,
1, \ldots, n \right\}  $, and each eigenvalue $2k$ appears with multiplicity
$\dbinom{n}{k}$.

\item[\textbf{(e)}] The \# of spanning trees of $Q_{n}$ is
\[
\dfrac{1}{2^{n}}\prod_{k=1}^{n}\left(  2k\right)  ^{\binom{n}{k}}.
\]

\end{enumerate}

[\textbf{Example:} As an example, here is the case $n=3$. In this case, the
graph $Q_{n}$ looks as follows:
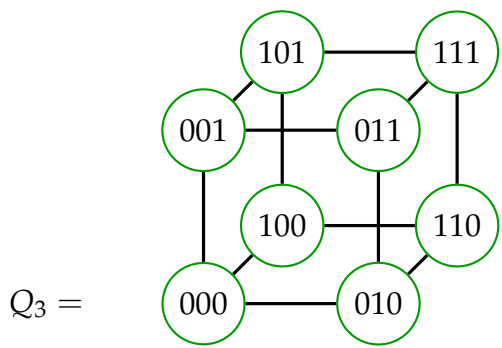
\[
\begin{tikzpicture}[scale=2.3]
\begin{scope}[every node/.style={circle,thick,draw=green!60!black}]
\node(000) at (0, 0) {$000$};
\node(001) at (0, 1) {$001$};
\node(010) at (1, 0) {$010$};
\node(011) at (1, 1) {$011$};
\node(100) at (0.45, 0.45) {$100$};
\node(101) at (0.45, 1.45) {$101$};
\node(110) at (1.45, 0.45) {$110$};
\node(111) at (1.45, 1.45) {$111$};
\end{scope}
\node(X) at (-0.9, 0) {$Q_3 = $};
\begin{scope}[every edge/.style={draw=black,very thick}]
\path[-] (000) edge (001) (001) edge (011) (011) edge (010) (010) edge (000);
\path[-] (100) edge (101) (101) edge (111) (111) edge (110) (110) edge (100);
\path[-] (100) edge (000) (101) edge (001) (111) edge (011) (110) edge (010);
\end{scope}
\end{tikzpicture}
\]
The matrices $L$, $G$ and $D$ are%
\[
L=%
\begin{pmatrix}
3 & -1 & -1 & 0 & -1 & 0 & 0 & 0\\
-1 & 3 & 0 & -1 & 0 & -1 & 0 & 0\\
-1 & 0 & 3 & -1 & 0 & 0 & -1 & 0\\
0 & -1 & -1 & 3 & 0 & 0 & 0 & -1\\
-1 & 0 & 0 & 0 & 3 & -1 & -1 & 0\\
0 & -1 & 0 & 0 & -1 & 3 & 0 & -1\\
0 & 0 & -1 & 0 & -1 & 0 & 3 & -1\\
0 & 0 & 0 & -1 & 0 & -1 & -1 & 3
\end{pmatrix}
,
\]
\[
G=%
\begin{pmatrix}
1 & 1 & 1 & 1 & 1 & 1 & 1 & 1\\
1 & -1 & 1 & -1 & 1 & -1 & 1 & -1\\
1 & 1 & -1 & -1 & 1 & 1 & -1 & -1\\
1 & -1 & -1 & 1 & 1 & -1 & -1 & 1\\
1 & 1 & 1 & 1 & -1 & -1 & -1 & -1\\
1 & -1 & 1 & -1 & -1 & 1 & -1 & 1\\
1 & 1 & -1 & -1 & -1 & -1 & 1 & 1\\
1 & -1 & -1 & 1 & -1 & 1 & 1 & -1
\end{pmatrix}
,
\]%
\[
D=%
\begin{pmatrix}
0 & 0 & 0 & 0 & 0 & 0 & 0 & 0\\
0 & 2 & 0 & 0 & 0 & 0 & 0 & 0\\
0 & 0 & 2 & 0 & 0 & 0 & 0 & 0\\
0 & 0 & 0 & 2 & 0 & 0 & 0 & 0\\
0 & 0 & 0 & 0 & 4 & 0 & 0 & 0\\
0 & 0 & 0 & 0 & 0 & 4 & 0 & 0\\
0 & 0 & 0 & 0 & 0 & 0 & 4 & 0\\
0 & 0 & 0 & 0 & 0 & 0 & 0 & 6
\end{pmatrix}
,
\]
where the rows and the columns are ordered by listing the eight bitstrings
$a\in V$ in the order $000,001,010,011,100,101,110,111$. ]
\end{exercise}

\bigskip

As we promised, let us make a few more remarks about Proposition
\ref{prop.sol.MTT.Knm.det0}. While this proposition can be proved by fairly
straightforward row transformations (first subtracting the first row from all
the other rows, then factoring an $x-a$ from all the latter rows, then
subtracting $a$ times each of the latter rows to the first row to obtain a
triangular matrix), it can also be viewed as a particular case of either of
the following two determinantal identities: \Needspace{9pc}

\begin{proposition}
\label{prop.sol.MTT.Knm.det1}Let $n\in\mathbb{N}$. Let $a_{1},a_{2}%
,\ldots,a_{n}$ be $n$ numbers, and let $x$ be a further number. Then,%
\[
\det\underbrace{\left(
\begin{array}
[c]{cccccc}%
x & a_{1} & a_{2} & \cdots & a_{n-1} & a_{n}\\
a_{1} & x & a_{2} & \cdots & a_{n-1} & a_{n}\\
a_{1} & a_{2} & x & \cdots & a_{n-1} & a_{n}\\
\vdots & \vdots & \vdots & \ddots & \vdots & \vdots\\
a_{1} & a_{2} & a_{3} & \cdots & x & a_{n}\\
a_{1} & a_{2} & a_{3} & \cdots & a_{n} & x
\end{array}
\right)  }_{\text{an }\left(  n+1\right)  \times\left(  n+1\right)
\text{-matrix}}=\left(  x+\sum_{i=1}^{n}a_{i}\right)  \prod_{i=1}^{n}\left(
x-a_{i}\right)  .
\]

\end{proposition}

\begin{proposition}
\label{prop.sol.MTT.Knm.det2}Let $n\in\mathbb{N}$. Let $x_{1},x_{2}%
,\ldots,x_{n}$ be $n$ numbers, and let $a$ be a further number. Then,%
\[
\det\left(
\begin{array}
[c]{ccccc}%
x_{1} & a & a & \cdots & a\\
a & x_{2} & a & \cdots & a\\
a & a & x_{3} & \cdots & a\\
\vdots & \vdots & \vdots & \ddots & \vdots\\
a & a & a & \cdots & x_{n}%
\end{array}
\right)  =\prod_{i=1}^{n}\left(  x_{i}-a\right)  +a\sum_{i=1}^{n}y_{i},
\]
where we set $y_{i}:=\prod_{\substack{k\in\left\{  1,2,\ldots,n\right\}
;\\k\neq i}}\left(  x_{k}-a\right)  $ for each $i\in\left\{  1,2,\ldots
,n\right\}  $.
\end{proposition}

Both of these propositions make good exercises in determinant evaluation.
(Proposition \ref{prop.sol.MTT.Knm.det1} is \cite[Exercise 6.21]{detnotes},
while Proposition \ref{prop.sol.MTT.Knm.det2} is
\url{https://math.stackexchange.com/a/2112473/} .)

See \cite{KleSta19} and \cite{Rubey00} for more applications of the
Matrix-Tree Theorem, and \cite{Holzer22} for many more related results.

\subsection{de Bruijn sequences}

\subsubsection{Definition}

Let me move on to a more intricate application of what we have learned about arborescences.

A little puzzle first: What is special about the periodic sequence%
\[
||\vphantom{a}:\ 0000\ 1111\ 0110\ 0101\ :{||}\ \ \ \ \ \ \ \ \ \text{?}%
\]
(This is an infinite sequence of $0$'s and $1$'s; the spaces between some of
them are only for readability. The $||\vphantom{a}:$ and $:||$ symbols are
\textquotedblleft repeat signs\textquotedblright\ -- they mean that everything
that stands between them should be repeated over and over. So the sequence
above is $0000\ 1111\ 0110\ 0101\ 0000\ 1111\ \ldots$.)

One nice property of this sequence is that if you slide a \textquotedblright
length-$4$ window\textquotedblright\ (i.e., a window that shows four
consecutive entries) along it, you get all $16$ possible bitstrings of length
$4$ depending on the position of the window, and these bitstrings do not
repeat until you move $16$ steps to the right. Just see:%
\begin{align*}
&  \fbox{$0000$}11110110010100001111\ldots\\
&  0\fbox{$0001$}1110110010100001111\ldots\\
&  00\fbox{$0011$}110110010100001111\ldots\\
&  000\fbox{$0111$}10110010100001111\ldots\\
&  0000\fbox{$1111$}0110010100001111\ldots\\
&  00001\fbox{$1110$}110010100001111\ldots\\
&  000011\fbox{$1101$}10010100001111\ldots\\
&  0000111\fbox{$1011$}0010100001111\ldots\\
&  00001111\fbox{$0110$}010100001111\ldots\\
&  000011110\fbox{$1100$}10100001111\ldots\\
&  0000111101\fbox{$1001$}0100001111\ldots\\
&  00001111011\fbox{$0010$}100001111\ldots\\
&  000011110110\fbox{$0101$}00001111\ldots\\
&  0000111101100\fbox{$1010$}0001111\ldots\\
&  00001111011001\fbox{$0100$}001111\ldots\\
&  000011110110010\fbox{$1000$}01111\ldots
\end{align*}
Note that, as you slide the window along the sequence, at each step, the first
bit is removed and a new bit is inserted at the end. Thus, by sliding a
length-$4$ window along the above sequence, you run through all $16$ possible
length-$4$ bitstrings in such a way that each bitstring is obtained from the
previous one by removing the first bit and inserting a new bit at the end.
This is nice and somewhat similar to Gray codes (in which you run through all
bitstrings of a given length in such a way that only a single bit is changed
at each step).

Can we find such nice sequences for any window length, not just $4$ ?

Here is an answer for window length $3$, for instance:%
\[
||\vphantom{a}:\ 000\ 111\ 01\ :{||}\ \ .
\]
What about higher window length?

Moreover, we can ask the same question with other alphabets. For instance,
instead of bits, here is a similar sequence for the alphabet $\left\{
0,1,2\right\}  $ (that is, we use the numbers $0,1,2$ instead of $0$ and $1$)
and window length $2$:%
\[
||\vphantom{a}:00\ 11\ 22\ 02\ 1\ :{||}\ .
\]

What about the general case? Let us give it a name:

\begin{definition}
\label{def.debr.debr}Let $n$ and $k$ be two positive integers, and let $K$ be
a $k$-element set.

A \textbf{de Bruijn sequence} of order $n$ on $K$ means a $k^{n}$-tuple
$\left(  c_{0},c_{1},\ldots,c_{k^{n}-1}\right)  $ of elements of $K$ such that

\begin{enumerate}
\item[\textbf{(A)}] for each $n$-tuple $\left(  a_{1},a_{2},\ldots
,a_{n}\right)  \in K^{n}$ of elements of $K$, there is a \textbf{unique}
$r\in\left\{  0,1,\ldots,k^{n}-1\right\}  $ such that%
\[
\left(  a_{1},a_{2},\ldots,a_{n}\right)  =\left(  c_{r},c_{r+1},\ldots
,c_{r+n-1}\right)  .
\]

\end{enumerate}

Here, the indices under the letter \textquotedblleft$c$\textquotedblright\ are
understood to be periodic modulo $k^{n}$; that is, we set $c_{q+k^{n}}=c_{q}$
for each $q\in\mathbb{Z}$ (so that $c_{k^{n}}=c_{0}$ and $c_{k^{n}+1}=c_{1}$
and so on).
\end{definition}

For example, for $n=2$ and $k=3$ and $K=\left\{  0,1,2\right\}  $, the
$9$-tuple%
\[
\left(  0,0,1,1,2,2,0,2,1\right)
\]
is a de Bruijn sequence of order $n$ on $K$, because if we label the entries
of this $9$-tuple as $c_{0},c_{1},\ldots,c_{8}$ (and extend the indices
periodically, so that $c_{9}=c_{0}$), then we have%
\begin{align*}
\left(  0,0\right)   &  =\left(  c_{0},c_{1}\right)  ;\qquad\left(
0,1\right)  =\left(  c_{1},c_{2}\right)  ;\qquad\left(  0,2\right)  =\left(
c_{6},c_{7}\right)  ;\\
\left(  1,0\right)   &  =\left(  c_{8},c_{9}\right)  ;\qquad\left(
1,1\right)  =\left(  c_{2},c_{3}\right)  ;\qquad\left(  1,2\right)  =\left(
c_{3},c_{4}\right)  ;\\
\left(  2,0\right)   &  =\left(  c_{5},c_{6}\right)  ;\qquad\left(
2,1\right)  =\left(  c_{7},c_{8}\right)  ;\qquad\left(  2,2\right)  =\left(
c_{4},c_{5}\right)  .
\end{align*}
This de Bruijn sequence $\left(  0,0,1,1,2,2,0,2,1\right)  $ corresponds to
the periodic sequence $||\vphantom{a}:00\ 11\ 22\ 02\ 1\ :{||}$ that we found above.

\subsubsection{Existence of de Bruijn sequences}

It turns out that de Bruijn sequences always exist:

\begin{theorem}
[de Bruijn, Sainte-Marie]\label{thm.debr.exist}Let $n$ and $k$ be positive
integers. Let $K$ be a $k$-element set. Then, a de Bruijn sequence of order
$n$ on $K$ exists.
\end{theorem}

\begin{proof}
It looks reasonable to approach this using a digraph. For example, we can
define a digraph whose vertices are the $n$-tuples in $K^{n}$, and that has an
arc from one $n$-tuple $i$ to another $n$-tuple $j$ if $j$ can be obtained
from $i$ by dropping the first entry and adding a new entry at the end. Then,
a de Bruijn sequence (of order $n$ on $K$) is the same as a Hamiltonian cycle
of this digraph.

Unfortunately, we don't have any useful criteria that would show that such a
cycle exists. So this idea seems to be a dead end.

However, let us do something counterintuitive: We try to reinterpret de Bruijn
sequences in terms of Eulerian circuits (rather than Hamiltonian cycles),
since we have a good criterion for the existence of Eulerian circuits (unlike
for that of Hamiltonian cycles)!

We need a different digraph for that. Namely, we let $D$ be the multidigraph
$\left(  K^{n-1},K^{n},\psi\right)  $, where the map $\psi:K^{n}\rightarrow
K^{n-1}\times K^{n-1}$ is given by the formula%
\[
\psi\left(  a_{1},a_{2},\ldots,a_{n}\right)  =\left(  \left(  a_{1}%
,a_{2},\ldots,a_{n-1}\right)  ,\ \ \left(  a_{2},a_{3},\ldots,a_{n}\right)
\right)  .
\]
Thus, the vertices of $D$ are the $\left(  n-1\right)  $-tuples (not the
$n$-tuples!) of elements of $K$, whereas the arcs are the $n$-tuples of
elements of $K$, and each such arc $\left(  a_{1},a_{2},\ldots,a_{n}\right)  $
has source $\left(  a_{1},a_{2},\ldots,a_{n-1}\right)  $ and target $\left(
a_{2},a_{3},\ldots,a_{n}\right)  $. Hence, there is an arc from each $\left(
n-1\right)  $-tuple $i\in K^{n-1}$ to each $\left(  n-1\right)  $-tuple $j\in
K^{n-1}$ that is obtained by dropping the first entry of $i$ and adding a new
entry at the end. (Be careful: If $n=1$, then $D$ has only one vertex but $n$
arcs. If this confuses you, just do the $n=1$ case by hand. For any $n>1$,
there are no parallel arcs in $D$.)

\begin{example}
For example, if $n=3$ and $k=2$ and $K=\left\{  0,1\right\}  $, then $D$ looks
as follows (we again write our tuples without commas and without parentheses):%
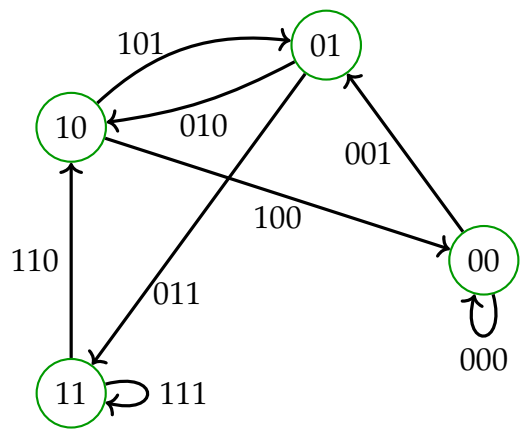
\[%
\begin{tikzpicture}[scale=1.5]
\begin{scope}[every node/.style={circle,thick,draw=green!60!black}]
\node(00) at (0:2) {$00$};
\node(01) at (360/5:2) {$01$};
\node(10) at (2*360/5:2) {$10$};
\node(11) at (3*360/5:2) {$11$};
\end{scope}
\begin{scope}[every edge/.style={draw=black,very thick}, every loop/.style={}]
\path[->] (00) edge node[left] {$001$}
(01) edge[loop below] node[below] {$000$} (00);
\path[->] (10) edge[bend left=25] node[pos=0.6, left, inner sep=1.2em] {$101$}
(01) edge node[below] {$100$} (00);
\path[->] (01) edge[bend left=10] node[below] {$010$}
(10) edge node[pos=0.6, below, inner sep=1em] {$011$}(11);
\path[->] (11) edge node[left] {$110$}
(10) edge[loop right] node[right] {$111$} (11);
\end{scope}
\end{tikzpicture}%
\]

\end{example}

Let us make a few observations about $D$:

\begin{itemize}
\item The multidigraph $D$ is strongly connected.

[\textit{Proof:} We need to show that for any two vertices $i$ and $j$ of $D$,
there is a walk from $i$ to $j$. But this is easy: Just insert the entries of
$j$ into $i$ one by one, pushing out the entries of $i$. In other words, using
the notation $k_{p}$ for the $p$-th entry of any tuple $k$, we have the walk%
\begin{align*}
i  &  =\left(  i_{1},i_{2},\ldots,i_{n-1}\right) \\
&  \rightarrow\left(  i_{2},i_{3},\ldots,i_{n-1},j_{1}\right) \\
&  \rightarrow\left(  i_{3},i_{4},\ldots,i_{n-1},j_{1},j_{2}\right) \\
&  \rightarrow\cdots\\
&  \rightarrow\left(  i_{n-1},j_{1},j_{2},\ldots,j_{n-2}\right) \\
&  \rightarrow\left(  j_{1},j_{2},\ldots,j_{n-1}\right)  =j.
\end{align*}
Note that this walk has length $n-1$, and is the unique walk from $i$ to $j$
that has length $n-1$. Thus, the \# of walks from $i$ to $j$ that have length
$n-1$ is $1$. This will come useful further below.]

\item Thus, the multidigraph $D$ is weakly connected (since any strongly
connected digraph is weakly connected).

\item The multidigraph $D$ is balanced, and in fact each vertex of $D$ has
outdegree $k$ and indegree $k$.

[\textit{Proof:} Let $i$ be a vertex of $D$. The arcs with source $i$ are the
$n$-tuples whose first $n-1$ entries form the $\left(  n-1\right)  $-tuple $i$
while the last, $n$-th entry is an arbitrary element of $K$. Thus, there are
$\left\vert K\right\vert $ many such arcs. In other words, $i$ has outdegree
$k$. A similar argument shows that $i$ has indegree $k$. This entails that
$\deg^{-}i=\deg^{+}i$. Since this holds for every vertex $i$, we conclude that
$D$ is balanced.]

\item The digraph $D$ has an Eulerian circuit.

[\textit{Proof:} This follows from the directed Euler--Hierholzer theorem
(Theorem \ref{thm.digr.euler-hier}), since $D$ is weakly connected and
balanced. Alternatively, we can derive this from the BEST theorem (Theorem
\ref{thm.BEST.from}) as follows: Pick an arbitrary arc $a$ of $D$, and let $r$
be its source. Then, $r$ is a from-root of $D$ (since $D$ is strongly
connected), and thus $D$ has a spanning arborescence rooted from $r$ (by
Theorem \ref{thm.spanning-arbor.exists}). In other words, using the notations
of the BEST theorem (Theorem \ref{thm.BEST.from}), we have $\tau\left(
D,r\right)  \neq0$. Moreover, each vertex of $D$ has indegree $k>0$. Thus, the
BEST theorem yields
\[
\varepsilon\left(  D,a\right)  =\underbrace{\tau\left(  D,r\right)  }_{\neq
0}\cdot\underbrace{\prod_{u\in V}\left(  \deg^{-}u-1\right)  !}_{\neq0}\neq0.
\]
But this shows that $D$ has an Eulerian circuit whose last arc is $a$.]
\end{itemize}

So we know that $D$ has an Eulerian circuit $\mathbf{c}$. This Eulerian
circuit leads to a de Bruijn sequence as follows:

Let $p_{0},p_{1},\ldots,p_{k^{n}-1}$ be the arcs of $\mathbf{c}$ (from first
to last). Extend the subscripts periodically modulo $k^{n}$ (that is, set
$p_{q+k^{n}}=p_{q}$ for all $q\in\mathbb{N}$). Thus, we obtain an infinite
walk\footnote{We have never formally defined infinite walks, but it should be
fairly clear what they are.} with arcs $p_{0},p_{1},p_{2},\ldots$ (since
$\mathbf{c}$ is a circuit). In other words, for each $i\in\mathbb{N}$, the
target of the arc $p_{i}$ is the source of the arc $p_{i+1}$.

In other words, for each $i\in\mathbb{N}$, the last $n-1$ entries of $p_{i}$
are the first $n-1$ entries of $p_{i+1}$ (since the target of $p_{i}$ is the
tuple consisting of the last $n-1$ entries of $p_{i}$, whereas the source of
$p_{i+1}$ is the tuple consisting of the first $n-1$ entries of $p_{i+1}$).
Therefore, for each $i\in\mathbb{N}$ and each $j\in\left\{  2,3,\ldots
,n\right\}  $, we have%
\begin{align}
&  \left(  \text{the }j\text{-th entry of }p_{i}\right) \nonumber\\
&  =\left(  \text{the }\left(  j-1\right)  \text{-st entry of }p_{i+1}\right)
. \label{pf.thm.debr.exist.j-1}%
\end{align}

Now, for each $i\in\mathbb{N}$, we let $x_{i}$ denote the first entry of the
$n$-tuple $p_{i}$. Then, $x_{q+k^{n}}=x_{q}$ for all $q\in\mathbb{N}$ (since
$p_{q+k^{n}}=p_{q}$ for all $q\in\mathbb{N}$). In other words, the sequence
$\left(  x_{0},x_{1},x_{2},\ldots\right)  $ repeats itself every $k^{n}$
terms. Note that the $k^{n}$-tuple $\left(  x_{0},x_{1},\ldots,x_{k^{n}%
-1}\right)  $ consists of the first entries of the arcs $p_{0},p_{1}%
,\ldots,p_{k^{n}-1}$ of $\mathbf{c}$ (by the definition of $x_{i}$).

For each $i\in\mathbb{N}$ and each $s\in\left\{  1,2,\ldots,n\right\}  $, we
have%
\begin{align*}
&  \left(  \text{the }s\text{-th entry of }p_{i}\right) \\
&  =\left(  \text{the }\left(  s-1\right)  \text{-st entry of }p_{i+1}\right)
\ \ \ \ \ \ \ \ \ \ \left(  \text{by (\ref{pf.thm.debr.exist.j-1})}\right) \\
&  =\left(  \text{the }\left(  s-2\right)  \text{-nd entry of }p_{i+2}\right)
\ \ \ \ \ \ \ \ \ \ \left(  \text{by (\ref{pf.thm.debr.exist.j-1})}\right) \\
&  =\left(  \text{the }\left(  s-3\right)  \text{-rd entry of }p_{i+3}\right)
\ \ \ \ \ \ \ \ \ \ \left(  \text{by (\ref{pf.thm.debr.exist.j-1})}\right) \\
&  =\cdots\\
&  =\left(  \text{the }1\text{-st entry of }p_{i+s-1}\right) \\
&  =x_{i+s-1}\ \ \ \ \ \ \ \ \ \ \left(  \text{since }x_{i+s-1}\text{ was
defined as the first entry of }p_{i+s-1}\right)  .
\end{align*}
In other words, for each $i\in\mathbb{N}$, the entries of $p_{i}$ (from first
to last) are \newline$x_{i},x_{i+1},\ldots,x_{i+n-1}$. In other words, for
each $i\in\mathbb{N}$, we have%
\begin{equation}
p_{i}=\left(  x_{i},x_{i+1},\ldots,x_{i+n-1}\right)  .
\label{pf.thm.debr.exist.p=x}%
\end{equation}

Now, recall that $\mathbf{c}$ is an Eulerian circuit. Thus, each arc of $D$
appears exactly once among its arcs $p_{0},p_{1},\ldots,p_{k^{n}-1}$. In other
words, each $n$-tuple in $K^{n}$ appears exactly once among $p_{0}%
,p_{1},\ldots,p_{k^{n}-1}$ (since the arcs of $D$ are the $n$-tuples in
$K^{n}$). In other words, as $i$ ranges from $0$ to $k^{n}-1$, the $n$-tuple
$p_{i}$ takes each possible value in $K^{n}$ exactly once.

In view of (\ref{pf.thm.debr.exist.p=x}), we can rewrite this as follows: As
$i$ ranges from $0$ to $k^{n}-1$, the $n$-tuple $\left(  x_{i},x_{i+1}%
,\ldots,x_{i+n-1}\right)  $ takes each possible value in $K^{n}$ exactly once
(since this $n$-tuple is precisely $p_{i}$, as we have shown in the previous
paragraph). In other words, for each $\left(  a_{1},a_{2},\ldots,a_{n}\right)
\in K^{n}$, there is a \textbf{unique} $r\in\left\{  0,1,\ldots,k^{n}%
-1\right\}  $ such that $\left(  a_{1},a_{2},\ldots,a_{n}\right)  =\left(
x_{r},x_{r+1},\ldots,x_{r+n-1}\right)  $.

Hence, the $k^{n}$-tuple $\left(  x_{0},x_{1},\ldots,x_{k^{n}-1}\right)  $ is
a de Bruijn sequence of order $n$ on $K$. This shows that a de Bruijn sequence
exists. Theorem \ref{thm.debr.exist} is thus proven.

\begin{example}
For $n=3$ and $k=2$ and $K=\left\{  0,1\right\}  $, one possible Eulerian
circuit $\mathbf{c}$ of $D$ is%
\[
\left(  00,\ \mathbf{001},\ 01,\ \mathbf{010},\ 10,\ \mathbf{101}%
,\ 01,\ \mathbf{011},\ 11,\ \mathbf{111},\ 11,\ \mathbf{110}%
,\ 10,\ \mathbf{100},\ 00\right)
\]
(where we have written the arcs in bold for readability). The first entries of
the arcs of this circuit form the sequence $0010111$, which is indeed a de
Bruijn sequence of order $3$ on $\left\{  0,1\right\}  $. Any $3$ consecutive
entries of this sequence (extended periodically to the infinite sequence
$||\vphantom{a}:0010111:||$) form the respective arc of $\mathbf{c}$.
\end{example}
\end{proof}

\bigskip

Theorem \ref{thm.debr.exist} is merely the starting point of a theory. Several
specific de Bruijn sequences are known, many of them having peculiar
properties. See \cite{Freder82} for a survey of various such
sequences\footnote{Some of these sequences (the \textquotedblleft
prefer-one\textquotedblright\ and \textquotedblleft
prefer-opposite\textquotedblright\ generators) are just disguised
implementations of the algorithm for finding an Eulerian circuit implicit in
our proof of the BEST theorem.} (note that they are called \textquotedblleft
full length nonlinear shift register sequences\textquotedblright\ in this
survey).\footnote{My favorite is the one obtained by concatenating all Lyndon
words whose length divides $n$ in lexicographically increasing order (assuming
that the set $K$ is totally ordered). See \cite{Moreno04} for the details of
that construction.} For the (somewhat confusing) history of de Bruijn
sequences (and of Theorem \ref{thm.debr.exist} and of the Theorem
\ref{thm.debr.num} to be stated later), see \cite{deBrui75}.

There are also several variations on de Bruijn sequences. For some of them,
see \cite{ChDiGr92}. (Note that some of the open questions in that paper are
still unsolved.) A variation that recently became quite popular is the notion
of a \textquotedblleft universal cycle for permutations\textquotedblright\ --
a string that contains all \textquotedblleft permutations\textquotedblright%
\ (more precisely, $n$-tuples of distinct elements of $K$) as factors. See
\cite{EngVat18} for some recent progress on minimizing the length of such a
string, including a contribution by a notorious hacker known as 4chan. (This
is no longer really about Eulerian circuits, since some amount of duplication
cannot be avoided in these strings.)

\subsubsection{Counting de Bruijn sequences}

Let us move in a different direction. Having proved the existence of de Bruijn
sequences in Theorem \ref{thm.debr.exist}, let us try to count them!

\textbf{Question.} Let $n$ and $k$ be two positive integers. Let $K$ be a
$k$-element set. How many de Bruijn sequences of order $n$ on $K$ are there?

To solve this, it makes sense to apply the BEST theorem to the digraph $D$ we
have constructed above. Alas, $D$ is not of the form $G^{\operatorname*{bidir}%
}$ for some undirected graph $G$, so we cannot apply the undirected MTT
(Matrix-Tree Theorem). However, $D$ is a balanced multidigraph, and for such
digraphs, a version of the undirected MTT still holds:

\begin{theorem}
[balanced Matrix-Tree Theorem]\label{thm.MTT.bal}Let $D=\left(  V,A,\psi
\right)  $ be a balanced multidigraph. Assume that $V=\left\{  1,2,\ldots
,n\right\}  $ for some positive integer $n$.

Let $L$ be the Laplacian of $D$. Then:

\begin{enumerate}
\item[\textbf{(a)}] For any vertex $r$ of $D$, we have%
\[
\left(  \text{\# of spanning arborescences of }D\text{ rooted to }r\right)
=\det\left(  L_{\sim r,\sim r}\right)  .
\]
Moreover, this number does not depend on $r$.

\item[\textbf{(b)}] Let $t$ be an indeterminate. Expand the determinant
$\det\left(  tI_{n}+L\right)  $ (here, $I_{n}$ denotes the $n\times n$
identity matrix) as a polynomial in $t$:%
\[
\det\left(  tI_{n}+L\right)  =c_{n}t^{n}+c_{n-1}t^{n-1}+\cdots+c_{1}%
t^{1}+c_{0}t^{0},
\]
where $c_{0},c_{1},\ldots,c_{n}$ are numbers. (Note that this is the
characteristic polynomial of $L$ up to substituting $-t$ for $t$ and
multiplying by a power of $-1$. Some of its coefficients are $c_{n}=1$ and
$c_{n-1}=\operatorname*{Tr}L$ and $c_{0}=\det L$.) Then, for any vertex $r$ of
$D$, we have%
\[
\left(  \text{\# of spanning arborescences of }D\text{ rooted to }r\right)
=\dfrac{1}{n}c_{1}.
\]

\item[\textbf{(c)}] Let $\lambda_{1},\lambda_{2},\ldots,\lambda_{n}$ be the
eigenvalues of $L$, listed in such a way that $\lambda_{n}=0$. Then, for any
vertex $r$ of $D$, we have%
\[
\left(  \text{\# of spanning arborescences of }D\text{ rooted to }r\right)
=\dfrac{1}{n}\cdot\lambda_{1}\lambda_{2}\cdots\lambda_{n-1}.
\]

\item[\textbf{(d)}] Let $\lambda_{1},\lambda_{2},\ldots,\lambda_{n}$ be the
eigenvalues of $L$, listed in such a way that $\lambda_{n}=0$. If all vertices
of $D$ have outdegree $>0$, then%
\[
\left(  \text{\# of Eulerian circuits of }D\right)  =\left\vert A\right\vert
\cdot\dfrac{1}{n}\cdot\lambda_{1}\lambda_{2}\cdots\lambda_{n-1}\cdot
\prod_{u\in V}\left(  \deg^{+}u-1\right)  !.
\]
(If you identify an Eulerian circuit with its cyclic rotations, then you
should drop the $\left\vert A\right\vert $ factor on the right hand side.)
\end{enumerate}
\end{theorem}

\begin{proof}
\textbf{(a)} The equality comes from the MTT (Theorem \ref{thm.MTT.MTT}). It
remains to prove that the \# of spanning arborescences of $D$ rooted to $r$
does not depend on $r$. But this is Corollary \ref{cor.BEST.tau=tau}.

\textbf{(b)} follows from \textbf{(a)} as in the undirected graph case (proof
of Theorem \ref{thm.MTT.undir} \textbf{(b)}).\footnote{In more detail: Just as
we proved in our above proof of Theorem \ref{thm.MTT.undir} (for the
undirected case), we have $c_{1}=\sum_{r=1}^{n}\det\left(  L_{\sim r,\sim
r}\right)  $. However, part \textbf{(a)} shows that the number $\det\left(
L_{\sim r,\sim r}\right)  $ does not depend on $r$. Thus, the sum $\sum
_{r=1}^{n}\det\left(  L_{\sim r,\sim r}\right)  $ consists of $n$ equal
addends, each of which can be written as $\det\left(  L_{\sim r,\sim
r}\right)  $ for any vertex $r$ of $D$. Therefore, this sum can be rewritten
as $n\cdot\det\left(  L_{\sim r,\sim r}\right)  $ for any vertex $r$ of $D$.
Hence, the equality $c_{1}=\sum_{r=1}^{n}\det\left(  L_{\sim r,\sim r}\right)
$ can be rewritten as $c_{1}=n\cdot\det\left(  L_{\sim r,\sim r}\right)  $ for
any vertex $r$ of $D$. Therefore, $\det\left(  L_{\sim r,\sim r}\right)
=\dfrac{1}{n}c_{1}$ for any vertex $r$ of $D$. But part \textbf{(a)} yields%
\[
\left(  \text{\# of spanning arborescences of }D\text{ rooted to }r\right)
=\det\left(  L_{\sim r,\sim r}\right)  =\dfrac{1}{n}c_{1}.
\]
}

\textbf{(c)} follows from \textbf{(b)} as in the undirected graph case (proof
of Theorem \ref{thm.MTT.undir} \textbf{(c)}).

\textbf{(d)} Assume that all vertices of $D$ have outdegree $>0$. Then,%
\begin{align*}
&  \left(  \text{\# of Eulerian circuits of }D\right) \\
&  =\sum_{a\in A}\left(  \text{\# of Eulerian circuits of }D\text{ whose first
arc is }a\right)  .
\end{align*}
However, if $a\in A$ is any arc, and if $r$ is the source of $a$, then%
\begin{align*}
&  \left(  \text{\# of Eulerian circuits of }D\text{ whose first arc is
}a\right) \\
&  =\left(  \text{\# of spanning arborescences of }D\text{ rooted to
}r\right)  \cdot\prod_{u\in V}\left(  \deg^{+}u-1\right)  !\\
&  \ \ \ \ \ \ \ \ \ \ \ \ \ \ \ \ \ \ \ \ \left(  \text{by the BEST' theorem
(Theorem \ref{thm.BEST.to})}\right) \\
&  =\dfrac{1}{n}\cdot\lambda_{1}\lambda_{2}\cdots\lambda_{n-1}\cdot\prod_{u\in
V}\left(  \deg^{+}u-1\right)  !\ \ \ \ \ \ \ \ \ \ \left(  \text{by part
\textbf{(c)}}\right)  .
\end{align*}
Hence,%
\begin{align*}
&  \left(  \text{\# of Eulerian circuits of }D\right) \\
&  =\sum_{a\in A}\underbrace{\left(  \text{\# of Eulerian circuits of }D\text{
whose first arc is }a\right)  }_{=\dfrac{1}{n}\cdot\lambda_{1}\lambda
_{2}\cdots\lambda_{n-1}\cdot\prod_{u\in V}\left(  \deg^{+}u-1\right)  !}\\
&  =\sum_{a\in A}\dfrac{1}{n}\cdot\lambda_{1}\lambda_{2}\cdots\lambda
_{n-1}\cdot\prod_{u\in V}\left(  \deg^{+}u-1\right)  !\\
&  =\left\vert A\right\vert \cdot\dfrac{1}{n}\cdot\lambda_{1}\lambda_{2}%
\cdots\lambda_{n-1}\cdot\prod_{u\in V}\left(  \deg^{+}u-1\right)  !.
\end{align*}
This proves part \textbf{(d)}.
\end{proof}

Now, let's try to solve our question -- i.e., let's count the de Bruijn
sequences of order $n$ on $K$.

Recall the digraph $D$ from our above proof of Theorem \ref{thm.debr.exist}.
We constructed a de Bruijn sequence of order $n$ on $K$ by finding an Eulerian
circuit of $D$. This actually works both ways: The map%
\begin{align*}
\left\{  \text{Eulerian circuits of }D\right\}   &  \rightarrow\left\{
\text{de Bruijn sequences of order }n\text{ on }K\right\}  ,\\
\mathbf{c}  &  \mapsto\left(  \text{the sequence of first entries of the arcs
of }\mathbf{c}\right)
\end{align*}
is a bijection (make sure you understand why!). Hence, by the bijection
principle, we have%
\begin{align}
&  \left(  \text{\# of de Bruijn sequences of order }n\text{ on }K\right)
\nonumber\\
&  =\left(  \text{\# of Eulerian circuits of }D\right)  .
\label{pf.thm.debr.num.1}%
\end{align}
By Theorem \ref{thm.MTT.bal} \textbf{(d)}, however, we have%
\begin{align}
&  \left(  \text{\# of Eulerian circuits of }D\right) \nonumber\\
&  =\left\vert K^{n}\right\vert \cdot\dfrac{1}{k^{n-1}}\cdot\lambda_{1}%
\lambda_{2}\cdots\lambda_{k^{n-1}-1}\cdot\prod_{u\in K^{n-1}}\left(  \deg
^{+}u-1\right)  !, \label{pf.thm.debr.num.2}%
\end{align}
where $\lambda_{1},\lambda_{2},\ldots,\lambda_{k^{n-1}}$ are the eigenvalues
of the Laplacian $L$ of $D$, indexed in such a way that $\lambda_{k^{n-1}}=0$.
(Note that the digraph $D=\left(  K^{n-1},K^{n},\psi\right)  $ has $k^{n-1}$
vertices, not $n$ vertices, so the \textquotedblleft$n$\textquotedblright\ in
Theorem \ref{thm.MTT.bal} is $k^{n-1}$ here.)

As we know, each vertex of $D$ has outdegree $k$. That is, we have $\deg
^{+}u=k$ for each $u\in K^{n-1}$. Thus,
\[
\prod_{u\in K^{n-1}}\left(  \deg^{+}u-1\right)  !=\prod_{u\in K^{n-1}}\left(
k-1\right)  !=\left(  \left(  k-1\right)  !\right)  ^{k^{n-1}}.
\]
Also,%
\[
\left\vert K^{n}\right\vert \cdot\dfrac{1}{k^{n-1}}=k^{n}\cdot\dfrac
{1}{k^{n-1}}=k.
\]
It remains to find $\lambda_{1}\lambda_{2}\cdots\lambda_{k^{n-1}-1}$. What are
the eigenvalues of $L$ ?

The Laplacian $L$ of our digraph $D$ is a $k^{n-1}\times k^{n-1}$-matrix whose
rows and columns are indexed by $\left(  n-1\right)  $-tuples in $K^{n-1}$.
Strictly speaking, we should relabel the vertices of $D$ as $1,2,\ldots
,k^{n-1}$ here, in order to have a \textquotedblleft proper
matrix\textquotedblright\ with a well-defined order on its rows and columns.
But let's not do this; instead, I trust you can do the relabeling yourself, or
just use the more general notion of matrices that allows for the rows and the
columns to be indexed by arbitrary things (see
\url{https://mathoverflow.net/questions/317105} for details).

Let $C$ be the adjacency matrix of the digraph $D$; this is the $k^{n-1}\times
k^{n-1}$-matrix (again with rows and columns indexed by $\left(  n-1\right)
$-tuples in $K^{n-1}$) whose $\left(  i,j\right)  $-th entry is the \# of arcs
with source $i$ and target $j$. In particular, the trace of $C$ is thus the \#
of loops of $D$. It is easy to see that the loops of $D$ are precisely the
arcs of the form $\left(  x,x,\ldots,x\right)  \in K^{n}$ for $x\in K$; thus,
$D$ has exactly $k$ loops. Hence, the trace of $C$ is $k$.

Recall the definition of the Laplacian matrix $L$. We can restate it as
follows:%
\begin{equation}
L=\Delta-C, \label{pf.thm.debr.num.L=1}%
\end{equation}
where $\Delta$ is the diagonal matrix whose diagonal entries are the
outdegrees of the vertices of $D$. Since each vertex of $D$ has outdegree $k$,
the latter diagonal matrix $\Delta$ is simply $k\cdot I$, where $I$ is the
identity matrix (of the appropriate size). Hence, (\ref{pf.thm.debr.num.L=1})
can be rewritten as%
\[
L=k\cdot I-C.
\]
Thus, if $\gamma_{1},\gamma_{2},\ldots,\gamma_{k^{n-1}}$ are the eigenvalues
of $C$, then $k-\gamma_{1},k-\gamma_{2},\ldots,k-\gamma_{k^{n-1}}$ are the
eigenvalues of $L$. Computing the former will thus help us find the latter.

Furthermore, let $J$ be the $k^{n-1}\times k^{n-1}$-matrix (again with rows
and columns indexed by $\left(  n-1\right)  $-tuples in $K^{n-1}$) whose all
entries are $1$. It is easy to see that the eigenvalues of $J$ are%
\[
\underbrace{0,0,\ldots,0}_{k^{n-1}-1\text{ many zeroes}},k^{n-1}.
\]
(The easiest way to see this is by noticing that $J$ has rank $1$ and trace
$k^{n-1}$.\ \ \ \ \footnote{Here are the details: The matrix $J$ has rank $1$
(since all its rows are the same); thus, all but one of its eigenvalues are
$0$. It remains to show that the remaining eigenvalue is $k^{n-1}$. However,
it is known that the sum of the eigenvalues of a square matrix equals its
trace. Thus, if all but one of the eigenvalues of a square matrix are $0$,
then the remaining eigenvalue equals its trace. Applying this to our matrix
$J$, we see that its remaining eigenvalue equals its trace, which is $k^{n-1}%
$.})

Now, here is something really underhanded: We observe that%
\[
C^{n-1}=J.
\]

[\textit{Proof:} We need to show that all entries of the matrix $C^{n-1}$ are
$1$. So let $i$ and $j$ be two vertices of $D$. We must then show that the
$\left(  i,j\right)  $-th entry of $C^{n-1}$ is $1$.

Recall the combinatorial interpretation of the powers of an adjacency matrix
(Theorem \ref{thm.adjmat.walks}): For any $\ell\in\mathbb{N}$, the $\left(
i,j\right)  $-th entry of $C^{\ell}$ is the \# of walks from $i$ to $j$ (in
$D$) that have length $\ell$. Thus, in particular, the $\left(  i,j\right)
$-th entry of $C^{n-1}$ is the \# of walks from $i$ to $j$ (in $D$) that have
length $n-1$. But this number is actually $1$, as we have already shown in our
above proof of Theorem \ref{thm.debr.exist}. This completes the proof of
$C^{n-1}=J$.] \medskip

How does this help us compute the eigenvalues of $C$ ? Well, let $\gamma
_{1},\gamma_{2},\ldots,\gamma_{k^{n-1}}$ be the eigenvalues of $C$. Then, for
any $\ell\in\mathbb{N}$, the eigenvalues of $C^{\ell}$ are $\gamma_{1}^{\ell
},\gamma_{2}^{\ell},\ldots,\gamma_{k^{n-1}}^{\ell}$ (this is a fact that holds
for any square matrix, and is probably easiest to prove using the Jordan
canonical form or triangularization\footnote{In fact this is a particular case
of the Spectral Mapping Theorem (\cite[Chapter 9, Theorem 2.1]{Treil17}).}).
Hence, in particular, $\gamma_{1}^{n-1},\gamma_{2}^{n-1},\ldots,\gamma
_{k^{n-1}}^{n-1}$ are the eigenvalues of $C^{n-1}=J$; but we know that the
latter eigenvalues are $\underbrace{0,0,\ldots,0}_{k^{n-1}-1\text{ many
zeroes}},k^{n-1}$. Hence, all but one of the $k^{n-1}$ numbers $\gamma
_{1}^{n-1},\gamma_{2}^{n-1},\ldots,\gamma_{k^{n-1}}^{n-1}$ equal $0$. Thus,
all but one of the $k^{n-1}$ numbers $\gamma_{1},\gamma_{2},\ldots
,\gamma_{k^{n-1}}$ equal $0$ (we don't know what the remaining number is,
since $\left(  n-1\right)  $-st roots are not uniquely determined in
$\mathbb{C}$). In other words, all but one of the eigenvalues of $C$ equal
$0$. The remaining eigenvalue must thus be the trace of $C$ (because the sum
of the eigenvalues of a square matrix is known to be the trace of that
matrix), and therefore equal $k$ (since we know that the trace of $C$ is $k$).

So we have shown that the eigenvalues of $C$ are $\underbrace{0,0,\ldots
,0}_{k^{n-1}-1\text{ many zeroes}},k$. Thus, the eigenvalues of $L$ are%
\[
\underbrace{k-0,\ \ k-0,\ \ \ldots,\ \ k-0}_{k^{n-1}-1\text{ many }\left(
k-0\right)  \text{'s}},\ \ k-k
\]
(because if $\gamma_{1},\gamma_{2},\ldots,\gamma_{k^{n-1}}$ are the
eigenvalues of $C$, then $k-\gamma_{1},k-\gamma_{2},\ldots,k-\gamma_{k^{n-1}}$
are the eigenvalues of $L$). In other words, the eigenvalues of $L$ are
\[
\underbrace{k,\ \ k,\ \ \ldots,\ \ k}_{k^{n-1}-1\text{ many }k\text{'s}%
},\ \ 0.
\]
Hence, the eigenvalues $\lambda_{1},\lambda_{2},\ldots,\lambda_{k^{n-1}-1}$ in
(\ref{pf.thm.debr.num.2}) all equal $k$. Thus, (\ref{pf.thm.debr.num.2})
simplifies to%
\begin{align*}
&  \left(  \text{\# of Eulerian circuits of }D\right) \\
&  =\underbrace{\left\vert K^{n}\right\vert \cdot\dfrac{1}{k^{n-1}}%
}_{\substack{=k^{n}\cdot\dfrac{1}{k^{n-1}}\\=k}}\cdot\underbrace{kk\cdots
k}_{k^{n-1}-1\text{ factors}}\cdot\underbrace{\prod_{u\in K^{n-1}}\left(
\deg^{+}u-1\right)  !}_{=\left(  \left(  k-1\right)  !\right)  ^{k^{n-1}}}\\
&  =\underbrace{k\cdot\underbrace{kk\cdots k}_{k^{n-1}-1\text{ factors}}%
}_{=k^{k^{n-1}}}\cdot\left(  \left(  k-1\right)  !\right)  ^{k^{n-1}%
}=k^{k^{n-1}}\cdot\left(  \left(  k-1\right)  !\right)  ^{k^{n-1}}\\
&  =\left(  \underbrace{k\cdot\left(  k-1\right)  !}_{=k!}\right)  ^{k^{n-1}%
}=k!^{k^{n-1}}.
\end{align*}
In view of this, we can rewrite (\ref{pf.thm.debr.num.1}) as%
\[
\left(  \text{\# of de Bruijn sequences of order }n\text{ on }K\right)
=k!^{k^{n-1}}.
\]
Thus, we have proved the following:

\begin{theorem}
\label{thm.debr.num}Let $n$ and $k$ be positive integers. Let $K$ be a
$k$-element set. Then,%
\[
\left(  \text{\# of de Bruijn sequences of order }n\text{ on }K\right)
=k!^{k^{n-1}}.
\]

\end{theorem}

What a nice (and huge) answer!

Our above proof of Theorem \ref{thm.debr.num} is essentially taken from
\cite[Chapter 10]{Stanley-AC}.

We note that a combinatorial proof of Theorem \ref{thm.debr.num} (avoiding any
use of linear algebra) has been recently given in \cite{BidKis02}.

\subsection{More on Laplacians}

Much more can be said about the Laplacian of a digraph. The study of matrices
associated to a graph or digraph is known as \textbf{spectral graph theory};
I'd say the Laplacian is probably the most prominent of these matrices (even
though the adjacency matrix is somewhat easier to define). The original form
of the matrix-tree theorem (actually a subtler variant of Theorem
\ref{thm.MTT.undir} \textbf{(a)}) was found by Gustav Kirchhoff in his study
of electricity \cite{Kirchh47} (see \cite[\S 2.1.1]{Holzer22} for a modern
exposition); the effective resistance between two nodes of an electrical
network is a ratio of spanning-tree counts and thus can be computed using the
Laplacian (see, e.g., \cite[\S 2 and \S 3]{Vos16}). To be more precise, this
relies on a \textquotedblleft weighted count\textquotedblright\ of spanning
trees, which is more general than the counting we have done so far; we will
learn about it in the next section.

Another application of Laplacians is to drawing graphs: see \textquotedblleft%
\href{https://en.wikipedia.org/wiki/Spectral_layout}{spectral layout}%
\textquotedblright\ or \textquotedblleft spectral graph
drawing\textquotedblright\ (e.g., \cite{Gallie13}).

\subsection{On the left nullspace of the Laplacian}

We shall now answer another natural question about Laplacians of digraphs.
Recall that the Laplacian $L$ of a digraph $D$ always satisfies $Le=0$, where
$e=\left(
\begin{array}
[c]{c}%
1\\
1\\
\vdots\\
1
\end{array}
\right)  $. Thus, the vector $e$ belongs to the right nullspace (= right
kernel) of $L$. It is not hard to see that if $D$ has a to-root and we are
working over a characteristic-$0$ field, then $e$ spans this nullspace, i.e.,
there are no vectors in that nullspace other than scalar multiples of $e$.
(This is actually an \textquotedblleft if and only if\textquotedblright.) What
about the left nullspace of $L$ ? Can we explicitly find a nonzero vector $f$
with $fL=0$ ? The answer is positive:

\begin{theorem}
[harmonic vector theorem for Laplacians]\label{thm.MTT.harm}Let $D=\left(
V,A,\psi\right)  $ be a multidigraph, where $V=\left\{  1,2,\ldots,n\right\}
$ for some $n\in\mathbb{N}$.

For each $r\in V$, let $\tau\left(  D,r\right)  $ be the \# of spanning
arborescences of $D$ rooted to $r$.

Let $f$ be the row vector $\left(  \tau\left(  D,1\right)  ,\ \tau\left(
D,2\right)  ,\ \ldots,\ \tau\left(  D,n\right)  \right)  $. Then, $fL=0$.
\end{theorem}

Before we prove this, some remarks are in order. Theorem \ref{thm.MTT.harm}
(or, more precisely, its weighted version, which we will see in the next
section) can be used to explicitly compute the steady state of a Markov chain
(see \cite{KrGrWi10} and also Corollary \ref{cor.MTT.markov} below); a similar
interpretation, but in economical terms (emergence of money in a barter
economy), appears in \cite[\S 1]{Sahi14}.

\begin{remark}
\label{rmk.MTT.harm.zero}The row vector $f$ in Theorem \ref{thm.MTT.harm} can
be the zero vector. Specifically, it will be the zero vector whenever $D$ has
no to-root. However, even in such cases, there exist nonzero vectors
$f^{\prime}$ such that $f^{\prime}L=0$, as long as $n>0$. Such vectors can be
constructed as follows: Pick a sink component $C$ of $D$ (see Definition
\ref{def.mdg.sinkcomp} for the meaning of \textquotedblleft sink
component\textquotedblright, and Theorem \ref{thm.mdg.sinkcomp} for its
existence). Then, the induced subdigraph $D\left[  C\right]  $ is strongly
connected (by Proposition \ref{prop.mdg.strong-component.ind-connected}). For
each $r\in C$, we let $\tau\left(  D\left[  C\right]  ,r\right)  $ be the \#
of spanning arborescences of $D\left[  C\right]  $ rooted to $r$; this is
nonzero because $D\left[  C\right]  $ is strongly connected. Also set
$\tau\left(  D\left[  C\right]  ,r\right)  :=0$ for all $r\in V\setminus C$.
Let $f^{\prime}$ be the row vector $\left(  \tau\left(  D\left[  C\right]
,1\right)  ,\ \tau\left(  D\left[  C\right]  ,2\right)  ,\ \ldots
,\ \tau\left(  D\left[  C\right]  ,n\right)  \right)  $. Then, $f^{\prime}%
L=0$. This is easy to show by applying Theorem \ref{thm.MTT.harm} to $D\left[
C\right]  $ instead of $D$.

Note that each sink component $C$ of $D$ yields a different nonzero row vector
$f^{\prime}$; and all these row vectors $f^{\prime}$ are linearly independent.
Hence, the dimension of the left nullspace of $L$ is at least the \# of sink
components of $D$. It can be proved that this is, in fact, an equality.
\end{remark}

We shall give a proof of Theorem \ref{thm.MTT.harm} based upon two lemmas. The
first lemma is a general linear-algebraic result:

\begin{lemma}
\label{lem.sol.7.1.matrix}Let $B$ be an $n\times n$-matrix over an arbitrary
commutative ring $\mathbb{K}$. (For example, $\mathbb{K}$ can be $\mathbb{R}$,
in which case $B$ is a real matrix.) Assume that the sum of all columns of $B$
is the zero vector. Then, for any $r,s,t\in\left\{  1,2,\ldots,n\right\}  $,
we have%
\[
\det\left(  B_{\sim r,\sim t}\right)  =\left(  -1\right)  ^{s-t}\det\left(
B_{\sim r,\sim s}\right)  .
\]

\end{lemma}

\begin{proof}
[Proof of Lemma \ref{lem.sol.7.1.matrix}.]There are various ways to prove
this, but here is probably the most elegant one:

We WLOG assume that $s\neq t$, since otherwise the claim is obvious. Let us
now change the $r$-th row of the matrix $B$ as follows:

\begin{itemize}
\item We replace the $s$-th entry of the $r$-th row by $1$.

\item We replace the $t$-th entry of the $r$-th row by $-1$.

\item We replace all other entries of the $r$-th row by $0$.
\end{itemize}

Let $C$ be the resulting $n\times n$-matrix.\footnote{For example, if $n=4$
and $B=\left(
\begin{array}
[c]{cccc}%
a & b & c & d\\
a^{\prime} & b^{\prime} & c^{\prime} & d^{\prime}\\
a^{\prime\prime} & b^{\prime\prime} & c^{\prime\prime} & d^{\prime\prime}\\
a^{\prime\prime\prime} & b^{\prime\prime\prime} & c^{\prime\prime\prime} &
d^{\prime\prime\prime}%
\end{array}
\right)  $ and $s=1$ and $t=3$ and $r=2$, then $C=\left(
\begin{array}
[c]{cccc}%
a & b & c & d\\
1 & 0 & -1 & 0\\
a^{\prime\prime} & b^{\prime\prime} & c^{\prime\prime} & d^{\prime\prime}\\
a^{\prime\prime\prime} & b^{\prime\prime\prime} & c^{\prime\prime\prime} &
d^{\prime\prime\prime}%
\end{array}
\right)  $.} Thus, $C$ agrees with $B$ in all rows other than the $r$-th one.
Hence, in particular,%
\begin{equation}
C_{\sim r,\sim k}=B_{\sim r,\sim k}\ \ \ \ \ \ \ \ \ \ \text{for each }%
k\in\left\{  1,2,\ldots,n\right\}  . \label{pf.lem.sol.7.1.matrix.2}%
\end{equation}
Note also that the only nonzero entries in the $r$-th row of $C$
are\footnote{We are using the notation $C_{r,k}$ for the entry of $C$ in the
$r$-th row and the $k$-th column.} $C_{r,s}=1$ and $C_{r,t}=-1$. Hence, the
entries in the $r$-th row of $C$ add up to $0$.

Recall that the sum of all columns of $B$ is the zero vector. In other words,
in each row of $B$, the entries add up to $0$. The matrix $C$ therefore also
has this property (because the only row of $C$ that differs from the
corresponding row of $B$ is the $r$-th row; however, we have shown above that
in the $r$-th row, the entries of $C$ also add up to $0$). In other words, the
sum of all columns of $C$ is the zero vector. This easily entails that $\det
C=0$\ \ \ \ \footnote{\textit{Proof.} It is well-known that the determinant of
a matrix does not change if we add a column to another. Hence, the determinant
of $C$ will not change if we add each column of $C$ other than the first one
to the first column of $C$. However, the result of this operation will be a
matrix whose first column is $0$ (since the sum of all columns of $C$ is the
zero vector), and therefore this matrix will have determinant $0$. Since the
operation did not change the determinant, we thus conclude that the
determinant of $C$ was $0$. In other words, $\det C=0$.}.

On the other hand, Laplace expansion along the $r$-th row yields%
\begin{align*}
\det C  &  =\sum_{k=1}^{n}\left(  -1\right)  ^{r+k}C_{r,k}\det\left(  C_{\sim
r,\sim k}\right) \\
&  =\left(  -1\right)  ^{r+s}1\det\left(  C_{\sim r,\sim s}\right)  +\left(
-1\right)  ^{r+t}\left(  -1\right)  \det\left(  C_{\sim r,\sim t}\right)
\end{align*}
(since the only nonzero entries $C_{r,k}$ in the $r$-th row of $C$ are
$C_{r,s}=1$ and $C_{r,t}=-1$). Comparing this with $\det C=0$, we obtain%
\begin{align*}
0  &  =\left(  -1\right)  ^{r+s}1\det\left(  C_{\sim r,\sim s}\right)
+\left(  -1\right)  ^{r+t}\left(  -1\right)  \det\left(  C_{\sim r,\sim
t}\right) \\
&  =\left(  -1\right)  ^{r+s}\det\underbrace{\left(  C_{\sim r,\sim s}\right)
}_{\substack{=B_{\sim r,\sim s}\\\text{(by (\ref{pf.lem.sol.7.1.matrix.2}))}%
}}-\left(  -1\right)  ^{r+t}\det\underbrace{\left(  C_{\sim r,\sim t}\right)
}_{\substack{=B_{\sim r,\sim t}\\\text{(by (\ref{pf.lem.sol.7.1.matrix.2}))}%
}}\\
&  =\left(  -1\right)  ^{r+s}\det\left(  B_{\sim r,\sim s}\right)  -\left(
-1\right)  ^{r+t}\det\left(  B_{\sim r,\sim t}\right)  .
\end{align*}
In other words, $\left(  -1\right)  ^{r+t}\det\left(  B_{\sim r,\sim
t}\right)  =\left(  -1\right)  ^{r+s}\det\left(  B_{\sim r,\sim s}\right)  $.
Dividing both sides of this by $\left(  -1\right)  ^{r+t}$, we obtain
$\det\left(  B_{\sim r,\sim t}\right)  =\left(  -1\right)  ^{s-t}\det\left(
B_{\sim r,\sim s}\right)  $. This proves Lemma \ref{lem.sol.7.1.matrix}.
\end{proof}

Our next lemma is the following generalization of Theorem \ref{thm.MTT.MTT}:

\begin{theorem}
[Matrix-Tree Theorem, off-diagonal version]\label{thm.MTT.MTTrs}Let $D=\left(
V,A,\psi\right)  $ be a multidigraph. Assume that $V=\left\{  1,2,\ldots
,n\right\}  $ for some positive integer $n$.

Let $L$ be the Laplacian of $D$. Let $r$ and $s$ be two vertices of $D$. Then,%
\[
\left(  \#\text{ of spanning arborescences of }D\text{ rooted to }r\right)
=\left(  -1\right)  ^{r+s}\det\left(  L_{\sim r,\sim s}\right)  .
\]

\end{theorem}

Note that Theorem \ref{thm.MTT.MTT} is the particular case of Theorem
\ref{thm.MTT.MTTrs} for $s=r$. Fortunately, using Lemma
\ref{lem.sol.7.1.matrix}, we can easily derive the general case from the particular:

\begin{proof}
[Proof of Theorem \ref{thm.MTT.MTTrs}.]We have seen (in the proof of
Proposition \ref{prop.MTT.detL=0}) that the sum of all columns of the
Laplacian $L$ is the zero vector. Hence, Lemma \ref{lem.sol.7.1.matrix}
(applied to $\mathbb{K}=\mathbb{Q}$ and $B=L$ and $t=r$) yields%
\[
\det\left(  L_{\sim r,\sim r}\right)  =\underbrace{\left(  -1\right)  ^{s-r}%
}_{=\left(  -1\right)  ^{r+s}}\det\left(  L_{\sim r,\sim s}\right)  =\left(
-1\right)  ^{r+s}\det\left(  L_{\sim r,\sim s}\right)  .
\]
However, the Matrix-Tree Theorem (Theorem \ref{thm.MTT.MTT}) yields%
\begin{align*}
\left(  \#\text{ of spanning arborescences of }D\text{ rooted to }r\right)
&  =\det\left(  L_{\sim r,\sim r}\right) \\
&  =\left(  -1\right)  ^{r+s}\det\left(  L_{\sim r,\sim s}\right)  .
\end{align*}
This proves Theorem \ref{thm.MTT.MTTrs}.
\end{proof}

We are now ready to prove Theorem \ref{thm.MTT.harm}:

\begin{proof}
[Proof of Theorem \ref{thm.MTT.harm}.]For each $r,s\in\left\{  1,2,\ldots
,n\right\}  $, we have%
\begin{align}
\tau\left(  D,r\right)   &  =\left(  \#\text{ of spanning arborescences of
}D\text{ rooted to }r\right) \nonumber\\
&  \ \ \ \ \ \ \ \ \ \ \ \ \ \ \ \ \ \ \ \ \left(  \text{by the definition of
}\tau\left(  D,r\right)  \right) \nonumber\\
&  =\left(  -1\right)  ^{r+s}\det\left(  L_{\sim r,\sim s}\right)
\label{sol.7.1.b.1}%
\end{align}
(by Theorem \ref{thm.MTT.MTTrs}).

However, we have $f=\left(  \tau\left(  D,1\right)  ,\tau\left(  D,2\right)
,\ldots,\tau\left(  D,n\right)  \right)  $. Thus, for each $s\in\left\{
1,2,\ldots,n\right\}  $, the $s$-th entry of the column vector $fL$
is\footnote{We are using the notation $L_{r,s}$ for the entry of the matrix
$L$ in the $r$-th row and the $s$-th column.}%
\begin{align*}
&  \sum_{r=1}^{n}\underbrace{\tau\left(  D,r\right)  }_{\substack{=\left(
-1\right)  ^{r+s}\det\left(  L_{\sim r,\sim s}\right)  \\\text{(by
(\ref{sol.7.1.b.1}))}}}L_{r,s}\\
&  =\sum_{r=1}^{n}\left(  -1\right)  ^{r+s}\det\left(  L_{\sim r,\sim
s}\right)  L_{r,s}\\
&  =\sum_{r=1}^{n}\left(  -1\right)  ^{r+s}L_{r,s}\det\left(  L_{\sim r,\sim
s}\right)  =\det L\\
&  \ \ \ \ \ \ \ \ \ \ \ \ \ \ \ \ \ \ \ \ \left(
\begin{array}
[c]{c}%
\text{since Laplace expansion along the }s\text{-th column}\\
\text{yields }\det L=\sum_{r=1}^{n}\left(  -1\right)  ^{r+s}L_{r,s}\det\left(
L_{\sim r,\sim s}\right)
\end{array}
\right) \\
&  =0
\end{align*}
(by Proposition \ref{prop.MTT.detL=0}). This shows that all entries of $fL$
are $0$. In other words, $fL=0$. Theorem \ref{thm.MTT.harm} is thus proved.
\end{proof}

Other proofs of Theorem \ref{thm.MTT.harm} exist. In particular, a
combinatorial proof is sketched in \cite[Theorem 1]{Sahi14}. (More precisely,
\cite[Theorem 1]{Sahi14} in this paper is the claim of Theorem
\ref{thm.MTT.harm} upon reversing all the arcs and replacing all matrices by
their transposes.)\footnote{I tried to explain this proof in more detail in
\href{https://www.cip.ifi.lmu.de/~grinberg/t/18s/mt3s.pdf}{the solutions to
Spring 2018 Math 4707 midterm \#3} -- see the proof of Theorem 0.7 in those
solutions; you be the judge if I succeeded.}

\subsection{A weighted Matrix-Tree Theorem}

\subsubsection{Definitions}

We have so far been \textbf{counting} arborescences. A natural generalization
of counting is \textbf{weighted counting} -- i.e., you assign a certain number
(a \textquotedblleft weight\textquotedblright) to each arborescence (or
whatever object you are interested in), and then you \textbf{sum} the weights
of all arborescences (instead of merely counting them). This generalizes
counting, because if all weights are $1$, then you get the \# of arborescences.

If you pick the weights to be completely random, then the sum won't usually be
particularly interesting. However, some choices of weights lead to good
behavior. Let us see what we get if we assign a weight to each \textbf{arc} of
our digraph, and then define the weight of an arborescence to be the
\textbf{product} of the weights of the arcs that appear in this arborescence.

\begin{definition}
\label{def.wMTT.wt-defs}Let $D=\left(  V,A,\psi\right)  $ be a multidigraph.

Let $\mathbb{K}$ be a commutative ring. Assume that an element $w_{a}%
\in\mathbb{K}$ is assigned to each arc $a\in A$. We call this $w_{a}$ the
\textbf{weight} of the arc $a$. (You can assume that $\mathbb{K}=\mathbb{R}$,
so that the weights are just numbers.)

\begin{enumerate}
\item[\textbf{(a)}] For any two vertices $i,j\in V$, we let $a_{i,j}^{w}$ be
the sum of the weights of all arcs of $D$ that have source $i$ and target $j$.

\item[\textbf{(b)}] For any vertex $i\in V$, we define the \textbf{weighted
outdegree} $\deg^{+w}i$ of $i$ to be the sum%
\[
\sum_{\substack{a\in A;\\\text{the source of }a\text{ is }i}}w_{a}.
\]

\item[\textbf{(c)}] If $B$ is a subdigraph of $D$, then the \textbf{weight}
$w\left(  B\right)  $ of $B$ is defined to be the product $\prod
\limits_{a\text{ is an arc of }B}w_{a}$. This is the product of the weights of
all arcs of $B$.

\item[\textbf{(d)}] Assume that $V=\left\{  1,2,\ldots,n\right\}  $ for some
$n\in\mathbb{N}$. The \textbf{weighted Laplacian} of $D$ (with respect to the
weights $w_{a}$) is defined to be the $n\times n$-matrix $L^{w}\in
\mathbb{K}^{n\times n}$ (note that the \textquotedblleft$w$\textquotedblright%
\ here is a superscript, not an exponent) whose entries are given by%
\[
L_{i,j}^{w}=\left(  \deg^{+w}i\right)  \cdot\left[  i=j\right]  -a_{i,j}%
^{w}\ \ \ \ \ \ \ \ \ \ \text{for all }i,j\in V.
\]

\end{enumerate}
\end{definition}

These definitions generalize analogous definitions in the \textquotedblleft
unweighted case\textquotedblright. Indeed, if we take all the arc weights
$w_{a}$ to be $1$, then the weighted outdegree $\deg^{+w}i$ of a vertex $i$
becomes its usual outdegree $\deg i$, and the weighted Laplacian $L^{w}$
becomes the usual Laplacian $L$. The weight $w\left(  B\right)  $ of a
subdigraph $B$ simply becomes $1$ in this case.

\subsubsection{The weighted Matrix-Tree Theorem}

We now can generalize the original MTT (= Matrix-Tree Theorem)\footnote{To
remind: The original MTT is Theorem \ref{thm.MTT.MTT}.} as follows:

\Needspace{20pc}

\begin{theorem}
[weighted Matrix-Tree Theorem]\label{thm.MTT.wMTT}Let $D=\left(
V,A,\psi\right)  $ be a multidigraph.

Let $\mathbb{K}$ be a commutative ring. Assume that an element $w_{a}%
\in\mathbb{K}$ is assigned to each arc $a\in A$. We call this $w_{a}$ the
\textbf{weight} of the arc $a$.

Assume that $V=\left\{  1,2,\ldots,n\right\}  $ for some $n\in\mathbb{N}$. Let
$L^{w}$ be the weighted Laplacian of $D$.

Let $r$ be a vertex of $D$. Then,%
\[
\sum_{\substack{B\text{ is a spanning}\\\text{arborescence}\\\text{of }D\text{
rooted to }r}}w\left(  B\right)  =\det\left(  L_{\sim r,\sim r}^{w}\right)  .
\]

\end{theorem}

\begin{example}
\label{exa.MTT.wMTT.1}Let $D$ be the following multidigraph:%
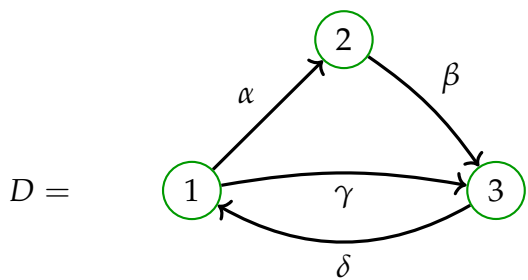
\[%
%
\ \ ,\ \ \ \ \ \ \ \ \ \ \text{and let }r=3.
\]
Then, $D$ has two spanning arborescences rooted to $r$. One of the two has
arcs $\alpha$ and $\beta$ (and thus has weight $w_{\alpha}w_{\beta}$); the
other has arcs $\gamma$ and $\beta$ (and thus has weight $w_{\gamma}w_{\beta}%
$). Hence,%
\begin{equation}
\sum_{\substack{B\text{ is a spanning}\\\text{arborescence}\\\text{of }D\text{
rooted to }r}}w\left(  B\right)  =w_{\alpha}w_{\beta}+w_{\gamma}w_{\beta},
\label{eq.exa.MTT.wMTT.1.1}%
\end{equation}
The weighted Laplacian $L^{w}$ is%
\[
L^{w}=\left(
\begin{array}
[c]{ccc}%
w_{\alpha}+w_{\gamma} & -w_{\alpha} & -w_{\gamma}\\
0 & w_{\beta} & -w_{\beta}\\
-w_{\delta} & 0 & w_{\delta}%
\end{array}
\right)
\]
(since, for example, $\deg^{+w}1=w_{\alpha}+w_{\gamma}$ and $a_{1,1}^{w}=0$
and $a_{1,2}^{w}=w_{\alpha}$). Thus,%
\begin{align*}
L_{\sim3,\sim3}^{w}  &  =\left(
\begin{array}
[c]{cc}%
w_{\alpha}+w_{\gamma} & -w_{\alpha}\\
0 & w_{\beta}%
\end{array}
\right)  \ \ \ \ \ \ \ \ \ \ \text{and therefore}\\
\det\left(  L_{\sim3,\sim3}^{w}\right)   &  =\left(  w_{\alpha}+w_{\gamma
}\right)  w_{\beta}=w_{\alpha}w_{\beta}+w_{\gamma}w_{\beta}.
\end{align*}
The right hand side of this agrees with that of (\ref{eq.exa.MTT.wMTT.1.1}).
This confirms the weighted MTT for our $D$ and $r$.
\end{example}

As we already said, the weighted MTT generalizes the original MTT, because if
we take all $w_{a}$'s to be $1$, we just recover the original MTT.

However, we can also go backwards: we can derive the weighted MTT from the
original MTT. Let us do this.

\subsubsection{The polynomial identity trick}

First, we recall a standard result in algebra, known as the \textbf{principle
of permanence of polynomial identities} or as the \textbf{polynomial identity
trick} (it also goes under several other names). Here is one incarnation of
this principle:

\begin{theorem}
[principle of permanence of polynomial identities]\label{thm.polyIDtrick}Let
$P\left(  x_{1},x_{2},\ldots,x_{m}\right)  $ and $Q\left(  x_{1},x_{2}%
,\ldots,x_{m}\right)  $ be two polynomials with integer coefficients in
several indeterminates $x_{1},x_{2},\ldots,x_{m}$. Assume that the equality%
\begin{equation}
P\left(  k_{1},k_{2},\ldots,k_{m}\right)  =Q\left(  k_{1},k_{2},\ldots
,k_{m}\right)  \label{eq.thm.polyIDtrick.eq}%
\end{equation}
holds for every $m$-tuple $\left(  k_{1},k_{2},\ldots,k_{m}\right)
\in\mathbb{N}^{m}$ of nonnegative integers. Then, $P\left(  x_{1},x_{2}%
,\ldots,x_{m}\right)  $ and $Q\left(  x_{1},x_{2},\ldots,x_{m}\right)  $ are
identical as polynomials (so that, in particular, the equality
(\ref{eq.thm.polyIDtrick.eq}) holds not only for every $\left(  k_{1}%
,k_{2},\ldots,k_{m}\right)  \in\mathbb{N}^{m}$, but also for every $\left(
k_{1},k_{2},\ldots,k_{m}\right)  \in\mathbb{C}^{m}$, and more generally, for
every $\left(  k_{1},k_{2},\ldots,k_{m}\right)  \in\mathbb{K}^{m}$ where
$\mathbb{K}$ is an arbitrary commutative ring).
\end{theorem}

Theorem \ref{thm.polyIDtrick} is often summarized as \textquotedblleft in
order to prove that two polynomials are equal, it suffices to show that they
are equal on all nonnegative integer points\textquotedblright\ (where a
\textquotedblleft nonnegative integer point\textquotedblright\ means a point
-- i.e., a tuple of inputs -- whose all entries are nonnegative integers).
Even shorter, one says that \textquotedblleft a polynomial identity (i.e., an
equality between two polynomials) needs only to be checked on nonnegative
integers\textquotedblright. For example, if you can prove the equality
\[
\left(  x+y\right)  ^{4}+\left(  x-y\right)  ^{4}=2x^{4}+12x^{2}y^{2}+2y^{4}%
\]
for all nonnegative integers $x$ and $y$, then you automatically conclude that
this equality holds as a polynomial identity, and thus is true for any
elements $x$ and $y$ of a commutative ring.

A typical application of Theorem \ref{thm.polyIDtrick} is to argue that a
polynomial identity you have proved for all nonnegative integers must
automatically hold for all inputs (because of Theorem \ref{thm.polyIDtrick}).
Some examples of such reasoning can be found in \cite[\S 2.6.3 and
\S 2.6.4]{19fco} and in \cite[\S 7.5.3]{20f}. A variant of Theorem
\ref{thm.polyIDtrick} is \cite[Theorem 2.6]{Conrad-univid}; actually, the
proof of \cite[Theorem 2.6]{Conrad-univid} can be trivially adapted to prove
Theorem \ref{thm.polyIDtrick} (just replace \textquotedblleft nonempty open
set in $\mathbb{C}^{k}$\textquotedblright\ by \textquotedblleft$\mathbb{N}%
^{k}$\textquotedblright). In truth, there is nothing special about nonnegative
integers and the set $\mathbb{N}$; you could replace $\mathbb{N}$ by any
infinite set of numbers (or even any sufficiently large set of numbers, where
\textquotedblleft sufficiently large\textquotedblright\ means
\textquotedblleft more than $\max\left\{  \deg P,\deg Q\right\}  $
many\textquotedblright). See \cite[Lemma 2.1]{Alon} or \cite[Theorem
2]{Lason10} for a fairly general version of Theorem \ref{thm.polyIDtrick} that
includes such cases\footnote{To be precise, \cite[Lemma 2.1]{Alon} and
\cite[Theorem 2]{Lason10} are not concerned with two polynomials being
identical, but rather with one polynomial being identically zero. But this is
an equivalent question: Two polynomials $P$ and $Q$ are identical if and only
if their difference $P-Q$ is identically zero.}.

\subsubsection{Proof of the weighted MTT}

We can now deduce the weighted MTT from the original MTT (Theorem
\ref{thm.MTT.MTT}):

\begin{proof}
[Proof of Theorem \ref{thm.MTT.wMTT}.]The claim of Theorem \ref{thm.MTT.wMTT}
(for fixed $D$ and $r$) is an equality between two polynomials in the arc
weights $w_{a}$. (For instance, in Example \ref{exa.MTT.wMTT.1}, this equality
is $w_{\alpha}w_{\beta}+w_{\gamma}w_{\beta}=\det\left(
\begin{array}
[c]{cc}%
w_{\alpha}+w_{\gamma} & -w_{\alpha}\\
0 & w_{\beta}%
\end{array}
\right)  $.)

Therefore, thanks to Theorem \ref{thm.polyIDtrick}, it suffices to prove this
equality in the case when all arc weights $w_{a}$ are nonnegative integers. So
let us WLOG assume that all arc weights $w_{a}$ are nonnegative integers.

Let us now replace each arc $a$ of $D$ by $w_{a}$ many copies of the arc $a$
(having the same source as $a$ and the same target as $a$). The result is a
new digraph $D^{\prime}$. Here is an example:

\begin{example}
Let $D$ be the digraph%
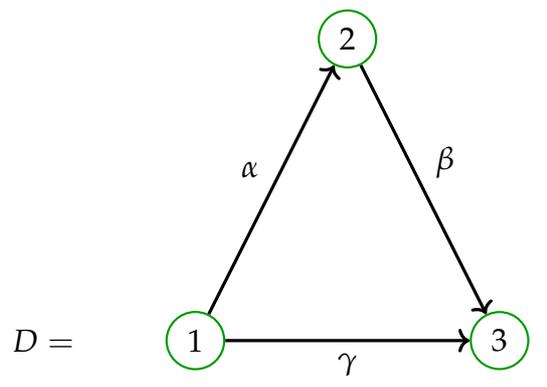
\[%
%
\ \ ,
\]
and let the arc weights be $w_{\alpha}=2$ and $w_{\beta}=3$ and $w_{\gamma}%
=2$. Then, $D^{\prime}$ looks as follows:%
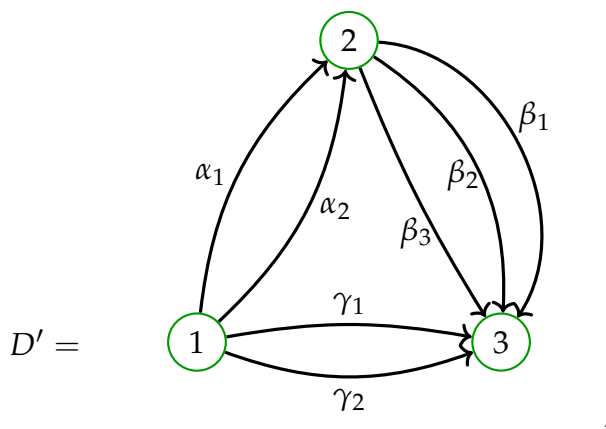
\[%
%
\ \ ,
\]
where $\alpha_{1},\alpha_{2}$ are the two arcs obtained from $\alpha$, and so on.
\end{example}

Now, recall that the digraph $D^{\prime}$ has the same vertices as $D$, but
each arc $a$ of $D$ has turned into $w_{a}$ arcs of $D^{\prime}$. Thus, the
weighted outdegree $\deg^{+w}i$ of a vertex $i$ of $D$ equals the (usual,
i.e., non-weighted) outdegree $\deg^{+}i$ of the same vertex $i$ of
$D^{\prime}$. Hence, the weighted Laplacian $L^{w}$ of $D$ is the (usual,
i.e., non-weighted) Laplacian of $D^{\prime}$.

Recall again that the digraph $D^{\prime}$ has the same vertices as $D$, but
each arc $a$ of $D$ has turned into $w_{a}$ arcs of $D^{\prime}$. Thus, each
subdigraph $B$ of $D$ gives rise to $w\left(  B\right)  $ many subdigraphs of
$D^{\prime}$ (because we can replace each arc $a$ of $B$ by any of the $w_{a}$
many copies of this arc in $D^{\prime}$). Moreover, this correspondence takes
spanning arborescences to spanning arborescences\footnote{More precisely: Let
$B$ be a subdigraph of $D$, and let $B^{\prime}$ be any of the $w\left(
B\right)  $ many subdigraphs of $D^{\prime}$ that are obtained from $B$
through this correspondence. Then, $B$ is a spanning arborescence of $D$
rooted to $r$ if and only if $B^{\prime}$ is a spanning arborescence of
$D^{\prime}$ rooted to $r$.}, and we can obtain any spanning arborescence of
$D^{\prime}$ in this way from exactly one $B$. Hence,%
\[
\sum_{\substack{B\text{ is a spanning}\\\text{arborescence}\\\text{of }D\text{
rooted to }r}}w\left(  B\right)  =\left(  \text{\# of spanning arborescences
of }D^{\prime}\text{ rooted to }r\right)  .
\]
Thus, applying the original MTT (Theorem \ref{thm.MTT.MTT}) to $D^{\prime}$
yields the weighted MTT for $D$ (since the weighted Laplacian $L^{w}$ of $D$
is the (usual, i.e., non-weighted) Laplacian of $D^{\prime}$). This completes
the proof of Theorem \ref{thm.MTT.wMTT}.

[\textit{Remark:} Alternatively, it is not hard to adapt our above proof of
the original MTT to the weighted case.]
\end{proof}

\subsubsection{Application: Counting trees by their degrees}

The weighted MTT has some applications that wouldn't be obvious from the
original MTT. Here is one:

\begin{exercise}
\label{exe.wMTT.degs}Let $n\geq2$ be an integer, and let $d_{1},d_{2}%
,\ldots,d_{n}$ be $n$ positive integers. An $n$\textbf{-tree} shall mean a
simple graph with vertex set $\left\{  1,2,\ldots,n\right\}  $ that is a tree.
We know from Corollary \ref{cor.cayley.n-trees} that there are $n^{n-2}$ many
$n$-trees. How many of these $n$-trees have the property that%
\[
\deg i=d_{i}\ \ \ \ \ \ \ \ \ \ \text{for each vertex }i\text{ ?}%
\]

\end{exercise}

\begin{proof}
[Solution.]The $n$-trees are just the spanning trees of the complete graph
$K_{n}$.

To incorporate the $\deg i=d_{i}$ condition into our count, we use a
generating function. So let us \textbf{not} fix the numbers $d_{1}%
,d_{2},\ldots,d_{n}$, but rather consider the polynomial%
\begin{equation}
P\left(  x_{1},x_{2},\ldots,x_{n}\right)  :=\sum_{T\text{ is an }%
n\text{-tree}}x_{1}^{\deg1}x_{2}^{\deg2}\cdots x_{n}^{\deg n}
\label{sol.wMTT.degs.1}%
\end{equation}
in $n$ indeterminates $x_{1},x_{2},\ldots,x_{n}$ (where $\deg i$ means the
degree of $i$ in $T$). Then, the $x_{1}^{d_{1}}x_{2}^{d_{2}}\cdots
x_{n}^{d_{n}}$-coefficient of this polynomial $P\left(  x_{1},x_{2}%
,\ldots,x_{n}\right)  $ is the \# of $n$-trees $T$ satisfying the property
that%
\[
\deg i=d_{i}\ \ \ \ \ \ \ \ \ \ \text{for each vertex }i
\]
(because each such $n$-tree $T$ contributes a monomial $x_{1}^{d_{1}}%
x_{2}^{d_{2}}\cdots x_{n}^{d_{n}}$ to the sum on the right hand side of
(\ref{sol.wMTT.degs.1}), whereas any other $n$-tree $T$ contributes a
different monomial to this sum).

Let us assign to each edge $ij$ of $K_{n}$ the weight $w_{ij}:=x_{i}x_{j}$.
Then, the definition of $P\left(  x_{1},x_{2},\ldots,x_{n}\right)  $ rewrites
as follows:%
\[
P\left(  x_{1},x_{2},\ldots,x_{n}\right)  =\sum_{T\text{ is an }n\text{-tree}%
}w\left(  T\right)  ,
\]
where $w\left(  T\right)  $ denotes the product of the weights of all edges of
$T$. (Indeed, for any subgraph $T$ of $K_{n}$, the weight $w\left(  T\right)
$ equals $x_{1}^{\deg1}x_{2}^{\deg2}\cdots x_{n}^{\deg n}$, where $\deg i$
means the degree of $i$ in $T$.)

We have assigned weights to the edges of the graph $K_{n}$; let us now assign
the same weights to the arcs of the digraph $K_{n}^{\operatorname*{bidir}}$.
That is, the two arcs $\left(  ij,1\right)  $ and $\left(  ij,2\right)  $
corresponding to an edge $ij$ of $K_{n}$ shall both have the weight%
\begin{equation}
w_{\left(  ij,1\right)  }=w_{\left(  ij,2\right)  }=w_{ij}=x_{i}x_{j}.
\label{sol.wMTT.degs.arcwt}%
\end{equation}

As we are already used to, we can replace spanning trees of $K_{n}$ by
spanning arborescences of $K_{n}^{\operatorname*{bidir}}$ rooted to $1$, since
the former are in bijection with the latter. Thus, we have
\begin{align*}
&  \left(  \text{\# of spanning trees of }K_{n}\right) \\
&  =\left(  \text{\# of spanning arborescences of }K_{n}%
^{\operatorname*{bidir}}\text{ rooted to }1\right)  .
\end{align*}
Moreover, since this bijection preserves weights (because of
(\ref{sol.wMTT.degs.arcwt})), we also have%
\[
\sum_{\substack{T\text{ is a spanning}\\\text{tree of }K_{n}}}w\left(
T\right)  =\sum_{\substack{B\text{ is a spanning}\\\text{arborescence of
}K_{n}^{\operatorname*{bidir}}\\\text{rooted to }1}}w\left(  B\right)  .
\]
In other words,
\[
\sum_{T\text{ is an }n\text{-tree}}w\left(  T\right)  =\sum_{\substack{B\text{
is a spanning}\\\text{arborescence of }K_{n}^{\operatorname*{bidir}%
}\\\text{rooted to }1}}w\left(  B\right)
\]
(since the spanning trees of $K_{n}$ are precisely the $n$-trees).

To compute the right hand side, we shall use the weighted Matrix-Tree Theorem.
The weighted Laplacian of $K_{n}^{\operatorname*{bidir}}$ (with the weights we
have just defined) is the $n\times n$-matrix $L^{w}$ with entries given by%
\begin{align*}
L_{i,j}^{w}  &  =\left(  \deg^{+w}i\right)  \cdot\left[  i=j\right]
-a_{i,j}^{w}\\
&  =%
\begin{cases}
\deg^{+w}i-a_{i,j}^{w}, & \text{if }i=j;\\
-a_{i,j}^{w}, & \text{if }i\neq j
\end{cases}
\\
&  =%
\begin{cases}
\deg^{+w}i, & \text{if }i=j;\\
-a_{i,j}^{w}, & \text{if }i\neq j
\end{cases}
\ \ \ \ \ \ \ \ \ \ \left(
\begin{array}
[c]{c}%
\text{since }a_{i,j}^{w}=0\text{ when }i=j\\
\text{(because }K_{n}^{\operatorname*{bidir}}\text{ has no loops)}%
\end{array}
\right) \\
&  =%
\begin{cases}
x_{i}\left(  x_{1}+x_{2}+\cdots+x_{n}\right)  -x_{i}x_{j}, & \text{if }i=j;\\
-x_{i}x_{j}, & \text{if }i\neq j
\end{cases}
\\
&  \ \ \ \ \ \ \ \ \ \ \ \ \ \ \ \ \ \ \ \ \left(
\begin{array}
[c]{c}%
\text{since }\deg^{+w}i=x_{i}x_{1}+x_{i}x_{2}+\cdots+x_{i}x_{i-1}+x_{i}%
x_{i+1}+\cdots+x_{i}x_{n}\\
=x_{i}\left(  x_{1}+x_{2}+\cdots+x_{i-1}+x_{i+1}+\cdots+x_{n}\right) \\
=x_{i}\left(  x_{1}+x_{2}+\cdots+x_{n}\right)  -x_{i}x_{i}\\
=x_{i}\left(  x_{1}+x_{2}+\cdots+x_{n}\right)  -x_{i}x_{j}\text{ whenever
}i=j\text{,}\\
\text{and since }a_{i,j}^{w}=x_{i}x_{j}\text{ whenever }i\neq j
\end{array}
\right) \\
&  =\left[  i=j\right]  x_{i}\left(  x_{1}+x_{2}+\cdots+x_{n}\right)
-x_{i}x_{j}\\
&  =x_{i}\left(  \left[  i=j\right]  \left(  x_{1}+x_{2}+\cdots+x_{n}\right)
-x_{j}\right)  .
\end{align*}
We can find its minor $\det\left(  L_{\sim1,\sim1}^{w}\right)  $ without too
much trouble (e.g., using row transformations similar to the ones we have done
back in the proof of Cayley's formula\footnote{The first step, of course, is
to factor an $x_{i}$ out of the $i$-th row for each $i$.}); the result is%
\[
\det\left(  L_{\sim1,\sim1}^{w}\right)  =x_{1}x_{2}\cdots x_{n}\left(
x_{1}+x_{2}+\cdots+x_{n}\right)  ^{n-2}.
\]

Summarizing what we have done so far,%
\begin{align}
P\left(  x_{1},x_{2},\ldots,x_{n}\right)   &  =\sum_{T\text{ is an
}n\text{-tree}}w\left(  T\right)  =\sum_{\substack{B\text{ is a spanning}%
\\\text{arborescence of }K_{n}^{\operatorname*{bidir}}\\\text{rooted to }%
1}}w\left(  B\right) \nonumber\\
&  =\det\left(  L_{\sim1,\sim1}^{w}\right)  \ \ \ \ \ \ \ \ \ \ \left(
\text{by the weighted Matrix-Tree Theorem}\right) \nonumber\\
&  =x_{1}x_{2}\cdots x_{n}\left(  x_{1}+x_{2}+\cdots+x_{n}\right)  ^{n-2}.
\label{sol.wMTT.degs.P=prod}%
\end{align}
As we recall, we are looking for the $x_{1}^{d_{1}}x_{2}^{d_{2}}\cdots
x_{n}^{d_{n}}$-coefficient in this polynomial. From
(\ref{sol.wMTT.degs.P=prod}), we see that%
\begin{align*}
&  \left(  \text{the }x_{1}^{d_{1}}x_{2}^{d_{2}}\cdots x_{n}^{d_{n}%
}\text{-coefficient of }P\left(  x_{1},x_{2},\ldots,x_{n}\right)  \right) \\
&  =\left(  \text{the }x_{1}^{d_{1}}x_{2}^{d_{2}}\cdots x_{n}^{d_{n}%
}\text{-coefficient of }x_{1}x_{2}\cdots x_{n}\left(  x_{1}+x_{2}+\cdots
+x_{n}\right)  ^{n-2}\right) \\
&  =\left(  \text{the }x_{1}^{d_{1}-1}x_{2}^{d_{2}-1}\cdots x_{n}^{d_{n}%
-1}\text{-coefficient of }\left(  x_{1}+x_{2}+\cdots+x_{n}\right)
^{n-2}\right)
\end{align*}
(because when we multiply a polynomial by $x_{1}x_{2}\cdots x_{n}$, all the
exponents in it get incremented by $1$, so its coefficients just shift by a
$1$ in each exponent).

Now, how can we describe the coefficients of $\left(  x_{1}+x_{2}+\cdots
+x_{n}\right)  ^{n-2}$, or, more generally, of $\left(  x_{1}+x_{2}%
+\cdots+x_{n}\right)  ^{m}$ for some $m\in\mathbb{N}$ ? These are the
so-called \textbf{multinomial coefficients} (named in analogy to the binomial
coefficients, which are their particular case for $n=2$). Their definition is
as follows: If $p_{1},p_{2},\ldots,p_{n},q$ are nonnegative integers with
$q=p_{1}+p_{2}+\cdots+p_{n}$, then the \textbf{multinomial coefficient}
$\dbinom{q}{p_{1},p_{2},\ldots,p_{n}}$ is defined to be $\dfrac{q!}%
{p_{1}!p_{2}!\cdots p_{n}!}$. If $q\neq p_{1}+p_{2}+\cdots+p_{n}$, then it is
defined to be $0$ instead. In either case, this coefficient is easily seen to
be an integer.\footnote{See \cite[Section 6.7]{24wd} for an introduction to
multinomial coefficients.} The
\textbf{\href{https://en.wikipedia.org/wiki/Multinomial_theorem}{\textbf{multinomial
formula} (aka \textbf{multinomial theorem})}} says that for each
$k\in\mathbb{N}$, we have%
\begin{align*}
\left(  x_{1}+x_{2}+\cdots+x_{n}\right)  ^{k}  &  =\sum_{\substack{i_{1}%
,i_{2},\ldots,i_{n}\in\mathbb{N};\\i_{1}+i_{2}+\cdots+i_{n}=k}}\dbinom
{k}{i_{1},i_{2},\ldots,i_{n}}x_{1}^{i_{1}}x_{2}^{i_{2}}\cdots x_{n}^{i_{n}}\\
&  =\sum_{i_{1},i_{2},\ldots,i_{n}\in\mathbb{N}}\dbinom{k}{i_{1},i_{2}%
,\ldots,i_{n}}x_{1}^{i_{1}}x_{2}^{i_{2}}\cdots x_{n}^{i_{n}}%
\end{align*}
(it does not matter whether we restrict the sum by the condition $i_{1}%
+i_{2}+\cdots+i_{n}=k$ or not, since the coefficient $\dbinom{k}{i_{1}%
,i_{2},\ldots,i_{n}}$ is defined to be $0$ when this condition is violated
anyway). Hence,%
\[
\left(  \text{the }x_{1}^{i_{1}}x_{2}^{i_{2}}\cdots x_{n}^{i_{n}%
}\text{-coefficient of }\left(  x_{1}+x_{2}+\cdots+x_{n}\right)  ^{k}\right)
=\dbinom{k}{i_{1},i_{2},\ldots,i_{n}}%
\]
for any $k\in\mathbb{N}$ and any $i_{1},i_{2},\ldots,i_{n}\in\mathbb{N}$. In
particular,%
\begin{align*}
&  \left(  \text{the }x_{1}^{d_{1}-1}x_{2}^{d_{2}-1}\cdots x_{n}^{d_{n}%
-1}\text{-coefficient of }\left(  x_{1}+x_{2}+\cdots+x_{n}\right)
^{n-2}\right) \\
&  =\dbinom{n-2}{d_{1}-1,\ d_{2}-1,\ \ldots,\ d_{n}-1}.
\end{align*}

Summarizing, we find%
\begin{align*}
&  \left(  \text{the }x_{1}^{d_{1}}x_{2}^{d_{2}}\cdots x_{n}^{d_{n}%
}\text{-coefficient of }P\left(  x_{1},x_{2},\ldots,x_{n}\right)  \right) \\
&  =\left(  \text{the }x_{1}^{d_{1}-1}x_{2}^{d_{2}-1}\cdots x_{n}^{d_{n}%
-1}\text{-coefficient of }\left(  x_{1}+x_{2}+\cdots+x_{n}\right)
^{n-2}\right) \\
&  =\dbinom{n-2}{d_{1}-1,\ d_{2}-1,\ \ldots,\ d_{n}-1}.
\end{align*}
However, the $x_{1}^{d_{1}}x_{2}^{d_{2}}\cdots x_{n}^{d_{n}}$-coefficient of
$P\left(  x_{1},x_{2},\ldots,x_{n}\right)  $ is the \# of $n$-trees $T$
satisfying the property that%
\[
\deg i=d_{i}\ \ \ \ \ \ \ \ \ \ \text{for each vertex }i
\]
(as we have seen above). Thus, we have proved the following:
\end{proof}

\begin{theorem}
[refined Cayley's formula]\label{thm.cayley.deg-refined}Let $n\geq2$ be an
integer, and let $d_{1},d_{2},\ldots,d_{n}$ be $n$ positive integers. Then,
the \# of $n$-trees with the property that
\[
\deg i=d_{i}\ \ \ \ \ \ \ \ \ \ \text{for each }i\in\left\{  1,2,\ldots
,n\right\}
\]
is the multinomial coefficient
\[
\dbinom{n-2}{d_{1}-1,\ d_{2}-1,\ \ldots,\ d_{n}-1}.
\]

\end{theorem}

This theorem can also be proved via a combinatorial bijection (see, e.g.,
\cite[\S 6]{Galvin}), or through an elementary induction argument using
Theorem \ref{thm.tree.leaf-ind.1} and Theorem \ref{thm.tree.leaf-ind.2} (see,
e.g., \cite[\S 5]{Galvin}).

\subsubsection{The weighted harmonic vector theorem}

The harmonic vector theorem for Laplacians (Theorem \ref{thm.MTT.harm}) also
has a weighted version:

\begin{theorem}
[harmonic vector theorem for weighted Laplacians]\label{thm.MTT.wharm}Let
$D=\left(  V,A,\psi\right)  $ be a multidigraph, where $V=\left\{
1,2,\ldots,n\right\}  $ for some $n\in\mathbb{N}$. Let $\mathbb{K}$ be a
commutative ring. Assume that an element $w_{a}\in\mathbb{K}$ is assigned to
each arc $a\in A$. For each $r\in V$, let $\tau^{w}\left(  D,r\right)  $ be
the sum of the weights of all the spanning arborescences of $D$ rooted to $r$.
Let $f^{w}$ be the row vector $\left(  \tau^{w}\left(  D,1\right)  ,\ \tau
^{w}\left(  D,2\right)  ,\ \ldots,\ \tau^{w}\left(  D,n\right)  \right)  $.
Let $L^{w}$ be the weighted Laplacian of $D$. Then, $f^{w}L^{w}=0$.
\end{theorem}

\begin{proof}
Similar to the unweighted case (Theorem \ref{thm.MTT.harm}).
\end{proof}

When the digraph $D$ is strongly connected, and the weights $w_{a}$ in Theorem
\ref{thm.MTT.wharm} are positive reals, we can make a slightly stronger claim:

\begin{corollary}
\label{cor.MTT.wharm.uniq}Let $D=\left(  V,A,\psi\right)  $ be a strongly
connected multidigraph, where $V=\left\{  1,2,\ldots,n\right\}  $ for some
$n\in\mathbb{N}$. Assume that a positive real number $w_{a}\in\mathbb{R}$ is
assigned to each arc $a\in A$. For each $r\in V$, let $\tau^{w}\left(
D,r\right)  $ be the sum of the weights of all the spanning arborescences of
$D$ rooted to $r$. Let $f^{w}$ be the row vector $\left(  \tau^{w}\left(
D,1\right)  ,\ \tau^{w}\left(  D,2\right)  ,\ \ldots,\ \tau^{w}\left(
D,n\right)  \right)  $. Let $L^{w}$ be the weighted Laplacian of $D$. Then:

\begin{enumerate}
\item[\textbf{(a)}] We have $f^{w}L^{w}=0$.

\item[\textbf{(b)}] All entries of the vector $f^{w}$ are positive reals.

\item[\textbf{(c)}] If all numbers $w_{a}$ are rational, then all entries of
the vector $f^{w}$ are rational.

\item[\textbf{(d)}] The matrix $L^{w}$ has rank $n-1$.

\item[\textbf{(e)}] Let $s_{w}$ be the sum of all entries of $f^{w}$. Then,
there is a unique row vector $x\in\mathbb{R}^{1\times n}$ such that $xL^{w}=0$
and such that the sum of all entries of $x$ is $1$. This vector $x$ is
$\dfrac{1}{s_{w}}f^{w}$.
\end{enumerate}
\end{corollary}

\begin{proof}
\textbf{(a)} This follows from Theorem \ref{thm.MTT.wharm} (applied to
$\mathbb{K}=\mathbb{R}$). \medskip

\textbf{(b)} We must show that $\tau^{w}\left(  D,r\right)  $ is positive for
each $r\in V$. So let us fix $r\in V$. Recall that $\tau^{w}\left(
D,r\right)  $ is the sum of the weights of all the spanning arborescences of
$D$ rooted to $r$. This sum is not an empty sum (since there exists a spanning
arborescence of $D$ rooted from $r$\ \ \ \ \footnote{This follows from the
analogue of Theorem \ref{thm.spanning-arbor.exists} for to-roots instead of
from-roots (since $D$ is strongly connected, and thus $r$ is a to-root of
$D$).}), and all its addends are positive (since the weight of a spanning
arborescence of $D$ is a product of some of the positive reals $w_{a}$, and
thus itself positive). Hence, this sum is positive. In other words, $\tau
^{w}\left(  D,r\right)  $ is positive. This proves part \textbf{(b)}. \medskip

\textbf{(c)} This follows from the construction of $f^{w}$. \medskip

\textbf{(d)} Pick any vertex $r\in V$. Theorem \ref{thm.MTT.wMTT} (applied to
$\mathbb{K}=\mathbb{R}$) shows that
\[
\sum_{\substack{B\text{ is a spanning}\\\text{arborescence}\\\text{of }D\text{
rooted to }r}}w\left(  B\right)  =\det\left(  L_{\sim r,\sim r}^{w}\right)  .
\]
The left hand side of this equality equals the number $\tau^{w}\left(
D,r\right)  $, which is positive (as we showed in the proof of part
\textbf{(b)} above). Thus, so is the right hand side. In other words,
$\det\left(  L_{\sim r,\sim r}^{w}\right)  >0$. Hence, in particular,
$\det\left(  L_{\sim r,\sim r}^{w}\right)  \neq0$. This shows that the
$n\times n$-matrix $L^{w}$ has an invertible $\left(  n-1\right)
\times\left(  n-1\right)  $-submatrix (namely, $L_{\sim r,\sim r}^{w}$).
Hence, its rank is at least $n-1$. In other words, $\operatorname*{rank}%
\left(  L^{w}\right)  \geq n-1$.

But the vector $f^{w}$ is nonzero (by part \textbf{(b)}) and lies in the left
nullspace\footnote{The \textbf{left nullspace} of a matrix $A$ means the
vector space of all row vectors $v$ satisfying $vA=0$.} of $L^{w}$ (since part
\textbf{(a)} says $f^{w}L^{w}=0$). Hence, the matrix $L^{w}$ is singular, and
thus its rank is $\operatorname*{rank}\left(  L^{w}\right)  <n$. Combined with
$\operatorname*{rank}\left(  L^{w}\right)  \geq n-1$, this yields
$\operatorname*{rank}\left(  L^{w}\right)  =n-1$. Thus, part \textbf{(d)} is
proved. \medskip

\textbf{(e)} First, we observe that $s_{w}$ is positive (by part
\textbf{(b)}), thus nonzero. The vector $\dfrac{1}{s_{w}}f^{w}$ really is a
row vector $x\in\mathbb{R}^{1\times n}$ such that $xL^{w}=0$ (since $\left(
\dfrac{1}{s_{w}}f^{w}\right)  L^{w}=\dfrac{1}{s_{w}}\underbrace{f^{w}L^{w}%
}_{\substack{=0\\\text{(by part \textbf{(a)})}}}=0$) and such that the sum of
all entries of $x$ is $1$ (since the sum of all entries of $f^{w}$ is $s_{w}%
$). It remains to show that it is the only such vector.

Indeed, assume the contrary. Thus, there are two distinct such vectors $x$.
Let them be $y$ and $z$. These two vectors $y$ and $z$ belong to the left
nullspace of $L^{w}$ (since they are vectors $x$ satisfying $xL^{w}=0$) but
cannot be linearly dependent (since each of them has its sum of entries equal
to $1$). Thus, they are linearly independent. This shows that the left
nullspace of $L^{w}$ has dimension $\geq2$ (since it contains the two linearly
independent vectors $y$ and $z$). In other words, the matrix $L^{w}$ has rank
$\leq n-2$. But this contradicts part \textbf{(d)} of the corollary. This
contradiction completes our proof of part \textbf{(e)}.
\end{proof}

A consequence of Corollary \ref{cor.MTT.wharm.uniq} is a graph-theoretical
proof of a famous fact from matrix analysis (often worded in the language of
Markov chains, cf. \cite[Theorem 11.10]{GriSne07}):

\begin{corollary}
\label{cor.MTT.markov}Let $n$ be a positive integer. Let $P=\left(
p_{i,j}\right)  _{1\leq i,j\leq n}\in\mathbb{R}^{n\times n}$ be a matrix whose
entries $p_{i,j}$ are nonnegative reals. Assume that%
\[
\sum_{j=1}^{n}p_{i,j}=1\ \ \ \ \ \ \ \ \ \ \text{for each }i\in\left\{
1,2,\ldots,n\right\}  .
\]
(Matrices $P$ satisfying these conditions are called \textbf{stochastic}.) Let
$D$ be the simple digraph $\left(  V,A\right)  $ with vertex set $V=\left\{
1,2,\ldots,n\right\}  $ and arc set $A=\left\{  \left(  i,j\right)  \in
V\times V\ \mid\ p_{i,j}>0\right\}  $. Assume that this digraph $D$ is
strongly connected.\ (Matrices $P$ satisfying this condition are called
\textbf{irreducible}.) Then:

\begin{enumerate}
\item[\textbf{(a)}] There is a unique row vector $x\in\mathbb{R}^{1\times n}$
such that $xP=x$ and such that the sum of all entries of $x$ is $1$. (In the
lingo of probabilists, this vector $x$ is called the \textbf{stationary
distribution} or the \textbf{steady state} of the Markov chain defined by $P$.)

\item[\textbf{(b)}] The entries of this vector $x$ are positive reals.

\item[\textbf{(c)}] If all the $p_{i,j}$ are rational, then the entries of
this vector $x$ are rational.
\end{enumerate}
\end{corollary}

\begin{proof}
We identify our simple digraph $D=\left(  V,A\right)  $ with the corresponding
multidigraph $D^{\operatorname*{mult}}=\left(  V,A,\psi\right)  $ (see
Definition \ref{def.sdg.Dmult}). This is a strongly connected multidigraph.

\begin{noncompile}
For each $\left(  i,j\right)  \in A$, we have $p_{i,j}>0$ (by the definition
of $A$). On the other hand, if $i$ and $j$ are two vertices of $D$ that do not
satisfy $\left(  i,j\right)  \in A$, then $p_{i,j}\leq0$ (since $p_{i,j}>0$
would mean $\left(  i,j\right)  \in A$) and therefore $p_{i,j}=0$ (since we
know that the entry $p_{i,j}$ is nonnegative).
\end{noncompile}

To each arc $\left(  i,j\right)  $ of $D$, we assign the positive real number
$w_{\left(  i,j\right)  }:=p_{i,j}$ (this is positive, because $\left(
i,j\right)  \in A$ guarantees that $p_{i,j}>0$). We call this number the
\textquotedblleft weight\textquotedblright\ of the arc $\left(  i,j\right)  $.
Then, Definition \ref{def.wMTT.wt-defs} \textbf{(a)} (applied to
$\mathbb{K}=\mathbb{R}$) defines a number $a_{i,j}^{w}$ for any two vertices
$i,j\in V=\left\{  1,2,\ldots,n\right\}  $. We claim that any two vertices
$i,j\in V=\left\{  1,2,\ldots,n\right\}  $ satisfy%
\begin{equation}
a_{i,j}^{w}=p_{i,j}. \label{pf.cor.MTT.markov.awij=}%
\end{equation}

[\textit{Proof of (\ref{pf.cor.MTT.markov.awij=}):} Let $i,j\in V=\left\{
1,2,\ldots,n\right\}  $ be two vertices. Then, $p_{i,j}\geq0$ (since all
entries of $P$ are nonnegative). Hence, we are in one of the following two cases:

\textit{Case 1:} We have $p_{i,j}>0$.

\textit{Case 2:} We have $p_{i,j}=0$.

Let us consider Case 1. In this case, we have $p_{i,j}>0$. Thus, $\left(
i,j\right)  \in A$ (by the definition of $A$). Hence, the digraph $D$ has
exactly one arc from $i$ to $j$, namely the arc $\left(  i,j\right)  $. The
definition of $a_{i,j}^{w}$ thus shows that $a_{i,j}^{w}=w_{\left(
i,j\right)  }=p_{i,j}$ (by the definition of $w_{\left(  i,j\right)  }$). This
proves (\ref{pf.cor.MTT.markov.awij=}) in Case 1.

Let us now consider Case 2. In this case, we have $p_{i,j}=0$. Thus, $\left(
i,j\right)  \notin A$ (by the definition of $A$). Hence, the digraph $D$ has
no arc from $i$ to $j$. The definition of $a_{i,j}^{w}$ thus shows that
$a_{i,j}^{w}=\left(  \text{empty sum}\right)  =0=p_{i,j}$. This proves
(\ref{pf.cor.MTT.markov.awij=}) in Case 2.

Thus, the proof of (\ref{pf.cor.MTT.markov.awij=}) is complete (since we have
covered both cases).] \medskip

Furthermore, Definition \ref{def.wMTT.wt-defs} \textbf{(b)} (applied to
$\mathbb{K}=\mathbb{R}$) yields that the weighted outdegree $\deg^{+w}i$ of
any vertex $i\in V=\left\{  1,2,\ldots,n\right\}  $ is%
\begin{align}
\deg^{+w}i  &  =\sum_{\substack{a\in A;\\\text{the source of }a\text{ is }%
i}}w_{a}=\underbrace{\sum_{j\in V}}_{=\sum_{j=1}^{n}}\ \ \underbrace{\sum
_{\substack{a\in A\text{ is an}\\\text{arc from }i\text{ to }j}}w_{a}%
}_{\substack{=a_{i,j}^{w}\\\text{(by the definition of }a_{i,j}^{w}\text{)}%
}}=\sum_{j=1}^{n}\underbrace{a_{i,j}^{w}}_{\substack{=p_{i,j}\\\text{(by
(\ref{pf.cor.MTT.markov.awij=}))}}}\nonumber\\
&  =\sum_{j=1}^{n}p_{i,j}=1 \label{pf.cor.MTT.markov.deg=1}%
\end{align}
(by one of the assumptions of the corollary).

Let $L^{w}$ be the weighted Laplacian of $D$. Then, for each $i,j\in
V=\left\{  1,2,\ldots,n\right\}  $, we have%
\begin{align}
L_{i,j}^{w}  &  =\underbrace{\left(  \deg^{+w}i\right)  }%
_{\substack{=1\\\text{(by (\ref{pf.cor.MTT.markov.deg=1}))}}}\cdot\left[
i=j\right]  -\underbrace{a_{i,j}^{w}}_{\substack{=p_{i,j}\\\text{(by
(\ref{pf.cor.MTT.markov.awij=}))}}}\ \ \ \ \ \ \ \ \ \ \left(  \text{by
Definition \ref{def.wMTT.wt-defs} \textbf{(d)}}\right) \nonumber\\
&  =\left[  i=j\right]  -p_{i,j}. \label{pf.cor.MTT.markov.Lwij=}%
\end{align}
Note that the numbers $\left[  i=j\right]  $ and $p_{i,j}$ on the right hand
side here are the $\left(  i,j\right)  $-th entries of the matrices $I_{n}$
and $P$, respectively. Thus, (\ref{pf.cor.MTT.markov.Lwij=}) shows that%
\begin{equation}
L^{w}=I_{n}-P. \label{pf.cor.MTT.markov.Lw=}%
\end{equation}

For each $r\in V$, let $\tau^{w}\left(  D,r\right)  $ be the sum of the
weights of all the spanning arborescences of $D$ rooted to $r$. Let $f^{w}$ be
the row vector \newline$\left(  \tau^{w}\left(  D,1\right)  ,\ \tau^{w}\left(
D,2\right)  ,\ \ldots,\ \tau^{w}\left(  D,n\right)  \right)  $. Let $s_{w}$ be
the sum of all entries of $f^{w}$.

Corollary \ref{cor.MTT.wharm.uniq} \textbf{(e)} shows that there is a unique
row vector $x\in\mathbb{R}^{1\times n}$ such that $xL^{w}=0$ and such that the
sum of all entries of $x$ is $1$. Since the equation $xL^{w}=0$ is equivalent
to $xP=x$ (because (\ref{pf.cor.MTT.markov.Lw=}) shows that $xL^{w}=x\left(
I_{n}-P\right)  =x-xP$, and thus we have the chain of logical equivalences
$\left(  xL^{w}=0\right)  \ \Longleftrightarrow\ \left(  x-xP=0\right)
\ \Longleftrightarrow\ \left(  xP=x\right)  $), we can rewrite this as
follows: There is a unique row vector $x\in\mathbb{R}^{1\times n}$ such that
$xP=x$ and such that the sum of all entries of $x$ is $1$. This proves
Corollary \ref{cor.MTT.markov} \textbf{(a)}. \medskip

\textbf{(b)} Consider the unique row vector $x\in\mathbb{R}^{1\times n}$ such
that $xP=x$ and such that the sum of all entries of $x$ is $1$. This is
precisely the unique vector $x\in\mathbb{R}^{1\times n}$ such that $xL^{w}=0$
and such that the sum of all entries of $x$ is $1$ (because, as we saw above,
the equation $xL^{w}=0$ is equivalent to $xP=x$). By Corollary
\ref{cor.MTT.wharm.uniq} \textbf{(e)}, this unique vector $x$ is therefore
$\dfrac{1}{s_{w}}f^{w}$. However, the entries of $f^{w}$ are positive reals
(by Corollary \ref{cor.MTT.wharm.uniq} \textbf{(b)}), and thus their sum
$s_{w}$ is a positive real as well. Therefore, the entries of $\dfrac{1}%
{s_{w}}f^{w}$ are positive reals. In other words, the entries of $x$ are
positive reals (since $x$ is $\dfrac{1}{s_{w}}f^{w}$). This proves Corollary
\ref{cor.MTT.markov} \textbf{(b)}. \medskip

\textbf{(c)} Assume that all the $p_{i,j}$ are rational. Hence, all entries of
the vector $f^{w}$ are rational (by Corollary \ref{cor.MTT.wharm.uniq}
\textbf{(c)}), whence their sum $s_{w}$ is rational as well. Therefore, the
entries of $\dfrac{1}{s_{w}}f^{w}$ are rational, too. In other words, the
entries of $x$ are rational (since $x$ is $\dfrac{1}{s_{w}}f^{w}$). This
proves Corollary \ref{cor.MTT.markov} \textbf{(c)}.
\end{proof}

Here ends our study of spanning trees and their enumeration. An interested
reader can learn more from \cite{Rubey00}, \cite{Holzer22}, \cite{Moon70} and
\cite{GrSaSu14}.

\section{Colorings}

Now to something different: Let's color the vertices of a graph!

\subsection{Definition}

This is a serious course, so our colors are positive integers. Coloring the
vertices thus means assigning a color (= a positive integer) to each vertex.
Here are the details:

\begin{definition}
\label{def.coloring.coloring}Let $G=\left(  V,E,\varphi\right)  $ be a
multigraph. Let $k\in\mathbb{N}$.

\begin{enumerate}
\item[\textbf{(a)}] A $k$\textbf{-coloring} of $G$ means a map $f:V\rightarrow
\left\{  1,2,\ldots,k\right\}  $. Given such a $k$-coloring $f$, we refer to
the numbers $1,2,\ldots,k$ as the \textbf{colors}, and we refer to each value
$f\left(  v\right)  $ as the \textbf{color} of the vertex $v$ in the
$k$-coloring $f$.

\item[\textbf{(b)}] A $k$-coloring $f$ of $G$ is said to be \textbf{proper} if
no two adjacent vertices of $G$ have the same color. (In other words, a
$k$-coloring $f$ of $G$ is proper if there exists no edge of $G$ whose
endpoints $u$ and $v$ satisfy $f\left(  u\right)  =f\left(  v\right)  $.)
\end{enumerate}
\end{definition}

\begin{example}
Here are two $7$-colorings of a graph:%
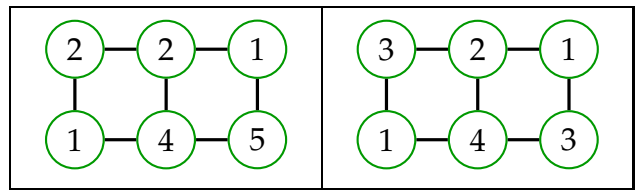
\[%
\begin{tabular}
[c]{|c|c|}\hline
$\ \
\begin{tikzpicture}[scale=1.2]
\begin{scope}[every node/.style={circle,thick,draw=green!60!black}]
\node(1) at (0,0) {$1$};
\node(2) at (0,1) {$2$};
\node(3) at (1,0) {$4$};
\node(4) at (1,1) {$2$};
\node(5) at (2,0) {$5$};
\node(6) at (2,1) {$1$};
\end{scope}
\begin{scope}[every edge/.style={draw=black,very thick}, every loop/.style={}]
\path
[-] (1) edge (2) edge (3) (4) edge (2) edge (3) edge (6) (5) edge (3) edge (6);
\end{scope}
\end{tikzpicture}%
\ \ $ & $\ \
\begin{tikzpicture}[scale=1.2]
\begin{scope}[every node/.style={circle,thick,draw=green!60!black}]
\node(1) at (0,0) {$1$};
\node(2) at (0,1) {$3$};
\node(3) at (1,0) {$4$};
\node(4) at (1,1) {$2$};
\node(5) at (2,0) {$3$};
\node(6) at (2,1) {$1$};
\end{scope}
\begin{scope}[every edge/.style={draw=black,very thick}, every loop/.style={}]
\path
[-] (1) edge (2) edge (3) (4) edge (2) edge (3) edge (6) (5) edge (3) edge (6);
\end{scope}
\end{tikzpicture}%
\ \ $\\\hline
\end{tabular}
\
\]
(where the numbers on the nodes are not the vertices, but rather the colors of
the vertices). The $7$-coloring on the left (yes, it is a $7$-coloring, even
though it does not actually use the colors $3$, $6$ and $7$) is not proper,
because the two adjacent vertices on the top left have the same color. The
$7$-coloring on the right, however, is proper.
\end{example}

\Needspace{15pc}

\begin{example}
Here is a bunch of graphs:%
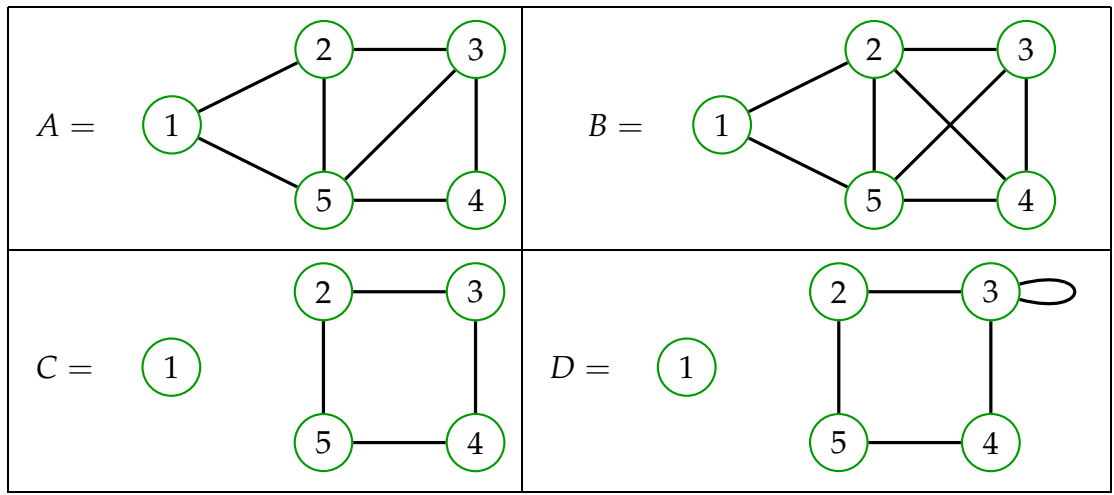
\[%
%
\ \ \ \ \ \ .
\]
Which of them have proper $3$-colorings?

\begin{itemize}
\item The graph $A$ has a proper $3$-coloring. For example, the map $f$ that
sends the vertices $1,2,3,4,5$ to the colors $1,2,1,2,3$ (respectively) is a
proper $3$-coloring.

\item The graph $B$ has no proper $3$-coloring. Indeed, the four vertices
$2,3,4,5$ are mutually adjacent, so they would have to have $4$ distinct
colors in a proper $k$-coloring; but this is not possible unless $k\geq4$.

\item The graph $C$ has a proper $3$-coloring and even a proper $2$-coloring
(e.g., assigning color $1$ to each odd vertex and color $2$ to each even vertex).

\item The graph $D$ has no proper $3$-coloring and, in fact, no proper
$k$-coloring for any $k\in\mathbb{N}$. The reason is that the vertex $3$ is
adjacent to itself, but obviously has the same color as itself no matter what
the $k$-coloring is. More generally, a graph with a loop cannot have a proper
$k$-coloring for any $k\in\mathbb{N}$.
\end{itemize}
\end{example}

\begin{example}
Let $n\in\mathbb{N}$. The $n$-hypercube $Q_{n}$ (introduced in Definition
\ref{def.hypercube}) has a proper $2$-coloring: Namely, the map
\begin{align*}
f:\left\{  0,1\right\}  ^{n}  &  \rightarrow\left\{  1,2\right\}  ,\\
\left(  a_{1},a_{2},\ldots,a_{n}\right)   &  \mapsto%
\begin{cases}
1, & \text{if }a_{1}+a_{2}+\cdots+a_{n}\text{ is odd};\\
2, & \text{if }a_{1}+a_{2}+\cdots+a_{n}\text{ is even}%
\end{cases}
\end{align*}
is a proper $2$-coloring of $Q_{n}$. (Check this! It boils down to the fact
that if two bitstrings $\left(  a_{1},a_{2},\ldots,a_{n}\right)  $ and
$\left(  b_{1},b_{2},\ldots,b_{n}\right)  $ differ in exactly one entry, then
the corresponding sums $a_{1}+a_{2}+\cdots+a_{n}$ and $b_{1}+b_{2}%
+\cdots+b_{n}$ differ by exactly $1$.)
\end{example}

\begin{example}
Let $n$ and $m$ be two positive integers. The Cartesian product $P_{n}\times
P_{m}$ of the $n$-th path graph $P_{n}$ and the $m$-th path graph $P_{m}$ is
known as the $\left(  n,m\right)  $\textbf{-grid graph}, as it looks as
follows:%
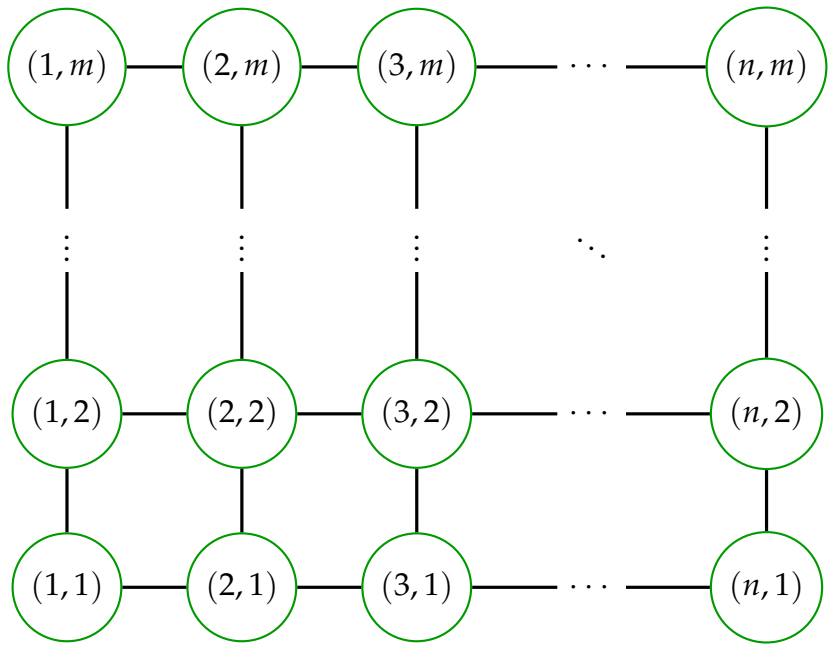
\[%
%
\ \ .
\]
This $\left(  n,m\right)  $-grid graph $P_{n}\times P_{m}$ has a proper
$2$-coloring: namely, the map that sends each vertex $\left(  i,j\right)  $ to
$%
\begin{cases}
1, & \text{if }i+j\text{ is even};\\
2, & \text{if }i+j\text{ is odd}.
\end{cases}
$

This $2$-coloring is called the \textquotedblleft chessboard
coloring\textquotedblright\ for a fairly obvious reason (view each vertex as a
square of a chessboard).

More generally, if $G$ and $H$ are two simple graphs each having a proper
$2$-coloring, then their Cartesian product $G\times H$ has a proper
$2$-coloring as well. (See Exercise \ref{exe.7.12} for the proof.)
\end{example}

\begin{example}
Here is the Petersen graph (as defined in Subsection
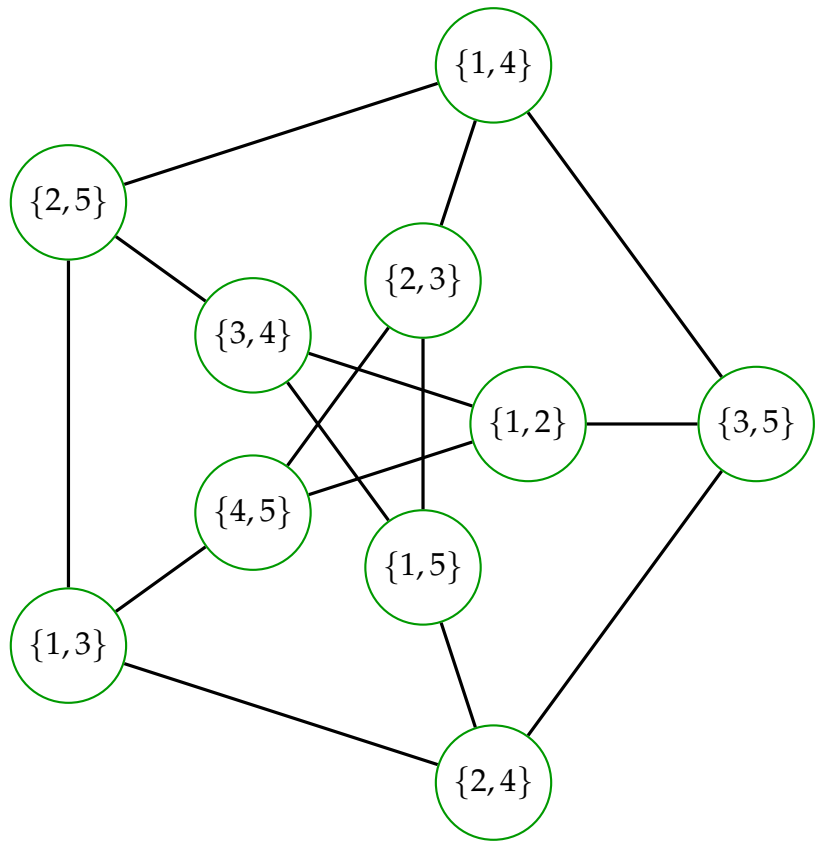
\ref{subsec.sg.complete.kneser}):%
\[%
\begin{tikzpicture}
\begin{scope}[every node/.style={circle,thick,draw=green!60!black}]
\node(A) at (0:2) {$\set{1,2}$};
\node(B) at (360/5:2) {$\set{2,3}$};
\node(C) at (2*360/5:2) {$\set{3,4}$};
\node(D) at (3*360/5:2) {$\set{4,5}$};
\node(E) at (4*360/5:2) {$\set{1,5}$};
\node(A2) at (0:5) {$\set{3,5}$};
\node(B2) at (360/5:5) {$\set{1,4}$};
\node(C2) at (2*360/5:5) {$\set{2,5}$};
\node(D2) at (3*360/5:5) {$\set{1,3}$};
\node(E2) at (4*360/5:5) {$\set{2,4}$};
\end{scope}
\begin{scope}[every edge/.style={draw=black,very thick}]
\path
[-] (A) edge (C) (C) edge (E) (E) edge (B) (B) edge (D) (D) edge (A) (A) edge (A2) (B) edge (B2) (C) edge (C2) (D) edge (D2) (E) edge (E2) (A2) edge (B2) (B2) edge (C2) (C2) edge (D2) (D2) edge (E2) (E2) edge (A2);
\end{scope}
\end{tikzpicture}%
\ \ .
\]
I claim that it has a proper $3$-coloring. Can you find it?
\end{example}

As we see, some graphs have proper $3$-colorings, while others don't. Clearly,
having $4$ mutually adjacent vertices makes a proper $3$-coloring impossible
(indeed, by the pigeonhole principle, two of them must have the same color),
but this is far from an \textquotedblleft if and only if\textquotedblright.
The question of determining whether a given graph has a proper $3$-coloring is
NP-complete (see, e.g., \cite[Proposition 4.11]{Goldre10}).

\subsection{2-colorings}

\subsubsection{The undirected case}

In contrast, the existence of proper $2$-colorings is a much simpler question.
The following is a nice criterion:

\begin{theorem}
[2-coloring equivalence theorem]\label{thm.coloring.2-color-eq}Let $G$ be a
multigraph. Then, the following three statements are equivalent:

\begin{itemize}
\item \textbf{Statement B1:} The graph $G$ has a proper $2$-coloring.

\item \textbf{Statement B2:} The graph $G$ has no cycles of odd length.

\item \textbf{Statement B3:} The graph $G$ has no circuits of odd length.
\end{itemize}
\end{theorem}

To prove this theorem, we will need a fact that is somewhat similar to
Proposition \ref{prop.mg.cyc.btf-walk-cyc}:

\begin{proposition}
\label{prop.coloring.2-color-odd-length-cyc}Let $G$ be a multigraph. Let $u$
and $v$ be two vertices of $G$. Let $\mathbf{w}$ be an odd-length walk from
$u$ to $v$. Then, $\mathbf{w}$ contains either an odd-length \textbf{path}
from $u$ to $v$ or an odd-length \textbf{cycle} (or both).
\end{proposition}

Here, we are using the following rather intuitive terminology:

\begin{itemize}
\item A walk is said to be \textbf{odd-length} if its length is odd.

\item A walk $\mathbf{w}$ is said to \textbf{contain} a walk $\mathbf{v}$ if
each edge of $\mathbf{v}$ is an edge of $\mathbf{w}$. (This does not
necessarily mean that $\mathbf{v}$ appears in $\mathbf{w}$ as a contiguous block.)

\item We remind the reader once again that a \textquotedblleft
circuit\textquotedblright\ just means a closed walk to us; we impose no
further requirements.
\end{itemize}

\begin{example}
Consider the following simple graph (which we treat as a multigraph):%
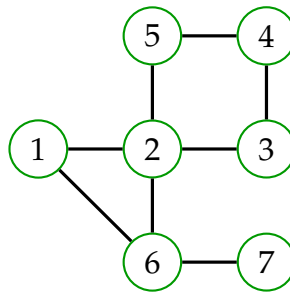
\[%
\begin{tikzpicture}[scale=1.5]
\begin{scope}[every node/.style={circle,thick,draw=green!60!black}]
\node(1) at (-1,0) {$1$};
\node(2) at (0,0) {$2$};
\node(3) at (1,0) {$3$};
\node(4) at (1,1) {$4$};
\node(5) at (0,1) {$5$};
\node(6) at (0,-1) {$6$};
\node(7) at (1,-1) {$7$};
\end{scope}
\begin{scope}[every edge/.style={draw=black,very thick}]
\path[-] (1) edge (2) (1) edge (6) (2) edge (3) (2) edge (5);
\path[-] (2) edge (6) (3) edge (4) (4) edge (5);
\path[-] (6) edge (7);
\end{scope}
\end{tikzpicture}%
\ \ .
\]

\textbf{(a)} The odd-length walk $\left(  1,\ast,2,\ast,3,\ast,4,\ast
,5,\ast,2,\ast,6,\ast,7\right)  $ (we are using asterisks for the edges, since
they can be trivially recovered from the vertices) contains the odd-length
path $\left(  1,\ast,2,\ast,6,\ast,7\right)  $ from $1$ to $7$.

\textbf{(b)} The odd-length walk $\left(  3,\ast,2,\ast,1,\ast,6,\ast
,2,\ast,3\right)  $ contains the odd-length cycle $\left(  2,\ast
,1,\ast,6,\ast,2\right)  $.
\end{example}

\begin{proof}
[Proof of Proposition \ref{prop.coloring.2-color-odd-length-cyc}.]We apply
strong induction on the length of $\mathbf{w}$.

Thus, we fix a $k\in\mathbb{N}$, and we assume (as the induction hypothesis)
that Proposition \ref{prop.coloring.2-color-odd-length-cyc} is already proved
for all odd-length walks of length $<k$. Now, we must prove it for an
odd-length walk $\mathbf{w}$ of length $k$.

Write this walk $\mathbf{w}$ as $\mathbf{w}=\left(  w_{0},\ast,w_{1}%
,\ast,w_{2},\ldots,\ast,w_{k}\right)  $. Hence, $k$ is the length of
$\mathbf{w}$, and thus is odd.

We must prove that $\mathbf{w}$ contains either an odd-length path from $u$ to
$v$ or an odd-length cycle.

If $\mathbf{w}$ itself is a path, then we are done. So WLOG assume that
$\mathbf{w}$ is not a path. Thus, two of the vertices $w_{0},w_{1}%
,\ldots,w_{k}$ of $\mathbf{w}$ are equal. In other words, there exists a pair
$\left(  i,j\right)  $ of integers $i$ and $j$ with $0\leq i<j\leq k$ and
$w_{i}=w_{j}$. Among all such pairs, we pick one with \textbf{minimum}
difference $j-i$. Then, the vertices $w_{i},w_{i+1},\ldots,w_{j-1}$ are
distinct (since $j-i$ is minimum).

Let $\mathbf{c}$ be the part of $\mathbf{w}$ between $w_{i}$ and $w_{j}$;
thus,\footnote{Here is an illustration (which, however, is a bit simplistic:
the walk $\mathbf{w}$ can intersect itself arbitrarily many times, not just
once as shown here):%
\[%
\begin{tikzpicture}[scale=2.1]
\begin{scope}[every node/.style={circle,thick,draw=green!60!black}]
\node(w0) at (0,0) {$w_0$};
\node(w1) at (1,0) {$w_1$};
\node(wi) at (3,0) {$w_i = w_j$};
\node(wj1) at (4,0) {$w_{j+1}$};
\node(wk) at (6,0) {$w_k$};
\node(wi1) at (2,1) {$w_{i+1}$};
\node(wjm1) at (4,1) {$w_{j-1}$};
\end{scope}
\node(w2) at (2,0) {$\cdots$};
\node(w5) at (5,0) {$\cdots$};
\node(ws) at (3,2) {$\cdots$};
\begin{scope}[every edge/.style={draw=black,very thick}, every loop/.style={}]
\path[-] (w0) edge (w1) (w1) edge (w2) (w2) edge (wi);
\path[-] (wi) edge (wj1) (wj1) edge (w5) (w5) edge (wk);
\end{scope}
\begin{scope}[every edge/.style={draw=blue,very thick}, every loop/.style={}]
\path[-] (wi) edge[bend left] (wi1) (wi1) edge[bend left] (ws);
\path[-] (ws) edge[bend left] (wjm1) (wjm1) edge[bend left] (wi);
\end{scope}
\end{tikzpicture}%
\ \ .
\]
The blue edges here form the walk $\mathbf{c}$.}%
\[
\mathbf{c}=\left(  w_{i},\ast,w_{i+1},\ast,\ldots,\ast,w_{j}\right)  .
\]
This $\mathbf{c}$ is clearly a closed walk (since $w_{i}=w_{j}$). If $j-i$ is
odd, then this closed walk $\mathbf{c}$ is a cycle (indeed, its vertices
$w_{i},w_{i+1},\ldots,w_{j-1}$ are distinct, and therefore its edges are
distinct as well\footnote{For the very skeptical, here is a \textit{proof} of
this:
\par
Assume (for the sake of contradiction) that the walk $\mathbf{c}$ has two
equal edges. Let the first of them be an edge between $w_{p}$ and $w_{p+1}$,
and let the second be an edge between $w_{q}$ and $w_{q+1}$, for some distinct
elements $p$ and $q$ of $\left\{  i,i+1,\ldots,j-1\right\}  $. Since equal
edges have equal endpoints, we thus have $\left\{  w_{p},w_{p+1}\right\}
=\left\{  w_{q},w_{q+1}\right\}  $, so that $w_{p}\in\left\{  w_{p}%
,w_{p+1}\right\}  =\left\{  w_{q},w_{q+1}\right\}  $. In other words, $w_{p}$
equals either $w_{q}$ or $w_{q+1}$. Since $w_{p}\neq w_{q}$ (because
$w_{i},w_{i+1},\ldots,w_{j-1}$ are distinct), this entails that $w_{p}%
=w_{q+1}$. Similarly, $w_{q}=w_{p+1}$.
\par
However, $p$ and $q$ are distinct. Thus, at least one of $p$ and $q$ is
distinct from $j-1$. We WLOG assume that $q\neq j-1$ (otherwise, we can swap
$p$ with $q$). Hence, $q+1\neq j$, so that $q+1\in\left\{  i,i+1,\ldots
,j-1\right\}  $. Thus, from $w_{p}=w_{q+1}$, we conclude that $p=q+1$ (since
$w_{i},w_{i+1},\ldots,w_{j-1}$ are distinct). Thus, $p=q+1>q$, so that
$p+1>p>q$ and therefore $p+1\neq q$. However, $w_{q}=w_{p+1}$. If $p+1$ was an
element of $\left\{  i,i+1,\ldots,j-1\right\}  $, then this would entail
$q=p+1$ (since $w_{i},w_{i+1},\ldots,w_{j-1}$ are distinct), which would
contradict $p+1\neq q$. Thus, $p+1$ cannot be an element of $\left\{
i,i+1,\ldots,j-1\right\}  $. Hence, $p+1=j$ (since $p+1$ clearly belongs to
$\left\{  i,i+1,\ldots,j\right\}  $). Thus, $w_{p+1}=w_{j}=w_{i}$, so that
$w_{i}=w_{p+1}=w_{q}$. This entails $i=q$ (since $w_{i},w_{i+1},\ldots
,w_{j-1}$ are distinct). Hence, $i=q=p-1$ (since $p=q+1$). Therefore,
$\underbrace{j}_{=p+1}-\underbrace{i}_{=p-1}=\left(  p+1\right)  -\left(
p-1\right)  =2$. This contradicts the fact that $j-i$ is odd.
\par
This contradiction shows that our assumption (that the walk $\mathbf{c}$ has
two equal edges) was false. Hence, the edges of $\mathbf{c}$ are distinct.}),
and thus we have found an odd-length cycle contained in $\mathbf{w}$ (namely,
$\mathbf{c}$ is such a cycle, since its length is $j-i$, which is odd). This
means that we are done if $j-i$ is odd.

Thus, we WLOG assume that $j-i$ is even. Hence, cutting out the closed walk
$\mathbf{c}$ from the original walk $\mathbf{w}$, we obtain a walk%
\[
\mathbf{w}^{\prime}:=\left(  w_{0},\ast,w_{1},\ast,\ldots,\ast,w_{i}%
=w_{j},\ast,w_{j+1},\ast,w_{j+2},\ldots,w_{k}\right)
\]
from $u$ to $v$. This new walk $\mathbf{w}^{\prime}$ has length $k-\left(
j-i\right)  $, which is odd (since $k$ is odd but $j-i$ is even) and smaller
than $k$ (since $i<j$). Hence, we can apply the induction hypothesis to this
walk $\mathbf{w}^{\prime}$. As a consequence, we conclude that this walk
$\mathbf{w}^{\prime}$ contains either an odd-length path from $u$ to $v$ or an
odd-length cycle. Therefore, the walk $\mathbf{w}$ also contains either an
odd-length path from $u$ to $v$ or an odd-length cycle (since anything
contained in $\mathbf{w}^{\prime}$ is automatically contained in $\mathbf{w}%
$). But this is precisely what we set out to prove. This completes the
induction step, and so we have proved Proposition
\ref{prop.coloring.2-color-odd-length-cyc}.
\end{proof}

Now, let us prove the 2-coloring equivalence theorem:

\begin{proof}
[Proof of Theorem \ref{thm.coloring.2-color-eq}.]Write the multigraph $G$ as
$G=\left(  V,E,\varphi\right)  $. We shall prove the implications B1
$\Longrightarrow$ B2 $\Longrightarrow$ B3 $\Longrightarrow$ B1. \medskip

\textit{Proof of the implication B1 }$\Longrightarrow$\textit{ B2:} Assume
that Statement B1 holds. We must prove that Statement B2 holds.

We have assumed that B1 holds. In other words, the graph $G$ has a proper
$2$-coloring. Let $f$ be this $2$-coloring. Thus, $f$ is a map from $V$ to
$\left\{  1,2\right\}  $ such that any two adjacent vertices $x$ and $y$ of
$G$ satisfy $f\left(  x\right)  \neq f\left(  y\right)  $.

Assume (for contradiction) that $G$ has a cycle of odd length. Let
\[
\left(  v_{0},\ast,v_{1},\ast,v_{2},\ast,\ldots,\ast,v_{k}\right)
\]
be this cycle. Thus, $k$ is odd, and we have $v_{k}=v_{0}$, so that $f\left(
v_{k}\right)  =f\left(  v_{0}\right)  $. Moreover, for each $i\in\left\{
1,2,\ldots,k\right\}  $, the vertex $v_{i}$ is adjacent to $v_{i-1}$ (since
$\left(  v_{0},\ast,v_{1},\ast,v_{2},\ast,\ldots,\ast,v_{k}\right)  $ is a
cycle) and therefore satisfies
\begin{equation}
f\left(  v_{i}\right)  \neq f\left(  v_{i-1}\right)
\label{pf.thm.coloring.2-color-eq.B1B2.1}%
\end{equation}
(since $f$ is a proper $2$-coloring).

We WLOG assume that $f\left(  v_{0}\right)  =1$ (otherwise, we
\textquotedblleft rename\textquotedblright\ the colors $1$ and $2$ so that the
color $f\left(  v_{0}\right)  $ becomes $1$). Then,
(\ref{pf.thm.coloring.2-color-eq.B1B2.1}) (applied to $i=1$) yields $f\left(
v_{1}\right)  \neq f\left(  v_{0}\right)  =1$, so that $f\left(  v_{1}\right)
=2$ (since $f\left(  v_{1}\right)  $ must be either $1$ or $2$). Hence,
(\ref{pf.thm.coloring.2-color-eq.B1B2.1}) (applied to $i=2$) yields $f\left(
v_{2}\right)  \neq f\left(  v_{1}\right)  =2$, so that $f\left(  v_{2}\right)
=1$ (since $f\left(  v_{2}\right)  $ must be either $1$ or $2$). For similar
reasons, we can successively obtain $f\left(  v_{3}\right)  =2$ and $f\left(
v_{4}\right)  =1$ and $f\left(  v_{5}\right)  =2$ and so on. The general
formula we obtain (strictly speaking, it needs to be proved by induction on
$i$) says that%
\[
f\left(  v_{i}\right)  =%
\begin{cases}
1, & \text{if }i\text{ is even;}\\
2, & \text{if }i\text{ is odd}%
\end{cases}
\ \ \ \ \ \ \ \ \ \ \text{for each }i\in\left\{  0,1,\ldots,k\right\}  .
\]
Applying this to $i=k$, we conclude that $f\left(  v_{k}\right)  =2$ (since
$k$ is odd). However, this contradicts $f\left(  v_{k}\right)  =f\left(
v_{0}\right)  =1\neq2$. This contradiction shows that our assumption was
false. Hence, $G$ has no cycle of odd length. In other words, Statement B2
holds. This proves the implication B1 $\Longrightarrow$ B2. \medskip

\textit{Proof of the implication B2 }$\Longrightarrow$\textit{ B3:} Assume
that Statement B2 holds. We must prove that Statement B3 holds. In other
words, we must show that $G$ has no odd-length circuits.

Assume the contrary. Thus, $G$ has an odd-length circuit $\mathbf{w}$. Let $u$
be the starting and ending point of $\mathbf{w}$. Thus, Proposition
\ref{prop.coloring.2-color-odd-length-cyc} (applied to $v=u$) shows that this
odd-length circuit $\mathbf{w}$ contains either an odd-length \textbf{path}
from $u$ to $u$ or an odd-length cycle. Since $G$ has no odd-length cycle
(because we assumed that Statement B2 holds), we thus concludes that
$\mathbf{w}$ contains an odd-length \textbf{path} from $u$ to $u$. However, an
odd-length path from $u$ to $u$ is impossible (since the only path from $u$ to
$u$ has length $0$). Thus, we obtain a contradiction, which shows that $G$ has
no odd-length circuits. This proves the implication B2 $\Longrightarrow$ B3.
\medskip

\textit{Proof of the implication B3 }$\Longrightarrow$\textit{ B1:} Assume
that Statement B3 holds. We must prove that Statement B1 holds.

We have assumed that Statement B3 holds. In other words, $G$ has no odd-length
circuits. We must find a proper $2$-coloring of $G$.

We WLOG assume that $G$ is connected (otherwise, let $C_{1},C_{2},\ldots
,C_{k}$ be the components of $G$, and apply the implication B3
$\Longrightarrow$ B1 to each of the smaller graphs $G\left[  C_{1}\right]
,\ G\left[  C_{2}\right]  ,\ \ldots,\ G\left[  C_{k}\right]  $, and then
combine the resulting proper $2$-colorings of these smaller graphs into a
single proper $2$-coloring of $G$). Fix any vertex $r$ of $G$. Define a map
$f:V\rightarrow\left\{  1,2\right\}  $ by setting%
\[
f\left(  v\right)  =%
\begin{cases}
1, & \text{if }d\left(  v,r\right)  \text{ is even};\\
2, & \text{if }d\left(  v,r\right)  \text{ is odd}%
\end{cases}
\ \ \ \ \ \ \ \ \ \ \text{for each }v\in V
\]
(where $d\left(  v,r\right)  $ denotes the distance from $v$ to $r$, that is,
the smallest length of a path from $v$ to $r$).

I claim that $f$ is a proper $2$-coloring.\footnote{Here is an illustrative
example:%
\[%
\begin{tikzpicture}[scale=1.2]
\begin{scope}[every node/.style={circle,thick,draw=green!60!black}]
\node(00) at (0, 0) {$1$};
\node(20) at (2, 0) {$2$};
\node(40) at (4, 0) {$2$};
\node(60) at (6, 0) {$1$};
\node(11) at (1, 1) {$2$};
\node(31) at (3, 1) {$1$};
\node(51) at (5, 1) {$2$};
\node(1-1) at (1, -1) {$2$};
\node(3-1) at (3, -1) {$1$};
\node(5-1) at (5, -1) {$2$};
\end{scope}
\begin{scope}[every edge/.style={draw=black,very thick}, every loop/.style={}]
\path[-] (00) edge (11) edge (20) edge (1-1);
\path[-] (31) edge (11) edge (20) edge (40) edge (51);
\path[-] (3-1) edge (1-1) edge (20) edge (40) edge (5-1);
\path[-] (60) edge (51) edge (5-1);
\end{scope}
\end{tikzpicture}%
\ \ .
\]
(Of course, the numbers on the nodes here are not the vertices, but rather the
colors of these vertices.)
\par
Note that all values of $f$ can be easily found by the following recursive
algorithm: Start by assigning the color $1$ to $r$. Then, assign the color $2$
to all neighbors of $r$. Then, assign the color $1$ to all neighbors of these
neighbors (unless they have already been colored). Then, assign the color $2$
to all neighbors of these neighbors of these neighbors, and so on.} Indeed,
assume the contrary. Thus, some two adjacent vertices $u$ and $v$ have the
same color $f\left(  u\right)  =f\left(  v\right)  $. Consider these $u$ and
$v$. Since $f\left(  u\right)  =f\left(  v\right)  $, we are in one of the
following two cases:

\textit{Case 1:} We have $f\left(  u\right)  =f\left(  v\right)  =1$.

\textit{Case 2:} We have $f\left(  u\right)  =f\left(  v\right)  =2$.

Let us consider Case 2. In this case, we have $f\left(  u\right)  =f\left(
v\right)  =2$. This means that $d\left(  u,r\right)  $ and $d\left(
v,r\right)  $ are both odd (by the definition of $f$). Hence, there is an
odd-length path $\mathbf{p}$ from $u$ to $r$ and an odd-length path
$\mathbf{q}$ from $v$ to $r$. Consider these $\mathbf{p}$ and $\mathbf{q}$.
Also, there is an edge $e$ that joins $u$ and $v$ (since $u$ and $v$ are
adjacent). Consider this edge $e$. By combining the paths $\mathbf{p}$ and
$\mathbf{q}$ and inserting the edge $e$ into the result, we obtain a circuit
from $r$ to $r$ (which starts by following the path $\mathbf{p}$ backwards to
$u$, then takes the edge $e$ to $v$, then follows the path $\mathbf{q}$ back
to $r$). This circuit has odd length (since $\mathbf{p}$ and $\mathbf{q}$ have
odd lengths, and since the edge $e$ adds $1$ to the length). Thus, we have
found an odd-length circuit of $G$. However, we assumed that $G$ has no
odd-length circuits. Contradiction!

Thus, we have found a contradiction in Case 2. Similarly, we can find a
contradiction in Case 1. Thus, we always get a contradiction. This shows that
$f$ is indeed a proper $2$-coloring. Thus, Statement B1 holds. This proves the
implication B3 $\Longrightarrow$ B1.\ \ \ \ \footnote{Note that this proof
provides a reasonably efficient algorithm for constructing a proper
$2$-coloring of $G$, as long as you know how to compute distances in a graph
(we have done this, e.g., in homework set \#4 exercise 5) and how to compute
the components of a graph (this is not hard).} \medskip

\begin{fineprint}
For aesthetical reasons, let me give a \textit{second proof of the implication
B3 }$\Longrightarrow$\textit{ B1}, which avoids the awkward \textquotedblleft
break $G$ up into components\textquotedblright\ step\textit{:}

Assume again that Statement B3 holds. We must prove that Statement B1 holds.

We assumed that Statement B3 holds. In other words, $G$ has no odd-length cycles.

Two vertices $u$ and $v$ of $G$ will be called \textbf{oddly connected} if $G$
has an odd-length path from $u$ to $v$. By Proposition
\ref{prop.coloring.2-color-odd-length-cyc}, this condition is equivalent to
\textquotedblleft$G$ has an odd-length walk from $u$ to $v$\textquotedblright,
since $G$ has no odd-length cycles. Moreover, a vertex $u$ cannot be oddly
connected to itself (since the only path from $u$ to $u$ is the trivial
length-$0$ path $\left(  u\right)  $, which is not odd-length).

A subset $A$ of $V$ will be called \textbf{odd-path-less} if no two vertices
in $A$ are oddly connected. (Note that \textquotedblleft two
vertices\textquotedblright\ doesn't mean \textquotedblleft two distinct
vertices\textquotedblright.)

Pick a maximum-size odd-path-less subset $A$ of $V$ (such an $A$ exists, since
$\varnothing$ is clearly odd-path-less). Now, let $f:V\rightarrow\left\{
1,2\right\}  $ be the $2$-coloring of $G$ that assigns the color $1$ to all
vertices in $A$ and assigns the color $2$ to all vertices not in $A$.

We shall show that this $2$-coloring $f$ is proper.

To prove this, we must show that no two adjacent vertices have color $1$ and
that no two adjacent vertices have color $2$. The first of these two claims is
obvious\footnote{\textit{Proof.} An edge always makes a walk of length $1$,
which is odd. Thus, two adjacent vertices are automatically oddly connected.
Hence, two adjacent vertices cannot both be contained in the odd-path-less
subset $A$. In other words, two adjacent vertices cannot both have color $1$%
.}. It thus remains to prove the second claim -- i.e., to prove that no two
adjacent vertices have color $2$.

Assume the contrary. Thus, there exist two adjacent vertices $u$ and $v$ that
both have color $2$. Consider these $u$ and $v$. These vertices $u$ and $v$
have color $2$; in other words, neither of them belongs to $A$.

The vertex $u$ is not oddly connected to itself (as we already saw). Hence,
the vertex $u$ is oddly connected to at least one vertex $a\in A$ (because
otherwise, we could insert $u$ into the odd-path-less set $A$ and obtain a
larger odd-path-less subset $A\cup\left\{  u\right\}  $ of $V$; but this would
contradict the fact that $A$ is a \textbf{maximum-size} odd-path-less subset
of $V$). For similar reasons, the vertex $v$ is oddly connected to at least
one vertex $b\in A$. Consider these vertices $a$ and $b$. Since $u$ is oddly
connected to $a$, there exists an odd-length walk $\mathbf{p}$ from $u$ to
$a$. Reversing this walk $\mathbf{p}$ yields an odd-length walk $\mathbf{p}%
^{\prime}$ from $a$ to $u$. Since $v$ is oddly connected to $b$, there exists
an odd-length walk $\mathbf{q}$ from $v$ to $b$. Finally, there is an edge $e$
with endpoints $u$ and $v$ (since $u$ and $v$ are adjacent). Combine the two
walks $\mathbf{p}^{\prime}$ and $\mathbf{q}$ and insert this edge $e$ between
them; this yields a walk from $a$ to $b$ (via $u$ and $v$) that has odd length
(since $\mathbf{p}^{\prime}$ and $\mathbf{q}$ have odd length each, and
inserting $e$ adds $1$ to the length). Thus, $G$ has an odd-length walk from
$a$ to $b$. In other words, the vertices $a$ and $b$ are oddly connected. This
contradicts the fact that the set $A$ is odd-path-less (since $a$ and $b$
belong to $A$).

This contradiction shows that our assumption was false. Thus, we have shown
that no two adjacent vertices have color $2$. This completes our proof that
$f$ is a proper $2$-coloring. Thus, Statement B1 holds. This proves the
implication B3 $\Longrightarrow$ B1 once again. \medskip
\end{fineprint}

Having proved all three implications B1 $\Longrightarrow$ B2 and B2
$\Longrightarrow$ B3 and B3 $\Longrightarrow$ B1, we now conclude that the
three statements B1, B2 and B3 are equivalent. This proves Theorem
\ref{thm.coloring.2-color-eq}.
\end{proof}

\begin{remark}
A graph $G$ that satisfies the three equivalent statements B1, B2, B3 of
Theorem \ref{thm.coloring.2-color-eq} is sometimes called a \textquotedblleft
bipartite graph\textquotedblright. This is slightly imprecise, since the
proper definition of a \textquotedblleft bipartite graph\textquotedblright\ is
(equivalent to) \textquotedblleft a graph \textbf{equipped with} a proper
$2$-coloring\textquotedblright. Thus, if we equip one and the same graph $G$
with different proper $2$-colorings, then we obtain different bipartite
graphs. We shall take a closer look at bipartite graphs in Sections
\ref{sec.matching.bip}, \ref{sec.matching.hall} and \ref{sec.matching.koenig}.
\end{remark}

A further simple property of proper $2$-colorings is the
following:\footnote{Recall that $\operatorname*{conn}G$ denotes the number of
components of a graph $G$.}

\begin{proposition}
\label{prop.coloring.2-color-count}Let $G$ be a multigraph that has a proper
$2$-coloring. Then, $G$ has exactly $2^{\operatorname*{conn}G}$ many proper
$2$-colorings.
\end{proposition}

\begin{proof}
[Proof sketch.]For each component $C$ of $G$, let us fix an arbitrary vertex
$r_{C}\in C$. When constructing a proper $2$-coloring $f$ of $G$, we can
freely choose the colors $f\left(  r_{C}\right)  $ of these vertices $r_{C}$;
the colors of all other vertices are then uniquely determined (see the first
proof of the implication B3 $\Longrightarrow$ B1 in our above proof of Theorem
\ref{thm.coloring.2-color-eq} for the details). Thus, we have
$2^{\operatorname*{conn}G}$ many options (since $G$ has $\operatorname*{conn}%
G$ many components). The proposition follows.
\end{proof}

\subsubsection{The directed case}

There is an interesting variant of Theorem \ref{thm.coloring.2-color-eq} for
directed graphs (\cite[Lemma 1.1.5]{Langlo23}):

\begin{theorem}
\label{thm.coloring.2-color-eq-dig}Let $D$ be a strongly connected
multidigraph. Then, the following three statements are equivalent:

\begin{itemize}
\item \textbf{Statement B'1:} The underlying undirected graph
$D^{\operatorname*{und}}$ has a proper $2$-coloring.

\item \textbf{Statement B'2:} The digraph $D$ has no cycles of odd length.

\item \textbf{Statement B'3:} The digraph $D$ has no circuits of odd length.
\end{itemize}
\end{theorem}

Keep in mind that cycles and circuits in a digraph are directed cycles,
traversing each arc from source to target; thus, a digraph $D$ will usually
have a lot fewer cycles than the underlying undirected graph
$D^{\operatorname*{und}}$. So the implication B'2 $\Longrightarrow$ B'1 in
Theorem \ref{thm.coloring.2-color-eq-dig} is much stronger than the
corresponding implication B2 $\Longrightarrow$ B1 in Theorem
\ref{thm.coloring.2-color-eq}, seeing that statement B'2 looks weaker than B2
for $G=D^{\operatorname*{und}}$. The \textquotedblleft strongly
connected\textquotedblright\ requirement in Theorem
\ref{thm.coloring.2-color-eq-dig} is the price to pay for this extra strength:
No connectedness was needed in Theorem \ref{thm.coloring.2-color-eq}, but we
do need to require $D$ to be strongly connected in Theorem
\ref{thm.coloring.2-color-eq-dig}, since otherwise the implication B'2
$\Longrightarrow$ B'1 would fail for the first digraph in Example
\ref{exa.tour.exas1} \textbf{(a)} (which has no cycles at all, but whose
underlying undirected graph $K_{3}$ certainly has no proper $2$-coloring).

The equivalence B'2 $\Longleftrightarrow$ B'3 in Theorem
\ref{thm.coloring.2-color-eq-dig} does indeed hold even when $D$ is not
strongly connected. This follows easily from the following directed analogue
of Proposition \ref{prop.coloring.2-color-odd-length-cyc}:

\begin{proposition}
\label{prop.coloring.2-color-odd-length-cyc-dg}Let $D$ be a multidigraph. Let
$u$ and $v$ be two vertices of $D$. Let $\mathbf{w}$ be an odd-length walk
from $u$ to $v$. Then, $\mathbf{w}$ contains either an odd-length
\textbf{path} from $u$ to $v$ or an odd-length \textbf{cycle} (or both).
\end{proposition}

\begin{proof}
Literally the same as for Proposition
\ref{prop.coloring.2-color-odd-length-cyc}. (It is even a bit easier, because
proving that a closed walk is a cycle in a multidigraph does not require
showing that its arcs are distinct.)
\end{proof}

\begin{corollary}
\label{cor.coloring.2-color-eq-dig-cyc-dg}Let $D$ be a multidigraph. Then, $D$
has no cycles of odd length if and only if $D$ has no circuits of odd length.
\end{corollary}

\begin{proof}
$\Longrightarrow:$ Assume that $D$ has no cycles of odd length. We must show
that $D$ has no circuits of odd length.

Let $\mathbf{c}$ be a circuit of odd length. Let $u$ be its starting point.
Then, $u$ is also its ending point (since $\mathbf{c}$ is a circuit). Hence,
$\mathbf{c}$ is an odd-length walk from $u$ to $u$. Thus, Proposition
\ref{prop.coloring.2-color-odd-length-cyc-dg} (applied to $\mathbf{w}%
=\mathbf{c}$ and $v=u$) shows that $\mathbf{c}$ contains either an odd-length
\textbf{path} from $u$ to $u$ or an odd-length \textbf{cycle} (or both). The
former option is impossible, since the only path from $u$ to $u$ is the
length-$0$ path $\left(  u\right)  $, which is certainly not odd-length. But
the latter option is also impossible, since we have assumed that $D$ has no
cycles of odd length. Thus, we obtain a contradiction in each case.

So we have found a contradiction whenever $\mathbf{c}$ is a circuit of odd
length. Hence, $D$ has no circuits of odd length. This proves the
\textquotedblleft$\Longrightarrow$\textquotedblright\ direction of Corollary
\ref{cor.coloring.2-color-eq-dig-cyc-dg}. \medskip

$\Longleftarrow:$ Obvious (since any cycle is a circuit).
\end{proof}

We can now prove Theorem \ref{thm.coloring.2-color-eq-dig}:

\begin{proof}
[Proof of Theorem \ref{thm.coloring.2-color-eq-dig}.]Corollary
\ref{cor.coloring.2-color-eq-dig-cyc-dg} shows that statements B'2 and B'3 are equivalent.

We shall now prove the implications B'1 $\Longrightarrow$ B'3 $\Longrightarrow
$ B'1. \medskip

\textit{Proof of the implication B'1 }$\Longrightarrow$\textit{ B'3:} Assume
that Statement B'1 holds. We must prove that Statement B'3 holds.

We have assumed that B'1 holds. In other words, the graph
$D^{\operatorname*{und}}$ has a proper $2$-coloring. Hence, the implication B1
$\Longrightarrow$ B3 of Theorem \ref{thm.coloring.2-color-eq} (applied to
$G=D^{\operatorname*{und}}$) shows that the graph $D^{\operatorname*{und}}$
has no circuits of odd length. Hence, the digraph $D$ has no circuits of odd
length either (since each circuit of $D$ is also a circuit of
$D^{\operatorname*{und}}$). In other words, Statement B'3 holds. This proves
the implication B'1 $\Longrightarrow$ B'3. \medskip

\textit{Proof of the implication B'3 }$\Longrightarrow$ \textit{B'1:} Assume
that Statement B'3 holds. We must prove that Statement B'1 holds.

We have assumed that B'3 holds. In other words, the digraph $D$ has no
circuits of odd length.

Let $V$ be the vertex set of $D$. For any two vertices $u,v\in V$, we let
$d\left(  u,v\right)  $ be the smallest length of a path from $u$ to $v$ in
$D$ (such a path exists, since $D$ is strongly connected). Now we claim the following:

\begin{statement}
\textit{Claim 1:} Let $u,v\in V$ be two vertices. Let $\mathbf{w}$ be any walk
from $v$ to $u$ in $D$. Then, the length of $\mathbf{w}$ is $\equiv d\left(
u,v\right)  \operatorname{mod}2$.
\end{statement}

\begin{proof}
[Proof of Claim 1.]Let $k$ be the length of $\mathbf{w}$. We must then show
that $k\equiv d\left(  u,v\right)  \operatorname{mod}2$.

The digraph $D$ has a path $\mathbf{p}$ from $u$ to $v$ whose length is
$d\left(  u,v\right)  $ (by the definition of $d\left(  u,v\right)  $).
Consider this $\mathbf{p}$. Splicing the path $\mathbf{p}$ (from $u$ to $v$)
with the walk $\mathbf{w}$ (from $v$ to $u$), we obtain a walk $\mathbf{p}%
\ast\mathbf{w}$ from $u$ to $u$. This walk $\mathbf{p}\ast\mathbf{w}$ is a
circuit (being a walk from $u$ to $u$), and thus cannot have odd length (since
$D$ has no circuits of odd length). Hence, $\mathbf{p}\ast\mathbf{w}$ has even
length. But the length of $\mathbf{p\ast w}$ is $d\left(  u,v\right)  +k$
(since $\mathbf{p}$ has length $d\left(  u,v\right)  $ while $\mathbf{w}$ has
length $k$). Combining the previous two sentences, we conclude that $d\left(
u,v\right)  +k$ is even. In other words, $d\left(  u,v\right)  \equiv-k\equiv
k\operatorname{mod}2$, so that $k\equiv d\left(  u,v\right)
\operatorname{mod}2$. This proves Claim 1.
\end{proof}

Now, pick any vertex $r$ of $D$ (this exists, since $D$ is strongly connected
and thus has at least $1$ vertex). Define a map $f:V\rightarrow\left\{
1,2\right\}  $ by setting%
\[
f\left(  v\right)  =%
\begin{cases}
1, & \text{if }d\left(  v,r\right)  \text{ is even};\\
2, & \text{if }d\left(  v,r\right)  \text{ is odd}%
\end{cases}
\ \ \ \ \ \ \ \ \ \ \text{for each }v\in V.
\]
I claim that $f$ is a proper $2$-coloring of $D^{\operatorname*{und}}$.
Indeed, assume the contrary. Thus, some two adjacent vertices $u$ and $v$ of
$D^{\operatorname*{und}}$ have the same color $f\left(  u\right)  =f\left(
v\right)  $. Consider these $u$ and $v$. Since $f\left(  u\right)  =f\left(
v\right)  $, we are in one of the following two cases:

\textit{Case 1:} We have $f\left(  u\right)  =f\left(  v\right)  =1$.

\textit{Case 2:} We have $f\left(  u\right)  =f\left(  v\right)  =2$.

Let us consider Case 2. In this case, we have $f\left(  u\right)  =f\left(
v\right)  =2$. This means that $d\left(  u,r\right)  $ and $d\left(
v,r\right)  $ are both odd (by the definition of $f$).

The digraph $D$ is strongly connected, and thus has a path $\mathbf{p}$ from
$r$ to $u$ as well as a path from $\mathbf{q}$ from $r$ to $v$. Consider these
paths $\mathbf{p}$ and $\mathbf{q}$. Claim 1 (applied to $\mathbf{p}$ and $r$
instead of $\mathbf{w}$ and $v$) yields that the length of the path
$\mathbf{p}$ is $\equiv d\left(  u,r\right)  \operatorname{mod}2$. In other
words, the length of the path $\mathbf{p}$ is odd (since $d\left(  u,r\right)
$ is odd). Likewise, the length of the path $\mathbf{q}$ is odd.

But the fact that $u$ and $v$ are adjacent in $D^{\operatorname*{und}}$ shows
that the digraph $D$ has either an arc from $u$ to $v$ or an arc from $v$ to
$u$ (or both). We WLOG assume that $D$ has an arc from $u$ to $v$ (otherwise,
we can just swap $u$ with $v$). Let $a$ be this arc. Attaching the arc $a$ to
the path $\mathbf{p}$ from $r$ to $u$ at its end, we obtain a walk
$\mathbf{w}$ from $r$ to $v$, and moreover this walk has even length (since
the path $\mathbf{p}$ has odd length, and we have obtained $\mathbf{w}$ by
attaching one more arc to it).

But Claim 1 (applied to $v$ and $r$ instead of $u$ and $v$) shows that the
length of $\mathbf{w}$ is $\equiv d\left(  v,r\right)  \operatorname{mod}2$.
In other words, the length of $\mathbf{w}$ is odd (since $d\left(  v,r\right)
$ is odd). But this contradicts the fact that $\mathbf{w}$ has even length.

Thus, we have found a contradiction in Case 2. Similarly, we can find a
contradiction in Case 1. Thus, we always get a contradiction. This shows that
$f$ is indeed a proper $2$-coloring of $D^{\operatorname*{und}}$. Thus,
Statement B'1 holds. This proves the implication B'3 $\Longrightarrow$ B'1.
\medskip

Having proved the implications B'1 $\Longrightarrow$ B'3 and B'3
$\Longrightarrow$ B'1 and the equivalence B'2 $\Longleftrightarrow$ B'3, we
conclude that all three statements B'1, B'2 and B'3 are equivalent. Thus,
Theorem \ref{thm.coloring.2-color-eq-dig} holds.
\end{proof}

\subsection{The Brooks theorems}

As we said, the existence of a proper $k$-coloring for a given graph $G$ is a
hard computational problem unless $k\leq2$. The same holds for theoretical
criteria: For $k>2$, I am not aware of any good criteria that are
simultaneously necessary and sufficient for the existence of a proper
$k$-coloring. However, some sufficient criteria are known. Here is
one:\footnote{Recall that a multigraph is called \textbf{loopless} if it has
no loops.}

\begin{theorem}
[Little Brooks theorem]\label{thm.coloring.little-brooks}Let $G=\left(
V,E,\varphi\right)  $ be a loopless multigraph with at least one vertex. Let%
\[
\alpha:=\max\left\{  \deg v\ \mid\ v\in V\right\}  .
\]
Then, $G$ has a proper $\left(  \alpha+1\right)  $-coloring.
\end{theorem}

\begin{proof}
[Proof sketch.]Let $v_{1},v_{2},\ldots,v_{n}$ be the vertices of $V$, listed
in some order (with no repetitions). We construct a proper $\left(
\alpha+1\right)  $-coloring $f:V\rightarrow\left\{  1,2,\ldots,\alpha
+1\right\}  $ of $G$ recursively as follows:

\begin{itemize}
\item First, we choose $f\left(  v_{1}\right)  $ arbitrarily.

\item Then, we choose $f\left(  v_{2}\right)  $ to be distinct from the colors
of all already-colored neighbors of $v_{2}$.

\item Then, we choose $f\left(  v_{3}\right)  $ to be distinct from the colors
of all already-colored neighbors of $v_{3}$.

\item Then, we choose $f\left(  v_{4}\right)  $ to be distinct from the colors
of all already-colored neighbors of $v_{4}$.

\item And so on, until all values $f\left(  v_{1}\right)  ,\ f\left(
v_{2}\right)  ,\ \ldots,\ f\left(  v_{n}\right)  $ have been chosen.
\end{itemize}

\noindent Why do we never run out of colors in this process? Well: When
choosing $f\left(  v_{i}\right)  $, we must choose a color distinct from the
colors of all already-colored neighbors of $v_{i}$. Since $v_{i}$ has at most
$\alpha$ neighbors (because $\deg\left(  v_{i}\right)  \leq\alpha$), this
means that we have at most $\alpha$ colors to avoid. Since there are
$\alpha+1$ colors in total, this leaves us at least $1$ color that we can
choose; therefore, we don't run out of colors.

The resulting $\left(  \alpha+1\right)  $-coloring $f:V\rightarrow\left\{
1,2,\ldots,\alpha+1\right\}  $ is called a \textbf{greedy coloring}. This
$\left(  \alpha+1\right)  $-coloring $f$ is indeed proper, because if an edge
has endpoints $v_{i}$ and $v_{j}$ with $i>j$, then the construction of
$f\left(  v_{i}\right)  $ ensures that $f\left(  v_{i}\right)  $ is distinct
from $f\left(  v_{j}\right)  $. (Note how we are using the fact that $G$ is
loopless here! If $G$ had a loop, then the endpoints of this loop could not be
written as $v_{i}$ and $v_{j}$ with $i>j$.)
\end{proof}

In general, the $\alpha+1$ in Theorem \ref{thm.coloring.little-brooks} cannot
be improved. Here are two examples:

\begin{itemize}
\item If $n\geq2$, then the cycle graph $C_{n}$\ \ \ \ \footnote{See
Definition \ref{def.mg.Cn} for the proper definition of $C_{n}$ when $n=2$.}
has maximum degree \newline$\alpha=\max\left\{  \deg v\ \mid\ v\in V\right\}
=2$. Thus, Theorem \ref{thm.coloring.little-brooks} shows that $C_{n}$ has a
proper $3$-coloring. When $n$ is even, $C_{n}$ has a proper $2$-coloring as
well, but this is not the case when $n$ is odd (by Theorem
\ref{thm.coloring.2-color-eq}).

\item If $n\geq1$, then the complete graph $K_{n}$ has maximum degree
\newline$\alpha=\max\left\{  \deg v\ \mid\ v\in V\right\}  =n-1$. Thus,
Theorem \ref{thm.coloring.little-brooks} shows that $K_{n}$ has a proper
$n$-coloring. By the pigeonhole principle, it is clear that $K_{n}$ has no
proper $\left(  n-1\right)  $-coloring.
\end{itemize}

Interestingly, these two examples are in fact the only cases when a connected
loopless multigraph with maximum degree $\alpha$ can fail to have a proper
$\alpha$-coloring. In all other cases, we can improve the $\alpha+1$ to
$\alpha$:

\begin{theorem}
[Brooks theorem]\label{thm.coloring.brooks}Let $G=\left(  V,E,\varphi\right)
$ be a connected loopless multigraph. Let%
\[
\alpha:=\max\left\{  \deg v\ \mid\ v\in V\right\}  .
\]
Assume that $G$ is neither a complete graph nor an odd-length cycle. Then, $G$
has a proper $\alpha$-coloring.
\end{theorem}

\begin{proof}
Despite the seemingly little difference, this is significantly harder to prove
than Theorem \ref{thm.coloring.little-brooks}. Various proofs can be found in
\cite{CraRab15} and in most serious textbooks on graph theory.
\end{proof}

\subsection{Exercises on proper colorings}

\begin{exercise}
\label{exe.7.9}Let $G$ be a simple graph with $n$ vertices. Let $k$ be a
positive integer.

Prove the following:

\begin{enumerate}
\item[\textbf{(a)}] If $G$ has a proper $k$-coloring, then $G$ has no subgraph
isomorphic to $K_{k+1}$.

\item[\textbf{(b)}] If $k\geq n-2$, then the converse to part \textbf{(a)}
also holds: If $G$ has no subgraph isomorphic to $K_{k+1}$, then $G$ has a
proper $k$-coloring.

\item[\textbf{(c)}] Does the converse to part \textbf{(a)} hold for $k<n-2$ as
well? Specifically, does it hold for $n=5$ and $k=2$ ?
\end{enumerate}
\end{exercise}

\begin{exercise}
\label{exe.7.10}Let $G$ be a connected loopless multigraph. Prove that $G$ has
a proper $2$-coloring if and only if every three vertices $u,v,w$ of $G$
satisfy
\[
2\mid d\left(  u,v\right)  +d\left(  v,w\right)  +d\left(  w,u\right)  .
\]

\end{exercise}

\begin{exercise}
\label{exe.7.7}Fix two positive integers $n$ and $k$ with $n\geq2k>0$. Let
$S=\left\{  1,2,\ldots,n\right\}  $. Consider the $k$-Kneser graph $K_{S,k}$
as defined in Subsection \ref{subsec.sg.complete.kneser}. Prove that $K_{S,k}$
has a proper $\left(  n-2k+2\right)  $-coloring. \medskip

[\textbf{Hint:} What can you say about the minima (i.e., smallest elements) of
two disjoint subsets of $S$? (Being distinct is a good first step.)] \medskip

[\textbf{Remark:} L\'{o}vasz has proved in 1978 (using topology!) that this
result is optimal -- in the sense that $n-2k+2$ is the smallest integer $q$
such that $K_{S,k}$ has a proper $q$-coloring.]
\end{exercise}

\begin{exercise}
\label{exe.7.12}Let $k\in\mathbb{N}$. Let $G$ and $H$ be two simple graphs.
Assume that each of $G$ and $H$ has a proper $k$-coloring. Prove that the
Cartesian product $G\times H$ (defined in Definition \ref{def.sg.cartprod})
has a proper $k$-coloring as well. \medskip

[\textbf{Remark:} It is easy to see that the converse holds as well (i.e., if
$G\times H$ has a proper $k$-coloring, then so do $G$ and $H$), provided that
the vertex sets $\operatorname*{V}\left(  G\right)  $ and $\operatorname*{V}%
\left(  H\right)  $ are both nonempty.]
\end{exercise}

\begin{exercise}
\label{exe.7.13}Let $n\in\mathbb{N}$. Let $G$ be the $n$-th coprimality graph
$\operatorname*{Cop}\nolimits_{n}$ defined in Example \ref{exa.sg.sg2}. Let
$k\in\mathbb{N}$. Let $m$ be the number of prime numbers in the set $\left\{
1,2,\ldots,n\right\}  $. Prove the following:

\begin{enumerate}
\item[\textbf{(a)}] The graph $G$ has a proper $k$-coloring if and only if
$k\geq m+1$.

\item[\textbf{(b)}] The graph $G$ has a subgraph isomorphic to $K_{k}$ if and
only if $k\leq m+1$.
\end{enumerate}
\end{exercise}

\begin{exercise}
\label{exe.7.8}Let $n$ and $k$ be two positive integers. Let $K$ be a set of
size $k$. Let $D$ be the de Bruijn digraph -- i.e., the multidigraph
constructed in the proof of Theorem \ref{thm.debr.exist}. Let $G$ be the
result of removing all loops from the undirected graph $D^{\operatorname{und}%
}$. Prove that $G$ has a proper $\left(  k+1\right)  $-coloring.
\end{exercise}

\begin{exercise}
\label{exe.7.14}Let $k\in\mathbb{N}$. Let $G$ be a simple graph with fewer
than $\dfrac{k\left(  k+1\right)  }{2}$ edges. Prove that $G$ has a proper $k$-coloring.
\end{exercise}

Colorings and edge-colorings can be interesting even when they are not proper.
Here is a sample result:

\begin{exercise}
Let $n$ be a positive integer. Let $E$ be the edge set of the complete graph
$K_{n}$. A map $f:E\rightarrow\mathbb{Z}$ is said to be an
\textbf{edge-coloring} of $K_{n}$. The values of this map $f$ are called the
\textbf{colors} of the edge-coloring $f$. (We imagine each value $f\left(
e\right)  $ to be a color assigned to the edge $e$.)

An edge-coloring $f:E\rightarrow\mathbb{Z}$ of $K_{n}$ is said to be
\textbf{Gallai} if for any three distinct vertices $i,j,k$, at least two of
the three numbers $f\left(  ij\right)  $, $f\left(  jk\right)  $ and $f\left(
ki\right)  $ are equal (where we use the shorthand $uv$ for the edge $\left\{
u,v\right\}  $, as usual). (Visually speaking, this means that each triangle
has at least two equally-colored sides.)

Prove that the maximum number of colors (i.e., the maximum size of $f\left(
E\right)  $) in a Gallai edge-coloring $f$ of $K_{n}$ is $n-1$.
\end{exercise}

\subsection{The chromatic polynomial}

Here is another surprise: The number of proper $k$-colorings of a given
multigraph $G$ turns out to be a polynomial function in $k$ (with integer
coefficients). More precisely:

\begin{theorem}
[Whitney's chromatic polynomial theorem]\label{thm.coloring.whitney}Let
$G=\left(  V,E,\varphi\right)  $ be a multigraph. Let $\chi_{G}$ be the
polynomial in the single indeterminate $x$ with coefficients in $\mathbb{Z}$
defined as follows:%
\[
\chi_{G}=\sum_{F\subseteq E}\left(  -1\right)  ^{\left\vert F\right\vert
}x^{\operatorname*{conn}\left(  V,F,\varphi\mid_{F}\right)  }=\sum
_{\substack{H\text{ is a spanning}\\\text{subgraph of }G}}\left(  -1\right)
^{\left\vert \operatorname*{E}\left(  H\right)  \right\vert }%
x^{\operatorname*{conn}H}.
\]
(The symbol ``$\sum_{F\subseteq E}$''\ means ``sum over all subsets $F$ of $E$''.)

Then, for any $k\in\mathbb{N}$, we have%
\[
\left(  \text{\# of proper }k\text{-colorings of }G\right)  =\chi_{G}\left(
k\right)  .
\]

\end{theorem}

The proper place for this theorem is probably a course on enumerative
combinatorics, but let us give here a proof for the sake of completeness
(optional material). The following proof is essentially due to Hassler Whitney
in 1930 (\cite[\S 6]{Whitney32}), and I am mostly copypasting it from my own
writeup \cite[\S 0.5]{17s-mt2s} (with some changes stemming from the fact that
we are here working with multigraphs rather than simple graphs).

We are going to use the \textbf{Iverson bracket notation}:

\begin{definition}
If $\mathcal{A}$ is any logical statement, then $\left[  \mathcal{A}\right]  $
shall denote the truth value of $\mathcal{A}$; this is the number $%
\begin{cases}
1, & \text{if }\mathcal{A}\text{ is true};\\
0, & \text{if }\mathcal{A}\text{ is false}.
\end{cases}
$

For instance, $\left[  2+2=4\right]  =1$ and $\left[  2+2=5\right]  =0$.
\end{definition}

We next recall a well-known combinatorial identity (see \cite[Lemma
3.3.5]{17s} or \cite[Proposition 2.9.10]{19fco}):

\begin{lemma}
\label{lem.dominating.heinrich-lemma1} Let $P$ be a finite set. Then,
\[
\sum_{\substack{A\subseteq P}}\left(  -1\right)  ^{\left\vert A\right\vert
}=\left[  P=\varnothing\right]  .
\]
(The symbol \textquotedblleft$\sum_{\substack{A\subseteq P}}$%
\textquotedblright\ means \textquotedblleft sum over all subsets $A$ of
$P$\textquotedblright.)
\end{lemma}

Next, we introduce a specific notation related to colorings:

\begin{definition}
\label{def.mt2.chrompoly.Ef} Let $G=\left(  V,E,\varphi\right)  $ be a
multigraph. Let $k\in\mathbb{N}$. Let $f:V\rightarrow\left\{  1,2,\ldots
,k\right\}  $ be a $k$-coloring. We then define a subset $E_{f}$ of $E$ by%
\[
E_{f}:=\left\{  e\in E\ \mid\ \text{the two endpoints of }e\text{ have the
same color in }f\right\}  .
\]
(Recall that the \textquotedblleft color in $f$\textquotedblright\ of a vertex
$v$ means the value $f\left(  v\right)  $. If an edge $e\in E$ is a loop, then
$e$ always belongs to $E_{f}$, since we think of the two endpoints of $e$ as
being equal.)

The elements of $E_{f}$ are called the $f$\textbf{-monochromatic} edges of
$G$. (\textquotedblleft Monochromatic\textquotedblright\ means
\textquotedblleft one-colored\textquotedblright, so no surprises here.)
\end{definition}

\begin{example}
Let $G=\left(  V,E,\varphi\right)  $ be the following multigraph:%
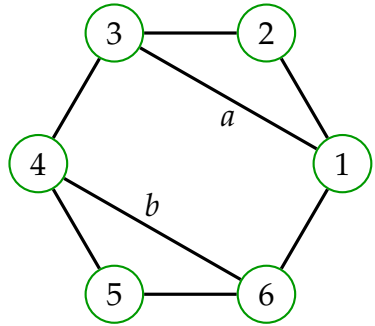
\[%
\begin{tikzpicture}
\begin{scope}[every node/.style={circle,thick,draw=green!60!black}]
\node(A) at (0:2) {$1$};
\node(B) at (60:2) {$2$};
\node(C) at (120:2) {$3$};
\node(D) at (180:2) {$4$};
\node(E) at (240:2) {$5$};
\node(F) at (300:2) {$6$};
\end{scope}
\begin{scope}[every edge/.style={draw=black,very thick}]
\path
[-] (A) edge (B) (B) edge (C) (C) edge (D) (D) edge (E) (E) edge (F) (F) edge (A);
\path[-] (C) edge node[below] {$a$} (A);
\path[-] (D) edge node[above] {$b$} (F);
\end{scope}
\end{tikzpicture}%
\ \ .
\]
Let $f:V\rightarrow\left\{  1,2\right\}  $ be the $2$-coloring of $G$ that
sends each odd vertex to $1$ and each even vertex to $2$. (Here, an
\textquotedblleft odd vertex\textquotedblright\ means a vertex that is odd as
an integer. Thus, the odd vertices are $1,3,5$. \textquotedblleft Even
vertices\textquotedblright\ are understood similarly.) Then, $E_{f}=\left\{
a,b\right\}  $.
\end{example}

Notice the following simple fact:

\begin{proposition}
\label{prop.mt2.chrompoly.EcapEf} Let $G=\left(  V,E,\varphi\right)  $ be a
multigraph. Let $k\in\mathbb{N}$. Let $f:V\rightarrow\left\{  1,2,\ldots
,k\right\}  $ be a $k$-coloring. Then, the $k$-coloring $f$ is proper if and
only if $E_{f}=\varnothing$.
\end{proposition}

\begin{proof}
[Proof of Proposition~\ref{prop.mt2.chrompoly.EcapEf}.]We have the following
chain of equivalences:
\begin{align*}
&  \ \left(  \text{the }k\text{-coloring }f\text{ is proper}\right) \\
&  \Longleftrightarrow\ \left(  \text{no two adjacent vertices have the same
color}\right) \\
&  \ \ \ \ \ \ \ \ \ \ \ \ \ \ \ \ \ \ \ \ \left(  \text{by the definition of
\textquotedblleft proper\textquotedblright}\right) \\
&  \Longleftrightarrow\ \left(  \text{there is no edge }e\in E\text{ such that
the two endpoints of }e\text{ have the same color}\right) \\
&  \ \ \ \ \ \ \ \ \ \ \ \ \ \ \ \ \ \ \ \ \left(
\begin{array}
[c]{c}%
\text{since adjacent vertices are vertices that}\\
\text{are the two endpoints of an edge}%
\end{array}
\right) \\
&  \Longleftrightarrow\ \left(  \text{there exists no element of }E_{f}\right)
\\
&  \ \ \ \ \ \ \ \ \ \ \ \ \ \ \ \ \ \ \ \ \left(
\begin{array}
[c]{c}%
\text{since the elements of }E_{f}\text{ are precisely the edges }e\in E\\
\text{such that the two endpoints of }e\text{ have the same color}\\
\text{(by the definition of }E_{f}\text{)}%
\end{array}
\right) \\
&  \Longleftrightarrow\ \left(  E_{f}=\varnothing\right)  .
\end{align*}
This proves Proposition~\ref{prop.mt2.chrompoly.EcapEf}.
\end{proof}

\begin{lemma}
\label{lem.mt2.chrompoly.kconn} Let $G=\left(  V,E,\varphi\right)  $ be a
multigraph. Let $B$ be a subset of $E$. Let $k\in\mathbb{N}$. Then, the number
of all $k$-colorings $f:V\rightarrow\left\{  1,2,\ldots,k\right\}  $
satisfying $B\subseteq E_{f}$ is $k^{\operatorname{conn}\left(  V,B,\varphi
\mid_{B}\right)  }$.
\end{lemma}

\begin{proof}
[Proof of Lemma~\ref{lem.mt2.chrompoly.kconn}.]If $C$ is a nonempty subset of
$V$, and if $f:V\rightarrow\left\{  1,2,\ldots,k\right\}  $ is any
$k$-coloring of $G$, then we shall say that $f$ is \textbf{constant on }$C$ if
the restriction $f\mid_{C}$ is a constant map (i.e., if the colors $f\left(
c\right)  $ for all $c\in C$ are equal). We shall show the following claim:

\begin{statement}
\textit{Claim 1:} Let $f:V\rightarrow\left\{  1,2,\ldots,k\right\}  $ be any
$k$-coloring of $G$. Then, we have $B\subseteq E_{f}$ if and only if $f$ is
constant on each component of the multigraph $\left(  V,B,\varphi\mid
_{B}\right)  $.
\end{statement}

[\textit{Proof of Claim 1:} This is an \textquotedblleft if and only
if\textquotedblright\ statement; we shall prove its \textquotedblleft%
$\Longrightarrow$\textquotedblright\ and \textquotedblleft$\Longleftarrow
$\textquotedblright\ directions separately:

$\Longrightarrow:$ Assume that $B\subseteq E_{f}$. We must prove that $f$ is
constant on each component of the multigraph $\left(  V,B,\varphi\mid
_{B}\right)  $.

Let $C$ be a component of $\left(  V,B,\varphi\mid_{B}\right)  $. We must
prove that $f$ is constant on $C$. In other words, we must prove that
$f\left(  c\right)  =f\left(  d\right)  $ for any $c,d\in C$.

So let us fix $c,d\in C$. Then, the vertices $c$ and $d$ belong to the same
component of the graph $\left(  V,B,\varphi\mid_{B}\right)  $ (namely, to
$C$). Hence, these vertices $c$ and $d$ are path-connected in this graph. In
other words, the graph $\left(  V,B,\varphi\mid_{B}\right)  $ has a path from
$c$ to $d$. Let%
\[
\mathbf{p}=\left(  v_{0},e_{1},v_{1},e_{2},v_{2},\ldots,e_{s},v_{s}\right)
\]
be this path. Hence, $v_{0}=c$ and $v_{s}=d$ and $e_{1},e_{2},\ldots,e_{s}\in
B$.

Let $i\in\left\{  1,2,\ldots,s\right\}  $. Then, the endpoints of the edge
$e_{i}$ are $v_{i-1}$ and $v_{i}$ (since $e_{i}$ is surrounded by $v_{i-1}$
and $v_{i}$ on the path $\mathbf{p}$). However, from $e_{1},e_{2},\ldots
,e_{s}\in B$, we obtain $e_{i}\in B\subseteq E_{f}$. Hence, the two endpoints
of $e_{i}$ have the same color in $f$ (by the definition of $E_{f}$). In other
words, $f\left(  v_{i-1}\right)  =f\left(  v_{i}\right)  $ (since the
endpoints of the edge $e_{i}$ are $v_{i-1}$ and $v_{i}$).

Forget that we fixed $i$. We thus have proved the equality $f\left(
v_{i-1}\right)  =f\left(  v_{i}\right)  $ for each $i\in\left\{
1,2,\ldots,s\right\}  $. Combining these equalities, we obtain%
\[
f\left(  v_{0}\right)  =f\left(  v_{1}\right)  =f\left(  v_{2}\right)
=\cdots=f\left(  v_{s}\right)  .
\]
Hence, $f\left(  v_{0}\right)  =f\left(  v_{s}\right)  $. In other words,
$f\left(  c\right)  =f\left(  d\right)  $ (since $v_{0}=c$ and $v_{s}=d$).

Forget that we fixed $c$ and $d$. We thus have shown that $f\left(  c\right)
=f\left(  d\right)  $ for any $c,d\in C$. In other words, $f$ is constant on
$C$. Since $C$ was allowed to be an arbitrary component of $\left(
V,B,\varphi\mid_{B}\right)  $, we thus conclude that $f$ is constant on each
component of the multigraph $\left(  V,B,\varphi\mid_{B}\right)  $. This
proves the \textquotedblleft$\Longrightarrow$\textquotedblright\ direction of
Claim 1.

$\Longleftarrow:$ Assume that $f$ is constant on each component of the
multigraph $\left(  V,B,\varphi\mid_{B}\right)  $. We must prove that
$B\subseteq E_{f}$.

Indeed, let $e\in B$. Let $u$ and $v$ be the two endpoints of $e$. Then,
$\left(  u,e,v\right)  $ is a walk from $u$ to $v$ in the multigraph $\left(
V,B,\varphi\mid_{B}\right)  $ (since $e\in B$). Hence, $u$ is path-connected
to $v$ in this multigraph. In other words, $u$ and $v$ belong to the same
component of the multigraph $\left(  V,B,\varphi\mid_{B}\right)  $. Therefore,
$f\left(  u\right)  =f\left(  v\right)  $ (since $f$ is constant on each
component of the multigraph $\left(  V,B,\varphi\mid_{B}\right)  $). This
means that the two endpoints of $e$ have the same color in $f$ (since $u$ and
$v$ are the endpoints of $e$). Combining this with the fact that $e\in E$
(because $e\in B\subseteq E$), we conclude that $e\in E_{f}$ (by the
definition of $E_{f}$).

Forget that we fixed $e$. We thus have shown that $e\in E_{f}$ for each $e\in
B$. In other words, $B\subseteq E_{f}$. This proves the \textquotedblleft%
$\Longleftarrow$\textquotedblright\ direction of Claim 1. The proof of Claim 1
is now complete.] \medskip

Now, Claim 1 shows that the $k$-colorings $f:V\rightarrow\left\{
1,2,\ldots,k\right\}  $ satisfying $B\subseteq E_{f}$ are precisely the
$k$-colorings $f:V\rightarrow\left\{  1,2,\ldots,k\right\}  $ that are
constant on each component of the graph $\left(  V,B,\varphi\mid_{B}\right)
$. Hence, all such $k$-colorings $f$ can be obtained by the following procedure:

\begin{itemize}
\item \textbf{For each} component $C$ of the graph $\left(  V,B,\varphi
\mid_{B}\right)  $, pick a color $c_{C}$ (that is, an element $c_{C}$ of
$\left\{  1,2,\ldots,k\right\}  $) and then assign this color $c_{C}$ to each
vertex in $C$ (that is, set $f\left(  v\right)  =c_{C}$ for each $v\in C$).
\end{itemize}

\medskip This procedure involves choices (because for each component $C$ of
$\left(  V,B,\varphi\mid_{B}\right)  $, we get to pick a color): Namely, for
each of the $\operatorname{conn}\left(  V,B,\varphi\mid_{B}\right)  $ many
components of the graph $\left(  V,B,\varphi\mid_{B}\right)  $, we must choose
a color from the set $\left\{  1,2,\ldots,k\right\}  $. Thus, we have a total
of $k^{\operatorname{conn}\left(  V,B,\varphi\mid_{B}\right)  }$ many options
(since we are choosing among $k$ colors for each of the $\operatorname{conn}%
\left(  V,B,\varphi\mid_{B}\right)  $ components). Each of these options gives
rise to a different $k$-coloring $f:V\rightarrow\left\{  1,2,\ldots,k\right\}
$. Therefore, the number of all $k$-colorings $f:V\rightarrow\left\{
1,2,\ldots,k\right\}  $ satisfying $B\subseteq E_{f}$ is
$k^{\operatorname{conn}\left(  V,B,\varphi\mid_{B}\right)  }$ (because all of
these $k$-colorings can be obtained by this procedure). This proves
Lemma~\ref{lem.mt2.chrompoly.kconn}.
\end{proof}

\begin{corollary}
\label{cor.mt2.chrompoly.kconn2} Let $\left(  V,E,\varphi\right)  $ be a
multigraph. Let $F$ be a subset of $E$. Let $k\in\mathbb{N}$. Then,
\[
k^{\operatorname{conn}\left(  V,F,\varphi\mid_{F}\right)  }=\sum
_{\substack{f:V\rightarrow\left\{  1,2,\ldots,k\right\}  ;\\F\subseteq E_{f}%
}}1.
\]

\end{corollary}

\begin{proof}
[Proof of Corollary~\ref{cor.mt2.chrompoly.kconn2}.]We have
\begin{align*}
\sum_{\substack{f:V\rightarrow\left\{  1,2,\ldots,k\right\}  ;\\F\subseteq
E_{f}}}1  &  =\left(  \text{the number of all }f:V\rightarrow\left\{
1,2,\ldots,k\right\}  \text{ satisfying }F\subseteq E_{f}\right)  \cdot1\\
&  =\left(  \text{the number of all }f:V\rightarrow\left\{  1,2,\ldots
,k\right\}  \text{ satisfying }F\subseteq E_{f}\right) \\
&  =k^{\operatorname{conn}\left(  V,F,\varphi\mid_{F}\right)  }%
\end{align*}
(because Lemma \ref{lem.mt2.chrompoly.kconn} (applied to $B=F$) shows that the
number of all $k$-colorings $f:V\rightarrow\left\{  1,2,\ldots,k\right\}  $
satisfying $F\subseteq E_{f}$ is $k^{\operatorname{conn}\left(  V,F,\varphi
\mid_{F}\right)  }$). This proves Corollary~\ref{cor.mt2.chrompoly.kconn2}.
\end{proof}

\begin{proof}
[Proof of Theorem \ref{thm.coloring.whitney}.]First of all, the equality%
\[
\sum_{F\subseteq E}\left(  -1\right)  ^{\left\vert F\right\vert }%
x^{\operatorname*{conn}\left(  V,F,\varphi\mid_{F}\right)  }=\sum
_{\substack{H\text{ is a spanning}\\\text{subgraph of }G}}\left(  -1\right)
^{\left\vert \operatorname*{E}\left(  H\right)  \right\vert }%
x^{\operatorname*{conn}H}%
\]
is clear, because the spanning subgraphs of $G$ are precisely the subgraphs of
the form $\left(  V,F,\varphi\mid_{F}\right)  $ for some $F\subseteq E$.

Now, let $k\in\mathbb{N}$. We must prove that $\left(  \text{\# of proper
}k\text{-colorings of }G\right)  =\chi_{G}\left(  k\right)  $.

Let us substitute $k$ for $x$ in the equality
\[
\chi_{G}=\sum_{F\subseteq E}\left(  -1\right)  ^{\left\vert F\right\vert
}x^{\operatorname{conn}\left(  V,F,\varphi\mid_{F}\right)  }.
\]
We thus obtain
\begin{align*}
&  \chi_{G}\left(  k\right) \\
&  =\sum_{F\subseteq E}\left(  -1\right)  ^{\left\vert F\right\vert
}\underbrace{k^{\operatorname{conn}\left(  V,F,\varphi\mid_{F}\right)  }%
}_{\substack{=\sum_{\substack{f:V\rightarrow\left\{  1,2,\ldots,k\right\}
;\\F\subseteq E_{f}}}1\\\text{(by Corollary \ref{cor.mt2.chrompoly.kconn2})}%
}}=\sum_{F\subseteq E}\left(  -1\right)  ^{\left\vert F\right\vert }%
\sum_{\substack{f:V\rightarrow\left\{  1,2,\ldots,k\right\}  ;\\F\subseteq
E_{f}}}1\\
&  =\underbrace{\sum_{F\subseteq E}\ \ \sum_{\substack{f:V\rightarrow\left\{
1,2,\ldots,k\right\}  ;\\F\subseteq E_{f}}}}_{=\sum_{f:V\rightarrow\left\{
1,2,\ldots,k\right\}  }\ \ \sum_{\substack{F\subseteq E;\\F\subseteq E_{f}}%
}}\underbrace{\left(  -1\right)  ^{\left\vert F\right\vert }1}_{=\left(
-1\right)  ^{\left\vert F\right\vert }}=\sum_{f:V\rightarrow\left\{
1,2,\ldots,k\right\}  }\ \ \underbrace{\sum_{\substack{F\subseteq
E;\\F\subseteq E_{f}}}}_{\substack{=\sum_{F\subseteq E_{f}}\\\text{(since
}E_{f}\subseteq E\text{)}}}\left(  -1\right)  ^{\left\vert F\right\vert }\\
&  =\sum_{f:V\rightarrow\left\{  1,2,\ldots,k\right\}  }\ \ \sum_{F\subseteq
E_{f}}\left(  -1\right)  ^{\left\vert F\right\vert }=\sum_{f:V\rightarrow
\left\{  1,2,\ldots,k\right\}  }\ \ \underbrace{\sum_{A\subseteq E_{f}}\left(
-1\right)  ^{\left\vert A\right\vert }}_{\substack{=\left[  E_{f}%
=\varnothing\right]  \\\text{(by Lemma \ref{lem.dominating.heinrich-lemma1}%
,}\\\text{applied to }P=E_{f}\text{)}}}\\
&  \ \ \ \ \ \ \ \ \ \ \ \ \ \ \ \ \ \ \ \ \left(
\begin{array}
[c]{c}%
\text{here, we have renamed the summation index }F\\
\text{in the inner sum as }A
\end{array}
\right) \\
&  =\sum_{f:V\rightarrow\left\{  1,2,\ldots,k\right\}  }\left[  E_{f}%
=\varnothing\right] \\
&  =\sum_{\substack{f:V\rightarrow\left\{  1,2,\ldots,k\right\}
;\\E_{f}=\varnothing}}\underbrace{\left[  E_{f}=\varnothing\right]
}_{\substack{=1\\\text{(since }E_{f}=\varnothing\text{ is true)}}%
}+\sum_{\substack{f:V\rightarrow\left\{  1,2,\ldots,k\right\}  ;\\\text{not
}E_{f}=\varnothing}}\underbrace{\left[  E_{f}=\varnothing\right]
}_{\substack{=0\\\text{(since }E_{f}=\varnothing\text{ is false)}}}\\
&  \ \ \ \ \ \ \ \ \ \ \ \ \ \ \ \ \ \ \ \ \left(
\begin{array}
[c]{c}%
\text{since each }f:V\rightarrow\left\{  1,2,\ldots,k\right\}  \text{ either
satisfies }E_{f}=\varnothing\\
\text{or does not}%
\end{array}
\right) \\
&  =\sum_{\substack{f:V\rightarrow\left\{  1,2,\ldots,k\right\}
;\\E_{f}=\varnothing}}1+\underbrace{\sum_{\substack{f:V\rightarrow\left\{
1,2,\ldots,k\right\}  ;\\\text{not }E_{f}=\varnothing}}0}_{\substack{=0}%
}=\sum_{\substack{f:V\rightarrow\left\{  1,2,\ldots,k\right\}  ;\\E_{f}%
=\varnothing}}1\\
&  =\left(  \text{the number of all }f:V\rightarrow\left\{  1,2,\ldots
,k\right\}  \text{ such that }E_{f}=\varnothing\right)  \cdot1\\
&  =\left(  \text{the number of all }f:V\rightarrow\left\{  1,2,\ldots
,k\right\}  \text{ such that }E_{f}=\varnothing\right) \\
&  =\left(  \text{the number of all }f:V\rightarrow\left\{  1,2,\ldots
,k\right\}  \text{ such that the }k\text{-coloring }f\text{ is proper}\right)
\\
&  \ \ \ \ \ \ \ \ \ \ \ \ \ \ \ \ \ \ \ \ \left(
\begin{array}
[c]{c}%
\text{since Proposition \ref{prop.mt2.chrompoly.EcapEf} shows that the
condition \textquotedblleft}E_{f}=\varnothing\text{\textquotedblright}\\
\text{is equivalent to \textquotedblleft the }k\text{-coloring }f\text{ is
proper\textquotedblright}%
\end{array}
\right) \\
&  =\left(  \text{the number of all proper }k\text{-colorings}\right)  .
\end{align*}
In other words, the number of proper $k$-colorings of $G$ is $\chi_{G}\left(
k\right)  $. This completes the proof of Theorem \ref{thm.coloring.whitney}.
\end{proof}

\begin{definition}
The polynomial $\chi_{G}$ in Theorem \ref{thm.coloring.whitney} is known as
the \textbf{chromatic polynomial} of $G$.
\end{definition}

Here are the chromatic polynomials of some graphs:

\begin{proposition}
\label{prop.coloring.chi-easy}Let $n\geq1$ be an integer.

\begin{enumerate}
\item[\textbf{(a)}] For the path graph $P_{n}$ with $n$ vertices, we have%
\[
\chi_{P_{n}}=x\left(  x-1\right)  ^{n-1}.
\]

\item[\textbf{(b)}] More generally, for any tree $T$ with $n$ vertices, we
have%
\[
\chi_{T}=x\left(  x-1\right)  ^{n-1}.
\]

\item[\textbf{(c)}] For the complete graph $K_{n}$ with $n$ vertices, we have%
\[
\chi_{K_{n}}=x\left(  x-1\right)  \left(  x-2\right)  \cdots\left(
x-n+1\right)  .
\]

\item[\textbf{(d)}] For the empty graph $E_{n}$ with $n$ vertices, we have%
\[
\chi_{E_{n}}=x^{n}.
\]

\item[\textbf{(e)}] Assume that $n\geq2$. For the cycle graph $C_{n}$ with $n$
vertices, we have%
\[
\chi_{C_{n}}=\left(  x-1\right)  ^{n}+\left(  -1\right)  ^{n}\left(
x-1\right)  .
\]

\end{enumerate}
\end{proposition}

\begin{fineprint}

\begin{proof}
[Proof sketch.]\textbf{(c)} In order to prove that two polynomials with real
coefficients are identical, it suffices to show that they agree on all
nonnegative integers (this is an instance of the \textquotedblleft principle
of permanence of polynomial identities\textquotedblright\ that we have already
stated as Theorem \ref{thm.polyIDtrick}). Thus, in order to prove that
$\chi_{K_{n}}=x\left(  x-1\right)  \left(  x-2\right)  \cdots\left(
x-n+1\right)  $, it suffices to show that $\chi_{K_{n}}\left(  k\right)
=k\left(  k-1\right)  \left(  k-2\right)  \cdots\left(  k-n+1\right)  $ for
each $k\in\mathbb{N}$.

So let us do this. Fix $k\in\mathbb{N}$. Theorem \ref{thm.coloring.whitney}
(applied to $G=K_{n}$) yields%
\begin{equation}
\left(  \text{\# of proper }k\text{-colorings of }K_{n}\right)  =\chi_{K_{n}%
}\left(  k\right)  . \label{pf.prop.coloring.chi-easy.c.1}%
\end{equation}
Now, how many proper $k$-colorings does $K_{n}$ have? We can construct such a
proper $k$-coloring as follows:

\begin{itemize}
\item First, choose the color of the vertex $1$. There are $k$ options for this.

\item Then, choose the color of the vertex $2$. There are $k-1$ options for
this, since it must differ from the color of $1$.

\item Then, choose the color of the vertex $3$. There are $k-2$ options for
this, since it must differ from the colors of $1$ and of $2$ (and the latter
two colors are distinct, so we must subtract $2$, not $1$).

\item And so on, until all $n$ vertices are colored.
\end{itemize}

The total number of options to perform this construction is therefore
\newline$k\left(  k-1\right)  \left(  k-2\right)  \cdots\left(  k-n+1\right)
$. Hence,
\[
\left(  \text{\# of proper }k\text{-colorings of }K_{n}\right)  =k\left(
k-1\right)  \left(  k-2\right)  \cdots\left(  k-n+1\right)  .
\]
Comparing this with (\ref{pf.prop.coloring.chi-easy.c.1}), we obtain
$\chi_{K_{n}}\left(  k\right)  =k\left(  k-1\right)  \left(  k-2\right)
\cdots\left(  k-n+1\right)  $. As we already explained, this completes the
proof of Proposition \ref{prop.coloring.chi-easy} \textbf{(c)}. \medskip

\textbf{(d)} This is similar to part \textbf{(c)}, but easier. We leave the
proof to the reader. Alternatively, it follows easily from the definition of
$\chi_{E_{n}}$, since $E_{n}$ has only one spanning subgraph (namely, $E_{n}$
itself). \medskip

\textbf{(b)} (This is an outline; see \cite[\S 0.6]{17s-mt2s} for details.)

We proceed by induction on $n$. If $n=1$, then this is easily checked by hand.
If $n>1$, then the tree $T$ has at least one leaf (by Theorem
\ref{thm.tree.leaves.2} \textbf{(a)}). Thus, we can fix a leaf $\ell$ of $T$.
The graph $T\setminus\ell$ then is a tree (by Theorem
\ref{thm.tree.leaf-ind.1}) and has $n-1$ vertices, and therefore (by the
induction hypothesis) its chromatic polynomial is $\chi_{T\setminus\ell
}=x\left(  x-1\right)  ^{n-2}$. However, for any given $k\in\mathbb{N}$, we
can construct a proper $k$-coloring of $T$ by first choosing a proper
$k$-coloring of $T\setminus\ell$ and then choosing the color of the remaining
leaf $\ell$ (there are $k-1$ choices for it, since it has to differ from the
color of the unique neighbor of $\ell$). Therefore, for each $k\in\mathbb{N}$,
we have%
\[
\left(  \text{\# of proper }k\text{-colorings of }T\right)  =\left(  \text{\#
of proper }k\text{-colorings of }T\setminus\ell\right)  \cdot\left(
k-1\right)  .
\]
In view of Theorem \ref{thm.coloring.whitney}, this equality can be rewritten
as
\[
\chi_{T}\left(  k\right)  =\chi_{T\setminus\ell}\left(  k\right)  \cdot\left(
k-1\right)  .
\]
Since this holds for all $k\in\mathbb{N}$, we thus conclude that%
\[
\chi_{T}=\underbrace{\chi_{T\setminus\ell}}_{=x\left(  x-1\right)  ^{n-2}%
}\cdot\left(  x-1\right)  =x\left(  x-1\right)  ^{n-2}\cdot\left(  x-1\right)
=x\left(  x-1\right)  ^{n-1}.
\]
This completes the induction step.

Alternatively, Proposition \ref{prop.coloring.chi-easy} \textbf{(b)} can also
be derived from the definition of $\chi_{T}$, using the fact that every
spanning subgraph $H$ of $T$ has no cycles and therefore satisfies
$\operatorname*{conn}H=n-\left\vert \operatorname*{E}\left(  H\right)
\right\vert $ (by Corollary \ref{cor.conn.eq}). \medskip

\textbf{(a)} This is a particular case of part \textbf{(b)}, since $P_{n}$ is
a tree with $n$ vertices. \medskip

\textbf{(e)} There are different ways to prove this; see \cite{LeeShi19} for
four different proofs. The simplest one is probably by induction on $n$: Let
$n\geq2$. Fix $k\in\mathbb{N}$. A proper $k$-coloring of $C_{n}$ is the same
as a proper $k$-coloring of $P_{n}$ that assigns different colors to the
vertices $1$ and $n$. Hence,%
\begin{align*}
&  \left(  \text{\# of proper }k\text{-colorings of }C_{n}\right) \\
&  =\left(  \text{\# of proper }k\text{-colorings of }P_{n}\text{ that assign
different colors to }1\text{ and }n\right) \\
&  =\underbrace{\left(  \text{\# of proper }k\text{-colorings of }%
P_{n}\right)  }_{\substack{=k\left(  k-1\right)  ^{n-1}\\\text{(by part
\textbf{(a)})}}}\\
&  \ \ \ \ \ \ \ \ \ \ -\underbrace{\left(  \text{\# of proper }%
k\text{-colorings of }P_{n}\text{ that assign the same color to }1\text{ and
}n\right)  }_{\substack{=\left(  \text{\# of proper }k\text{-colorings of
}C_{n-1}\right)  \\\text{(why?)}}}\\
&  =k\left(  k-1\right)  ^{n-1}-\left(  \text{\# of proper }k\text{-colorings
of }C_{n-1}\right)  .
\end{align*}
In view of Theorem \ref{thm.coloring.whitney}, this equality can be rewritten
as
\[
\chi_{C_{n}}\left(  k\right)  =k\left(  k-1\right)  ^{n-1}-\chi_{C_{n-1}%
}\left(  k\right)  .
\]
Since this holds for all $k\in\mathbb{N}$, we thus obtain%
\[
\chi_{C_{n}}=x\left(  x-1\right)  ^{n-1}-\chi_{C_{n-1}}.
\]
This is a recursion that is easily solved for $\chi_{C_{n}}$, yielding the
claim of part \textbf{(e)}.

(Proposition \ref{prop.coloring.chi-easy} \textbf{(e)} also appeared as
Exercise 2 \textbf{(a)} on midterm \#3 in my Spring 2017 course; see
\href{https://www.cip.ifi.lmu.de/~grinberg/t/17s/}{the course website} for solutions.)
\end{proof}
\end{fineprint}

\begin{exercise}
\label{exe.mt3.chromatic-example-b}Let $g\in\mathbb{N}$. Let $G$ be the simple
graph whose vertices are the $2g+1$ integers $-g,-g+1,\ldots,g-1,g$, and whose
edges are
\begin{align*}
&  \left\{  0,i\right\}  \qquad\text{ for all }i\in\left\{  1,2,\ldots
,g\right\}  ;\\
&  \left\{  0,-i\right\}  \qquad\text{ for all }i\in\left\{  1,2,\ldots
,g\right\}  ;\\
&  \left\{  i,-i\right\}  \qquad\text{ for all }i\in\left\{  1,2,\ldots
,g\right\}
\end{align*}
(these are $3g$ edges in total).

Compute the chromatic polynomial $\chi_{G}$ of $G$.

[Here is how $G$ looks like in the case when $g=4$:
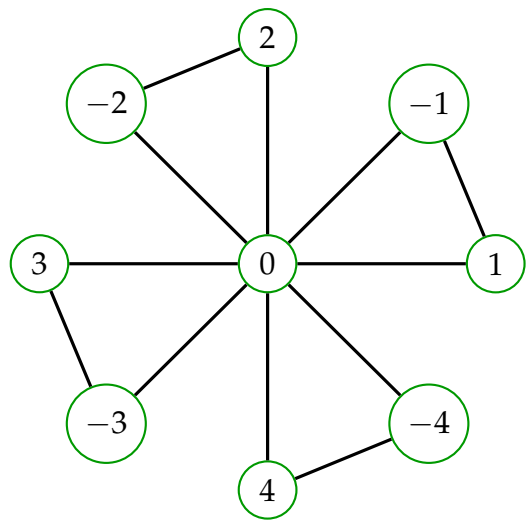
\[%
\begin{tikzpicture}[scale=1.5]
\begin{scope}[every node/.style={circle,thick,draw=green!60!black}]
\node(O) at (0,0) {$0$};
\node(A) at (0:2) {$1$};
\node(B) at (45:2) {$-1$};
\node(C) at (90:2) {$2$};
\node(D) at (135:2) {$-2$};
\node(E) at (180:2) {$3$};
\node(F) at (225:2) {$-3$};
\node(G) at (270:2) {$4$};
\node(H) at (315:2) {$-4$};
\end{scope}
\begin{scope}[every edge/.style={draw=black,very thick}]
\path[-] (A) edge (B) (C) edge (D) (E) edge (F) (G) edge (H);
\path
[-] (O) edge (A) (O) edge (B) (O) edge (C) (O) edge (D) (O) edge (E) (O) edge (F) (O) edge (G) (O) edge (H);
\end{scope}
\end{tikzpicture}%
\]
] \medskip

[\textbf{Solution:} This is Exercise 2 \textbf{(b)} on midterm \#3 from my
Spring 2017 course; see \href{https://www.cip.ifi.lmu.de/~grinberg/t/17s/}{the
course page} for solutions.]
\end{exercise}

\begin{exercise}
Let $G$ be a multigraph. Let $\chi_{G}$ be its chromatic polynomial. For each
$i\in\mathbb{N}$, let $a_{i}$ be the coefficient of $x^{i}$ in $\chi_{G}$.
Prove that $a_{1}$ is odd if and only if the graph $G$ is connected and has a
proper $2$-coloring.
\end{exercise}

\subsection{Vizing's theorem}

So far we have been coloring the vertices of a graph. We can also color the edges:

\begin{definition}
Let $G=\left(  V,E,\varphi\right)  $ be a multigraph. Let $k\in\mathbb{N}$.

A $k$\textbf{-edge-coloring} of $G$ means a map $f:E\rightarrow\left\{
1,2,\ldots,k\right\}  $.

Such a $k$-edge-coloring $f$ is called \textbf{proper} if no two distinct
edges that have a common endpoint have the same color.
\end{definition}

The most prominent fact about edge-colorings is the following theorem:

\begin{theorem}
[Vizing's theorem]Let $G$ be a simple graph with at least one vertex. Let%
\[
\alpha:=\max\left\{  \deg v\ \mid\ v\in V\right\}  .
\]
Then, $G$ has a proper $\left(  \alpha+1\right)  $-edge-coloring.
\end{theorem}

\begin{proof}
See, e.g., \cite{Schrij04} or various textbooks on graph theory.\footnote{Note
that \cite{Schrij04} uses some standard graph-theoretical notations: What we
call $\alpha$ is denoted by $\Delta\left(  G\right)  $ in \cite{Schrij04},
whereas $\chi^{\prime}\left(  G\right)  $ denotes the minimum $k\in\mathbb{N}$
for which $G$ has a proper $k$-edge-coloring.}
\end{proof}

Two remarks:

\begin{itemize}
\item The $\alpha+1$ in Vizing's theorem cannot be improved in general (e.g.,
take $G$ to be an odd-length cycle graph $C_{n}$).

\item Vizing's theorem can be adapted to work for multigraphs instead of
simple graphs. However, this requires replacing the $\alpha+1$ by $\alpha+m$,
where $m$ is the maximum number of distinct mutually parallel edges in $G$
(since otherwise, the multigraph $\left(  K_{3}^{\operatorname*{bidir}%
}\right)  ^{\operatorname*{und}}$ would be a counterexample, as it has
$\alpha=4$ but has no proper $5$-edge-coloring). For a proof of this, see
\cite[Corollary 2]{BerFou91}.
\end{itemize}

\subsection{Further exercises}

Some interesting things can be said about colorings of graphs, even about
non-proper colorings:

\begin{exercise}
\label{exe.hw0.enmity}Let $G=\left(  V,E\right)  $ be a simple graph.

Prove that there exists a $2$-coloring $f$ of $G$ with the following property:
For each vertex $v\in V$, at most $\dfrac{1}{2}\deg v$ among the neighbors of
$v$ have the same color as $v$. \medskip

[\textbf{Remark:} This problem is often restated as follows: You are given a
(finite) set of politicians; some politicians are mutual enemies. (No
politician is his own enemy. If $u$ is an enemy of $v$, then $v$ is an enemy
of $u$. An enemy of an enemy is not necessarily a friend. So this is just a
simple graph.) Prove that it is possible to subdivide this set into two
(disjoint) parties such that no politician has more than half of his enemies
in his own party.] \medskip

[\textbf{Hint:} First, pick an arbitrary $2$-coloring $f$ of $G$. Then,
gradually improve it until it satisfies the required property.] \medskip

[\textbf{Solution:} This is Exercise 1 on homework set \#0 from my Spring 2017
course; see \href{https://www.cip.ifi.lmu.de/~grinberg/t/17s/}{the course
page} for solutions.]
\end{exercise}

Exercise \ref{exe.hw0.enmity} can be generalized to multiple colors:

\begin{exercise}
\label{exe.mt1.enmity} Let $k\in\mathbb{N}$. Let $p_{1},p_{2},\ldots,p_{k}$ be
$k$ nonnegative real numbers such that $p_{1}+p_{2}+\cdots+p_{k}\geq1$.

Let $G=\left(  V,E\right)  $ be a simple graph.

Prove that there exists a $k$-coloring $f$ of $G$ with the following property:
For each vertex $v\in V$, at most $p_{f\left(  v\right)  }\deg v$ neighbors of
$v$ have the same color as $v$. \medskip

[\textbf{Solution:} This is Exercise 5 on midterm \#1 from my Spring 2017
course; see \href{https://www.cip.ifi.lmu.de/~grinberg/t/17s/}{the course
page} for solutions.]
\end{exercise}

\subsection{Some recent results}

The theory of proper $k$-colorings can be surprisingly deep and difficult. The
following theorem was conjectured by Ioan Tomescu in 1971, and only proved in
2019 by Fox, He and Manners using sophisticated probabilistic arguments
\cite[Theorem 1]{FoHeMa19}:

\begin{theorem}
[Tomescu, Fox, He, Manners]Let $G$ be a graph with $n$ vertices, and let $k$
be a positive integer. Assume that $G$ has no proper $\left(  k-1\right)
$-coloring. Then:

\begin{enumerate}
\item[\textbf{(a)}] The number of proper $k$-colorings of $G$ is at most
$k!\cdot k^{n-k}$.

\item[\textbf{(b)}] If $G$ is connected and $k\geq4$, then the number of
proper $k$-colorings of $G$ is at most $k!\cdot\left(  k-1\right)  ^{n-k}$.
\end{enumerate}
\end{theorem}

\noindent(Part \textbf{(a)} is \cite[problem 10.22]{Tomesc85}; see
\url{https://mathoverflow.net/a/492188/} for a clearer version of the proof.)
An even more general conjecture of Tomescu says that if $G$ is a connected
graph with $n$ vertices, and if $s$ and $k$ are two integers such that $s\geq
k\geq4$ and such that $G$ has no proper $\left(  k-1\right)  $-coloring, then
the number of proper $s$-colorings of $G$ is at most $s\left(  s-1\right)
\left(  s-2\right)  \cdots\left(  s-k+1\right)  \cdot\left(  s-1\right)
^{n-k}$; this is still unsolved (\cite[Conjecture 12]{FoHeMa19}). \medskip

The coefficients of the chromatic polynomial $\chi_{G}$ defined in Theorem
\ref{thm.coloring.whitney} have also been the topic of much research, from
which I shall only mention a few highlights:

\begin{theorem}
[Whitney, Huh]\label{thm.coloring.chrom-ineqs}Let $G$ be a graph with $n$
vertices. Write its chromatic polynomial $\chi_{G}$ as $\sum_{i\in\mathbb{N}%
}a_{i}x^{i}$ with $a_{i}\in\mathbb{Z}$ (so that $a_{0},a_{1},a_{2},\ldots$ are
its coefficients). Then:

\begin{enumerate}
\item[\textbf{(a)}] We have $a_{i}=0$ for all $i>n$.

\item[\textbf{(b)}] If $n>0$, then $a_{0}=0$.

\item[\textbf{(c)}] For each $i\in\mathbb{N}$, we have $\left(  -1\right)
^{n-i}a_{i}\geq0$. (Thus, the coefficients $a_{0},a_{1},a_{2},\ldots$
alternate in sign; in other words, all coefficients of the polynomial $\left(
-1\right)  ^{n}\chi_{G}\left(  -x\right)  $ are nonnegative.)

\item[\textbf{(d)}] For each $i>0$, we have $a_{i}^{2}\geq a_{i-1}a_{i+1}$.
\end{enumerate}
\end{theorem}

Parts \textbf{(a)} and \textbf{(b)} of Theorem \ref{thm.coloring.chrom-ineqs}
are very easy consequences of Theorem \ref{thm.coloring.whitney}, and are only
mentioned for completeness's sake. Part \textbf{(c)} was proved by Whitney
\cite[\S 7]{Whitney32} combinatorially (i.e., by constructing a set with size
$\left(  -1\right)  ^{n-i}a_{i}$). Part \textbf{(d)} is a famous conjecture,
stated by Hoggar in 1974, and only proved in 2012 by June Huh \cite{Huh12}
using deep algebraic geometry. A later proof by Br\"{a}nden and Leake
\cite[Theorem 4.1]{BraLea21} relies on multivariate analysis, which is still
far from elementary. Huh's work has brought the subject to the fore of
research again, and many similar results have since been proven using his
methods. More about the chromatic polynomial can be found in the recent text
\cite{Tittma25}.

Another famous result on proper colorings is
\href{https://en.wikipedia.org/wiki/Four_color_theorem}{the \textbf{four-color
theorem}}, which says that every planar graph (i.e., any graph that can be
embedded in the real plane in such a way that no two edges intersect; cf.
Example \ref{exa.intro.draw}) has a proper $4$-coloring. This was conjectured
by Guthrie in 1852 and only proved by Appel and Haken in 1976 using a long and
complex computer-aided case analysis; see \cite{FriFri98} for the details and
the history of this result.

\section{\label{chp.indsets}Independent sets}

\subsection{\label{sec.indsets.carowei}Definition and the Caro--Wei theorem}

Next, we define one of the most fundamental notions in graph theory:

\begin{definition}
\label{def.indep.indep}An \textbf{independent set} of a multigraph $G$ means a
subset $S$ of $\operatorname*{V}\left(  G\right)  $ such that no two elements
of $S$ are adjacent.
\end{definition}

In other words, an independent set of $G$ means an induced subgraph of $G$
that has no edges\footnote{This is a somewhat sloppy statement. Of course, an
independent set is not literally an induced subgraph, since the former is just
a set, while the latter is a graph. What I mean is that a subset $S$ of
$\operatorname*{V}\left(  G\right)  $ is independent if and only if the
induced subgraph $G\left[  S\right]  $ has no edges.}. Note that
\textquotedblleft no two elements of $S$\textquotedblright\ doesn't mean
\textquotedblleft no two distinct elements of $S$\textquotedblright.

Thus, for example, what we called an \textquotedblleft
anti-triangle\textquotedblright\ (back in Definition \ref{def.sg.triangle}) is
an independent set of size $3$.

\begin{example}
Let $G$ be the graph $G$ from Example \ref{exa.sg.K4-24}. Then, the
independent sets of $G$ are the sets $\varnothing$, $\left\{  1\right\}  $,
$\left\{  2\right\}  $, $\left\{  3\right\}  $, $\left\{  4\right\}  $ and
$\left\{  2,4\right\}  $.
\end{example}

\begin{remark}
Independent sets are closely related to proper colorings. Indeed, let $G$ be a
graph, and let $k\in\mathbb{N}$. Let $f:V\rightarrow\left\{  1,2,\ldots
,k\right\}  $ be a $k$-coloring. For each $i\in\left\{  1,2,\ldots,k\right\}
$, let
\begin{align*}
V_{i}:=  &  \left\{  v\in V\ \mid\ f\left(  v\right)  =i\right\} \\
=  &  \left\{  \text{all vertices of }G\text{ that have color }i\right\}  .
\end{align*}
Then, the $k$-coloring $f$ is proper if and only if the $k$ sets $V_{1}%
,V_{2},\ldots,V_{k}$ are independent sets of $G$. (Proving this is a matter of
unraveling the definitions of \textquotedblleft independent
sets\textquotedblright\ and \textquotedblleft proper $k$%
-colorings\textquotedblright.)
\end{remark}

One classical computational problem in graph theory is to find a maximum-size
independent set of a given graph. This problem is NP-hard\footnote{And, in
fact, NP-complete; see \cite[Proposition 4.10]{Goldre10}.}, so don't expect a
quick algorithm or even a good formula for the maximum size of an independent
set. However, there are some lower bounds for this maximum size. Here is one,
known as the \textbf{Caro--Wei theorem} (\cite[Chapter 6, Probabilistic
Lens]{AloSpe16}):

\begin{theorem}
[Caro--Wei theorem]\label{thm.indep.lower-bd}Let $G=\left(  V,E,\varphi
\right)  $ be a loopless multigraph. Then, $G$ has an independent set of size
\[
\geq\sum_{v\in V}\dfrac{1}{1+\deg v}.
\]

\end{theorem}

\begin{example}
Let $G$ be the following loopless multigraph:%
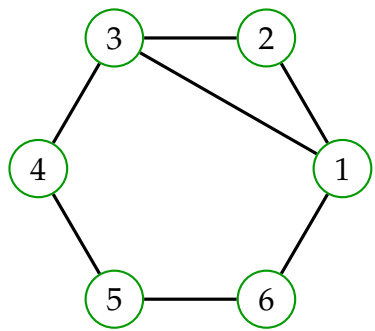
\[%
\begin{tikzpicture}
\begin{scope}[every node/.style={circle,thick,draw=green!60!black}]
\node(A) at (0:2) {$1$};
\node(B) at (60:2) {$2$};
\node(C) at (120:2) {$3$};
\node(D) at (180:2) {$4$};
\node(E) at (240:2) {$5$};
\node(F) at (300:2) {$6$};
\end{scope}
\begin{scope}[every edge/.style={draw=black,very thick}]
\path
[-] (A) edge (B) (B) edge (C) (C) edge (A) (C) edge (D) (D) edge (E) (E) edge (F) (F) edge (A);
\end{scope}
\end{tikzpicture}%
\ \ .
\]
Then, the degrees of the vertices of $G$ are $3,2,3,2,2,2$. Hence, Theorem
\ref{thm.indep.lower-bd} yields that $G$ has an independent set of size%
\[
\geq\dfrac{1}{1+3}+\dfrac{1}{1+2}+\dfrac{1}{1+3}+\dfrac{1}{1+2}+\dfrac{1}%
{1+2}+\dfrac{1}{1+2}=\dfrac{11}{6}\approx1.83.
\]
Since the size of an independent set is always an integer, we can round this
up and conclude that $G$ has an independent set of size $\geq2$. In truth, $G$
actually has an independent set of size $3$ (namely, $\left\{  2,4,6\right\}
$), but there is no way to tell this from the degrees of its vertices alone.
For example, the vertices of the graph%
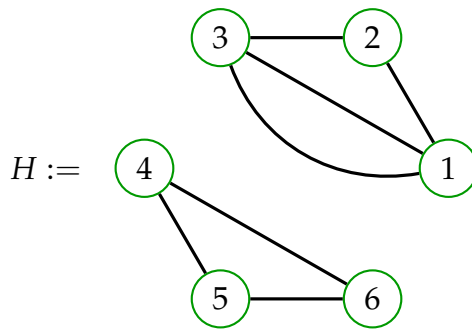
\[%
\begin{tikzpicture}
\begin{scope}[every node/.style={circle,thick,draw=green!60!black}]
\node(A) at (0:2) {$1$};
\node(B) at (60:2) {$2$};
\node(C) at (120:2) {$3$};
\node(D) at (180:2) {$4$};
\node(E) at (240:2) {$5$};
\node(F) at (300:2) {$6$};
\end{scope}
\node(H) at (-3.3, 0) {$H:=$};
\begin{scope}[every edge/.style={draw=black,very thick}]
\path
[-] (A) edge (B) (B) edge (C) (C) edge (A) (C) edge[bend right=40] (A) (D) edge (E) (E) edge (F) (F) edge (D);
\end{scope}
\end{tikzpicture}%
\]
have the same degrees as those of $G$, but $H$ has no independent set of size
$3$.
\end{example}

We shall give two proofs of Theorem \ref{thm.indep.lower-bd}, both of them
illustrating useful techniques.\footnote{Note that the looplessness
requirement in Theorem \ref{thm.indep.lower-bd} is important: If $G$ has a
loop at each vertex, then the only independent set of $G$ is $\varnothing$.}

\begin{proof}
[First proof of Theorem \ref{thm.indep.lower-bd}.]Assume the contrary. Thus,
each independent set $S$ of $G$ has size%
\begin{equation}
\left\vert S\right\vert <\sum_{v\in V}\dfrac{1}{1+\deg v}.
\label{pf.thm.indep.lower-bd.1}%
\end{equation}

A $V$\textbf{-listing} shall mean a list of all vertices in $V$, with each
vertex occurring exactly once in the list. If $\sigma$ is a $V$-listing, then
we define a subset $J_{\sigma}$ of $V$ as follows:%
\[
J_{\sigma}:=\left\{  v\in V\ \mid\ v\text{ occurs \textbf{before} all
neighbors of }v\text{ in }\sigma\right\}  .
\]

[\textbf{Example:} Let $G$ be the following graph:%
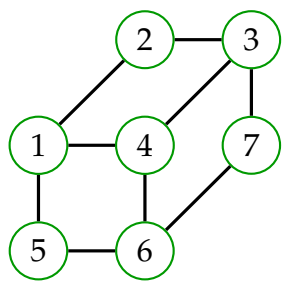
\[%
\begin{tikzpicture}[scale=1.4]
\begin{scope}[every node/.style={circle,thick,draw=green!60!black}]
\node(5) at (0, 0) {$5$};
\node(1) at (0, 1) {$1$};
\node(6) at (1, 0) {$6$};
\node(4) at (1, 1) {$4$};
\node(2) at (1, 2) {$2$};
\node(3) at (2, 2) {$3$};
\node(7) at (2, 1) {$7$};
\end{scope}
\begin{scope}[every edge/.style={draw=black,very thick}]
\path[-] (5) edge (1) edge (6);
\path[-] (4) edge (1) edge (6) edge (3);
\path[-] (3) edge (2) edge (7);
\path[-] (1) edge (2) (6) edge (7);
\end{scope}
\end{tikzpicture}%
\ \ .
\]
Let $\sigma$ be the $V$-listing $\left(  1,2,7,5,3,4,6\right)  $. Then, the
vertex $1$ occurs before all its neighbors ($2$, $4$ and $5$) in $\sigma$, and
thus we have $1\in J_{\sigma}$. Likewise, the vertex $7$ occurs before all its
neighbors ($3$ and $6$) in $\sigma$, so that we have $7\in J_{\sigma}$. But
the vertex $2$ does not occur before all its neighbors in $\sigma$ (indeed, it
occurs after its neighbor $1$), so that we have $2\notin J_{\sigma}$.
Likewise, the vertices $5,3,4,6$ don't belong to $J_{\sigma}$. Altogether, we
thus obtain $J_{\sigma}=\left\{  1,7\right\}  $.] \medskip

The set $J_{\sigma}$ is an independent set of $G$ (because if two vertices $u$
and $v$ in $J_{\sigma}$ were adjacent, then $u$ would have to occur before $v$
in $\sigma$, but $v$ would have to occur before $u$ in $\sigma$; but these two
statements clearly contradict each other). Thus,
(\ref{pf.thm.indep.lower-bd.1}) (applied to $S=J_{\sigma}$) yields%
\[
\left\vert J_{\sigma}\right\vert <\sum_{v\in V}\dfrac{1}{1+\deg v}.
\]

This inequality holds for \textbf{each} $V$-listing $\sigma$. Thus, summing
this inequality over all $V$-listings $\sigma$, we obtain%
\begin{align}
\sum_{\sigma\text{ is a }V\text{-listing}}\left\vert J_{\sigma}\right\vert  &
<\sum_{\sigma\text{ is a }V\text{-listing}}\ \ \sum_{v\in V}\dfrac{1}{1+\deg
v}\nonumber\\
&  =\left(  \text{\# of all }V\text{-listings}\right)  \cdot\sum_{v\in
V}\dfrac{1}{1+\deg v}. \label{pf.thm.indep.lower-bd.3}%
\end{align}

On the other hand, I claim the following:

\begin{statement}
\textit{Claim 1:} For each $v\in V$, we have%
\[
\left(  \text{\# of all }V\text{-listings }\sigma\text{ satisfying }v\in
J_{\sigma}\right)  \geq\dfrac{\left(  \text{\# of all }V\text{-listings}%
\right)  }{1+\deg v}.
\]

\end{statement}

[\textit{Proof of Claim 1:} Fix a vertex $v\in V$. Define $\deg^{\prime}v$ to
be the \# of all neighbors of $v$. Clearly, $\deg^{\prime}v\leq\deg v$.

We shall call a $V$-listing $\sigma$ \textbf{good} if the vertex $v$ occurs in
it before all its neighbors. In other words, a $V$-listing $\sigma$ is good if
and only if it satisfies $v\in J_{\sigma}$ (because $v\in J_{\sigma}$ means
that the vertex $v$ occurs in $\sigma$ before all its neighbors\footnote{This
follows straight from the definition of $J_{\sigma}$.}). Thus, we must show
that
\[
\left(  \text{\# of all good }V\text{-listings}\right)  \geq\dfrac{\left(
\text{\# of all }V\text{-listings}\right)  }{1+\deg v}.
\]

We define a map%
\[
\Gamma:\left\{  \text{all }V\text{-listings}\right\}  \rightarrow\left\{
\text{all good }V\text{-listings}\right\}
\]
as follows: Whenever $\tau$ is a $V$-listing, we let $\Gamma\left(
\tau\right)  $ be the $V$-listing obtained from $\tau$ by swapping $v$ with
the first neighbor of $v$ that occurs in $\tau$ (or, if $\tau$ is already
good, then we just do nothing, i.e., we set $\Gamma\left(  \tau\right)  =\tau
$). This map $\Gamma$ is a $\left(  1+\deg^{\prime}v\right)  $-to-$1$
correspondence -- i.e., for each good $V$-listing $\sigma$, there are exactly
$1+\deg^{\prime}v$ many $V$-listings $\tau$ that satisfy $\Gamma\left(
\tau\right)  =\sigma$ (in fact, one of these $\tau$'s is $\sigma$ itself, and
the remaining $\deg^{\prime}v$ many of these $\tau$'s are obtained from
$\sigma$ by switching $v$ with some neighbor of $v$). Hence, by the
multijection principle\footnote{See a footnote in the proof of Theorem
\ref{thm.BEST.to} for the statement of the multijection principle.}, we
conclude that
\[
\left\vert \left\{  \text{all }V\text{-listings}\right\}  \right\vert =\left(
1+\deg^{\prime}v\right)  \cdot\left\vert \left\{  \text{all good
}V\text{-listings}\right\}  \right\vert .
\]
In other words,%
\[
\left(  \text{\# of all }V\text{-listings}\right)  =\left(  1+\deg^{\prime
}v\right)  \cdot\left(  \text{\# of all good }V\text{-listings}\right)  .
\]
Hence,%
\[
\left(  \text{\# of all good }V\text{-listings}\right)  =\dfrac{\left(
\text{\# of all }V\text{-listings}\right)  }{1+\deg^{\prime}v}\geq
\dfrac{\left(  \text{\# of all }V\text{-listings}\right)  }{1+\deg v}%
\]
(since $\deg^{\prime}v\leq\deg v$). This proves Claim 1 (since the good
$V$-listings are precisely the $V$-listings $\sigma$ satisfying $v\in
J_{\sigma}$).] \medskip

Next, we recall a basic property of the Iverson bracket notation\footnote{See,
e.g., Definition \ref{def.iverson} for the definition of the Iverson bracket
notation.}: If $T$ is a subset of a finite set $S$, then
\begin{equation}
\left\vert T\right\vert =\sum_{v\in S}\left[  v\in T\right]  .
\label{pf.thm.indep.lower-bd.4}%
\end{equation}
(Indeed, the sum $\sum_{v\in S}\left[  v\in T\right]  $ contains an addend
equal to $1$ for each $v\in T$, and an addend equal to $0$ for each $v\in
S\setminus T$. Thus, this sum amounts to $\left\vert T\right\vert
\cdot1+\left\vert S\setminus T\right\vert \cdot0=\left\vert T\right\vert $.)

Now, (\ref{pf.thm.indep.lower-bd.3}) yields%
\begin{align*}
&  \left(  \text{\# of all }V\text{-listings}\right)  \cdot\sum_{v\in V}%
\dfrac{1}{1+\deg v}\\
&  >\sum_{\sigma\text{ is a }V\text{-listing}}\underbrace{\left\vert
J_{\sigma}\right\vert }_{\substack{=\sum\limits_{v\in V}\left[  v\in
J_{\sigma}\right]  \\\text{(by (\ref{pf.thm.indep.lower-bd.4}))}%
}}=\underbrace{\sum_{\sigma\text{ is a }V\text{-listing}}\ \ \sum\limits_{v\in
V}}_{=\sum\limits_{v\in V}\ \ \sum_{\sigma\text{ is a }V\text{-listing}}%
}\left[  v\in J_{\sigma}\right] \\
&  =\sum\limits_{v\in V}\ \ \underbrace{\sum_{\sigma\text{ is a }%
V\text{-listing}}\left[  v\in J_{\sigma}\right]  }%
_{\substack{_{\substack{=\left(  \text{\# of all }V\text{-listings }%
\sigma\text{ satisfying }v\in J_{\sigma}\right)  }}\\\text{(because the sum
}\sum_{\sigma\text{ is a }V\text{-listing}}\left[  v\in J_{\sigma}\right]
\\\text{contains an addend equal to }1\text{ for each }V\text{-listing }%
\sigma\text{ satisfying }v\in J_{\sigma}\text{,}\\\text{and an addend equal to
}0\text{ for each other }V\text{-listing }\sigma\text{)}}}\\
&  =\sum\limits_{v\in V}\underbrace{\left(  \text{\# of all }V\text{-listings
}\sigma\text{ satisfying }v\in J_{\sigma}\right)  }_{\substack{\geq
\dfrac{\left(  \text{\# of all }V\text{-listings}\right)  }{1+\deg
v}\\\text{(by Claim 1)}}}\\
&  \geq\sum\limits_{v\in V}\dfrac{\left(  \text{\# of all }V\text{-listings}%
\right)  }{1+\deg v}=\left(  \text{\# of all }V\text{-listings}\right)
\cdot\sum_{v\in V}\dfrac{1}{1+\deg v}.
\end{align*}
This is absurd (since no real number $x$ can satisfy $x>x$). So we got a
contradiction, and our proof of Theorem \ref{thm.indep.lower-bd} is complete.
\end{proof}

\begin{remark}
This proof is an example of a \textbf{probabilistic proof}. Why? We have been
manipulating sums, but we could easily replace these sums by averages. Claim 1
then would say the following: For any given vertex $v\in V$, the
\textbf{probability} that a (uniformly random) $V$-listing $\sigma$ satisfies
$v\in J_{\sigma}$ is $\geq\dfrac{1}{1+\deg v}$. Thus, the expectation of
$\left\vert J_{\sigma}\right\vert $ is $\geq\sum\limits_{v\in V}\dfrac
{1}{1+\deg v}$ (by linearity of expectation). Therefore, at least one
$V$-listing $\sigma$ actually satisfies $\left\vert J_{\sigma}\right\vert
\geq\sum\limits_{v\in V}\dfrac{1}{1+\deg v}$. So the whole proof can be
restated in terms of probabilities and expectations.

Note that this proof (as it stands) is fairly useless as it comes to actually
\textbf{finding} an independent set of size $\geq\sum\limits_{v\in V}\dfrac
{1}{1+\deg v}$. It does not give any better algorithm than \textquotedblleft
try the subsets $J_{\sigma}$ for all possible $V$-listings $\sigma$; one of
them will work\textquotedblright, which is even slower than trying all subsets
of $V$.

Note also that the proof does \textbf{not} entail that at least half of the
$V$-listings $\sigma$ will satisfy $\left\vert J_{\sigma}\right\vert \geq
\sum\limits_{v\in V}\dfrac{1}{1+\deg v}$. The mean is not the median!
\end{remark}

Let us now give a second proof of the theorem, which does provide a good algorithm:

\begin{proof}
[Second proof of Theorem \ref{thm.indep.lower-bd}.]We proceed by strong
induction on $\left\vert V\right\vert $. Thus, we fix $p\in\mathbb{N}$, and we
assume (as the induction hypothesis) that Theorem \ref{thm.indep.lower-bd} is
already proved for all loopless multigraphs $G$ with $<p$ vertices. We must
now prove it for a loopless multigraph $G=\left(  V,E,\varphi\right)  $ with
$p$ vertices.

If $\left\vert V\right\vert =0$, then this is clear (since $\varnothing$ is an
independent set of appropriate size). Thus, we WLOG assume that $\left\vert
V\right\vert \neq0$. We furthermore assume WLOG that $G$ is a simple graph
(because otherwise, we can replace $G$ by $G^{\operatorname*{simp}}$; this can
only decrease the degrees $\deg v$ of the vertices $v\in V$, and thus our
claim only becomes stronger).

Since $\left\vert V\right\vert \neq0$, there exists a vertex $u\in V$ with
$\deg_{G}u$ minimum\footnote{Here, the notation $\deg_{H}u$\ means the degree
of a vertex $u$ in a graph $H$.}. Pick such a $u$. Thus,%
\begin{equation}
\deg_{G}v\geq\deg_{G}u\ \ \ \ \ \ \ \ \ \ \text{for each }v\in V.
\label{pf.thm.indep.lower-bd.2nd.mindeg}%
\end{equation}

Let $U:=\left\{  u\right\}  \cup\left\{  \text{all neighbors of }u\right\}  $.
Thus, $U\subseteq V$ and $\left\vert U\right\vert =1+\deg_{G}u$ (this is a
honest equality, since $G$ is a simple graph).

Let $G^{\prime}$ be the induced subgraph of $G$ on the set $V\setminus U$.
This is the simple graph obtained from $G$ by removing all vertices belonging
to $U$ (that is, removing the vertex $u$ along with all its neighbors) and
removing all edges that require these vertices. Then, $G^{\prime}$ has fewer
vertices than $G$. Hence, $G^{\prime}$ has $<p$ vertices (since $G$ has $p$
vertices). Hence, by the induction hypothesis, Theorem
\ref{thm.indep.lower-bd} is already proved for $G^{\prime}$. In other words,
$G^{\prime}$ has an independent set of size
\[
\geq\sum_{v\in V\setminus U}\dfrac{1}{1+\deg_{G^{\prime}}v}.
\]
Let $T$ be such an independent set. Set $S:=\left\{  u\right\}  \cup T$. Then,
$S$ is an independent set of $G$ (since $T\subseteq V\setminus U$, so that $T$
contains no neighbors of $u$). Moreover, I claim that $\left\vert S\right\vert
\geq\sum\limits_{v\in V}\dfrac{1}{1+\deg_{G}v}$. Indeed, this follows from%
\begin{align*}
\sum\limits_{v\in V}\dfrac{1}{1+\deg_{G}v}  &  =\sum\limits_{v\in
U}\underbrace{\dfrac{1}{1+\deg_{G}v}}_{\substack{\leq\dfrac{1}{1+\deg_{G}%
u}\\\text{(since }\deg_{G}v\geq\deg_{G}u\\\text{(by
(\ref{pf.thm.indep.lower-bd.2nd.mindeg})))}}}+\sum\limits_{v\in V\setminus
U}\underbrace{\dfrac{1}{1+\deg_{G}v}}_{\substack{\leq\dfrac{1}{1+\deg
_{G^{\prime}}v}\\\text{(since }\deg_{G}v\geq\deg_{G^{\prime}}v\\\text{(because
}G^{\prime}\text{ is a subgraph of }G\text{))}}}\\
&  \leq\underbrace{\sum\limits_{v\in U}\dfrac{1}{1+\deg_{G}u}}%
_{\substack{=\left\vert U\right\vert \cdot\dfrac{1}{1+\deg_{G}u}%
\\=1\\\text{(since }\left\vert U\right\vert =1+\deg_{G}u\text{)}%
}}+\underbrace{\sum\limits_{v\in V\setminus U}\dfrac{1}{1+\deg_{G^{\prime}}v}%
}_{\substack{\leq\left\vert T\right\vert \\\text{(since }T\text{ has size
}\geq\sum_{v\in V\setminus U}\dfrac{1}{1+\deg_{G^{\prime}}v}\text{)}}}\\
&  \leq1+\left\vert T\right\vert =\left\vert S\right\vert
\ \ \ \ \ \ \ \ \ \ \left(  \text{since }S=\left\{  u\right\}  \cup T\right)
.
\end{align*}
So we have found an independent set of $G$ having size $\geq\sum\limits_{v\in
V}\dfrac{1}{1+\deg_{G}v}$ (namely, $S$). This means that Theorem
\ref{thm.indep.lower-bd} holds for our $G$. This completes the induction step,
and Theorem \ref{thm.indep.lower-bd} is proved.
\end{proof}

\begin{remark}
The second proof of Theorem \ref{thm.indep.lower-bd} (unlike the first one)
does give a fairly efficient algorithm for finding an independent set of the
appropriate size. However, the second proof is actually not that much
different from the first proof; it can in fact be recovered from the first
proof by \textbf{derandomization}, specifically using the
\textbf{\href{https://en.wikipedia.org/wiki/Method_of_conditional_probabilities}{\textbf{method
of conditional probabilities}}}. (This is a general technique for
\textquotedblleft derandomizing\textquotedblright\ probabilistic proofs, i.e.,
turning them into algorithmic ones. It often requires some ingenuity and is
not guaranteed to always work, but the above is an example where it can be
applied. See \cite[Chapter 13]{Aspnes23} for much more about derandomization.)
\end{remark}

See also \cite{Chen14} and \cite{AloSpe16} for more about probabilistic proofs
in combinatorics and in general. Here are two more applications of
probabilistic proofs:

\begin{exercise}
\label{exe.8.3}Let $G=\left(  V,E\right)  $ be a simple graph such that each
vertex of $G$ has degree $\geq1$. Prove that there exists a subset $S$ of $V$
having size $\geq\sum_{v\in V}\dfrac{2}{1+\deg v}$ and with the property that
the induced subgraph $G\left[  S\right]  $ is a forest. \medskip

[\textbf{Hint:} As the example of $\begin{tikzpicture}[scale=1.4]
\begin{scope}[every node/.style={circle,thick,draw=green!60!black}]
\node(1) at (0, 0) {$1$};
\node(2) at (1, 0) {$2$};
\node(3) at (2, 0) {$3$};
\end{scope}
\begin{scope}[every edge/.style={draw=black,very thick}, every loop/.style={}]
\path[-] (1) edge[bend left=20] (2) edge[bend right=20] (2) (2) edge (3);
\end{scope}
\end{tikzpicture}
$ shows, this claim is not true for loopless multigraphs (unlike the similar
Theorem \ref{thm.indep.lower-bd}).]
\end{exercise}

\begin{exercise}
\label{exe.8.4}Let $n$ be a positive integer. Prove that there exists a
tournament with $n$ vertices and at least $\dfrac{n!}{2^{n-1}}$ Hamiltonian paths.
\end{exercise}

\subsection{A weaker (but simpler) lower bound}

Let us now weaken Theorem \ref{thm.indep.lower-bd} a bit:

\begin{corollary}
\label{cor.indep.anti-turan}Let $G$ be a loopless multigraph with $n$ vertices
and $m$ edges. Then, $G$ has an independent set of size%
\[
\geq\dfrac{n^{2}}{n+2m}.
\]

\end{corollary}

In order to prove this, we will need the following inequality:

\begin{lemma}
\label{lem.indep.cauchy-lame}Let $a_{1},a_{2},\ldots,a_{n}$ be $n$ positive
reals. Then,%
\[
\dfrac{1}{a_{1}}+\dfrac{1}{a_{2}}+\cdots+\dfrac{1}{a_{n}}\geq\dfrac{n^{2}%
}{a_{1}+a_{2}+\cdots+a_{n}}.
\]

\end{lemma}

\begin{proof}
[Proof of Lemma \ref{lem.indep.cauchy-lame}.]There are several ways to prove
this:\footnote{For unexplained terminology used in the bullet points below,
see any textbook on inequalities, such as \cite{Steele04}. (That said,
notation is not completely standardized; what I call \textquotedblleft AM-HM
inequality\textquotedblright\ is dubbed \textquotedblleft HM-AM
inequality\textquotedblright\ in \cite{Steele04}.)}

\begin{itemize}
\item Apply \href{https://en.wikipedia.org/wiki/Jensen_inequality}{Jensen's
inequality} to the convex function $\mathbb{R}^{+}\rightarrow\mathbb{R}%
^{+},\ x\mapsto\dfrac{1}{x}$.

\item Apply \href{https://en.wikipedia.org/wiki/Cauchy-Schwarz_inequality}{the
Cauchy-Schwarz inequality} to get
\begin{align*}
&  \left(  a_{1}+a_{2}+\cdots+a_{n}\right)  \left(  \dfrac{1}{a_{1}}+\dfrac
{1}{a_{2}}+\cdots+\dfrac{1}{a_{n}}\right) \\
&  \geq\left(  \underbrace{\sqrt{a_{1}\dfrac{1}{a_{1}}}+\sqrt{a_{2}\dfrac
{1}{a_{2}}}+\cdots+\sqrt{a_{n}\dfrac{1}{a_{n}}}}_{=n}\right)  ^{2}=n^{2}.
\end{align*}

\item Apply \href{https://en.wikipedia.org/wiki/HM-GM-AM-QM_inequalities}{the
AM-HM inequality}.

\item Apply \href{https://en.wikipedia.org/wiki/HM-GM-AM-QM_inequalities}{the
AM-GM inequality} twice, then multiply.

\item There is a direct proof, too: First, recall the famous inequality
\begin{equation}
\dfrac{u}{v}+\dfrac{v}{u}\geq2, \label{pf.lem.indep.cauchy-lame.uvvu}%
\end{equation}
which holds for any two positive reals $u$ and $v$. (This follows by observing
that $\dfrac{u}{v}+\dfrac{v}{u}-2=\dfrac{\left(  u-v\right)  ^{2}}{uv}\geq0$.)
Now,
\begin{align*}
&  \left(  a_{1}+a_{2}+\cdots+a_{n}\right)  \left(  \dfrac{1}{a_{1}}+\dfrac
{1}{a_{2}}+\cdots+\dfrac{1}{a_{n}}\right) \\
&  =\left(  \sum_{i=1}^{n}a_{i}\right)  \left(  \sum_{j=1}^{n}\dfrac{1}{a_{j}%
}\right)  =\sum_{i=1}^{n}\ \ \sum_{j=1}^{n}a_{i}\dfrac{1}{a_{j}}=\sum
_{i=1}^{n}\ \ \sum_{j=1}^{n}\dfrac{a_{i}}{a_{j}}\\
&  =\dfrac{1}{2}\left(  \sum_{i=1}^{n}\ \ \sum_{j=1}^{n}\dfrac{a_{i}}{a_{j}%
}+\sum_{i=1}^{n}\ \ \sum_{j=1}^{n}\dfrac{a_{i}}{a_{j}}\right)
\ \ \ \ \ \ \ \ \ \ \left(  \text{since }x=\dfrac{1}{2}\left(  x+x\right)
\text{ for any }x\in\mathbb{R}\right) \\
&  =\dfrac{1}{2}\left(  \sum_{i=1}^{n}\ \ \sum_{j=1}^{n}\dfrac{a_{i}}{a_{j}%
}+\sum_{j=1}^{n}\ \ \sum_{i=1}^{n}\dfrac{a_{j}}{a_{i}}\right)
\ \ \ \ \ \ \ \ \ \ \left(
\begin{array}
[c]{c}%
\text{here, we renamed }i\text{ and }j\text{ as }j\text{ and }i\\
\text{in the second double sum}%
\end{array}
\right) \\
&  =\dfrac{1}{2}\left(  \sum_{i=1}^{n}\ \ \sum_{j=1}^{n}\dfrac{a_{i}}{a_{j}%
}+\sum_{i=1}^{n}\ \ \sum_{j=1}^{n}\dfrac{a_{j}}{a_{i}}\right)
\ \ \ \ \ \ \ \ \ \ \left(
\begin{array}
[c]{c}%
\text{here, we swapped the two}\\
\text{summation signs in the}\\
\text{second double sum}%
\end{array}
\right) \\
&  =\dfrac{1}{2}\sum_{i=1}^{n}\ \ \sum_{j=1}^{n}\underbrace{\left(
\dfrac{a_{i}}{a_{j}}+\dfrac{a_{j}}{a_{i}}\right)  }_{\substack{\geq
2\\\text{(by (\ref{pf.lem.indep.cauchy-lame.uvvu}))}}}\geq\dfrac{1}%
{2}\underbrace{\sum_{i=1}^{n}\ \ \sum_{j=1}^{n}2}_{=n^{2}\cdot2}=\dfrac{1}%
{2}n^{2}\cdot2=n^{2},
\end{align*}
from which the claim of Lemma \ref{lem.indep.cauchy-lame} follows.
\end{itemize}
\end{proof}

\begin{proof}
[Proof of Corollary \ref{cor.indep.anti-turan}.]Write the multigraph $G$ as
$G=\left(  V,E,\varphi\right)  $. Thus, $\left\vert V\right\vert =n$ and
$\left\vert E\right\vert =m$. We WLOG assume that $V=\left\{  1,2,\ldots
,n\right\}  $ (since $\left\vert V\right\vert =n$). Hence,%
\begin{align*}
\sum_{v=1}^{n}\deg v  &  =\sum_{v\in V}\deg v=2\cdot\underbrace{\left\vert
E\right\vert }_{=m}\ \ \ \ \ \ \ \ \ \ \left(  \text{by Proposition
\ref{prop.deg.euler}}\right) \\
&  =2m.
\end{align*}
However, Theorem \ref{thm.indep.lower-bd} yields that $G$ has an independent
set of size%
\begin{align*}
&  \geq\sum_{v\in V}\dfrac{1}{1+\deg v}=\sum_{v=1}^{n}\dfrac{1}{1+\deg
v}\ \ \ \ \ \ \ \ \ \ \left(  \text{since }V=\left\{  1,2,\ldots,n\right\}
\right) \\
&  \geq\dfrac{n^{2}}{\sum_{v=1}^{n}\left(  1+\deg v\right)  }%
\ \ \ \ \ \ \ \ \ \ \left(
\begin{array}
[c]{c}%
\text{by Lemma \ref{lem.indep.cauchy-lame}, applied to the }n\text{
positive}\\
\text{reals }a_{v}=1+\deg v\text{ for all }v\in\left\{  1,2,\ldots,n\right\}
\end{array}
\right) \\
&  =\dfrac{n^{2}}{n+2m}\ \ \ \ \ \ \ \ \ \ \left(  \text{since }\sum_{v=1}%
^{n}\left(  1+\deg v\right)  =n+\underbrace{\sum_{v\in V}\deg v}%
_{=2m}=n+2m\right)  .
\end{align*}
This proves Corollary \ref{cor.indep.anti-turan}.
\end{proof}

\subsection{A proof of Turan's theorem}

Recall Turan's theorem (Theorem \ref{thm.sg.turan}), whose proof we have not
given so far. Now is the time. For the sake of convenience, let me repeat the
statement of the theorem:

\begin{theorem}
[Turan's theorem]\label{thm.indep.turan}Let $r$ be a positive integer. Let $G$
be a simple graph with $n$ vertices and $e$ edges. Assume that%
\[
e>\dfrac{r-1}{r}\cdot\dfrac{n^{2}}{2}.
\]
Then, there exist $r+1$ distinct vertices of $G$ that are mutually adjacent
(i.e., any two distinct vertices among these $r+1$ vertices are adjacent).
\end{theorem}

We can now easily derive it from Corollary \ref{cor.indep.anti-turan}:

\begin{proof}
[Proof of Theorem \ref{thm.indep.turan}.]Write the simple graph $G$ as
$G=\left(  V,E\right)  $. Thus, $\left\vert V\right\vert =n$ and $\left\vert
E\right\vert =e$ and $E\subseteq\mathcal{P}_{2}\left(  V\right)  $.

Let $\overline{E}:=\mathcal{P}_{2}\left(  V\right)  \setminus E$. Thus, the
set $\overline{E}$ consists of all \textquotedblleft
non-edges\textquotedblright\ of $G$ -- that is, of all $2$-element subsets of
$V$ that are not edges of $G$. Clearly,%
\[
\left\vert \overline{E}\right\vert =\left\vert \mathcal{P}_{2}\left(
V\right)  \setminus E\right\vert =\underbrace{\left\vert \mathcal{P}%
_{2}\left(  V\right)  \right\vert }_{=\dbinom{n}{2}}-\underbrace{\left\vert
E\right\vert }_{=e}=\dbinom{n}{2}-e.
\]

Now, let $\overline{G}$ be the simple graph $\left(  V,\overline{E}\right)  $.
This simple graph $\overline{G}$ is called the \textbf{complementary graph} of
$G$; it has $n$ vertices and $\left\vert \overline{E}\right\vert =\dbinom
{n}{2}-e$ edges.\footnote{For example, if $%
\begin{tikzpicture}[scale=0.8]
\begin{scope}[every node/.style={circle,thick,draw=green!60!black}]
\node(A) at (0:2) {$1$};
\node(B) at (360/5:2) {$2$};
\node(C) at (2*360/5:2) {$3$};
\node(D) at (3*360/5:2) {$4$};
\node(E) at (4*360/5:2) {$5$};
\end{scope}
\node(X) at (-3.2, 0) {$G = $};
\begin{scope}[every edge/.style={draw=black,very thick}]
\path[-] (A) edge (B) (B) edge (C) (C) edge (D) (D) edge (E);
\end{scope}
\end{tikzpicture}%
$ , then $%
\begin{tikzpicture}[scale=0.8]
\begin{scope}[every node/.style={circle,thick,draw=green!60!black}]
\node(A) at (0:2) {$1$};
\node(B) at (360/5:2) {$2$};
\node(C) at (2*360/5:2) {$3$};
\node(D) at (3*360/5:2) {$4$};
\node(E) at (4*360/5:2) {$5$};
\end{scope}
\node(X) at (-3.2, 0) {$\overline{G} = $};
\begin{scope}[every edge/.style={draw=black,very thick}]
\path
[-] (A) edge (C) (C) edge (E) (E) edge (B) (B) edge (D) (D) edge (A) (A) edge (E);
\end{scope}
\end{tikzpicture}%
$ .} Hence, Corollary \ref{cor.indep.anti-turan} (applied to $\overline{G}$
and $\dbinom{n}{2}-e$ instead of $G$ and $m$) yields that $\overline{G}$ has
an independent set of size%
\[
\geq\dfrac{n^{2}}{n+2\cdot\left(  \dbinom{n}{2}-e\right)  }.
\]
Let $S$ be this independent set. Its size is%
\[
\left\vert S\right\vert \geq\dfrac{n^{2}}{n+2\cdot\left(  \dbinom{n}%
{2}-e\right)  }=\dfrac{n^{2}}{n+n\left(  n-1\right)  -2e}=\dfrac{n^{2}}%
{n^{2}-2e}>r
\]
(this inequality follows by high-school algebra from $e>\dfrac{r-1}{r}%
\cdot\dfrac{n^{2}}{2}$). Hence, $\left\vert S\right\vert \geq r+1$ (since
$\left\vert S\right\vert $ and $r$ are integers). However, $S$ is an
independent set of $\overline{G}$. Thus, any two distinct vertices in $S$ are
non-adjacent in $\overline{G}$ and therefore adjacent in $G$ (by the
definition of $\overline{G}$). Since $\left\vert S\right\vert \geq r+1$, we
have thus found $r+1$ (or more) distinct vertices of $G$ that are mutually
adjacent in $G$. This proves Theorem \ref{thm.indep.turan}.
\end{proof}

Several other beautiful proofs of Theorem \ref{thm.indep.turan} can be found
in \cite[Chapter 41]{AigZie} and \cite[\S 1.2]{Zhao23}.

\subsection{\label{sec.ker}A brief introduction to digraph kernels}

\subsubsection{Definitions and main results}

We shall now curtly touch on the topic of independent sets in directed graphs,
and particularly on the so-called \textbf{kernels}, which combine independence
with \textbf{absorptiveness} (a digraph analogue of domination). The relevant
definitions are very simple:

\begin{definition}
Two vertices $u$ and $v$ of a multidigraph $D$ are said to be
\textbf{adjacent} if they are adjacent in the undirected graph
$D^{\operatorname*{und}}$. (In other words, they are adjacent if and only if
$D$ has an arc with source $u$ and target $v$ or an arc with source $v$ and
target $u$.)
\end{definition}

\begin{definition}
Let $u$ be a vertex of a multidigraph $D$. An \textbf{outneighbor} of $u$
means a vertex $v$ of $D$ such that $D$ has an arc with source $u$ and target
$v$.
\end{definition}

\begin{definition}
\label{def.ker.ker}Let $D$ be a multidigraph with vertex set $V$. Let $S$ be a
subset of $V$. Then:

\begin{enumerate}
\item[\textbf{(a)}] The set $S$ is called \textbf{independent} (or an
\textbf{independent set} of $D$) if no two elements of $S$ are adjacent. In
other words, it is called independent if it is an independent set of the
undirected graph $D^{\operatorname*{und}}$.

\item[\textbf{(b)}] The set $S$ is called \textbf{absorbing} (or an
\textbf{absorbing set} of $D$) if it has the following property: Each vertex
$v\in V\setminus S$ has at least one outneighbor in $S$. (Informally, this
property is saying that each vertex of $D$ that does not itself belong to $S$
is just \textquotedblleft one step away\textquotedblright\ from $S$.)

\item[\textbf{(c)}] The set $S$ is called a \textbf{kernel of }$D$ if it is
both independent and absorbing.
\end{enumerate}
\end{definition}

The concept of a kernel originates in the game-theoretical work of von Neumann
and Morgenstern in \cite[\S 65]{NeuMor04}.

\begin{example}
\label{exa.ker.ker1}\ \ 

\begin{enumerate}
\item[\textbf{(a)}] Let $D$ be the simple digraph%
\[%
\begin{tikzpicture}[scale=2]
\begin{scope}[every node/.style={circle,thick,draw=green!60!black}]
\node(1) at (0,1) {$1$};
\node(2) at (1,1) {$2$};
\node(3) at (2,1) {$3$};
\node(4) at (2,0) {$4$};
\end{scope}
\begin{scope}[every edge/.style={draw=black,very thick}]
\path[->] (1) edge (2) (2) edge (3) (3) edge (4) (2) edge (4);
\end{scope}
\end{tikzpicture}%
\ \ .
\]
Then, the set $\left\{  1,3\right\}  $ is independent but not absorbing (since
the vertex $4$ has no outneighbor in $\left\{  1,3\right\}  $). On the other
hand, the set $\left\{  2,4\right\}  $ is absorbing but not independent (since
its two elements $2$ and $4$ are adjacent). However, the set $\left\{
1,4\right\}  $ is both independent and absorbing, thus a kernel of $D$. It is
easy to see that this is the only kernel of $D$.

\item[\textbf{(b)}] Consider the $5$-cycle digraph $\overrightarrow{C}_{5}$,
as defined (and drawn) in Example \ref{exa.spanning-arbor.Cn}. Then, the set
$\left\{  1,3\right\}  $ is independent but not absorbing (since the vertex
$4$ has no outneighbor in $\left\{  1,3\right\}  $). On the other hand, the
set $\left\{  1,3,5\right\}  $ is absorbing but not independent (since the
vertices $1$ and $5$ are adjacent). It is easy to see that each independent
set of $\overrightarrow{C}_{5}$ has at most $2$ elements, while each absorbing
set of $\overrightarrow{C}_{5}$ has at least $3$ elements. Hence,
$\overrightarrow{C}_{5}$ has no kernel.

Likewise, the $n$-th cycle digraph $\overrightarrow{C}_{n}$ for any odd $n$
has no kernel.

\item[\textbf{(c)}] On the other hand, the $n$-th cycle digraph
$\overrightarrow{C}_{n}$ for any even $n$ has two kernels: $\left\{
1,3,5,\ldots,n-1\right\}  $ and $\left\{  2,4,6,\ldots,n\right\}  $.

\item[\textbf{(d)}] Now let $D$ be the simple digraph%
\[%
\begin{tikzpicture}[scale=2]
\begin{scope}[every node/.style={circle,thick,draw=green!60!black}]
\node(1) at (0,1) {$1$};
\node(2) at (1,1) {$2$};
\node(3) at (2,1) {$3$};
\node(4) at (2,0) {$4$};
\node(5) at (1,0) {$5$};
\end{scope}
\begin{scope}[every edge/.style={draw=black,very thick}]
\path[->] (1) edge (2) (2) edge (3) (3) edge (4) (4) edge (5) (5) edge (2);
\end{scope}
\end{tikzpicture}%
\ \ .
\]
Then, $D$ has two kernels: $\left\{  1,3,5\right\}  $ and $\left\{
2,4\right\}  $.
\end{enumerate}
\end{example}

As this example shows, a digraph $D$ can have multiple kernels or just one or
none at all. Thus it is natural to ask for necessary and sufficient criteria
for the existence of kernels. Example \ref{exa.ker.ker1} \textbf{(b)} might
suggest that an odd-length cycle would always prevent a kernel from existing
(just like an odd cycle in an undirected graph makes a proper $2$-coloring
impossible), but this is not the case: e.g., the digraph%
\[%
\begin{tikzpicture}
\begin{scope}[every node/.style={circle,thick,draw=green!60!black}]
\node(A) at (0:2) {$1$};
\node(B) at (360/5:2) {$2$};
\node(C) at (2*360/5:2) {$3$};
\node(D) at (3*360/5:2) {$4$};
\node(E) at (4*360/5:2) {$5$};
\node(O) at (0,0) {$0$};
\end{scope}
\begin{scope}[every edge/.style={draw=black,very thick}]
\path[->] (A) edge (B) (B) edge (C) (C) edge (D) (D) edge (E) (E) edge (A);
\path[->] (A) edge (O) (B) edge (O) (C) edge (O) (D) edge (O) (E) edge (O);
\end{scope}
\end{tikzpicture}%
\]
has a kernel $\left\{  0\right\}  $ despite having a length-$5$ cycle.
However, the converse is true, as Moses Richardson proved in 1946
\cite[Theorem]{Richar46}:

\begin{theorem}
[Richardson]\label{thm.ker.rich}Let $D$ be a multidigraph. If $D$ has no
odd-length cycles, then $D$ has a kernel.
\end{theorem}

Moreover, if $D$ has no cycles at all, then the kernel is unique, as shown by
von Neumann and Morgenstern in 1944 \cite[\S 65]{NeuMor04}:

\begin{theorem}
[von Neumann, Morgenstern]\label{thm.ker.vnm}Let $D$ be a multidigraph. If $D$
has no cycles, then $D$ has a unique kernel.
\end{theorem}

In fact, there is a \textquotedblleft dual\textquotedblright\ to Theorem
\ref{thm.ker.rich} for uniqueness instead of existence:

\begin{theorem}
\label{thm.ker.uni}Let $D$ be a multidigraph. If $D$ has no even-length
cycles\footnotemark, then $D$ has at most one kernel.
\end{theorem}

\footnotetext{An \textquotedblleft even-length cycle\textquotedblright\ means
a cycle that has even length.}

We shall prove these three theorems soon. The interested reader can consult
Langlois's thesis \cite{Langlo23} for a comprehensive survey of recent and
older work on kernels, including several generalizations of Theorem
\ref{thm.ker.rich}.

\subsubsection{Proofs}

\begin{fineprint}
\begin{proof}
[Proof of Theorem \ref{thm.ker.rich}.] We proceed by strong induction on
$\left\vert \operatorname*{V}\left(  D\right)  \right\vert $ (that is, the \#
of vertices of $D$).

\textit{Base case:} Theorem \ref{thm.ker.rich} is obvious if $\left\vert
\operatorname*{V}\left(  D\right)  \right\vert =0$, because in this case the
empty set $\varnothing$ is a kernel of $D$.

\textit{Induction step:} Let $N$ be a positive integer. Assume (as the
induction hypothesis) that Theorem \ref{thm.ker.rich} holds whenever
$\left\vert \operatorname*{V}\left(  D\right)  \right\vert <N$. We must now
prove that Theorem \ref{thm.ker.rich} holds whenever $\left\vert
\operatorname*{V}\left(  D\right)  \right\vert =N$.

So let $D$ be a multidigraph with $\left\vert \operatorname*{V}\left(
D\right)  \right\vert =N$ that has no odd-length cycles. We must show that $D$
has a kernel.

The digraph $D$ has at least one vertex (since $\left\vert \operatorname*{V}%
\left(  D\right)  \right\vert =N$ is positive). Thus, Theorem
\ref{thm.mdg.sinkcomp} shows that $D$ has at least one sink component. Pick
such a sink component, and call it $C$. Then, $C$ is a strong component of
$D$. Hence, the induced subdigraph $D\left[  C\right]  $ is strongly connected
(by Proposition \ref{prop.mdg.strong-component.ind-connected}). Moreover, this
subdigraph $D\left[  C\right]  $ has no odd-length cycles (since it is a
subdigraph of $D$, which has no odd-length cycles). In other words, the
statement B'2 of Theorem \ref{thm.coloring.2-color-eq-dig} is satisfied for
$D\left[  C\right]  $ instead of $D$. Hence, the implication
B'2$\Longrightarrow$B'1 in Theorem \ref{thm.coloring.2-color-eq-dig} (applied
to $D\left[  C\right]  $ instead of $D$) shows that the underlying undirected
graph $\left(  D\left[  C\right]  \right)  ^{\operatorname*{und}}$ has a
proper $2$-coloring. Let $f:C\rightarrow\left\{  1,2\right\}  $ be this proper
$2$-coloring.

The map $f:C\rightarrow\left\{  1,2\right\}  $ allows us to decompose the set
$C$ into two disjoint subsets%
\begin{align*}
C_{1}  & :=\left\{  v\in C\ \mid\ f\left(  v\right)  =1\right\}
\ \ \ \ \ \ \ \ \ \ \text{and}\\
C_{2}  & :=\left\{  v\in C\ \mid\ f\left(  v\right)  =2\right\}  ,
\end{align*}
which satisfy $C_{1}\cup C_{2}=C$. Since $C_{1}\cup C_{2}=C$ is nonempty
(because $C$ is a strong component of $D$), we conclude that at least one of
these two subsets $C_{1}$ and $C_{2}$ is nonempty. We WLOG assume that $C_{1}$
is nonempty (indeed, if $C_{2}$ is nonempty, then we can swap the two colors
$1$ and $2$ in the proper $2$-coloring $f$, and thus make $C_{1}$ nonempty instead).

Now, let $L:=\operatorname*{V}\left(  D\right)  \setminus C$ be the complement
of $C$ in $\operatorname*{V}\left(  D\right)  $. Thus, $C\cap L=\varnothing$
and $C\cup L=\operatorname*{V}\left(  D\right)  $.

Define a subset $P$ of $L$ by%
\[
P:=\left\{  v\in L\ \mid\ v\text{ has \textbf{no} outneighbor in }%
C_{1}\right\}  .
\]

\begin{noncompile}
Let $L_{2}:=L\setminus L_{1}$ be the complement of $L_{1}$ in $L$. Thus,
$L_{1}\cap L_{2}=\varnothing$ and $L_{1}\cup L_{2}=L$.
\end{noncompile}

Now consider the induced subdigraph $D\left[  P\right]  $ of $D$. Its number
of vertices is%
\begin{align*}
\left\vert \operatorname*{V}\left(  D\left[  P\right]  \right)  \right\vert  &
=\left\vert P\right\vert \leq\left\vert L\right\vert
\ \ \ \ \ \ \ \ \ \ \left(  \text{since }P\subseteq L\right)  \\
& =\underbrace{\left\vert \operatorname*{V}\left(  D\right)  \right\vert
}_{=N}-\underbrace{\left\vert C\right\vert }_{\substack{>0\\\text{(since
}C\text{ is nonempty)}}}\ \ \ \ \ \ \ \ \ \ \left(  \text{since }%
L=\operatorname*{V}\left(  D\right)  \setminus C\right)  \\
& <N.
\end{align*}
Moreover, this digraph $D\left[  P\right]  $ has no odd-length cycles (since
it is a subdigraph of $D$, which has no odd-length cycles). Hence, our
induction hypothesis ensures that Theorem \ref{thm.ker.rich} can be applied to
$D\left[  P\right]  $ instead of $D$. As a result, we conclude that $D\left[
P\right]  $ has a kernel. Let $K$ be this kernel of $D\left[  P\right]  $.
Thus, $K\subseteq P$ is an independent set of $D\left[  P\right]  $ and
simultaneously an absorbing set of $D\left[  P\right]  $.

We now claim that $K\cup C_{1}$ is a kernel of $D$. In order to prove this, we
must verify the following two claims:

\begin{statement}
\textit{Claim 1:} The set $K\cup C_{1}$ is independent (with respect to $D$).
\end{statement}

\begin{statement}
\textit{Claim 2:} The set $K\cup C_{1}$ is absorbing (with respect to $D$).
\end{statement}

\begin{proof}
[Proof of Claim 1.] Assume the contrary. Thus, $K\cup C_{1}$ contains two
adjacent vertices $u$ and $v$. Consider these $u$ and $v$.

Since $u$ and $v$ are adjacent, we know that $D$ has an arc from $u$ to $v$ or
an arc from $v$ to $u$. We WLOG assume that $D$ has an arc from $u$ to $v$
(because in the other case, we can just swap $u$ with $v$). Let $a$ be this
arc. The existence of this arc shows that $v$ is an outneighbor of $u$.

Each of the vertices $u$ and $v$ belongs to $K\cup C_{1}$, hence belongs to
$K$ or to $C_{1}$. Thus, we are in one of the following four cases:

\textit{Case 1:} We have $u\in K$ and $v\in K$.

\textit{Case 2:} We have $u\in K$ and $v\in C_{1}$.

\textit{Case 3:} We have $u\in C_{1}$ and $v\in K$.

\textit{Case 4:} We have $u\in C_{1}$ and $v\in C_{1}$.

Let us first consider Case 1. In this case, we have $u\in K$ and $v\in K$.
Hence, $u\in K\subseteq P$ and $v\in K\subseteq P$. Thus, $u$ and $v$ are
adjacent vertices of $D$ both belonging to $P$. In other words, $u$ and $v$
are adjacent vertices of $D\left[  P\right]  $. These two vertices both belong
to $K$ (since $u\in K$ and $v\in K$). But this impossible, since $K$ is an
independent set of $D\left[  P\right]  $ and thus contains no two adjacent
vertices. Thus, we obtained a contradiction in Case 1.

Let us next consider Case 2. In this case, we have $u\in K$ and $v\in C_{1}$.
But $v$ is an outneighbor of $u$. Hence, $u$ has an outneighbor in $C_{1}$
(namely, $v$). However, $u\in K\subseteq P$ shows that $u$ has no outneighbor
in $C_{1}$ (by the definition of $P$). These two sentences contradict each
other. Thus, we have found a contradiction in Case 2.

Next, let us consider Case 3. In this case, we have $u\in C_{1}$ and $v\in K$.
We have $u\in C_{1}\subseteq C$. In other words, the source of the arc $a$
lies in $C$ (since this source is $u$). But $C$ is a sink component of $D$.
Hence, if the source of an arc of $D$ lies in $C$, then its target also lies
in $C$ (by the definition of \textquotedblleft sink
component\textquotedblright). Applying this to the arc $a$, we conclude that
the target of $a$ lies in $C$ (since the source of $a$ lies in $C$). But this
is absurd, since the target of $a$ is $v\in K\subseteq P\subseteq
L=\operatorname*{V}\left(  D\right)  \setminus C$. Thus, we have found a
contradiction in Case 3.

Finally, let us consider Case 4. In this case, we have  $u\in C_{1}$ and $v\in
C_{1}$. By the definition of $C_{1}$, this entails that $f\left(  u\right)
=1$ and $f\left(  v\right)  =1$. But $u$ and $v$ are two vertices in $C$
(since $u\in C_{1}\subseteq C$ and $v\in C_{1}\subseteq C$) that are adjacent
in $D$, and thus are two adjacent vertices in the induced subdigraph $D\left[
C\right]  $. In other words, $u$ and $v$ are two adjacent vertices in the
undirected graph $\left(  D\left[  C\right]  \right)  ^{\operatorname*{und}}$.
Since $f:C\rightarrow\left\{  1,2\right\}  $ is a proper $2$-coloring of this
graph, we thus obtain $f\left(  u\right)  \neq f\left(  v\right)  $. But this
contradicts $f\left(  u\right)  =1=f\left(  v\right)  $. Hence, we have found
a contradiction in Case 4.

We now have found contradictions in all four cases. Hence, our assumption was
wrong, and Claim 1 is proved.
\end{proof}

\begin{proof}
[Proof of Claim 2.] We must prove that each vertex $v\in\operatorname*{V}%
\left(  D\right)  \setminus\left(  K\cup C_{1}\right)  $ has at least one
outneighbor in  $K\cup C_{1}$. So let $v\in\operatorname*{V}\left(  D\right)
\setminus\left(  K\cup C_{1}\right)  $ be a vertex. We must show that $v$ has
an outneighbor in  $K\cup C_{1}$. This is the claim that we shall now prove.

We are in one of the following two cases:

\textit{Case 1:} We have $v\in C$.

\textit{Case 2:} We have $v\notin C$.

We first consider Case 1. In this case, we have $v\in C$. But $v\in
\operatorname*{V}\left(  D\right)  \setminus\left(  K\cup C_{1}\right)  $ and
thus $v\notin K\cup C_{1}$, so that $v\notin C_{1}$. Combining this with $v\in
C=C_{1}\cup C_{2}$, we obtain $v\in\left(  C_{1}\cup C_{2}\right)  \setminus
C_{1}\subseteq C_{2}$. Therefore, $f\left(  v\right)  =2$ (by the definition
of $C_{2}$).

The set $C_{1}$ is nonempty, thus contains some vertex $w\in C_{1}$. Consider
this $w$. Thus, $w\in C_{1}\subseteq C$. But the digraph $D\left[  C\right]  $
is strongly connected, and thus has a path from $v$ to $w$ (since both $v$ and
$w$ belong to $C$). This path cannot have length $0$ (since $v\notin C_{1}$
and $w\in C_{1}$ entail $v\neq w$), and thus has a first arc. This first arc
has source $v$; let $u$ be its target. Thus, $u$ is an outneighbor of $v$ in
$D\left[  C\right]  $. Hence, $u$ and $v$ are adjacent vertices of the digraph
$D\left[  C\right]  $, hence adjacent vertices of the undirected graph
$\left(  D\left[  C\right]  \right)  ^{\operatorname*{und}}$. Since $f$ is a
proper $2$-coloring of $\left(  D\left[  C\right]  \right)
^{\operatorname*{und}}$, we thus obtain $f\left(  u\right)  \neq f\left(
v\right)  =2$. Hence, $f\left(  u\right)  =1$ (since $f\left(  u\right)  $ is
either $1$ or $2$). Therefore, $u\in C_{1}$ (by the definition of $C_{1}$), so
that $u\in C_{1}\subseteq K\cup C_{1}$. Thus, $v$ has an outneighbor in $K\cup
C_{1}$ (namely, $u$). This proves our claim in Case 1.

Let us now consider Case 2. In this case, we have $v\notin C$. Hence,
$v\in\operatorname*{V}\left(  D\right)  \setminus C=L$. On the other hand,
$v\in\operatorname*{V}\left(  D\right)  \setminus\left(  K\cup C_{1}\right)  $
and thus $v\notin K\cup C_{1}$, so that $v\notin K$.

If $v\notin P$, then $v$ has an outneighbor in $C_{1}$ (by the definition of
$P$) and thus has an outneighbor in $K\cup C_{1}$ (since $C_{1}\subseteq K\cup
C_{1}$). Thus, we are done if $v\notin P$. Hence, for the rest of this proof,
we WLOG assume that $v\in P$. Thus, $v\in P\setminus K$ (since $v\notin K$).

Recall that $K$ is an absorbing set of $D\left[  P\right]  $. Hence, from
$v\in P\setminus K$, we conclude that $v$ has at least one outneighbor in $K$.
Thus, $v$ has at least one outneighbor in $K\cup C_{1}$ (since $K\subseteq
K\cup C_{1}$). Thus, we have proved our claim in Case 2.

Thus, in both cases, we have proved our claim that $v$ has an outneighbor in
$K\cup C_{1}$. This completes the proof of Claim 2.
\end{proof}

Now, the subset $K\cup C_{1}$ of $\operatorname*{V}\left(  D\right)  $ is
independent (by Claim 1) and absorbing (by Claim 2); hence, it is a kernel of
$D$ (by the definition of a kernel). Thus, $D$ has a kernel. This proves
Theorem \ref{thm.ker.rich} for our $D$, and thus completes the induction step.
\end{proof}

\begin{proof}
[Proof of Theorem \ref{thm.ker.uni}.] Assume that $D$ has no even-length
cycles. We must show that $D$ has at most one kernel.

Assume the contrary. Thus, $D$ has two distinct kernels $K$ and $L$.

A vertex will be called \textbf{disputed} if it belongs to $K\setminus L$ or
to $L\setminus K$. Let $\Delta=\left(  K\setminus L\right)  \cup\left(
L\setminus K\right)  $ be the set of all disputed vertices of $D$. This set
$\Delta$ is nonempty (since $K$ and $L$ are distinct). Now we shall show the following:

\begin{statement}
\textit{Claim 1:} Each vertex $v\in K\setminus L$ has at least one outneighbor
in $L\setminus K$.
\end{statement}

\begin{proof}
[Proof of Claim 1.] Let $v\in K\setminus L$ be a vertex. We must show that $v$
has at least one outneighbor in $L\setminus K$.

Indeed, $v\in K\setminus L\subseteq\operatorname*{V}\left(  D\right)
\setminus L$. But $L$ is a kernel of $D$, and thus an absorbing set. Hence,
from $v\in\operatorname*{V}\left(  D\right)  \setminus L$, we conclude that
$v$ has an outneighbor in $L$. Let $w\in L$ be this outneighbor. Then, $v$ and
$w$ are two adjacent vertices. But $K$ is a kernel of $D$, and thus an
independent set. Hence, $K$ contains no two adjacent vertices. If we had $w\in
K$, then $v$ and $w$ would be two adjacent vertices in $K$ (since $v\in
K\setminus L\subseteq K$ and $w\in K$), which would contradict the previous
sentence. Thus, we cannot have $w\in K$. Combining this with $w\in L$, we
obtain $w\in L\setminus K$. Hence, $v$ has an outneighbor in $L\setminus K$
(namely, $w$). This proves Claim 1.
\end{proof}

\begin{statement}
\textit{Claim 2:} Each vertex $v\in L\setminus K$ has at least one outneighbor
in $K\setminus L$.
\end{statement}

\begin{proof}
[Proof of Claim 2.] Analogous to Claim 1 (just swap $K$ with $L$).
\end{proof}

\begin{statement}
\textit{Claim 3:} Each vertex $v\in\Delta$ has at least one outneighbor in
$\Delta$.
\end{statement}

\begin{proof}
[Proof of Claim 3.] This follows from Claim 1 and Claim 2 (since
$\Delta=\left(  K\setminus L\right)  \cup\left(  L\setminus K\right)  $).
\end{proof}

Now, consider the induced subdigraph $D\left[  \Delta\right]  $ of $D$ on the
set $\Delta$. This subdigraph $D\left[  \Delta\right]  $ has at least one
vertex (since $\Delta$ is nonempty). Moreover, each vertex $v$ of this
subdigraph $D\left[  \Delta\right]  $ has at least one outneighbor in this
subdigraph (by Claim 3), i.e., has outdegree $\deg_{D\left[  \Delta\right]
}^{+}v>0$ (where $\deg_{D\left[  \Delta\right]  }^{+}v$ denotes the outdegree
of $v$ with respect to the digraph $D\left[  \Delta\right]  $). Hence,
Exercise \ref{exe.4.2} \textbf{(a)} (applied to $D\left[  \Delta\right]  $
instead of $D$) shows that the digraph $D\left[  \Delta\right]  $ has a cycle.

On the other hand, each vertex $v\in\Delta$ belongs either to $K\setminus L$
or to $L\setminus K$ (since $\Delta=\left(  K\setminus L\right)  \cup\left(
L\setminus K\right)  $), and these two possibilities are mutually exclusive
(since the sets $K\setminus L$ and $L\setminus K$ are disjoint). Hence, we can
define a map $f:\Delta\rightarrow\left\{  1,2\right\}  $ by%
\[
f\left(  v\right)  =%
\begin{cases}
1, & \text{if }v\in K\setminus L;\\
2, & \text{if }v\in L\setminus K
\end{cases}
\ \ \ \ \ \ \ \ \ \ \text{for each }v\in\Delta.
\]
Consider this map $f$. We now claim the following:

\begin{statement}
\textit{Claim 4:} The map $f$ is a proper $2$-coloring of the graph $\left(
D\left[  \Delta\right]  \right)  ^{\operatorname*{und}}$.
\end{statement}

\begin{proof}
[Proof of Claim 4.] Let $u$ and $v$ be two adjacent vertices of $\left(
D\left[  \Delta\right]  \right)  ^{\operatorname*{und}}$. We must show that
$f\left(  u\right)  \neq f\left(  v\right)  $.

The vertices $u$ and $v$ are adjacent in the graph $\left(  D\left[
\Delta\right]  \right)  ^{\operatorname*{und}}$, hence also in the digraph
$D\left[  \Delta\right]  $, thus also in the digraph $D$. But $K$ is a kernel
of $D$, thus an independent set of $D$. Hence, $K$ contains no two adjacent
vertices. Therefore, the two adjacent vertices $u$ and $v$ cannot both belong
to $K$. Hence, we have either $u\notin K$ or $v\notin K$ (or both). We WLOG
assume that $u\notin K$ (otherwise, swap $u$ with $v$). Hence, $u\notin
K\setminus L$. But $u$ is a vertex of $\left(  D\left[  \Delta\right]
\right)  ^{\operatorname*{und}}$; thus, $u\in\Delta=\left(  K\setminus
L\right)  \cup\left(  L\setminus K\right)  $. Combining this with $u\notin
K\setminus L$, we obtain $u\in\left(  \left(  K\setminus L\right)  \cup\left(
L\setminus K\right)  \right)  \setminus\left(  K\setminus L\right)  \subseteq
L\setminus K$. Therefore, the definition of $f$ yields $f\left(  u\right)  =2$.

We have shown that $u\notin K$ or $v\notin K$. Likewise, we can show that
$u\notin L$ or $v\notin L$ (since $L$ is also a kernel). Since $u\notin L$ is
impossible (because $u\in L\setminus K\subseteq L$), we thus obtain $v\notin
L$. Hence, $v\notin L\setminus K$. But $v$ is a vertex of $\left(  D\left[
\Delta\right]  \right)  ^{\operatorname*{und}}$; thus, $v\in\Delta=\left(
K\setminus L\right)  \cup\left(  L\setminus K\right)  $. Combining this with
$v\notin L\setminus K$, we obtain $v\in\left(  \left(  K\setminus L\right)
\cup\left(  L\setminus K\right)  \right)  \setminus\left(  L\setminus
K\right)  \subseteq K\setminus L$. Therefore, the definition of $f$ yields
$f\left(  v\right)  =1$. Thus, $f\left(  u\right)  =2\neq1=f\left(  v\right)
$. Claim 4 is thus proved.
\end{proof}

Claim 4 shows that the graph $\left(  D\left[  \Delta\right]  \right)
^{\operatorname*{und}}$ has a proper $2$-coloring. Hence, the implication
B1$\Longrightarrow$B2 in Theorem \ref{thm.coloring.2-color-eq} (applied to
$G=\left(  D\left[  \Delta\right]  \right)  ^{\operatorname*{und}}$) shows
that $\left(  D\left[  \Delta\right]  \right)  ^{\operatorname*{und}}$ has no
cycles of odd length. Thus, the digraph $D\left[  \Delta\right]  $ has no
cycles of odd length either (since any such cycle would also be a cycle of
$\left(  D\left[  \Delta\right]  \right)  ^{\operatorname*{und}}$). But we
know that the digraph $D\left[  \Delta\right]  $ has a cycle. Thus, this cycle
must have even length (since $D\left[  \Delta\right]  $ has no cycles of odd
length). So we have shown that $D\left[  \Delta\right]  $ has an even-length
cycle. Thus, $D$ has an even-length cycle (since $D\left[  \Delta\right]  $ is
a subgraph of $D$). This contradicts the assumption that $D$ has no
even-length cycles. This contradiction shows that our assumption was false.
Hence, Theorem \ref{thm.ker.uni} is proved.
\end{proof}

\begin{proof}
[Proof of Theorem \ref{thm.ker.vnm}.] The digraph $D$ has no cycles, hence no
odd-length cycles and no even-length cycles. Thus, Theorem \ref{thm.ker.rich}
shows that $D$ has a kernel, whereas Theorem \ref{thm.ker.uni} shows that $D$
has at most one kernel. Combining these, we conclude that  $D$ has a unique
kernel. This proves Theorem \ref{thm.ker.vnm}. (Note that there are much
easier proofs of Theorem \ref{thm.ker.vnm} around.)
\end{proof}
\end{fineprint}

\section{Matchings}

\subsection{Introduction}

Independent sets of a graph consist of vertices that \textquotedblleft have no
edges in common\textquotedblright\ (i.e., no two belong to the same edge).

In a sense, \textbf{matchings} are the dual notion to this: they consist of
edges that \textquotedblleft have no vertices in common\textquotedblright%
\ (i.e., no two contain the same vertex). Here is the formal definition:

\begin{definition}
\label{def.match.match}Let $G=\left(  V,E,\varphi\right)  $ be a loopless multigraph.

\begin{enumerate}
\item[\textbf{(a)}] A \textbf{matching} of $G$ means a subset $M$ of $E$ such
that no two distinct edges in $M$ have a common endpoint.

\item[\textbf{(b)}] If $M$ is a matching of $G$, then an $M$\textbf{-edge}
shall mean an edge that belongs to $M$.

\item[\textbf{(c)}] If $M$ is a matching of $G$, and if $v\in V$ is any
vertex, then we say that $v$ is \textbf{matched} in $M$ (or \textbf{saturated}
in $M$) if $v$ is an endpoint of an $M$-edge. In this case, this latter
$M$-edge is necessarily unique (since $M$ is a matching), and is called the
$M$\textbf{-edge} of $v$. The other endpoint of this $M$-edge (i.e., its
endpoint different from $v$) is called the $M$\textbf{-partner} of $v$.

\item[\textbf{(d)}] A matching $M$ of $G$ is said to be \textbf{perfect} if
each vertex of $G$ is matched in $M$.

\item[\textbf{(e)}] Let $A$ be a subset of $V$. A matching $M$ of $G$ is said
to be $A$\textbf{-complete} if each vertex in $A$ is matched in $M$.
\end{enumerate}
\end{definition}

Thus, a matching $M$ of a multigraph $G=\left(  V,E,\varphi\right)  $ is
perfect if and only if it is $V$-complete.

\begin{example}
Let $G$ be the following simple graph:%
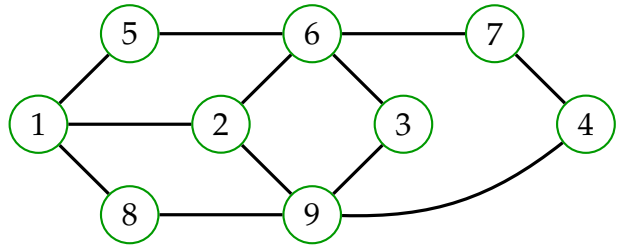
\[%
\begin{tikzpicture}[scale=1.2]
\begin{scope}[every node/.style={circle,thick,draw=green!60!black}]
\node(00) at (0, 0) {$1$};
\node(20) at (2, 0) {$2$};
\node(40) at (4, 0) {$3$};
\node(60) at (6, 0) {$4$};
\node(11) at (1, 1) {$5$};
\node(31) at (3, 1) {$6$};
\node(51) at (5, 1) {$7$};
\node(1-1) at (1, -1) {$8$};
\node(3-1) at (3, -1) {$9$};
\end{scope}
\begin{scope}[every edge/.style={draw=black,very thick}, every loop/.style={}]
\path[-] (00) edge (11) edge (20) edge (1-1);
\path[-] (31) edge (11) edge (20) edge (40) edge (51);
\path[-] (3-1) edge (1-1) edge (20) edge (40);
\path[-] (60) edge (51) edge[bend left=20] (3-1);
\end{scope}
\end{tikzpicture}%
\ \ .
\]
Then:

\begin{itemize}
\item The set $\left\{  12,\ 36,\ 47\right\}  $ is a matching of $G$. If we
call this set $M$, then the vertices matched in $M$ are $1,2,3,4,6,7$, and
their respective $M$-partners are $2,1,6,7,3,4$. This matching is not perfect,
but it is (for example) $\left\{  1,3,4\right\}  $-complete and $\left\{
1,2,3,4,6,7\right\}  $-complete.

\item The set $\left\{  12,\ 36,\ 67\right\}  $ is not a matching of $G$,
since the two distinct edges $36$ and $67$ from this set have a common endpoint.

\item The sets $\varnothing$, $\left\{  36\right\}  $, $\left\{
15,\ 29,\ 36,\ 47\right\}  $ are matchings of $G$ as well.
\end{itemize}
\end{example}

We see that any matching \textquotedblleft pairs up\textquotedblright\ some
vertices using the existing edges of the graph. Clearly, the $M$-partner of
the $M$-partner of a vertex $v$ is $v$ itself. Also, no two distinct vertices
have the same $M$-partner (since otherwise, their $M$-edges would have a
common endpoint).

\begin{remark}
\label{rmk.match.degrees}A matching of a loopless multigraph $G=\left(
V,E,\varphi\right)  $ can also be characterized as a subset $M$ of its edge
set $E$ such that all vertices of the spanning subgraph $\left(
V,M,\varphi\mid_{M}\right)  $ have degree $\leq1$.
\end{remark}

\begin{warning}
If a multigraph $G$ has loops, then most authors additionally require that a
matching must not contain any loops. This ensures that Remark
\ref{rmk.match.degrees} remains valid.
\end{warning}

Here are some natural questions:

\begin{itemize}
\item Does a given graph $G$ have a perfect matching?

\item If not, can we find a maximum-size matching?

\item What about an $A$-complete matching for a given $A\subseteq V$ ?
\end{itemize}

Here are some examples:

\begin{example}
\label{exa.gridgraph}Let $n$ and $m$ be two positive integers. The Cartesian
product $P_{n}\times P_{m}$ of the $n$-th path graph $P_{n}$ and the $m$-th
path graph $P_{m}$ is known as the $\left(  n,m\right)  $\textbf{-grid graph},
as it looks as follows:%
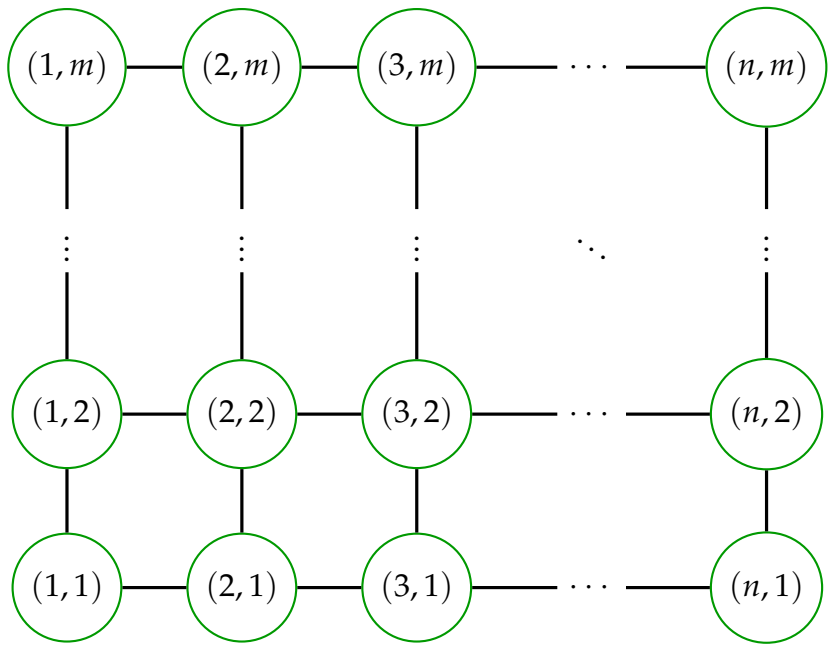
\[%
%
\ \ .
\]

\begin{enumerate}
\item[\textbf{(a)}] If $n$ is even, then
\[
\left\{  \left\{  \left(  i,j\right)  ,\ \left(  i+1,j\right)  \right\}
\ \mid\ i\text{ is odd, while }j\text{ is arbitrary}\right\}
\]
is a perfect matching of $P_{n}\times P_{m}$. For example, here is this
perfect matching for $n=4$ and $m=3$ (we have drawn all edges that do
\textbf{not} belong to this matching as dotted lines):%
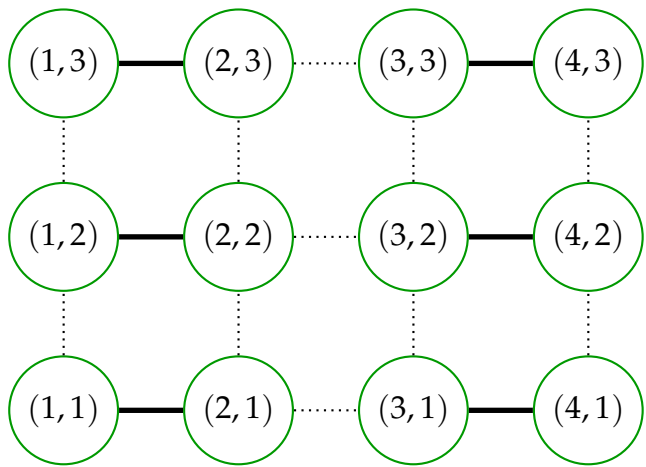
\[%
%
\]

\item[\textbf{(b)}] Likewise, if $m$ is even, then
\[
\left\{  \left\{  \left(  i,j\right)  ,\ \left(  i,j+1\right)  \right\}
\ \mid\ j\text{ is odd, while }i\text{ is arbitrary}\right\}
\]
is a perfect matching of $P_{n}\times P_{m}$.

\item[\textbf{(c)}] If $n$ and $m$ are both odd, then $P_{n}\times P_{m}$ has
no perfect matching. Indeed, any loopless multigraph $G$ with an odd number of
vertices cannot have a perfect matching, since each edge of the matching
covers exactly $2$ vertices.
\end{enumerate}
\end{example}

\begin{example}
\label{exa.match.penta-antlers}The \textquotedblleft pentagon with two
antlers\textquotedblright\ $C_{5}^{\prime\prime}$ (this is my notation,
hopefully sufficiently natural) is the following graph:%
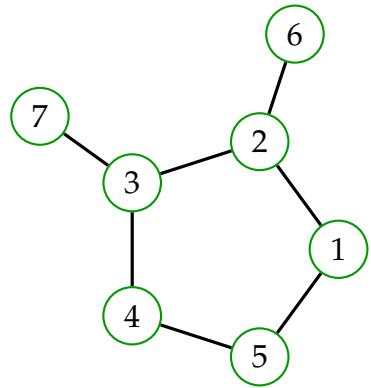
\[%
\begin{tikzpicture}
\begin{scope}[every node/.style={circle,thick,draw=green!60!black}]
\node(1) at (0*360/5 : 1.5) {$1$};
\node(2) at (1*360/5 : 1.5) {$2$};
\node(3) at (2*360/5 : 1.5) {$3$};
\node(4) at (3*360/5 : 1.5) {$4$};
\node(5) at (4*360/5 : 1.5) {$5$};
\node(6) at (1*360/5 : 3) {$6$};
\node(7) at (2*360/5 : 3) {$7$};
\end{scope}
\begin{scope}[every edge/.style={draw=black,very thick}]
\path[-] (1) edge (2) (2) edge (3) (3) edge (4) (4) edge (5) (5) edge (1);
\path[-] (2) edge (6) (3) edge (7);
\end{scope}
\end{tikzpicture}%
\ \ .
\]
It has no perfect matching. This is easiest to see as follows: The graph
$C_{5}^{\prime\prime}$ is loopless, so each edge contains exactly two
vertices. Thus, any matching $M$ of $C_{5}^{\prime\prime}$ matches exactly
$2\cdot\left\vert M\right\vert $ vertices. In particular, any matching of
$C_{5}^{\prime\prime}$ matches an even number of vertices. Since the total
number of vertices $C_{5}^{\prime\prime}$ is odd, this entails that
$C_{5}^{\prime\prime}$ has no perfect matching.

What is the maximum size of a matching of $C_{5}^{\prime\prime}$ ? The
matching $\left\{  12,\ 34\right\}  $ of $C_{5}^{\prime\prime}$ has size $2$
and cannot be improved by adding any new edges. Thus, one is tempted to
believe that the maximum size of a matching is $2$. However, this is not the
case. Indeed, the matching $\left\{  12,\ 37,\ 45\right\}  $ has size $3$.
This latter matching is actually maximum-size.
\end{example}

Example \ref{exa.match.penta-antlers} shows that when searching for a
maximum-size matching, it is not sufficient to just keep adding edges until no
further edges can be added; this strategy may lead to a non-improvable but
non-maximum matching. This suggests that finding a maximum-size matching may
be one of those hard problems like finding a maximum-size independent set. But
no -- there is a polynomial-time algorithm! It's known as the
\href{https://en.wikipedia.org/wiki/Blossom_algorithm}{Edmonds blossom
algorithm}, and it has a running time of $O\left(  \left\vert E\right\vert
\cdot\left\vert V\right\vert ^{2}\right)  $; however, it is too complicated to
be covered in this course. We shall here focus on a simple case of the problem
that is already interesting enough and almost as useful as the general case.

Namely, we shall study matchings of \textbf{bipartite graphs}.

\subsection{\label{sec.matching.bip}Bipartite graphs}

\subsubsection{Definition and examples}

\begin{definition}
\label{def.match.bip}A \textbf{bipartite graph} means a triple $\left(
G,X,Y\right)  $, where

\begin{itemize}
\item $G=\left(  V,E,\varphi\right)  $ is a multigraph, and

\item $X$ and $Y$ are two disjoint subsets of $V$ such that $X\cup Y=V$ and
such that each edge of $G$ has one endpoint in $X$ and one endpoint in $Y$.
\end{itemize}
\end{definition}

\begin{example}
\label{exa.match.bipartite-lr}Consider the $6$-th cycle graph $C_{6}$:%
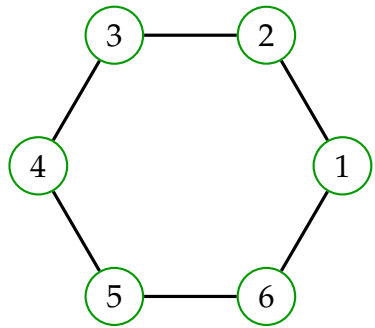
\[%
\begin{tikzpicture}
\begin{scope}[every node/.style={circle,thick,draw=green!60!black}]
\node(A) at (0:2) {$1$};
\node(B) at (60:2) {$2$};
\node(C) at (120:2) {$3$};
\node(D) at (180:2) {$4$};
\node(E) at (240:2) {$5$};
\node(F) at (300:2) {$6$};
\end{scope}
\begin{scope}[every edge/.style={draw=black,very thick}]
\path
[-] (A) edge (B) (B) edge (C) (C) edge (D) (D) edge (E) (E) edge (F) (F) edge (A);
\end{scope}
\end{tikzpicture}%
\ \ .
\]
Then, $\left(  C_{6},\ \left\{  1,3,5\right\}  ,\ \left\{  2,4,6\right\}
\right)  $ is a bipartite graph, since each edge of $G$ has one endpoint in
$\left\{  1,3,5\right\}  $ and one endpoint in $\left\{  2,4,6\right\}  $.
Also, $\left(  C_{6},\ \left\{  2,4,6\right\}  ,\ \left\{  1,3,5\right\}
\right)  $ is a bipartite graph.

Note that a bipartite graph $\left(  G,X,Y\right)  $ is not just the graph $G$
but rather the whole package consisting of the graph $G$ and the subsets $X$
and $Y$. Two different bipartite graphs can have the same underlying graph $G$
but different choices of $X$ and $Y$. For example, the two bipartite graphs
$\left(  C_{6},\ \left\{  1,3,5\right\}  ,\ \left\{  2,4,6\right\}  \right)  $
and $\left(  C_{6},\ \left\{  2,4,6\right\}  ,\ \left\{  1,3,5\right\}
\right)  $ are different.

We typically draw a bipartite graph $\left(  G,X,Y\right)  $ by drawing the
graph $G$ in such a way that the vertices in $X$ are aligned along one
vertical line and the vertices $Y$ are aligned along another, with the former
line being left of the latter. Thus, for example, the bipartite graph $\left(
C_{6},\ \left\{  1,3,5\right\}  ,\ \left\{  2,4,6\right\}  \right)  $ can be
drawn as follows:%
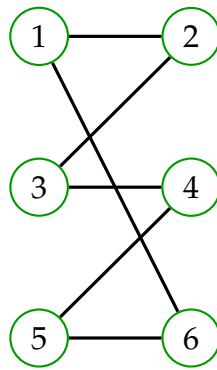
\[%
%
\ \ .
\]
Similarly, the bipartite graph $\left(  C_{6},\ \left\{  2,4,6\right\}
,\ \left\{  1,3,5\right\}  \right)  $ can be drawn as follows:%
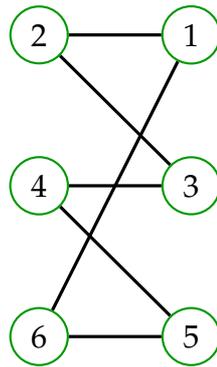
\[%
%
\ \ .
\]

\end{example}

This example suggests the following terminology:

\begin{definition}
\label{def.match.lr}Let $\left(  G,X,Y\right)  $ be a bipartite graph. We
shall refer to the vertices in $X$ as the \textbf{left vertices} of this
bipartite graph. We shall refer to the vertices in $Y$ as the \textbf{right
vertices} of this bipartite graph. Moreover, the edges of $G$ will be called
the \textbf{edges} of this bipartite graph.
\end{definition}

Thus, each edge of a bipartite graph joins one left vertex with one right
vertex.\footnote{Another example of a bipartite graph is $\left(
K_{n,m},\ \left\{  1,2,\ldots,n\right\}  ,\ \left\{  -1,-2,\ldots,-m\right\}
\right)  $, where $K_{n,m}$ is as in Exercise \ref{exe.MTT.Knm}.}

\subsubsection{Bipartite graphs as graphs with a proper $2$-coloring}

Bipartite graphs are \textquotedblleft the same as\textquotedblright%
\ multigraphs with a proper $2$-coloring. To wit:

\begin{proposition}
\label{prop.match.bip-is-2col}Let $G=\left(  V,E,\varphi\right)  $ be a multigraph.

\begin{enumerate}
\item[\textbf{(a)}] If $\left(  G,X,Y\right)  $ is a bipartite graph, then the
map
\begin{align*}
f:V  &  \rightarrow\left\{  1,2\right\}  ,\\
v  &  \mapsto%
\begin{cases}
1, & \text{if }v\in X;\\
2, & \text{if }v\in Y
\end{cases}
\end{align*}
is a proper $2$-coloring of $G$.

\item[\textbf{(b)}] Conversely, if $f:V\rightarrow\left\{  1,2\right\}  $ is a
proper $2$-coloring of $G$, then $\left(  G,V_{1},V_{2}\right)  $ is a
bipartite graph, where we set%
\[
V_{i}:=\left\{  \text{all vertices with color }i\right\}
\ \ \ \ \ \ \ \ \ \ \text{for each }i\in\left\{  1,2\right\}  .
\]

\item[\textbf{(c)}] These constructions are mutually inverse. (That is, going
from a bipartite graph to a proper $2$-coloring and back again results in the
original bipartite graph, whereas going from a proper $2$-coloring to a
bipartite graph and back again results in the original $2$-coloring.)
\end{enumerate}
\end{proposition}

\begin{proof}
An exercise in understanding the definitions.
\end{proof}

\begin{proposition}
\label{prop.match.bip-has-no-odd-circs}Let $\left(  G,X,Y\right)  $ be a
bipartite graph. Then, the graph $G$ has no circuits of odd length. In
particular, $G$ has no loops or triangles.
\end{proposition}

\begin{proof}
By Proposition \ref{prop.match.bip-is-2col} \textbf{(a)}, we know that $G$ has
a proper $2$-coloring. Hence, the $2$-coloring equivalence theorem (Theorem
\ref{thm.coloring.2-color-eq}) shows that $G$ has no circuits of odd length.
In particular, $G$ has no loops or triangles (since these would yield circuits
of length $1$ or $3$, respectively).
\end{proof}

\subsubsection{Neighbor sets}

We need another piece of notation:

\begin{definition}
Let $G=\left(  V,E,\varphi\right)  $ be any multigraph. Let $U$ be a subset of
$V$. Then, we define%
\[
N\left(  U\right)  :=\left\{  v\in V\ \mid\ v\text{ has a neighbor in
}U\right\}  .
\]
This is called the \textbf{neighbor set} of $U$.
\end{definition}

\begin{example}
If $G$ is the \textquotedblleft pentagon with antlers\textquotedblright%
\ $C_{5}^{\prime\prime}$ from Example \ref{exa.match.penta-antlers}, then%
\begin{align*}
N\left(  \left\{  1,5,6\right\}  \right)   &  =\left\{  1,2,4,5\right\}  ;\\
N\left(  \left\{  1\right\}  \right)   &  =\left\{  2,5\right\}  ;\\
N\left(  \varnothing\right)   &  =\varnothing.
\end{align*}

\end{example}

For bipartite graphs, the neighbor set has a nice property:

\begin{proposition}
\label{prop.bipartite.NAY} Let $\left(  G,X,Y\right)  $ be a bipartite graph.
Let $A\subseteq X$. Then,%
\[
N\left(  A\right)  \subseteq Y.
\]

\end{proposition}

\begin{proof}
Let $v\in N\left(  A\right)  $. Thus, the vertex $v$ has a neighbor in $A$ (by
definition of $N\left(  A\right)  $). Let $w$ be this neighbor. Then, $w\in
A\subseteq X$, so that $w\notin Y$ (since the bipartiteness of $\left(
G,X,Y\right)  $ shows that the sets $X$ and $Y$ are disjoint).

There exists some edge that has endpoints $v$ and $w$ (since $w$ is a neighbor
of $v$). This edge must have an endpoint in $Y$ (since the bipartiteness of
$\left(  G,X,Y\right)  $ shows that each edge of $G$ has one endpoint in $Y$).
In other words, one of $v$ and $w$ must belong to $Y$ (since the endpoints of
this edge are $v$ and $w$). Since $w\notin Y$, we thus conclude that $v\in Y$.

Thus, we have shown that $v\in Y$ for each $v\in N\left(  A\right)  $. In
other words, $N\left(  A\right)  \subseteq Y$.
\end{proof}

\begin{exercise}
\label{exe.8.6}Let $\left(  G,X,Y\right)  $ be a bipartite graph. Prove that
\[
\sum_{A\subseteq X}\left(  -1\right)  ^{\left\vert A\right\vert }\left[
N\left(  A\right)  =Y\right]  =\sum_{B\subseteq Y}\left(  -1\right)
^{\left\vert B\right\vert }\left[  N\left(  B\right)  =X\right]
\]
(where we are using the Iverson bracket notation).
\end{exercise}

\subsection{\label{sec.matching.hall}Hall's marriage theorem}

How can we tell whether a bipartite graph has a perfect matching? an
$X$-complete matching?

\subsubsection{Generalities}

First, to keep the suspense, let us prove some trivialities:

\begin{proposition}
\label{prop.match.bipar-match-basics}Let $\left(  G,X,Y\right)  $ be a
bipartite graph. Let $M$ be a matching of $G$. Then:

\begin{enumerate}
\item[\textbf{(a)}] The $M$-partner of a vertex $x\in X$ (if it exists)
belongs to $Y$.

The $M$-partner of a vertex $y\in Y$ (if it exists) belongs to $X$.

\item[\textbf{(b)}] We have $\left\vert M\right\vert \leq\left\vert
X\right\vert $ and $\left\vert M\right\vert \leq\left\vert Y\right\vert $.

\item[\textbf{(c)}] If $M$ is $X$-complete, then $\left\vert X\right\vert
\leq\left\vert Y\right\vert $.

\item[\textbf{(d)}] If $M$ is perfect, then $\left\vert X\right\vert
=\left\vert Y\right\vert $.

\item[\textbf{(e)}] If $\left\vert M\right\vert \geq\left\vert X\right\vert $,
then $M$ is $X$-complete.

\item[\textbf{(f)}] If $M$ is $X$-complete and we have $\left\vert
X\right\vert =\left\vert Y\right\vert $, then $M$ is perfect.
\end{enumerate}
\end{proposition}

\begin{proof}
Each edge of $G$ has an endpoint in $X$ and an endpoint in $Y$ (since $\left(
G,X,Y\right)  $ is a bipartite graph). Thus, in particular, each $M$-edge has
an endpoint in $X$ and an endpoint in $Y$. Moreover, no two $M$-edges share a
common endpoint (since $M$ is a matching). \medskip

\textbf{(a)} This follows from the fact that each $M$-edge has an endpoint in
$X$ and an endpoint in $Y$. \medskip

\textbf{(b)} Recall that each $M$-edge has an endpoint in $X$. Since no two
$M$-edges share a common endpoint, we thus have found at least $\left\vert
M\right\vert $ many endpoints in $X$. This entails $\left\vert M\right\vert
\leq\left\vert X\right\vert $. Similarly, $\left\vert M\right\vert
\leq\left\vert Y\right\vert $. \medskip

\textbf{(c)} Assume that $M$ is $X$-complete. Hence, each vertex in $X$ is
matched in $M$ and therefore has an $M$-edge that contains it. In other words,
for each vertex $x\in X$, there exists an $M$-edge $m$ such that $x$ is an
endpoint of $m$. Since no two $M$-edges share an endpoint, this yields that
there are at least $\left\vert X\right\vert $ many $M$-edges. In other words,
$\left\vert M\right\vert \geq\left\vert X\right\vert $. Hence, $\left\vert
X\right\vert \leq\left\vert M\right\vert \leq\left\vert Y\right\vert $ (by
part \textbf{(b)}). \medskip

\textbf{(d)} Assume that $M$ is perfect. Then, $M$ is both $X$-complete and
$Y$-complete. Hence, part \textbf{(c)} yields $\left\vert X\right\vert
\leq\left\vert Y\right\vert $; similarly, $\left\vert Y\right\vert
\leq\left\vert X\right\vert $. Combining these two inequalities, we obtain
$\left\vert X\right\vert =\left\vert Y\right\vert $. \medskip

\textbf{(e)} Assume that $\left\vert M\right\vert \geq\left\vert X\right\vert
$.

However, each $M$-edge has an endpoint in $X$. These endpoints are all
distinct (since no two $M$-edges share a common endpoint), and there are at
least $\left\vert X\right\vert $ many of them (since there are $\left\vert
M\right\vert $ many of them, but we have $\left\vert M\right\vert
\geq\left\vert X\right\vert $). Therefore, these endpoints must cover
\textbf{all} the vertices in $X$ (because the only way to choose $\left\vert
X\right\vert $ many distinct vertices in $X$ is to choose \textbf{all}
vertices in $X$). In other words, all the vertices in $X$ must be matched in
$M$. In other words, the matching $M$ is $X$-complete. \medskip

\textbf{(f)} Assume that $M$ is $X$-complete and that we have $\left\vert
X\right\vert =\left\vert Y\right\vert $.

The matching $M$ is $X$-complete; thus, all vertices $x\in X$ are matched in
$M$. The $M$-partners of all these vertices $x\in X$ belong to $Y$ (by
Proposition \ref{prop.match.bipar-match-basics} \textbf{(a)}), and are also
matched in $M$. Hence, at least $\left\vert X\right\vert $ many vertices in
$Y$ must be matched in $M$ (since these $M$-partners are all
distinct\footnote{because the $M$-partners of distinct vertices are
distinct}). In other words, at least $\left\vert Y\right\vert $ many vertices
in $Y$ must be matched in $M$ (since $\left\vert X\right\vert =\left\vert
Y\right\vert $). This means that \textbf{all} vertices in $Y$ are matched in
$M$ (since \textquotedblleft at least $\left\vert Y\right\vert $ many vertices
in $Y$\textquotedblright\ means \textquotedblleft all vertices in
$Y$\textquotedblright). Since we also know that all vertices $x\in X$ are
matched in $M$, we thus conclude that all vertices of $G$ are matched in $M$.
In other words, the matching $M$ is perfect.
\end{proof}

\begin{example}
\label{exa.match.VV}Consider the bipartite graph%
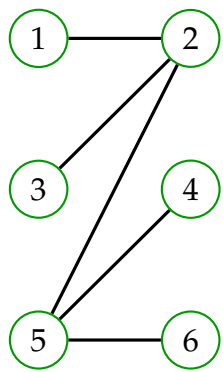
\[%
%
\]
(drawn as explained in Example \ref{exa.match.bipartite-lr}). Does this graph
have a perfect matching? No, because the two left vertices $1$ and $3$ would
necessarily have the same partner in such a matching (since their only
possible partner is $2$).

Similarly, the bipartite graph%
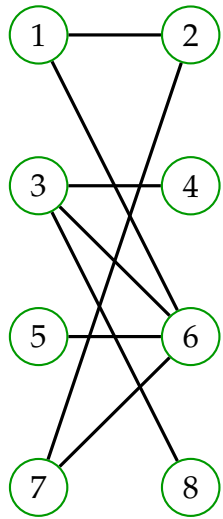
\[%
%
\]
has no perfect matching, since the three left vertices $1$, $5$ and $7$ have
only two potential partners (viz., $2$ and $6$).
\end{example}

So we see that a subset $A\subseteq X$ satisfying $\left\vert N\left(
A\right)  \right\vert <\left\vert A\right\vert $ is an obstruction to the
existence of an $X$-complete matching. Let us state this in a positive way:

\begin{proposition}
\label{prop.match.HMT-easy}Let $\left(  G,X,Y\right)  $ be a bipartite graph.
Let $A$ be a subset of $X$. Assume that $G$ has an $X$-complete matching.
Then, $\left\vert N\left(  A\right)  \right\vert \geq\left\vert A\right\vert $.
\end{proposition}

\begin{proof}
Let $V$ be the vertex set of $G$. We assumed that $G$ has an $X$-complete
matching. Let $M$ be such a matching. Thus, each $x\in X$ has an $M$-partner.
The map%
\begin{align*}
\mathbf{p}:X  &  \rightarrow V,\\
x  &  \mapsto\left(  \text{the }M\text{-partner of }x\right)
\end{align*}
is injective (since two distinct vertices cannot have the same $M$-partner).
Thus, $\left\vert \mathbf{p}\left(  A\right)  \right\vert =\left\vert
A\right\vert $ (because any injective map preserves the size of a subset).
However, $\mathbf{p}\left(  A\right)  \subseteq N\left(  A\right)  $, because
the $M$-partner of an element of $A$ will always belong to $N\left(  A\right)
$. Hence, $\left\vert \mathbf{p}\left(  A\right)  \right\vert \leq\left\vert
N\left(  A\right)  \right\vert $. Thus, $\left\vert N\left(  A\right)
\right\vert \geq\left\vert \mathbf{p}\left(  A\right)  \right\vert =\left\vert
A\right\vert $, qed.
\end{proof}

\subsubsection{Hall's marriage theorem}

Proposition \ref{prop.match.HMT-easy} gives a necessary condition for the
existence of an $X$-complete matching in a bipartite graph $\left(
G,X,Y\right)  $. Interestingly, this condition is also sufficient:

\begin{theorem}
[Hall's marriage theorem, short: HMT]\label{thm.match.HMT}Let $\left(
G,X,Y\right)  $ be a bipartite graph. Assume that each subset $A$ of $X$
satisfies $\left\vert N\left(  A\right)  \right\vert \geq\left\vert
A\right\vert $. (This assumption is called the \textquotedblleft\textbf{Hall
condition}\textquotedblright.)

Then, $G$ has an $X$-complete matching.
\end{theorem}

This is called \textquotedblleft marriage theorem\textquotedblright\ because
one can interpret a bipartite graph as a dating scene, with $X$ being the guys
and $Y$ the ladies. A guy $x$ and a lady $y$ are adjacent if and only if they
are interested in one another. Thus, an $X$-complete matching is a way of
marrying off each guy to some lady he is mutually interested in (without
allowing polygamy). This is a classical model for bipartite graphs and appears
all across the combinatorics literature; to my knowledge, however, no
real-life applications have been found along these lines. Nevertheless, Hall's
marriage theorem can be applied in many other situations, such as logistics
(although its generalizations, which we will soon see, are even more useful in
that). Philip Hall has originally invented the theorem in 1935 (in a somewhat
obfuscated form), motivated (I believe) by a problem about finite groups. So
did Wilhelm Maak, also in 1935, for use in analysis (defining a notion of
integrals for almost-periodic functions).

There are many proofs of Hall's marriage theorem, some pretty easy. Two short
and self-contained proofs can be found in \cite[\S 12.5.2]{LeLeMe} and in
\cite[Theorem 3.9]{Harju14}. I will tease you by leaving the theorem unproved
for several pages, while exploring some of its many consequences. Afterwards,
I will give two proofs of Hall's marriage theorem:

\begin{itemize}
\item one proof using the theory of \textbf{network flows} (Section
\ref{sec.flows.hk}) -- an elegant theory created for use in
logistics\footnote{and, more generally, operations research} in the 1950s that
has proved to be quite useful in combinatorics. Among other consequences, this
proof will also provide a polynomial-time algorithm for actually finding a
maximum matching in a bipartite graph (Theorem \ref{thm.match.HMT} by itself
does not help here).

\item another proof using the \textbf{Gallai--Milgram theorem} (Subsection
\ref{subsec.paths.gm.apps}) -- an elegant and surprising property of paths in digraphs.
\end{itemize}

\subsection{\label{sec.matching.koenig}K\"{o}nig and Hall--K\"{o}nig}

Hall's marriage theorem is famous for its many forms and versions, most of
which are \textquotedblleft secretly\textquotedblright\ equivalent to it
(i.e., can be derived from it and conversely can be used to derive it without
too much trouble). We will start with one that is known as \textbf{K\"{o}nig's
theorem} (discovered independently by D\'{e}nes K\H{o}nig and Jen\H{o}
Egerv\'{a}ry in 1931). This relies on the notion of a \textbf{vertex cover}.
Here is its definition:

\begin{definition}
\label{def.match.vercover}Let $G=\left(  V,E,\varphi\right)  $ be a
multigraph. A \textbf{vertex cover} of $G$ means a subset $C$ of $V$ such that
each edge of $G$ contains at least one vertex in $C$.
\end{definition}

\begin{example}
Let $n\geq1$. What are the vertex covers of the complete graph $K_{n}$ ?

A quick thought reveals that any subset $S$ of $\left\{  1,2,\ldots,n\right\}
$ that has at least $n-1$ elements is a vertex cover of $K_{n}$. (In fact,
$K_{n}$ has no loops, so that each edge of $K_{n}$ contains two different
vertices, and thus at least one of these two vertices belongs to $S$.) On the
other hand, a subset $S$ with fewer than $n-1$ vertices will never be a vertex
cover of $K_{n}$ (since there will be at least two distinct vertices that
don't belong to $S$, and the edge that joins these two vertices contains no
vertex in $S$).
\end{example}

\begin{example}
Let $G=\left(  V,E,\varphi\right)  $ be the graph from Example
\ref{exa.match.VV}. Then, the set $\left\{  2,5\right\}  $ is a vertex cover
of $G$. Of course, any subset of $V$ that contains $\left\{  2,5\right\}  $ as
a subset will thus also be a vertex cover of $G$.
\end{example}

Note that the notion of a vertex cover is (in some sense) \textquotedblleft
dual\textquotedblright\ to the notion of an edge cover, which we defined in
Exercise \ref{exe.1.6}. For those getting confused, here is a convenient table
(courtesy of
\href{https://canvas.dartmouth.edu/courses/46201/files/folder/Notes}{Nadia
Lafreni\`{e}re, Math 38, Spring 2021}):%
\[%
\begin{tabular}
[c]{|c|c|c|}\hline
\textbf{a ...} & \textbf{is a set of ...} & \textbf{that contains
...}\\\hline\hline
matching & edges & at most one edge per vertex\\\hline
edge cover & edges & at least one edge per vertex\\\hline
independent set & vertices & at most one vertex per edge\\\hline
vertex cover & vertices & at least one vertex per edge\\\hline
\end{tabular}
\]

The notion of vertex covers is also somewhat reminiscent of the notion of
dominating sets; here is the precise relation:

\begin{remark}
Each vertex cover of a multigraph $G$ is a dominating set (as long as $G$ has
no vertices of degree $0$). But the converse is not true.
\end{remark}

\begin{proposition}
\label{prop.match.konig-easy}Let $G$ be a loopless multigraph.

Let $m$ be the largest size of a matching of $G$.

Let $c$ be the smallest size of a vertex cover of $G$.

Then, $m\leq c$.
\end{proposition}

\begin{proof}
By the definition of $m$, we know that $G$ has a matching $M$ of size $m$.

By the definition of $c$, we know that $G$ has a vertex cover $C$ of size $c$.

Consider these $M$ and $C$. Every $M$-edge $e\in M$ contains at least one
vertex in $C$ (since $C$ is a vertex cover). Thus, we can define a map
$f:M\rightarrow C$ that sends each $M$-edge $e$ to some vertex in $C$ that is
contained in $e$. (If there are two such vertices, then we just pick one of
them at random.) This map $f$ is injective, because no two $M$-edges contain
the same vertex (after all, $M$ is a matching). Thus, we have found an
injective map from $M$ to $C$ (namely, $f$). Therefore, $\left\vert
M\right\vert \leq\left\vert C\right\vert $. But the definitions of $M$ and $C$
show that $\left\vert M\right\vert =m$ and $\left\vert C\right\vert =c$. Thus,
$m=\left\vert M\right\vert \leq\left\vert C\right\vert =c$, and Proposition
\ref{prop.match.konig-easy} is proved.
\end{proof}

In general, we can have $m<c$ in Proposition \ref{prop.match.konig-easy}.
However, for a \textbf{bipartite} graph, equality reigns:

\begin{theorem}
[K\"{o}nig's theorem]\label{thm.match.konig}Let $\left(  G,X,Y\right)  $ be a
bipartite graph.

Let $m$ be the largest size of a matching of $G$.

Let $c$ be the smallest size of a vertex cover of $G$.

Then, $m=c$.
\end{theorem}

Both Hall's and K\"{o}nig's theorems easily follow from the following theorem:

\begin{theorem}
[Hall--K\"{o}nig matching theorem]\label{thm.match.HKMT}Let $\left(
G,X,Y\right)  $ be a bipartite graph. Then, there exist a matching $M$ of $G$
and a subset $U$ of $X$ such that%
\[
\left\vert M\right\vert \geq\left\vert N\left(  U\right)  \right\vert
+\left\vert X\right\vert -\left\vert U\right\vert .
\]

\end{theorem}

We will prove this theorem in Section \ref{sec.flows.hk} and again in
Subsection \ref{subsec.paths.gm.apps}. For now, let us show that Hall's
marriage theorem (Theorem \ref{thm.match.HMT}), K\"{o}nig's theorem (Theorem
\ref{thm.match.konig}) and the Hall--K\"{o}nig matching theorem (Theorem
\ref{thm.match.HKMT}) are mutually equivalent. More precisely, we will explain
how to derive the first two from the third, and outline the reverse derivations.

\begin{proof}
[Proof of Theorem \ref{thm.match.HMT} using Theorem \ref{thm.match.HKMT}%
.]Assume that Theorem \ref{thm.match.HKMT} has already been proved.

Theorem \ref{thm.match.HKMT} yields that there exist a matching $M$ of $G$ and
a subset $U$ of $X$ such that%
\[
\left\vert M\right\vert \geq\left\vert N\left(  U\right)  \right\vert
+\left\vert X\right\vert -\left\vert U\right\vert .
\]
Consider these $M$ and $U$. The Hall condition shows that each subset $A$ of
$X$ satisfies $\left\vert N\left(  A\right)  \right\vert \geq\left\vert
A\right\vert $. Applying this to $A=U$, we obtain $\left\vert N\left(
U\right)  \right\vert \geq\left\vert U\right\vert $. Thus,
\[
\left\vert M\right\vert \geq\underbrace{\left\vert N\left(  U\right)
\right\vert }_{\geq\left\vert U\right\vert }+\left\vert X\right\vert
-\left\vert U\right\vert \geq\left\vert X\right\vert .
\]
Hence, the matching $M$ is $X$-complete (by Proposition
\ref{prop.match.bipar-match-basics} \textbf{(e)}). Thus, we have found an
$X$-complete matching. This proves Theorem \ref{thm.match.HMT} (assuming that
Theorem \ref{thm.match.HKMT} is true).
\end{proof}

\begin{proof}
[Proof of Theorem \ref{thm.match.konig} using Theorem \ref{thm.match.HKMT}%
.]Assume that Theorem \ref{thm.match.HKMT} has already been proved.

Write the multigraph $G$ as $G=\left(  V,E,\varphi\right)  $. Theorem
\ref{thm.match.HKMT} yields that there exist a matching $M$ of $G$ and a
subset $U$ of $X$ such that%
\begin{equation}
\left\vert M\right\vert \geq\left\vert N\left(  U\right)  \right\vert
+\left\vert X\right\vert -\left\vert U\right\vert .
\label{pf.thm.match.konig.1}%
\end{equation}
Consider these $M$ and $U$. Clearly, $\left\vert M\right\vert \leq m$ (since
$m$ is the largest size of a matching of $G$).

Let $C:=\left(  X\setminus U\right)  \cup N\left(  U\right)  $. This is a
subset of $V$. Moreover, each edge of $G$ has at least one endpoint in $C$
(this is easy to see\footnote{\textit{Proof.} Let $e$ be an edge of $G$. We
must show that $e$ has at least one endpoint in $C$.
\par
Clearly, the edge $e$ has an endpoint in $X$ (since $\left(  G,X,Y\right)  $
is a bipartite graph). Let $x$ be this endpoint. This $x$ either belongs to
$U$ or doesn't.
\par
\begin{itemize}
\item If $x$ belongs to $U$, then the other endpoint of $e$ (that is, the
endpoint distinct from $x$) belongs to $N\left(  U\right)  $ (since its
neighbor $x$ belongs to $U$) and therefore to $C$ (since $N\left(  U\right)
\subseteq\left(  X\setminus U\right)  \cup N\left(  U\right)  =C$).
\par
\item If $x$ does not belong to $U$, then $x$ belongs to $X\setminus U$ (since
$x\in X$) and therefore to $C$ (since $X\setminus U\subseteq\left(  X\setminus
U\right)  \cup N\left(  U\right)  =C$).
\end{itemize}
\par
In either of these two cases, we have found an endpoint of $e$ that belongs to
$C$. Thus, $e$ has at least one endpoint in $C$, qed.}). Hence, $C$ is a
vertex cover of $G$. Therefore, $\left\vert C\right\vert \geq c$ (since $c$ is
the smallest size of a vertex cover of $G$). The definition of $C$ yields
\begin{align*}
\left\vert C\right\vert  &  =\left\vert \left(  X\setminus U\right)  \cup
N\left(  U\right)  \right\vert \\
&  \leq\underbrace{\left\vert X\setminus U\right\vert }_{\substack{=\left\vert
X\right\vert -\left\vert U\right\vert \\\text{(since }U\subseteq X\text{)}%
}}+\left\vert N\left(  U\right)  \right\vert \ \ \ \ \ \ \ \ \ \ \left(
\text{actually an equality, but we don't care}\right) \\
&  =\left\vert X\right\vert -\left\vert U\right\vert +\left\vert N\left(
U\right)  \right\vert =\left\vert N\left(  U\right)  \right\vert +\left\vert
X\right\vert -\left\vert U\right\vert \leq\left\vert M\right\vert
\ \ \ \ \ \ \ \ \ \ \left(  \text{by (\ref{pf.thm.match.konig.1})}\right) \\
&  \leq m.
\end{align*}
Hence, $m\geq\left\vert C\right\vert \geq c$. Combining this with $m\leq c$
(which follows from Proposition \ref{prop.match.konig-easy}), we obtain $m=c$.
Thus, Theorem \ref{thm.match.konig} follows.
\end{proof}

\begin{fineprint}
Conversely, it is not hard to derive the HKMT from either Hall or K\"{o}nig:

\begin{proof}
[Proof of Theorem \ref{thm.match.HKMT} using Theorem \ref{thm.match.HMT}
(sketched).]Assume that Theorem \ref{thm.match.HMT} has already been proved.

Add a bunch of \textquotedblleft dummy vertices\textquotedblright\ to $Y$ and
join each of these \textquotedblleft dummy vertices\textquotedblright\ by a
new edge to each vertex in $X$. How many \textquotedblleft dummy
vertices\textquotedblright\ should we add? As many as it takes to ensure that
every subset $A$ of $X$ satisfies the Hall condition -- i.e., exactly
$\max\left\{  \left\vert A\right\vert -\left\vert N\left(  A\right)
\right\vert \ \mid\ A\text{ is a subset of }X\right\}  $ many.

Let $G^{\prime}$ be the resulting graph. Let also $D$ be the set of all dummy
vertices that were added to $Y$, and let $Y^{\prime}=Y\cup D$ be the set of
all right vertices of $G^{\prime}$. (The set of left vertices of $G^{\prime}$
is still $X$.) Then, the bipartite graph $\left(  G^{\prime},X,Y^{\prime
}\right)  $ satisfies the Hall condition, and therefore we can apply Theorem
\ref{thm.match.HMT} to $\left(  G^{\prime},X,Y^{\prime}\right)  $ instead of
$\left(  G,X,Y\right)  $, and conclude that the graph $G^{\prime}$ has an
$X$-complete matching. Let $M^{\prime}$ be this matching. By removing from
$M^{\prime}$ all edges that contain dummy vertices, we obtain a matching $M$
of $G$. This matching $M$ has size%
\begin{align}
\left\vert M\right\vert  &  =\left\vert M^{\prime}\right\vert
-\underbrace{\left(  \text{the number of edges that were removed from
}M^{\prime}\right)  }_{\substack{\leq\left(  \text{the number of dummy
vertices}\right)  \\\text{(since each dummy vertex is contained in at most one
}M^{\prime}\text{-edge)}}}\nonumber\\
&  \geq\left\vert M^{\prime}\right\vert -\underbrace{\left(  \text{the number
of dummy vertices}\right)  }_{\substack{=\max\left\{  \left\vert A\right\vert
-\left\vert N\left(  A\right)  \right\vert \ \mid\ A\text{ is a subset of
}X\right\}  \\\text{(by the construction of the dummy vertices)}}}\nonumber\\
&  =\left\vert M^{\prime}\right\vert -\max\left\{  \left\vert A\right\vert
-\left\vert N\left(  A\right)  \right\vert \ \mid\ A\text{ is a subset of
}X\right\}  . \label{pf.thm.match.HKMT.via.HMT.1}%
\end{align}

However, the maximum of a set is always an element of this set. Hence, there
exists a subset $U$ of $X$ such that%
\[
\max\left\{  \left\vert A\right\vert -\left\vert N\left(  A\right)
\right\vert \ \mid\ A\text{ is a subset of }X\right\}  =\left\vert
U\right\vert -\left\vert N\left(  U\right)  \right\vert .
\]
Consider this $U$. Then, (\ref{pf.thm.match.HKMT.via.HMT.1}) becomes%
\begin{align*}
\left\vert M\right\vert  &  \geq\underbrace{\left\vert M^{\prime}\right\vert
}_{\substack{\geq\left\vert X\right\vert \\\text{(since }M^{\prime}\text{ is
}X\text{-complete,}\\\text{and thus each }x\in X\text{ has}\\\text{an
}M^{\prime}\text{-edge (and these}\\\text{edges are distinct))}}%
}-\underbrace{\max\left\{  \left\vert A\right\vert -\left\vert N\left(
A\right)  \right\vert \ \mid\ A\text{ is a subset of }X\right\}
}_{=\left\vert U\right\vert -\left\vert N\left(  U\right)  \right\vert }\\
&  \geq\left\vert X\right\vert -\left(  \left\vert U\right\vert -\left\vert
N\left(  U\right)  \right\vert \right)  =\left\vert N\left(  U\right)
\right\vert +\left\vert X\right\vert -\left\vert U\right\vert .
\end{align*}
Hence, we have found a matching $M$ of $G$ and a subset $U$ of $X$ such that
$\left\vert M\right\vert \geq\left\vert N\left(  U\right)  \right\vert
+\left\vert X\right\vert -\left\vert U\right\vert $. This proves Theorem
\ref{thm.match.HKMT} (assuming that Theorem \ref{thm.match.HMT} is true).
\end{proof}

\begin{proof}
[Proof of Theorem \ref{thm.match.HKMT} from Theorem \ref{thm.match.konig}
(sketched).]Assume that Theorem \ref{thm.match.konig} has already been proved.

Let $M$ be a maximum-size matching of $G$. Let $C$ be a minimum-size vertex
cover of $G$. Then, Theorem \ref{thm.match.konig} says that $\left\vert
M\right\vert =\left\vert C\right\vert $.

Let $U:=X\setminus C$. Then, $N\left(  U\right)  \subseteq C\setminus X$
(why?). Hence, $\left\vert N\left(  U\right)  \right\vert \leq\left\vert
C\setminus X\right\vert $, so that
\[
\underbrace{\left\vert N\left(  U\right)  \right\vert }_{\leq\left\vert
C\setminus X\right\vert }+\left\vert X\right\vert -\left\vert \underbrace{U}%
_{=X\setminus C}\right\vert \leq\left\vert C\setminus X\right\vert
+\underbrace{\left\vert X\right\vert -\left\vert X\setminus C\right\vert
}_{=\left\vert C\cap X\right\vert }=\left\vert C\setminus X\right\vert
+\left\vert C\cap X\right\vert =\left\vert C\right\vert =\left\vert
M\right\vert .
\]
Hence, $\left\vert M\right\vert \geq\left\vert N\left(  U\right)  \right\vert
+\left\vert X\right\vert -\left\vert U\right\vert $. This proves Theorem
\ref{thm.match.HKMT} (assuming that Theorem \ref{thm.match.konig} is true).
\end{proof}
\end{fineprint}

Theorem \ref{thm.match.HKMT} thus occupies a convenient \textquotedblleft high
ground\textquotedblright\ between the Hall and K\"{o}nig theorems, allowing
easy access to both of them. We shall prove Theorem \ref{thm.match.HKMT} in
Section \ref{sec.flows.hk} and again in Subsection \ref{subsec.paths.gm.apps}.

\subsection{Systems of representatives}

There are two more equivalent form of the HMT that have the \textquotedblleft
advantage\textquotedblright\ that they do not rely on the notion of a graph.
When non-combinatorialists use the HMT, they often use it in one of these
forms. Here is the first form:

\begin{theorem}
[existence of SDR]\label{thm.match.SDR}Let $A_{1},A_{2},\ldots,A_{n}$ be any
$n$ sets. Assume that the union of any $p$ of these sets has size $\geq p$,
for all $p\in\left\{  0,1,\ldots,n\right\}  $. (In other words, assume that%
\[
\left\vert A_{i_{1}}\cup A_{i_{2}}\cup\cdots\cup A_{i_{p}}\right\vert \geq
p\ \ \ \ \ \ \ \ \ \ \text{for any }1\leq i_{1}<i_{2}<\cdots<i_{p}\leq n.
\]
)

Then, we can find $n$ \textbf{distinct} elements%
\[
a_{1}\in A_{1},\ \ \ \ \ \ \ \ \ \ a_{2}\in A_{2},\ \ \ \ \ \ \ \ \ \ \ldots
,\ \ \ \ \ \ \ \ \ \ a_{n}\in A_{n}.
\]

\end{theorem}

\begin{remark}
An $n$-tuple $\left(  a_{1},a_{2},\ldots,a_{n}\right)  $ of $n$ distinct
elements like this is called a \textbf{system of distinct representatives} for
our $n$ sets $A_{1},A_{2},\ldots,A_{n}$. (This is often abbreviated
\textquotedblleft SDR\textquotedblright.)
\end{remark}

\begin{example}
Take \href{https://en.wikipedia.org/wiki/Standard_52-card_deck}{a standard
deck of cards}, and deal them out into $13$ piles of $4$ cards each -- e.g.,
as follows:%
\begin{align*}
&  \left\{  2\spadesuit,\ 2\heartsuit,\ 9\diamondsuit,\ K\diamondsuit\right\}
,\ \ \ \ \ \ \ \ \ \ \left\{  A\spadesuit,\ A\heartsuit,\ 3\spadesuit
,\ 3\diamondsuit\right\}  ,\ \ \ \ \ \ \ \ \ \ \left\{  A\diamondsuit
,\ 4\clubsuit,\ 5\clubsuit,\ Q\clubsuit\right\}  ,\\
&  \left\{  2\diamondsuit,\ 4\heartsuit,\ 5\heartsuit,\ 5\spadesuit\right\}
,\ \ \ \ \ \ \ \ \ \ \left\{  A\clubsuit,\ 7\clubsuit,\ 7\spadesuit
,\ 7\heartsuit\right\}  ,\ \ \ \ \ \ \ \ \ \ \left\{  4\spadesuit
,\ 6\spadesuit,\ 6\diamondsuit,\ 6\clubsuit\right\}  ,\\
&  \left\{  3\heartsuit,\ 3\clubsuit,\ 8\spadesuit,\ 8\heartsuit\right\}
,\ \ \ \ \ \ \ \ \ \ \left\{  2\clubsuit,\ K\clubsuit,\ K\heartsuit
,\ 10\heartsuit\right\}  ,\ \ \ \ \ \ \ \ \ \ \left\{  4\diamondsuit
,\ 5\diamondsuit,\ 9\spadesuit,\ 9\heartsuit\right\}  ,\\
&  \left\{  Q\spadesuit,\ Q\heartsuit,\ Q\diamondsuit,\ Q\clubsuit\right\}
,\ \ \ \ \ \ \ \ \ \ \left\{  6\heartsuit,\ J\spadesuit,\ J\diamondsuit
,\ J\clubsuit\right\}  ,\ \ \ \ \ \ \ \ \ \ \left\{  7\diamondsuit
,\ 8\diamondsuit,\ 8\clubsuit,\ 9\clubsuit\right\}  ,\\
&  \left\{  10\spadesuit,\ J\heartsuit,\ 10\diamondsuit,\ 10\clubsuit\right\}
\end{align*}
(you can distribute the cards among the piles randomly; this is just one
example). Then, I claim that it is possible to select exactly $1$ card from
each pile so that the $13$ selected cards contain exactly $1$ card of each
rank (i.e., exactly one ace, exactly one $2$, exactly one $3$, and so on).

Indeed, this follows from Theorem \ref{thm.match.SDR} (applied to
$A_{i}=\left\{  \text{ranks that occur in the }i\text{-th pile}\right\}  $)
because any $p$ piles contain cards of at least $p$ different ranks.
\end{example}

\begin{proof}
[Proof of Theorem \ref{thm.match.SDR}.]First, we WLOG assume that all $n$ sets
$A_{1},A_{2},\ldots,A_{n}$ are finite. (If not, then we can just replace each
infinite one by an $n$-element subset thereof. The assumption $\left\vert
A_{i_{1}}\cup A_{i_{2}}\cup\cdots\cup A_{i_{p}}\right\vert \geq p$ will not be
disturbed by this change -- make sure you understand why!)

Furthermore, we WLOG assume that no integer belongs to any of the $n$ sets
$A_{1},A_{2},\ldots,A_{n}$ (otherwise, we just rename the elements of these
sets so that they aren't integers any more).

Now, let $X=\left\{  1,2,\ldots,n\right\}  $ and $Y=A_{1}\cup A_{2}\cup
\cdots\cup A_{n}$. Both sets $X$ and $Y$ are finite, and are disjoint.

We define a simple graph $G$ as follows:

\begin{itemize}
\item The vertices of $G$ are the elements of $X\cup Y$.

\item A vertex $x\in X$ is adjacent to a vertex $y\in Y$ if and only if $y\in
A_{x}$. There are no further adjacencies.
\end{itemize}

Thus, $\left(  G,X,Y\right)  $ is a bipartite graph. The assumption
$\left\vert A_{i_{1}}\cup A_{i_{2}}\cup\cdots\cup A_{i_{p}}\right\vert \geq p$
ensures that it satisfies the Hall condition. Hence, by the HMT (Theorem
\ref{thm.match.HMT}), we conclude that this graph $G$ has an $X$-complete
matching. This matching must have the form%
\[
\left\{  \left\{  1,a_{1}\right\}  ,\ \left\{  2,a_{2}\right\}  ,\ \ldots
,\ \left\{  n,a_{n}\right\}  \right\}
\]
for some $a_{1},a_{2},\ldots,a_{n}\in Y$ (since $\left(  G,X,Y\right)  $ is
bipartite, so that the partners of the vertices $1,2,\ldots,n\in X$ must
belong to $Y$). These elements $a_{1},a_{2},\ldots,a_{n}\in Y$ are distinct
(since two edges in a matching cannot have a common endpoint), and each
$i\in\left\{  1,2,\ldots,n\right\}  $ satisfies $a_{i}\in A_{i}$ (since the
vertex $a_{i}$ is adjacent to $i$ in $G$). Thus, these $a_{1},a_{2}%
,\ldots,a_{n}$ are precisely the $n$ distinct elements we are looking for.
This proves Theorem \ref{thm.match.SDR}.
\end{proof}

Conversely, it is not hard to derive the HMT from Theorem \ref{thm.match.SDR}.
Thus, Theorem \ref{thm.match.SDR} is an equivalent version of the HMT. It is
Theorem \ref{thm.match.SDR} that Hall originally discovered (\cite[Theorem
1]{Hall35}).

\bigskip

Here is the second set-theoretical restatement of the HMT:

\begin{theorem}
[existence of SCR]Let $A_{1},A_{2},\ldots,A_{n}$ be $n$ sets. Let $B_{1}%
,B_{2},\ldots,B_{m}$ be $m$ sets. Assume that for any numbers $1\leq
i_{1}<i_{2}<\cdots<i_{p}\leq n$, there exist at least $p$ elements
$j\in\left\{  1,2,\ldots,m\right\}  $ such that the union $A_{i_{1}}\cup
A_{i_{2}}\cup\cdots\cup A_{i_{p}}$ has nonempty intersection with $B_{j}$.
Then, there exists an injective map $\sigma:\left\{  1,2,\ldots,n\right\}
\rightarrow\left\{  1,2,\ldots,m\right\}  $ such that all $i\in\left\{
1,2,\ldots,n\right\}  $ satisfy $A_{i}\cap B_{\sigma\left(  i\right)  }%
\neq\varnothing$.
\end{theorem}

\begin{proof}
We leave this to the reader. Again, construct an appropriate bipartite graph
and apply the HMT.
\end{proof}

(The \textquotedblleft SCR\textquotedblright\ in the name of the theorem is
short for \textquotedblleft system of common representatives\textquotedblright.)

See \cite{MirPer66} for much more about systems of representatives.

\subsection{Regular bipartite graphs}

The HMT gives a necessary and sufficient criterion for the existence of an
$X$-complete matching in an arbitrary bipartite graph. In the more restrictive
setting of \textbf{regular} bipartite graphs -- i.e., bipartite graphs where
each vertex has the same degree --, there is a simpler sufficient condition:
such a matching always exists! We shall soon prove this surprising fact (which
is not hard using the HMT), but first let us get the definition in order:

\begin{definition}
Let $k\in\mathbb{N}$. A multigraph $G$ is said to be $k$\textbf{-regular} if
all its vertices have degree $k$.
\end{definition}

\begin{example}
A $1$-regular graph is a graph whose entire edge set is a perfect matching. In
other words, a $1$-regular graph is a graph that is a disjoint union of copies
of the $2$-nd path graph $P_{2}$. Here is an example of such a graph:%
\[%
%
\]

\end{example}

\begin{example}
A $2$-regular graph is a graph that is a disjoint union of cycle graphs. Here
is an example of such a graph:%
\[%
%
\]
(yes, $C_{1}$ and $C_{2}$ are allowed).
\end{example}

\begin{example}
The $3$-regular graphs are known as \textbf{cubic graphs} or \textbf{trivalent
graphs}. An example is the Petersen graph (defined in Subsection
\ref{subsec.sg.complete.kneser}). Here is another example (known as
\href{https://en.wikipedia.org/wiki/Frucht_graph}{the \textbf{Frucht graph}}):%
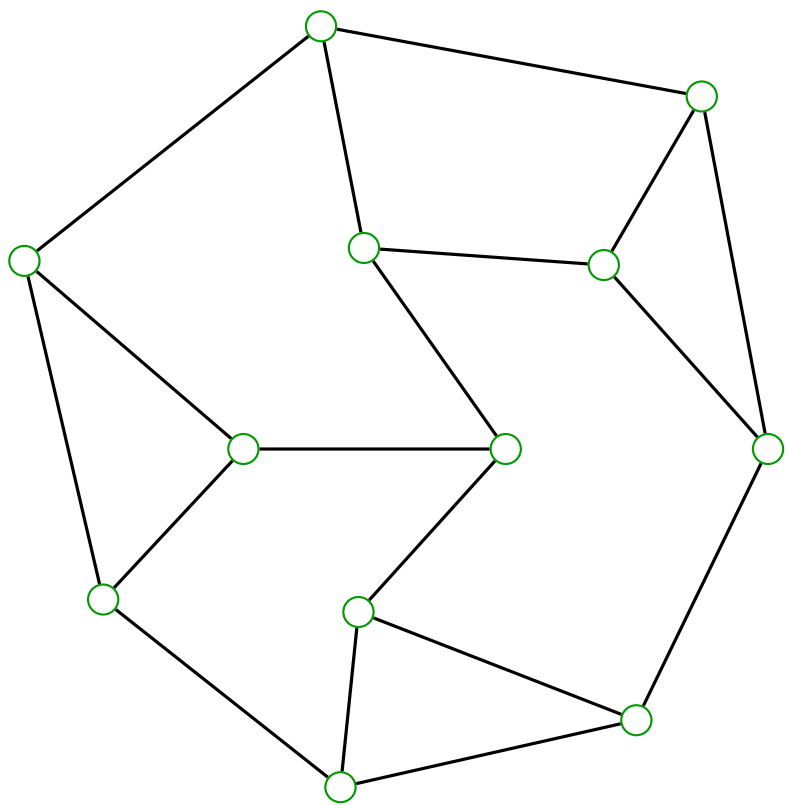
\[%
%
\ \ .
\]
More examples of cubic graphs can be found on
\href{https://en.wikipedia.org/wiki/Cubic_graph}{the Wikipedia page}. There is
no hope of describing them all.
\end{example}

Recall the Kneser graphs defined in Subsection \ref{subsec.sg.complete.kneser}%
. They are all regular:

\begin{example}
Any Kneser graph $K_{S,k}$ is $\dbinom{\left\vert S\right\vert -k}{k}$-regular.
\end{example}

\begin{proof}
This is saying that if $A$ is a $k$-element subset of a finite set $S$, then
there are precisely $\dbinom{\left\vert S\right\vert -k}{k}$ many $k$-element
subsets of $S$ that are disjoint from $A$. But this is clear, since the latter
subsets are just the $k$-element subsets of the $\left(  \left\vert
S\right\vert -k\right)  $-element set $S\setminus A$.
\end{proof}

\begin{proposition}
\label{prop.match.frob-equal}Let $k>0$. Let $\left(  G,X,Y\right)  $ be a
$k$-regular bipartite graph (i.e., a bipartite graph such that $G$ is
$k$-regular). Then, $\left\vert X\right\vert =\left\vert Y\right\vert $.
\end{proposition}

\begin{proof}
Write the multigraph $G$ as $G=\left(  V,E,\varphi\right)  $. Each edge $e\in
E$ contains exactly one vertex $x\in X$ (since $\left(  G,X,Y\right)  $ is a
bipartite graph). Hence,
\begin{align*}
\left\vert E\right\vert  &  =\sum_{x\in X}\underbrace{\left(  \text{\# of
edges that contain the vertex }x\right)  }_{=\deg x}=\sum_{x\in X}%
\underbrace{\deg x}_{\substack{=k\\\text{(since }G\text{ is }k\text{-regular)}%
}}\\
&  =\sum_{x\in X}k=k\cdot\left\vert X\right\vert .
\end{align*}
Similarly, $\left\vert E\right\vert =k\cdot\left\vert Y\right\vert $.
Comparing these two equalities, we obtain $k\cdot\left\vert X\right\vert
=k\cdot\left\vert Y\right\vert $. Since $k>0$, we can divide this by $k$, and
conclude $\left\vert X\right\vert =\left\vert Y\right\vert $.
\end{proof}

\begin{theorem}
[Frobenius matching theorem]\label{thm.match.frob}Let $k>0$. Let $\left(
G,X,Y\right)  $ be a $k$-regular bipartite graph (i.e., a bipartite graph such
that $G$ is $k$-regular). Then, $G$ has a perfect matching.
\end{theorem}

\begin{proof}
First, we claim that each subset $A$ of $X$ satisfies $\left\vert N\left(
A\right)  \right\vert \geq\left\vert A\right\vert $.

Indeed, let $A$ be a subset of $X$. Consider the edges of $G$ that have at
least one endpoint in $A$. We shall call such edges \textquotedblleft%
$A$-edges\textquotedblright. How many $A$-edges are there?

On the one hand, each $A$-edge contains exactly one vertex in $A$
(why?\footnote{Here we are using the fact that $A\subseteq X$, so that no two
vertices in $A$ can be adjacent.}). Thus,%
\begin{align*}
\left(  \text{\# of }A\text{-edges}\right)   &  =\sum_{x\in A}%
\underbrace{\left(  \text{\# of }A\text{-edges containing the vertex
}x\right)  }_{\substack{=\deg x\\\text{(since each edge that contains the
vertex }x\\\text{is an }A\text{-edge)}}}\\
&  =\sum_{x\in A}\underbrace{\deg x}_{\substack{=k\\\text{(since }G\text{ is
}k\text{-regular)}}}=\sum_{x\in A}k=k\cdot\left\vert A\right\vert .
\end{align*}

On the other hand, each $A$-edge contains exactly one vertex in $N\left(
A\right)  $ (why?\footnote{Here we are using the fact that $N\left(  A\right)
\subseteq Y$ (which follows from $A\subseteq X$ using Proposition
\ref{prop.bipartite.NAY}), so that no two vertices in $N\left(  A\right)  $
can be adjacent.}). Thus,%
\begin{align*}
\left(  \text{\# of }A\text{-edges}\right)   &  =\sum_{y\in N\left(  A\right)
}\underbrace{\left(  \text{\# of }A\text{-edges containing the vertex
}y\right)  }_{\substack{\leq\deg y}}\\
&  \leq\sum_{y\in N\left(  A\right)  }\underbrace{\deg y}%
_{\substack{=k\\\text{(since }G\text{ is }k\text{-regular)}}}=\sum_{y\in
N\left(  A\right)  }k=k\cdot\left\vert N\left(  A\right)  \right\vert .
\end{align*}
Hence,%
\[
k\cdot\left\vert N\left(  A\right)  \right\vert \geq\left(  \text{\# of
}A\text{-edges}\right)  =k\cdot\left\vert A\right\vert .
\]
Since $k>0$, we can divide this inequality by $k$, and thus find $\left\vert
N\left(  A\right)  \right\vert \geq\left\vert A\right\vert $.

Forget that we fixed $A$. We thus have proved $\left\vert N\left(  A\right)
\right\vert \geq\left\vert A\right\vert $ for each subset $A$ of $X$. Hence,
the HMT (Theorem \ref{thm.match.HMT}) yields that the graph $G$ has an
$X$-complete matching $M$. Consider this $M$.

However, Proposition \ref{prop.match.frob-equal} yields $\left\vert
X\right\vert =\left\vert Y\right\vert $. Hence, Proposition
\ref{prop.match.bipar-match-basics} \textbf{(f)} shows that the matching $M$
is perfect (since $M$ is $X$-complete). Therefore, $G$ has a perfect matching.
This proves Theorem \ref{thm.match.frob}.
\end{proof}

\begin{fineprint}
For a rather surprising alternative proof of Theorem \ref{thm.match.frob},
using Eulerian circuits instead of Hall's marriage theorem, see example 2 in
the MathOverflow question \url{https://mathoverflow.net/q/271608/} .
\end{fineprint}

\subsection{Latin squares}

One of many applications of Theorem \ref{thm.match.frob} is to the study of
Latin squares. Here is the definition of this concept:

\begin{definition}
\label{def.match.latinsq}Let $n\in\mathbb{N}$. A \textbf{Latin square} of
order $n$ is an $n\times n$-matrix $M$ that satisfies the following conditions:

\begin{enumerate}
\item The entries of $M$ are the numbers $1,2,\ldots,n$, each appearing
exactly $n$ times.

\item In each row of $M$, the entries are distinct.

\item In each column of $M$, the entries are distinct.
\end{enumerate}
\end{definition}

\Needspace{10pc}

\begin{example}
\label{exa.match.latinsq-triv}Here is a Latin square of order $5$:%
\[
\left(
\begin{array}
[c]{ccccc}%
1 & 2 & 3 & 4 & 5\\
2 & 3 & 4 & 5 & 1\\
3 & 4 & 5 & 1 & 2\\
4 & 5 & 1 & 2 & 3\\
5 & 1 & 2 & 3 & 4
\end{array}
\right)  .
\]
Similarly, for each $n\in\mathbb{N}$, the matrix $\left(  c_{i+j-1}\right)
_{1\leq i\leq n,\ 1\leq j\leq n}$, where
\[
c_{k}=%
\begin{cases}
k, & \text{if }k\leq n;\\
k-n, & \text{else,}%
\end{cases}
\]
is a Latin square of order $n$.
\end{example}

A popular example of Latin squares of order $9$ are Sudokus (but they have to
satisfy an additional requirement, concerning certain $3\times3$ subsquares).
See \href{https://en.wikipedia.org/wiki/Latin_square}{the Wikipedia page} and
the book \cite{LayMul98} for much more about Latin squares.

\bigskip

The Latin squares in Example \ref{exa.match.latinsq-triv} are rather boring.
What would be a good algorithm to construct general Latin squares?

Here is an attempt at a recursive algorithm: We just start by filling in the
first row, then the second row, then the third row, and so on, making sure at
each step that the distinctness conditions (Conditions 2 and 3 in Definition
\ref{def.match.latinsq}) are satisfied.

\begin{example}
Let us construct a Latin square of order $5$ by this algorithm. We begin
(e.g.) with the first row%
\[
\left(
\begin{array}
[c]{ccccc}%
3 & 1 & 4 & 2 & 5
\end{array}
\right)  .
\]
Then, we append a second row $\left(
\begin{array}
[c]{ccccc}%
2 & 4 & 1 & 5 & 3
\end{array}
\right)  $ to it, chosen in such a way that its five entries are distinct and
also each entry is distinct from the entry above (again, there are many
possibilities; we have just picked one). Thus, we have our first two rows:%
\[
\left(
\begin{array}
[c]{ccccc}%
3 & 1 & 4 & 2 & 5\\
2 & 4 & 1 & 5 & 3
\end{array}
\right)  .
\]
We continue along the same lines, ending up with the Latin square
\[
\left(
\begin{array}
[c]{ccccc}%
3 & 1 & 4 & 2 & 5\\
2 & 4 & 1 & 5 & 3\\
1 & 5 & 2 & 3 & 4\\
5 & 2 & 3 & 4 & 1\\
4 & 3 & 5 & 1 & 2
\end{array}
\right)
\]
(or another, depending on the choices we have made).
\end{example}

Does this algorithm always work?

To be fully honest, it's not a fully specified algorithm, since I haven't
explained how to fill a row (it's not straightforward). But let's assume that
we know how to do this, if it is at all possible. The natural question is:
Will we always be able to produce a complete Latin square using this
algorithm, or will we get stuck somewhere (having constructed $k$ rows for
some $k<n$, but being unable to produce a $\left(  k+1\right)  $-st row)?

It turns out that we won't get stuck this way. In other words, the following holds:

\begin{proposition}
\label{prop.match.latin}Let $n\in\mathbb{N}$ and $k\in\left\{  0,1,\ldots
,n-1\right\}  $. Then, any $k\times n$ \textbf{Latin rectangle} (i.e., any
$k\times n$-matrix that contains the entries $1,2,\ldots,n$, each appearing
exactly $k$ times, and satisfies the Conditions 2 and 3 from Definition
\ref{def.match.latinsq}) can be extended to a $\left(  k+1\right)  \times n$
Latin rectangle by adding an appropriately chosen extra row at the bottom.
\end{proposition}

\begin{proof}
Let $M$ be a $k\times n$ Latin rectangle\footnote{For example, if $n=5$ and
$k=3$, then $M$ can be $\left(
\begin{array}
[c]{ccccc}%
3 & 1 & 4 & 2 & 5\\
2 & 4 & 1 & 5 & 3\\
1 & 5 & 2 & 3 & 4
\end{array}
\right)  $.}. We want to find a new row that we can append to $M$ at the
bottom, such that the result will be a $\left(  k+1\right)  \times n$ Latin rectangle.

This new row should contain the numbers $1,2,\ldots,n$ in some order.
Moreover, for each $i\in\left\{  1,2,\ldots,n\right\}  $, its $i$-th entry
should be distinct from all entries of the $i$-th column of $M$. How do we
find such a new row?

Let $X=\left\{  1,2,\ldots,n\right\}  $ and $Y=\left\{  -1,-2,\ldots
,-n\right\}  $.

Let $G$ be the simple graph with vertex set $X\cup Y$, where we let a vertex
$i\in X$ be adjacent to a vertex $-j\in Y$ if and only if the number $j$ does
not appear in the $i$-th column of $M$. There should be no further adjacencies.

Thus, $\left(  G,X,Y\right)  $ is a bipartite graph. Moreover, the graph $G$
is $\left(  n-k\right)  $-regular (this is not hard to
see\footnote{\textit{Proof.} Each vertex $i\in X$ has degree $n-k$ (after all,
there are $k$ numbers in $\left\{  1,2,\ldots,n\right\}  $ that appear in the
$i$-th column of $M$, thus $n-k$ numbers in $\left\{  1,2,\ldots,n\right\}  $
that \textbf{do not} appear in this column). It remains to show that each
vertex $-j\in Y$ has degree $n-k$ as well. To see this, consider some vertex
$-j\in Y$. Then, the number $j$ appears exactly once in each row of $M$ (since
Condition 2 forces each row to contain the numbers $1,2,\ldots,n$ in some
order). Hence, the number $j$ appears a total of $k$ times in $M$. These $k$
appearances of $j$ must be in $k$ distinct columns (since having two of them
in the same column would conflict with Condition 3). Thus, there are $k$
columns of $M$ that contain $j$, and therefore $n-k$ columns that don't. In
other words, the vertex $-j\in Y$ has degree $n-k$.}). Thus, by the Frobenius
matching theorem (Theorem \ref{thm.match.frob}), the graph $G$ has a perfect
matching. Let
\[
\left\{  \left\{  1,\ -a_{1}\right\}  ,\ \left\{  2,\ -a_{2}\right\}
,\ \ldots,\ \left\{  n,\ -a_{n}\right\}  \right\}
\]
be this perfect matching. Then, the numbers $a_{1},a_{2},\ldots,a_{n}$ are
distinct (since two edges in a matching cannot have a common endpoint), and
the number $a_{i}$ does not appear in the $i$-th column of $M$ (since
$\left\{  i,\ -a_{i}\right\}  $ is an edge of $G$). Thus, we can append the
row%
\[
\left(
\begin{array}
[c]{cccc}%
a_{1} & a_{2} & \cdots & a_{n}%
\end{array}
\right)
\]
to $M$ at the bottom and obtain a $\left(  k+1\right)  \times n$ Latin
rectangle. This proves Proposition \ref{prop.match.latin}.
\end{proof}

Proposition \ref{prop.match.latin} is a result of Marshall Hall (no relation
to Philip Hall) from 1945 (see \cite{Hall45}), and the proof given above is
exactly his.

\subsection{Magic matrices and the Birkhoff--von Neumann theorem}

Let us now apply the HMT to linear algebra.

Recall that $\mathbb{N}=\left\{  0,1,2,\ldots\right\}  $. We also set
$\mathbb{R}_{+}:=\left\{  \text{all nonnegative reals}\right\}  $.

Here are three very similar definitions:

\begin{definition}
\label{def.match.bvn.Nmag}An $\mathbb{N}$\textbf{-magic matrix} means an
$n\times n$-matrix $M$ that satisfies the following three conditions:

\begin{enumerate}
\item All entries of $M$ are nonnegative integers.

\item The sum of the entries in each row of $M$ is equal.

\item The sum of the entries in each column of $M$ is equal.
\end{enumerate}
\end{definition}

\begin{definition}
\label{def.match.bvn.Rmag}An $\mathbb{R}_{+}$\textbf{-magic matrix} means an
$n\times n$-matrix $M$ that satisfies the following three conditions:

\begin{enumerate}
\item All entries of $M$ are nonnegative reals.

\item The sum of the entries in each row of $M$ is equal.

\item The sum of the entries in each column of $M$ is equal.
\end{enumerate}
\end{definition}

\begin{definition}
\label{def.match.bvn.ds}A \textbf{doubly stochastic matrix} means an $n\times
n$-matrix $M$ that satisfies the following three conditions:

\begin{enumerate}
\item All entries of $M$ are nonnegative reals.

\item The sum of the entries in each row of $M$ is $1$.

\item The sum of the entries in each column of $M$ is $1$.
\end{enumerate}
\end{definition}

Clearly, these three concepts are closely related (in particular, all
$\mathbb{N}$-magic matrices and all doubly stochastic matrices are
$\mathbb{R}_{+}$-magic). The most important of them is the last; in
particular, majorization theory (one of the main methods for proving
inequalities) is deeply connected to the properties of doubly stochastic
matrices (see \cite[Chapter 2]{MaOlAr11}). See \cite[Chapter 2]{BapRag97} for
a chapter-length treatment of doubly stochastic matrices. We shall only prove
some of their most basic properties. First, some examples:

\begin{example}
For any $n>0$, the $n\times n$-matrix%
\[
\left(
\begin{array}
[c]{cccc}%
1 & 1 & \cdots & 1\\
1 & 1 & \cdots & 1\\
\vdots & \vdots & \ddots & \vdots\\
1 & 1 & \cdots & 1
\end{array}
\right)
\]
is $\mathbb{N}$-magic and also $\mathbb{R}_{+}$-magic. This matrix is not
doubly stochastic (unless $n=1$), since the sum of the entries in a row or
column is $n$, not $1$. However, if we divide this matrix by $n$, it becomes
doubly stochastic.
\end{example}

\begin{example}
Here is an $\mathbb{N}$-magic $3\times3$-matrix:%
\[
\left(
\begin{array}
[c]{ccc}%
7 & 0 & 5\\
2 & 6 & 4\\
3 & 6 & 3
\end{array}
\right)  .
\]
Dividing this matrix by $12$ gives a doubly stochastic matrix.
\end{example}

\begin{example}
A \textbf{permutation matrix} is an $n\times n$-matrix whose entries are $0$'s
and $1$'s, and which has exactly one $1$ in each row and exactly one $1$ in
each column. For example, $\left(
\begin{array}
[c]{cccc}%
0 & 0 & 1 & 0\\
1 & 0 & 0 & 0\\
0 & 1 & 0 & 0\\
0 & 0 & 0 & 1
\end{array}
\right)  $ is a permutation matrix of size $4$.

For any $n\in\mathbb{N}$, there are $n!$ many permutation matrices (of size
$n$), since they are in bijection with the permutations of $\left\{
1,2,\ldots,n\right\}  $. Namely, if $\sigma$ is a permutation of $\left\{
1,2,\ldots,n\right\}  $, then the corresponding permutation matrix $P\left(
\sigma\right)  $ has its $\left(  i,\sigma\left(  i\right)  \right)  $-th
entries equal to $1$ for all $i\in\left\{  1,2,\ldots,n\right\}  $, while its
remaining $n^{2}-n$ entries are $0$. For example, if $\sigma$ is the
permutation of $\left\{  1,2,3\right\}  $ sending $1,2,3$ to $2,3,1$, then the
corresponding permutation matrix $P\left(  \sigma\right)  $ is $\left(
\begin{array}
[c]{ccc}%
0 & 1 & 0\\
0 & 0 & 1\\
1 & 0 & 0
\end{array}
\right)  $.

Any permutation matrix is $\mathbb{N}$-magic, $\mathbb{R}_{+}$-magic and
doubly stochastic.
\end{example}

It turns out that these permutation matrices are (in a sense) the
\textquotedblleft building blocks\textquotedblright\ of all magic (and doubly
stochastic) matrices! Namely, the following holds:

\begin{theorem}
[Birkhoff--von Neumann theorem]\label{thm.match.bvn.bvn}Let $n\in\mathbb{N}$. Then:

\begin{enumerate}
\item[\textbf{(a)}] Any $\mathbb{N}$-magic $n\times n$-matrix can be expressed
as a finite sum of permutation matrices.

\item[\textbf{(b)}] Any $\mathbb{R}_{+}$-magic $n\times n$-matrix can be
expressed as an $\mathbb{R}_{+}$-linear combination of permutation matrices
(i.e., it can be expressed in the form $\lambda_{1}P_{1}+\lambda_{2}%
P_{2}+\cdots+\lambda_{k}P_{k}$, where $\lambda_{1},\lambda_{2},\ldots
,\lambda_{k}\in\mathbb{R}_{+}$ are numbers and where $P_{1},P_{2},\ldots
,P_{k}$ are permutation matrices).

\item[\textbf{(c)}] Let $n>0$. Any doubly stochastic $n\times n$-matrix can be
expressed as a \href{https://en.wikipedia.org/wiki/Convex_combination}{convex
combination} of permutation matrices (i.e., it can be expressed in the form
$\lambda_{1}P_{1}+\lambda_{2}P_{2}+\cdots+\lambda_{k}P_{k}$, where
$\lambda_{1},\lambda_{2},\ldots,\lambda_{k}\in\mathbb{R}_{+}$ are numbers
satisfying $\lambda_{1}+\lambda_{2}+\cdots+\lambda_{k}=1$ and where
$P_{1},P_{2},\ldots,P_{k}$ are permutation matrices).
\end{enumerate}
\end{theorem}

Soon we will sketch a proof of this theorem using the HMT. First, two simple
results that will be used in the proof.

\begin{proposition}
\label{prop.match.bvn.samenum}Let $A$ be an $\mathbb{N}$-magic or
$\mathbb{R}_{+}$-magic $n\times n$-matrix. Then, the sum of all entries in a
row of $A$ equals the sum of all entries in a column of $A$.
\end{proposition}

\begin{proof}
Both sums equal $\dfrac{1}{n}$ times the sum of all entries of $A$ (since $A$
has $n$ rows and $n$ columns).
\end{proof}

\begin{lemma}
\label{lem.match.bvn.nonzero}Let $M$ be an $\mathbb{N}$-magic or
$\mathbb{R}_{+}$-magic matrix that is not the zero matrix. Then, there exists
a permutation $\sigma$ of $\left\{  1,2,\ldots,n\right\}  $ such that all
entries $M_{1,\sigma\left(  1\right)  },\ M_{2,\sigma\left(  2\right)
},\ \ldots,\ M_{n,\sigma\left(  n\right)  }$ are nonzero.
\end{lemma}

\begin{example}
If $n=3$ and $M=\left(
\begin{array}
[c]{ccc}%
2 & 7 & 1\\
0 & 1 & 9\\
8 & 2 & 0
\end{array}
\right)  $, then the permutation $\sigma$ that sends $1,2,3$ to $3,2,1$ has
this property.
\end{example}

\begin{proof}
[Proof of Lemma \ref{lem.match.bvn.nonzero}.]Let $s$ denote the sum of the
entries in any given row of $M$ (it doesn't matter which row we take, since
$M$ is magic). Then, $s$ is also the sum of the entries in any given column of
$M$ (by Proposition \ref{prop.match.bvn.samenum}). Also, the sum of all
entries of $M$ is $ns$. Hence, $ns>0$ (since $M$ has nonnegative entries and
is not the zero matrix). Thus, $s>0$.

Let $X=\left\{  1,2,\ldots,n\right\}  $ and $Y=\left\{  -1,-2,\ldots
,-n\right\}  $.

Let $G$ be the simple graph with vertex set $X\cup Y$ and with edges defined
as follows: A vertex $i\in X$ shall be adjacent to a vertex $-j\in Y$ if and
only if $M_{i,j}>0$ (here, $M_{i,j}$ denotes the $\left(  i,j\right)  $-th
entry of $M$). There shall be no further adjacencies.

Thus, $\left(  G,X,Y\right)  $ is a bipartite graph.

We shall now prove that it satisfies the Hall condition. That is, we shall
prove that every subset $A$ of $\left\{  1,2,\ldots,n\right\}  $ satisfies
$\left\vert N\left(  A\right)  \right\vert \geq\left\vert A\right\vert $.

Assume the contrary. Thus, there exists a subset $A$ of $\left\{
1,2,\ldots,n\right\}  $ that satisfies $\left\vert N\left(  A\right)
\right\vert <\left\vert A\right\vert $. Consider this $A$. WLOG assume that
$A=\left\{  1,2,\ldots,k\right\}  $ for some $k\in\left\{  0,1,\ldots
,n\right\}  $ (otherwise, we permute the rows of $M$). Thus, all positive
entries in the first $k$ rows of $A$ are concentrated in fewer than $k$
columns (since the columns in which they lie are the $j$-th columns for $j\in
N\left(  A\right)  $, but we have $\left\vert N\left(  A\right)  \right\vert
<\left\vert A\right\vert =k$). Therefore, the sum of these entries is smaller
than $ks$ (since the sum of all entries in any given column is $s$). On the
other hand, however, the sum of these entries equals $ks$, because they are
all the positive entries in the first $k$ rows of $A$ (and the sum of all
positive entries in a given row equals the sum of \textbf{all} entries in this
row, which is $s$). The two preceding sentences clearly contradict each other.
This contradiction shows that our assumption was false.

Hence, the Hall condition is satisfied. Thus, the HMT yields that $G$ has a
perfect matching. Let
\[
\left\{  \left\{  1,-a_{1}\right\}  ,\ \left\{  2,-a_{2}\right\}
,\ \ldots,\ \left\{  n,-a_{n}\right\}  \right\}
\]
be this perfect matching. Then, $a_{1},a_{2},\ldots,a_{n}$ are distinct, so we
can find a permutation $\sigma$ of $\left\{  1,2,\ldots,n\right\}  $ such that
$a_{i}=\sigma\left(  i\right)  $ for all $i\in\left\{  1,2,\ldots,n\right\}
$. This permutation $\sigma$ then satisfies $M_{i,\sigma\left(  i\right)  }>0$
for all $i\in\left\{  1,2,\ldots,n\right\}  $, which is what we wanted. Thus,
Lemma \ref{lem.match.bvn.nonzero} is proved.
\end{proof}

\begin{proof}
[Proof of Theorem \ref{thm.match.bvn.bvn} (sketched).]\textbf{(a)} Let $M$ be
an $\mathbb{N}$-magic $n\times n$-matrix. How can we express $M$ as a sum of
permutation matrices?

We can try the following method: Try to subtract a permutation matrix from $M$
in such a way that the result will still be an $\mathbb{N}$-magic matrix. Then
do this again, and again and again... until we reach the zero matrix. Once we
have arrived at the zero matrix, the sum of all the permutation matrices that
we have subtracted along the way must be $M$.

Let us experience this method on an example: Let $n=3$ and\footnote{We are
here omitting zero entries from matrices. Thus, $\left(
\begin{array}
[c]{ccc}%
2 & 7 & 1\\
& 1 & 9\\
8 & 2 &
\end{array}
\right)  $ means the matrix $\left(
\begin{array}
[c]{ccc}%
2 & 7 & 1\\
0 & 1 & 9\\
8 & 2 & 0
\end{array}
\right)  $.} $M=\left(
\begin{array}
[c]{ccc}%
2 & 7 & 1\\
& 1 & 9\\
8 & 2 &
\end{array}
\right)  $. If we subtract a permutation matrix from $M$, then the resulting
matrix will still satisfy Conditions 2 and 3 of Definition
\ref{def.match.bvn.Nmag} (since the sum of the entries in any row has been
decreased by $1$, and the sum of the entries in any column has also been
decreased by $1$); however, Condition 1 is not guaranteed, since the
subtraction may turn an entry of $M$ negative (which is not allowed). For
example, this would happen if we tried to subtract the permutation matrix
$\left(
\begin{array}
[c]{ccc}%
1 &  & \\
& 1 & \\
&  & 1
\end{array}
\right)  $ from $M$. Fortunately, Lemma \ref{lem.match.bvn.nonzero} tells us
that there is a permutation $\sigma$ of $\left\{  1,2,\ldots,n\right\}  $ such
that all entries $M_{1,\sigma\left(  1\right)  },\ M_{2,\sigma\left(
2\right)  },\ \ldots,\ M_{n,\sigma\left(  n\right)  }$ are nonzero. If we
choose such a $\sigma$, and subtract the corresponding permutation matrix
$P\left(  \sigma\right)  $ from $M$, then we obtain an $\mathbb{N}$-magic
matrix, because subtracting $1$ from the nonzero entries $M_{1,\sigma\left(
1\right)  },\ M_{2,\sigma\left(  2\right)  },\ \ldots,\ M_{n,\sigma\left(
n\right)  }$ cannot render any of these entries negative. In our example, we
can pick $\sigma$ to be the permutation that sends $1,2,3$ to $3,2,1$. The
corresponding permutation matrix $P\left(  \sigma\right)  $ is $\left(
\begin{array}
[c]{ccc}
&  & 1\\
& 1 & \\
1 &  &
\end{array}
\right)  $. Subtracting this matrix from $M$, we find
\[
\left(
\begin{array}
[c]{ccc}%
2 & 7 & 1\\
& 1 & 9\\
8 & 2 &
\end{array}
\right)  -\left(
\begin{array}
[c]{ccc}
&  & 1\\
& 1 & \\
1 &  &
\end{array}
\right)  =\left(
\begin{array}
[c]{ccc}%
2 & 7 & \\
&  & 9\\
7 & 2 &
\end{array}
\right)  .
\]
This is again an $\mathbb{N}$-magic matrix. Thus, let us do the same to it
that we did to $M$: We again subtract a permutation matrix.

This time, we can actually do better: We can subtract the permutation matrix
$\left(
\begin{array}
[c]{ccc}
& 1 & \\
&  & 1\\
1 &  &
\end{array}
\right)  $ from $\left(
\begin{array}
[c]{ccc}%
2 & 7 & \\
&  & 9\\
7 & 2 &
\end{array}
\right)  $ not just once, but $7$ times, without rendering any entry negative,
because the relevant entries $7,9,7$ are all $\geq7$. The result is%
\[
\left(
\begin{array}
[c]{ccc}%
2 & 7 & \\
&  & 9\\
7 & 2 &
\end{array}
\right)  -7\cdot\left(
\begin{array}
[c]{ccc}
& 1 & \\
&  & 1\\
1 &  &
\end{array}
\right)  =\left(
\begin{array}
[c]{ccc}%
2 &  & \\
&  & 2\\
& 2 &
\end{array}
\right)  .
\]

Now, we follow the same recipe and again subtract a permutation matrix. This
time, we can do it $2$ times, and obtain%
\[
\left(
\begin{array}
[c]{ccc}%
2 &  & \\
&  & 2\\
& 2 &
\end{array}
\right)  -2\cdot\left(
\begin{array}
[c]{ccc}%
1 &  & \\
&  & 1\\
& 1 &
\end{array}
\right)  =\left(
\begin{array}
[c]{ccc}%
\phantom{1} &  & \\
& \phantom{1} & \\
&  & \phantom{1}
\end{array}
\right)  =0_{3\times3}%
\]
(the zero matrix, in case you're wondering).

Thus, we have arrived at the zero matrix by successively subtracting
permutation matrices from $M$. Hence, $M$ is the sum of all the permutation
matrices that have been subtracted: namely,%
\[
M=\left(
\begin{array}
[c]{ccc}
&  & 1\\
& 1 & \\
1 &  &
\end{array}
\right)  +7\cdot\left(
\begin{array}
[c]{ccc}
& 1 & \\
&  & 1\\
1 &  &
\end{array}
\right)  +2\cdot\left(
\begin{array}
[c]{ccc}%
1 &  & \\
&  & 1\\
& 1 &
\end{array}
\right)  ,
\]
which is a sum of $1+7+2$ permutation matrices.

This method works in general, because:

\begin{itemize}
\item If $M$ is an $\mathbb{N}$-magic matrix that is not the zero matrix, then
Lemma \ref{lem.match.bvn.nonzero} tells us that there is a permutation
$\sigma$ of $\left\{  1,2,\ldots,n\right\}  $ such that all entries
$M_{1,\sigma\left(  1\right)  },\ M_{2,\sigma\left(  2\right)  }%
,\ \ldots,\ M_{n,\sigma\left(  n\right)  }$ are nonzero. We can then choose
such a $\sigma$ and subtract the corresponding permutation matrix $P\left(
\sigma\right)  $ from $M$.

\item Better yet, we can subtract $m\cdot P\left(  \sigma\right)  $ from $M$,
where%
\[
m=\min\left\{  M_{1,\sigma\left(  1\right)  },\ M_{2,\sigma\left(  2\right)
},\ \ldots,\ M_{n,\sigma\left(  n\right)  }\right\}  .
\]
This results in an $\mathbb{N}$-magic matrix (since the sum of the entries
decreases by $m$ in each row and by $m$ in each column, and since we are only
subtracting $m$ from a bunch of entries that are $\geq m$) that has at least
one fewer nonzero entry than $M$ (since at least one of the nonzero entries
$M_{1,\sigma\left(  1\right)  },\ M_{2,\sigma\left(  2\right)  }%
,\ \ldots,\ M_{n,\sigma\left(  n\right)  }$ becomes $0$ when $m$ is subtracted
from it).

\item This way, in each step of our process, the number of nonzero entries of
our matrix decreases by at least $1$ (but the matrix remains an $\mathbb{N}%
$-magic matrix throughout the process). Hence, we eventually (after at most
$n^{2}$ steps) will end up with the zero matrix.
\end{itemize}

This proves Theorem \ref{thm.match.bvn.bvn} \textbf{(a)}. \medskip

\textbf{(b)} This is analogous to the proof of part \textbf{(a)} (but this
time, we \textbf{have} to subtract $m\cdot P\left(  \sigma\right)  $ rather
than $P\left(  \sigma\right)  $ in our procedure, since the nonzero entries
$M_{1,\sigma\left(  1\right)  },\ M_{2,\sigma\left(  2\right)  }%
,\ \ldots,\ M_{n,\sigma\left(  n\right)  }$ are not necessarily $\geq1$).
\medskip

\textbf{(c)} Let $M$ be a doubly stochastic $n\times n$-matrix. Then, $M$ is
also $\mathbb{R}_{+}$-magic. Hence, part \textbf{(b)} shows that $M$ can be
expressed in the form $\lambda_{1}P_{1}+\lambda_{2}P_{2}+\cdots+\lambda
_{k}P_{k}$, where $\lambda_{1},\lambda_{2},\ldots,\lambda_{k}\in\mathbb{R}%
_{+}$ are numbers and where $P_{1},P_{2},\ldots,P_{k}$ are permutation
matrices. Consider these $\lambda_{1},\lambda_{2},\ldots,\lambda_{k}$ and
these $P_{1},P_{2},\ldots,P_{k}$.

Now, consider the sum of all entries in the first row of $M$. It is easy to
see that this sum is $\lambda_{1}+\lambda_{2}+\cdots+\lambda_{k}$ (because
$M=\lambda_{1}P_{1}+\lambda_{2}P_{2}+\cdots+\lambda_{k}P_{k}$, but each
permutation matrix $P_{i}$ contributes a $1$ to the sum of all entries in the
first row). But we know that this sum is $1$, since $M$ is doubly stochastic.
Comparing these, we conclude that $\lambda_{1}+\lambda_{2}+\cdots+\lambda
_{k}=1$. Thus, we have expressed $M$ in the form $\lambda_{1}P_{1}+\lambda
_{2}P_{2}+\cdots+\lambda_{k}P_{k}$, where $\lambda_{1},\lambda_{2}%
,\ldots,\lambda_{k}\in\mathbb{R}_{+}$ are numbers satisfying $\lambda
_{1}+\lambda_{2}+\cdots+\lambda_{k}=1$ and where $P_{1},P_{2},\ldots,P_{k}$
are permutation matrices. This proves Theorem \ref{thm.match.bvn.bvn}
\textbf{(c)}.
\end{proof}

\subsection{Further uses of Hall's marriage theorem}

The following few exercises illustrate other applications of Hall's marriage theorem:

\begin{exercise}
\label{exe.9.1}Let $X$ and $Y$ be two finite sets such that $\left\vert
X\right\vert \leq\left\vert Y\right\vert $. Let $f:X\rightarrow Y$ be a map
that is not constant. (A map is said to be \textbf{constant} if all its values
are equal.) Prove that there exists an injective map $g:X\rightarrow Y$ such
that each $x\in X$ satisfies $g\left(  x\right)  \neq f\left(  x\right)  $.
\end{exercise}

\begin{exercise}
\label{exe.8.7}Let $A$ and $B$ be two finite sets such that $\left\vert
B\right\vert \geq\left\vert A\right\vert $. Let $d_{i,j}$ be a real number for
each $\left(  i,j\right)  \in A\times B$. Let
\[
m_{1}=\min\limits_{\sigma:A\rightarrow B\text{ injective}}\ \max\limits_{i\in
A}d_{i,\sigma\left(  i\right)  }%
\]
and
\[
m_{2}=\max\limits_{\substack{I\subseteq A;\ J\subseteq B;\\\left\vert
I\right\vert +\left\vert J\right\vert =\left\vert B\right\vert +1}%
}\ \min\limits_{\left(  i,j\right)  \in I\times J}d_{i,j}.
\]
(The notation \textquotedblleft$\min_{\text{some kind of objects}}\text{some
kind of value}$\textquotedblright\ means the minimum of the given value over
all objects of the given kind. An analogous notation is used for a maximum.)
Prove that $m_{1}=m_{2}$.
\end{exercise}

\begin{exercise}
\label{exe.mt2.verts-to-edges-incident}Let $G=\left(  V,E\right)  $ be a
simple graph such that $\left\vert E\right\vert \geq\left\vert V\right\vert $.
Show that there exists an injective map $f:V\rightarrow E$ such that for each
vertex $v\in V$, the edge $f\left(  v\right)  $ does not contain $v$.

(In other words, show that we can assign to each vertex an edge that does not
contain this vertex in such a way that no edge is assigned twice.) \medskip

[\textbf{Remark:} This is, in some sense, an \textquotedblleft evil
twin\textquotedblright\ to Exercise \ref{exe.5.4}. However, it requires a
simple graph, not a multigraph, since a multigraph with a single vertex and a
single loop would constitute a counterexample. Incidentally, Exercise
\ref{exe.5.4} can also be solved using Hall's marriage theorem.] \medskip

[\textbf{Solution:} This is Exercise 1 on midterm \#2 from my Spring 2017
course; see \href{https://www.cip.ifi.lmu.de/~grinberg/t/17s/}{the course
page} for solutions.]
\end{exercise}

\begin{exercise}
\label{exe.hw4.parabolic-derangements} Let $S$ be a finite set, and let $k
\in\mathbb{N}$. Let $A_{1}, A_{2}, \ldots, A_{k}$ be $k$ subsets of $S$ such
that each element of $S$ lies in exactly one of these $k$ subsets. Prove that
the following statements are equivalent:

\begin{itemize}
\item \textit{Statement 1:} There exists a bijection $\sigma: S \to S$ such
that each $i \in\left\{  1, 2, \ldots, k \right\}  $ satisfies $\sigma\left(
A_{i} \right)  \cap A_{i} = \varnothing$.

\item \textit{Statement 2:} Each $i\in\left\{  1,2,\ldots,k\right\}  $
satisfies $\left\vert A_{i}\right\vert \leq\left\vert S\right\vert /2$.
\end{itemize}

[\textbf{Solution:} This is Exercise 5 on homework set \#4 from my Spring 2017
course; see \href{https://www.cip.ifi.lmu.de/~grinberg/t/17s/}{the course
page} for solutions.]
\end{exercise}

\begin{exercise}
\label{exe.hw4.subsets} Let $S$ be a finite set. Let $k\in\mathbb{N}$ be such
that $\left\vert S\right\vert \geq2k+1$. Prove that there exists an injective
map $f:\mathcal{P}_{k}\left(  S\right)  \rightarrow\mathcal{P}_{k+1}\left(
S\right)  $ such that each $X\in\mathcal{P}_{k}\left(  {S}\right)  $ satisfies
$f\left(  X\right)  \supseteq X$.

(In other words, prove that we can add to each $k$-element subset $X$ of $S$
an additional element from $S\setminus X$ such that the resulting $\left(
k+1\right)  $-element subsets are distinct.)\medskip

[\textbf{Example:} For $S=\left\{  1,2,3,4,5\right\}  $ and $k=2$, we can (for
instance) have the map $f$ send
\begin{align*}
\left\{  1,2\right\}   &  \mapsto\left\{  1,2,4\right\}
,\ \ \ \ \ \ \ \ \ \ \left\{  1,3\right\}  \mapsto\left\{  1,3,4\right\}
,\ \ \ \ \ \ \ \ \ \ \left\{  1,4\right\}  \mapsto\left\{  1,4,5\right\}  ,\\
\left\{  1,5\right\}   &  \mapsto\left\{  1,3,5\right\}
,\ \ \ \ \ \ \ \ \ \ \left\{  2,3\right\}  \mapsto\left\{  1,2,3\right\}
,\ \ \ \ \ \ \ \ \ \ \left\{  2,4\right\}  \mapsto\left\{  2,4,5\right\}  ,\\
\left\{  2,5\right\}   &  \mapsto\left\{  1,2,5\right\}
,\ \ \ \ \ \ \ \ \ \ \left\{  3,4\right\}  \mapsto\left\{  2,3,4\right\}
,\ \ \ \ \ \ \ \ \ \ \left\{  3,5\right\}  \mapsto\left\{  2,3,5\right\}  ,\\
\left\{  4,5\right\}   &  \mapsto\left\{  3,4,5\right\}  .
\end{align*}
Do you see any pattern behind these values? (I don't).]\medskip

[\textbf{Hint:} First, reduce the problem to the case when $\left\vert
S\right\vert =2k+1$. Then, in that case, restate it as a claim about matchings
in a certain bipartite graph.]\medskip

[\textbf{Solution:} This is Exercise 4 on homework set \#4 from my Spring 2017
course; see \href{https://www.cip.ifi.lmu.de/~grinberg/t/17s/}{the course
page} for solutions.]
\end{exercise}

\begin{exercise}
\label{exe.hw4.lattice} Let $\left(  G,X,Y\right)  $ be a bipartite graph.
Assume that each subset $A$ of $X$ satisfies $\left\vert N\left(  A\right)
\right\vert \geq\left\vert A\right\vert $. (Thus, Theorem \ref{thm.match.HMT}
shows that $G$ has an $X$-complete matching.)

A subset $A$ of $X$ will be called \textbf{neighbor-critical} if $\left\vert
N\left(  A\right)  \right\vert =\left\vert A\right\vert $.

Let $A$ and $B$ be two neighbor-critical subsets of $X$. Prove that the
subsets $A\cup B$ and $A\cap B$ are also neighbor-critical.\medskip

[\textbf{Solution:} This is Exercise 6 on homework set \#4 from my Spring 2017
course; see \href{https://www.cip.ifi.lmu.de/~grinberg/t/17s/}{the course
page} for solutions.]
\end{exercise}

For more about applications, restatements and variants of Hall's marriage
theorem, see the recent survey \cite{Camero25}.

\subsection{Further exercises on matchings}

\begin{exercise}
\label{exe.8.5}Let $G=\left(  V,E,\varphi\right)  $ be a multigraph. Let $M$
be a matching of $G$.

An \textbf{augmenting path} for $M$ shall mean a path $\left(  v_{0}%
,e_{1},v_{1},e_{2},v_{2},\ldots,e_{k},v_{k}\right)  $ of $G$ such that $k$ is
odd (note that $k=1$ is allowed) and such that

\begin{itemize}
\item the even-indexed edges $e_{2}, e_{4}, \ldots, e_{k-1}$ belong to $M$
(note that this condition is vacuously true if $k = 1$);

\item the odd-indexed edges $e_{1}, e_{3}, \ldots, e_{k}$ belong to $E
\setminus M$;

\item neither the starting point $v_{0}$ nor the ending point $v_{k}$ is
matched in $M$.
\end{itemize}

Prove that $M$ has maximum size among all matchings of $G$ if and only if
there exists no augmenting path for $M$. \medskip

[\textbf{Hint:} If $M$ and $M^{\prime}$ are two matchings of $G$, what can you
say about the symmetric difference $\left(  M \cup M^{\prime}\right)
\setminus\left(  M \cap M^{\prime}\right)  $ ?]
\end{exercise}

\begin{exercise}
\label{exe.9.2}Let $\left(  G,X,Y\right)  $ be a bipartite graph. Let $A$ be a
subset of $X$, and let $B$ be a subset of $Y$. Assume that $G$ has an
$A$-complete matching, and that $G$ has a $B$-complete matching. Prove that
$G$ has an $A\cup B$-complete matching.
\end{exercise}

\begin{exercise}
\label{exe.9.3}Let $\left(  G,X,Y\right)  $ be a bipartite graph with
$X\neq\varnothing$. Assume that $G$ has an $X$-complete matching.

An edge $e$ of $G$ will be called \textbf{useless} if $G$ has no $X$-complete
matching that contains $e$.

Prove that there exists a vertex $x\in X$ such that no edge that contains $x$
is useless.
\end{exercise}

\begin{exercise}
\label{exe.mt3.lopsided-bipartite} Let $\left(  G,X,Y\right)  $ be a bipartite
graph such that $\left\vert Y\right\vert \geq2\left\vert X\right\vert -1$.
Prove that there exists an injective map $f:X\rightarrow Y$ such that each
$x\in X$ satisfies one of the following two statements:

\begin{itemize}
\item \textit{Statement 1:} The vertices $x$ and $f \left(  x \right)  $ of
$G$ are adjacent.

\item \textit{Statement 2:} There exists no $x^{\prime}\in X$ such that the
vertices $x$ and $f\left(  x^{\prime}\right)  $ of $G$ are adjacent.\medskip
\end{itemize}

[\textbf{Remark:} In vaguely matrimonial terminology, this is saying that in a
group of $m$ men and $w$ women satisfying $w\geq2m-1$, we can always marry
each man (monogamously) to a woman in such a way that either he likes his
partner or all women he likes are unmarried.]\medskip

[\textbf{Solution:} This is Exercise 3 on midterm \#3 from my Spring 2017
course; see \href{https://www.cip.ifi.lmu.de/~grinberg/t/17s/}{the course
page} for solutions.]
\end{exercise}

\begin{exercise}
\label{exe.hw4.matching-product} Let $\left(  G,X,Y\right)  $ and $\left(
H,U,V\right)  $ be bipartite graphs.

Assume that $G$ is a simple graph and has an $X$-complete matching.

Assume that $H$ is a simple graph and has a $U$-complete matching.

Consider the Cartesian product $G\times H$ of $G$ and $H$ defined in
Definition \ref{def.sg.cartprod}. (Note that we required $G$ and $H$ to be
simple graphs only to avoid having to define $G\times H$ for multigraphs.)

\begin{enumerate}
\item[\textbf{(a)}] Show that $\left(  G\times H,\ \left(  X\times V\right)
\cup\left(  Y\times U\right)  ,\ \left(  X\times U\right)  \cup\left(  Y\times
V\right)  \right)  $ is a bipartite graph.

\item[\textbf{(b)}] Prove that the graph $G\times H$ has an $\left(  X\times
V\right)  \cup\left(  Y\times U\right)  $-complete matching.
\end{enumerate}

[\textbf{Solution:} This is Exercise 3 on homework set \#4 from my Spring 2017
course; see \href{https://www.cip.ifi.lmu.de/~grinberg/t/17s/}{the course
page} for solutions.]
\end{exercise}

\section{\label{chp.flows}Networks and flows}

In this chapter, I will give an introduction to \textbf{network flows} and
their optimization. This is a topic of great interest to logisticians, as even
the simplest results have obvious applications to scheduling trains and
trucks. It also has lots of purely mathematical consequences; in particular,
we will use network flows to finally prove the Hall--K\"{o}nig matching
theorem (and thus the HMT, K\"{o}nig's theorem, and their many consequences).

I will follow my notes \cite{17s-lec16}, which are a good place to look up the
details of some proofs that I will only sketch. That said, I will be using
multidigraphs instead of simple digraphs, so some adaptations will be
necessary (since \cite{17s-lec16} only works with simple digraphs). These
adaptations are generally easy.

I will only cover the very basics of network flow optimization, leading to a
proof of the max-flow-min-cut theorem (for integer-valued flows) and to a
proof of the Hall--K\"{o}nig matching theorem. For the deeper reaches of the
theory, see \cite{ForFul74} (a classical textbook written by the inventors of
the subject), \cite[Chapter 4]{Schrij-ACO} and \cite[Part I]{Schrij-CO1}.

\subsection{\label{sec.flows.def}Definitions}

\subsubsection{Networks}

Recall that we use the notation $\mathbb{N}=\left\{  0,1,2,\ldots\right\}  $.

\begin{definition}
\label{def.net.network}A \textbf{network} consists of

\begin{itemize}
\item a multidigraph $D=\left(  V,A,\psi\right)  $;

\item two distinct vertices $s\in V$ and $t\in V$, called the \textbf{source}
and the \textbf{sink}, respectively;

\item a function $c:A\rightarrow\mathbb{N}$, called the \textbf{capacity
function}.
\end{itemize}
\end{definition}

\begin{example}
\label{exa.net.network.1}Here is an example of a network:%
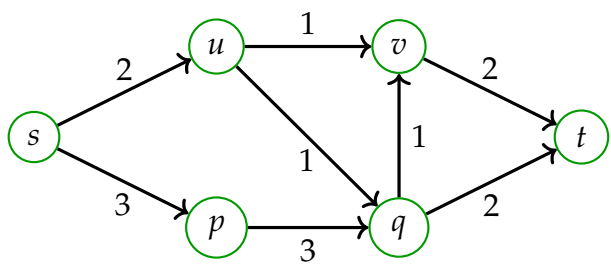
\[%
\begin{tikzpicture}[scale=1.2]
\begin{scope}[every node/.style={circle,thick,draw=green!60!black}]
\node(s) at (-1,0) {$s$};
\node(u) at (1,1) {$u$};
\node(p) at (1,-1) {$p$};
\node(v) at (3,1) {$v$};
\node(q) at (3,-1) {$q$};
\node(t) at (5,0) {$t$};
\end{scope}
\begin{scope}[every edge/.style={draw=black,very thick}, every loop/.style={}]
\path[->] (s) edge node[above] {$2$} (u) edge node[below] {$3$} (p);
\path[->] (p) edge node[below] {$3$} (q);
\path[->] (u) edge node[above] {$1$} (v) edge node[below] {$1$} (q);
\path[->] (q) edge node[right] {$1$} (v) edge node[below] {$2$} (t);
\path[->] (v) edge node[above] {$2$} (t);
\end{scope}
\end{tikzpicture}%
\ \ .
\]
Here, the multidigraph $D$ is the one we drew (it is a simple digraph, so we
have not labeled its arcs); the vertices $s$ and $t$ are the vertices labeled
$s$ and $t$; the values of the function $c$ on the arcs of $D$ are written on
top of these respective arcs (e.g., we have $c\left(  \left(  s,p\right)
\right)  =3$ and $c\left(  \left(  u,q\right)  \right)  =1$).
\end{example}

\begin{remark}
The digraph $D$ in Example \ref{exa.net.network.1} has no cycles and satisfies
$\deg^{-}s=\deg^{+}t=0$. This is not required in the definition of a network,
although it is satisfied in many basic applications.

Also, all capacities $c\left(  a\right)  $ in Example \ref{exa.net.network.1}
were positive. This, too, is not required; however, arcs with capacity $0$ do
not contribute anything useful to the situation, so they could just as well be absent.
\end{remark}

\begin{remark}
The notion of \textquotedblleft network\textquotedblright\ we just defined is
just one of a myriad notions of \textquotedblleft network\textquotedblright%
\ that can be found all over mathematics. Most of them can be regarded as
graphs with \textquotedblleft some extra structures\textquotedblright; apart
from this, they don't have much in common.
\end{remark}

\subsubsection{The notations $\overline{S}$, $\left[  P,Q\right]  $ and
$d\left(  P,Q\right)  $}

\begin{definition}
\label{def.net.capacity}Let $N$ be a network consisting of a multidigraph
$D=\left(  V,A,\psi\right)  $, a source $s\in V$, a sink $t\in V$ and a
capacity function $c:A\rightarrow\mathbb{N}$. Then:

\begin{enumerate}
\item[\textbf{(a)}] For any arc $a\in A$, we call the number $c\left(
a\right)  \in\mathbb{N}$ the \textbf{capacity} of the arc $a$.

\item[\textbf{(b)}] For any subset $S$ of $V$, we let $\overline{S}$ denote
the subset $V\setminus S$ of $V$.

\item[\textbf{(c)}] If $P$ and $Q$ are two subsets of $V$, then $\left[
P,Q\right]  $ shall mean the set of all arcs of $D$ whose source belongs to
$P$ and whose target belongs to $Q$. That is,%
\[
\left[  P,Q\right]  :=\left\{  a\in A\ \mid\ \psi\left(  a\right)  \in P\times
Q\right\}  .
\]

\item[\textbf{(d)}] If $P$ and $Q$ are two subsets of $V$, and if
$d:A\rightarrow\mathbb{N}$ is any function, then the number $d\left(
P,Q\right)  \in\mathbb{N}$ is defined by%
\[
d\left(  P,Q\right)  :=\sum_{a\in\left[  P,Q\right]  }d\left(  a\right)  .
\]
(In particular, we can apply this to $d=c$, and then get $c\left(  P,Q\right)
=\sum\limits_{a\in\left[  P,Q\right]  }c\left(  a\right)  $.)
\end{enumerate}
\end{definition}

\begin{example}
Let us again consider the network from Example \ref{exa.net.network.1}. For
the subset $\left\{  s,u\right\}  $ of $V$, we have $\overline{\left\{
s,u\right\}  }=\left\{  p,v,q,t\right\}  $ and%
\[
\left[  \left\{  s,u\right\}  ,\overline{\left\{  s,u\right\}  }\right]
=\left\{  sp,\ uv,\ uq\right\}
\]
(recall that our $D$ is a simple digraph, so an arc is just a pair of two
vertices) and%
\begin{align*}
c\left(  \left\{  s,u\right\}  ,\overline{\left\{  s,u\right\}  }\right)   &
=\sum_{a\in\left[  \left\{  s,u\right\}  ,\overline{\left\{  s,u\right\}
}\right]  }c\left(  a\right)  =\underbrace{c\left(  sp\right)  }%
_{=3}+\underbrace{c\left(  uv\right)  }_{=1}+\underbrace{c\left(  uq\right)
}_{=1}\\
&  =3+1+1=5.
\end{align*}
We can make this visually clearer if we draw a \textquotedblleft
border\textquotedblright\ between the sets $\left\{  s,u\right\}  $ and
$\overline{\left\{  s,u\right\}  }$:%
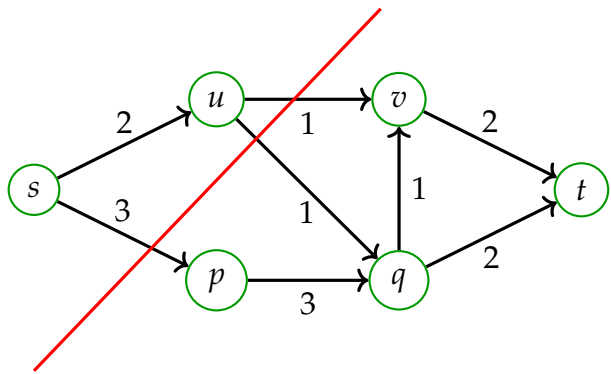
\[%
\begin{tikzpicture}[scale=1.2]
\begin{scope}[every node/.style={circle,thick,draw=green!60!black}]
\node(s) at (-1,0) {$s$};
\node(u) at (1,1) {$u$};
\node(p) at (1,-1) {$p$};
\node(v) at (3,1) {$v$};
\node(q) at (3,-1) {$q$};
\node(t) at (5,0) {$t$};
\end{scope}
\begin{scope}[every edge/.style={draw=black,very thick}, every loop/.style={}]
\path[->] (s) edge node[above] {$2$} (u) edge node[above] {$3$} (p);
\path[->] (p) edge node[below] {$3$} (q);
\path[->] (u) edge node[below] {$1$} (v) edge node[below] {$1$} (q);
\path[->] (q) edge node[right] {$1$} (v) edge node[below] {$2$} (t);
\path[->] (v) edge node[above] {$2$} (t);
\end{scope}
\begin{scope}[every edge/.style={draw=red, very thick}]
\path[-] (2.8,2) edge (-1, -2);
\end{scope}
\end{tikzpicture}%
\ \ .
\]
Then, $\left[  \left\{  s,u\right\}  ,\overline{\left\{  s,u\right\}
}\right]  $ is the set of all arcs that cross this border from $\left\{
s,u\right\}  $ to $\overline{\left\{  s,u\right\}  }$. (Of course, this
visualization works only for sums of the form $d\left(  P,\overline{P}\right)
$, not for the more general case of $d\left(  P,Q\right)  $ where $P$ and $Q$
can have elements in common. But the $d\left(  P,\overline{P}\right)  $ are
the most useful sums.)
\end{example}

\begin{exercise}
\label{exe.9.5}Let $D=\left(  V,A,\psi\right)  $ be a balanced multidigraph.
For any subset $S$ of $V$, we set $\overline{S}:=V\setminus S$. For any two
subsets $S$ and $T$ of $V$, we set
\begin{align*}
\left[  S,T\right]   &  :=\left\{  a\in A\ \mid\ \text{the source of }a\text{
belongs to }S\text{,}\right. \\
&  \ \ \ \ \ \ \ \ \ \ \ \ \ \ \ \ \ \ \ \ \left.  \text{and the target of
}a\text{ belongs to }T\right\}  .
\end{align*}
Prove that $\left\vert \left[  S,\overline{S}\right]  \right\vert =\left\vert
\left[  \overline{S},S\right]  \right\vert $ for any subset $S$ of $V$.
\end{exercise}

\subsubsection{Flows}

Let us now define flows on a network:

\begin{definition}
\label{def.net.flow}Let $N$ be a network consisting of a multidigraph
$D=\left(  V,A,\psi\right)  $, a source $s\in V$, a sink $t\in V$ and a
capacity function $c:A\rightarrow\mathbb{N}$.

A \textbf{flow} (on the network $N$) means a function $f:A\rightarrow
\mathbb{N}$ with the following properties:

\begin{itemize}
\item We have $0\leq f\left(  a\right)  \leq c\left(  a\right)  $ for each arc
$a\in A$. This condition is called the \textbf{capacity constraints} (we are
using the plural form, since there is one constraint for each arc $a\in A$).

\item For any vertex $v\in V\setminus\left\{  s,t\right\}  $, we have%
\[
f^{-}\left(  v\right)  =f^{+}\left(  v\right)  ,
\]
where we set%
\[
f^{-}\left(  v\right)  :=\sum_{\substack{a\in A\text{ is an arc}\\\text{with
target }v}}f\left(  a\right)  \ \ \ \ \ \ \ \ \ \ \text{and}%
\ \ \ \ \ \ \ \ \ \ f^{+}\left(  v\right)  :=\sum_{\substack{a\in A\text{ is
an arc}\\\text{with source }v}}f\left(  a\right)  .
\]
This is called the \textbf{conservation constraints}.
\end{itemize}

If $f:A\rightarrow\mathbb{N}$ is a flow and $a\in A$ is an arc, then the
nonnegative integer $f\left(  a\right)  $ will be called the \textbf{arc flow}
of $f$ on $a$.
\end{definition}

\begin{example}
\label{exa.net.flow.1}To draw a flow $f$ on a network $N$, we draw the network
$N$, with one little tweak: Instead of writing the capacity $c\left(
a\right)  $ atop each arc $a\in A$, we write \textquotedblleft$f\left(
a\right)  $ of $c\left(  a\right)  $\textquotedblright\ atop each arc $a\in
A$. For example, here is a flow $f$ on the network $N$ from Example
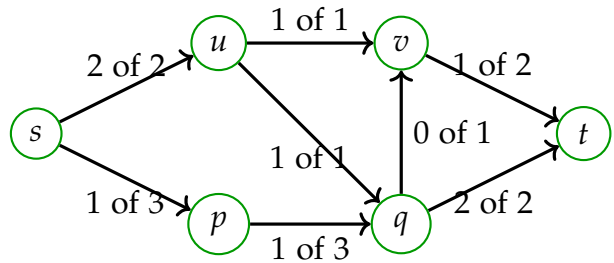
\ref{exa.net.network.1}:%
\[%
%
\]
(so, for example, $f\left(  su\right)  =2$, $f\left(  pq\right)  =1$ and
$f\left(  qv\right)  =0$).

For another example, here is a different flow $g$ on the same network $N$:%
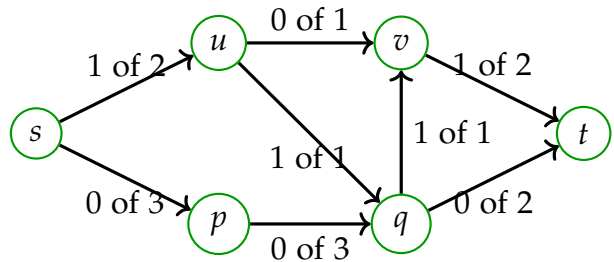
\[%
%
\]

\end{example}

There are several intuitive ways to think of a network $N$ and of a flow on it:

\begin{itemize}
\item We can visualize $N$ as a collection of one-way roads: Each arc $a\in A$
is a one-way road, and its capacity $c\left(  a\right)  $ is how much traffic
it can (maximally) handle per hour. A flow $f$ on $N$ can then be understood
as traffic flowing through these roads, where $f\left(  a\right)  $ is the
amount of traffic that travels through the arc $a$ in an hour. The
conservation constraints say that the traffic out of a given vertex $v$ equals
the traffic into $v$ unless $v$ is one of $s$ and $t$. (We imagine that
traffic can arbitrarily materialize or dematerialize at $s$ and $t$.)

\item We can visualize $N$ as a collection of pipes: Each arc $a\in A$ is a
pipe, and its capacity $c\left(  a\right)  $ is how much water it can
maximally transport in a second. A flow $f$ on $N$ can then be viewed as water
flowing through the pipes, where $f\left(  a\right)  $ is the amount of water
traveling through a pipe $a$ in a second. The capacity constraints say that no
pipe is over its capacity or carries a negative amount of water. The
conservation constraints say that at every vertex $v$ other than $s$ and $t$,
the amount of water coming in (that is, $f^{-}\left(  v\right)  $) equals the
amount of water moving out (that is, $f^{+}\left(  v\right)  $); that is,
there are no leaks and no water being injected into the system other than at
$s$ and $t$. This is why $s$ is called the \textquotedblleft
source\textquotedblright\ and $t$ is the \textquotedblleft
sink\textquotedblright. A slightly counterintuitive aspect of this
visualization is that each pipe has a direction, and water can only flow in
that one direction (from source to target). That said, you can always model an
undirected pipe by having two pipes of opposite directions.

\item We can regard $N$ as a money transfer scheme: Each vertex $v\in V$ is a
bank account, and the goal is to transfer some money from $s$ to $t$. All
other vertices $v$ act as middlemen. Each arc $a\in A$ corresponds to a
possibility of transfer from its source to its target; the maximum amount that
can be transferred on this arc is $c\left(  a\right)  $. A flow describes a
way in which money is transferred such that each middleman vertex $v\in
V\setminus\left\{  s,t\right\}  $ ends up receiving exactly as much money as
it gives away.
\end{itemize}

Needless to say, these visualizations have been chosen for their intuitive
grasp; the real-life applications of network flows are somewhat different.

\begin{remark}
\label{rmk.net.pathflow}Flows on a network $N$ can be viewed as a
generalization of paths on the underlying digraph $D$. Indeed, if $\mathbf{p}$
is a path from $s$ to $t$ on the digraph $D=\left(  V,A,\psi\right)  $
underlying a network $N$, then we can define a flow $f_{\mathbf{p}}$ on $N$ as
follows:%
\[
f_{\mathbf{p}}\left(  a\right)  =%
\begin{cases}
1, & \text{if }a\text{ is an arc of }\mathbf{p};\\
0, & \text{otherwise}%
\end{cases}
\ \ \ \ \ \ \ \ \ \ \text{for each }a\in A,
\]
provided that all arcs of $\mathbf{p}$ have capacity $\geq1$. An example of
such a flow is the flow $g$ in Example \ref{exa.net.flow.1}.
\end{remark}

\subsubsection{Inflow, outflow and value of a flow}

Next, we define certain numbers related to any flow on a network:

\begin{definition}
\label{def.net.value}Let $N$ be a network consisting of a multidigraph
$D=\left(  V,A,\psi\right)  $, a source $s\in V$, a sink $t\in V$ and a
capacity function $c:A\rightarrow\mathbb{N}$. Let $f:A\rightarrow\mathbb{N}$
be an arbitrary map (e.g., a flow on $N$). Then:

\begin{enumerate}
\item[\textbf{(a)}] For each vertex $v\in V$, we set
\[
f^{-}\left(  v\right)  :=\sum_{\substack{a\in A\text{ is an arc}\\\text{with
target }v}}f\left(  a\right)  \ \ \ \ \ \ \ \ \ \ \text{and}%
\ \ \ \ \ \ \ \ \ \ f^{+}\left(  v\right)  :=\sum_{\substack{a\in A\text{ is
an arc}\\\text{with source }v}}f\left(  a\right)  .
\]
We call $f^{-}\left(  v\right)  $ the \textbf{inflow} of $f$ into $v$, and we
call $f^{+}\left(  v\right)  $ the \textbf{outflow} of $f$ from $v$.

Note that if $f$ is a flow, then each vertex $v\in V\setminus\left\{
s,t\right\}  $ satisfies $f^{-}\left(  v\right)  =f^{+}\left(  v\right)  $;
this is the conservation constraints.

\item[\textbf{(b)}] We define the \textbf{value} of the map $f$ to be the
number $f^{+}\left(  s\right)  -f^{-}\left(  s\right)  $. This value is
denoted by $\left\vert f\right\vert $.
\end{enumerate}
\end{definition}

\begin{example}
The flow $f$ in Example \ref{exa.net.flow.1} satisfies%
\begin{align*}
f^{+}\left(  s\right)   &  =3,\ \ \ \ \ \ \ \ \ \ f^{-}\left(  s\right)  =0,\\
f^{+}\left(  u\right)   &  =f^{-}\left(  u\right)  =2,\\
f^{+}\left(  p\right)   &  =f^{-}\left(  p\right)  =1,\\
f^{+}\left(  v\right)   &  =f^{-}\left(  v\right)  =1,\\
f^{+}\left(  q\right)   &  =f^{-}\left(  q\right)  =2,,\\
f^{+}\left(  t\right)   &  =0,\ \ \ \ \ \ \ \ \ \ f^{-}\left(  t\right)  =3
\end{align*}
and has value $\left\vert f\right\vert =3$. The flow $g$ in Example
\ref{exa.net.flow.1} has value $\left\vert g\right\vert =1$. More generally,
the flow $f_{\mathbf{p}}$ in Remark \ref{rmk.net.pathflow} always has value
$\left\vert f_{\mathbf{p}}\right\vert =1$.
\end{example}

\begin{example}
\label{exa.net.zeroflow}For any network $N$, we can define the \textbf{zero
flow} on $N$. This is the flow $0_{A}:A\rightarrow\mathbb{N}$ that sends each
arc $a\in A$ to $0$. This flow has value $\left\vert 0_{A}\right\vert =0$.
\end{example}

\subsection{The maximum flow problem and bipartite graphs}

Now we can state an important optimization problem, known as the
\textbf{maximum flow problem}: Given a network $N$, how can we find a flow of
maximum possible value?

\begin{example}
\label{exa.net.bip-as-network}Finding a maximum matching in a bipartite graph
is a particular case of the maximum flow problem.

Indeed, let $\left(  G,X,Y\right)  $ be a bipartite graph. Then, we can
transform this graph into a network $N$ as follows:

\begin{itemize}
\item Add two new vertices $s$ and $t$.

\item Turn each edge $e$ of $G$ into an arc $\overrightarrow{e}$ whose source
is the $X$-endpoint of $e$ (that is, the endpoint of $e$ that belongs to $X$)
and whose target is the $Y$-endpoint of $e$ (that is, the endpoint of $e$ that
belongs to $Y$).

\item Add an arc from $s$ to each vertex in $X$.

\item Add an arc from each vertex in $Y$ to $t$.

\item Assign to each arc the capacity $1$.
\end{itemize}

Here is an example of a bipartite graph $\left(  G,X,Y\right)  $ (as usual,
drawn with the $X$-vertices on the left and with the $Y$-vertices on the
right) and the corresponding network $N$:
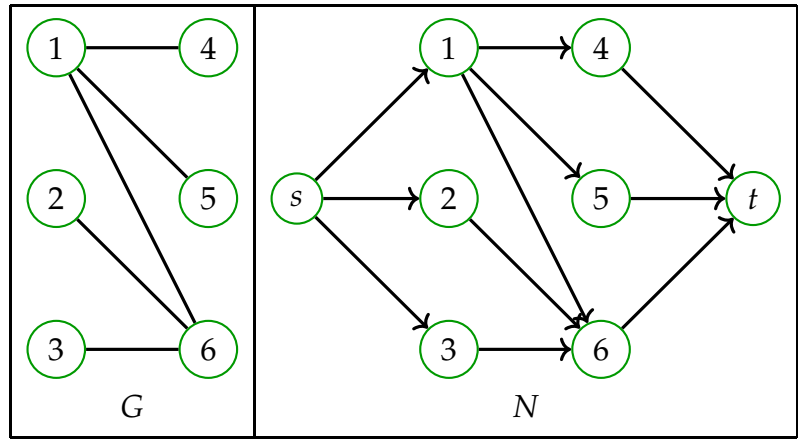
\[%

\
\]
(we are not showing the capacities of the arcs, since they are all equal to
$1$).

The flows of this network $N$ are in bijection with the matchings of $G$.
Namely, if $f$ is a flow on $N$, then the set%
\[
\left\{  e\in\operatorname*{E}\left(  G\right)  \ \mid\ f\left(
\overrightarrow{e}\right)  =1\right\}
\]
is a matching of $G$. Conversely, if $M$ is a matching of $G$, then we obtain
a flow $f$ on $N$ by assigning the arc flow $1$ to all arcs of the form
$\overrightarrow{e}$ where $e\in M$, as well as assigning the arc flow $1$ to
every new arc that joins $s$ or $t$ to a vertex matched in $M$. All other arcs
are assigned the arc flow $0$. For instance, in our above example, the
matching $\left\{  15,\ 36\right\}  $ corresponds to the following flow:%
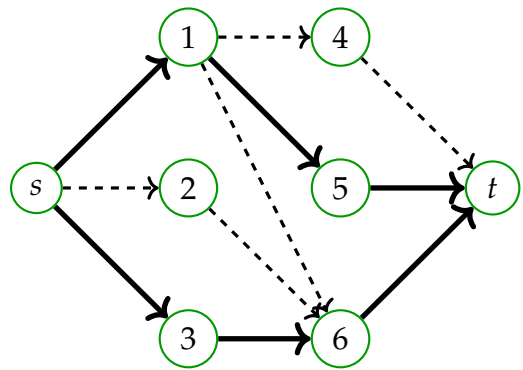
\[%
\begin{tikzpicture}[scale=2]
\begin{scope}[every node/.style={circle,thick,draw=green!60!black}]
\node(1) at (0, 1) {$1$};
\node(2) at (0, 0) {$2$};
\node(3) at (0, -1) {$3$};
\node(4) at (1, 1) {$4$};
\node(5) at (1, 0) {$5$};
\node(6) at (1, -1) {$6$};
\node(s) at (-1, 0) {$s$};
\node(t) at (2, 0) {$t$};
\end{scope}
\begin{scope}[every edge/.style={draw=black,dashed,very thick}]
\path[->] (1) edge (4) edge (6);
\path[->] (2) edge (6);
\path[->] (s) edge (2);
\path[->] (4) edge (t);
\end{scope}
\begin{scope}[every edge/.style={draw=black,line width=2pt}]
\path[->] (1) edge (5);
\path[->] (3) edge (6);
\path[->] (s) edge (1) edge (3);
\path[->] (5) edge (t) (6) edge (t);
\end{scope}
\end{tikzpicture}%
\ \ ,
\]
where we are using the convention that an arc $a$ with $f\left(  a\right)  =0$
is drawn dashed whereas an arc $a$ with $f\left(  a\right)  =1$ is drawn
boldfaced (thankfully, the only possibilities for $f\left(  a\right)  $ are
$0$ and $1$, because all capacities are $1$).

One nice property of this bijection is that if a flow $f$ corresponds to a
matching $M$, then $\left\vert f\right\vert =\left\vert M\right\vert $. Thus,
finding a flow of maximum value means finding a matching of maximum size.

(See \cite[Proposition 1.36 till Proposition 1.40]{17s-lec16} for details and
proofs; that said, the proofs are straightforward and you will probably
\textquotedblleft see\textquotedblright\ them just by starting at an example.)
\end{example}

\subsection{\label{sec.flows.basics}Basic properties of flows}

Before we approach the maximum flow problem, let us prove some simple
observations about flows:

\begin{proposition}
\label{prop.net.flow-t}Let $N$ be a network consisting of a multidigraph
$D=\left(  V,A,\psi\right)  $, a source $s\in V$, a sink $t\in V$ and a
capacity function $c:A\rightarrow\mathbb{N}$. Let $f:A\rightarrow\mathbb{N}$
be a flow on $N$. Then,%
\begin{align*}
\left\vert f\right\vert  &  =f^{+}\left(  s\right)  -f^{-}\left(  s\right) \\
&  =f^{-}\left(  t\right)  -f^{+}\left(  t\right)  .
\end{align*}

\end{proposition}

\begin{proof}
Each $v\in V$ satisfies $f^{+}\left(  v\right)  =\sum_{\substack{a\in A\text{
is an arc}\\\text{with source }v}}f\left(  a\right)  $ (by the definition of
$f^{+}\left(  v\right)  $). Summing this equality over all $v\in V$, we obtain%
\[
\sum\limits_{v\in V}f^{+}\left(  v\right)  =\underbrace{\sum\limits_{v\in
V}\ \ \sum_{\substack{a\in A\text{ is an arc}\\\text{with source }v}}}%
_{=\sum_{a\in A}}f\left(  a\right)  =\sum_{a\in A}f\left(  a\right)
\]
(note that this is a generalization of the familiar fact that $\sum
\limits_{v\in V}\deg^{+}v=\left\vert A\right\vert $). Similarly,
$\sum\limits_{v\in V}f^{-}\left(  v\right)  =\sum\limits_{a\in A}f\left(
a\right)  $. Hence,
\begin{equation}
\sum\limits_{v\in V}\left(  f^{-}\left(  v\right)  -f^{+}\left(  v\right)
\right)  =\underbrace{\sum\limits_{v\in V}f^{-}\left(  v\right)  }%
_{=\sum\limits_{a\in A}f\left(  a\right)  }-\underbrace{\sum_{v\in V}%
f^{+}\left(  v\right)  }_{=\sum\limits_{a\in A}f\left(  a\right)  }=0.
\label{pf.prop.net.flow-t.1}%
\end{equation}

However, by the conservation constraints, we have $f^{-}\left(  v\right)
=f^{+}\left(  v\right)  $ for each $v\in V\setminus\left\{  s,t\right\}  $. In
other words, $f^{-}\left(  v\right)  -f^{+}\left(  v\right)  =0$ for each
$v\in V\setminus\left\{  s,t\right\}  $. Thus, in the sum $\sum\limits_{v\in
V}\left(  f^{-}\left(  v\right)  -f^{+}\left(  v\right)  \right)  $, all
addends are $0$ except for the addends for $v=s$ and for $v=t$. Hence, the sum
boils down to these two addends:%
\[
\sum\limits_{v\in V}\left(  f^{-}\left(  v\right)  -f^{+}\left(  v\right)
\right)  =\left(  f^{-}\left(  s\right)  -f^{+}\left(  s\right)  \right)
+\left(  f^{-}\left(  t\right)  -f^{+}\left(  t\right)  \right)  .
\]
Comparing this with (\ref{pf.prop.net.flow-t.1}), we obtain%
\[
\left(  f^{-}\left(  s\right)  -f^{+}\left(  s\right)  \right)  +\left(
f^{-}\left(  t\right)  -f^{+}\left(  t\right)  \right)  =0,
\]
so that%
\[
f^{-}\left(  t\right)  -f^{+}\left(  t\right)  =-\left(  f^{-}\left(
s\right)  -f^{+}\left(  s\right)  \right)  =f^{+}\left(  s\right)
-f^{-}\left(  s\right)  =\left\vert f\right\vert
\]
(by the definition of $\left\vert f\right\vert $). This proves Proposition
\ref{prop.net.flow-t}.
\end{proof}

\begin{proposition}
\label{prop.net.flow-setS}Let $N$ be a network consisting of a multidigraph
$D=\left(  V,A,\psi\right)  $, a source $s\in V$, a sink $t\in V$ and a
capacity function $c:A\rightarrow\mathbb{N}$. Let $f:A\rightarrow\mathbb{N}$
be a flow on $N$. Let $S$ be a subset of $V$. Then:

\begin{enumerate}
\item[\textbf{(a)}] We have%
\[
f\left(  S,\overline{S}\right)  -f\left(  \overline{S},S\right)  =\sum_{v\in
S}\left(  f^{+}\left(  v\right)  -f^{-}\left(  v\right)  \right)  .
\]
(Recall that we are using Definition \ref{def.net.capacity} here, so that
$f\left(  P,Q\right)  $ means $\sum_{a\in\left[  P,Q\right]  }f\left(
a\right)  $.)

\item[\textbf{(b)}] Assume that $s\in S$ and $t\notin S$. Then,%
\[
\left\vert f\right\vert =f\left(  S,\overline{S}\right)  -f\left(
\overline{S},S\right)  .
\]

\item[\textbf{(c)}] Assume that $s\in S$ and $t\notin S$. Then,%
\[
\left\vert f\right\vert \leq c\left(  S,\overline{S}\right)  .
\]

\item[\textbf{(d)}] Assume that $s\in S$ and $t\notin S$. Then, $\left\vert
f\right\vert =c\left(  S,\overline{S}\right)  $ if and only if%
\[
\left(  f\left(  a\right)  =0\text{ for all }a\in\left[  \overline
{S},S\right]  \right)
\]
and%
\[
\left(  f\left(  a\right)  =c\left(  a\right)  \text{ for all }a\in\left[
S,\overline{S}\right]  \right)  .
\]

\end{enumerate}
\end{proposition}

\begin{proof}
Let me first make these claims intuitive in terms of the \textquotedblleft
money transfer scheme\textquotedblright\ model for our network. Consider $S$
as a country. Then, $f\left(  S,\overline{S}\right)  $ is the
\textquotedblleft export\textquotedblright\ from this country $S$ (that is,
the total wealth exported from $S$), whereas $f\left(  \overline{S},S\right)
$ is the \textquotedblleft import\textquotedblright\ into this country $S$
(that is, the total wealth imported into $S$). Thus, part \textbf{(a)} of the
proposition is saying that the \textquotedblleft net export\textquotedblright%
\ of $S$ (that is, the export from $S$ minus the import into $S$) can be
computed by summing the \textquotedblleft outflow minus
inflow\textquotedblright\ values of all accounts in $S$. This should match the
intuition for exports and imports (in particularly, any transfers that happen
within $S$ should cancel out when we sum the \textquotedblleft outflow minus
inflow\textquotedblright\ values of all accounts in $S$). Part \textbf{(b)}
says that if the country $S$ contains the source $s$ but not the sink $t$
(that is, the goal of the network is to transfer money out of the country),
then the total value transferred is actually the net export of $S$. Part
\textbf{(c)} claims that this total value is no larger than the total
\textquotedblleft export capacity\textquotedblright\ $c\left(  S,\overline
{S}\right)  $ (that is, the total capacity of the \textquotedblleft export
arcs\textquotedblright\ $a\in\left[  S,\overline{S}\right]  $). Part
\textbf{(d)} says that if equality holds in this inequality (i.e., if the
total value equals the total export capacity), then each \textquotedblleft
import arc\textquotedblright\ $a\in\left[  \overline{S},S\right]  $ is unused
(i.e., nothing is imported into $S$), whereas each \textquotedblleft export
arc\textquotedblright\ $a\in\left[  S,\overline{S}\right]  $ is used to its
full capacity. \medskip

I hope this demystifies all claims of the proposition. But for the sake of
completeness, here are rigorous proofs (though rather terse ones, since I
assume you have seen enough manipulations of sum to fill in the details):
\medskip

\begin{fineprint}
\textbf{(a)} This follows from
\begin{align*}
\sum_{v\in S}\left(  f^{+}\left(  v\right)  -f^{-}\left(  v\right)  \right)
&  =\underbrace{\sum_{v\in S}f^{+}\left(  v\right)  }_{\substack{=f\left(
S,V\right)  \\\text{(why?)}}}-\underbrace{\sum_{v\in S}f^{-}\left(  v\right)
}_{\substack{=f\left(  V,S\right)  \\\text{(why?)}}}\\
&  =\underbrace{f\left(  S,V\right)  }_{\substack{=f\left(  S,S\right)
+f\left(  S,\overline{S}\right)  \\\text{(since }V\text{ is the union of
the}\\\text{two disjoint sets }S\text{ and }\overline{S}\text{)}%
}}-\underbrace{f\left(  V,S\right)  }_{\substack{=f\left(  S,S\right)
+f\left(  \overline{S},S\right)  \\\text{(since }V\text{ is the union of
the}\\\text{two disjoint sets }S\text{ and }\overline{S}\text{)}}}\\
&  =f\left(  S,S\right)  +f\left(  S,\overline{S}\right)  -\left(  f\left(
S,S\right)  +f\left(  \overline{S},S\right)  \right)  =f\left(  S,\overline
{S}\right)  -f\left(  \overline{S},S\right)  .
\end{align*}

\textbf{(b)} We have $S\setminus\left\{  s\right\}  \subseteq V\setminus
\left\{  s,t\right\}  $ (since $t\notin S$). From part \textbf{(a)}, we obtain%
\begin{align*}
f\left(  S,\overline{S}\right)  -f\left(  \overline{S},S\right)   &
=\sum_{v\in S}\left(  f^{+}\left(  v\right)  -f^{-}\left(  v\right)  \right)
\\
&  =\underbrace{\left(  f^{+}\left(  s\right)  -f^{-}\left(  s\right)
\right)  }_{\substack{=\left\vert f\right\vert \\\text{(by the definition of
}\left\vert f\right\vert \text{)}}}+\sum_{v\in S\setminus\left\{  s\right\}
}\underbrace{\left(  f^{+}\left(  v\right)  -f^{-}\left(  v\right)  \right)
}_{\substack{=0\\\text{(by the conservation constraints,}\\\text{since }v\in
S\setminus\left\{  s\right\}  \subseteq V\setminus\left\{  s,t\right\}
\text{)}}}\ \ \ \ \ \ \ \ \ \ \left(  \text{since }s\in S\right) \\
&  =\left\vert f\right\vert .
\end{align*}
This proves part \textbf{(b)}. \medskip

\textbf{(c)} The capacity constraints yield that $f\left(  a\right)  \leq
c\left(  a\right)  $ for each arc $a\in A$. Summing up these inequalities over
all $a\in\left[  S,\overline{S}\right]  $, we obtain $f\left(  S,\overline
{S}\right)  \leq c\left(  S,\overline{S}\right)  $. The capacity constraints
furthermore yield that $f\left(  a\right)  \geq0$ for each arc $a\in A$.
Summing up these inequalities over all $a\in\left[  \overline{S},S\right]  $,
we obtain $f\left(  \overline{S},S\right)  \geq0$. Hence, part \textbf{(b)}
yields%
\[
\left\vert f\right\vert =\underbrace{f\left(  S,\overline{S}\right)  }_{\leq
c\left(  S,\overline{S}\right)  }-\underbrace{f\left(  \overline{S},S\right)
}_{\geq0}\leq c\left(  S,\overline{S}\right)  .
\]
This proves part \textbf{(c)}. \medskip

\textbf{(d)} We must characterize the equality case in part \textbf{(c)}.
However, recall the proof of part \textbf{(c)}: We obtained the inequality
$\left\vert f\right\vert \leq c\left(  S,\overline{S}\right)  $ by summing up
the inequalities $f\left(  a\right)  \leq c\left(  a\right)  $ over all arcs
$a\in\left[  S,\overline{S}\right]  $ and subtracting the sum of the
inequalities $f\left(  a\right)  \geq0$ over all arcs $a\in\left[
\overline{S},S\right]  $. Hence, in order for the inequality $\left\vert
f\right\vert \leq c\left(  S,\overline{S}\right)  $ to become an equality, it
is necessary and sufficient that all the inequalities involved -- i.e., the
inequalities $f\left(  a\right)  \leq c\left(  a\right)  $ for all arcs
$a\in\left[  S,\overline{S}\right]  $ as well as the inequalities $f\left(
a\right)  \geq0$ for all arcs $a\in\left[  \overline{S},S\right]  $ -- become
equalities. In other words, it is necessary and sufficient that we have%
\[
\left(  f\left(  a\right)  =0\text{ for all }a\in\left[  \overline
{S},S\right]  \right)
\]
and%
\[
\left(  f\left(  a\right)  =c\left(  a\right)  \text{ for all }a\in\left[
S,\overline{S}\right]  \right)  .
\]
This proves Proposition \ref{prop.net.flow-setS} \textbf{(d)}.
\end{fineprint}
\end{proof}

\subsection{\label{sec.flows.mfmc}The max-flow-min-cut theorem}

\subsubsection{Cuts and their capacities}

One more definition, before we get to the hero of this story:

\begin{definition}
\label{def.net.cut}Let $N$ be a network consisting of a multidigraph
$D=\left(  V,A,\psi\right)  $, a source $s\in V$, a sink $t\in V$ and a
capacity function $c:A\rightarrow\mathbb{N}$. Then:

\begin{enumerate}
\item[\textbf{(a)}] A \textbf{cut} of $N$ shall mean a subset of $A$ that has
the form $\left[  S,\overline{S}\right]  $, where $S$ is a subset of $V$
satisfying $s\in S$ and $t\notin S$.

\item[\textbf{(b)}] The \textbf{capacity} of a cut $\left[  S,\overline
{S}\right]  $ is defined to be the number $c\left(  S,\overline{S}\right)
=\sum\limits_{a\in\left[  S,\overline{S}\right]  }c\left(  a\right)  $.
\end{enumerate}
\end{definition}

\begin{example}
Let us again consider the network from Example \ref{exa.net.network.1}. Then,
$\left[  \left\{  s,u\right\}  ,\overline{\left\{  s,u\right\}  }\right]
=\left\{  sp,\ uv,\ uq\right\}  $ is a cut of this network, and its capacity
is $c\left(  \left\{  s,u\right\}  ,\overline{\left\{  s,u\right\}  }\right)
=5$.
\end{example}

\subsubsection{The max-flow-min-cut theorem: statement}

Now, Proposition \ref{prop.net.flow-setS} \textbf{(c)} says that the value of
any flow $f$ can never be larger than the capacity of any cut $\left[
S,\overline{S}\right]  $. Thus, in particular, the maximum value of a flow is
$\leq$ to the minimum capacity of a cut.

Furthermore, Proposition \ref{prop.net.flow-setS} \textbf{(d)} says that if
this inequality is an equality -- i.e., if the value of some flow $f$ equals
the capacity of some cut $\left[  S,\overline{S}\right]  $ --, then the flow
$f$ must use each arc that crosses the cut in the right direction (from $S$ to
$\overline{S}$) to its full capacity and must not use any of the arcs that
cross the cut in the wrong direction (from $\overline{S}$ to $S$).

It turns out that this inequality actually \textbf{is} an equality for any
maximum flow and any minimum cut:

\begin{theorem}
[max-flow-min-cut theorem]\label{thm.net.MFMC}Let $N$ be a network consisting
of a multidigraph $D=\left(  V,A,\psi\right)  $, a source $s\in V$, a sink
$t\in V$ and a capacity function $c:A\rightarrow\mathbb{N}$. Then,%
\[
\max\left\{  \left\vert f\right\vert \ \mid\ f\text{ is a flow}\right\}
=\min\left\{  c\left(  S,\overline{S}\right)  \ \mid\ S\subseteq V;\ s\in
S;\ t\notin S\right\}  .
\]
In other words, the maximum value of a flow equals the minimum capacity of a cut.
\end{theorem}

We shall soon sketch a proof of this theorem that doubles as a fairly
efficient (polynomial-time) algorithm for finding both a maximum flow (i.e., a
flow of maximum value) and a minimum cut (i.e., a cut of minimum capacity).
The algorithm is known as the \textbf{Ford-Fulkerson algorithm}, and is
sufficiently fast to be useful in practice.

\subsubsection{How to augment a flow}

The idea of this algorithm is to start by having $f$ be the zero flow\ (i.e.,
the flow from Example \ref{exa.net.zeroflow}), and then gradually increase its
value $\left\vert f\right\vert $ by making changes to some of its arc flows
$f\left(  a\right)  $.

Of course, we cannot unilaterally change the arc flow $f\left(  a\right)  $ on
a single arc, since this will (usually) mess up the conservation constraints.
Thus, if we change $f\left(  a\right)  $, then we will also have to change
$f\left(  b\right)  $ for some other arcs $b\in A$ to make the result a flow
again. One way to do this is to increase all arc flows $f\left(  a\right)  $
along some path from $s$ to $t$. Here is an example of such an increase:

\begin{example}
Consider the flow $f$ from Example \ref{exa.net.flow.1}. We can increase the
arc flows $f\left(  sp\right)  ,\ f\left(  pq\right)  ,\ f\left(  qv\right)
,\ f\left(  vt\right)  $ of $f$ on all the arcs of the path $\left(
s,p,q,v,t\right)  $ (since neither of these arcs is used to its full
capacity). As a result, we obtain the following flow $h$:%
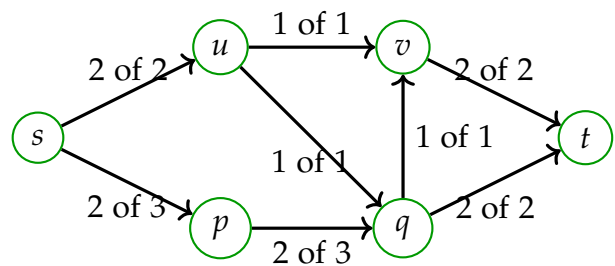
\[%
\begin{tikzpicture}[scale=1.2]
\begin{scope}[every node/.style={circle,thick,draw=green!60!black}]
\node(s) at (-1,0) {$s$};
\node(u) at (1,1) {$u$};
\node(p) at (1,-1) {$p$};
\node(v) at (3,1) {$v$};
\node(q) at (3,-1) {$q$};
\node(t) at (5,0) {$t$};
\end{scope}
\begin{scope}[every edge/.style={draw=black,very thick}, every loop/.style={}]
\path[->] (s) edge node[left=8pt, above] {$2 \of2$}
(u) edge node[left=7pt, below] {$2 \of3$} (p);
\path[->] (p) edge node[below] {$2 \of3$} (q);
\path[->] (u) edge node[above] {$1 \of1$} (v) edge node[pos=0.6, left] {$1 \of
1$} (q);
\path[->] (q) edge node[right] {$1 \of1$}
(v) edge node[below=3pt, pos=0.6] {$2 \of2$} (t);
\path[->] (v) edge node[above=5pt, pos=0.6] {$2 \of2$} (t);
\end{scope}
\end{tikzpicture}%
\ \ ,
\]
whose value $\left\vert h\right\vert $ is $4$. It is easy to see that this is
actually the maximum value of a flow on our network (since $\left\vert
h\right\vert =4$ equals the capacity $c\left(  \overline{\left\{  t\right\}
},\left\{  t\right\}  \right)  $ of the cut $\left[  \overline{\left\{
t\right\}  },\left\{  t\right\}  \right]  $, but Proposition
\ref{prop.net.flow-setS} \textbf{(c)} tells us that the value of any flow is
$\leq$ to the capacity of any cut).
\end{example}

However, simple increases like the one we just did are not always enough to
find a maximum flow. They can leave us stuck at a \textquotedblleft local
maximum\textquotedblright\ -- i.e., at a flow which does not have any more
paths from $s$ to $t$ that can be used for any further increases (i.e., any
path from $s$ to $t$ contains an arc that is already used to its full
capacity), yet is not a maximum flow. Here is an example:

\begin{example}
\label{exa.net.ff.stuck}Consider the following network and flow:%
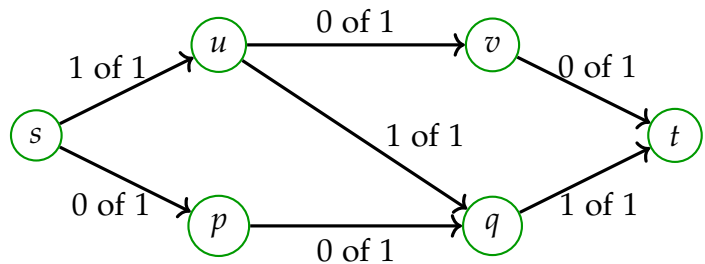
\[%
\begin{tikzpicture}[scale=1.2]
\begin{scope}[every node/.style={circle,thick,draw=green!60!black}]
\node(s) at (-1,0) {$s$};
\node(u) at (1,1) {$u$};
\node(p) at (1,-1) {$p$};
\node(v) at (4,1) {$v$};
\node(q) at (4,-1) {$q$};
\node(t) at (6,0) {$t$};
\end{scope}
\begin{scope}[every edge/.style={draw=black,very thick}, every loop/.style={}]
\path[->] (s) edge node[above] {$1 \of1\ \ \ \ $} (u) edge node[below] {$0 \of
1\ \ \ $} (p);
\path[->] (p) edge node[below] {$0 \of1$} (q);
\path[->] (u) edge node[above] {$0 \of1$} (v) edge node[right] {$\ \ 1 \of1$}
(q);
\path[->] (q) edge node[below] {$\ \ \ 1 \of1$} (t);
\path[->] (v) edge node[above] {$\ \ \ 0 \of1$} (t);
\end{scope}
\end{tikzpicture}%
\ \ .
\]
This flow is not maximum, but each path from $s$ to $t$ has at least one arc
that is used to its full capacity. Thus, we cannot improve this flow by
increasing all its arc flows on any given path from $s$ to $t$.
\end{example}

The trick to get past this hurdle is to use a \textquotedblleft zig-zag
path\textquotedblright\ -- i.e., not a literal path, but rather a sequence
$\left(  v_{0},a_{1},v_{1},a_{2},v_{2},\ldots,a_{k},v_{k}\right)  $ of
vertices and arcs that can use arcs both in the forward and backward
directions (i.e., any $i\in\left\{  1,2,\ldots,k\right\}  $ has to satisfy
either $\psi\left(  a_{i}\right)  =\left(  v_{i-1},v_{i}\right)  $ or
$\psi\left(  a_{i}\right)  =\left(  v_{i},v_{i-1}\right)  $). Instead of
increasing the flow on all arcs of this \textquotedblleft
path\textquotedblright, we do something slightly subtler: On the forward arcs,
we increase the flow; on the backward arcs, we decrease it (all by the same
amount). This, too, preserves the conservation constraints (think about why;
we will soon see a rigorous proof), so it is a valid way of increasing the
value of a flow. Here is an example:

\begin{example}
\label{exa.net.ff.unstuck}Consider the flow in Example \ref{exa.net.ff.stuck}.
The underlying digraph has a \textquotedblleft zig-zag path\textquotedblright%
\ $\left(  s,p,q,u,v,t\right)  $, which uses the arc $uq$ in the backward
direction. We can decrease the arc flows of $f$ on all forward arcs $sp$,
$pq$, $uv$ and $vt$ of this \textquotedblleft zig-zag path\textquotedblright,
and decrease it on the backward arc $uq$. As a result, we obtain the flow%
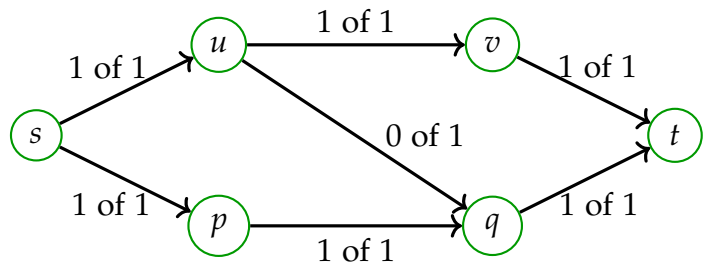
\[%
\begin{tikzpicture}[scale=1.2]
\begin{scope}[every node/.style={circle,thick,draw=green!60!black}]
\node(s) at (-1,0) {$s$};
\node(u) at (1,1) {$u$};
\node(p) at (1,-1) {$p$};
\node(v) at (4,1) {$v$};
\node(q) at (4,-1) {$q$};
\node(t) at (6,0) {$t$};
\end{scope}
\begin{scope}[every edge/.style={draw=black,very thick}, every loop/.style={}]
\path[->] (s) edge node[above] {$1 \of1\ \ \ \ $} (u) edge node[below] {$1 \of
1\ \ \ $} (p);
\path[->] (p) edge node[below] {$1 \of1$} (q);
\path[->] (u) edge node[above] {$1 \of1$} (v) edge node[right] {$\ \ 0 \of1$}
(q);
\path[->] (q) edge node[below] {$\ \ \ 1 \of1$} (t);
\path[->] (v) edge node[above] {$\ \ \ 1 \of1$} (t);
\end{scope}
\end{tikzpicture}%
\ \ .
\]
This new flow has value $2$, and can easily be seen to be a maximum flow.
\end{example}

Good news: Allowing ourselves to use \textquotedblleft zig-zag
paths\textquotedblright\ like this (rather than literal paths only), we never
get stuck at a non-maximum flow; we can always increase the value further and
further until we eventually arrive at a maximum flow.

In order to prove this, we introduce some convenient notations. We prefer not
to talk about \textquotedblleft zig-zag paths\textquotedblright, but rather
reinterpret these \textquotedblleft zig-zag paths\textquotedblright\ as
(literal) paths of an appropriately chosen digraph (not $D$). This has the
advantage of allowing us to use known properties of paths without having to
first generalize them to \textquotedblleft zig-zag paths\textquotedblright.

\subsubsection{The residual digraph}

The appropriately chosen digraph is the so-called \textbf{residual digraph} of
a flow; it is defined as follows:

\begin{definition}
\label{def.net.Df}Let $N$ be a network consisting of a multidigraph $D=\left(
V,A,\psi\right)  $, a source $s\in V$, a sink $t\in V$ and a capacity function
$c:A\rightarrow\mathbb{N}$.

\begin{enumerate}
\item[\textbf{(a)}] For each arc $a\in A$, we introduce a new arc $a^{-1}$,
which should act like a reversal of the arc $a$ (that is, its source should be
the target of $a$, and its target should be the source of $a$). We don't add
these new arcs $a^{-1}$ to our digraph $D$, but we keep them ready for use in
a different digraph (which we will define below).

Here is what this means in rigorous terms: For each arc $a\in A$, we introduce
a new object, which we call $a^{-1}$. We let $A^{-1}$ be the set of these new
objects $a^{-1}$ for $a\in A$. We extend the map $\psi:A\rightarrow V\times V$
to a map $\widehat{\psi}:A\cup A^{-1}\rightarrow V\times V$ as follows: For
each $a\in A$, we let
\[
\widehat{\psi}\left(  a\right)  =\left(  u,v\right)
\ \ \ \ \ \ \ \ \ \ \text{and}\ \ \ \ \ \ \ \ \ \ \widehat{\psi}\left(
a^{-1}\right)  =\left(  v,u\right)  ,
\]
where $u$ and $v$ are defined by $\left(  u,v\right)  =\psi\left(  a\right)  $.

For each arc $a\in A$, we shall refer to the new arc $a^{-1}$ as the
\textbf{reversal} of $a$, and conversely, we shall refer to the original arc
$a$ as the \textbf{reversal} of $a^{-1}$. We set $\left(  a^{-1}\right)
^{-1}:=a$ for each $a\in A$.

We shall refer to the arcs $a\in A$ as \textbf{forward arcs}, and to their
reversals $a^{-1}$ as \textbf{backward arcs}.

\item[\textbf{(b)}] Let $f:A\rightarrow\mathbb{N}$ be any flow on $N$. We
define the \textbf{residual digraph }$D_{f}$ of this flow $f$ to be the
multidigraph $\left(  V,A_{f},\psi_{f}\right)  $, where%
\[
A_{f}=\left\{  a\in A\ \mid\ f\left(  a\right)  <c\left(  a\right)  \right\}
\cup\left\{  a^{-1}\ \mid\ a\in A;\ f\left(  a\right)  >0\right\}
\]
and $\psi_{f}:=\widehat{\psi}\mid_{A_{f}}$. (This is usually not a subdigraph
of $D$.) Thus, the residual digraph $D_{f}$ has the same vertices as $V$, but
its arcs are those arcs of $D$ that are not used to their full capacity by $f$
as well as the reversals of all arcs of $D$ that are used by $f$.
\end{enumerate}
\end{definition}

\begin{example}
\label{exa.net.Df.1}Let $f$ be the flow $f$ from Example \ref{exa.net.flow.1}.
Then, the residual digraph $D_{f}$ is%
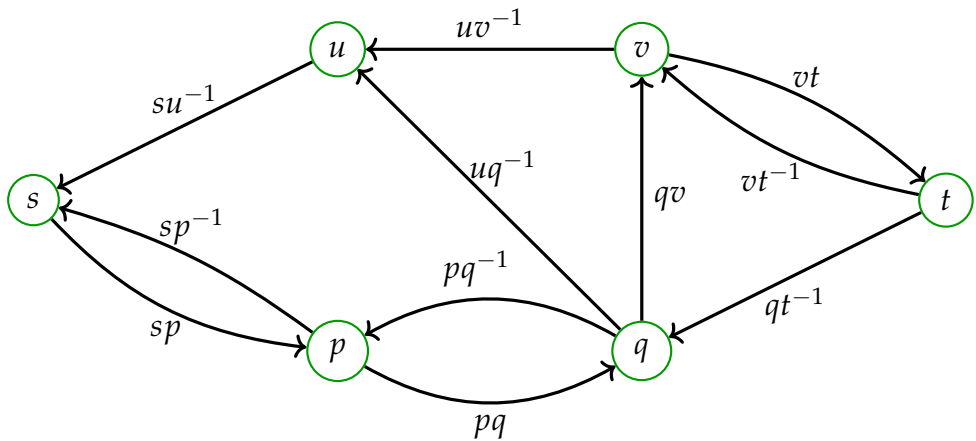
\[%
%
\ \ .
\]
Notice that the digraph $D_{f}$ has cycles even though $D$ has none!
\end{example}

\begin{example}
Let $f$ be the non-maximum flow from Example \ref{exa.net.ff.stuck}. Then, the
residual digraph $D_{f}$ is%
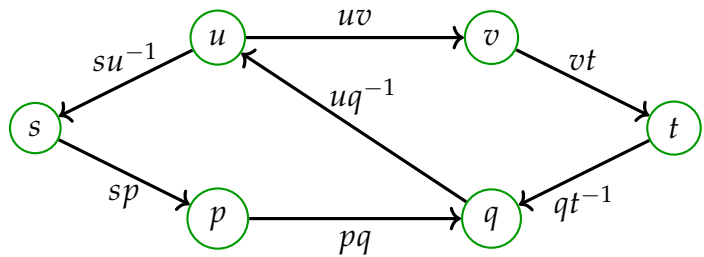
\[%
%
\ \ .
\]
This digraph $D_{f}$ has a path from $s$ to $t$, which corresponds precisely
to the \textquotedblleft zig-zag path\textquotedblright\ $\left(
s,p,q,u,v,t\right)  $ we found in Example \ref{exa.net.ff.unstuck}.
\end{example}

You can think of the residual digraph $D_{f}$ as follows: Each arc of $D_{f}$
corresponds to an opportunity to change an arc flow $f\left(  a\right)  $;
namely, a forward arc $a$ of $D_{f}$ means that $f\left(  a\right)  $ can be
increased, whereas a backward arc $a^{-1}$ of $D_{f}$ means that $f\left(
a\right)  $ can be decreased. Hence, the paths of the residual digraph $D_{f}$
are the \textquotedblleft zig-zag paths\textquotedblright\ of $D$ that allow
the flow $f$ to be increased (on forward arcs) or decreased (on backward arcs)
as in Example \ref{exa.net.ff.unstuck}. Thus, using $D_{f}$, we can avoid
talking about \textquotedblleft zig-zag paths\textquotedblright.

\subsubsection{The augmenting path lemma}

The following crucial lemma tells us that such \textquotedblleft zig-zag path
increases\textquotedblright\ are valid (i.e., turn flows into flows), and are
sufficient to find a maximum flow (i.e., if no more \textquotedblleft zig-zag
path increases\textquotedblright\ are possible, then our flow is already maximal):

\begin{lemma}
[augmenting path lemma]\label{lem.net.augpath}Let $N$ be a network consisting
of a multidigraph $D=\left(  V,A,\psi\right)  $, a source $s\in V$, a sink
$t\in V$ and a capacity function $c:A\rightarrow\mathbb{N}$. Let
$f:A\rightarrow\mathbb{N}$ be a flow.

\begin{enumerate}
\item[\textbf{(a)}] If the digraph $D_{f}$ has a path from $s$ to $t$, then
the network $N$ has a flow $f^{\prime}$ with a larger value than $f$.

\item[\textbf{(b)}] If the digraph $D_{f}$ has no path from $s$ to $t$, then
the flow $f$ has maximum value (among all flows on $N$), and there exists a
subset $S$ of $V$ satisfying $s\in S$ and $t\notin S$ and $\left\vert
f\right\vert =c\left(  S,\overline{S}\right)  $.
\end{enumerate}
\end{lemma}

\begin{proof}
\textbf{(a)} Assume that the digraph $D_{f}$ has a path from $s$ to $t$. Pick
such a path, and call it $\mathbf{p}$. Each arc of $\mathbf{p}$ is an arc of
$D_{f}$.

For each forward arc $a\in A$ that appears in $\mathbf{p}$, we have $f\left(
a\right)  <c\left(  a\right)  $ (since $a$ is an arc of $D_{f}$), and thus we
can increase the arc flow $f\left(  a\right)  $ by some positive
$\varepsilon\in\mathbb{N}$ (namely, by any $\varepsilon\leq c\left(  a\right)
-f\left(  a\right)  $) without violating the capacity constraints.\footnote{Of
course, such a unilateral increase will likely violate the conservation
constraints.}

For each backward arc $a^{-1}\in A^{-1}$ that appears in $\mathbf{p}$, we have
$f\left(  a\right)  >0$ (since $a^{-1}$ is an arc of $D_{f}$), and thus we can
decrease the arc flow $f\left(  a\right)  $ by some positive $\varepsilon
\in\mathbb{N}$ (namely, by any $\varepsilon\leq f\left(  a\right)  $) without
violating the capacity constraints.

Let now
\begin{align*}
\varepsilon &  :=\min\Big(\left\{  c\left(  a\right)  -f\left(  a\right)
\ \mid\ a\in A\text{ is a forward arc that appears in }\mathbf{p}\right\} \\
&  \ \ \ \ \ \ \ \ \ \ \ \ \ \ \ \ \ \ \ \ \cup\left\{  f\left(  a\right)
\ \mid\ a^{-1}\in A^{-1}\text{ is a backward arc that appears in }%
\mathbf{p}\right\}  \Big).
\end{align*}
This $\varepsilon$ is a positive integer (since it is a minimum of a set of
positive integers\footnote{because
\par
\begin{itemize}
\item for each forward arc $a\in A$ that appears in $\mathbf{p}$, we have
$f\left(  a\right)  <c\left(  a\right)  $ and thus $c\left(  a\right)
-f\left(  a\right)  >0$;
\par
\item for each backward arc $a^{-1}\in A^{-1}$ that appears in $\mathbf{p}$,
we have $f\left(  a\right)  >0$.
\end{itemize}
}). Let $f^{\prime}:A\rightarrow\mathbb{N}$ be the map obtained from $f$ as follows:

\begin{itemize}
\item For each forward arc $a\in A$ that appears in $\mathbf{p}$, we increase
the arc flow $f\left(  a\right)  $ by $\varepsilon$ (that is, we set
$f^{\prime}\left(  a\right)  :=f\left(  a\right)  +\varepsilon$).

\item For each backward arc $a^{-1}\in A^{-1}$ that appears in $\mathbf{p}$,
we decrease the arc flow $f\left(  a\right)  $ by $\varepsilon$ (that is, we
set $f^{\prime}\left(  a\right)  :=f\left(  a\right)  -\varepsilon$).

\item For all other arcs $a$ of $D$, we keep the arc flow $f\left(  a\right)
$ unchanged (i.e., we set $f^{\prime}\left(  a\right)  :=f\left(  a\right)  $).
\end{itemize}

This new map $f^{\prime}$ still satisfies the capacity
constraints\footnote{since the definition of $\varepsilon$ shows that
\par
\begin{itemize}
\item for each forward arc $a$ that appears in $\mathbf{p}$, we have
$\varepsilon\leq c\left(  a\right)  -f\left(  a\right)  $ and thus $f\left(
a\right)  +\varepsilon\leq c\left(  a\right)  $;
\par
\item for each backward arc $a^{-1}\in A^{-1}$ that appears in $\mathbf{p}$,
we have $\varepsilon\leq f\left(  a\right)  $ and thus $f\left(  a\right)
-\varepsilon\geq0$.
\end{itemize}
}. We claim that it also satisfies the conservation constraints. To check
this, we have to verify that $\left(  f^{\prime}\right)  ^{-}\left(  v\right)
=\left(  f^{\prime}\right)  ^{+}\left(  v\right)  $ for each vertex $v\in
V\setminus\left\{  s,t\right\}  $. So let us do this.

Let $v\in V\setminus\left\{  s,t\right\}  $ be a vertex. We know that
$f^{-}\left(  v\right)  =f^{+}\left(  v\right)  $ (since $f$ is a flow). We
must prove that $\left(  f^{\prime}\right)  ^{-}\left(  v\right)  =\left(
f^{\prime}\right)  ^{+}\left(  v\right)  $.

The path $\mathbf{p}$ is a path from $s$ to $t$. Thus, it neither starts nor
ends at $v$ (since $v\in V\setminus\left\{  s,t\right\}  $). Hence, if $v$ is
a vertex of $\mathbf{p}$, then the path $\mathbf{p}$ enters $v$ by some arc
and exits $v$ by another. Hence, we are in one of the following five cases:

\textit{Case 1:} The vertex $v$ is not a vertex of the path $\mathbf{p}$.

\textit{Case 2:} The path $\mathbf{p}$ enters $v$ by a forward arc and exits
$v$ by a forward arc.

\textit{Case 3:} The path $\mathbf{p}$ enters $v$ by a forward arc and exits
$v$ by a backward arc.

\textit{Case 4:} The path $\mathbf{p}$ enters $v$ by a backward arc and exits
$v$ by a forward arc.

\textit{Case 5:} The path $\mathbf{p}$ enters $v$ by a backward arc and exits
$v$ by a backward arc.

Now, we can prove $\left(  f^{\prime}\right)  ^{-}\left(  v\right)  =\left(
f^{\prime}\right)  ^{+}\left(  v\right)  $ in each of these five cases by
hand. Here is how this can be done in the first three cases:

First, we consider Case 1. In this case, $v$ is not a vertex of the path
$\mathbf{p}$. Hence, each arc $a\in A$ with target $v$ satisfies $f^{\prime
}\left(  a\right)  =f\left(  a\right)  $ (because neither $a$ nor $a^{-1}$
appears in $\mathbf{p}$). Therefore, $\left(  f^{\prime}\right)  ^{-}\left(
v\right)  =f^{-}\left(  v\right)  $. Similarly, $\left(  f^{\prime}\right)
^{+}\left(  v\right)  =f^{+}\left(  v\right)  $. Hence, $\left(  f^{\prime
}\right)  ^{-}\left(  v\right)  =f^{-}\left(  v\right)  =f^{+}\left(
v\right)  =\left(  f^{\prime}\right)  ^{+}\left(  v\right)  $. Thus, we have
proved $\left(  f^{\prime}\right)  ^{-}\left(  v\right)  =\left(  f^{\prime
}\right)  ^{+}\left(  v\right)  $ in Case 1.

Let us now consider Case 2. In this case, the path $\mathbf{p}$ enters $v$ by
a forward arc and exits $v$ by a forward arc. Let $b$ be the former arc, and
$c$ the latter. Then, both $b$ and $c$ are arcs of $D$, and the vertex $v$ is
the target of $b$ and the source of $c$. The definition of $f^{\prime}$ yields
that $f^{\prime}\left(  b\right)  =f\left(  b\right)  +\varepsilon$, whereas
each other arc $a\in A$ with target $v$ satisfies $f^{\prime}\left(  a\right)
=f\left(  a\right)  $. Hence, $\left(  f^{\prime}\right)  ^{-}\left(
v\right)  =f^{-}\left(  v\right)  +\varepsilon$. Similarly, using the arc $c$,
we see that $\left(  f^{\prime}\right)  ^{+}\left(  v\right)  =f^{+}\left(
v\right)  +\varepsilon$. Hence, $\left(  f^{\prime}\right)  ^{-}\left(
v\right)  =\underbrace{f^{-}\left(  v\right)  }_{=f^{+}\left(  v\right)
}+\varepsilon=f^{+}\left(  v\right)  +\varepsilon=\left(  f^{\prime}\right)
^{+}\left(  v\right)  $. Thus, we have proved $\left(  f^{\prime}\right)
^{-}\left(  v\right)  =\left(  f^{\prime}\right)  ^{+}\left(  v\right)  $ in
Case 2.

Let us next consider Case 3. In this case, the path $\mathbf{p}$ enters $v$ by
a forward arc and exits $v$ by a backward arc. Let $b$ be the former arc, and
$c^{-1}$ the latter. Then, both $b$ and $c$ are arcs of $D$, and the vertex
$v$ is the target of both $b$ and $c$. The definition of $f^{\prime}$ yields
that $f^{\prime}\left(  b\right)  =f\left(  b\right)  +\varepsilon$ (since
$\mathbf{p}$ uses the forward arc $b$) and $f^{\prime}\left(  c\right)
=f\left(  c\right)  -\varepsilon$ (since $\mathbf{p}$ uses the backward arc
$c^{-1}$), whereas each other arc $a\in A$ with target $v$ satisfies
$f^{\prime}\left(  a\right)  =f\left(  a\right)  $. Hence, $\left(  f^{\prime
}\right)  ^{-}\left(  v\right)  =f^{-}\left(  v\right)  +\varepsilon
-\varepsilon=f^{-}\left(  v\right)  $. Moreover, $\left(  f^{\prime}\right)
^{+}\left(  v\right)  =f^{+}\left(  v\right)  $ (since none of the arcs of $D$
with source $v$ appears in $\mathbf{p}$, nor does its reversal). Hence,
$\left(  f^{\prime}\right)  ^{-}\left(  v\right)  =f^{-}\left(  v\right)
=f^{+}\left(  v\right)  =\left(  f^{\prime}\right)  ^{+}\left(  v\right)  $.
Thus, we have proved $\left(  f^{\prime}\right)  ^{-}\left(  v\right)
=\left(  f^{\prime}\right)  ^{+}\left(  v\right)  $ in Case 3.

The other two cases are similar (Case 4 is analogous to Case 3, while Case 5
is analogous to Case 2). Thus, altogether, we have proved $\left(  f^{\prime
}\right)  ^{-}\left(  v\right)  =\left(  f^{\prime}\right)  ^{+}\left(
v\right)  $ in all five cases.

Forget that we fixed $v$. We thus have shown that each vertex $v\in
V\setminus\left\{  s,t\right\}  $ satisfies $\left(  f^{\prime}\right)
^{-}\left(  v\right)  =\left(  f^{\prime}\right)  ^{+}\left(  v\right)  $. In
other words, the map $f^{\prime}$ satisfies the conservation constraints.
Since $f^{\prime}$ also satisfies the capacity constraints, we thus conclude
that $f^{\prime}$ is a flow.

What is the value $\left\vert f^{\prime}\right\vert $ of this flow? The path
$\mathbf{p}$ starts at $s$, so it exits $s$ by some arc $\gamma$ (it must have
at least one arc, since $s\neq t$) and never comes back to $s$ again. If this
arc $\gamma$ is a forward arc $b$, then $f^{\prime}\left(  b\right)  =f\left(
b\right)  +\varepsilon$ and therefore $\left(  f^{\prime}\right)  ^{+}\left(
s\right)  =f^{+}\left(  s\right)  +\varepsilon$ and $\left(  f^{\prime
}\right)  ^{-}\left(  s\right)  =f^{-}\left(  s\right)  $. If this arc
$\gamma$ is a backward arc $c^{-1}$, then $f^{\prime}\left(  c\right)
=f\left(  c\right)  -\varepsilon$ and therefore $\left(  f^{\prime}\right)
^{-}\left(  s\right)  =f^{-}\left(  s\right)  -\varepsilon$ and $\left(
f^{\prime}\right)  ^{+}\left(  s\right)  =f^{+}\left(  s\right)  $. Thus,%
\begin{align*}
\left(  f^{\prime}\right)  ^{+}\left(  s\right)  -\left(  f^{\prime}\right)
^{-}\left(  s\right)   &  =%
\begin{cases}
\left(  f^{+}\left(  s\right)  +\varepsilon\right)  -f^{-}\left(  s\right)
, & \text{if }\gamma\text{ is a forward arc};\\
f^{+}\left(  s\right)  -\left(  f^{-}\left(  s\right)  -\varepsilon\right)
, & \text{if }\gamma\text{ is a backward arc}%
\end{cases}
\\
&  =%
\begin{cases}
f^{+}\left(  s\right)  -f^{-}\left(  s\right)  +\varepsilon, & \text{if
}\gamma\text{ is a forward arc};\\
f^{+}\left(  s\right)  -f^{-}\left(  s\right)  +\varepsilon, & \text{if
}\gamma\text{ is a backward arc}%
\end{cases}
\\
&  =\underbrace{f^{+}\left(  s\right)  -f^{-}\left(  s\right)  }%
_{\substack{=\left\vert f\right\vert \\\text{(by the definition of }\left\vert
f\right\vert \text{)}}}+\,\varepsilon=\left\vert f\right\vert +\varepsilon.
\end{align*}
However, the definition of the value $\left\vert f^{\prime}\right\vert $
yields%
\[
\left\vert f^{\prime}\right\vert =\left(  f^{\prime}\right)  ^{+}\left(
s\right)  -\left(  f^{\prime}\right)  ^{-}\left(  s\right)  =\left\vert
f\right\vert +\varepsilon>\left\vert f\right\vert \ \ \ \ \ \ \ \ \ \ \left(
\text{since }\varepsilon>0\right)  .
\]
In other words, the flow $f^{\prime}$ has a larger value than $f$. Thus, we
have found a flow $f^{\prime}$ with a larger value than $f$. This proves Lemma
\ref{lem.net.augpath} \textbf{(a)}. \medskip

\textbf{(b)} Assume that the digraph $D_{f}$ has no path from $s$ to $t$.
Define a subset $S$ of $V$ by%
\[
S=\left\{  v\in V\ \mid\ \text{the digraph }D_{f}\text{ has a path from
}s\text{ to }v\right\}  .
\]
Then, $s\in S$ (because the trivial path $\left(  s\right)  $ is a path from
$s$ to $s$) and $t\notin S$ (since we assumed that $D_{f}$ has no path from
$s$ to $t$). We shall next show that $\left\vert f\right\vert =c\left(
S,\overline{S}\right)  $.

Indeed, we shall obtain this from Proposition \ref{prop.net.flow-setS}
\textbf{(d)}. To do so, we will first show that
\begin{equation}
\left(  f\left(  a\right)  =0\text{ for all }a\in\left[  \overline
{S},S\right]  \right)  \label{pf.lem.net.augpath.b.1}%
\end{equation}
and%
\begin{equation}
\left(  f\left(  a\right)  =c\left(  a\right)  \text{ for all }a\in\left[
S,\overline{S}\right]  \right)  . \label{pf.lem.net.augpath.b.2}%
\end{equation}

[\textit{Proof of (\ref{pf.lem.net.augpath.b.1}):} Let $a\in\left[
\overline{S},S\right]  $. Assume that $f\left(  a\right)  \neq0$. Thus,
$f\left(  a\right)  >0$ (since the capacity constraints yield $f\left(
a\right)  \geq0$). Hence, the backward arc $a^{-1}$ is an arc of the residual
digraph $D_{f}$. Let $u$ be the source of $a$, and let $v$ be the target of
$a$. Since $a\in\left[  \overline{S},S\right]  $, we thus have $u\in
\overline{S}$ and $v\in S$. From $v\in S$, we see that the digraph $D_{f}$ has
a path from $s$ to $v$. Let $\mathbf{q}$ be this path. Appending the backward
arc $a^{-1}$ (which is an arc from $v$ to $u$) and the vertex $u$ to this path
$\mathbf{q}$ (at the end), we obtain a walk from $s$ to $u$ in $D_{f}$. Hence,
$D_{f}$ has a walk from $s$ to $u$, thus also a path from $s$ to $u$ (by
Corollary \ref{cor.mdg.walk-thus-path}). This entails $u\in S$ (by the
definition of $S$). However, this contradicts $u\in\overline{S}=V\setminus S$.
This contradiction shows that our assumption (that $f\left(  a\right)  \neq0$)
was wrong. Therefore, $f\left(  a\right)  =0$. This proves
(\ref{pf.lem.net.augpath.b.1}).]

[\textit{Proof of (\ref{pf.lem.net.augpath.b.2}):} Let $a\in\left[
S,\overline{S}\right]  $. Assume that $f\left(  a\right)  \neq c\left(
a\right)  $. Thus, $f\left(  a\right)  <c\left(  a\right)  $ (since the
capacity constraints yield $f\left(  a\right)  \leq c\left(  a\right)  $).
Hence, the forward arc $a$ is an arc of the residual digraph $D_{f}$. Let $u$
be the source of $a$, and let $v$ be the target of $a$. Since $a\in\left[
S,\overline{S}\right]  $, we thus have $u\in S$ and $v\in\overline{S}$. From
$u\in S$, we see that the digraph $D_{f}$ has a path from $s$ to $u$. Let
$\mathbf{q}$ be this path. Appending the forward arc $a$ (which is an arc from
$u$ to $v$) and the vertex $v$ to this path $\mathbf{q}$ (at the end), we
obtain a walk from $s$ to $v$ in $D_{f}$. Hence, $D_{f}$ has a walk from $s$
to $v$, thus also a path from $s$ to $v$ (by Corollary
\ref{cor.mdg.walk-thus-path}). This entails $v\in S$ (by the definition of
$S$). However, this contradicts $v\in\overline{S}=V\setminus S$. This
contradiction shows that our assumption (that $f\left(  a\right)  \neq
c\left(  a\right)  $) was wrong. Therefore, $f\left(  a\right)  =c\left(
a\right)  $. This proves (\ref{pf.lem.net.augpath.b.2}).]

Now, Proposition \ref{prop.net.flow-setS} \textbf{(d)} yields that $\left\vert
f\right\vert =c\left(  S,\overline{S}\right)  $ holds (since
(\ref{pf.lem.net.augpath.b.1}) and (\ref{pf.lem.net.augpath.b.2}) hold).

We have now found a subset $S$ of $V$ satisfying $s\in S$ and $t\notin S$ and
$\left\vert f\right\vert =c\left(  S,\overline{S}\right)  $. In order to prove
Lemma \ref{lem.net.augpath} \textbf{(b)}, it suffices to show that the flow
$f$ has maximum value (among all flows on $N$). However, this is now easy: Any
flow $g$ on $N$ has value $\left\vert g\right\vert \leq c\left(
S,\overline{S}\right)  $ (by Proposition \ref{prop.net.flow-setS}
\textbf{(c)}, applied to $g$ instead of $f$). In other words, any flow $g$ on
$N$ has value $\left\vert g\right\vert \leq\left\vert f\right\vert $ (since
$\left\vert f\right\vert =c\left(  S,\overline{S}\right)  $). Thus, the flow
$f$ has maximum value. This completes the proof of Lemma \ref{lem.net.augpath}
\textbf{(b)}.
\end{proof}

\subsubsection{Proof of max-flow-min-cut}

We are now ready to prove the max-flow-min-cut theorem (Theorem
\ref{thm.net.MFMC}):

\begin{proof}
[Proof of Theorem \ref{thm.net.MFMC}.]We let $f:A\rightarrow\mathbb{N}$ be the
zero flow on $N$ (see Example \ref{exa.net.zeroflow} for its definition). Now,
we shall incrementally increase the value $\left\vert f\right\vert $ of this
flow by the following algorithm (known as the \textbf{Ford-Fulkerson
algorithm}):

\begin{enumerate}
\item Construct the residual digraph $D_{f}$.

\item If the digraph $D_{f}$ has a path from $s$ to $t$, then Lemma
\ref{lem.net.augpath} \textbf{(a)} shows that the network $N$ has a flow
$f^{\prime}$ with a larger value than $f$ (and furthermore, the proof of Lemma
\ref{lem.net.augpath} \textbf{(a)} shows how to find such an $f^{\prime}$
efficiently\footnote{Of course, this requires an algorithm for finding a path
from $s$ to $t$ in $D_{f}$. But there are many efficient algorithms for this
(see, e.g., homework set \#4 exercise 5).}). Fix such an $f^{\prime}$, and
replace $f$ by $f^{\prime}$. Then, go back to step 1.

\item If the digraph $D_{f}$ has no path from $s$ to $t$, then we end the algorithm.
\end{enumerate}

The replacement of $f$ by $f^{\prime}$ in Step 2 of this algorithm will be
called an \textbf{augmentation}. Thus, the algorithm proceeds by repeatedly
performing augmentations until this is no longer possible.

I claim that the algorithm will eventually end -- i.e., it cannot keep
performing augmentations forever. Indeed, each augmentation increases the
value $\left\vert f\right\vert $ of the flow $f$, and therefore it increases
this value $\left\vert f\right\vert $ by at least $1$ (because increasing an
integer always means increasing it by at least $1$). However, the value
$\left\vert f\right\vert $ is bounded from above by the capacity $c\left(
S,\overline{S}\right)  $ of an arbitrary cut $\left[  S,\overline{S}\right]  $
(by Proposition \ref{prop.net.flow-setS} \textbf{(c)}), and thus cannot get
increased by $1$ more than $c\left(  S,\overline{S}\right)  $ many times
(since its initial value is $0$). Therefore, we cannot perform more than
$c\left(  S,\overline{S}\right)  $ many augmentations in sequence.

Thus, the algorithm eventually ends. Let us consider the flow $f$ that is
obtained once the algorithm has ended. This flow $f$ has the property that the
digraph $D_{f}$ has no path from $s$ to $t$. Thus, Lemma \ref{lem.net.augpath}
\textbf{(b)} shows that the flow $f$ has maximum value (among all flows on
$N$), and there exists a subset $S$ of $V$ satisfying $s\in S$ and $t\notin S$
and $\left\vert f\right\vert =c\left(  S,\overline{S}\right)  $. Consider this
$S$.

Since the flow $f$ has maximum value, we have%
\[
\left\vert f\right\vert =\max\left\{  \left\vert g\right\vert \ \mid\ g\text{
is a flow}\right\}  .
\]

On the other hand, for each subset $T$ of $V$ satisfying $s\in T$ and $t\notin
T$, we have%
\[
c\left(  S,\overline{S}\right)  =\left\vert f\right\vert \leq c\left(
T,\overline{T}\right)
\]
(by Proposition \ref{prop.net.flow-setS} \textbf{(c)}, applied to $T$ instead
of $S$). Hence,
\[
c\left(  S,\overline{S}\right)  =\min\left\{  c\left(  T,\overline{T}\right)
\ \mid\ T\subseteq V;\ s\in T;\ t\notin T\right\}  .
\]
Comparing this with%
\[
c\left(  S,\overline{S}\right)  =\left\vert f\right\vert =\max\left\{
\left\vert g\right\vert \ \mid\ g\text{ is a flow}\right\}  ,
\]
we obtain%
\[
\max\left\{  \left\vert g\right\vert \ \mid\ g\text{ is a flow}\right\}
=\min\left\{  c\left(  T,\overline{T}\right)  \ \mid\ T\subseteq V;\ s\in
T;\ t\notin T\right\}  .
\]
In other words, the maximum value of a flow equals the minimum capacity of a
cut. This proves Theorem \ref{thm.net.MFMC}. (Of course, we cannot use the
letters $f$ and $S$ for the bound variables in $\max\left\{  \left\vert
g\right\vert \ \mid\ g\text{ is a flow}\right\}  $ and \newline$\min\left\{
c\left(  T,\overline{T}\right)  \ \mid\ T\subseteq V;\ s\in T;\ t\notin
T\right\}  $, since $f$ and $S$ already stand for a specific flow and a
specific set.)
\end{proof}

\begin{remark}
All the theorems, propositions and lemmas we proved in this chapter still hold
if we replace the set $\mathbb{N}$ by the set $\mathbb{Q}_{+}:=\left\{
\text{nonnegative rational numbers}\right\}  $ or the set $\mathbb{R}%
_{+}:=\left\{  \text{nonnegative real numbers}\right\}  $. However, their
proofs get more complicated. The problem is that if the arc flows of $f$
belong to $\mathbb{Q}_{+}$ or $\mathbb{R}_{+}$ rather than $\mathbb{N}$, it is
possible for $\left\vert f\right\vert $ to increase endlessly (cf.
\href{https://en.wikipedia.org/wiki/Zenos paradoxes}{Zeno's paradox of
Achilles and the tortoise}), as we make smaller and smaller improvements to
our flow but never achieve (or even approach!) the maximum value.

With rational values, this fortunately cannot happen, since the lowest common
denominator of all arc flows $f\left(  a\right)  $ does not change when we
perform an augmentation. (To put it differently: The case of rational values
can be reduced to the case of integer values by multiplying through with the
lowest common denominator.) With real values, however, this misbehavior can
occur (see \cite[\S I.8]{ForFul74} for an example). Fortunately, there is a
way to avoid it by choosing a \textbf{shortest} path from $s$ to $t$ in
$D_{f}$ at each step. This is known as the \textbf{Edmonds-Karp version of the
Ford-Fulkerson algorithm} (or, for short, the \textbf{Edmonds-Karp
algorithm}). Proving that it works takes a bit more work, which we won't do
here (see, e.g., \cite[Theorem 4.4]{Schrij-ACO}). Incidentally, this technique
also helps keep the algorithm fast for integer-valued flows (running time
$O\left(  \left\vert V\right\vert \cdot\left\vert A\right\vert ^{2}\right)
$). Even faster algorithms exist (see \cite[\S 5.3]{Even12}).
\end{remark}

\subsection{\label{sec.flows.hk}Application: Deriving Hall--K\"{o}nig}

Now, let us apply the max-flow-min-cut theorem to prove the Hall--K\"{o}nig
matching theorem (\ref{thm.match.HKMT}):

\begin{proof}
[Proof of Theorem \ref{thm.match.HKMT} (sketched).](This is an outline; see
\cite[proof of Lemma 1.42]{17s-lec16} for details.\footnote{Note that
\cite[Lemma 1.42]{17s-lec16} is stated only for a simple graph $G$, not for a
multigraph $G$. However, this really makes no difference here: If $\left(
G,X,Y\right)  $ is a bipartite graph with $G$ being a multigraph, then
$\left(  G^{\operatorname*{simp}},X,Y\right)  $ is a bipartite graph as well,
and clearly any matching of $G^{\operatorname*{simp}}$ yields a matching of
$G$ having the same size (and the set $N\left(  U\right)  $ does not change
from $G$ to $G^{\operatorname*{simp}}$ either). Thus, in proving Theorem
\ref{thm.match.HKMT}, we can WLOG assume that $G$ is a simple graph.}) As
explained in Example \ref{exa.net.bip-as-network}, we can turn the bipartite
graph $\left(  G,X,Y\right)  $ into a network so that the matchings of $G$
become the flows $f$ of this network. The max-flow-min-cut theorem (Theorem
\ref{thm.net.MFMC}) yields that
\[
\max\left\{  \left\vert f\right\vert \ \mid\ f\text{ is a flow}\right\}
=\min\left\{  c\left(  S,\overline{S}\right)  \ \mid\ S\subseteq V;\ s\in
S;\ t\notin S\right\}  ,
\]
where $V$ is the vertex set of the digraph that underlies our network. Thus,
there exist a flow $f$ and a cut $\left[  S,\overline{S}\right]  $ of this
network such that $\left\vert f\right\vert =c\left(  S,\overline{S}\right)  $.
Consider these $f$ and $S$. Thus, $S$ is a subset of $V$ such that $s\in S$
and $t\notin S$.

Let $M$ be the matching of $G$ corresponding to the flow $f$ (that is, we let
$M$ be the set of all edges $e$ of $G$ such that $f\left(  \overrightarrow{e}%
\right)  =1$). Thus, $\left\vert M\right\vert =\left\vert f\right\vert $.

Let $U:=X\cap S$. Then, $U$ is a subset of $X$. Here is an illustration of the
cut $\left[  S,\overline{S}\right]  $ on a simple example (the flow $f$ is not
shown):%
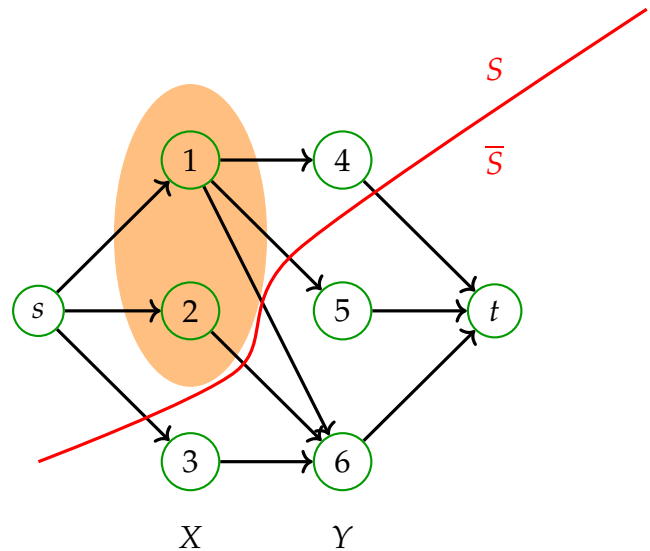
\[%
\begin{tikzpicture}[scale=2]
\filldraw[orange!50!white] (0, 0.5) ellipse (0.5 and 1);
\begin{scope}[every node/.style={circle,thick,draw=green!60!black}]
\node(1) at (0, 1) {$1$};
\node(2) at (0, 0) {$2$};
\node(3) at (0, -1) {$3$};
\node(4) at (1, 1) {$4$};
\node(5) at (1, 0) {$5$};
\node(6) at (1, -1) {$6$};
\node(s) at (-1, 0) {$s$};
\node(t) at (2, 0) {$t$};
\end{scope}
\node(x) at (0, -1.5) {$X$};
\node(x) at (1, -1.5) {$Y$};
\begin{scope}[every edge/.style={draw=black,very thick}]
\path[->] (1) edge (4) edge (5) edge (6);
\path[->] (2) edge (6);
\path[->] (3) edge (6);
\path[->] (s) edge (1) edge (2) edge (3);
\path[->] (4) edge (t) (5) edge (t) (6) edge (t);
\end{scope}
\draw
[red, very thick] plot [smooth] coordinates {(3,2) (0.7,0.4) (0.3,-0.4) (-1, -1)}%
;
\node(X) at (2,1.6) {$\color{red}{S}$};
\node(X) at (2,1) {$\color{red}{\overline{S}}$};
\end{tikzpicture}%
\]
(the orange oval is the set $U$).

Now, we have%
\begin{align*}
\left\vert M\right\vert  &  =\left\vert f\right\vert =c\left(  S,\ \overline
{S}\right)  =\underbrace{c\left(  \left\{  s\right\}  ,\ \overline{S}\right)
}_{\substack{=\left\vert X\setminus U\right\vert \\\text{(why?)}}}+c\left(
\underbrace{X\cap S}_{=U},\ \overline{S}\right)  +\underbrace{c\left(  Y\cap
S,\ \overline{S}\right)  }_{\substack{=\left\vert Y\cap S\right\vert
\\\text{(why?)}}}\\
&  \ \ \ \ \ \ \ \ \ \ \ \ \ \ \ \ \ \ \ \ \left(  \text{since }S\text{ is the
union of the disjoint sets }\left\{  s\right\}  \text{, }X\cap S\text{ and
}Y\cap S\right) \\
&  =\underbrace{\left\vert X\setminus U\right\vert }_{=\left\vert X\right\vert
-\left\vert U\right\vert }+\underbrace{c\left(  U,\ \overline{S}\right)
+\left\vert Y\cap S\right\vert }_{\substack{\geq\left\vert N\left(  U\right)
\right\vert \\\text{(since each vertex }y\in N\left(  U\right)  \text{ either
belongs to }Y\cap S\\\text{and thus contributes to }\left\vert Y\cap
S\right\vert \text{, or belongs to }\overline{S}\\\text{and thus contributes
to }c\left(  U,\ \overline{S}\right)  \text{)}}}\\
&  \geq\left\vert X\right\vert -\left\vert U\right\vert +\left\vert N\left(
U\right)  \right\vert =\left\vert N\left(  U\right)  \right\vert +\left\vert
X\right\vert -\left\vert U\right\vert .
\end{align*}
This proves Theorem \ref{thm.match.HKMT}.
\end{proof}

Having proved the Hall--K\"{o}nig matching theorem (Theorem
\ref{thm.match.HKMT}), we have thus completed the proofs of Hall's marriage
theorem (Theorem \ref{thm.match.HMT}) and of K\"{o}nig's theorem (Theorem
\ref{thm.match.konig}) as well, because we already know how to derive the
latter two theorems from the former.

\subsection{Other applications}

Further applications of the max-flow-min-cut theorem include:

\begin{itemize}
\item A curious fact about rounding matrix entries (stated in terms of a
digraph in \cite[Exercise 4.13]{Schrij-ACO}): Let $A$ be an $m\times n$-matrix
with real entries. Assume that all row sums\footnote{A \textbf{row sum} of a
matrix means the sum of all entries in some row of this matrix. Thus, an
$m\times n$-matrix has $m$ row sums.} of $A$ and all column sums\footnote{A
\textbf{column sum} of a matrix means the sum of all entries in some column of
this matrix. Thus, an $m\times n$-matrix has $n$ column sums.} of $A$ are
integers. Then, we can round each non-integer entry of $A$ (that is, replace
it either by the next-smaller integer or the next-larger integer) in such a
way that the resulting matrix has the same row sums as $A$ and the same column
sums as $A$.

\item An Euler-Hierholzer-like criterion for the existence of an Eulerian
circuit in a \textquotedblleft mixed graph\textquotedblright\ (a general
notion of a graph that can contain both undirected edges and directed arcs)
\cite[\S II.7]{ForFul74}.

\item A proof \cite[\S 6.3]{Berge91} of the Erd\"{o}s--Gallai theorem, which
states that for a given weakly decreasing $n$-tuple $\left(  d_{1}\geq
d_{2}\geq\cdots\geq d_{n}\right)  $ of nonnegative integers, there exists a
simple graph with $n$ vertices whose $n$ vertices have degrees $d_{1}%
,d_{2},\ldots,d_{n}$ if and only if the sum $d_{1}+d_{2}+\cdots+d_{n}$ is even
and each $k\in\left\{  1,2,\ldots,n\right\}  $ satisfies
\[
\sum_{i=1}^{k}d_{i}\leq k\left(  k-1\right)  +\sum_{i=k+1}^{n}\min\left\{
d_{i},k\right\}  .
\]
(The \textquotedblleft only if\textquotedblright\ part of this theorem was
Exercise \ref{exe.2.6} = Exercise 6 on homework set \#2.)
\end{itemize}

The following exercise can be solved both with and without using the
max-flow-min-cut theorem; it should make good practice to solve it in both ways.

\begin{exercise}
\label{exe.flows-cuts.cut-lattice} Consider a network consisting of a
multidigraph $D=\left(  V,A,\psi\right)  $, a source $s\in V$ and a sink $t\in
V$, and a capacity function $c:A\rightarrow\mathbb{N}$ such that $s\neq t$.
(You can replace $\mathbb{N}$ by $\mathbb{Q}_{+}$ or $\mathbb{R}_{+}$ here.)

An $s$\textbf{-}$t$\textbf{-cutting subset} shall mean a subset $S$ of $V$
satisfying $s\in S$ and $t\notin S$.

Let $m$ denote the minimum possible value of $c\left(  S,\overline{S}\right)
$ where $S$ ranges over the $s$-$t$-cutting subsets. (Recall that this is the
maximum value of a flow, according to Theorem \ref{thm.net.MFMC}.)

An $s$-$t$-cutting subset $S$ is said to be \textbf{cut-minimal} if it
satisfies $c\left(  S,\overline{S}\right)  =m$.

Let $X$ and $Y$ be two cut-minimal $s$-$t$-cutting subsets. Prove that $X\cap
Y$ and $X\cup Y$ also are cut-minimal $s$-$t$-cutting subsets.\medskip

[\textbf{Solution:} This is Exercise 7 on homework set \#5 from my Spring 2017
course (except that the simple digraph has been replaced by a multidigraph);
see \href{https://www.cip.ifi.lmu.de/~grinberg/t/17s/}{the course page} for solutions.]
\end{exercise}

\section{More about paths}

In this chapter, we will learn a few more things about paths in graphs and digraphs.

\subsection{Menger's theorems}

We begin with a series of fundamental results known as \textbf{Menger's
theorems} (named after Karl Menger, who discovered one of them in 1927 as an
auxiliary result in a topological study of curves\footnote{See \cite[\S 9.6e]%
{Schrij-CO1} for more about its history.}).

Imagine you have $4$ different ways to get from Philadelphia to NYC, all using
different roads (i.e., no piece of road is used by more than one of your $4$
ways). Then, if $3$ arbitrary roads get blocked, then you still have a way to
get to NYC.

This is obvious (indeed, each blocked road destroys at most one of your $4$
paths, so you still have at least one path left undisturbed after $3$ roads
have been blocked). A more interesting question is the converse: If the road
network is sufficiently robust that blocking $3$ arbitrary roads will not
disconnect you from NYC, does this mean that you can find $4$ different ways
to NYC all using different roads?

Menger's theorems answer this question (and various questions of this kind) in
the positive, in several different setups. Each of these theorems can be
roughly described as \textquotedblleft the maximum number of pairwise
independent paths from some place to another place equals the minimum size of
a bottleneck that separates the former from the latter\textquotedblright.
Here, the \textquotedblleft places\textquotedblright\ can be vertices or sets
of vertices; the word \textquotedblleft independent\textquotedblright\ can
mean \textquotedblleft having no arcs in common\textquotedblright\ or
\textquotedblleft having no intermediate vertices in common\textquotedblright%
\ or \textquotedblleft having no vertices at all in common\textquotedblright;
and the word \textquotedblleft bottleneck\textquotedblright\ can mean a set of
arcs or of vertices whose removal would disconnect the former place from the
latter. Here is a quick overview of all Menger's theorems that we will
prove:\footnote{All undefined terminology used here will be defined further
below.}

\begin{itemize}
\item for directed graphs:%
\[
\hspace{-3.5263pc}%
\begin{tabular}
[c]{|c|c|c|c|}\hline
Theorem ... & the places are ... & the paths must be ... & the bottleneck
consists of ...\\\hline\hline
\ref{thm.menger.da1} & vertices & arc-disjoint & arcs\\\hline
\ref{thm.menger.da2} & vertices & arc-disjoint & arcs of a cut\\\hline
\ref{thm.menger.dam} & sets of vertices & arc-disjoint & arcs of a cut\\\hline
\ref{thm.menger.dv2} & vertices & internally vertex-disjoint & vertices $\in
V\setminus\left\{  s,t\right\}  $\\\hline
\ref{thm.menger.dvm} & sets of vertices & internally vertex-disjoint &
vertices $\in V\setminus\left(  X\cup Y\right)  $\\\hline
\ref{thm.menger.dvm-vd} & sets of vertices & vertex-disjoint & vertices $\in
V$\\\hline
\end{tabular}
\ \
\]

\item for undirected graphs:%
\[
\hspace{-3.5263pc}%
\begin{tabular}
[c]{|c|c|c|c|}\hline
Theorem ... & the places are ... & the paths must be ... & the bottleneck
consists of ...\\\hline\hline
\ref{thm.menger.ue1} & vertices & edge-disjoint & edges\\\hline
\ref{thm.menger.ue2} & vertices & edge-disjoint & edges of a cut\\\hline
\ref{thm.menger.uv2} & vertices & internally vertex-disjoint & vertices $\in
V\setminus\left\{  s,t\right\}  $\\\hline
\ref{thm.menger.uvm} & sets of vertices & internally vertex-disjoint &
vertices $\in V\setminus\left(  X\cup Y\right)  $\\\hline
\ref{thm.menger.uvm-vd} & sets of vertices & vertex-disjoint & vertices $\in
V$\\\hline
\end{tabular}
\ \
\]
(I could state more, but I don't want this to go on forever.)
\end{itemize}

\subsubsection{The arc-Menger theorem for directed graphs}

We begin with the most natural setup: a directed graph (one-way roads) with
roads being arcs. The following definitions will help keep the theorems short:

\begin{definition}
\label{def.menger.adis}Two walks $\mathbf{p}$ and $\mathbf{q}$ in a digraph
are said to be \textbf{arc-disjoint} if they have no arc in common.
\end{definition}

\begin{example}
\label{exa.menger.adis.1}The following picture shows two arc-disjoint paths
$\mathbf{p}$ and $\mathbf{q}$ (they can be told apart by their labels: each
arc of $\mathbf{p}$ is labelled with a \textquotedblleft$\mathbf{p}%
$\textquotedblright, and likewise for $\mathbf{q}$):%
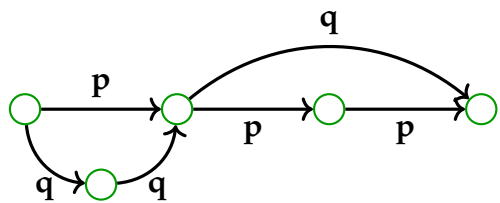
\[%
%
\ \ .
\]

The following picture shows two paths $\mathbf{r}$ and $\mathbf{s}$ that are
\textbf{not} arc-disjoint (the common arc is marked with \textquotedblleft%
$\mathbf{r},\mathbf{s}$\textquotedblright):%
\[%
%
\ \ .
\]

\end{example}

\begin{definition}
\label{def.menger.arc-sep}Let $D=\left(  V,A,\psi\right)  $ be a multidigraph,
and let $s$ and $t$ be two vertices of $D$. A subset $B$ of $A$ is said to be
an $s$\textbf{-}$t$\textbf{-arc-separator} if each path from $s$ to $t$
contains at least one arc from $B$. Equivalently, a subset $B$ of $A$ is said
to be an $s$\textbf{-}$t$\textbf{-arc-separator} if the multidigraph $\left(
V,\ A\setminus B,\ \psi\mid_{A\setminus B}\right)  $ has no path from $s$ to
$t$ (in other words, removing from $D$ all arcs contained in $B$ destroys all
paths from $s$ to $t$).
\end{definition}

\begin{example}
\label{exa.menger.arc-sep.1}Let $D=\left(  V,A,\psi\right)  $ be the following
multidigraph:%
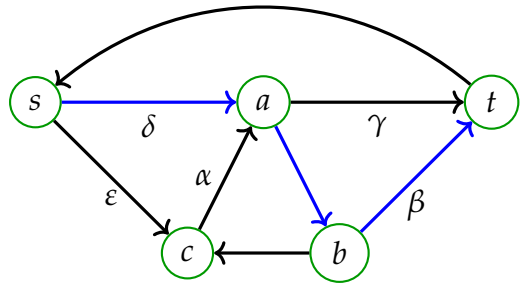
\[%
\begin{tikzpicture}[scale=2]
\begin{scope}[every node/.style={circle,thick,draw=green!60!black}]
\node(s) at (-1,0) {$s$};
\node(a) at (0.5, 0) {$a$};
\node(c) at (0,-1) {$c$};
\node(b) at (1,-1) {$b$};
\node(t) at (2,0) {$t$};
\end{scope}
\begin{scope}[every edge/.style={draw=black,very thick}, every loop/.style={}]
\path[->] (s) edge node[below] {$\varepsilon$} (c);
\path[->] (b) edge (c);
\path[->] (c) edge node[left] {$\alpha$} (a);
\path[->] (a) edge node[below] {$\gamma$} (t);
\path[->] (t) edge[bend right=40] (s);
\end{scope}
\begin{scope}[every edge/.style={draw=blue,very thick}, every loop/.style={}]
\path[->] (a) edge (b);
\path[->] (s) edge node[below] {$\delta$} (a);
\path[->] (b) edge node[below] {$\beta$} (t);
\end{scope}
\end{tikzpicture}%
\ \ .
\]
Then, the set $\left\{  \alpha,\gamma\right\}  $ is not an $s$-$t$%
-arc-separator (since the path drawn in blue contains no arc from this set).
However, the set $\left\{  \beta,\gamma\right\}  $ is an $s$-$t$%
-arc-separator, and so is the set $\left\{  \delta,\varepsilon\right\}  $. Of
course, any set that contains any of $\left\{  \beta,\gamma\right\}  $ and
$\left\{  \delta,\varepsilon\right\}  $ as a subset is therefore an $s$%
-$t$-arc-separator as well.
\end{example}

\begin{example}
Let $D$ be a multidigraph. Let $s$ and $t$ be two vertices of $D$. Then, the
empty set $\varnothing$ is an $s$-$t$-arc-separator if and only if $D$ has no
path from $s$ to $t$. This degenerate case should not be forgotten!
\end{example}

We can now state the first Menger's theorem:

\begin{theorem}
[arc-Menger theorem for directed graphs, version 1]\label{thm.menger.da1}Let
$D=\left(  V,A,\psi\right)  $ be a multidigraph, and let $s$ and $t$ be two
distinct vertices of $D$. Then, the maximum number of pairwise arc-disjoint
paths from $s$ to $t$ equals the minimum size of an $s$-$t$-arc-separator.
\end{theorem}

\begin{example}
Let $D$ be the multidigraph from Example \ref{exa.menger.arc-sep.1}. Then, the
minimum size of an $s$-$t$-arc-separator is $2$ (indeed, $\left\{
\beta,\gamma\right\}  $ is an $s$-$t$-arc-separator of size $2$, and it is
easy to see that there are no $s$-$t$-arc-separators of smaller size). Hence,
Theorem \ref{thm.menger.da1} yields that the maximum number of pairwise
arc-disjoint paths from $s$ to $t$ is $2$ as well. And indeed, we can easily
find $2$ arc-disjoint paths from $s$ to $t$, namely the red and the blue paths
in the following figure:%
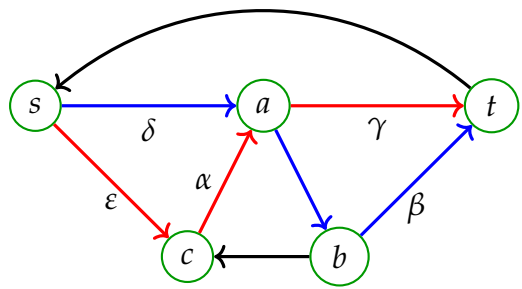
\[%
\begin{tikzpicture}[scale=2]
\begin{scope}[every node/.style={circle,thick,draw=green!60!black}]
\node(s) at (-1,0) {$s$};
\node(a) at (0.5, 0) {$a$};
\node(c) at (0,-1) {$c$};
\node(b) at (1,-1) {$b$};
\node(t) at (2,0) {$t$};
\end{scope}
\begin{scope}[every edge/.style={draw=black,very thick}, every loop/.style={}]
\path[->] (b) edge (c);
\path[->] (t) edge[bend right=40] (s);
\end{scope}
\begin{scope}[every edge/.style={draw=blue,very thick}, every loop/.style={}]
\path[->] (s) edge node[below] {$\delta$} (a);
\path[->] (a) edge (b);
\path[->] (b) edge node[below] {$\beta$} (t);
\end{scope}
\begin{scope}[every edge/.style={draw=red,very thick}, every loop/.style={}]
\path[->] (s) edge node[below] {$\varepsilon$} (c);
\path[->] (c) edge node[left] {$\alpha$} (a);
\path[->] (a) edge node[below] {$\gamma$} (t);
\end{scope}
\end{tikzpicture}%
\ \ .
\]

\end{example}

Before proving Theorem \ref{thm.menger.da1}, let me state another variant of
this theorem, which is closer to the proof. First, some notations:

\begin{definition}
\label{def.menger.da-cuts}Let $D=\left(  V,A,\psi\right)  $ be a multidigraph,
and let $s$ and $t$ be two distinct vertices of $D$.

\begin{enumerate}
\item[\textbf{(a)}] For each subset $S$ of $V$, we set $\overline
{S}:=V\setminus S$ and
\begin{align*}
\left[  S,\overline{S}\right]   &  :=\left\{  a\in A\ \mid\ \text{the source
of }a\text{ belongs to }S\text{,}\right. \\
&  \ \ \ \ \ \ \ \ \ \ \ \ \ \ \ \ \ \ \ \ \left.  \text{and the target of
}a\text{ belongs to }\overline{S}\right\}  .
\end{align*}
(These are the same definitions that we introduced for networks in Definition
\ref{def.net.capacity}.)

\item[\textbf{(b)}] An $s$\textbf{-}$t$\textbf{-cut} means a subset of $A$
that has the form $\left[  S,\overline{S}\right]  $, where $S$ is a subset of
$V$ that satisfies $s\in S$ and $t\notin S$. (This was just called a
\textquotedblleft cut\textquotedblright\ back in Definition \ref{def.net.cut}
\textbf{(a)}.)
\end{enumerate}
\end{definition}

An $s$-$t$-cut is called this way because its removal would cut the vertex $s$
from the vertex $t$. More precisely:

\begin{remark}
\label{rmk.menger.cut}Let $D=\left(  V,A,\psi\right)  $ be a multidigraph, and
let $s$ and $t$ be two distinct vertices of $D$. Then, any $s$-$t$-cut is an
$s$-$t$-arc-separator.
\end{remark}

\begin{proof}
Let $B$ be an $s$-$t$-cut. We must prove that $B$ is an $s$-$t$-arc-separator.
In other words, we must prove that each path from $s$ to $t$ contains at least
one arc from $B$.

We know that $B$ is an $s$-$t$-cut. In other words, $B=\left[  S,\overline
{S}\right]  $, where $S$ is a subset of $V$ that satisfies $s\in S$ and
$t\notin S$. Consider this subset $S$.

Each path from $s$ to $t$ starts at a vertex in $S$ (since $s\in S$) and ends
at a vertex outside of $S$ (since $t\notin S$). Thus, each such path has to
escape the set $S$ at some point -- i.e., it must contain an arc whose source
is in $S$ and whose target is outside of $S$. But such an arc must necessarily
belong to $\left[  S,\overline{S}\right]  $ (by the definition of $\left[
S,\overline{S}\right]  $). Thus, each path from $s$ to $t$ must contain an arc
from $\left[  S,\overline{S}\right]  $. In other words, each path from $s$ to
$t$ must contain an arc from $B$ (since $B=\left[  S,\overline{S}\right]  $).
In other words, $B$ is an $s$-$t$-arc-separator (by the definition of an
$s$-$t$-arc-separator). This proves Remark \ref{rmk.menger.cut}.
\end{proof}

\begin{theorem}
[arc-Menger theorem for directed graphs, version 2]\label{thm.menger.da2}Let
$D=\left(  V,A,\psi\right)  $ be a multidigraph, and let $s$ and $t$ be two
distinct vertices of $D$. Then, the maximum number of pairwise arc-disjoint
paths from $s$ to $t$ equals the minimum size of an $s$-$t$-cut.
\end{theorem}

\begin{example}
\label{exa.menger.da2.1}Let $D$ be the following multidigraph:%
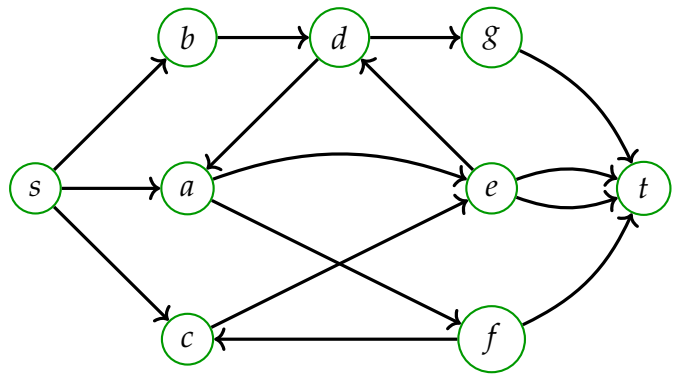
\[%
%
\ \ .
\]
Then, the maximum number of pairwise arc-disjoint paths from $s$ to $t$ is
$3$. Indeed, the following picture shows $3$ such paths in red, blue and
brown, respectively:%
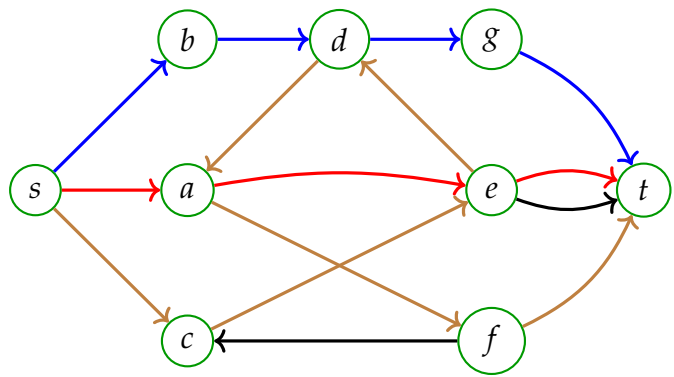
\[%
%
\ \ .
\]
More than $3$ pairwise arc-disjoint paths from $s$ to $t$ cannot exist in $D$,
since (e.g.) there are only $3$ arcs outgoing from $s$.

By Theorem \ref{thm.menger.da2}, this shows that the minimum size of an
$s$-$t$-cut in $D$ is $3$ as well. There are many $s$-$t$-cuts of size $3$
(for instance, the \textquotedblleft obvious\textquotedblright\ cut $\left[
\left\{  s\right\}  ,\overline{\left\{  s\right\}  }\right]  $ has this
property, as does the $s$-$t$-cut $\left[  \left\{  s,a,f\right\}
,\overline{\left\{  s,a,f\right\}  }\right]  $).

Let us now reverse of the direction of the arc from $c$ to $e$ in $D$ (thus
destroying the brown path). The resulting multidigraph $D^{\prime}$ looks as
follows:%
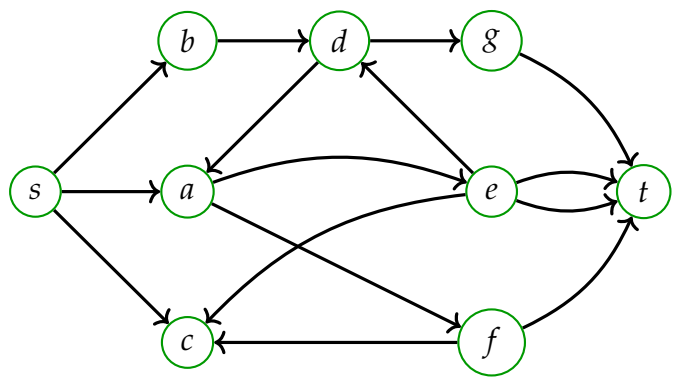
\[%
\begin{tikzpicture}[scale=2]
\begin{scope}[every node/.style={circle,thick,draw=green!60!black}]
\node(s) at (-1,0) {$s$};
\node(a) at (0,0) {$a$};
\node(b) at (0,1) {$b$};
\node(c) at (0,-1) {$c$};
\node(d) at (1,1) {$d$};
\node(e) at (2,0) {$e$};
\node(f) at (2,-1) {$f$};
\node(g) at (2,1) {$g$};
\node(t) at (3,0) {$t$};
\end{scope}
\begin{scope}[every edge/.style={draw=black,very thick}, every loop/.style={}]
\path[->] (s) edge (a);
\path[->] (s) edge (b);
\path[->] (s) edge (c);
\path[->] (b) edge (d);
\path[->] (a) edge[bend left=20] (e);
\path[->] (a) edge (f);
\path[->] (e) edge[bend right=20] (c);
\path[->] (d) edge (a);
\path[->] (e) edge (d);
\path[->] (e) edge[bend left=20] (t) edge[bend right=20] (t);
\path[->] (f) edge (c);
\path[->] (f) edge[bend right=20] (t);
\path[->] (d) edge (g);
\path[->] (g) edge[bend left=20] (t);
\end{scope}
\end{tikzpicture}%
\ \ .
\]
This digraph $D^{\prime}$ has no more than $2$ pairwise arc-disjoint paths
from $s$ to $t$. This can be seen by observing that the $s$-$t$-cut $\left[
\left\{  s,c\right\}  ,\overline{\left\{  s,c\right\}  }\right]  $ has size
$2$ (it consists of the arc from $s$ to $a$ and the arc from $s$ to $b$), so
that the minimum size of an $s$-$t$-cut is at most $2$, and therefore (by
Theorem \ref{thm.menger.da2}) the maximum number of pairwise arc-disjoint
paths from $s$ to $t$ is at most $2$ as well. It is easy to see that the
latter number is exactly $2$ (since our red and blue paths still exist in
$D^{\prime}$).
\end{example}

To prove the above two arc-Menger theorems, we need one more lemma about
networks. We recall the notations from Section \ref{sec.flows.def} and from
Definition \ref{def.net.cut}, and introduce a couple more:

\begin{definition}
Let $D=\left(  V,A,\psi\right)  $ be a multidigraph. Let $f,g:A\rightarrow
\mathbb{N}$ be two maps. Then:

\begin{enumerate}
\item[\textbf{(a)}] We let $f+g$ denote the map from $A$ to $\mathbb{N}$ that
sends each arc $a\in A$ to $f\left(  a\right)  +g\left(  a\right)  $. (This is
the pointwise sum of $f$ and $g$.)

\item[\textbf{(b)}] We write $g\leq f$ if and only if each arc $a\in A$
satisfies $g\left(  a\right)  \leq f\left(  a\right)  $.

\item[\textbf{(c)}] If $g\leq f$, then we let $f-g$ denote the map from $A$ to
$\mathbb{N}$ that sends each arc $a\in A$ to $f\left(  a\right)  -g\left(
a\right)  $. (This is really a map to $\mathbb{N}$, since $g\leq f$ entails
$g\left(  a\right)  \leq f\left(  a\right)  $.)
\end{enumerate}
\end{definition}

These notations satisfy the properties that you'd expect: e.g., the pointwise
sum of maps from $A$ to $\mathbb{N}$ is associative (meaning that $\left(
f+g\right)  +h=f+\left(  g+h\right)  $, so that you can write $f+g+h$ for both
sides); inequalities can be manipulated in the usual way (e.g., we have
$f-g\leq h$ if and only if $f\leq g+h$). Verifying this all is straightforward.

The following definition codifies the flows that we constructed in Remark
\ref{rmk.net.pathflow}:

\begin{definition}
\label{def.net.pathflow}Let $N$ be a network consisting of a multidigraph
$D=\left(  V,A,\psi\right)  $, a source $s\in V$, a sink $t\in V$ and a
capacity function $c:A\rightarrow\mathbb{N}$. Let $\mathbf{p}$ be a path from
$s$ to $t$ in $D$. Then, we define a map $f_{\mathbf{p}}:A\rightarrow
\mathbb{N}$ by setting%
\[
f_{\mathbf{p}}\left(  a\right)  =%
\begin{cases}
1, & \text{if }a\text{ is an arc of }\mathbf{p};\\
0, & \text{otherwise}%
\end{cases}
\ \ \ \ \ \ \ \ \ \ \text{for each }a\in A.
\]
We call this map $f_{\mathbf{p}}$ the \textbf{path flow} of $\mathbf{p}$. It
is an actual flow of value $1$ if all the arcs of $\mathbf{p}$ have capacity
$\geq1$.
\end{definition}

\begin{example}
Consider the following network:%
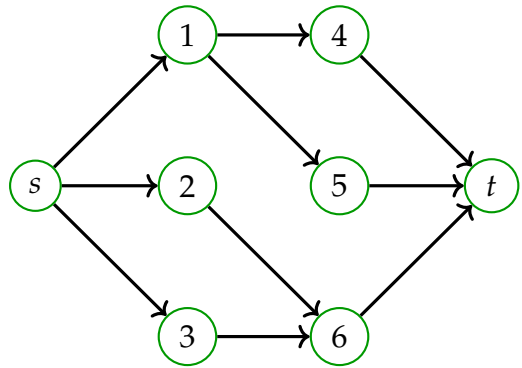
\[%
%
\ \ ,
\]
where each arc has capacity $1$. Then, the path $\mathbf{p}=\left(
s,2,6,t\right)  $ leads to the following path flow $f_{\mathbf{p}}$:%
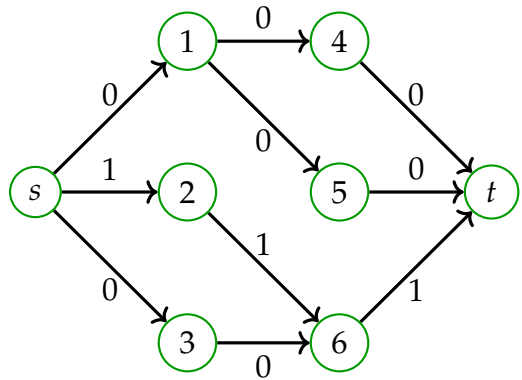
\[%
%
\ \ .
\]
Here, in order not to crowd the picture, we have left out the
\textquotedblleft of $1$\textquotedblright\ part of the label of each arc (so
you should read the \textquotedblleft$0$\textquotedblright s and the
\textquotedblleft$1$\textquotedblright s atop the arcs as \textquotedblleft$0$
of $1$\textquotedblright\ and \textquotedblleft$1$ of $1$\textquotedblright, respectively).
\end{example}

The path flow thus turns any path from $s$ to $t$ in a network into a flow,
provided that the arcs have enough capacity to carry this flow. If we have $m$
paths $\mathbf{p}_{1},\mathbf{p}_{2},\ldots,\mathbf{p}_{m}$ from $s$ to $t$,
then we can add their path flows together, and obtain a flow $f_{\mathbf{p}%
_{1}}+f_{\mathbf{p}_{2}}+\cdots+f_{\mathbf{p}_{m}}$ of value $m$, provided
(again) that the arcs have enough capacity for it. (In general, we cannot
uniquely reconstruct $\mathbf{p}_{1},\mathbf{p}_{2},\ldots,\mathbf{p}_{m}$
back from this latter flow, as they might have gotten \textquotedblleft mixed
together\textquotedblright.)

Our next lemma can be viewed as a (partial) converse of this observation: Any
flow $f$ of value $m$ \textquotedblleft contains\textquotedblright\ a sum
$f_{\mathbf{p}_{1}}+f_{\mathbf{p}_{2}}+\cdots+f_{\mathbf{p}_{m}}$ of $m$ path
flows $f_{\mathbf{p}_{1}},f_{\mathbf{p}_{2}},\ldots,f_{\mathbf{p}_{m}}$
corresponding to $m$ (not necessarily distinct) paths $\mathbf{p}%
_{1},\mathbf{p}_{2},\ldots,\mathbf{p}_{m}$ from $s$ to $t$. Here, the word
\textquotedblleft contains\textquotedblright\ signals that $f$ is not
necessarily equal to $f_{\mathbf{p}_{1}}+f_{\mathbf{p}_{2}}+\cdots
+f_{\mathbf{p}_{m}}$, but only satisfies $f_{\mathbf{p}_{1}}+f_{\mathbf{p}%
_{2}}+\cdots+f_{\mathbf{p}_{m}}\leq f$ in general. So here is the lemma:

\begin{lemma}
[flow path decomposition lemma]\label{lem.net.decomp1}Let $N$ be a network
consisting of a multidigraph $D=\left(  V,A,\psi\right)  $, a source $s\in V$,
a sink $t\in V$ and a capacity function $c:A\rightarrow\mathbb{N}$. Let $f$ be
a flow on $N$ that has value $m$. Then, there exist $m$ paths $\mathbf{p}%
_{1},\mathbf{p}_{2},\ldots,\mathbf{p}_{m}$ from $s$ to $t$ in $D$ such that%
\[
f_{\mathbf{p}_{1}}+f_{\mathbf{p}_{2}}+\cdots+f_{\mathbf{p}_{m}}\leq f.
\]

\end{lemma}

\begin{proof}
We induct on $m$.

The \textit{base case} ($m=0$) is obvious (since the empty sum $f_{\mathbf{p}%
_{1}}+f_{\mathbf{p}_{2}}+\cdots+f_{\mathbf{p}_{m}}$ is the zero flow, and thus
is $\leq f$ because of the capacity constraints).

\textit{Induction step:} Let $m$ be a positive integer. Assume (as the
induction hypothesis) that the lemma holds for $m-1$. We must prove the lemma
for $m$.

So we consider a flow $f$ on $N$ that has value $m$. We need to show that
there exist $m$ paths $\mathbf{p}_{1},\mathbf{p}_{2},\ldots,\mathbf{p}_{m}$
from $s$ to $t$ in $D$ such that $f_{\mathbf{p}_{1}}+f_{\mathbf{p}_{2}}%
+\cdots+f_{\mathbf{p}_{m}}\leq f$.

We shall first find some path $\mathbf{p}$ from $s$ to $t$ such that
$f_{\mathbf{p}}\leq f$.

We shall refer to the arcs $a\in A$ satisfying $f\left(  a\right)  >0$ as the
\textbf{active} arcs. Let $A^{\prime}:=\left\{  a\in A\ \mid\ f\left(
a\right)  >0\right\}  $ be the set of these active arcs. Consider the spanning
subdigraph $D^{\prime}:=\left(  V,A^{\prime},\psi\mid_{A^{\prime}}\right)  $
of $D$.

Let $S$ be the set of all vertices $v\in V$ such that $D^{\prime}$ has a path
from $s$ to $v$. Then, $s\in S$ (since the trivial path $\left(  s\right)  $
is a path of $D^{\prime}$).

We next claim that each arc $b\in\left[  S,\overline{S}\right]  $ satisfies
$f\left(  b\right)  =0$.

[\textit{Proof:} Assume the contrary. Thus, some arc $b\in\left[
S,\overline{S}\right]  $ satisfies $f\left(  b\right)  \neq0$. Consider this
$b$. From $f\left(  b\right)  \neq0$, we obtain $f\left(  b\right)  >0$ (since
$f$ is a flow), thus $b\in A^{\prime}$ (by the definition of $A^{\prime}$).
Hence, $b$ is an arc of $D^{\prime}$ (by the definition of $D^{\prime}$).

Let $u$ be the source of the arc $b$, and $v$ its target. Since $b\in\left[
S,\overline{S}\right]  $, we therefore have $u\in S$ and $v\in\overline{S}$.
Since $u\in S$, the digraph $D^{\prime}$ has a path $\mathbf{p}$ from $s$ to
$u$ (by the definition of $S$). Consider this path $\mathbf{p}$. Appending the
arc $b$ and the vertex $v$ at the end of this path $\mathbf{p}$, we obtain a
walk from $s$ to $v$ in $D^{\prime}$ (since $b$ is an arc of $D^{\prime}$ with
source $u$ and target $v$). Hence, the digraph $D^{\prime}$ has a walk from
$s$ to $v$, thus also a path from $s$ to $v$ (by Corollary
\ref{cor.mdg.walk-thus-path}). This means that $v\in S$ (by the definition of
$S$). But this contradicts $v\in\overline{S}=V\setminus S$. This contradiction
shows that our assumption was wrong, qed.]

We thus have proved that each $b\in\left[  S,\overline{S}\right]  $ satisfies
$f\left(  b\right)  =0$. Therefore, $f\left(  S,\overline{S}\right)  =0$
(using the notations of Definition \ref{def.net.capacity} \textbf{(d)}).
However, recall that $s\in S$. Thus, if we had $t\notin S$, then Proposition
\ref{prop.net.flow-setS} \textbf{(b)} would yield
\[
\left\vert f\right\vert =\underbrace{f\left(  S,\overline{S}\right)  }%
_{=0}-\underbrace{f\left(  \overline{S},S\right)  }_{\geq0}\leq0-0=0,
\]
which would contradict $\left\vert f\right\vert =m>0$. Hence, we must have
$t\in S$. In other words, the digraph $D^{\prime}$ has a path from $s$ to $t$
(by the definition of $S$). Let $\mathbf{p}$ be this path. Then, $\mathbf{p}$
is also a path in $D$ and satisfies $f_{\mathbf{p}}\leq f$%
\ \ \ \ \footnote{\textit{Proof.} We need to prove that each arc $a\in A$
satisfies $f_{\mathbf{p}}\left(  a\right)  \leq f\left(  a\right)  $.
\par
So let $a\in A$ be an arc. If $a$ is not an arc of $\mathbf{p}$, then the
definition of $f_{\mathbf{p}}$ yields $f_{\mathbf{p}}\left(  a\right)  =0\leq
f\left(  a\right)  $ (since $f$ is a flow), so we are done in this case.
Hence, assume WLOG that $a$ is an arc of $\mathbf{p}$. Thus, $a$ is an arc of
$D^{\prime}$ (since $\mathbf{p}$ is a path of $D^{\prime}$). In other words,
$a\in A^{\prime}$. By the definition of $A^{\prime}$, this means that
$f\left(  a\right)  >0$. Since $f\left(  a\right)  $ is an integer, we thus
have $f\left(  a\right)  \geq1=f_{\mathbf{p}}\left(  a\right)  $ (since $a$ is
an arc of $\mathbf{p}$). In other words, $f_{\mathbf{p}}\left(  a\right)  \leq
f\left(  a\right)  $. This is precisely what we wanted to prove.}. Therefore,
$f-f_{\mathbf{p}}$ is a map from $A$ to $\mathbb{N}$. Moreover,
$f-f_{\mathbf{p}}$ is again a flow\footnote{Here, we are using the fact (which
is straightforward to prove) that if $g$ and $h$ are two flows with $h\leq g$,
then $g-h$ is again a flow.}, and has value $\left\vert f-f_{\mathbf{p}%
}\right\vert =m-1$\ \ \ \ \footnote{Here, we are using the fact (which is
straightforward to prove) that if $g$ and $h$ are two flows satisfying $h\leq
g$, then $\left\vert g-h\right\vert =\left\vert g\right\vert -\left\vert
h\right\vert $. Applying this fact to $g=f$ and $h=f_{\mathbf{p}}$, we obtain
$\left\vert f-f_{\mathbf{p}}\right\vert =\underbrace{\left\vert f\right\vert
}_{=m}-\underbrace{\left\vert f_{\mathbf{p}}\right\vert }_{=1}=m-1$.}. Thus,
by the induction hypothesis, we can apply Lemma \ref{lem.net.decomp1} to $m-1$
and $f-f_{\mathbf{p}}$ instead of $m$ and $f$. As a result, we conclude that
there exist $m-1$ paths $\mathbf{p}_{1},\mathbf{p}_{2},\ldots,\mathbf{p}%
_{m-1}$ from $s$ to $t$ in $D$ such that $f_{\mathbf{p}_{1}}+f_{\mathbf{p}%
_{2}}+\cdots+f_{\mathbf{p}_{m-1}}\leq f-f_{\mathbf{p}}$. Consider these $m-1$
paths $\mathbf{p}_{1},\mathbf{p}_{2},\ldots,\mathbf{p}_{m-1}$, and set
$\mathbf{p}_{m}:=\mathbf{p}$. Then, $f_{\mathbf{p}_{1}}+f_{\mathbf{p}_{2}%
}+\cdots+f_{\mathbf{p}_{m-1}}\leq f-f_{\mathbf{p}}=f-f_{\mathbf{p}_{m}}$
(since $\mathbf{p}=\mathbf{p}_{m}$), so that $f_{\mathbf{p}_{1}}%
+f_{\mathbf{p}_{2}}+\cdots+f_{\mathbf{p}_{m}}\leq f$.

Thus, we have found $m$ paths $\mathbf{p}_{1},\mathbf{p}_{2},\ldots
,\mathbf{p}_{m}$ from $s$ to $t$ in $D$ such that $f_{\mathbf{p}_{1}%
}+f_{\mathbf{p}_{2}}+\cdots+f_{\mathbf{p}_{m}}\leq f$. But this is precisely
what we wanted. Thus, the induction step is complete, and Lemma
\ref{lem.net.decomp1} is proved. \medskip
\end{proof}

\begin{remark}
\label{rmk.net.decomp1.euler-proof}There exists an alternative proof of Lemma
\ref{lem.net.decomp1}, which is too nice to leave unmentioned. Here is a quick
outline: Consider a new multidigraph that is obtained from $D$ by replacing
each arc $a$ by $f\left(  a\right)  $ many parallel arcs (if $f\left(
a\right)  =0$, this means that $a$ is simply removed). Add $m$ many arcs from
$t$ to $s$ to this new multidigraph. The resulting digraph is balanced
(because of the conservation constraints for $f$). It may fail to be weakly
connected; however, the vertices $s$ and $t$ belong to the same weak component
of it (as long as $m>0$). Hence, applying the directed Euler-Hierholzer
theorem (Theorem \ref{thm.digr.euler-hier} \textbf{(a)}) to this component, we
see that this component has an Eulerian circuit. Cutting the $m$ arcs from $t$
to $s$ out of this circuit, we obtain $m$ arc-disjoint walks from $s$ to $t$.
Each of these $m$ walks contains some path from $s$ to $t$, and thus we obtain
$m$ paths $\mathbf{p}_{1},\mathbf{p}_{2},\ldots,\mathbf{p}_{m}$ from $s$ to
$t$ in $D$ such that $f_{\mathbf{p}_{1}}+f_{\mathbf{p}_{2}}+\cdots
+f_{\mathbf{p}_{m}}\leq f$.
\end{remark}

\begin{remark}
Let $N$ be a network consisting of a multidigraph $D=\left(  V,A,\psi\right)
$, a source $s\in V$, a sink $t\in V$ and a capacity function $c:A\rightarrow
\mathbb{N}$. If $\mathbf{c}$ is a cycle of $D$, then we can define a map
$f_{\mathbf{c}}:A\rightarrow\mathbb{N}$ by setting%
\[
f_{\mathbf{c}}\left(  a\right)  =%
\begin{cases}
1, & \text{if }a\text{ is an arc of }\mathbf{c};\\
0, & \text{otherwise}%
\end{cases}
\ \ \ \ \ \ \ \ \ \ \text{for each }a\in A.
\]
We call this map $f_{\mathbf{c}}$ the \textbf{cycle flow} of $\mathbf{c}$. It
is an actual flow of value $0$ if all the arcs of $\mathbf{c}$ have capacity
$\geq1$.

Now, the conclusion of Lemma \ref{lem.net.decomp1} can be improved as follows:
There exist $m$ paths $\mathbf{p}_{1},\mathbf{p}_{2},\ldots,\mathbf{p}_{m}$
from $s$ to $t$ in $D$ as well as a (possibly empty) collection of cycles
$\mathbf{c}_{1},\mathbf{c}_{2},\ldots,\mathbf{c}_{k}$ of $D$ such that%
\[
f=\left(  f_{\mathbf{p}_{1}}+f_{\mathbf{p}_{2}}+\cdots+f_{\mathbf{p}_{m}%
}\right)  +\left(  f_{\mathbf{c}_{1}}+f_{\mathbf{c}_{2}}+\cdots+f_{\mathbf{c}%
_{k}}\right)  .
\]
Proving this improved claim is a bit harder than proving Lemma
\ref{lem.net.decomp1}, but not by too much (in particular, the argument in
Remark \ref{rmk.net.decomp1.euler-proof} can be adapted, since a walk becomes
a path if we successively remove all cycles from it).
\end{remark}

\begin{proof}
[Proof of Theorem \ref{thm.menger.da2}.]We make $D$ into a network $N$ (with
source $s$ and sink $t$) by assigning the capacity $1$ to each arc $a\in A$.
Clearly, a cut of this network is the same as what we call an $s$-$t$-cut.
Moreover, the capacity $c\left(  S,\overline{S}\right)  $ of a cut $\left[
S,\overline{S}\right]  $ is simply the size of this cut (since each arc has
capacity $1$).

The max-flow-min-cut theorem (Theorem \ref{thm.net.MFMC}) tells us that the
maximum value of a flow equals the minimum capacity of a cut, i.e., the
minimum size of an $s$-$t$-cut (because, as we just explained, a cut is the
same as an $s$-$t$-cut, and its capacity is simply its size). It thus remains
to show that the maximum value of a flow is the maximum number of pairwise
arc-disjoint paths from $s$ to $t$. But this is easy by now:

\begin{itemize}
\item If you have a flow $f$ of value $m$, then you can find $m$ pairwise
arc-disjoint paths from $s$ to $t$ (because Lemma \ref{lem.net.decomp1} gives
you $m$ paths $\mathbf{p}_{1},\mathbf{p}_{2},\ldots,\mathbf{p}_{m}$ such that
$f_{\mathbf{p}_{1}}+f_{\mathbf{p}_{2}}+\cdots+f_{\mathbf{p}_{m}}\leq f$, and
the latter inequality tells you that these $m$ paths $\mathbf{p}%
_{1},\mathbf{p}_{2},\ldots,\mathbf{p}_{m}$ are
arc-disjoint\footnote{\textit{Proof.} Assume the contrary. Thus, these $m$
paths are not arc-disjoint. In other words, there exists an arc $a$ that is
used by two paths $\mathbf{p}_{i}$ and $\mathbf{p}_{j}$ with $i\neq j$.
Consider this arc $a$ and the corresponding indices $i$ and $j$. Since $a$ is
used by $\mathbf{p}_{i}$, we have $f_{\mathbf{p}_{i}}\left(  a\right)  =1$.
Likewise, $f_{\mathbf{p}_{j}}\left(  a\right)  =1$. However, $f_{\mathbf{p}%
_{1}}+f_{\mathbf{p}_{2}}+\cdots+f_{\mathbf{p}_{m}}\leq f$, so that%
\begin{align*}
\left(  f_{\mathbf{p}_{1}}+f_{\mathbf{p}_{2}}+\cdots+f_{\mathbf{p}_{m}%
}\right)  \left(  a\right)   &  \leq f\left(  a\right)  \leq c\left(
a\right)  \ \ \ \ \ \ \ \ \ \ \left(  \text{by the capacity constraints}%
\right) \\
&  =1\ \ \ \ \ \ \ \ \ \ \left(  \text{since each arc has capacity }1\right)
.
\end{align*}
Thus,%
\begin{align*}
1  &  \geq\left(  f_{\mathbf{p}_{1}}+f_{\mathbf{p}_{2}}+\cdots+f_{\mathbf{p}%
_{m}}\right)  \left(  a\right)  =f_{\mathbf{p}_{1}}\left(  a\right)
+f_{\mathbf{p}_{2}}\left(  a\right)  +\cdots+f_{\mathbf{p}_{m}}\left(
a\right) \\
&  \geq\underbrace{f_{\mathbf{p}_{i}}\left(  a\right)  }_{=1}%
+\underbrace{f_{\mathbf{p}_{j}}\left(  a\right)  }_{=1}%
\ \ \ \ \ \ \ \ \ \ \left(
\begin{array}
[c]{c}%
\text{since }f_{\mathbf{p}_{i}}\left(  a\right)  \text{ and }f_{\mathbf{p}%
_{j}}\left(  a\right)  \text{ are two distinct addends}\\
\text{of the sum }f_{\mathbf{p}_{1}}\left(  a\right)  +f_{\mathbf{p}_{2}%
}\left(  a\right)  +\cdots+f_{\mathbf{p}_{m}}\left(  a\right) \\
\text{(because }i\neq j\text{), and since all the remaining}\\
\text{addends are }\geq0\text{ (since }f_{\mathbf{p}}\left(  a\right)
\geq0\text{ for each path }\mathbf{p}\text{)}%
\end{array}
\right) \\
&  =1+1>1,
\end{align*}
which is absurd. This contradiction shows that our assumption was false,
qed.}). Thus,
\begin{align}
&  \left(  \text{the maximum number of pairwise arc-disjoint paths from
}s\text{ to }t\right) \nonumber\\
&  \geq\left(  \text{the maximum value of a flow}\right)  .
\label{pf.thm.menger.da2.1}%
\end{align}

\item Conversely, if you have $m$ pairwise arc-disjoint paths $\mathbf{p}%
_{1},\mathbf{p}_{2},\ldots,\mathbf{p}_{m}$ from $s$ to $t$, then you obtain a
flow of value $m$ (namely, $f_{\mathbf{p}_{1}}+f_{\mathbf{p}_{2}}%
+\cdots+f_{\mathbf{p}_{m}}$ is such a flow\footnote{\textit{Proof.} First, we
observe that the map $f_{\mathbf{p}_{1}}+f_{\mathbf{p}_{2}}+\cdots
+f_{\mathbf{p}_{m}}$ satisfies the conservation constraints (because it is the
sum of the functions $f_{\mathbf{p}_{1}},f_{\mathbf{p}_{2}},\ldots
,f_{\mathbf{p}_{m}}$, each of which satisfies the conservation constraints).
Let us now check that it satisfies the capacity constraints.
\par
Indeed, let $a\in A$ be an arc. Then, $a$ belongs to at most one of the $m$
paths $\mathbf{p}_{1},\mathbf{p}_{2},\ldots,\mathbf{p}_{m}$ (since these $m$
paths are arc-disjoint). In other words, at most one of the $m$ numbers
$f_{\mathbf{p}_{1}}\left(  a\right)  ,\ f_{\mathbf{p}_{2}}\left(  a\right)
,\ \ldots,\ f_{\mathbf{p}_{m}}\left(  a\right)  $ equals $1$; all the
remaining numbers equal $0$. Hence, the sum $f_{\mathbf{p}_{1}}\left(
a\right)  +f_{\mathbf{p}_{2}}\left(  a\right)  +\cdots+f_{\mathbf{p}_{m}%
}\left(  a\right)  $ of these $m$ numbers equals either $1$ or $0$; in either
case, we thus have $f_{\mathbf{p}_{1}}\left(  a\right)  +f_{\mathbf{p}_{2}%
}\left(  a\right)  +\cdots+f_{\mathbf{p}_{m}}\left(  a\right)  \in\left\{
0,1\right\}  $. Now,%
\[
\left(  f_{\mathbf{p}_{1}}+f_{\mathbf{p}_{2}}+\cdots+f_{\mathbf{p}_{m}%
}\right)  \left(  a\right)  =f_{\mathbf{p}_{1}}\left(  a\right)
+f_{\mathbf{p}_{2}}\left(  a\right)  +\cdots+f_{\mathbf{p}_{m}}\left(
a\right)  \in\left\{  0,1\right\}  ,
\]
so that%
\[
0\leq\left(  f_{\mathbf{p}_{1}}+f_{\mathbf{p}_{2}}+\cdots+f_{\mathbf{p}_{m}%
}\right)  \left(  a\right)  \leq1=c\left(  a\right)
\]
(since each arc has capacity $1$). Since we have proved this for each arc
$a\in A$, we thus have shown that the map $f_{\mathbf{p}_{1}}+f_{\mathbf{p}%
_{2}}+\cdots+f_{\mathbf{p}_{m}}$ satisfies the capacity constraints. Hence,
this map is a flow (since it also satisfies the conservation constraints).
\par
It remains to show that the value of this flow is $m$. But this is easy: For
any flows $g_{1},g_{2},\ldots,g_{k}$, we have $\left\vert g_{1}+g_{2}%
+\cdots+g_{k}\right\vert =\left\vert g_{1}\right\vert +\left\vert
g_{2}\right\vert +\cdots+\left\vert g_{k}\right\vert $ (this is
straightforward to see from the definition of value). Thus,
\[
\left\vert f_{\mathbf{p}_{1}}+f_{\mathbf{p}_{2}}+\cdots+f_{\mathbf{p}_{m}%
}\right\vert =\left\vert f_{\mathbf{p}_{1}}\right\vert +\left\vert
f_{\mathbf{p}_{2}}\right\vert +\cdots+\left\vert f_{\mathbf{p}_{m}}\right\vert
=\sum_{k=1}^{m}\underbrace{\left\vert f_{\mathbf{p}_{k}}\right\vert }%
_{=1}=\sum_{k=1}^{m}1=m.
\]
In other words, the value of the flow $f_{\mathbf{p}_{1}}+f_{\mathbf{p}_{2}%
}+\cdots+f_{\mathbf{p}_{m}}$ is $m$.}). Thus,%
\begin{align*}
&  \left(  \text{the maximum value of a flow}\right) \\
&  \geq\left(  \text{the maximum number of pairwise arc-disjoint paths from
}s\text{ to }t\right)  .
\end{align*}

\end{itemize}

Combining this last inequality with (\ref{pf.thm.menger.da2.1}), we obtain%
\begin{align*}
&  \left(  \text{the maximum number of pairwise arc-disjoint paths from
}s\text{ to }t\right) \\
&  =\left(  \text{the maximum value of a flow}\right) \\
&  =\left(  \text{the minimum size of an }s\text{-}t\text{-cut}\right)
\ \ \ \ \ \ \ \ \ \ \left(  \text{as we have proved before}\right)  .
\end{align*}
Thus, Theorem \ref{thm.menger.da2} is proved.
\end{proof}

Theorem \ref{thm.menger.da2} can also be proved without using network flows
(see, e.g., \cite[Corollary 4.1b]{Schrij-ACO} for such a proof).

\begin{proof}
[Proof of Theorem \ref{thm.menger.da1}.]Let $x$ denote the maximum number of
pairwise arc-disjoint paths from $s$ to $t$.

Let $n_{c}$ denote the minimum size of an $s$-$t$-cut.

Let $n_{s}$ denote the minimum size of an $s$-$t$-arc-separator.\footnote{If
you are wondering why we chose the baroque notations \textquotedblleft%
$x$\textquotedblright, \textquotedblleft$n_{c}$\textquotedblright\ and
\textquotedblleft$n_{s}$\textquotedblright\ for these three numbers: The
letter \textquotedblleft$x$\textquotedblright\ appears in \textquotedblleft
maximum\textquotedblright, whereas the letter \textquotedblleft$n$%
\textquotedblright\ appears in \textquotedblleft minimum\textquotedblright.
The subscripts \textquotedblleft$c$\textquotedblright\ and \textquotedblleft%
$s$\textquotedblright\ should be reasonably clear.}

Theorem \ref{thm.menger.da2} says that $x=n_{c}$. Our goal is to prove that
$x=n_{s}$.

Remark \ref{rmk.menger.cut} shows that any $s$-$t$-cut is an $s$%
-$t$-arc-separator. Thus, $n_{s}\leq n_{c}$.

The inequality $x\leq n_{s}$ follows easily from the pigeonhole
principle\footnote{\textit{Proof.} We know that there exist $x$ pairwise
arc-disjoint paths from $s$ to $t$ (by the definition of $x$). Let
$\mathbf{p}_{1},\mathbf{p}_{2},\ldots,\mathbf{p}_{x}$ be these $x$ paths.
\par
We know that there exists an $s$-$t$-arc-separator of size $n_{s}$ (by the
definition of $n_{s}$). Let $B$ be this $s$-$t$-arc-separator. Thus, each path
from $s$ to $t$ contains at least one arc from $B$ (by the definition of an
$s$-$t$-arc-separator). Hence, in particular, each of the $x$ paths
$\mathbf{p}_{1},\mathbf{p}_{2},\ldots,\mathbf{p}_{x}$ contains at least one
arc from $B$. These altogether $x$ arcs must be distinct (since the $x$ paths
$\mathbf{p}_{1},\mathbf{p}_{2},\ldots,\mathbf{p}_{x}$ are arc-disjoint); thus,
we have found at least $x$ arcs that belong to $B$. This shows that
$\left\vert B\right\vert \geq x$. However, $B$ has size $n_{s}$; in other
words, we have $\left\vert B\right\vert =n_{s}$. Thus, $n_{s}=\left\vert
B\right\vert \geq x$, so that $x\leq n_{s}$.}. Combining this with $n_{s}\leq
n_{c}=x$ (since $x=n_{c}$), we obtain $x=n_{s}$. Thus, Theorem
\ref{thm.menger.da1} is proved.
\end{proof}

\begin{exercise}
\label{exe.9.6}Let $D$ be a balanced multidigraph. Let $s$ and $t$ be two
vertices of $D$. Let $k\in\mathbb{N}$. Assume that $D$ has $k$ pairwise
arc-disjoint paths from $s$ to $t$. Show that $D$ has $k$ pairwise
arc-disjoint paths from $t$ to $s$.
\end{exercise}

\begin{exercise}
\label{exe.mt2.menger-postnikov} Let $D$ be a multidigraph. Let $k\in
\mathbb{N}$. Let $u$, $v$ and $w$ be three vertices of $D$. Assume that there
exist $k$ arc-disjoint paths from $u$ to $v$. Assume furthermore that there
exist $k$ arc-disjoint paths from $v$ to $w$.

Prove that there exist $k$ arc-disjoint paths from $u$ to $w$.

[\textbf{Note:} If $u=w$, then the trivial path $\left(  u\right)  $ counts as
being arc-disjoint from itself (so in this case, there exist arbitrarily many
arc-disjoint paths from $u$ to $w$).] \medskip

[\textbf{Solution:} This is Exercise 3 on midterm \#2 from my Spring 2017
course (except that it is stated for multidigraphs instead of simple
digraphs); see \href{https://www.cip.ifi.lmu.de/~grinberg/t/17s/}{the course
page} for solutions.]
\end{exercise}

\bigskip

We can also extend the arc-Menger theorem to paths between different pairs of vertices:

\begin{theorem}
[arc-Menger theorem for directed graphs, multi-terminal version]%
\label{thm.menger.dam}Let $D=\left(  V,A,\psi\right)  $ be a multidigraph, and
let $X$ and $Y$ be two disjoint subsets of $V$.

A \textbf{path from }$X$\textbf{ to }$Y$ shall mean a path whose starting
point belongs to $X$ and whose ending point belongs to $Y$.

An $X$\textbf{-}$Y$\textbf{-cut} shall mean a subset of $A$ that has the form
$\left[  S,\overline{S}\right]  $, where $S$ is a subset of $V$ that satisfies
$X\subseteq S$ and $Y\subseteq\overline{S}$.

Then, the maximum number of pairwise arc-disjoint paths from $X$ to $Y$ equals
the minimum size of an $X$-$Y$-cut.
\end{theorem}

\begin{example}
\label{exa.menger.dam.1}Here is an example of a digraph $D=\left(
V,A,\psi\right)  $, with two disjoint subsets $X$ and $Y$ of $V$ drawn as
ovals:%
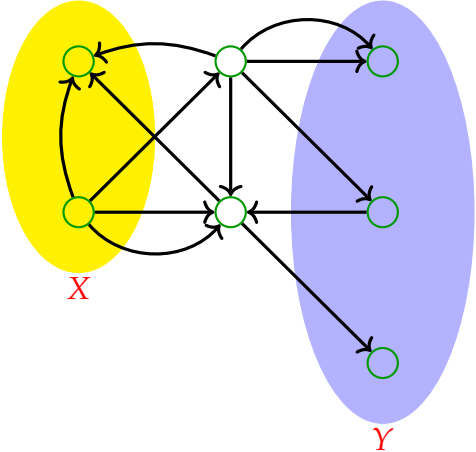
\[%
%
\ \ .
\]
In this digraph $D$, the maximum number of pairwise arc-disjoint paths from
$X$ to $Y$ is $2$; here are two such paths (marked in red and blue):%
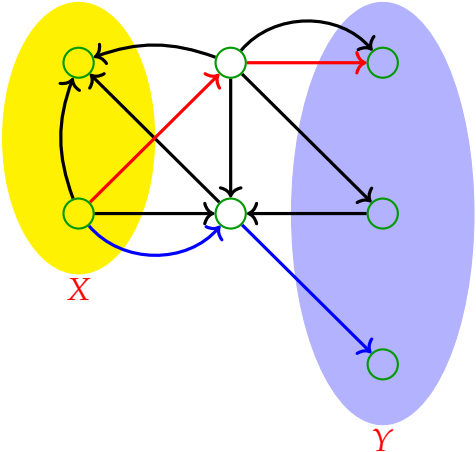
\[%
%
\ \ .
\]
According to Theorem \ref{thm.menger.dam}, the minimum size of an $X$-$Y$-cut
must thus also be $2$. And indeed, here is such an $X$-$Y$-cut:%
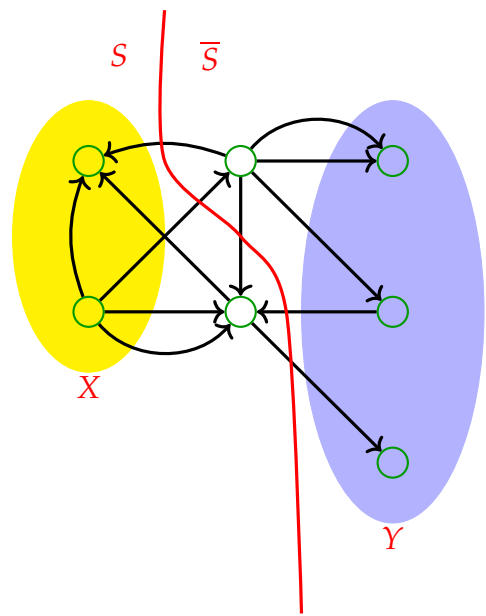
\[%
%
\ \ .
\]

\end{example}

\begin{proof}
[Proof of Theorem \ref{thm.menger.dam}.]We transform our digraph $D=\left(
V,A,\psi\right)  $ into a new multidigraph $D^{\prime}=\left(  V^{\prime
},A^{\prime},\psi^{\prime}\right)  $ as follows:

\begin{itemize}
\item We replace all the vertices in $X$ by a single (new) vertex $s$, and
replace all the vertices in $Y$ by a single (new) vertex $t$. (Thus, formally
speaking, we set $V^{\prime}=\left(  V\setminus\left(  X\cup Y\right)
\right)  \cup\left\{  s,t\right\}  $, where $s$ and $t$ are two objects not in
$V$.)

For any vertex $p\in V$, we define a vertex $p^{\prime}\in V^{\prime}$ by%
\[
p^{\prime}=%
\begin{cases}
s, & \text{if }p\in X;\\
t, & \text{if }p\in Y;\\
p, & \text{otherwise.}%
\end{cases}
\]
We refer to this vertex $p^{\prime}$ as the \textbf{projection} of $p$.

\item We keep all the arcs of $D$ around, but we replace all their endpoints
(i.e., sources and targets) by their projections (thus, any endpoint in $X$
gets replaced by $s$, and any endpoint in $Y$ gets replaced by $t$, while an
endpoint that belongs neither to $X$ nor to $Y$ stays unchanged). For example,
an arc with source in $X$ becomes an arc with source in $s$. (Formally
speaking, this means the following: We set $A^{\prime}=A$ and we define the
map $\psi^{\prime}:A^{\prime}\rightarrow V^{\prime}\times V^{\prime}$ as
follows: For any $a\in A^{\prime}=A$, we set $\psi^{\prime}\left(  a\right)
=\left(  u^{\prime},v^{\prime}\right)  $, where $\left(  u,v\right)
=\psi\left(  a\right)  $.)
\end{itemize}

For instance, if $D$ is the digraph from Example \ref{exa.menger.dam.1}, then
$D^{\prime}$ looks as follows:%
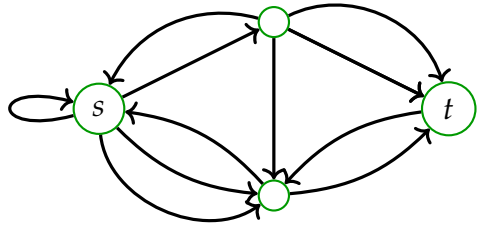
\[%
\begin{tikzpicture}[scale=2.3]
\begin{scope}[every node/.style={circle,thick,draw=green!60!black}]
\node(s) at (-1,0.5) {$s$};
\node(3) at (0,1) {};
\node(4) at (0,0) {};
\node(t) at (1,0.5) {$t$};
\end{scope}
\begin{scope}[every edge/.style={draw=black,very thick}, every loop/.style={}]
\path[->] (s) edge[loop left] (s);
\path[->] (3) edge[bend right=40] (s);
\path[->] (4) edge[bend right=20] (s);
\path[->] (s) edge (3);
\path[->] (s) edge[bend right=20] (4);
\path[->] (s) edge[bend right=60] (4);
\path[->] (3) edge (4);
\path[->] (3) edge (t);
\path[->] (3) edge[bend left=50] (t);
\path[->] (3) edge (t);
\path[->] (4) edge[bend right=20] (t);
\path[->] (t) edge[bend right=20] (4);
\end{scope}
\end{tikzpicture}%
\ \ .
\]

Now, Theorem \ref{thm.menger.da2} (applied to $D^{\prime}=\left(  V^{\prime
},A^{\prime},\psi^{\prime}\right)  $ instead of $D=\left(  V,A,\psi\right)  $)
shows that the maximum number of pairwise arc-disjoint paths from $s$ to $t$
in $D^{\prime}$ equals the minimum size of an $s$-$t$-cut in $D^{\prime}$.

Let us now connect this with the claim that we want to prove. It is easy to
see that the minimum size of an $s$-$t$-cut in $D^{\prime}$ equals the minimum
size of an $X$-$Y$-cut in $D$ (indeed, the $s$-$t$-cuts in $D^{\prime}$ are
precisely the $X$-$Y$-cuts in $D$\ \ \ \ \footnote{In more detail:
\par
\begin{itemize}
\item Any $s$-$t$-cut in $D^{\prime}$ has the form $\left[  S,\overline
{S}\right]  $ for some subset $S$ of $V^{\prime}$ satisfying $s\in S$ and
$t\notin S$; it is therefore equal to the set $\left[  S^{\prime}%
,\overline{S^{\prime}}\right]  $, where $S^{\prime}$ is the subset of $V$
given by $S^{\prime}:=\left(  S\setminus\left\{  s\right\}  \right)  \cup X$.
Therefore, it is an $X$-$Y$-cut in $D$.
\par
\item Conversely, any $X$-$Y$-cut in $D$ has the form $\left[  S,\overline
{S}\right]  $ for some subset $S$ of $V$ satisfying $X\subseteq S$ and
$Y\subseteq\overline{S}$; it is therefore equal to the set $\left[  S^{\prime
},\overline{S^{\prime}}\right]  $, where $S^{\prime}$ is the subset of
$V^{\prime}$ given by $S^{\prime}:=\left(  S\setminus X\right)  \cup\left\{
s\right\}  $. Therefore, it is an $s$-$t$-cut in $D^{\prime}$.
\end{itemize}
}). If we can also show that the maximum number of pairwise arc-disjoint paths
from $s$ to $t$ in $D^{\prime}$ equals the maximum number of pairwise
arc-disjoint paths from $X$ to $Y$ in $D$, then the result of the preceding
paragraph will thus become the claim of Theorem \ref{thm.menger.dam}, so we
will be done.

So how can we show that the maximum number of pairwise arc-disjoint paths from
$s$ to $t$ in $D^{\prime}$ equals the maximum number of pairwise arc-disjoint
paths from $X$ to $Y$ in $D$ ? It would be easy if there was a well-behaved
bijection between the former paths and the latter paths that preserves the
arcs of any path, but this is not quite the case. Each path from $X$ to $Y$ in
$D$ becomes a \textbf{walk} from $s$ to $t$ in $D^{\prime}$ if we replace each
of its vertices by its projection.\ However, the latter walk is not
necessarily a path, since different vertices can have the same projection.

Fortunately, this is easy to fix. If we have $k$ pairwise arc-disjoint paths
from $X$ to $Y$ in $D$, then we can turn them into $k$ pairwise arc-disjoint
\textbf{walks} from $s$ to $t$ in $D^{\prime}$, and then we also obtain $k$
pairwise arc-disjoint \textbf{paths} from $s$ to $t$ in $D^{\prime}$ (since
any walk from $s$ to $t$ contains a path from $s$ to $t$). Thus,%
\begin{align*}
&  \left(  \text{the maximum number of pairwise arc-disjoint paths from
}s\text{ to }t\text{ in }D^{\prime}\right) \\
&  \geq\left(  \text{the maximum number of pairwise arc-disjoint paths from
}X\text{ to }Y\text{ in }D\right)  .
\end{align*}
Conversely, if we have $k$ pairwise arc-disjoint paths from $s$ to $t$ in
$D^{\prime}$, then we can \textquotedblleft lift\textquotedblright\ these $k$
paths back to the digraph $D$ (preserving the arcs, and replacing the vertices
$s$ and $t$ by appropriate vertices in $X$ and $Y$ to make them belong to the
right arcs), and thus obtain $k$ pairwise arc-disjoint paths from $X$ to $Y$
in $D$. Thus,%
\begin{align*}
&  \left(  \text{the maximum number of pairwise arc-disjoint paths from
}X\text{ to }Y\text{ in }D\right) \\
&  \geq\left(  \text{the maximum number of pairwise arc-disjoint paths from
}s\text{ to }t\text{ in }D^{\prime}\right)  .
\end{align*}
Combining these two inequalities, we obtain%
\begin{align*}
&  \left(  \text{the maximum number of pairwise arc-disjoint paths from
}s\text{ to }t\text{ in }D^{\prime}\right) \\
&  =\left(  \text{the maximum number of pairwise arc-disjoint paths from
}X\text{ to }Y\text{ in }D\right)  .
\end{align*}
As explained above, this completes the proof of Theorem \ref{thm.menger.dam}.
\end{proof}

\subsubsection{The edge-Menger theorem for undirected graphs}

We shall now state analogues of Theorem \ref{thm.menger.da1} and Theorem
\ref{thm.menger.da2} for undirected graphs. First, the unsurprising definitions:

\begin{definition}
Two walks $\mathbf{p}$ and $\mathbf{q}$ in a graph are said to be
\textbf{edge-disjoint} if they have no edge in common.
\end{definition}

\begin{definition}
Let $G=\left(  V,E,\varphi\right)  $ be a multigraph, and let $s$ and $t$ be
two vertices of $G$. A subset $B$ of $E$ is said to be an $s$\textbf{-}%
$t$\textbf{-edge-separator} if each path from $s$ to $t$ contains at least one
edge from $B$. Equivalently, a subset $B$ of $E$ is said to be an
$s$\textbf{-}$t$\textbf{-edge-separator} if the multigraph $\left(
V,\ E\setminus B,\ \varphi\mid_{E\setminus B}\right)  $ has no path from $s$
to $t$ (in other words, removing from $G$ all edges contained in $B$ destroys
all paths from $s$ to $t$).
\end{definition}

Now comes the analogue of Theorem \ref{thm.menger.da1}:

\begin{theorem}
[edge-Menger theorem for undirected graphs, version 1]\label{thm.menger.ue1}%
Let $G=\left(  V,E,\varphi\right)  $ be a multigraph, and let $s$ and $t$ be
two distinct vertices of $G$. Then, the maximum number of pairwise
edge-disjoint paths from $s$ to $t$ equals the minimum size of an $s$-$t$-edge-separator.
\end{theorem}

To state the analogue of Theorem \ref{thm.menger.da2}, we need to first adopt
Definition \ref{def.menger.da-cuts} to undirected graphs:

\begin{definition}
\label{def.menger.ue-cuts}Let $G=\left(  V,E,\varphi\right)  $ be a
multigraph, and let $s$ and $t$ be two distinct vertices of $G$.

\begin{enumerate}
\item[\textbf{(a)}] For each subset $S$ of $V$, we set $\overline
{S}:=V\setminus S$ and
\begin{align*}
\left[  S,\overline{S}\right]  _{\operatorname*{und}}  &  :=\left\{  e\in
E\ \mid\ \text{one endpoint of }e\text{ belongs to }S\text{,}\right. \\
&  \ \ \ \ \ \ \ \ \ \ \ \ \ \ \ \ \ \ \ \ \left.  \text{while the other
belongs to }\overline{S}\right\}  .
\end{align*}

\item[\textbf{(b)}] An (\textbf{undirected) }$s$\textbf{-}$t$\textbf{-cut}
means a subset of $E$ that has the form $\left[  S,\overline{S}\right]
_{\operatorname*{und}}$, where $S$ is a subset of $V$ that satisfies $s\in S$
and $t\notin S$.
\end{enumerate}
\end{definition}

The following remark is an analogue of Remark \ref{rmk.menger.cut}:

\begin{remark}
\label{rmk.menger.udcut}Let $G=\left(  V,E,\varphi\right)  $ be a multigraph,
and let $s$ and $t$ be two distinct vertices of $G$. Then, any (undirected)
$s$-$t$-cut is an $s$-$t$-edge-separator.
\end{remark}

\begin{proof}
Analogous to the proof of Remark \ref{rmk.menger.cut}.
\end{proof}

And here is the analogue of Theorem \ref{thm.menger.da2}:

\begin{theorem}
[edge-Menger theorem for undirected graphs, version 2]\label{thm.menger.ue2}%
Let $G=\left(  V,E,\varphi\right)  $ be a multigraph, and let $s$ and $t$ be
two distinct vertices of $G$. Then, the maximum number of pairwise
edge-disjoint paths from $s$ to $t$ equals the minimum size of an (undirected)
$s$-$t$-cut.
\end{theorem}

\begin{proof}
[Proof of Theorem \ref{thm.menger.ue2}.]We shall not prove this from scratch,
but rather derive this from the directed version (Theorem \ref{thm.menger.da2}).

Namely, we apply Theorem \ref{thm.menger.da2} to\footnote{Recall that
$G^{\operatorname*{bidir}}$ is the multidigraph obtained from $G$ by replacing
each edge by two arcs in opposite directions. (If the edge has endpoints $u$
and $v$, then one of the two arcs has source $u$ and target $v$, while the
other has source $v$ and target $u$.) See Definition \ref{def.mg.bidir} for a
formal definition.} $D=G^{\operatorname*{bidir}}$. We thus see that the
maximum number of pairwise arc-disjoint paths from $s$ to $t$ (in
$G^{\operatorname*{bidir}}$) equals the minimum size of an $s$-$t$-cut (in
$G^{\operatorname*{bidir}}$). This is similar to the claim that we want to
prove, but not quite the same statement, because $G^{\operatorname*{bidir}}$
is not $G$. To obtain the claim that we want to prove, we must prove the
following two claims:

\begin{statement}
\textit{Claim 1:} The maximum number of pairwise arc-disjoint paths from $s$
to $t$ (in $G^{\operatorname*{bidir}}$) equals the maximum number of pairwise
edge-disjoint paths from $s$ to $t$ (in $G$).
\end{statement}

\begin{statement}
\textit{Claim 2:} The minimum size of a directed $s$-$t$-cut\footnote{A
\textquotedblleft directed $s$-$t$-cut\textquotedblright\ here simply means an
$s$-$t$-cut in a digraph.} (in $G^{\operatorname*{bidir}}$) equals the minimum
size of an (undirected) $s$-$t$-cut (in $G$).
\end{statement}

Claim 2 is very easy to verify, since the directed $s$-$t$-cuts in
$G^{\operatorname*{bidir}}$ are essentially the same as the undirected $s$%
-$t$-cuts in $G$\ \ \ \ \footnote{In more detail: If $S$ is a subset of $V$
that satisfies $s\in S$ and $t\notin S$, then the directed $s$-$t$-cut
$\left[  S,\overline{S}\right]  $ in $G^{\operatorname*{bidir}}$ and the
undirected $s$-$t$-cut $\left[  S,\overline{S}\right]  _{\operatorname*{und}}$
in $G$ have the same size (because each edge in $\left[  S,\overline
{S}\right]  _{\operatorname*{und}}$ corresponds to exactly one arc in $\left[
S,\overline{S}\right]  $). Thus, the sizes of the directed $s$-$t$-cuts in
$G^{\operatorname*{bidir}}$ are exactly the sizes of the undirected $s$%
-$t$-cuts in $G$. In particular, the minimum size of a former cut equals the
minimum size of a latter cut. This proves Claim 2.}.

It remains to verify Claim 1. The simplest approach is to argue that each path
from $s$ to $t$ in $G^{\operatorname*{bidir}}$ becomes a path from $s$ to $t$
in $G$ (just replace each arc of the path by the corresponding undirected
edge). Unfortunately, this alone does not suffice, since two arc-disjoint
paths in $G^{\operatorname*{bidir}}$ won't necessarily become edge-disjoint
paths in $G$. Here is an example of how this can go wrong (imagine that the
two arcs between $u$ and $v$ come from the same edge of $G$, and the two paths
are marked red and blue):%
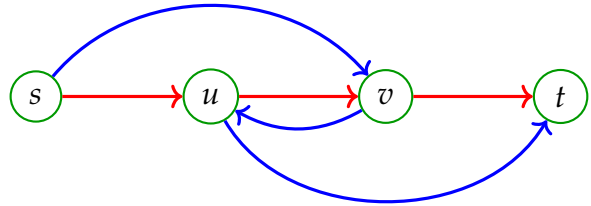
\begin{equation}%
\begin{tikzpicture}[scale=2.3]
\begin{scope}[every node/.style={circle,thick,draw=green!60!black}]
\node(s) at (-1,0) {$s$};
\node(u) at (0,0) {$u$};
\node(v) at (1,0) {$v$};
\node(t) at (2,0) {$t$};
\end{scope}
\begin{scope}[every edge/.style={draw=red,very thick}, every loop/.style={}]
\path[->] (s) edge (u);
\path[->] (u) edge (v);
\path[->] (v) edge (t);
\end{scope}
\begin{scope}[every edge/.style={draw=blue,very thick}, every loop/.style={}]
\path[->] (s) edge[bend left=50] (v);
\path[->] (v) edge[bend left=30] (u);
\path[->] (u) edge[bend right=60] (t);
\end{scope}
\end{tikzpicture}%
\ \ . \label{eq.thm.menger.ue2.badcase}%
\end{equation}
If we replace each arc by the corresponding edge here, then the two paths will
no longer be edge-disjoint (since the edge between $u$ and $v$ will be used by
both paths).

However, this kind of situation can be averted. To do this, we let $k$ be the
maximum number of pairwise arc-disjoint paths from $s$ to $t$ in
$G^{\operatorname*{bidir}}$. We now choose $k$ pairwise arc-disjoint paths
$\mathbf{p}_{1},\mathbf{p}_{2},\ldots,\mathbf{p}_{k}$ from $s$ to $t$ in
$G^{\operatorname*{bidir}}$ in such a way that their \textbf{total length}
(i.e., the sum of the lengths of $\mathbf{p}_{1},\mathbf{p}_{2},\ldots
,\mathbf{p}_{k}$) is \textbf{as small as possible}. Then, it is easy to see
that these paths $\mathbf{p}_{1},\mathbf{p}_{2},\ldots,\mathbf{p}_{k}$ become
edge-disjoint paths in $G$ when we replace each arc by the corresponding edge.
\medskip

\begin{fineprint}
[\textit{Proof:} Assume the contrary. Thus, two of these paths $\mathbf{p}%
_{1},\mathbf{p}_{2},\ldots,\mathbf{p}_{k}$ end up sharing an edge when we
replace each arc by the corresponding edge. Let $\mathbf{p}_{i}$ and
$\mathbf{p}_{j}$ be these two paths (where $i\neq j$, of course). Let $e$ be
the edge that they end up sharing, and let $u$ and $v$ be the two endpoints of
$e$, in the order in which they appear on $\mathbf{p}_{i}$. Hence, the path
$\mathbf{p}_{i}$ uses the edge $e$ (or, more precisely, one of the two arcs of
$G^{\operatorname*{bidir}}$ corresponding to $e$) to get from $u$ to $v$.

Since the paths $\mathbf{p}_{i}$ and $\mathbf{p}_{j}$ are arc-disjoint, they
cannot both use the edge $e$ in the same direction (because this would mean
that $\mathbf{p}_{i}$ and $\mathbf{p}_{j}$ share the same arc of
$G^{\operatorname*{bidir}}$). Hence, the path $\mathbf{p}_{j}$ uses the edge
$e$ to get from $v$ to $u$ (since the path $\mathbf{p}_{i}$ uses the edge $e$
to get from $u$ to $v$). Hence, the paths $\mathbf{p}_{i}$ and $\mathbf{p}%
_{j}$ have the following forms:%
\begin{align*}
\mathbf{p}_{i}  &  =\left(  \ldots,u,e_{1},v,\ldots\right)  ;\\
\mathbf{p}_{j}  &  =\left(  \ldots,v,e_{2},u,\ldots\right)  ,
\end{align*}
where $e_{1}$ and $e_{2}$ are the two arcs of $G^{\operatorname*{bidir}}$ that
correspond to the edge $e$. Now, let us replace the two paths $\mathbf{p}_{i}$
and $\mathbf{p}_{j}$ by two new walks\footnote{Here is an illustration:
\[%

\ \ \ .
\]
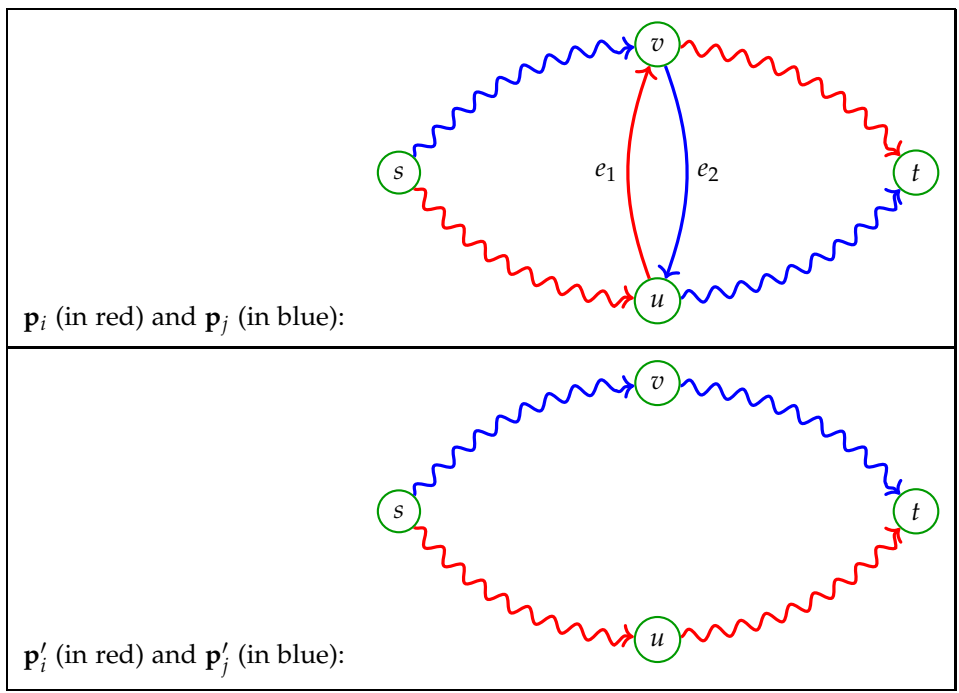
(The wavy arrows stand not for single arcs, but for sequences of multiple
arcs.)}%
\begin{align*}
\mathbf{p}_{i}^{\prime}  &  =\left(  \underbrace{\ldots}_{\substack{\text{the
part of }\mathbf{p}_{i}\\\text{before }u}},u,\underbrace{\ldots}%
_{\substack{\text{the part of }\mathbf{p}_{j}\\\text{after }u}}\right)
\ \ \ \ \ \ \ \ \ \ \text{and}\\
\mathbf{p}_{j}^{\prime}  &  =\left(  \underbrace{\ldots}_{\substack{\text{the
part of }\mathbf{p}_{j}\\\text{before }v}},v,\underbrace{\ldots}%
_{\substack{\text{the part of }\mathbf{p}_{i}\\\text{after }v}}\right)  .
\end{align*}
These walks $\mathbf{p}_{i}^{\prime}$ and $\mathbf{p}_{j}^{\prime}$ are two
walks from $s$ to $t$, and they don't use any arcs that were not already used
by $\mathbf{p}_{i}$ or $\mathbf{p}_{j}$. Thus, they are arc-disjoint from all
of the paths $\mathbf{p}_{1},\mathbf{p}_{2},\ldots,\mathbf{p}_{k}$ except for
$\mathbf{p}_{i}$ and $\mathbf{p}_{j}$. Moreover, they are arc-disjoint from
each other (since $\mathbf{p}_{i}$ and $\mathbf{p}_{j}$ were arc-disjoint, and
since the arcs of any path are distinct). Furthermore, their total length is
smaller by $2$ than the total length of $\mathbf{p}_{i}$ and $\mathbf{p}_{j}$
(since they use all the arcs of $\mathbf{p}_{i}$ and $\mathbf{p}_{j}$ except
for $e_{1}$ and $e_{2}$). They are not necessarily paths, but we can turn them
into paths from $s$ to $t$ by successively removing cycles (as in the proof of
Corollary \ref{cor.mdg.walk-thus-path}). If we do this, we end up with two
paths $\mathbf{p}_{i}^{\prime\prime}$ and $\mathbf{p}_{j}^{\prime\prime}$ from
$s$ to $t$ that are arc-disjoint from each other and from all of the paths
$\mathbf{p}_{1},\mathbf{p}_{2},\ldots,\mathbf{p}_{k}$ except for
$\mathbf{p}_{i}$ and $\mathbf{p}_{j}$, and whose total length is smaller by at
least $2$ than the total length of $\mathbf{p}_{i}$ and $\mathbf{p}_{j}$.

Thus, if we replace $\mathbf{p}_{i}$ and $\mathbf{p}_{j}$ by these two paths
$\mathbf{p}_{i}^{\prime\prime}$ and $\mathbf{p}_{j}^{\prime\prime}$ (while
leaving the remaining $k-2$ of our $k$ paths $\mathbf{p}_{1},\mathbf{p}%
_{2},\ldots,\mathbf{p}_{k}$ unchanged), then we obtain $k$ mutually
arc-disjoint paths from $s$ to $t$ whose total length is smaller than the
total length of our original $k$ paths $\mathbf{p}_{1},\mathbf{p}_{2}%
,\ldots,\mathbf{p}_{k}$. However, this is absurd, because we chose our
original $k$ pairwise arc-disjoint paths $\mathbf{p}_{1},\mathbf{p}_{2}%
,\ldots,\mathbf{p}_{k}$ from $s$ to $t$ in such a way that their total length
is as small as possible. This contradiction shows that our assumption was
wrong. Thus, we have proved that the paths $\mathbf{p}_{1},\mathbf{p}%
_{2},\ldots,\mathbf{p}_{k}$ become edge-disjoint paths in $G$ when we replace
each arc by the corresponding edge.] \medskip
\end{fineprint}

Hence, we have found $k$ pairwise edge-disjoint paths from $s$ to $t$ in $G$
(namely, the $k$ paths that are obtained from the paths $\mathbf{p}%
_{1},\mathbf{p}_{2},\ldots,\mathbf{p}_{k}$ when we replace each arc by the
corresponding edge). This shows that%
\begin{align*}
&  \left(  \text{the maximum number of pairwise edge-disjoint paths from
}s\text{ to }t\text{ in }G\right) \\
&  \geq k\\
&  =\left(  \text{the maximum number of pairwise arc-disjoint paths from
}s\text{ to }t\text{ in }G^{\operatorname*{bidir}}\right)
\end{align*}
(by the definition of $k$). Conversely, we can easily see that%
\begin{align*}
&  \left(  \text{the maximum number of pairwise arc-disjoint paths from
}s\text{ to }t\text{ in }G^{\operatorname*{bidir}}\right) \\
&  \geq\left(  \text{the maximum number of pairwise edge-disjoint paths from
}s\text{ to }t\text{ in }G\right)
\end{align*}
(since there is an obvious way to transform paths in $G$ into paths in
$G^{\operatorname*{bidir}}$ (just replace each edge by one of the two
corresponding arcs of $G^{\operatorname*{bidir}}$), and applying this
transformation to edge-disjoint paths of $G$ yields arc-disjoint paths of
$G^{\operatorname*{bidir}}$). Combining these two inequalities, we obtain%
\begin{align*}
&  \left(  \text{the maximum number of pairwise arc-disjoint paths from
}s\text{ to }t\text{ in }G^{\operatorname*{bidir}}\right) \\
&  =\left(  \text{the maximum number of pairwise edge-disjoint paths from
}s\text{ to }t\text{ in }G\right)  .
\end{align*}
This proves Claim 1. As we explained, this concludes the proof of Theorem
\ref{thm.menger.ue2}.
\end{proof}

\begin{proof}
[Proof of Theorem \ref{thm.menger.ue1}.]This can be derived from Theorem
\ref{thm.menger.ue2} and Remark \ref{rmk.menger.udcut} in the same way as we
derived Theorem \ref{thm.menger.da1} from Theorem \ref{thm.menger.da2} and
Remark \ref{rmk.menger.cut}.
\end{proof}

\begin{exercise}
\label{exe.9.7}Let $G$ be a multigraph such that each vertex of $G$ has even
degree. Let $s$ and $t$ be two distinct vertices of $G$. Prove that the
maximum number of pairwise edge-disjoint paths from $s$ to $t$ is even.
\end{exercise}

\subsubsection{The vertex-Menger theorem for directed graphs}

The Menger theorems we have seen so far have been concerned with paths not
having arcs in common. What if we want to avoid common vertices too?

\begin{definition}
Let $\mathbf{p}$ be a path of some graph or digraph. Then, an
\textbf{intermediate vertex} of $\mathbf{p}$ shall mean a vertex of
$\mathbf{p}$ that is neither the starting point nor the ending point of
$\mathbf{p}$.
\end{definition}

\begin{definition}
Two paths $\mathbf{p}$ and $\mathbf{q}$ in a graph or digraph are said to be
\textbf{internally vertex-disjoint} if they have no common intermediate vertices.
\end{definition}

\begin{example}
The two paths $\mathbf{p}$ and $\mathbf{q}$ in Example \ref{exa.menger.adis.1}
are arc-disjoint, but not internally vertex-disjoint.

Here are two internally vertex-disjoint paths $\mathbf{p}$ and $\mathbf{q}$:%
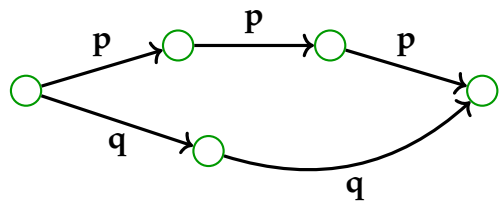
\[%
\begin{tikzpicture}[scale=2]
\begin{scope}[every node/.style={circle,thick,draw=green!60!black}]
\node(s) at (-1,0) {};
\node(1) at (0,0.3) {};
\node(2) at (1,0.3) {};
\node(3) at (0.2,-0.4) {};
\node(t) at (2,0) {};
\end{scope}
\begin{scope}[every edge/.style={draw=black,very thick}, every loop/.style={}]
\path[->] (s) edge node[above] {$\mathbf{p}$} (1);
\path[->] (1) edge node[above] {$\mathbf{p}$} (2);
\path[->] (2) edge node[above] {$\mathbf{p}$} (t);
\path[->] (s) edge node[below] {$\mathbf{q}$} (3);
\path[->] (3) edge[bend right=30] node[below] {$\mathbf{q}$} (t);
\end{scope}
\end{tikzpicture}%
\ \ .
\]

One trivial case of internally vertex-disjoint paths is a path of length
$\leq1$: Namely, a path of length $\leq1$ is internally vertex-disjoint from
any path, including itself (since it has no intermediate vertices).
\end{example}

\begin{definition}
\label{def.menger.ver-sep}Let $D=\left(  V,A,\psi\right)  $ be a multidigraph,
and let $s$ and $t$ be two vertices of $D$. A subset $W$ of $V\setminus
\left\{  s,t\right\}  $ is said to be an \textbf{internal }$s$\textbf{-}%
$t$\textbf{-vertex-separator} if each path from $s$ to $t$ contains at least
one vertex from $W$. Equivalently, a subset $W$ of $V\setminus\left\{
s,t\right\}  $ is said to be an \textbf{internal }$s$\textbf{-}$t$%
\textbf{-vertex-separator} if the induced subdigraph of $D$ on the set
$V\setminus W$ has no path from $s$ to $t$ (in other words, removing from $D$
all vertices contained in $W$ destroys all paths from $s$ to $t$).
\end{definition}

\begin{example}
Let $D=\left(  V,A,\psi\right)  $ be the following multidigraph:%
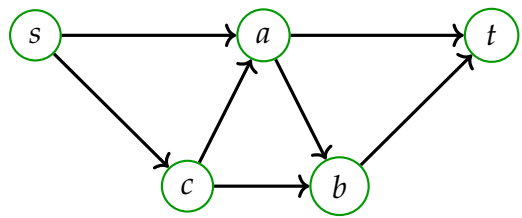
\[%
\begin{tikzpicture}[scale=2]
\begin{scope}[every node/.style={circle,thick,draw=green!60!black}]
\node(s) at (-1,0) {$s$};
\node(a) at (0.5, 0) {$a$};
\node(c) at (0,-1) {$c$};
\node(b) at (1,-1) {$b$};
\node(t) at (2,0) {$t$};
\end{scope}
\begin{scope}[every edge/.style={draw=black,very thick}, every loop/.style={}]
\path[->] (s) edge (c);
\path[->] (c) edge (b);
\path[->] (c) edge (a);
\path[->] (a) edge (t);
\path[->] (a) edge (b);
\path[->] (s) edge (a);
\path[->] (b) edge (t);
\end{scope}
\end{tikzpicture}%
\ \ .
\]
Then, the sets $\left\{  a,b\right\}  $ and $\left\{  a,c\right\}  $ are
internal $s$-$t$-vertex-separators (indeed, removing the vertices $a$ and $b$
cuts off $s$ from the rest of the digraph, whereas removing the vertices $a$
and $c$ does the same to $t$), but the sets $\left\{  a\right\}  $ and
$\left\{  b,c\right\}  $ are not (since the path from $s$ to $t$ via $c$ and
$b$ avoids $a$, whereas the path from $s$ to $t$ via $a$ avoids $b$ and $c$).
\end{example}

\begin{example}
Let $D=\left(  V,A,\psi\right)  $ be a multidigraph. Let $s$ and $t$ be two
distinct vertices of $D$. Then:

\begin{enumerate}
\item[\textbf{(a)}] The empty set $\varnothing$ is an internal $s$%
-$t$-vertex-separator if and only if $D$ has no path from $s$ to $t$.

\item[\textbf{(b)}] If $D$ has no arc with source $s$ and target $t$, then the
set $V\setminus\left\{  s,t\right\}  $ is an internal $s$-$t$-vertex-separator
(since any path from $s$ to $t$ contains at least one intermediate vertex, and
such a vertex must belong to $V\setminus\left\{  s,t\right\}  $).

\item[\textbf{(c)}] If $D$ has an arc with source $s$ and target $t$, then
there exists no internal $s$-$t$-vertex-separator (since the \textquotedblleft
direct\textquotedblright\ length-$1$ path from $s$ to $t$ contains no vertices
besides $s$ and $t$).
\end{enumerate}
\end{example}

Now, we state the analogue of Theorem \ref{thm.menger.da1} and Theorem
\ref{thm.menger.da2} for internally vertex-disjoint paths:

\begin{theorem}
[vertex-Menger theorem for directed graphs]\label{thm.menger.dv2}Let
$D=\left(  V,A,\psi\right)  $ be a multidigraph, and let $s$ and $t$ be two
distinct vertices of $D$. Assume that $D$ has no arc with source $s$ and
target $t$. Then, the maximum number of pairwise internally vertex-disjoint
paths from $s$ to $t$ equals the minimum size of an internal $s$-$t$-vertex-separator.
\end{theorem}

\begin{example}
Let $D$ be the following multidigraph:%
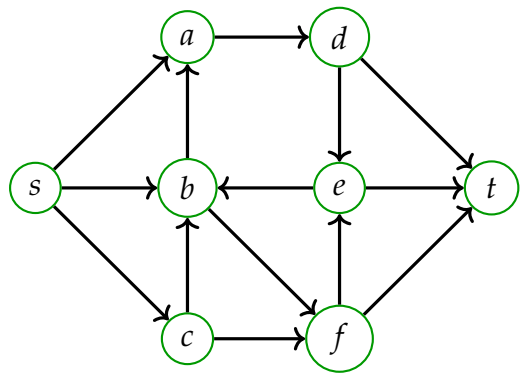
\[%
%
\ \ .
\]
Then, the maximum number of pairwise internally vertex-disjoint paths from $s$
to $t$ is $2$. Indeed, the following picture shows $2$ such paths in red and
blue, respectively:%
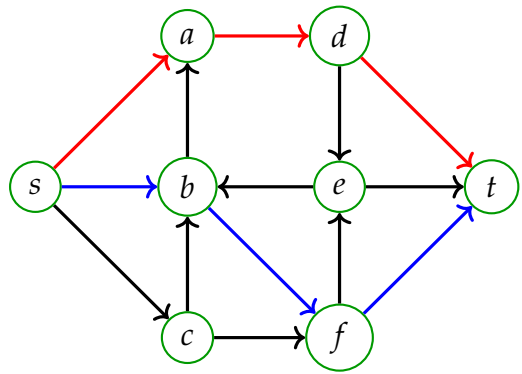
\[%
%
\ \ .
\]
Why can there be no $3$ such paths? This is not obvious from a quick look, but
can be easily derived from Theorem \ref{thm.menger.dv2}. Indeed, Theorem
\ref{thm.menger.dv2} yields that the maximum number of pairwise internally
vertex-disjoint paths from $s$ to $t$ equals the minimum size of an internal
$s$-$t$-vertex-separator. Thus, if the former number was larger than $2$, then
so would be the latter number. But this cannot be the case, since the
$2$-element set $\left\{  a,f\right\}  $ is easily checked to be an internal
$s$-$t$-vertex-separator. Hence, we see that both of these numbers are $2$.
\end{example}

\begin{example}
Consider again the digraph $D$ from Example \ref{exa.menger.da2.1}. In that
example, we found $3$ pairwise arc-disjoint paths from $s$ to $t$. These $3$
paths are not internally vertex-disjoint (in fact, the brown path has
non-starting and non-ending vertices in common with both the red and the blue
path). However, there do exist $3$ pairwise internally vertex-disjoint paths
from $s$ to $t$. Can you find them?
\end{example}

\begin{proof}
[Proof of Theorem \ref{thm.menger.dv2}.]We will again derive this from the
arc-Menger theorem (Theorem \ref{thm.menger.da2}), applied to an appropriate
multidigraph $D^{\prime}=\left(  V^{\prime},A^{\prime},\psi^{\prime}\right)  $.

What is this multidigraph $D^{\prime}$ ? The idea is to modify the digraph $D$
in such a way that paths having a common vertex become paths having a common
arc. The most natural way to achieve this is to \textquotedblleft stretch
out\textquotedblright\ each vertex $v$ of $D$ into a little arc. In order to
do this in a systematic manner, we replace each vertex $v$ of $D$ by two
distinct vertices $v^{i}$ and $v^{o}$ (the notations stand for
\textquotedblleft$v$-in\textquotedblright\ and \textquotedblleft%
$v$-out\textquotedblright, and we can think of $v^{i}$ as the
\textquotedblleft entrance\textquotedblright\ to $v$ while $v^{o}$ is the
\textquotedblleft exit\textquotedblright\ from $v$) and an arc $v^{io}$ that
goes from $v^{i}$ to $v^{o}$. Any existing arc $a$ of $D$ becomes a new arc
$a^{oi}$ of $D$, whose source and target are specified as follows: If $a$ has
source $u$ and target $v$, then $a^{oi}$ will have source $u^{o}$ and target
$v^{i}$.

Here is an example: If
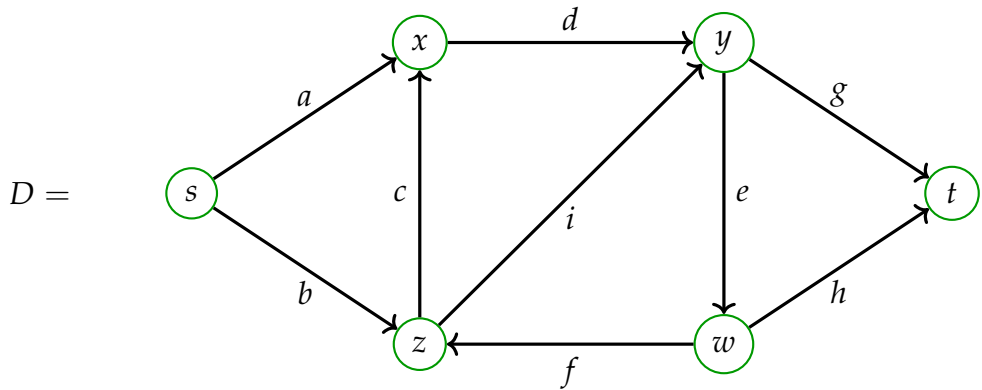
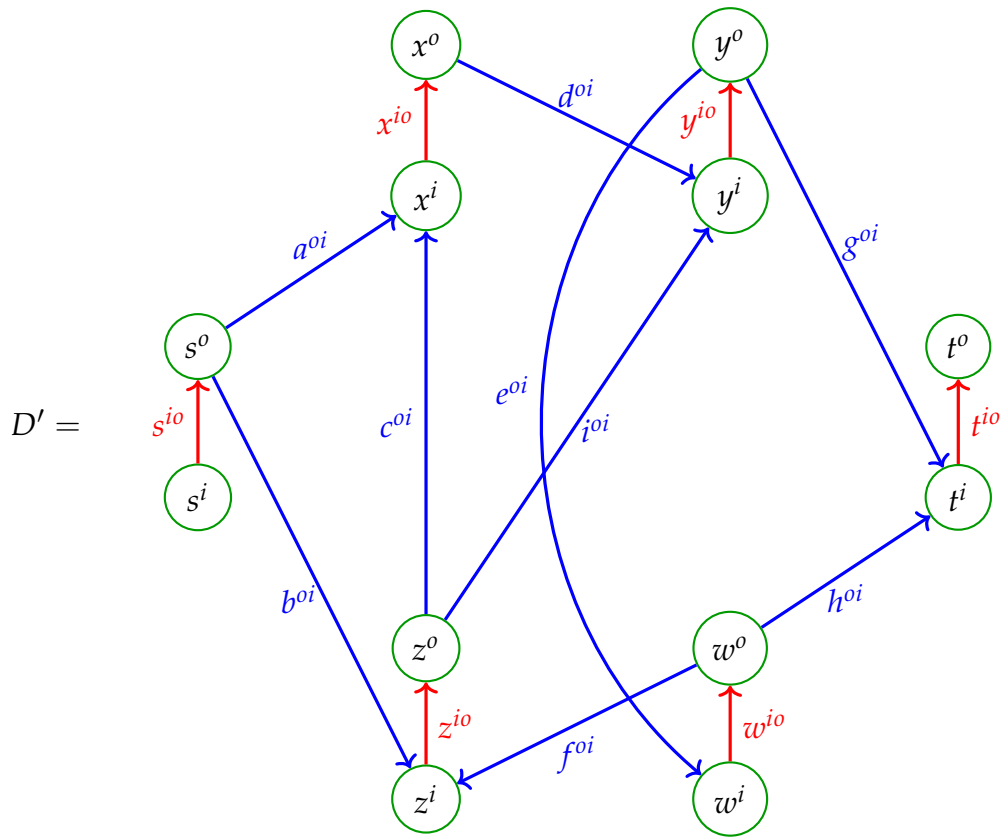
\[%
%
\]
(where all arcs of the form $a^{oi}$ for $a\in A$ have been colored blue,
whereas all arcs of the form $v^{io}$ for $v\in V$ have been colored red).
This $D^{\prime}$ satisfies the property that we want it to satisfy: For
instance, the two paths%
\begin{align*}
&  \left(  s,a,x,d,y,g,t\right)  \ \ \ \ \ \ \ \ \ \ \text{and}\\
&  \left(  s,b,z,i,y,e,w,h,t\right)
\end{align*}
of $D$ have the vertex $y$ in common, so the corresponding two paths%
\begin{align*}
&  \left(  s^{o},a^{oi},x^{i},x^{io},x^{o},d^{oi},y^{i},y^{io},y^{o}%
,g^{oi},t^{i}\right)  \ \ \ \ \ \ \ \ \ \ \text{and}\\
&  \left(  s^{o},b^{oi},z^{i},z^{io},z^{o},i^{oi},y^{i},y^{io},y^{o}%
,e^{oi},w^{i},w^{io},w^{o},h^{oi},t^{i}\right)
\end{align*}
of $D^{\prime}$ have the arc $y^{io}$ in common. If you think of $D$ as a
railway network with the vertices being train stations and the arcs being
train rides, then $D^{\prime}$ is a more detailed version of this network that
records a change of platforms as an arc as well.

Here is a formal definition of the multidigraph $D^{\prime}=\left(  V^{\prime
},A^{\prime},\psi^{\prime}\right)  $ in full generality:

\begin{itemize}
\item We replace each vertex $v$ of $D$ by two new vertices $v^{i}$ and
$v^{o}$. We call $v^{i}$ an \textquotedblleft\textbf{in-vertex}%
\textquotedblright\ and $v^{o}$ an \textquotedblleft\textbf{out-vertex}%
\textquotedblright. The vertex set of $D^{\prime}$ will be the set%
\[
V^{\prime}:=\underbrace{\left\{  v^{i}\ \mid\ v\in V\right\}  }%
_{\text{in-vertices}}\cup\underbrace{\left\{  v^{o}\ \mid\ v\in V\right\}
}_{\text{out-vertices}}.
\]

\item Each arc $a\in A$ is replaced by a new arc $a^{oi}$, which is defined as
follows: If the arc $a\in A$ has source $u$ and target $v$, then we replace it
by a new arc $a^{oi}$, which has source $u^{o}$ and target $v^{i}$. This arc
$a^{oi}$ will be called an \textquotedblleft\textbf{arc-arc}\textquotedblright%
\ of $D^{\prime}$ (since it originates from an arc of $D$).

\item For any vertex $v\in V$ of $D$, we introduce a new arc $v^{io}$, which
has source $v^{i}$ and target $v^{o}$. This arc $v^{io}$ will be called a
\textquotedblleft\textbf{vertex-arc}\textquotedblright\ of $D^{\prime}$ (since
it originates from a vertex of $D$).

\item The arc set of $D^{\prime}$ will be the set%
\[
A^{\prime}:=\underbrace{\left\{  a^{oi}\ \mid\ a\in A\right\}  }_{\text{the
arc-arcs}}\cup\underbrace{\left\{  v^{io}\ \mid\ v\in V\right\}  }_{\text{the
vertex-arcs}}.
\]
The map $\psi^{\prime}:A^{\prime}\rightarrow V^{\prime}\times V^{\prime}$ is
defined as we already explained:

\begin{itemize}
\item For any arc-arc $a^{oi}\in A$, we let $\psi^{\prime}\left(
a^{oi}\right)  :=\left(  u^{o},v^{i}\right)  $, where $\left(  u,v\right)
=\psi\left(  a\right)  $.

\item For any vertex-arc $v^{io}\in V$, we let $\psi^{\prime}\left(
v^{io}\right)  :=\left(  v^{i},v^{o}\right)  $.
\end{itemize}
\end{itemize}

Note that $D^{\prime}$ is something like a \textquotedblleft bipartite
digraph\textquotedblright: Each of its arcs goes either from an out-vertex to
an in-vertex or vice versa. Namely, each arc-arc goes from an out-vertex to an
in-vertex, whereas each vertex-arc goes from an in-vertex to an out-vertex.
Thus, on any walk of $D^{\prime}$, the arc-arcs and the vertex-arcs have to alternate.

If $\mathbf{p}=\left(  v_{0},a_{1},v_{1},a_{2},v_{2},\ldots,a_{k}%
,v_{k}\right)  $ is any nontrivial\footnote{We say that a path is
\textbf{nontrivial} if it has length $>0$.} path of $D$, then we can define a
corresponding path $\mathbf{p}^{oi}$ of $D^{\prime}$ by%
\[
\mathbf{p}^{oi}:=\left(  v_{0}^{o},a_{1}^{oi},v_{1}^{i},v_{1}^{io},v_{1}%
^{o},a_{2}^{oi},v_{2}^{i},v_{2}^{io},v_{2}^{o},\ldots,a_{k}^{oi},v_{k}%
^{i}\right)  .
\]
This path $\mathbf{p}^{oi}$ is obtained from $\mathbf{p}$ by

\begin{itemize}
\item replacing the starting point $v_{0}$ by $v_{0}^{o}$;

\item replacing the ending point $v_{k}$ by $v_{k}^{i}$;

\item replacing each other vertex $v_{j}$ by the sequence $v_{j}^{i}%
,v_{j}^{io},v_{j}^{o}$;

\item replacing each arc $a_{j}$ by $a_{j}^{oi}$.
\end{itemize}

Informally speaking, this simply means that we stretch out each intermediate
vertex of $\mathbf{p}$ to the corresponding arc.

If $\mathbf{p}$ is a path from $s$ to $t$ in $D$, then $\mathbf{p}^{oi}$ is a
path from $s^{o}$ to $t^{i}$ in $D^{\prime}$. Conversely, any path from
$s^{o}$ to $t^{i}$ in $D^{\prime}$ must have the form $\mathbf{p}^{oi}$, where
$\mathbf{p}$ is some path from $s$ to $t$ in $D$ (because on any walk of
$D^{\prime}$, the arc-arcs and the vertex-arcs have to alternate). Therefore,
the map%
\begin{align}
\left\{  \text{paths from }s\text{ to }t\text{ in }D\right\}   &
\rightarrow\left\{  \text{paths from }s^{o}\text{ to }t^{i}\text{ in
}D^{\prime}\right\}  ,\nonumber\\
\mathbf{p}  &  \mapsto\mathbf{p}^{oi} \label{pf.thm.menger.dv2.bij}%
\end{align}
is a bijection. Moreover, two paths $\mathbf{p}$ and $\mathbf{q}$ of $D$ are
internally vertex-disjoint if and only if the paths $\mathbf{p}^{oi}$ and
$\mathbf{q}^{oi}$ are arc-disjoint (since each vertex of a path $\mathbf{p}$
except for its starting and ending points is represented by an arc in
$\mathbf{p}^{oi}$).

Now, let $k$ be the maximum number of pairwise arc-disjoint paths from $s^{o}$
to $t^{i}$ in $D^{\prime}$. Thus, $D^{\prime}$ has $k$ pairwise arc-disjoint
paths from $s^{o}$ to $t^{i}$. Applying the inverse of the bijection
(\ref{pf.thm.menger.dv2.bij}) to these $k$ paths, we obtain $k$ pairwise
internally vertex-disjoint paths from $s$ to $t$ in $D$ (because two paths
$\mathbf{p}$ and $\mathbf{q}$ of $D$ are internally vertex-disjoint if and
only if the paths $\mathbf{p}^{oi}$ and $\mathbf{q}^{oi}$ are arc-disjoint).
Hence,%
\begin{align}
&  \left(  \text{the maximum number of pairwise internally vertex-disjoint
}\right. \nonumber\\
&  \ \ \ \ \ \ \ \ \ \ \ \ \ \ \ \ \ \ \ \ \left.  \text{paths from }s\text{
to }t\text{ in }D\right) \nonumber\\
&  \geq k. \label{pf.thm.menger.dv2.max-geq-k}%
\end{align}

Our next goal is to find an internal $s$-$t$-vertex-separator $W\subseteq
V\setminus\left\{  s,t\right\}  $ of size $\left\vert W\right\vert \leq k$.

First, we simplify our setting a bit.

A path from $s$ to $t$ cannot contain a loop; nor can it contain an arc with
source $t$ and target $s$ (since the vertices of a path must be distinct).
Hence, we can remove such arcs (i.e., loops as well as arcs with source $t$
and target $s$) from $D$ without affecting the meaning of the claim we are
proving. Thus, we WLOG assume that the digraph $D$ has no such arcs. Since we
also know (by assumption) that $D$ has no arc with source $s$ and target $t$,
we thus conclude that $D$ has no arc with source $\in\left\{  s,t\right\}  $
and target $\in\left\{  s,t\right\}  $ (because each such arc would either
have source $s$ and target $t$, or have source $t$ and target $s$, or be a
loop). In other words, each arc of $D$ has at least one endpoint\footnote{An
\textbf{endpoint} of an arc means a vertex that is either the source or the
target of this arc.} distinct from both $s$ and $t$.

However, $k$ is the maximum number of pairwise arc-disjoint paths from $s^{o}$
to $t^{i}$ in $D^{\prime}$. Therefore, by Theorem \ref{thm.menger.da2}
(applied to $D^{\prime}=\left(  V^{\prime},A^{\prime},\psi^{\prime}\right)  $,
$s^{o}$ and $t^{i}$ instead of $D=\left(  V,A,\psi\right)  $, $s$ and $t$),
this number $k$ equals the minimum size of an $s^{o}$-$t^{i}$-cut in
$D^{\prime}$. Hence, there exists an $s^{o}$-$t^{i}$-cut $\left[
S,\overline{S}\right]  $ in $D^{\prime}$ such that $\left\vert \left[
S,\overline{S}\right]  \right\vert =k$. Consider this $s^{o}$-$t^{i}$-cut
$\left[  S,\overline{S}\right]  $. Since $\left[  S,\overline{S}\right]  $ is
an $s^{o}$-$t^{i}$-cut, we have $S\subseteq V^{\prime}$ and $s^{o}\in S$ and
$t^{i}\notin S$.

Let $B:=\left[  S,\overline{S}\right]  $. Then, $\left\vert B\right\vert
=\left\vert \left[  S,\overline{S}\right]  \right\vert =k$. Moreover, it is
easy to see that $s^{io}\notin B$\ \ \ \ \footnote{\textit{Proof.} If we had
$s^{io}\in\left[  S,\overline{S}\right]  $, then we would have $s^{i}\in S$
and $s^{o}\in\overline{S}$; however, $s^{o}\in\overline{S}$ would contradict
$s^{o}\in S$. Thus, we cannot have $s^{io}\in\left[  S,\overline{S}\right]  $.
In other words, we cannot have $s^{io}\in B$ (since $B=\left[  S,\overline
{S}\right]  $). Hence, $s^{io}\notin B$.} and $t^{io}\notin B$%
\ \ \ \ \footnote{\textit{Proof.} If we had $t^{io}\in\left[  S,\overline
{S}\right]  $, then we would have $t^{i}\in S$ and $t^{o}\in\overline{S}$;
however, $t^{i}\in S$ would contradict $t^{i}\in\overline{S}$. Thus, we cannot
have $t^{io}\in\left[  S,\overline{S}\right]  $. In other words, we cannot
have $t^{io}\in B$ (since $B=\left[  S,\overline{S}\right]  $). Hence,
$t^{io}\notin B$.}.

To each vertex $w\in V^{\prime}$ of $D^{\prime}$, we assign a vertex
$\beta\left(  w\right)  \in V$ of $D$ as follows: If $w=v^{i}$ or $w=v^{o}$
for some $v\in V$, then we set $\beta\left(  w\right)  :=v$. In other words,
$\beta\left(  w\right)  $ is the vertex $v$ such that $w\in\left\{
v^{i},v^{o}\right\}  $. We shall call $\beta\left(  w\right)  $ the
\textbf{base} of the vertex $w$.

For each arc $b\in B$, there exists at least one endpoint $w$ of $b$ such that
$\beta\left(  w\right)  \in V\setminus\left\{  s,t\right\}  $%
\ \ \ \ \footnote{\textit{Proof:} Let $b\in B$ be an arc. We must prove that
there exists at least one endpoint $w$ of $b$ such that $\beta\left(
w\right)  \in V\setminus\left\{  s,t\right\}  $.
\par
The arc $b$ is either a vertex-arc or an arc-arc. Thus, we are in one of the
following two cases:
\par
\textit{Case 1:} The arc $b$ is a vertex-arc.
\par
\textit{Case 2:} The arc $b$ is an arc-arc.
\par
Let us first consider Case 1. In this case, the arc $b$ is a vertex-arc. In
other words, $b=v^{io}$ for some $v\in V$. Consider this $v$. Then,
$v^{io}=b\in B$. Hence, $v\neq s$ (since $v=s$ would entail $v^{io}%
=s^{io}\notin B$, which would contradict $v^{io}\in B$) and $v\neq t$ (since
$v=t$ would entail $v^{io}=t^{io}\notin B$, which would contradict $v^{io}\in
B$). Therefore, $v\in V\setminus\left\{  s,t\right\}  $. Also, clearly,
$v^{i}$ is an endpoint of $b$ and satisfies $\beta\left(  v^{i}\right)  =v\in
V\setminus\left\{  s,t\right\}  $. Hence, there exists at least one endpoint
$w$ of $b$ such that $\beta\left(  w\right)  \in V\setminus\left\{
s,t\right\}  $ (namely, $w=v^{i}$). Thus, our proof is complete in Case 1.
\par
Let us now consider Case 2. In this case, the arc $b$ is an arc-arc. In other
words, $b=a^{oi}$ for some $a\in A$. Consider this $a$. Now, $a$ is an arc of
$D$ (since $a\in A$), and thus has at least one endpoint distinct from both
$s$ and $t$ (since we have shown above that each arc of $D$ has at least one
endpoint distinct from both $s$ and $t$). Let $x$ be this endpoint. Then,
$x\in V\setminus\left\{  s,t\right\}  $ (since $x$ is distinct from both $s$
and $t$).
\par
But $x$ is an endpoint of $a$. In other words, $x$ is either the source or the
target of $a$. Hence, the arc $a^{oi}$ of $D^{\prime}$ either has source
$x^{o}$ or has target $x^{i}$ (by the definition of $a^{oi}$). In other words,
the arc $b$ of $D^{\prime}$ either has source $x^{o}$ or has target $x^{i}$
(since $b=a^{oi}$). Since $\beta\left(  x^{o}\right)  =x\in V\setminus\left\{
s,t\right\}  $ and $\beta\left(  x^{i}\right)  =x\in V\setminus\left\{
s,t\right\}  $, we thus conclude that the arc $b$ of $D^{\prime}$ has at least
one endpoint $w$ such that $\beta\left(  w\right)  \in V\setminus\left\{
s,t\right\}  $ (namely, $w=x^{o}$ if $b$ has source $x^{o}$, and $w=x^{i}$ if
$b$ has target $x^{i}$). This completes our proof in Case 2.
\par
Thus, we are done in both Cases 1 and 2, so that our proof is complete.}. We
choose such an endpoint $w$ arbitrarily, and we denote its base $\beta\left(
w\right)  $ by $\beta\left(  b\right)  $. We shall call $\beta\left(
b\right)  $ the \textbf{basepoint} of the arc $b$. Thus, by definition, we
have%
\begin{equation}
\beta\left(  b\right)  \in V\setminus\left\{  s,t\right\}
\ \ \ \ \ \ \ \ \ \ \text{for each }b\in B. \label{pf.thm.menger.dv2.betabin}%
\end{equation}

We let $\beta\left(  B\right)  $ denote the set $\left\{  \beta\left(
b\right)  \ \mid\ b\in B\right\}  $. Clearly, $\left\vert \beta\left(
B\right)  \right\vert \leq\left\vert B\right\vert =k$.

Now, we claim that
\begin{equation}
\text{every path from }s\text{ to }t\text{ (in }D\text{) contains a vertex in
}\beta\left(  B\right)  . \label{pf.thm.menger.dv2.base-path}%
\end{equation}

\begin{fineprint}
[\textit{Proof of (\ref{pf.thm.menger.dv2.base-path}):} Let $\mathbf{p}$ be a
path from $s$ to $t$ (in $D$). We must prove that $\mathbf{p}$ contains a
vertex in $\beta\left(  B\right)  $.

Recall that we have assigned a path $\mathbf{p}^{oi}$ of $D^{\prime}$ to the
path $\mathbf{p}$ of $D$. The definition of $\mathbf{p}^{oi}$ shows that the
base of any vertex of $\mathbf{p}^{oi}$ is a vertex of $\mathbf{p}$ (indeed,
if $v_{0},v_{1},\ldots,v_{k}$ are the vertices of $\mathbf{p}$, then the
vertices of $\mathbf{p}^{oi}$ are $v_{0}^{o},v_{1}^{i},v_{1}^{o},v_{2}%
^{i},v_{2}^{o},\ldots,v_{k-1}^{i},v_{k-1}^{o},v_{k}^{i}$, and their respective
bases are $v_{0},v_{1},v_{1},v_{2},v_{2},\ldots,v_{k-1},v_{k-1},v_{k}$).

The path $\mathbf{p}^{oi}$ is a path from $s^{o}$ to $t^{i}$ (since
$\mathbf{p}$ is a path from $s$ to $t$). Hence, it starts at a vertex in $S$
(since $s^{o}\in S$) and ends at a vertex that is not in $S$ (since
$t^{i}\notin S$). Thus, this path $\mathbf{p}^{oi}$ must cross from $S$ into
$\overline{S}$ somewhere. In other words, there exists an arc $b$ of
$\mathbf{p}^{oi}$ such that the source of $b$ belongs to $S$ but the target of
$b$ belongs to $\overline{S}$. Consider this arc $b$. Thus, $b\in\left[
S,\overline{S}\right]  =B$, so that $\beta\left(  b\right)  \in\beta\left(
B\right)  $ (by the definition of $\beta\left(  B\right)  $). Both endpoints
of $b$ are vertices of $\mathbf{p}^{oi}$ (since $b$ is an arc of
$\mathbf{p}^{oi}$).

Now, consider the basepoint $\beta\left(  b\right)  $ of this arc $b$. This
basepoint $\beta\left(  b\right)  $ is the base of an endpoint of $b$ (by the
definition of $\beta\left(  b\right)  $). Thus, $\beta\left(  b\right)  $ is
the base of a vertex of $\mathbf{p}^{oi}$ (since both endpoints of $b$ are
vertices of $\mathbf{p}^{oi}$). Hence, $\beta\left(  b\right)  $ is a vertex
of $\mathbf{p}$ (since the base of any vertex of $\mathbf{p}^{oi}$ is a vertex
of $\mathbf{p}$). In other words, the path $\mathbf{p}$ contains the vertex
$\beta\left(  b\right)  $. Since $\beta\left(  b\right)  \in\beta\left(
B\right)  $, we thus conclude that $\mathbf{p}$ contains a vertex in
$\beta\left(  B\right)  $. This proves (\ref{pf.thm.menger.dv2.base-path}).]
\medskip
\end{fineprint}

Now, the set $\beta\left(  B\right)  $ is a subset of $V\setminus\left\{
s,t\right\}  $ (since $\beta\left(  b\right)  \in V\setminus\left\{
s,t\right\}  $ for each $b\in B$) and has the property that every path from
$s$ to $t$ contains a vertex in $\beta\left(  B\right)  $ (by
(\ref{pf.thm.menger.dv2.base-path})). In other words, $\beta\left(  B\right)
$ is a subset $W\subseteq V\setminus\left\{  s,t\right\}  $ such that every
path from $s$ to $t$ contains a vertex in $W$. In other words, $\beta\left(
B\right)  $ is an internal $s$-$t$-vertex-separator (by the definition of an
\textquotedblleft internal $s$-$t$-vertex-separator\textquotedblright). Thus,
\begin{align*}
&  \left(  \text{the minimum size of an internal }s\text{-}%
t\text{-vertex-separator}\right) \\
&  \leq\left\vert \beta\left(  B\right)  \right\vert =k\\
&  \leq\left(  \text{the maximum number of pairwise internally vertex-disjoint
}\right. \\
&  \ \ \ \ \ \ \ \ \ \ \ \ \ \ \ \ \ \ \ \ \left.  \text{paths from }s\text{
to }t\text{ in }D\right)
\end{align*}
(by (\ref{pf.thm.menger.dv2.max-geq-k})).

On the other hand, we have
\begin{align*}
&  \left(  \text{the minimum size of an internal }s\text{-}%
t\text{-vertex-separator}\right) \\
&  \geq\left(  \text{the maximum number of pairwise internally
vertex-disjoint}\right. \\
&  \ \ \ \ \ \ \ \ \ \ \ \ \ \ \ \ \ \ \ \ \left.  \text{paths from }s\text{
to }t\text{ in }D\right)
\end{align*}
(by the pigeonhole principle\footnote{\textit{Proof} in more detail\textit{:}
Let $n$ be the minimum size of an internal $s$-$t$-vertex-separator. Let $x$
be the maximum number of pairwise internally vertex-disjoint paths from $s$ to
$t$ in $D$. We must show that $n\geq x$.
\par
Assume the contrary. Thus, $n<x$.
\par
The definition of $n$ shows that there exists an internal $s$-$t$%
-vertex-separator $W$ that has size $n$.
\par
The set $W$ is an internal $s$-$t$-vertex-separator. In other words, $W$ is a
subset of $V\setminus\left\{  s,t\right\}  $ such that every path from $s$ to
$t$ contains a vertex in $W$. Moreover, $W$ has size $n$; thus, $\left\vert
W\right\vert =n<x$.
\par
The definition of $x$ shows that there exist $x$ pairwise internally
vertex-disjoint paths from $s$ to $t$ in $D$. Let $\mathbf{p}_{1}%
,\mathbf{p}_{2},\ldots,\mathbf{p}_{x}$ be these $x$ paths. Each of these $x$
paths $\mathbf{p}_{1},\mathbf{p}_{2},\ldots,\mathbf{p}_{x}$ must contain at
least one vertex in $W$ (since every path from $s$ to $t$ contains a vertex in
$W$). Since $\left\vert W\right\vert <x$, we thus conclude by the pigeonhole
principle that at least two of the $x$ paths $\mathbf{p}_{1},\mathbf{p}%
_{2},\ldots,\mathbf{p}_{x}$ must contain the same vertex in $W$. In other
words, there exist two distinct elements $i,j\in\left\{  1,2,\ldots,x\right\}
$ such that $\mathbf{p}_{i}$ and $\mathbf{p}_{j}$ contain the same vertex in
$W$. Let $w$ be the latter vertex. Thus, $w\in W\subseteq V\setminus\left\{
s,t\right\}  $. Hence, $w$ is distinct from both $s$ and $t$. Therefore, $w$
is an intermediate vertex of $\mathbf{p}_{i}$ (since the path $\mathbf{p}_{i}$
has starting point $s$ and ending point $t$). Likewise, $w$ is an intermediate
vertex of $\mathbf{p}_{j}$.
\par
However, the paths $\mathbf{p}_{i}$ and $\mathbf{p}_{j}$ are internally
vertex-disjoint, and thus have no common intermediate vertex. This contradicts
the fact that $w$ is an intermediate vertex of both paths $\mathbf{p}_{i}$ and
$\mathbf{p}_{j}$. This contradiction shows that our assumption was false.
Hence, $n\geq x$ is proved, qed.}). Combining this inequality with the
preceding one, we obtain
\begin{align*}
&  \left(  \text{the minimum size of an internal }s\text{-}%
t\text{-vertex-separator}\right) \\
&  =\left(  \text{the maximum number of pairwise internally vertex-disjoint}%
\right. \\
&  \ \ \ \ \ \ \ \ \ \ \ \ \ \ \ \ \ \ \ \ \left.  \text{paths from }s\text{
to }t\text{ in }D\right)  .
\end{align*}
This proves Theorem \ref{thm.menger.dv2}.
\end{proof}

\bigskip

There is also a variant of the vertex-Menger theorem similar to what Theorem
\ref{thm.menger.dam} did for the arc-Menger theorem. Again, we need some
notations first:

\begin{definition}
Let $D=\left(  V,A,\psi\right)  $ be a multidigraph, and let $X$ and $Y$ be
two subsets of $V$.

\begin{enumerate}
\item[\textbf{(a)}] A \textbf{path from }$X$\textbf{ to }$Y$ shall mean a path
whose starting point belongs to $X$ and whose ending point belongs to $Y$.

\item[\textbf{(b)}] A subset $W$ of $V$ is said to be an $X$\textbf{-}%
$Y$\textbf{-vertex-separator} if each path from $X$ to $Y$ contains at least
one vertex from $W$. Equivalently, a subset $W$ of $V$ is said to be an
$X$\textbf{-}$Y$\textbf{-vertex-separator} if the induced subdigraph of $D$ on
the set $V\setminus W$ has no path from $X$ to $Y$ (in other words, removing
from $D$ all vertices contained in $W$ destroys all paths from $X$ to $Y$).

\item[\textbf{(c)}] An $X$-$Y$-vertex-separator $W$ is said to be
\textbf{internal} if it is a subset of $V\setminus\left(  X\cup Y\right)  $
(that is, if it is disjoint from $X$ and from $Y$).
\end{enumerate}
\end{definition}

\begin{theorem}
[vertex-Menger theorem for directed graphs, multi-terminal version
1]\label{thm.menger.dvm}Let $D=\left(  V,A,\psi\right)  $ be a multidigraph,
and let $X$ and $Y$ be two disjoint subsets of $V$. Assume that $D$ has no arc
with source in $X$ and target in $Y$.

Then, the maximum number of pairwise internally vertex-disjoint paths from $X$
to $Y$ equals the minimum size of an internal $X$-$Y$-vertex-separator.
\end{theorem}

\begin{example}
Let $D$ be the following multidigraph:%
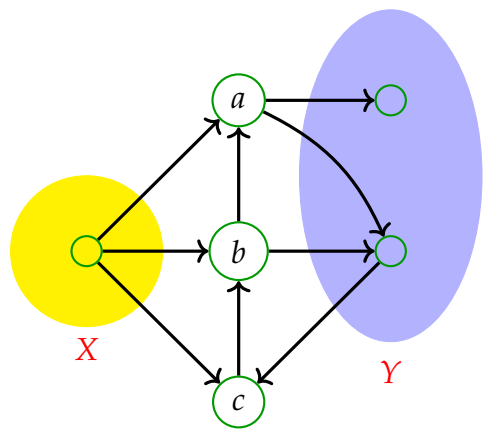
\[%
%
\ \ .
\]
Then, the maximum number of pairwise internally vertex-disjoint paths from $X$
to $Y$ is $2$; here are two such paths (drawn in red and blue):%
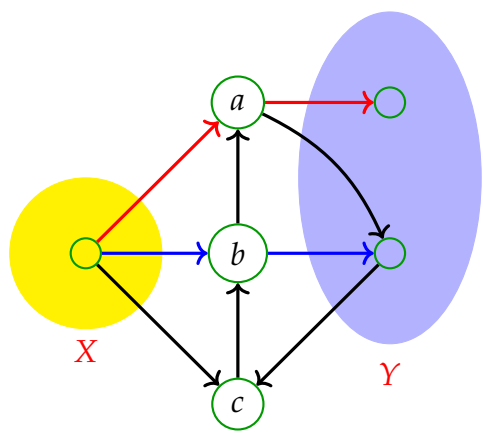
\[%
%
\]
(there are other choices, of course). The minimum size of an internal $X$%
-$Y$-vertex-separator is $2$ as well; indeed, $\left\{  a,b\right\}  $ is such
an internal $X$-$Y$-vertex-separator. These two numbers are equal, just as
Theorem \ref{thm.menger.dvm} predicts.
\end{example}

\begin{proof}
[Proof of Theorem \ref{thm.menger.dvm}.]We define a new multidigraph
$D^{\prime}=\left(  V^{\prime},A^{\prime},\psi^{\prime}\right)  $ as in the
proof of Theorem \ref{thm.menger.dam}. Then, $D^{\prime}$ has no arc with
source $s$ and target $t$ (since $D$ has no arc with source in $X$ and target
in $Y$).

Hence, Theorem \ref{thm.menger.dv2} (applied to $D^{\prime}=\left(  V^{\prime
},A^{\prime},\psi^{\prime}\right)  $ instead of $D=\left(  V,A,\psi\right)  $)
shows that the maximum number of pairwise internally vertex-disjoint paths
from $s$ to $t$ in $D^{\prime}$ equals the minimum size of an internal $s$%
-$t$-vertex-separator in $D^{\prime}$.

Let us now see what this result means for our original digraph $D$. Indeed:

\begin{itemize}
\item The minimum size of an internal $s$-$t$-vertex-separator in $D^{\prime}$
equals the minimum size of an internal $X$-$Y$-vertex-separator in $D$
(indeed, the internal $s$-$t$-vertex-separators in $D^{\prime}$ are precisely
the internal $X$-$Y$-vertex-separators in $D$\ \ \ \ \footnote{\textit{Proof.}
We recall the definitions of internal $s$-$t$-vertex-separators in $D^{\prime
}$ and of internal $X$-$Y$-vertex-separators in $D$:
\par
\begin{itemize}
\item An internal $s$-$t$-vertex-separator in $D^{\prime}$ is a subset $W$ of
$V^{\prime}\setminus\left\{  s,t\right\}  $ such that each path from $s$ to
$t$ contains at least one vertex from $W$.
\par
\item An internal $X$-$Y$-vertex-separator in $D$ is a subset $W$ of
$V\setminus\left(  X\cup Y\right)  $ such that each path from $X$ to $Y$
contains at least one vertex from $W$.
\end{itemize}
\par
These two definitions describe the same object, because of the following two
reasons:
\par
\begin{itemize}
\item We have $V^{\prime}\setminus\left\{  s,t\right\}  =V\setminus\left(
X\cup Y\right)  $.
\par
\item The paths from $s$ to $t$ are in bijection with the paths from $X$ to
$Y$ (indeed, any path of the latter kind can be transformed into a path of the
former kind by replacing the starting point by $s$ and replacing the ending
point by $t$). This bijection preserves the intermediate vertices (i.e., the
vertices other than the starting point and the ending point). Thus, a path
$\mathbf{p}$ from $s$ to $t$ contains at least one vertex from $W$ if and only
if the corresponding path from $X$ to $Y$ (that is, the image of $\mathbf{p}$
under our bijection) contains at least one vertex from $W$.
\end{itemize}
\par
Thus, the internal $s$-$t$-vertex-separators in $D^{\prime}$ are precisely the
internal $X$-$Y$-vertex-separators in $D$.}).

\item The maximum number of pairwise internally vertex-disjoint paths from $s$
to $t$ in $D^{\prime}$ equals the maximum number of pairwise internally
vertex-disjoint paths from $X$ to $Y$ in $D$\ \ \ \ \footnote{\textit{Proof.}
We make the following two observations:
\par
\begin{statement}
\textit{Observation 1:} Let $k\in\mathbb{N}$. If $D^{\prime}$ has $k$ pairwise
internally vertex-disjoint paths from $s$ to $t$, then $D$ has $k$ pairwise
internally vertex-disjoint paths from $X$ to $Y$.
\end{statement}
\par
[\textit{Proof of Observation 1:} Assume that $D^{\prime}$ has $k$ pairwise
internally vertex-disjoint paths from $s$ to $t$. We can \textquotedblleft
lift\textquotedblright\ these $k$ paths to $k$ paths from $X$ to $Y$ in $D$
(preserving the arcs, and replacing the vertices $s$ and $t$ by appropriate
vertices in $X$ and $Y$ to make them belong to the right arcs). The resulting
$k$ paths from $X$ to $Y$ in $D$ are still pairwise internally vertex-disjoint
(since our \textquotedblleft lifting\textquotedblright\ operation has not
changed the intermediate vertices of our paths). Thus, $D$ has $k$ pairwise
internally vertex-disjoint paths from $X$ to $Y$. This proves Observation 1.]
\par
\begin{statement}
\textit{Observation 2:} Let $k\in\mathbb{N}$. If $D$ has $k$ pairwise
internally vertex-disjoint paths from $X$ to $Y$, then $D^{\prime}$ has $k$
pairwise internally vertex-disjoint paths from $s$ to $t$.
\end{statement}
\par
[\textit{Proof of Observation 2:} Assume that $D$ has $k$ pairwise internally
vertex-disjoint paths from $X$ to $Y$. We can replace these $k$ paths by $k$
pairwise internally vertex-disjoint \textbf{walks} from $s$ to $t$ in
$D^{\prime}$ (by replacing their starting points with $s$ and replacing their
ending points with $t$). Thus, $D^{\prime}$ has $k$ pairwise internally
vertex-disjoint \textbf{walks} from $s$ to $t$. Therefore, $D^{\prime}$ has
$k$ pairwise internally vertex-disjoint \textbf{paths} from $s$ to $t$ as well
(since each walk contains a path, and of course we don't lose internal
vertex-disjointness if we restrict our walk to a path contained in it). This
proves Observation 2.]
\par
Observation 2 shows that the maximum number of pairwise internally
vertex-disjoint paths from $s$ to $t$ in $D^{\prime}$ is $\geq$ to the maximum
number of pairwise internally vertex-disjoint paths from $X$ to $Y$ in $D$.
But Observation 1 shows the reverse inequality (i.e., it shows that the former
number is $\leq$ to the latter number). Thus, the inequality is an equality,
i.e., the two numbers are equal. Qed.}.
\end{itemize}

\noindent Hence, the result of the preceding paragraph is precisely the claim
of Theorem \ref{thm.menger.dvm}, and our proof is thus complete.
\end{proof}

Another variant of this result can be stated for vertex-disjoint (as opposed
to internally vertex-disjoint) paths. These are even easier to define:

\begin{definition}
Two paths $\mathbf{p}$ and $\mathbf{q}$ in a graph or digraph are said to be
\textbf{vertex-disjoint} if they have no common vertices.
\end{definition}

\begin{theorem}
[vertex-Menger theorem for directed graphs, multi-terminal version
2]\label{thm.menger.dvm-vd}Let $D=\left(  V,A,\psi\right)  $ be a
multidigraph, and let $X$ and $Y$ be two subsets of $V$.

Then, the maximum number of pairwise vertex-disjoint paths from $X$ to $Y$
equals the minimum size of an $X$-$Y$-vertex-separator.
\end{theorem}

\begin{example}
\label{exa.menger.dvm-vd.1}Let $D$ be the following multidigraph:%
\[%
%
\ \ .
\]
Then, the maximum number of pairwise vertex-disjoint paths from $X$ to $Y$ is
$2$. Here are two such paths (drawn in red and blue):%
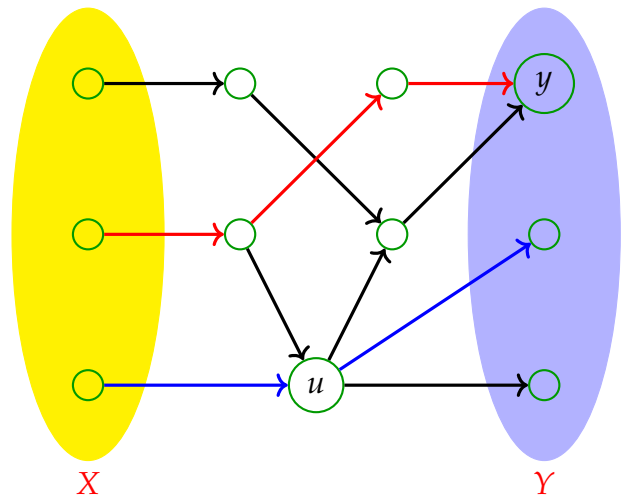
\[%
%
\ \ .
\]
If we were only looking for \textbf{internally} vertex-disjoint paths, then we
could add a third path to these two (namely, the path that starts at the
topmost vertex of $X$ and ends at the topmost vertex of $Y$). However, this
path and our red paths are only \textbf{internally} vertex-disjoint, not
vertex-disjoint. A little bit of thought shows that $D$ has no more than $2$
vertex-disjoint paths from $X$ to $Y$.

The minimum size of an $X$-$Y$-vertex-separator is $2$ as well; indeed,
$\left\{  u,y\right\}  $ is such an $X$-$Y$-vertex-separator. This number
equals the maximum number of pairwise vertex-disjoint paths from $X$ to $Y$,
just as Theorem \ref{thm.menger.dvm-vd} predicts.
\end{example}

\begin{proof}
[Proof of Theorem \ref{thm.menger.dvm-vd}.]We will reduce this to Theorem
\ref{thm.menger.dv2}, again by tweaking our digraph appropriately. This time,
the tweak is pretty simple: We add two new vertices $s$ and $t$ to $D$, and we
furthermore add an arc from $s$ to each $x\in X$ and an arc from each $y\in Y$
to $t$ (thus, we add a total of $\left\vert X\right\vert +\left\vert
Y\right\vert $ new arcs). We denote the resulting digraph by $D^{\prime}$. In
more detail, the definition of $D^{\prime}$ is as follows:

\begin{itemize}
\item We introduce two new vertices $s$ and $t$, and we set $V^{\prime}%
:=V\cup\left\{  s,t\right\}  $. This set $V^{\prime}$ will be the vertex set
of $D^{\prime}$.

\item For each $x\in X$, we introduce a new arc $a_{x}$, which shall have
source $s$ and target $x$.

\item For each $y\in Y$, we introduce a new arc $b_{y}$, which shall have
source $y$ and target $t$.

\item We let $A^{\prime}:=A\cup\left\{  a_{x}\ \mid\ x\in X\right\}
\cup\left\{  b_{y}\ \mid\ y\in Y\right\}  $. This set $A^{\prime}$ will be the
arc set of $D^{\prime}$.

\item We extend our map $\psi:A\rightarrow V\times V$ to a map $\psi^{\prime
}:A^{\prime}\rightarrow V^{\prime}\times V^{\prime}$ by setting%
\[
\psi^{\prime}\left(  a_{x}\right)  =\left(  s,x\right)
\ \ \ \ \ \ \ \ \ \ \text{for each }x\in X
\]
and%
\[
\psi^{\prime}\left(  b_{y}\right)  =\left(  y,t\right)
\ \ \ \ \ \ \ \ \ \ \text{for each }y\in Y
\]
(and, of course, $\psi^{\prime}\left(  c\right)  =\psi\left(  c\right)  $ for
each $c\in A$).

\item We define $D^{\prime}$ to be the multidigraph $\left(  V^{\prime
},A^{\prime},\psi^{\prime}\right)  $.
\end{itemize}

For instance, if $D$ is the multidigraph from Example
\ref{exa.menger.dvm-vd.1}, then $D^{\prime}$ looks as follows:%
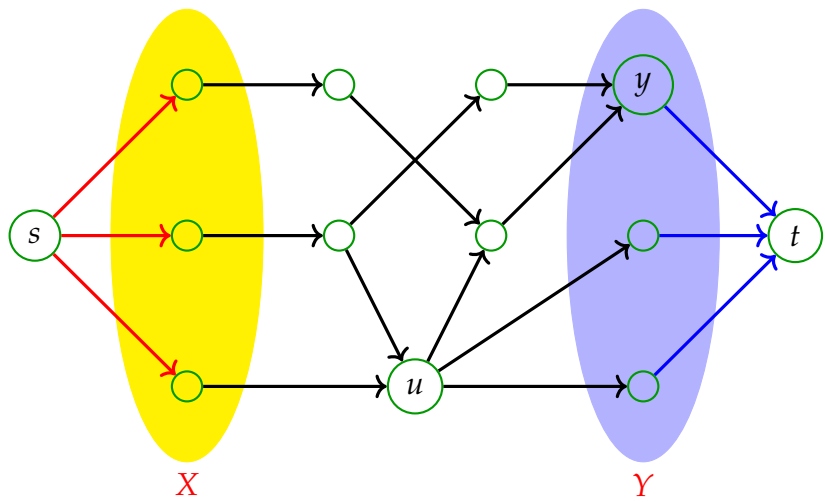
\[%
%
\ \ .
\]
(The arcs $a_{x}$ are drawn in red; the arcs $b_{y}$ are drawn in blue.)

By its construction, the digraph $D^{\prime}$ has no arc with source $s$ and
target $t$. Hence, Theorem \ref{thm.menger.dv2} (applied to $D^{\prime
}=\left(  V^{\prime},A^{\prime},\psi^{\prime}\right)  $ instead of $D=\left(
V,A,\psi\right)  $) yields that the maximum number of pairwise internally
vertex-disjoint paths from $s$ to $t$ equals the minimum size of an internal
$s$-$t$-vertex-separator. However, it is easy to see the following two claims:

\begin{statement}
\textit{Claim 1:} The maximum number of pairwise internally vertex-disjoint
paths from $s$ to $t$ equals the maximum number of pairwise vertex-disjoint
paths from $X$ to $Y$ (in $D$).
\end{statement}

\begin{statement}
\textit{Claim 2:} The minimum size of an internal $s$-$t$-vertex-separator
equals the minimum size of an $X$-$Y$-vertex-separator (in $D$).
\end{statement}

[\textit{Proof of Claim 1 (sketched):} Given any path $\mathbf{p}$ from $s$ to
$t$, we can remove the starting point and the ending point of this path; the
result will always be a path from $X$ to $Y$ (in $D$). Let us denote the
latter path by $\overline{\mathbf{p}}$. Thus, we obtain a map%
\begin{align*}
\left\{  \text{paths from }s\text{ to }t\right\}   &  \rightarrow\left\{
\text{paths from }X\text{ to }Y\text{ (in }D\text{)}\right\}  ,\\
\mathbf{p}  &  \mapsto\overline{\mathbf{p}}.
\end{align*}
This map is easily seen to be a bijection (indeed, if $\mathbf{q}$ is a path
from $X$ to $Y$ (in $D$), then we can easily extend it to a path $\mathbf{p}$
from $s$ to $t$ by inserting the appropriate arc $a_{x}$ at its beginning and
the appropriate arc $b_{y}$ at its end; this latter path $\mathbf{p}$ will
then satisfy $\overline{\mathbf{p}}=\mathbf{q}$). Moreover, two paths
$\mathbf{p}$ and $\mathbf{q}$ from $s$ to $t$ are internally vertex-disjoint
if and only if the corresponding paths $\overline{\mathbf{p}}$ and
$\overline{\mathbf{q}}$ are vertex-disjoint (because the intermediate vertices
of $\mathbf{p}$ are the vertices of $\overline{\mathbf{p}}$, whereas the
intermediate vertices of $\mathbf{q}$ are the vertices of $\overline
{\mathbf{q}}$). This proves Claim 1.] \medskip

[\textit{Proof of Claim 2 (sketched):} It is easy to see that the internal
$s$-$t$-vertex-separators are precisely the $X$-$Y$-vertex-separators (in
$D$). (To show this, compare the definitions of these two objects using the
bijection from the proof of Claim 1, and observe that $V^{\prime}%
\setminus\left\{  s,t\right\}  =V$.) From this, Claim 2 follows.] \medskip

Recall that the maximum number of pairwise internally vertex-disjoint paths
from $s$ to $t$ equals the minimum size of an internal $s$-$t$%
-vertex-separator. In view of Claim 1 and Claim 2, we can rewrite this as
follows: The maximum number of pairwise vertex-disjoint paths from $X$ to $Y$
equals the minimum size of an $X$-$Y$-vertex-separator. Thus, Theorem
\ref{thm.menger.dvm-vd} is proved.
\end{proof}

We note that Hall's Marriage Theorem (Theorem \ref{thm.match.HMT}) can be
easily derived from any of the directed Menger theorems (exercise!). I have
heard that this can also be done in reverse. This places the Menger theorems
in the cluster of theorems equivalent to Hall's Marriage Theorem (such as
K\"{o}nig's theorem).

\subsubsection{The vertex-Menger theorem for undirected graphs}

Vertex-Menger theorems also exist for undirected graphs. Here are the
undirected analogues of Theorem \ref{thm.menger.dv2}, Theorem
\ref{thm.menger.dvm} and Theorem \ref{thm.menger.dvm-vd}, along with the
definitions they rely on:

\begin{definition}
Let $G=\left(  V,E,\varphi\right)  $ be a multigraph, and let $s$ and $t$ be
two vertices of $G$. A subset $W$ of $V\setminus\left\{  s,t\right\}  $ is
said to be an \textbf{internal }$s$\textbf{-}$t$\textbf{-vertex-separator} if
each path from $s$ to $t$ contains at least one vertex from $W$. Equivalently,
a subset $W$ of $V\setminus\left\{  s,t\right\}  $ is said to be an
\textbf{internal }$s$\textbf{-}$t$\textbf{-vertex-separator} if the induced
subgraph of $G$ on the set $V\setminus W$ has no path from $s$ to $t$ (in
other words, removing from $G$ all vertices contained in $W$ destroys all
paths from $s$ to $t$).
\end{definition}

\begin{theorem}
[vertex-Menger theorem for undirected graphs]\label{thm.menger.uv2}Let
$G=\left(  V,E,\varphi\right)  $ be a multigraph, and let $s$ and $t$ be two
distinct vertices of $G$. Assume that $G$ has no edge with endpoints $s$ and
$t$. Then, the maximum number of pairwise internally vertex-disjoint paths
from $s$ to $t$ equals the minimum size of an internal $s$-$t$-vertex-separator.
\end{theorem}

\begin{definition}
Let $G=\left(  V,E,\varphi\right)  $ be a multigraph, and let $X$ and $Y$ be
two subsets of $V$.

\begin{enumerate}
\item[\textbf{(a)}] A \textbf{path from }$X$\textbf{ to }$Y$ shall mean a path
whose starting point belongs to $X$ and whose ending point belongs to $Y$.

\item[\textbf{(b)}] A subset $W$ of $V$ is said to be an $X$\textbf{-}%
$Y$\textbf{-vertex-separator} if each path from $X$ to $Y$ contains at least
one vertex from $W$. Equivalently, a subset $W$ of $V$ is said to be an
$X$\textbf{-}$Y$\textbf{-vertex-separator} if the induced subgraph of $G$ on
the set $V\setminus W$ has no path from $X$ to $Y$ (in other words, removing
from $G$ all vertices contained in $W$ destroys all paths from $X$ to $Y$).

\item[\textbf{(c)}] An $X$-$Y$-vertex-separator $W$ is said to be
\textbf{internal} if it is a subset of $V\setminus\left(  X\cup Y\right)  $
(that is, if it is disjoint from $X$ and from $Y$).
\end{enumerate}
\end{definition}

\begin{theorem}
[vertex-Menger theorem for undirected graphs, multi-terminal version
1]\label{thm.menger.uvm}Let $G=\left(  V,E,\varphi\right)  $ be a multigraph,
and let $X$ and $Y$ be two disjoint subsets of $V$. Assume that $G$ has no
edge with one endpoint in $X$ and the other endpoint in $Y$.

Then, the maximum number of pairwise internally vertex-disjoint paths from $X$
to $Y$ equals the minimum size of an internal $X$-$Y$-vertex-separator.
\end{theorem}

\begin{theorem}
[vertex-Menger theorem for undirected graphs, multi-terminal version
2]\label{thm.menger.uvm-vd}Let $G=\left(  V,E,\varphi\right)  $ be a
multigraph, and let $X$ and $Y$ be two subsets of $V$.

Then, the maximum number of pairwise vertex-disjoint paths from $X$ to $Y$
equals the minimum size of an $X$-$Y$-vertex-separator.
\end{theorem}

Theorem \ref{thm.menger.uv2}, Theorem \ref{thm.menger.uvm} and Theorem
\ref{thm.menger.uvm-vd} follow immediately by applying the analogous theorems
for directed graphs (i.e., Theorem \ref{thm.menger.dv2}, Theorem
\ref{thm.menger.dvm} and Theorem \ref{thm.menger.dvm-vd}) to the digraph
$G^{\operatorname*{bidir}}$ instead of $D$ (since the paths of $G$ are in
bijection with the paths of $G^{\operatorname*{bidir}}$).

\subsection{\label{sec.paths.gm}The Gallai--Milgram theorem}

Next, we proceed to some more obscure properties of paths in digraphs and graphs.

\subsubsection{Definitions}

In order to state the first of these properties, we need the following three
definitions (the first two of which were already made in Section \ref{sec.ker}):

\begin{definition}
Two vertices $u$ and $v$ of a multidigraph $D$ are said to be
\textbf{adjacent} if they are adjacent in the undirected graph
$D^{\operatorname*{und}}$. (In other words, they are adjacent if and only if
$D$ has an arc with source $u$ and target $v$ or an arc with source $v$ and
target $u$.)
\end{definition}

\begin{definition}
An \textbf{independent set} of a multidigraph $D$ means a subset $S$ of
$\operatorname*{V}\left(  D\right)  $ such that no two elements of $S$ are
adjacent. In other words, it means an independent set of the undirected graph
$D^{\operatorname*{und}}$.
\end{definition}

\begin{definition}
A \textbf{path cover} of a multidigraph $D$ means a set of paths of $D$ such
that each vertex of $D$ is contained in exactly one of these paths.
\end{definition}

\begin{example}
\label{exa.gm.pc}Let $D$ be the following digraph:%
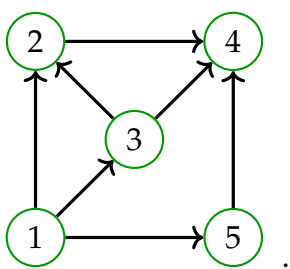
\[%
\begin{tikzpicture}[scale=1.3]
\begin{scope}[every node/.style={circle,thick,draw=green!60!black}]
\node(1) at (-1, -1) {$1$};
\node(2) at (-1, 1) {$2$};
\node(3) at (0, 0) {$3$};
\node(4) at (1, 1) {$4$};
\node(5) at (1, -1) {$5$};
\end{scope}
\begin{scope}[every edge/.style={draw=black,very thick}, every loop/.style={}]
\path[->] (1) edge (2) edge (3) edge (5);
\path[->] (2) edge (4);
\path[->] (3) edge (2) edge (4);
\path[->] (5) edge (4);
\end{scope}
\end{tikzpicture}%
\ \ .
\]
Then, $\left\{  \left(  1,\ast,5,\ast,4\right)  ,\ \left(  3,\ast,2\right)
\right\}  $ is a path cover of $D$ (we are again writing asterisks for the
arcs, since the arcs of $D$ are uniquely determined by their sources and their
targets). Another path cover of $D$ is $\left\{  \left(  1,\ast,3,\ast
,4\right)  ,\ \left(  2\right)  ,\ \left(  5\right)  \right\}  $. Yet another
path cover of $D$ is $\left\{  \left(  1\right)  ,\ \left(  2\right)
,\ \left(  3\right)  ,\ \left(  4\right)  ,\ \left(  5\right)  \right\}  $.
There are many more.

Note that the set $\left\{  \left(  1,\ast,5,\ast,4\right)  ,\ \left(
3,\ast,2,\ast,4\right)  \right\}  $ is not a path cover of $D$, since the
vertex $4$ is contained in two (not one) of its paths.

Let us draw the three path covers we have mentioned (by simply drawing the
arcs of the paths they contain, while omitting all other arcs of $D$):%
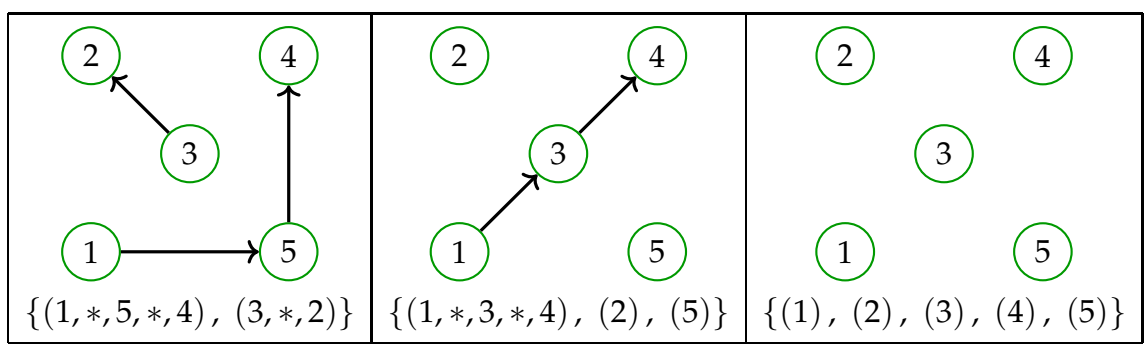
\[%

\]

\end{example}

(Note that we have already seen path covers of a \textquotedblleft
complete\textquotedblright\ simple digraph $\left(  V,\ V\times V\right)  $ in
the proof of Theorem \ref{thm.hamp.Dbar}; we called them \textquotedblleft
path covers of $V$\textquotedblright.)

\begin{remark}
Let $D$ be a digraph. A path cover of $D$ consisting of only $1$ path is the
same as a Hamiltonian path of $D$. (More precisely: A single path $\mathbf{p}$
forms a path cover $\left\{  \mathbf{p}\right\}  $ of $D$ if and only if
$\mathbf{p}$ is a Hamiltonian path.)
\end{remark}

\subsubsection{The Gallai--Milgram theorem}

Now, the \textbf{Gallai--Milgram theorem} states the following:

\begin{theorem}
[Gallai--Milgram theorem]\label{thm.gm.gm}Let $D$ be a loopless digraph. Then,
there exist a path cover $\mathcal{P}$ of $D$ and an independent set $S$ of
$D$ such that $S$ has exactly one vertex from each path in $\mathcal{P}$ (in
other words, for each path $\mathbf{p}\in\mathcal{P}$, exactly one vertex of
$\mathbf{p}$ belongs to $S$).
\end{theorem}

\begin{example}
Let $D$ be the digraph from Example \ref{exa.gm.pc}. Then, Theorem
\ref{thm.gm.gm} tells us that there exist a path cover $\mathcal{P}$ of $D$
and an independent set $S$ of $D$ such that $S$ has exactly one vertex from
each path in $\mathcal{P}$. For example, we can take $\mathcal{P}=\left\{
\left(  1,\ast,5,\ast,4\right)  ,\ \left(  3,\ast,2\right)  \right\}  $ and
$S=\left\{  5,3\right\}  $.
\end{example}

We will now prove Theorem \ref{thm.gm.gm}, following Diestel's book
\cite[Theorem 2.5.1]{Dieste17}:

\begin{noncompile}
See also https://math.stackexchange.com/questions/897278/ for details missing
from Diestel.
\end{noncompile}

\begin{proof}
[Proof of Theorem \ref{thm.gm.gm}.]Write the multidigraph $D$ as $D=\left(
V,A,\psi\right)  $. We introduce a notation:

\begin{itemize}
\item If $\mathcal{P}$ is a path cover of $D$, then a \textbf{cross-cut} of
$\mathcal{P}$ means a subset $S$ of $V$ that contains exactly one vertex from
each path in $\mathcal{P}$.
\end{itemize}

Thus, the claim we must prove is saying that there exist a path cover
$\mathcal{P}$ of $D$ and an independent cross-cut of $\mathcal{P}$.

We will show something stronger:

\begin{statement}
\textit{Claim 1:} Any minimum-size path cover of $D$ has an independent cross-cut.
\end{statement}

Note that the size of a path cover means the number of paths in it. Thus, a
minimum-size path cover means a path cover with the smallest possible number
of paths.

We will show something even stronger than Claim 1. To state this stronger
claim, we need more notations:

\begin{itemize}
\item If $\mathcal{P}$ is a path cover, then $\operatorname*{Ends}\mathcal{P}$
means the set of the ending points of all paths in $\mathcal{P}$. Note that
$\left\vert \operatorname*{Ends}\mathcal{P}\right\vert =\left\vert
\mathcal{P}\right\vert $.

\item A path cover $\mathcal{P}$ is said to be \textbf{end-minimal} if no
proper subset of $\operatorname*{Ends}\mathcal{P}$ can be written as
$\operatorname*{Ends}\mathcal{Q}$ for a path cover $\mathcal{Q}$.
\end{itemize}

\begin{example}
For instance, if $D$ is as in Example \ref{exa.gm.pc}, and if
\begin{align*}
\mathcal{P}  &  =\left\{  \left(  1,\ast,5,\ast,4\right)  ,\ \left(
3,\ast,2\right)  \right\}  ,\\
\mathcal{Q}  &  =\left\{  \left(  1,\ast,3,\ast,4\right)  ,\ \left(  2\right)
,\ \left(  5\right)  \right\}  ,\\
\mathcal{R}  &  =\left\{  \left(  1\right)  ,\ \left(  2\right)  ,\ \left(
3\right)  ,\ \left(  4\right)  ,\ \left(  5\right)  \right\}
\end{align*}
are the three path covers from Example \ref{exa.gm.pc}, then
\[
\operatorname*{Ends}\mathcal{P}=\left\{  4,2\right\}
,\ \ \ \ \ \ \ \ \ \ \operatorname*{Ends}\mathcal{Q}=\left\{  4,2,5\right\}
,\ \ \ \ \ \ \ \ \ \ \operatorname*{Ends}\mathcal{R}=\left\{
1,2,3,4,5\right\}  ,
\]
which shows immediately that neither $\mathcal{Q}$ nor $\mathcal{R}$ is
end-minimal (since $\operatorname*{Ends}\mathcal{P}$ is a proper subset of
each of $\operatorname*{Ends}\mathcal{Q}$ and $\operatorname*{Ends}%
\mathcal{R}$). It is easy to see that $\mathcal{P}$ is end-minimal (and also minimum-size).
\end{example}

Back to the general case. Clearly, any minimum-size path cover of $D$ is also
end-minimal\footnote{\textit{Proof.} Let $\mathcal{P}$ be a minimum-size path
cover of $D$. If $\mathcal{P}$ was not end-minimal, then there would be a path
cover $\mathcal{Q}$ with $\left\vert \operatorname*{Ends}\mathcal{Q}%
\right\vert <\left\vert \operatorname*{Ends}\mathcal{P}\right\vert $ and
therefore $\left\vert \mathcal{Q}\right\vert =\left\vert \operatorname*{Ends}%
\mathcal{Q}\right\vert <\left\vert \operatorname*{Ends}\mathcal{P}\right\vert
=\left\vert \mathcal{P}\right\vert $; but this would contradict the fact that
$\mathcal{P}$ is minimum-size. Hence, $\mathcal{P}$ is end-minimal.}. Thus,
the following claim is stronger than Claim 1:

\begin{statement}
\textit{Claim 2:} Any end-minimal path cover of $D$ has an independent cross-cut.
\end{statement}

It is Claim 2 that we will be proving.\footnote{On a sidenote: Is Claim 2
really stronger than Claim 1? Yes, because it can happen that some end-minimal
path cover fails to be minimum-size. For example, the path cover $\left\{
\left(  1,\ast,2,\ast,3\right)  ,\ \left(  4\right)  \right\}  $ in the
digraph%
\[%
\begin{tikzpicture}[scale=2]
\begin{scope}[every node/.style={circle,thick,draw=green!60!black}]
\node(1) at (0:1) {$1$};
\node(2) at (90:1) {$2$};
\node(3) at (180:1) {$3$};
\node(4) at (270:1) {$4$};
\end{scope}
\begin{scope}[every edge/.style={draw=black,very thick}, every loop/.style={}]
\path[->] (1) edge (2) (2) edge[bend right=20] (3);
\path[->] (3) edge[bend right=20] (2);
\path[->] (1) edge (4) (4) edge (3);
\end{scope}
\end{tikzpicture}%
\]
has this property.}

[\textit{Proof of Claim 2:} We proceed by induction on $\left\vert
V\right\vert $.

\textit{Base case:} Claim 2 is obvious when $\left\vert V\right\vert =0$
(since $\varnothing$ is an independent cross-cut in this case).

\textit{Induction step:} Consider a multidigraph $D=\left(  V,A,\psi\right)  $
with $\left\vert V\right\vert =N$. Assume (as the induction hypothesis) that
Claim 2 is already proved for any multidigraph with $N-1$ vertices.

Let $\mathcal{P}$ be an end-minimal path cover of $D$. We must show that
$\mathcal{P}$ has an independent cross-cut.

Let $\mathbf{p}_{1},\mathbf{p}_{2},\ldots,\mathbf{p}_{k}$ be the paths in
$\mathcal{P}$ (listed without repetitions), and let $v_{1},v_{2},\ldots,v_{k}$
be their respective ending points. Thus, $\left\{  v_{1},v_{2},\ldots
,v_{k}\right\}  =\operatorname*{Ends}\mathcal{P}$ and $k=\left\vert
\operatorname*{Ends}\mathcal{P}\right\vert $.

Recall that we must find an independent cross-cut of $\mathcal{P}$. If the set
$\left\{  v_{1},v_{2},\ldots,v_{k}\right\}  $ is independent, then we are done
(since this set $\left\{  v_{1},v_{2},\ldots,v_{k}\right\}  $ is clearly a
cross-cut of $\mathcal{P}$). Thus, we WLOG assume that this is not the case.
Hence, there is an arc from some vertex $v_{i}$ to some vertex $v_{j}$. These
two vertices $v_{i}$ and $v_{j}$ are distinct (because $D$ is loopless). Since
we can swap our paths $\mathbf{p}_{1},\mathbf{p}_{2},\ldots,\mathbf{p}_{k}$
(and thus their ending points $v_{1},v_{2},\ldots,v_{k}$) at will, we can thus
WLOG assume that $i=2$ and $j=1$. Assume this. Thus, there is an arc from
$v_{2}$ to $v_{1}$. We shall refer to this arc as the \emph{blue arc}, and we
will draw it accordingly:\footnote{This picture illustrates just one
representative case, with $k=4$. The four columns (from left to right) are the
four paths $\mathbf{p}_{1},\mathbf{p}_{2},\mathbf{p}_{3},\mathbf{p}_{4}$. Of
course, the digraph $D$ can have many more arcs than we have drawn on this
picture, but we are not interested in them right now.}%
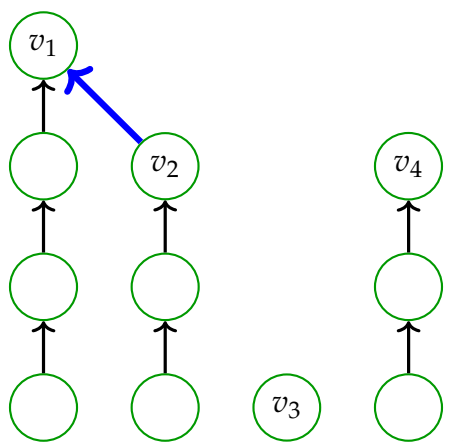
\[%
\begin{tikzpicture}[scale=1.6]
\begin{scope}[every node/.style={circle,thick,draw=green!60!black}]
\node(11) at (0, 0) {$\phantom{v_1}$};
\node(12) at (0, 1) {$\phantom{v_1}$};
\node(13) at (0, 2) {$\phantom{v_1}$};
\node(14) at (0, 3) {$v_1$};
\node(21) at (1, 0) {$\phantom{v_1}$};
\node(22) at (1, 1) {$\phantom{v_1}$};
\node(23) at (1, 2) {$v_2$};
\node(31) at (2, 0) {$v_3$};
\node(41) at (3, 0) {$\phantom{v_4}$};
\node(42) at (3, 1) {$\phantom{v_4}$};
\node(43) at (3, 2) {$v_4$};
\end{scope}
\begin{scope}[every edge/.style={draw=black,very thick}, every loop/.style={}]
\path[->] (11) edge (12) (12) edge (13) (13) edge (14);
\path[->] (21) edge (22) (22) edge (23);
\path[->] (41) edge (42) (42) edge (43);
\end{scope}
\begin{scope}[every edge/.style={draw=blue,line width=0.2pc}%
, every loop/.style={}]
\path[->] (23) edge (14);
\end{scope}
\end{tikzpicture}%
\ \ .
\]

We can extend the path $\mathbf{p}_{2}$ beyond its ending point $v_{2}$ by
inserting the blue arc and the vertex $v_{1}$ at its end. This results in a
new path, which we denote by $\mathbf{p}_{2}+v_{1}$; this path has ending
point $v_{1}$.

If $v_{1}$ is the only vertex on the path $\mathbf{p}_{1}$ (that is, if the
path $\mathbf{p}_{1}$ has length $0$), then we can therefore replace the path
$\mathbf{p}_{2}$ by $\mathbf{p}_{2}+v_{1}$ and remove the length-$0$ path
$\mathbf{p}_{1}$ from our path cover $\mathcal{P}$, and we thus obtain a new
path cover $\mathcal{Q}$ such that $\operatorname*{Ends}\mathcal{Q}$ is a
proper subset of $\operatorname*{Ends}\mathcal{P}$. But this is impossible,
since we assumed that $\mathcal{P}$ is end-minimal. Therefore, $v_{1}$ is not
the only vertex on $\mathbf{p}_{1}$.

Thus, let $v$ be the second-to-last vertex on $\mathbf{p}_{1}$ (that is, the
vertex that is immediately followed by $v_{1}$). Then, the path $\mathbf{p}%
_{1}$ contains an arc from $v$ to $v_{1}$. We shall refer to this arc as the
\emph{red arc}, and we will draw it accordingly:%
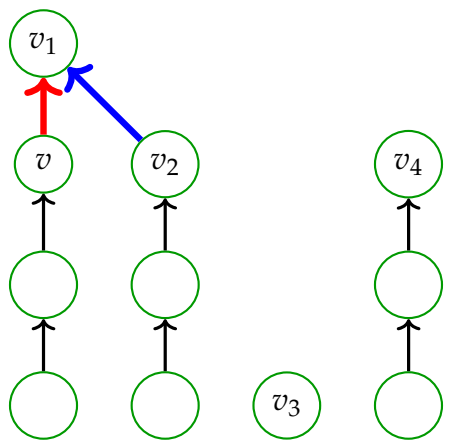
\[%
%
\ \ .
\]

Let $D^{\prime}$ be the digraph $D\setminus v_{1}$ (that is, the digraph
obtained from $D$ by removing the vertex $v_{1}$ and all arcs that have
$v_{1}$ as source or target). Let $\mathbf{p}_{1}^{\prime}$ be the result of
removing the vertex $v_{1}$ and the red arc from the path $\mathbf{p}_{1}$.
Then, $\mathcal{P}^{\prime}:=\left\{  \mathbf{p}_{1}^{\prime},\mathbf{p}%
_{2},\mathbf{p}_{3},\ldots,\mathbf{p}_{k}\right\}  $ is a path cover of
$D^{\prime}$. Note that the path $\mathbf{p}_{1}^{\prime}$ has ending point
$v$ (since it is obtained from $\mathbf{p}_{1}$ by removing the last vertex
and the last arc, but we know that the second-to-last vertex on $\mathbf{p}%
_{1}$ is $v$), whereas the paths $\mathbf{p}_{2},\mathbf{p}_{3},\ldots
,\mathbf{p}_{k}$ have ending points $v_{2},v_{3},\ldots,v_{k}$. Thus,
$\operatorname*{Ends}\left(  \mathcal{P}^{\prime}\right)  =\left\{
v,v_{2},v_{3},\ldots,v_{k}\right\}  $. Here is an illustration of the digraph
$D^{\prime}=D\setminus v_{1}$ and its path cover $\mathcal{P}^{\prime}$:%
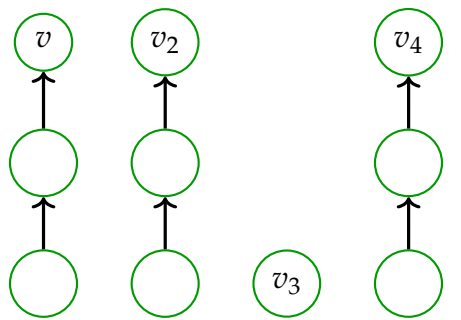
\[%
\begin{tikzpicture}[scale=1.6]
\begin{scope}[every node/.style={circle,thick,draw=green!60!black}]
\node(11) at (0, 0) {$\phantom{v_1}$};
\node(12) at (0, 1) {$\phantom{v_1}$};
\node(13) at (0, 2) {$v\vphantom{v_1}$};
\node(21) at (1, 0) {$\phantom{v_1}$};
\node(22) at (1, 1) {$\phantom{v_1}$};
\node(23) at (1, 2) {$v_2$};
\node(31) at (2, 0) {$v_3$};
\node(41) at (3, 0) {$\phantom{v_4}$};
\node(42) at (3, 1) {$\phantom{v_4}$};
\node(43) at (3, 2) {$v_4$};
\end{scope}
\begin{scope}[every edge/.style={draw=black,very thick}, every loop/.style={}]
\path[->] (11) edge (12) (12) edge (13);
\path[->] (21) edge (22) (22) edge (23);
\path[->] (41) edge (42) (42) edge (43);
\end{scope}
\end{tikzpicture}%
\ \ .
\]

Consider the path cover $\mathcal{P}^{\prime}$ of $D^{\prime}$. If we can find
an independent cross-cut of $\mathcal{P}^{\prime}$, then we will be done,
because any such cross-cut will also be an independent cross-cut of our
original path cover $\left\{  \mathbf{p}_{1},\mathbf{p}_{2},\ldots
,\mathbf{p}_{k}\right\}  =\mathcal{P}$. Since the digraph $D\setminus v_{1}$
has $N-1$ vertices\footnote{because the digraph $D$ has $\left\vert
V\right\vert =N$ vertices}, we can find such an independent cross-cut by our
induction hypothesis if we can prove that the path cover $\mathcal{P}^{\prime
}$ is end-minimal (as a path cover of $D^{\prime}$).

So let us prove this now. Indeed, assume the contrary. Thus, $D^{\prime}$ has
a path cover $\mathcal{Q}^{\prime}$ such that $\operatorname*{Ends}\left(
\mathcal{Q}^{\prime}\right)  $ is a proper subset of $\operatorname*{Ends}%
\left(  \mathcal{P}^{\prime}\right)  $. Consider this $\mathcal{Q}^{\prime}$.
Note that\footnote{The symbol \textquotedblleft$\subsetneq$\textquotedblright%
\ (note that the stroke only crosses the straight line, not the curved one)
means \textquotedblleft proper subset of\textquotedblright.}
\[
\operatorname*{Ends}\left(  \mathcal{Q}^{\prime}\right)  \subsetneq
\operatorname*{Ends}\left(  \mathcal{P}^{\prime}\right)  =\left\{
v,v_{2},v_{3},\ldots,v_{k}\right\}  .
\]
As a consequence, $\left\vert \operatorname*{Ends}\left(  \mathcal{Q}^{\prime
}\right)  \right\vert <\left\vert \left\{  v,v_{2},v_{3},\ldots,v_{k}\right\}
\right\vert =k$.

Now, we are in one of the following three cases:

\textit{Case 1:} We have $v\in\operatorname*{Ends}\left(  \mathcal{Q}^{\prime
}\right)  $.

\textit{Case 2:} We have $v\notin\operatorname*{Ends}\left(  \mathcal{Q}%
^{\prime}\right)  $ but $v_{2}\in\operatorname*{Ends}\left(  \mathcal{Q}%
^{\prime}\right)  $.

\textit{Case 3:} We have $v\notin\operatorname*{Ends}\left(  \mathcal{Q}%
^{\prime}\right)  $ and $v_{2}\notin\operatorname*{Ends}\left(  \mathcal{Q}%
^{\prime}\right)  $.

Let us consider these cases one by one:

\begin{itemize}
\item We first consider Case 1. In this case, we have $v\in
\operatorname*{Ends}\left(  \mathcal{Q}^{\prime}\right)  $. In other words,
some path $\mathbf{p}\in\mathcal{Q}^{\prime}$ ends at $v$. Let us extend this
path $\mathbf{p}$ beyond $v$ by inserting the red arc and the vertex $v_{1}$
at its end. Thus, we obtain a path of $D$, which we call $\mathbf{p}+v_{1}$.
Replacing $\mathbf{p}$ by $\mathbf{p}+v_{1}$ in $\mathcal{Q}^{\prime}$, we
obtain a path cover $\mathcal{Q}$ of $D$ such that $\operatorname*{Ends}%
\mathcal{Q}$ is a proper subset of $\operatorname*{Ends}\mathcal{P}%
$\ \ \ \ \footnote{\textit{Proof.} We obtained $\mathcal{Q}$ from
$\mathcal{Q}^{\prime}$ by replacing $\mathbf{p}$ by $\mathbf{p}+v_{1}$. As a
consequence of this replacement, the ending point $v$ of $\mathbf{p}$ has been
replaced by the ending point $v_{1}$ of $\mathbf{p}+v_{1}$. Thus,
\begin{align*}
\operatorname*{Ends}\mathcal{Q}  &  =\underbrace{\left(  \operatorname*{Ends}%
\left(  \mathcal{Q}^{\prime}\right)  \setminus\left\{  v\right\}  \right)
}_{\substack{\subseteq\left\{  v_{2},v_{3},\ldots,v_{k}\right\}
\\\text{(since }\operatorname*{Ends}\left(  \mathcal{Q}^{\prime}\right)
\subseteq\left\{  v,v_{2},v_{3},\ldots,v_{k}\right\}  \text{)}}}\cup\left\{
v_{1}\right\} \\
&  \subseteq\left\{  v_{2},v_{3},\ldots,v_{k}\right\}  \cup\left\{
v_{1}\right\}  =\left\{  v_{1},v_{2},\ldots,v_{k}\right\}
=\operatorname*{Ends}\mathcal{P}.
\end{align*}
For the same reason, we have $\left\vert \operatorname*{Ends}\mathcal{Q}%
\right\vert =\left\vert \operatorname*{Ends}\left(  \mathcal{Q}^{\prime
}\right)  \right\vert <k=\left\vert \operatorname*{Ends}\mathcal{P}\right\vert
$, so that $\operatorname*{Ends}\mathcal{Q}\neq\operatorname*{Ends}%
\mathcal{P}$. Combining this with $\operatorname*{Ends}\mathcal{Q}%
\subseteq\operatorname*{Ends}\mathcal{P}$, we conclude that
$\operatorname*{Ends}\mathcal{Q}$ is a proper subset of $\operatorname*{Ends}%
\mathcal{P}$.}. But this contradicts the fact that $\mathcal{P}$ is
end-minimal. Thus, we have obtained a contradiction in Case 1.

\item Next, we consider Case 2. In this case, we have $v\notin%
\operatorname*{Ends}\left(  \mathcal{Q}^{\prime}\right)  $ but $v_{2}%
\in\operatorname*{Ends}\left(  \mathcal{Q}^{\prime}\right)  $. Combining
$\operatorname*{Ends}\left(  \mathcal{Q}^{\prime}\right)  \subseteq\left\{
v,v_{2},v_{3},\ldots,v_{k}\right\}  $ with $v\notin\operatorname*{Ends}\left(
\mathcal{Q}^{\prime}\right)  $, we obtain
\[
\operatorname*{Ends}\left(  \mathcal{Q}^{\prime}\right)  \subseteq\left\{
v,v_{2},v_{3},\ldots,v_{k}\right\}  \setminus\left\{  v\right\}  =\left\{
v_{2},v_{3},\ldots,v_{k}\right\}  .
\]
From $v_{2}\in\operatorname*{Ends}\left(  \mathcal{Q}^{\prime}\right)  $, we
see that some path $\mathbf{p}\in\mathcal{Q}^{\prime}$ ends at $v_{2}$. Let us
extend this path $\mathbf{p}$ beyond $v_{2}$ by inserting the blue arc and the
vertex $v_{1}$ at its end. Thus, we obtain a path of $D$, which we call
$\mathbf{p}+v_{1}$. Replacing $\mathbf{p}$ by $\mathbf{p}+v_{1}$ in
$\mathcal{Q}^{\prime}$, we obtain a path cover $\mathcal{Q}$ of $D$ such that
$\operatorname*{Ends}\mathcal{Q}$ is a proper subset of $\operatorname*{Ends}%
\mathcal{P}$\ \ \ \ \footnote{\textit{Proof.} We obtained $\mathcal{Q}$ from
$\mathcal{Q}^{\prime}$ by replacing $\mathbf{p}$ by $\mathbf{p}+v_{1}$. As a
consequence of this replacement, the ending point $v_{2}$ of $\mathbf{p}$ has
been replaced by the ending point $v_{1}$ of $\mathbf{p}+v_{1}$. Thus,
\begin{align*}
\operatorname*{Ends}\mathcal{Q}  &  =\underbrace{\left(  \operatorname*{Ends}%
\left(  \mathcal{Q}^{\prime}\right)  \setminus\left\{  v_{2}\right\}  \right)
}_{\substack{\subseteq\left\{  v_{3},v_{4},\ldots,v_{k}\right\}
\\\text{(since }\operatorname*{Ends}\left(  \mathcal{Q}^{\prime}\right)
\subseteq\left\{  v_{2},v_{3},\ldots,v_{k}\right\}  \text{)}}}\cup\left\{
v_{1}\right\} \\
&  \subseteq\left\{  v_{3},v_{4},\ldots,v_{k}\right\}  \cup\left\{
v_{1}\right\}  =\left\{  v_{1},v_{3},v_{4},\ldots,v_{k}\right\} \\
&  \subsetneq\left\{  v_{1},v_{2},\ldots,v_{k}\right\}  =\operatorname*{Ends}%
\mathcal{P}.
\end{align*}
In other words, $\operatorname*{Ends}\mathcal{Q}$ is a proper subset of
$\operatorname*{Ends}\mathcal{P}$.}. But this contradicts the fact that
$\mathcal{P}$ is end-minimal. Thus, we have obtained a contradiction in Case 2.

\item Finally, we consider Case 3. In this case, we have $v\notin%
\operatorname*{Ends}\left(  \mathcal{Q}^{\prime}\right)  $ and $v_{2}%
\notin\operatorname*{Ends}\left(  \mathcal{Q}^{\prime}\right)  $. Combining
this with $\operatorname*{Ends}\left(  \mathcal{Q}^{\prime}\right)
\subseteq\left\{  v,v_{2},v_{3},\ldots,v_{k}\right\}  $, we obtain%
\[
\operatorname*{Ends}\left(  \mathcal{Q}^{\prime}\right)  \subseteq\left\{
v,v_{2},v_{3},\ldots,v_{k}\right\}  \setminus\left\{  v,v_{2}\right\}
=\left\{  v_{3},v_{4},\ldots,v_{k}\right\}  ,
\]
so that $\left\vert \operatorname*{Ends}\left(  \mathcal{Q}^{\prime}\right)
\right\vert \leq\left\vert \left\{  v_{3},v_{4},\ldots,v_{k}\right\}
\right\vert =k-2$. Now, adding the trivial path $\left(  v_{1}\right)  $ to
$\mathcal{Q}^{\prime}$ yields a path cover $\mathcal{Q}$ of $D$ such that
$\operatorname*{Ends}\mathcal{Q}$ is a proper subset of $\operatorname*{Ends}%
\mathcal{P}$\ \ \ \ \footnote{\textit{Proof.} We obtained $\mathcal{Q}$ from
$\mathcal{Q}^{\prime}$ by adding the trivial path $\left(  v_{1}\right)  $,
whose ending point is $v_{1}$. Thus,
\begin{align*}
\operatorname*{Ends}\mathcal{Q}  &  =\underbrace{\operatorname*{Ends}\left(
\mathcal{Q}^{\prime}\right)  }_{\subseteq\left\{  v_{3},v_{4},\ldots
,v_{k}\right\}  }\cup\left\{  v_{1}\right\}  \subseteq\left\{  v_{3}%
,v_{4},\ldots,v_{k}\right\}  \cup\left\{  v_{1}\right\}  =\left\{  v_{1}%
,v_{3},v_{4},\ldots,v_{k}\right\} \\
&  \subsetneq\left\{  v_{1},v_{2},\ldots,v_{k}\right\}  =\operatorname*{Ends}%
\mathcal{P}.
\end{align*}
In other words, $\operatorname*{Ends}\mathcal{Q}$ is a proper subset of
$\operatorname*{Ends}\mathcal{P}$.}. But this contradicts the fact that
$\mathcal{P}$ is end-minimal. Thus, we have found a contradiction in Case 3.
\end{itemize}

So we have obtained a contradiction in each case. Thus, our assumption was
false. This shows that the path cover $\mathcal{P}^{\prime}$ is end-minimal.
As we already said above, this allows us to apply the induction hypothesis to
$D^{\prime}$ instead of $D$, and conclude that the end-minimal path cover
$\mathcal{P}^{\prime}$ of $D^{\prime}$ has an independent cross-cut. This
independent cross-cut is clearly an independent cross-cut of $\mathcal{P}$ as
well, and thus we have shown that $\mathcal{P}$ has an independent cross-cut.
This proves Claim 2.] \medskip

As explained above, this completes the proof of Theorem \ref{thm.gm.gm}.
\end{proof}

\subsubsection{\label{subsec.paths.gm.apps}Applications}

Here are two simple applications of the Gallai--Milgram theorem:

\begin{itemize}
\item Remember the Easy R\'{e}dei theorem (Theorem \ref{thm.tour.redei-easy}),
which we proved long ago. It says that each tournament has a Hamiltonian path.

We can now prove it again using the Gallai--Milgram theorem: \medskip

\textit{New proof of the Easy R\'{e}dei theorem:} Indeed, let $D$ be a
tournament. The Gallai--Milgram theorem shows that $D$ has a path cover with
an independent cross-cut\footnote{See the above proof of Theorem
\ref{thm.gm.gm} for the definition of a \textquotedblleft
cross-cut\textquotedblright.}. Consider this path cover and this cross-cut.
But since $D$ is a tournament, any independent set of $D$ has size $\leq1$.
Thus, our independent cross-cut must have size $\leq1$. Hence, our path cover
must consist of $1$ path only (because the size of the path cover equals the
size of its cross-cut). But this means that it is a Hamiltonian path (or, more
precisely, it consists of a single path, which is necessarily a Hamiltonian
path). Hence, $D$ has a Hamiltonian path. So we have proved the Easy R\'{e}dei
theorem (Theorem \ref{thm.tour.redei-easy}) again.

\item Less obviously, Hall's Marriage Theorem (Theorem \ref{thm.match.HMT})
and the Hall--K\"{o}nig matching theorem (Theorem \ref{thm.match.HKMT}) can be
proved again using Gallai--Milgram. Here is how: \medskip

\textit{New proof of the Hall--K\"{o}nig matching theorem:} Let $\left(
G,X,Y\right)  $ be a bipartite graph.

Let $D$ be the digraph obtained from $G$ by directing each edge so that it
goes from $X$ to $Y$ (in other words, each edge with endpoints $x\in X$ and
$y\in Y$ becomes an arc with source $x$ and target $y$). Thus, in the digraph
$D$, no vertex can simultaneously be the source of some arc and the target of
some arc. Thus, any path of $D$ has length $\leq1$. Here is an illustration of
a bipartite graph $\left(  G,X,Y\right)  $ (drawn as agreed in Example
\ref{exa.match.bipartite-lr}) and the corresponding digraph $D$:%
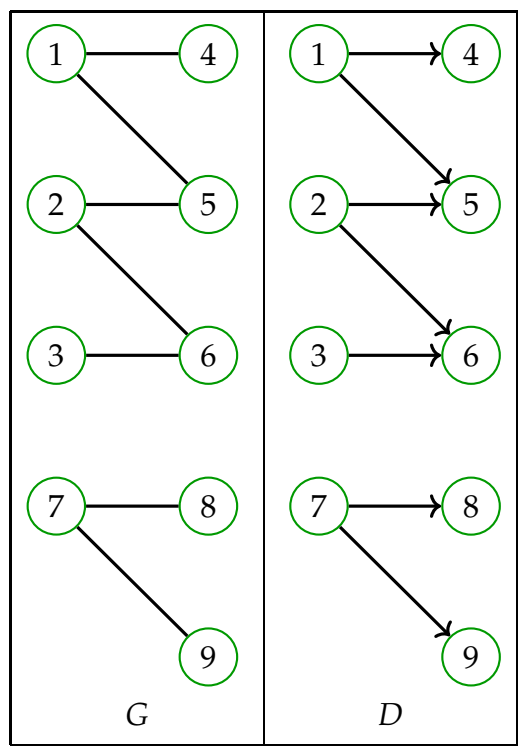
\[%

\ \ \ \ \ .
\]

As we said, any path of $D$ has length $\leq1$. Thus, any path of $D$
corresponds either to a vertex of $G$ or to an edge of $G$ (depending on
whether its length is $0$ or $1$). Hence, any path cover $\mathcal{P}$ of $D$
necessarily consists of length-$0$ paths (corresponding to vertices of $G$)
and length-$1$ paths (corresponding to edges of $G$); moreover, the edges of
$\mathcal{P}$ (that is, the edges corresponding to the length-$1$ paths in
$\mathcal{P}$) form a matching of $G$, and the vertices of $\mathcal{P}$ (that
is, the vertices corresponding to length-$0$ paths in $\mathcal{P}$) are
precisely the vertices that are not matched in this matching.

Now, Theorem \ref{thm.gm.gm} shows that there exist a path cover $\mathcal{P}$
of $D$ and an independent cross-cut $S$ of $\mathcal{P}$. Consider these
$\mathcal{P}$ and $S$. For the purpose of illustration, let us draw a path
cover $\mathcal{P}$ (by marking the arcs in red) and an independent cross-cut
$S$ of $\mathcal{P}$ (by drawing each vertex $s\in S$ as a blue diamond
instead of a green circle):%
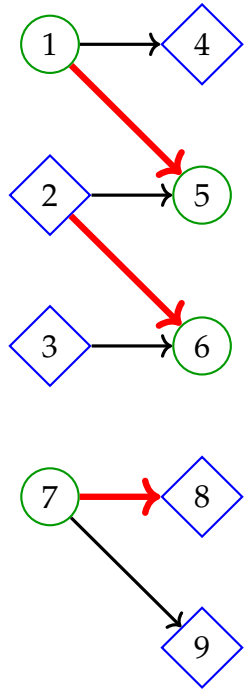
\[%
\begin{tikzpicture}[scale=2]
\begin{scope}[every node/.style={circle,thick,draw=green!60!black}]
\node(5) at (1, 0) {$5$};
\node(6) at (1, -1) {$6$};
\node(7) at (0, -2) {$7$};
\end{scope}
\begin{scope}[every node/.style={diamond,thick,draw=blue}]
\node(1) at (0, 1) {$1$};
\node(2) at (0, 0) {$2$};
\node(3) at (0, -1) {$3$};
\node(4) at (1, 1) {$4$};
\node(8) at (1, -2) {$8$};
\node(9) at (1, -3) {$9$};
\end{scope}
\begin{scope}[every edge/.style={draw=black,very thick}]
\path[->] (2) edge (5);
\path[->] (3) edge (6);
\path[->] (7) edge (9);
\end{scope}
\begin{scope}[every edge/.style={draw=red, line width=0.2pc}]
\path[->] (1) edge (5);
\path[->] (2) edge (6);
\path[->] (7) edge (8);
\end{scope}
\end{tikzpicture}%
\ \ .
\]

We have $\left\vert S\right\vert =\left\vert X\cap S\right\vert +\left\vert
Y\cap S\right\vert $ (since the set $S$ is the union of its two disjoint
subsets $X\cap S$ and $Y\cap S$).

The set $S$ is an independent set of the digraph $D$, thus also an independent
set of the graph $D^{\operatorname*{und}}=G$. From this, we easily obtain
$N\left(  X\cap S\right)  \subseteq Y\setminus S$ (since $\left(
G,X,Y\right)  $ is a bipartite graph)\footnote{\textit{Proof.} Let $v\in
N\left(  X\cap S\right)  $. Thus, $v$ is a vertex with a neighbor in $X\cap
S$. Let $x$ be this neighbor. Then, $x\in X\cap S\subseteq X$, so that the
vertex $v$ has a neighbor in $X$ (namely, $x$). Since $\left(  G,X,Y\right)  $
is a bipartite graph, this entails that $v\in Y$. Furthermore, we have $x\in
X\cap S\subseteq S$. If we had $v\in S$, then the set $S$ would contain two
adjacent vertices (namely, $v$ and $x$), which would contradict the fact that
$S$ is an independent set of $G$. Thus, we have $v\notin S$. Combining $v\in
Y$ with $v\notin S$, we obtain $v\in Y\setminus S$.
\par
Forget that we fixed $v$. We thus have shown that $v\in Y\setminus S$ for each
$v\in N\left(  X\cap S\right)  $. In other words, $N\left(  X\cap S\right)
\subseteq Y\setminus S$.}. Therefore, $\left\vert N\left(  X\cap S\right)
\right\vert \leq\left\vert Y\setminus S\right\vert $, so that $\left\vert
Y\setminus S\right\vert \geq\left\vert N\left(  X\cap S\right)  \right\vert $.
Hence,%
\begin{align}
\left\vert Y\right\vert  &  =\underbrace{\left\vert Y\setminus S\right\vert
}_{\geq\left\vert N\left(  X\cap S\right)  \right\vert }+\left\vert Y\cap
S\right\vert \geq\left\vert N\left(  X\cap S\right)  \right\vert
+\underbrace{\left\vert Y\cap S\right\vert }_{\substack{=\left\vert
S\right\vert -\left\vert X\cap S\right\vert \\\text{(since }\left\vert
S\right\vert =\left\vert X\cap S\right\vert +\left\vert Y\cap S\right\vert
\text{)}}}\nonumber\\
&  =\left\vert N\left(  X\cap S\right)  \right\vert +\left\vert S\right\vert
-\left\vert X\cap S\right\vert . \label{pf.HKMT.proof-using-gm.1}%
\end{align}

Now, let $M$ be the set of edges of $G$ corresponding to the length-$1$ paths
in our path cover $\mathcal{P}$. As we already mentioned, this set $M$ is a
matching of $G$ (since two paths in $\mathcal{P}$ cannot have a vertex in
common). The vertices that are not matched in $M$ are precisely the vertices
that don't belong to any of the length-$1$ paths in $\mathcal{P}$; in other
words, they are the vertices that belong to length-$0$ paths in $\mathcal{P}$
(since $\mathcal{P}$ is a path cover, and any path has length $\leq1$). We let
$p$ be the number of such vertices that lie in $X$, and we let $q$ be the
number of such vertices that lie in $Y$.

Thus, our path cover $\mathcal{P}$ contains exactly $p+q$ length-$0$ paths:
namely, $p$ length-$0$ paths consisting of a vertex in $X$ and $q$ length-$0$
paths consisting of a vertex in $Y$. Hence, the path cover $\mathcal{P}$
consists of $\left\vert M\right\vert +p+q$ paths in total (since it contains
$\left\vert M\right\vert $ many length-$1$ paths). The set $S$ contains
exactly one vertex from each of these $\left\vert M\right\vert +p+q$ paths
(since $S$ is a cross-cut of $\mathcal{P}$); therefore,
\begin{equation}
\left\vert S\right\vert =\left\vert M\right\vert +p+q.
\label{pf.HKMT.proof-using-gm.S=}%
\end{equation}

Each vertex $y\in Y$ that is matched in $M$ belongs to exactly one $M$-edge
(namely, to its $M$-edge), and conversely, each $M$-edge contains exactly one
vertex in $Y$ (which, of course, is matched in $M$). Thus, the map%
\begin{align*}
\left\{  \text{vertices in }Y\text{ that are matched in }M\right\}   &
\rightarrow M,\\
y  &  \mapsto\left(  \text{the }M\text{-edge of }y\right)
\end{align*}
is a bijection. Hence, the bijection principle yields%
\begin{equation}
\left(  \text{\# of vertices in }Y\text{ that are matched in }M\right)
=\left\vert M\right\vert . \label{pf.HKMT.proof-using-gm.2}%
\end{equation}

On the other hand, the set $Y$ contains exactly $q$ vertices that are not
matched in $M$ (by the definition of $q$). Therefore, $Y$ contains exactly
$\left\vert Y\right\vert -q$ vertices that are matched in $M$. In other words,%
\[
\left(  \text{\# of vertices in }Y\text{ that are matched in }M\right)
=\left\vert Y\right\vert -q.
\]

Comparing this with (\ref{pf.HKMT.proof-using-gm.2}), we obtain $\left\vert
M\right\vert =\left\vert Y\right\vert -q$. In other words,%
\begin{equation}
\left\vert M\right\vert +q=\left\vert Y\right\vert .
\label{pf.HKMT.proof-using-gm.3}%
\end{equation}
The same argument (but applied to $X$ and $p$ instead of $Y$ and $q$) yields%
\begin{equation}
\left\vert M\right\vert +p=\left\vert X\right\vert .
\label{pf.HKMT.proof-using-gm.4}%
\end{equation}
Hence, (\ref{pf.HKMT.proof-using-gm.S=}) becomes%
\begin{equation}
\left\vert S\right\vert =\underbrace{\left\vert M\right\vert +p}%
_{\substack{=\left\vert X\right\vert \\\text{(by
(\ref{pf.HKMT.proof-using-gm.4}))}}}+\,q=\left\vert X\right\vert +q.
\label{pf.HKMT.proof-using-gm.4b}%
\end{equation}

Now, from (\ref{pf.HKMT.proof-using-gm.3}), we obtain
\begin{align*}
\left\vert M\right\vert +q  &  =\left\vert Y\right\vert \\
&  \geq\left\vert N\left(  X\cap S\right)  \right\vert +\underbrace{\left\vert
S\right\vert }_{\substack{=\left\vert X\right\vert +q\\\text{(by
(\ref{pf.HKMT.proof-using-gm.4b}))}}}-\left\vert X\cap S\right\vert
\ \ \ \ \ \ \ \ \ \ \left(  \text{by (\ref{pf.HKMT.proof-using-gm.1})}\right)
\\
&  =\left\vert N\left(  X\cap S\right)  \right\vert +\left\vert X\right\vert
+q-\left\vert X\cap S\right\vert \\
&  =\left\vert N\left(  X\cap S\right)  \right\vert +\left\vert X\right\vert
-\left\vert X\cap S\right\vert +q.
\end{align*}
Cancelling $q$, we obtain%
\begin{equation}
\left\vert M\right\vert \geq\left\vert N\left(  X\cap S\right)  \right\vert
+\left\vert X\right\vert -\left\vert X\cap S\right\vert .
\label{pf.HKMT.proof-using-gm.5}%
\end{equation}
Thus, we have found a matching $M$ of $G$ and a subset $U$ of $X$ (namely,
$U=X\cap S$) such that $\left\vert M\right\vert \geq\left\vert N\left(
U\right)  \right\vert +\left\vert X\right\vert -\left\vert U\right\vert $.
This proves the Hall--K\"{o}nig matching theorem (once again). \medskip

\textit{New proof of Hall's Marriage Theorem:} Proceed as in the proof of the
Hall--K\"{o}nig matching theorem that we just gave. But now assume that our
bipartite graph $\left(  G,X,Y\right)  $ satisfies the Hall condition (i.e.,
we have $\left\vert N\left(  A\right)  \right\vert \geq\left\vert A\right\vert
$ for each subset $A$ of $X$). Hence, in particular, $\left\vert N\left(
X\cap S\right)  \right\vert \geq\left\vert X\cap S\right\vert $. Therefore,
(\ref{pf.HKMT.proof-using-gm.5}) becomes%
\[
\left\vert M\right\vert \geq\underbrace{\left\vert N\left(  X\cap S\right)
\right\vert }_{\geq\left\vert X\cap S\right\vert }+\left\vert X\right\vert
-\left\vert X\cap S\right\vert \geq\left\vert X\right\vert .
\]
Hence, Proposition \ref{prop.match.bipar-match-basics} \textbf{(e)} shows that
the matching $M$ is $X$-complete. Thus, $G$ has an $X$-complete matching
(namely, $M$). This proves Hall's Marriage Theorem (once again).
\end{itemize}

\begin{exercise}
\label{exe.9.4}Let $c$ and $r$ be two positive integers. Let $T$ be a
tournament with more than $r^{c}$ vertices. Each arc of $T$ is colored with
one of the $c$ colors $1,2,\ldots,c$. Prove that $T$ has a monochromatic path
of length $r$.

(A path is said to be \textbf{monochromatic} if all its arcs have the same
color.) \medskip

[\textbf{Hint:} Induct on $c$, and apply Gallai-Milgram to a certain digraph
in the induction step.]
\end{exercise}

\begin{remark}
If we apply Exercise \ref{exe.9.4} to $c=1$, then we recover the easy
R\'{e}dei theorem (Theorem \ref{thm.tour.redei-easy}). Indeed, if $T$ is any
tournament, then we can color all its arcs with the color $1$, and then use
Exercise \ref{exe.9.4} (applied to $c=1$ and $r=\left\vert \operatorname*{V}%
\left(  T\right)  \right\vert -1$) to conclude that $T$ has a monochromatic
path of length $\left\vert \operatorname*{V}\left(  T\right)  \right\vert -1$.
But such a path must necessarily be a Hamiltonian path (since its length
forces it to contain all vertices of $T$).
\end{remark}

\subsection{Path-missing sets}

We move on to less well-trodden ground.

Menger's theorem (one of the many) is from 1927; the Gallai--Milgram theorem
is from 1960. One might think that everything that can be said about paths in
graphs has been said long ago.

Apparently, this is not the case. In 2017, when trying to come up with a
homework exercise for a previous iteration of this course, I was experimenting
with paths in Python. Specifically, I was looking at digraphs $D=\left(
V,A,\psi\right)  $ with two distinct vertices $s$ and $t$ selected. Inspired
by the arc-Menger theorems, I was looking at the subsets $B$ of $A$ that could
be removed without disconnecting $s$ from $t$ (more precisely, without
destroying all paths from $s$ to $t$). I noticed that the number of such
subsets $B$ seemed to be even whenever $D$ has a cycle or a \textquotedblleft
useless arc\textquotedblright\ (i.e., an arc that is used by no path from $s$
to $t$)\ \ \ \footnote{With one exception: If $A=\varnothing$, then it is
odd.}, and odd otherwise.

I could not prove this observation. Soon after, Joel Brewster Lewis and Lukas
Katth\"{a}n came up with a proof and multiple stronger results. The proofs can
now be found in a joint preprint \cite{GrKaLe21}.

The first way to strengthen the observation is to replace the parity claim
(i.e., the claim that the number is even or odd depending on cycles and
useless arcs) by a stronger claim about an alternating sum. This is an
instance of a general phenomenon, in which a statement of the form
\textquotedblleft the number of some class of things is even\textquotedblright%
\ can often be replaced by a stronger statement of the form \textquotedblleft
we can assign a plus or minus sign to each of these things, and then the total
number of plus signs equals the total number of minus signs\textquotedblright.
The stronger statement is as follows:

\begin{theorem}
[Grinberg--Lewis--Katth\"{a}n]\label{thm.grleka.altsum}Let $D=\left(
V,A,\psi\right)  $ be a multidigraph. Let $s$ and $t$ be two distinct vertices
of $D$. A subset $B$ of $A$ will be called \textbf{path-missing} if $D$ has a
path from $s$ to $t$ that does not use any of the arcs in $B$ (that is, a path
from $s$ to $t$ that would not be destroyed if we remove all arcs in $B$ from
$D$). (In the terminology of Definition \ref{def.menger.arc-sep}, this is the
same as saying that $B$ is \textbf{not} an $s$-$t$-arc-separator.)

Let $\mathbf{M}$ be the set of all path-missing subsets of $A$.

\begin{enumerate}
\item[\textbf{(a)}] If $D$ has an arc that is not used by any path from $s$ to
$t$ (this is what we call a \textquotedblleft useless arc\textquotedblright),
then%
\[
\sum_{B\in\mathbf{M}}\left(  -1\right)  ^{\left\vert B\right\vert }=0
\]
(and thus $\left\vert \mathbf{M}\right\vert $ is even).

\item[\textbf{(b)}] If $D$ has a cycle, then%
\[
\sum_{B\in\mathbf{M}}\left(  -1\right)  ^{\left\vert B\right\vert }=0
\]
(and thus $\left\vert \mathbf{M}\right\vert $ is even).

\item[\textbf{(c)}] If $A=\varnothing$, then%
\[
\sum_{B\in\mathbf{M}}\left(  -1\right)  ^{\left\vert B\right\vert }=0
\]
(and thus $\left\vert \mathbf{M}\right\vert $ is even).

\item[\textbf{(d)}] In all other cases, we have%
\[
\sum_{B\in\mathbf{M}}\left(  -1\right)  ^{\left\vert B\right\vert }=\left(
-1\right)  ^{\left\vert A\right\vert -\left\vert V^{\prime}\right\vert },
\]
where $V^{\prime}$ is the set of all vertices of $D$ that have outdegree $>0$
(and thus $\left\vert \mathbf{M}\right\vert $ is odd).
\end{enumerate}
\end{theorem}

\begin{example}
Let $D=\left(  V,A,\psi\right)  $ be the following digraph:%
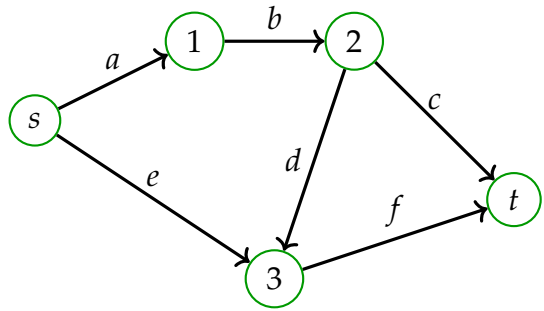
\[%
\begin{tikzpicture}[scale=2.1]
\begin{scope}[every node/.style={circle,thick,draw=green!60!black}]
\node(s) at (-1, 0.5) {$s$};
\node(t) at (2, 0) {$t$};
\node(1) at (0, 1) {$1$};
\node(2) at (1, 1) {$2$};
\node(3) at (0, -0.2) {$3$};
\end{scope}
\begin{scope}[every edge/.style={draw=black,very thick}, every loop/.style={}]
\path[->] (s) edge node[above] {$a$} (1);
\path[->] (s) edge node[above] {$e$} (3);
\path[->] (2) edge node[left=1pt] {$d$} (3);
\path[->] (2) edge node[above] {$c$} (t);
\path[->] (3) edge node[above] {$f$} (t);
\path[->] (1) edge node[above] {$b$} (2);
\end{scope}
\end{tikzpicture}%
\ \ .
\]
Let $s$ and $t$ be the vertices labelled $s$ and $t$ here. Then, $D$ has
neither a cycle nor a \textquotedblleft useless arc\textquotedblright, and its
arc set $A$ is nonempty; thus, Theorem \ref{thm.grleka.altsum} \textbf{(d)}
applies. The path-missing subsets of $A$ are the three sets $\left\{
a,b,c,d\right\}  $, $\left\{  c,e\right\}  $ and $\left\{  d,e,f\right\}  $ as
well as all their subsets (such as $\left\{  b,c,d\right\}  $). In other
words,%
\begin{align*}
\mathbf{M}  &  =\left\{  \text{all subsets of }\left\{  a,b,c,d\right\}
\right\}  \cup\left\{  \text{all subsets of }\left\{  c,e\right\}  \right\} \\
&  \ \ \ \ \ \ \ \ \ \ \cup\left\{  \text{all subsets of }\left\{
d,e,f\right\}  \right\} \\
&  =\{\varnothing,\ \left\{  a\right\}  ,\ \left\{  b\right\}  ,\ \left\{
c\right\}  ,\ \left\{  d\right\}  ,\ \left\{  a,b\right\}  ,\ \left\{
a,c\right\}  ,\ \left\{  a,d\right\}  ,\ \left\{  b,c\right\}  ,\ \left\{
b,d\right\}  ,\\
&  \ \ \ \ \ \ \ \ \ \ \left\{  c,d\right\}  ,\ \left\{  a,b,c\right\}
,\ \left\{  a,b,d\right\}  ,\ \left\{  a,c,d\right\}  ,\ \left\{
b,c,d\right\}  ,\ \left\{  a,b,c,d\right\}  ,\\
&  \ \ \ \ \ \ \ \ \ \ \left\{  e\right\}  ,\ \left\{  c,e\right\}
,\ \left\{  f\right\}  ,\ \left\{  d,e\right\}  ,\ \left\{  e,f\right\}
,\ \left\{  d,f\right\}  ,\ \left\{  d,e,f\right\}  \}.
\end{align*}
Hence, the sum $\sum_{B\in\mathbf{M}}\left(  -1\right)  ^{\left\vert
B\right\vert }$ has $11$ addends equal to $-1$ and $12$ addends equal to $1$;
thus, this sum equals to $1$. This is precisely the value $\left(  -1\right)
^{\left\vert A\right\vert -\left\vert V^{\prime}\right\vert }=\left(
-1\right)  ^{6-4}=1$ predicted by Theorem \ref{thm.grleka.altsum} \textbf{(d)}.
\end{example}

\begin{proof}
[Proof of Theorem \ref{thm.grleka.altsum}.]See \cite[Theorem 1.3]{GrKaLe21}
(where $\mathbf{M}$ is denoted by $\operatorname*{PM}\left(  D\right)  $, and
where arcs are called \textquotedblleft edges\textquotedblright). Of course,
part \textbf{(c)} is obvious, and part \textbf{(a)} is easy (since inserting a
useless arc into a set $B\in\mathbf{M}$ or removing it from a set
$B\in\mathbf{M}$ always results in a set in $\mathbf{M}$). Parts \textbf{(b)}
and \textbf{(d)} are the interesting ones. The proof in \cite[Theorem
1.3]{GrKaLe21} relies on a recursive argument (\textquotedblleft
deletion-contraction\textquotedblright) in which we pick an arc with source
$s$ and consider the two smaller digraphs $D\setminus a$ and $D\diagup a$
obtained (respectively) by deleting the arc $a$ from $D$ and by
\textquotedblleft contracting\textquotedblright\ $a$ \textquotedblleft to a
point\textquotedblright.
\end{proof}

Further levels of strength can be reached by treating $\mathbf{M}$ as a
topological space. Indeed, $\mathbf{M}$ is not just a random collection of
sets of arcs, but actually a
\textbf{\href{https://en.wikipedia.org/wiki/Simplicial_complex}{\textbf{simplicial
complex}}} (since any subset of a path-missing subset of $A$ is again
path-missing). Simplicial complexes are known to be a combinatorial model for
topological spaces, and in particular they have homology groups, homotopy
types, etc.. Thus, in particular, we can ask ourselves how the topological
space corresponding to the simplicial complex $\mathbf{M}$ looks like. This,
too, has been answered in \cite[Theorem 1.3]{GrKaLe21}: It is homotopic to a
sphere or a ball (depending on the existence of cycles or \textquotedblleft
useless arcs\textquotedblright); its dimension can also be determined
explicitly. (The sum $\sum_{B\in\mathbf{M}}\left(  -1\right)  ^{\left\vert
B\right\vert }$ discussed above is, of course, its reduced Euler
characteristic.) Even stronger statements can be made using \textbf{discrete
Morse theory} (see again \cite{GrKaLe21} for details).

\subsection{Elser's sums}

We now return to undirected (multi)graphs. Here is a result found by Veit
Elser in 1984 (\cite[Lemma 1]{Elser84}), as a lemma for his work in
statistical mechanics:\footnote{I have restated the result beyond recognition;
see \cite[Remark 1.4]{Grinbe21} for why Theorem \ref{thm.elser.1} actually
implies \cite[Lemma 1]{Elser84}.}

\begin{theorem}
[Elser's theorem, in my version]\label{thm.elser.1}Let $G=\left(
V,E,\varphi\right)  $ be a multigraph with at least one edge. Fix a vertex
$v\in V$.

If $F\subseteq E$, then an $F$\textbf{-path} shall mean a path of $G$ such
that all edges of this path belong to $F$. In other words, it means a path of
the spanning subgraph $\left(  V,F,\varphi\mid_{F}\right)  $.

If $e\in E$ is an edge and $F\subseteq E$ is a subset, then we say that $F$
\textbf{infects} $e$ if there exists an $F$-path from $v$ to some endpoint of
$e$. (The terminology is inspired by the idea that some infectious disease
starts at $v$ and spreads along the $F$-edges.)

(Note that if an edge $e$ contains the vertex $v$, then \textbf{any} subset
$F$ of $E$ (even the empty set) infects $e$, because $\left(  v\right)  $ is a
trivial $F$-path from $v$ to $v$.)

Then,%
\[
\sum_{\substack{F\subseteq E\text{ infects}\\\text{every edge }e\in E}}\left(
-1\right)  ^{\left\vert F\right\vert }=0.
\]

\end{theorem}

\begin{example}
\label{exa.elser.1}Let $G=\left(  V,E,\varphi\right)  $ be the following
graph:%
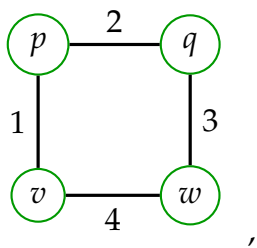
\[%
\begin{tikzpicture}[scale=2]
\begin{scope}[every node/.style={circle,thick,draw=green!60!black}]
\node(1) at (0, 0) {$v$};
\node(2) at (0, 1) {$p$};
\node(3) at (1, 1) {$q$};
\node(4) at (1, 0) {$w$};
\end{scope}
\begin{scope}[every edge/.style={draw=black,very thick}, every loop/.style={}]
\path[-] (1) edge node[left] {$1$} (2) (2) edge node[above] {$2$} (3);
\path[-] (3) edge node[right] {$3$} (4) (4) edge node[below] {$4$} (1);
\end{scope}
\end{tikzpicture}
\ \ ,
\]
and let $v$ be the vertex labelled $v$. Then, the subsets of $E$ that infect
every edge are%
\[
\left\{  1,2\right\}  ,\ \ \left\{  1,4\right\}  ,\ \ \left\{  3,4\right\}
,\ \ \left\{  1,2,3\right\}  ,\ \ \left\{  1,3,4\right\}  ,\ \ \left\{
1,2,4\right\}  ,\ \ \left\{  2,3,4\right\}  ,\ \ \left\{  1,2,3,4\right\}  .
\]
Thus,%
\begin{align*}
&  \sum_{\substack{F\subseteq E\text{ infects}\\\text{every edge }e\in
E}}\left(  -1\right)  ^{\left\vert F\right\vert }\\
&  =\left(  -1\right)  ^{2}+\left(  -1\right)  ^{2}+\left(  -1\right)
^{2}+\left(  -1\right)  ^{3}+\left(  -1\right)  ^{3}+\left(  -1\right)
^{3}+\left(  -1\right)  ^{3}+\left(  -1\right)  ^{4}\\
&  =0,
\end{align*}
exactly as predicted by Theorem \ref{thm.elser.1}.
\end{example}

\begin{remark}
It might appear more natural to study subsets $F\subseteq E$ infecting
vertices rather than edges. However, Theorem \ref{thm.elser.1} would be false
if we replaced \textquotedblleft every edge $e\in E$\textquotedblright\ by
\textquotedblleft every vertex $v\in V$\textquotedblright. The graph in
Example \ref{exa.elser.1} provides a counterexample.

However, if we go further and replace $F\subseteq E$ by $W\subseteq V$, then
we get something true again -- see Theorem \ref{thm.elser.2} below.
\end{remark}

\begin{proof}
[Proof of Theorem \ref{thm.elser.1}.]Elser's proof is somewhat complicated. I
give a different proof in \cite[Theorem 1.2]{Grinbe21}, which is elementary
and nice if I may say so myself.

My proof should also be not very hard to discover, once you have the following
hint: It suffices to prove the equality%
\[
\sum_{\substack{F\subseteq E\text{ does \textbf{not} infect}\\\text{every edge
}e\in E}}\left(  -1\right)  ^{\left\vert F\right\vert }=0
\]
(because the total sum $\sum_{F\subseteq E}\left(  -1\right)  ^{\left\vert
F\right\vert }$ is known to be $0$). In order to prove this equality, we equip
the set $E$ with some total order (it doesn't matter how; we can just rank the
edges arbitrarily), and we make the following definition: If $F\subseteq E$ is
a subset that does \textbf{not} infect every edge $e\in E$, then we let
$\varepsilon\left(  F\right)  $ be the smallest (with respect to our chosen
total order) edge that is not infected by $F$. Now, you can show that if
$F\subseteq E$ is a subset that does \textbf{not} infect every edge $e\in E$,
then the set\footnote{The symbol $\bigtriangleup$ stands for the symmetric
difference of two sets. Recall its definition: If $X$ and $Y$ are two sets,
then their \textbf{symmetric difference} $X\bigtriangleup Y$ is defined to be
the set
\begin{align*}
\left(  X\cup Y\right)  \setminus\left(  X\cap Y\right)   &  =\left(
X\setminus Y\right)  \cup\left(  Y\setminus X\right) \\
&  =\left\{  \text{all elements that belong to \textbf{exactly} one of
}X\text{ and }Y\right\}  .
\end{align*}
} $F^{\prime}:=F\bigtriangleup\left\{  \varepsilon\left(  F\right)  \right\}
$ (that is, the set obtained from $F$ by inserting $\varepsilon\left(
F\right)  $ if $\varepsilon\left(  F\right)  \notin F$ and by removing
$\varepsilon\left(  F\right)  $ if $\varepsilon\left(  F\right)  \in F$) has
the same property (viz., it does not infect every edge $e\in E$) and satisfies
$\varepsilon\left(  F^{\prime}\right)  =\varepsilon\left(  F\right)  $. This
entails that the addends in the sum $\sum_{\substack{F\subseteq E\text{ does
\textbf{not} infect}\\\text{every edge }e\in E}}\left(  -1\right)
^{\left\vert F\right\vert }$ cancel each other in pairs (namely, the addend
for a given set $F$ cancels the addend for the set $F^{\prime}=F\bigtriangleup
\left\{  \varepsilon\left(  F\right)  \right\}  $), and thus the whole sum is
$0$.
\end{proof}

Elser's theorem, too, can be generalized and strengthened. The strengthening
is similar to what we did with Theorem \ref{thm.grleka.altsum}: We treat the
set of all \textquotedblleft non-pandemic-causing subsets\textquotedblright%
\ (i.e., of all subsets $F\subseteq E$ that \textbf{don't} infect every edge)
as a simplicial complex (since a subset of a non-pandemic-causing subset is
again non-pandemic-causing), and analyze this complex as a topological space.
The claim of Theorem \ref{thm.elser.1} then says that the reduced Euler
characteristic of this space is $0$; but we can actually show that this space
is contractible (i.e., homotopy-equivalent to a point). Even better, we can
prove that the simplicial complex of all non-pandemic-causing subsets is
\textbf{collapsible} (a combinatorial property that is stronger than
contractibility of the corresponding space). See \cite[\S 5]{Grinbe21} for
definitions and proofs.

We can furthermore generalize the theorem. One way to do so is to replace our
\textquotedblleft patient zero\textquotedblright\ $v$ by a set of vertices.
This leads to a much less trivial situation. The recent paper \cite{DHLetc19}
by Dorpalen-Barry, Hettle, Livingston, Martin, Nasr, Vega and Whitlatch proves
some results and asks some questions (that are still open as of 2022).

A different direction in which Elser's theorem can be generalized is more
fundamental: It turns out that the theorem is not really about graphs and
edges. Instead, there is a general structure that I call a \textquotedblleft
shade map\textquotedblright, which always leads to a certain sum being $0$.
See \cite[\S 4]{Grinbe21} for the details of this generalization. I will not
explain it here, but I will state one more particular case of it
(\cite[Theorem 3.2]{Grinbe21}), which replaces edges by vertices throughout
Theorem \ref{thm.elser.1}:

\begin{theorem}
[vertex-Elser's theorem]\label{thm.elser.2}Let $G=\left(  V,E,\varphi\right)
$ be a multigraph with at least two vertices. Fix a vertex $v\in V$.

If $W\subseteq V$, then a $W$\textbf{-vertex-path} shall mean a path
$\mathbf{p}$ such that all intermediate vertices of $\mathbf{p}$ belong to
$W$. (Recall that the \textquotedblleft intermediate vertices of $\mathbf{p}%
$\textquotedblright\ mean all vertices of $\mathbf{p}$ except for the starting
and ending points of $\mathbf{p}$.) (Note that any path of length $\leq1$ is
automatically a $W$-vertex-path, since it has no intermediate vertices.)

If $w\in V\setminus\left\{  v\right\}  $ is any vertex, and $W\subseteq
V\setminus\left\{  v\right\}  $ is any subset, then we say that $W$
\textbf{vertex-infects} $w$ if there exists a $W$-vertex-path from $v$ to $w$.
(This is always true when $w$ is a neighbor of $v$.)

Then,%
\[
\sum_{\substack{W\subseteq V\setminus\left\{  v\right\}  \text{
vertex-infects}\\\text{every vertex }w\in V\setminus\left\{  v\right\}
}}\left(  -1\right)  ^{\left\vert W\right\vert }=0.
\]

\end{theorem}

\end{document}